\documentclass[12pt,leqno]{article}

\usepackage{amsmath,amssymb,amsfonts,theorem,latexsym,epsfig}
\usepackage{psfrag}
\usepackage{color}
\usepackage{float}
\usepackage{setspace}
\usepackage{relsize}

\usepackage{indentfirst}
\usepackage{xspace}
\usepackage{multirow}
\usepackage{hyperref}
\usepackage{xcolor}
\hypersetup{colorlinks=true,urlcolor=blue,linkcolor=blue,citecolor=blue}

\newtheorem{thm}{Theorem}
\newtheorem{theo}{Theorem}[section]
\newtheorem{lem}[theo]{Lemma}

\newtheorem{prop}[theo]{Proposition}
\newtheorem{defi}[theo]{Definition}
\newtheorem{claim}[theo]{Claim}

{\theorembodyfont{\rmfamily} \newtheorem{rem}[theo]{Remark}}
{\theorembodyfont{\rmfamily} \newtheorem{ex}[theo]{Example}}
{\theorembodyfont{\rmfamily} }
{\theorembodyfont{\rmfamily} }
{\theorembodyfont{\rmfamily} \newtheorem{problem}[theo]{Problem}}

\numberwithin{equation}{section}

\usepackage[all]{xy}
\usepackage{multirow}

\newcommand{\Appendix}[1]{%
  \refstepcounter{section}%
  \addtocontents{toc}{\protect\setcounter{tocdepth}{1}}
  \addcontentsline{toc}{section}%
    {\bfseries\appendixname~\thesection:\ #1}%
    {\medskip\noindent \Large\bfseries\appendixname\ \thesection:\ #1}%
\sectionmark{#1}\smallskip\noindent
\renewcommand{\theequation}{{\bf 
{{\thesection}}.{\arabic{equation}}}}
}
\DeclareFontFamily{U}{rsf}{}
\DeclareFontShape{U}{rsf}{m}{n}{
  <5> <6> rsfs5 <7> <8> <9> rsfs7 <10->  rsfs10}{}
\DeclareMathAlphabet{\mathscr}{U}{rsf}{m}{n}
\newcommand{\mycal}[1]{\mathscr{#1}}
\usepackage{mathtools}

\DeclarePairedDelimiter\fl{\lfloor}{\rfloor}

\newcounter{num}

\newcommand{\BIGOP}[1]{\mathop{\mathchoice%
{\raise-0.22em\hbox{\huge $#1$}}%
{\raise-0.05em\hbox{\Large $#1$}}{\hbox{\large $#1$}}{#1}}}
\newcommand{\bigtimes}{\BIGOP{\times}}

\newcommand{\op}[1]{\operatorname{#1}}
\newcommand{\lan}[1]{{}^{L}{#1}}

\newcommand{\Higgs}{\op{{\bf Higgs}}}
\newcommand{\gHiggs}{\op{\boldsymbol{\mathcal{H}iggs}}}

\newcommand{\gLoc}{\op{\boldsymbol{\mathcal{F}lat}}}
\newcommand{\Heck}{\op{{\bf Hecke}}}
\newcommand{\sHeck}{\op{{\boldsymbol{\mathcal{H}}{\bf ecke}}}}
\newcommand{\sH}{\op{\boldsymbol{\mathcal{H}}}}

\newcommand{\bD}{\boldsymbol{\sf{D}}}

\newcommand{\chr}{\op{{\sf{char}}}}
\newcommand{\cchr}{\op{{\sf{cochar}}}}
\newcommand{\aj}{\mathfrak{aj}}

\newcommand{\Bun}{\op{{\bf Bun}}}
\newcommand{\sBun}{\op{\boldsymbol{\mathcal{B}un}}}
\newcommand{\sN}{\op{\boldsymbol{\mathcal{N}}}}
\newcommand{\Pic}{\op{{\bf Pic}}}
\newcommand{\gPic}{\op{\boldsymbol{\mathcal{P}ic}}}
\newcommand{\bh}{\boldsymbol{h}}
\newcommand{\bH}{\boldsymbol{H}}

\newcommand{\bb}{\boldsymbol{\sf{b}}}

\newcommand{\ba}{\boldsymbol{\sf{a}}}
\newcommand{\bF}{\boldsymbol{\sf{F}}}
\newcommand{\bG}{\boldsymbol{\sf{G}}}

\newcommand{\pardeg}{\op{{\sf{pardeg}}}}
\newcommand{\parch}{\op{{\sf{parch}}}}

\newcommand{\kk}{{\sf{k}}}
\newcommand{\be}{\boldsymbol{\sf{e}}}
\newcommand{\bT}{\boldsymbol{\sf{T}}}
\newcommand{\bt}{\boldsymbol{\sf{t}}}
\newcommand{\bws}{\boldsymbol{\sf{ws}}}
\newcommand{\bps}{\boldsymbol{\sf{ps}}}
\newcommand{\pt}{\boldsymbol{\sf{pt}}}
\newcommand{\bdelta}{\boldsymbol{\delta}}
\newcommand{\bepsilon}{\boldsymbol{\varepsilon}}
\newcommand{\bzeta}{\boldsymbol{\zeta}}
\newcommand{\bp}{\boldsymbol{\sf{p}}}
\newcommand{\br}{\boldsymbol{\sf{r}}}

\newcommand{\blambda}{\boldsymbol{\lambda}}
\newcommand{\bLambda}{\boldsymbol{\Lambda}}
\newcommand{\beps}{\boldsymbol{\epsilon}}
\newcommand{\sq}{\boldsymbol{\sf{sq}}}
\newcommand{\bP}{\boldsymbol{\sf{P}}}
\newcommand{\bbH}{\boldsymbol{\sf{H}}}
\newcommand{\balpha}{\boldsymbol{\alpha}}
\newcommand{\bpi}{\boldsymbol{\pi}}
\newcommand{\bB}{\boldsymbol{B}}
\newcommand{\sC}{\widetilde{C}}
\newcommand{\sP}{\widetilde{P}}
\newcommand{\bff}{\boldsymbol{f}}
\newcommand{\bmu}{\boldsymbol{\mu}}
\newcommand{\oddL}{\boldsymbol{{\sf Odd}}}
\newcommand{\evenL}{\boldsymbol{{\sf Even}}}
\newcommand{\btheta}{\boldsymbol{\theta}}
\newcommand{\Par}{\op{\sf{Par}}}
\newcommand{\res}{\op{{\mathsf{res}}}}
\newcommand{\bWob}{\op{\sf{Wob}}}
\newcommand{\lhs}{\op{\sf{LHS}}}
\newcommand{\rhs}{\op{\sf{RHS}}}
\newcommand{\Dol}{{\op{DOL}^{\op{par}}}}
\newcommand{\hor}{\op{{\sf hor}}}
\newcommand{\ver}{\op{{\sf ver}}}
\newcommand{\lift}{\mathbf{\sf up}}
\newcommand{\btau}{\boldsymbol{\tau}}
\newcommand{\abH}{\widetilde{H}}
\newcommand{\abp}{\tilde{p}}
\newcommand{\abq}{\tilde{q}}
\newcommand{\Jac}{\mathbf{\mathsf J}}
\newcommand{\AJ}{\mathsf{aj}}
\newcommand{\bQ}{\mathfrak{M}}
\newcommand{\bbF}{\mathfrak{F}}
\newcommand{\verpsi}{\psi_{q}}
\newcommand{\mybeta}{{\alpha_{q}}}
\newcommand{\hatYpi}{\widehat{\mathsf{Yp}}_{i}}
\newcommand{\hatECI}{\widehat{\mathsf{EC}}_{I}}
\newcommand{\hatGCI}{\widehat{\mathsf{GC}}_{I}}
\newcommand{\hatAJ}{\mathsf{AJ}}
\newcommand{\bI}{{\mathbb I}}
\newcommand{\one}{\textsf{\textbf{1}}}
\newcommand{\LL}{M_{T,T}}
\newcommand{\OO}{{\mathcal O}}
\newcommand{\cube}{\textsf{\textbf{C}}}
\newcommand{\sumij}{\sum\limits_{\substack{I,J \\ \bI_{IJ} = 1}}}
\newcommand{\bR}{{\mathbb R}}
\newcommand{\bZ}{{\mathbb Z}}

\newcommand{\weights}{\op{{\sf weights}}}
\newcommand{\Aom}{\op{{\sf Oka}}}
\newcommand{\shriek}{!}
\newcommand{\bba}{\boldsymbol{\mathfrak{e}}}
\newcommand{\bbb}{\boldsymbol{\mathfrak{d}}}
\newcommand{\bvartheta}{\boldsymbol{\vartheta}}
\newcommand{\bFF}{\boldsymbol{\mathsf{f}}}
\newcommand{\ubf}{\underline{\bFF}}
\newcommand{\sM}{\boldsymbol{\mathcal{M}}}
\newcommand{\bxi}{\boldsymbol{\xi}}
\newcommand{\unone}{\mathbb{J}}
\newcommand{\incdelta}{\mathbb{D}}
\newcommand{\lat}{\mathsf{Q}}
\newcommand{\bimath}{\boldsymbol{\mathlarger{\mathlarger{\imath}}}}
\newcommand{\acan}{\boldsymbol{\mathfrak{u}}}
\newcommand{\bfri}{\boldsymbol{\mathfrak{i}}}

\usepackage{multind}

\makeindex{terms}
\makeindex{notations}

\begin{document}

\title{\textbf{Parabolic Hecke eigensheaves}}
\author{R. Donagi \and T. Pantev}
\date{}
\maketitle

\begin{abstract}
  We study the Geometric Langlands Conjecture (GLC) for rank two flat
  bundles on the projective line $C$ with tame ramification at five
  points $\{p_{1}, p_{2}, p_{3}, p_{4}, p_{5} \}$.  In particular we
  construct the automorphic $\mathcal{D}$-modules predicted by GLC on
  the moduli space of rank two parabolic bundles on $(C, \{p_{1},
  p_{2}, p_{3}, p_{4}, p_{5} \})$. The construction uses non-abelian
  Hodge theory and a Fourier-Mukai transform along the fibers of the
  Hitchin fibration to reduce the problem to one in classical
  projective geometry on the intersection of two quadrics in
  $\mathbb{P}^{4}$.
\end{abstract}

\

\begin{center}
{\LARGE\textbf{Faisceux propre Hecke paraboliques}}

\

{\large R. Donagi \qquad \qquad  T. Pantev}
\end{center}

\renewcommand{\abstractname}{R\'{e}sum\'{e}}

\begin{abstract}
Nous \'{e}tudions la conjecture g\'{e}om\'{e}trique de Langlands (CGL)
pour les fibr\'{e}s plats de rang deux sur la ligne projective $C$
avec une ramification mod\'{e}r\'{e}e en cinq points $\{p_{1}, p_{2},
p_{3}, p_{4}, p_{5} \}$. En particulier, nous construisons les
$\mathcal{D}$-modules automorphes pr\'{e}dits par CGL sur l'espace des
modules de fibr\'{e}s paraboliques de rang deux sur \linebreak $(C, \{p_{1},
p_{2}, p_{3}, p_{4}, p_{5} \})$. La construction utilise la
th\'{e}orie de Hodge non ab\'{e}lienne et une transform\'{e} de
Fourier-Mukai le long des fibres de la fibration de Hitchin pour
r\'{e}duire le probl\`{e}me \`{a} une question en g\'{e}om\'{e}trie
projective classique \`{a} propos l'intersection de deux quadriques
dans $\mathbb{P}^{4}$.
 \end{abstract}

\newpage


\smallskip


\smallskip

\noindent
{\bfseries Mots-cl\'{e}s:} th\'{e}orie de Hodge
non-ab\'{e}lienne, fibration de Hitchin, $\mathcal{D}$-modules, fibr\'{e}s
\linebreak paraboliques, correspondance de Langlands g\'{e}om\'{e}trique,
propri\'{e}t\'{e} de Hecke, rev\^{e}tement spectrale,
ab\'{e}lianisation

\

\newpage

\tableofcontents

\

\newpage

\section{Introduction} \label{sec:introduction}
%\addcontentsline{toc}{section}{Introduction}

The goal of this work is to provide proof-of-concept for the authors'
approach \cite{dp-jdg, dp-Langlands} to the Geometric Langlands
Conjecture (GLC) via abelianization and non-abelian Hodge theory
(NAHT). We do this by carrying out all details in one non trivial
case: rank $2$ bundles on $C=\mathbb{P}^1$ with parabolic structure at
$5$ points. In this introduction we explain the GLC, its parabolic
version, and our approach to solving them, and we try to convince the
reader that the special case we consider here already exhibits many of
the features that should occur in the full version. 

\index{terms}{Geometric Langlands Conjecture}
\index{terms}{GLC|see{Geometric Langlands Conjecture}}
\index{terms}{Non-abelian Hodge theory}
\index{terms}{NAHT|see{Non-abelian Hodge theory}}
\index{notations}{C@$C$}

\subsection{The conjecture}

The input for the Geometric Langlands Conjecture consists of a smooth
compact Riemann surface $C$ and a pair of Langlands dual complex
reductive groups $G$, $\lan{G}$
\index{notations}{G@$G$}
\index{notations}{LG@$\lan{G}$}.
With this data we can associate
moduli stacks:
\begin{description}
\item[$\sBun$, $\lan{\sBun}$:] the moduli stacks of principal $G$,
  $\lan{G}$ bundles $V$ on $C$.
  \index{terms}{principal!bundle}
  \index{terms}{moduli stack!of principal bundles}
  \index{terms}{moduli stack!of flat  bundles}
  \index{terms}{moduli stack!of Higgs bundles}
  \index{notations}{Buns@$\sBun$}
  \index{notations}{LBuns@$\lan{\sBun}$}
\item[$\gLoc$, $\lan{\gLoc}$:] the moduli stacks of $G$,
  $\lan{G}$ flat bundles{\footnote{unfortunately some references, including our own
    \cite{dp-Langlands}, call this moduli stack
    $\boldsymbol{\mathcal{L}oc}$. We prefer to reserve the term local
    systems \index{terms}{local system}
    \index{notations}{Locs@$\boldsymbol{\mathcal{L}oc}$}
    \index{terms}{Betti!version} for the Betti incarnation,
    parametrizing representations of the fundamental group, rather
    than the de Rham incarnation involving bundles with flat
    connections.}} $\mathbb{V} = (V,\nabla)$ on $C$.
  \index{terms}{principal!flat bundle}
  \index{terms}{flat!bundle} \index{terms}{flat!connection}
\index{notations}{Vbb@$\mathbb{V}$}
\index{notations}{Vnabla@$(V,\nabla)$}
\index{notations}{Flats@$\gLoc$}
\index{notations}{LFlats@$\lan{\gLoc}$}
\end{description}

The GLC predicts the existence of a canonical equivalence of derived
categories, relating coherent $\mathcal{O}$-modules on $\lan{\gLoc}$ to
coherent $\mathcal{D}$-modules on $\sBun$:
\index{terms}{Dmodule@$\mathcal{D}$-module}
\index{notations}{LFlats@$\lan{\gLoc}$}
  \index{notations}{Buns@$\sBun$}
\[
\tag{{\bfseries\sf GLC}} \label{eq:glc}
\boxed{
\mathfrak{c} : D_{\op{coh}}(\lan{\gLoc},\mathcal{O})
\stackrel{\cong}{\longrightarrow}
D_{\op{coh}}(\sBun,\mathcal{D}).}
\]
This $\mathfrak{c}$ \index{notations}{c@$\mathfrak{c}$} needs to
intertwine the action of the tensorization functors on
$D_{\op{coh}}(\lan{\gLoc},\mathcal{O})$
\index{notations}{DcohLFlat@$D_{\op{coh}}(\lan{\gLoc},\mathcal{O})$}
with the action of the Hecke
functors on $D_{\op{coh}}(\sBun,\mathcal{D})$,
\index{notations}{DcohBun@$D_{\op{coh}}(\sBun,\mathcal{D})$}
uniformly over
$C$. \index{notations}{C@$C$}
Explaining this last sentence is unfortunately a bit technical:

\subsubsection{Tensorization, Hecke, Intertwining}

\

\noindent
The tensorization functors $\lan{W}^{\mu,x} :
D_{\op{coh}}(\lan{\gLoc},\mathcal{O}) \to D_{\op{coh}}(\lan{\gLoc},\mathcal{O})$,
and the Hecke functors $H^{\mu,x} :
D_{\op{coh}}(\sBun,\mathcal{D})\to
D_{\op{coh}}(\sBun,\mathcal{D})$, are endofunctors of the
\index{notations}{LWmux@$\lan{W}^{\mu,x}$} 
\index{notations}{Hmux@$H^{\mu,x}$}
\index{terms}{tensorization!functors}
\index{terms}{Hecke!functors}
respective categories of sheaves. They are labeled by the same data:
pairs $(x,\mu)$, where $x \in C$ is a closed point and $\mu \in
\cchr^{+}(G) = \chr^{+}(\lan{G})$ is a dominant cocharacter for
$G$, or equivalently a dominant character for
$\lan{G}$. \index{terms}{dominant!character}
\index{terms}{dominant!cocharacter}
\index{notations}{char@$\chr^{+}(G)$, $\chr^{+}(\lan{G})$}
\index{notations}{cocharG@$\cchr^{+}(G)$, $\cchr^{+}(\lan{G})$}

Given such a pair $(x,\mu)$, one defines the {\bfseries\em
  tensorization functor} $\lan{W}^{\mu,x}$ as
\[
\xymatrix@R-1.5pc{
\lan{W}^{\mu,x} : & \hspace{-0.3in} D_{\op{coh}}(\lan{\gLoc},\mathcal{O}) \ar[r] &
D_{\op{coh}}(\lan{\gLoc},\mathcal{O})  \\
& \mycal{F} \ar@{|->}[r] & \mycal{F} {\otimes}
\rho^{\mu}\left(\mycal{V}_{|\lan{\gLoc}\times \{ x\}}\right),
}
\]
where $\mycal{V} \to \lan{\gLoc}\times C$ is the principal $\lan{G}$-bundle
underlying the universal flat bundle
$\mathcal{V} = (\mycal{V},\nabla)$,
\index{notations}{Vuniv@$\mycal{V}$} \index{terms}{universal!bundle}
\index{terms}{universal!flat bundle}
\index{notations}{Vunivf@$\mathcal{V}$}
$\rho^{\mu}$ is the irreducible
representation of $\lan{G}$ with highest weight $\mu$,
\index{notations}{mu@$\mu$} \index{terms}{dominant!character}
\index{terms}{dominant!cocharacter}
\index{notations}{rhomu@$\rho^{\mu}$} \index{terms}{representation!
  of highest weight $\mu$} and
$\rho^{\mu}\left(\mycal{V}_{|\lan{\gLoc}\times \{ x\}}\right)$ is the vector
bundle on $\lan{\gLoc}$ associated with $\mycal{V}_{|\lan{\gLoc}\times \{ x\}}$
via the representation $\rho^{\mu}$. \index{terms}{associated vector
  bundle}

  \

\noindent
To define the Hecke functors, we need first to construct the Hecke
correspondences: \index{terms}{Hecke!correspondence}
\index{terms}{Hecke!stack} 

\begin{itemize}
\item The {\bfseries\em Hecke stack} is the moduli
  stack $\sHeck$ of quadruples $(V,V',x,\mathsf{elm})$, where
  \index{notations}{Heck@$\sHeck$}
	\begin{itemize}
		\item $V$, $V'$ are principal $G$-bundles on $C$,
		\item $x \in C$ is a point of the curve,
		\item $\xymatrix@1{\mathsf{elm} : \hspace{-1.8pc} &
                  V_{|C -\{x\}} \ar[r]^-{\cong} & V'_{|C -\{x\}}}$ 
                  (an "elementary modification" of  $V$ at $x$)  is
                  an isomorphism of principal bundles away from the
                  point. \index{notations}{elm@$\mathsf{elm}$}
	\end{itemize}
\item $\sHeck^{\mu}$ is the closed substack of $\sHeck$ of quadruples
  $(V,V',x,\mathsf{elm})$ where the order of the pole of
  $\mathsf{elm}$ at $x$ is bounded in terms of $\mu$. The precise
  condition is that if \index{notations}{Heckmu@$\sHeck^{\mu}$}
  $\lambda \in \chr^{+}(G)$ is a dominant cocharacter and if
  $\rho^{\lambda}$ is the irreducible representation of $G$ with
  highest weight $\lambda$, then $\mathsf{elm}$ induces an inclusion
  of locally free sheaves \linebreak $\rho^{\lambda}(\mathsf{elm}) :
  \rho^{\lambda}(V) \hookrightarrow \rho^{\lambda}(V')\otimes
  \mathcal{O}_{C}(\langle \mu, \lambda\rangle x)$.
\item $\sHeck^{\mu,x}$ is the closed substack of
  $\sHeck^{\mu}$ of quadruples $(V,V',x,\mathsf{elm})$ with the
  specified $x$. \index{terms}{Hecke!stack}
  \index{terms}{Hecke!correspondence}
  \index{notations}{Heckmux@$\sHeck^{\mu,x}$}  
\item   These stacks are equipped with natural projections
  \[
{\small
\xymatrix@C-1.9pc@R-0.5pc{& \sHeck \ar[dl]_-{p} \ar[dr]^-{q} & \\
\sBun & & \sBun\times C} \quad
\xymatrix@C-1.9pc@R-0.5pc{&
  \sHeck^{\mu} \ar[dl]_-{p^{\mu}} \ar[dr]^-{q^{\mu}} & \\ 
\sBun & & \sBun\times C} \quad
\xymatrix@C-1.9pc@R-0.5pc{
& \sHeck^{\mu,x} \ar[dl]_-{p^{\mu,x}} \ar[dr]^-{q^{\mu,x}} & \\
\sBun & & \sBun\times C,}
}
\]
where $p(V,V',x,\mathsf{elm}) := V$, $q(V,V',x,\mathsf{elm}) := V'$,
while $p^{\mu}, q^{\mu}, p^{\mu,x}$ and $q^{\mu,x}$
are the restrictions of $p$ and $q$ to
$\sHeck^{\mu}$ and $\sHeck^{\mu,x}$
respectively. \index{notations}{p@$p$} \index{notations}{q@$q$}
\index{notations}{pmu@$p^{\mu}$} \index{notations}{qmu@$q^{\mu}$}
\index{notations}{pmux@$p^{\mu,x}$} \index{notations}{qmux@$q^{\mu,x}$}
\item The {\bfseries\em Hecke kernel} \ $\mathfrak{I}^{\mu,x}$ is
  the Deligne-Goresky-MacPherson middle perversity
  \index{terms}{DeligneGG@Deligne-Goresky-MacPherson!extension}
  \index{terms}{middle perversity
    extension|see{\\ Deligne-Goresky-MacPherson \linebreak extension}}
  \index{terms}{minimal extension|see{\\ Deligne-Goresky-MacPherson
      \linebreak extension}} \index{terms}{Hecke!kernel}
  \index{notations}{Imux@$\mathfrak{I}^{\mu,x}$} extension
  $j_{!*}\left((\mathcal{O},d)\left[\dim
    \sHeck^{\mu,x}\right]\right)$ of the trivial rank one flat
  bundle on the smooth part \linebreak $j :
  \left(\sHeck^{\mu,x}\right)^{\op{smooth}} \hookrightarrow
  \sHeck^{\mu,x}$ of the Hecke stack. $\mathfrak{I}^{\mu}$
  is defined similarly.
\end{itemize}

\

\noindent The {\bfseries\em Hecke functor} $H^{\mu}$ is defined
as the integral transform \index{terms}{Hecke!functor}
\index{notations}{Hmu@$H^{\mu}$}
\index{notations}{Hmux@$H^{\mu,x}$}
\[
\xymatrix@R-1.5pc{ H^{\mu} : & \hspace{-0.3in}
D_{\op{coh}}(\sBun,\mathcal{D}) \ar[r] &
D_{\op{coh}}(\sBun \times C  ,\mathcal{D})\\ & \mycal{M} \ar@{|->}[r] &
q^{\mu}_{!}\left((p^{\mu})^{*}\mycal{M}\otimes
\mathfrak{I}^{\mu}\right). }
\]
The Hecke functor $H^{\mu,x}$ is defined analogously, replacing
the curve $C$ by its point $x$:
\[
\xymatrix@R-1.5pc{ H^{\mu,x} : & \hspace{-0.3in}
D_{\op{coh}}(\sBun,\mathcal{D}) \ar[r] &
D_{\op{coh}}(\sBun,\mathcal{D})\\ & \mycal{M} \ar@{|->}[r] &
q^{\mu,x}_{!}\left((p^{\mu,x})^{*}\mycal{M}\otimes
\mathfrak{I}^{\mu,x} \right). }
\]
Finally, the intertwining property
\index{terms}{intertwining property}
can and should be stated over $C$
in terms of the $H^{\mu}$ rather than point-by-point in terms of
the $H^{\mu,x}$:
\index{notations}{C@$C$}

Using the universal relative flat $G$-bundle
\index{notations}{Vunivf@$\mathcal{V}$} 
$\mathcal{V} \to \lan{\gLoc} \times C$, 
the tensorization functors $\lan{W}^{\mu,x}$ extend to
\index{notations}{Wmu@$W^{\mu}$} 
\[
\lan{W}^{\mu} : D_{\op{coh}}(\lan{\gLoc},\mathcal{O}) \to
D_{\op{coh}}(\lan{\gLoc} \times C,  p_{C}^{*}\mathcal{D}),
\]
and any functor
\[
\mathfrak{c} : D_{\op{coh}}(\lan{\gLoc},\mathcal{O})
\to
D_{\op{coh}}(\sBun,\mathcal{D}).
\]
extends naturally to a functor
(which should be denoted $\mathfrak{c} \boxtimes 1$,
but we still call it $\mathfrak{c}$ to avoid clumsy notation):
\[
\mathfrak{c} : D_{\op{coh}}(\lan{\gLoc} \times C, p_{C}^{*}\mathcal{D})
\to
D_{\op{coh}}(\sBun \times C, \mathcal{D}).
\]
Now the full intertwining condition is:
\[
\mathfrak{c} (\lan{W}^{\mu}(F)) = H^{\mu}(\mathfrak{c}(F))
\]
for all $\mu$ and all coherent sheaves $F$ on $\lan{\gLoc}$.
\index{notation}{DcohLFlat@$D_{\op{coh}}(\lan{\gLoc},\mathcal{O})$}
\index{notation}{DcohBun@$D_{\op{coh}}(\sBun,\mathcal{D})$}

\subsubsection{Hecke eigensheaves}

The GLC implies that $\mathfrak{c}$ sends the structure sheaves of
points $\mathbb{V}$ in $\lan{\gLoc}$ to Hecke eigen sheaves (or more
precisely, Hecke eigen $\mathcal{D}$-modules) on $\sBun$.  This
is simply the special case of the intertwining property when $F$ is
the structure sheaf $\mathcal{O}_{\mathbb{V}}$ of a point:
\index{terms}{Hecke!eigensheaf} \index{terms}{HES|see{Hecke eigensheaf}}
\[ \tag{{\bfseries\sf HES}} \label{eq:hes}
\boxed{H^{\mu}\left(\mathfrak{c}(\mathcal{O}_{\mathbb{V}})\right) =
\mathfrak{c}(\mathcal{O}_{\mathbb{V}})\boxtimes \rho^{\mu}(\mathbb{V}).}
\]
In fact, the GLC is uniquely characterized by this property. Due to
some subtleties (see \cite{dp-Langlands,dp-jdg,ag-sing,bznp}) arising
from choices of connected components, twists by gerbes, and singular
supports of coherent sheaves, the statement of GLC needs to be
modified in several ways.  These modifications have only a minor
effect on the existence of Hecke eigen sheaves.  In the present work
we focus on the explicit construction of Hecke eigen sheaves for
neutral components and generic inputs which suppresses all such
modifications. \index{terms}{singular support conditions}

\subsubsection{$GL_n$}

Consider the case $G = GL_{n}(\mathbb{C})$. Then 
$\lan{G} = GL_{n}(\mathbb{C})$ and $\lan{\gLoc}$ can be identified with the
  stack of rank $n$ vector bundles on $C$ equipped with an integrable
  connection. \index{terms}{flat!vector bundle}
Among the Hecke correspondences ${\sHeck}^{\mu}$, those specified by 
the {\em fundamental} weights
$\mu =i \in \{ 0 \dots n \}$ of $GL_{n}(\mathbb{C})$
are particularly accessible:
\[
\sHeck^{i} := \left\{ (V,V',x) \left| \text{\begin{minipage}[c]{3in}
$V$ and $V'$ are locally free sheaves of rank $n$ such that 
$V \subset V' \subset V(x)$ and $\op{length}(V'/V) = i$
\end{minipage}} \right.\right\}.
\]
The Hecke operators $H^{i}$ given by the correspondences $\sHeck^{i}$
generate (the commutative algebra of) all Hecke operators $H^{\mu,x}$.
The fibers of the projection $q^{i} : \sHeck^{i} \to \sBun\times C$
are all isomorphic to the Grassmanian $Gr(i,n)$ of $i$-dimensional
subspaces in an $n$-dimensional vector space.
\index{terms}{Hecke!fibers} \index{notations}{Heckei@$\sHeck^{i}$}
\index{notations}{Hi@$H^{i}$}

\subsubsection{The Abelian case}\label{abelian}

The abelian case
$G \cong GL_{1}(\mathbb{C}) \cong {\mathbb{C}}^{\times} \cong \lan{G}$
of GLC is the geometric analogue of Class Field
Theory. \index{terms}{Geometric Langlands Conjecture!abelian case of}
\index{terms}{Class Field Theory} This is essentially the well-known
theory of curves and their Jacobians.  In this case $\sBun = \gPic(C)$
is the Picard stack of $C$. To avoid cumbersome technical statements
arising from derived structures and singular support conditions on
$\gLoc$ we will focus on the part of the $GL_{1}(\mathbb{C})$ Geometric Langlands Conjecture that
takes place on the Picard variety $\Pic(C)$ of
$C$. \index{terms}{Picard!variety} \index{notations}{PicC@$\Pic(C)$}
\index{terms}{Jacobian} \index{terms}{singular support conditions}

Here there is only one interesting Hecke operator
\index{terms}{Picard!stack} \index{notations}{PicsC@$\gPic(C)$}
\[
H^{1} : D_{\op{coh}}(\Pic(C),\mathcal{D}) \to D_{\op{coh}}(C\times
\Pic(C),\mathcal{D}) 
\]
which is simply the pull-back $H^{1} := \AJ^{*}$ 
via the classical Abel-Jacobi map \index{terms}{Abel-Jacobi map}
\index{notations}{aj@$\AJ$}
\[
\xymatrix@R-2pc{
\AJ : \hspace{-2pc} & C\times \Pic^{d}(C) \ar[r] &  \Pic^{d+1}(C) \\
& (x,L) \ar@{|->}[r] & L(x).
}
\]
In this case the geometric Langlands correspondence $\mathfrak{c}$ can
be described explicitly. Let $\mathbb{L} = (L,\nabla)$ be a rank one
flat bundle on $C$. Since the fundamental group $\pi_{1}(\Pic^{d}(C))$
of the $d$-th component of the Picard variety $\Pic(C)$ is the
abelianization of $\pi_{1}(C)$ and the monodromy representation of
$\mathbb{L}$ is abelian, it follows that we can view $\mathbb{L}$ as a
flat bundle  on each component $\Pic^{d}(C)$ of $\Pic(C)$. 
 \index{terms}{Picard!variety!component of}
 We get: \index{terms}{Picard!variety}
\index{notations}{PicC@$\Pic(C)$}\index{notations}{PicCcomp@$\Pic^{d}(C)$} 
\begin{spacing}{1.0}
\[
  \mathfrak{c}\left(\mathcal{O}_{\mathbb{L}}\right) :=
  \left(\text{\begin{minipage}[c]{2.5in} the unique
    rank one flat bundle on $\Pic(C)$ whose restriction on each
    component $\Pic^{d}(C)$ has the same monodromy as $\mathbb{L}$
    \end{minipage}}\right).
\]
\end{spacing}

\subsubsection{Literature}

\noindent
Many cases of these conjectures have been proven through different
approaches and constructions.
\cite{drinfeld-icm,drinfeld-ajm,drinfeld-jsm},
\cite{laumon-langlands}, \cite{beilinson-drinfeld-langlands},
\cite{lafforgue}, \cite{fgkv}, \cite{fgv}, \cite{dennis-nearby},
\cite{laumon-fgv}. These include profs of versions of the conjecture
for $GL_{2}$ \cite{drinfeld-ajm} and later, using Lafforgue's
spectacular work \cite{lafforgue}, also for $GL_{n}$
\cite{fgv,dennis-nearby}. In contrast only special cases of the
conjecture are understood for other groups, or in the parabolic case,
but see \cite{beilinson-drinfeld-langlands}, \cite{hny},
\cite{vlafforgue}, \cite{yun-rigidity} and \cite{arinkin-lysenko},
\cite{arinkin-4points}, \cite{heinloth}. Even for $GL_n$, the proof
is indirect: no construction of non-abelian Hecke eigensheaves is
known, except perhaps for the works of Bezrukavnikov-Braverman
\cite{bb}, Shen \cite{shen}, and Travkin \cite{travkin} which for
curves over finite fields implement an approach similar to ours.

\subsection{The program}

The program outlined in \cite{dp-jdg} aims to convert the GLC for
general group $G$ to the abelian case.  Two key tools are Non-Abelian
Hodge Theory (NAHT) and Hitchin's integrable system. This approach can
be viewed as a mathematical incarnation of important physics ideas
introduced in the fundamental works of Kapustin-Witten \cite{kw},
Gukov-Witten \cite{gukov-witten}, and Frenkel-Witten
\cite{frenkel-witten} where the geometric Langlands correspondence is
interpretted as a duality in quantum field theory.

\index{terms}{Non-abelian Hodge theory}
\index{terms}{Hitchin!integrable system}

NAHT sets up an equivalence between moduli spaces of three types of
objects on a projective variety $X$:
\begin{itemize}
\item $\gHiggs_{X,G}$, parametrizing $G$-Higgs bundles on $X$
  (Dolbeault incarnation),
\item $\gLoc_{X,G}$, parametrizing $G$-bundles with a flat connection
  (de Rham incarnation), and
\item $\op{\boldsymbol{\mathcal{L}oc}}_{X,G}$, parametrizing $G$-local
  systems (Betti incarnation).
\end{itemize}
\index{terms}{Betti!version}
\index{terms}{de Rham!version}
\index{terms}{Dolbeault!version}
The objects are subject to restrictions on their Chern classes, and
the equivalence can be interpreted as a categorical equivalence or as
a diffeomorphism of underlying moduli spaces
$\Higgs, \op{\bf Flat}, \op{\bf Loc}$ of semistable objects. 
\index{terms}{Higgs!bundle}
\index{terms}{principal!Higgs bundle}
\index{terms}{principal!flat bundle}
\index{terms}{flat!bundle}
\index{terms}{local system}
\index{notations}{LocsXG@$\op{\boldsymbol{\mathcal{L}oc}}_{X,G}$}
\index{notations}{HiggssXG@$\gHiggs_{X,G}$}
\index{notations}{FlatsXG@$\gLoc_{X,G}$}
\index{notations}{Loc@$\op{\bf Loc}$}
\index{notations}{Higgs@$\Higgs$}
\index{notations}{Flat@$\op{\bf Flat}$}

{\bfseries\em Hitchin's integrable system} for a curve $C$ and a group
$G$ is a natural proper morphism $\bh: \Higgs_{C,G} \to \bB$ called the
{\bfseries\em Hitchin map}.  Here the {\bfseries\em Hitchin base} $\bB$
is an appropriate vector space of pluridifferentials on
$C$. $\Higgs_{C,G}$ is the moduli space of semistable Higgs
$G$-bundles on $C$. It is a quasi projective variety which is
birational to the cotangent bundle $T^{\vee}\Bun$ of $\Bun$, and has a
natural symplectic structure for which the map $\bh$ is Lagrangian. For
generic $b \in \bB$, the Hitchin fiber $\bh^{-1}(b)$ is an abelian variety
that can be described explicitly as a generalized Prym variety.
\index{terms}{Hitchin!map} \index{terms}{Hitchin!fiber}
\index{terms}{Hitchin!base} \index{terms}{Prym!variety!generalized}
\index{notations}{Bun@$\Bun$}
\index{notations}{HiggsCG@$\Higgs_{C,G}$}
\index{notations}{h@$\bh$} \index{notations}{B@$\bB$}

The main result of \cite{dp-Langlands} is a duality of Hitchin systems
for a Langlands-dual pair of groups: the Hitchin bases for the two
systems are identified, and the non-singular fibers over matching base
points are dual abelian varieties: a point on one determines a degree
0 line bundle on the other, and vice versa. As explained there, this
result can be viewed as an abelianized version of GLC.
\index{terms}{Geometric Langlands Conjecture}
\index{terms}{Non-abelian Hodge theory}

Our approach, here as in \cite{dp-jdg}, is based on using NAHT to
reduce GLC to the abelianized version proved in \cite{dp-Langlands}.
The idea is to use NAHT twice: first for flat $\lan{G}$-bundles on
$C$, then for flat vector bundles on open subsets of the moduli space
$\Bun$ of semistable $G$-bundles on $C$.
\index{terms}{Geometric Langlands Conjecture}
\index{terms}{Non-abelian Hodge theory} The
input is a $\lan{G}$-bundle with flat connection on $C$.  The first
application of NAHT converts this to a $\lan{G}$-Higgs bundle on $C$.
This is interpreted as a point of $\lan{\Higgs} = \Higgs_{C,\lan{G}}$.
The duality for the Hitchin system transforms this to a line bundle on
a fiber of the Langlands-dual Hitchin system $\Higgs =
\Higgs_{C,G}$. Via the rational map $\Higgs \dashrightarrow \Bun$, we
view this fiber as a (meromorphic) spectral cover of $\Bun$, and refer
to it as the {\bfseries\em Modular Spectral Cover}.
\index{terms}{spectral cover!modular} \index{terms}{Higgs!bundle}
\index{terms}{flat!bundle} Pushing down by the rational map $\Higgs
\dashrightarrow \Bun$ gives us a Higgs sheaf on
\index{terms}{Higgs!sheaf} $\Bun$.  This turns out to be a Higgs
bundle (i.e. locally free) away from a particular divisor that we call
the {\bfseries\em Wobbly locus}, $\bWob \subset \Bun$.
\index{terms}{wobbly!locus} \index{terms}{wobbly!divisor}
\index{notations}{Wob@$\bWob$} The second application of NAHT then
converts this Higgs bundle back to a flat bundle on $\Bun
- \bWob$.  By taking the Deligne-Goresky-MacPherson middle
perversity extension we obtain the desired Hecke eigen $\mathcal{D}$-module on
$\lan{\Bun}$.
\index{terms}{DeligneGG@Deligne-Goresky-MacPherson!extension}
\index{notations}{LBun@$\lan{\Bun}$} \index{terms}{Dmodule@$\mathcal{D}$-module}

It is crucial to note that the second application of NAHT takes place
on the non-compact space $\Bun - \bWob$.  To handle this open case,
one replaces objects on the open space by parabolic objects on the
compactification $\Bun$, with parabolic structure on the boundary.  We
discuss various kinds of parabolic objects in section
\ref{sec:parabolic.objects}.  A major result is Mochizuki's open NAHT
which is reclled in section \ref{Mo}, cf.  also \cite{mochizuki-kh1}
and \cite[Theorem~1.1 and Corollary~1.5]{mochizuki-kh2}.
\index{terms}{Mochizuki's!open NAHT} \index{terms}{parabolic!objects}

Much of the hard calculations in the present work is needed in order
to verify that the conditions in Mochizuki's theorem, especially the
stability and the vanishing of parabolic Chern classes, hold in our
situation.  \index{terms}{Mochizuki's!conditions}
\index{terms}{parabolic!Chern classes}

A key feature is the appearance of the {\bfseries\em wobbly} locus.  This is
the locus of bundles on $C$ that are stable, but not {\em very
  stable}.  According to Laumon \cite{laumon-nilpotent, laumon-gln}, a
vector bundle $V$ on a curve $C$ is {\bfseries\em very stable} if the
only nilpotent Higgs field on $V$ is the zero one.  A very stable
bundle is automatically stable.  A stable bundle that is not very
stable is called {\bfseries\em wobbly}.  The wobbly locus
$\bWob \subset \Bun$ is a divisor.  As we see in section
\ref{sssec:wobbly}, it is the locus over which the Hecke eigen sheaves
will have singularities.  In our application of Mochizuki's theorem we
therefore consider Higgs and flat bundles with parabolic structure
along $\bWob$. \index{notations}{Wob@$\bWob$}
\index{terms}{bundle!verry stable} \index{terms}{bundle!stable}
\index{terms}{bundle!wobbly} \index{terms}{Higgs!field!nilpotent}
\index{terms}{Laumon}

\subsection{Parabolic version}

So, starting with flat bundles on $C$, we are led to Higgs bundles on
$\Bun$ with poles along the wobbly locus.  Working a little more
carefully, we see that these actually have the structure of parabolic
Higgs bundles on $(\Bun, \bWob).$ By a theorem of Mochizuki, these are
equivalent to $\mathcal{D}$-modules on $\Bun$.  The notion of parabolic bundle,
and the basic results we need, are reviewed in section
\ref{sec:parabolic.objects}. \index{terms}{Mochizuki's!theorem}
\index{terms}{parabolic!Higgs bundle} \index{terms}{parabolic!bundle}
\index{terms}{Dmodule@$\mathcal{D}$-module}

Since parabolic bundles arise naturally, one may ask whether it is
possible to extend the GLC to include parabolic bundles in the input
as well.  In other words, we want to allow the input $\mathbb{V}$ to
be ramified along the parabolic divisor $\Par_C$, a reduced divisor on
our $C$ consisting of $k$ distinct points.  We take $\mathbb{V}$ to be
a tame parabolic flat bundle $\mathbb{V} = (V,\nabla)$ on
$(C, \Par_C)$.  From this we want to produce a Hecke eigensheaf whose
eigenvalue is $\mathbb{V}$.  This is a rather drastic extension, and
it exhibits several new features.
\index{terms}{parabolic!flat bundle!tame} \index{terms}{Hecke!eigensheaf}
\index{terms}{parabolic!divisor} \index{notations}{ParC@$\Par_{C}$}
\index{notations}{C@$C$} \index{notations}{Vbb@$\mathbb{V}$}
\index{notations}{Vnabla@$(V,\nabla)$}

\subsubsection{The parabolic abelian case}\label{parabelian}
To get a sense for this, start again with the abelian case.  The
first, obvious, observation is that a ramified rank one  flat bundle
$\mathbb{L}$ on $C  -  \Par_C$ cannot be the pullback of a
similar object on the Jacobian via the $\AJ$ map: the map on
fundamental groups $\pi_1(C  -  \Par_C) \to \pi_1(\Pic^d(C))$
sends loops around points of $\Par_C$ to $0$, so pulling back gives only
unramified flat bundles. \index{terms}{flat!bundle!ramified}
\index{terms}{flat!bundle!unramified}

Nevertheless, this is easy to fix: we simply replace $\Jac:=\Pic^0(C)$
by a commutative algebraic group $\underline{\Jac}$ which has an
Abel-Jacobi map $\AJ : C  -  \Par_C \to \underline{\Jac}$
inducing an isomorphism on first
homologies. \index{terms}{Jacobian!generalized}
\index{notations}{J@$\Jac$} \index{notations}{Jnc@$\underline{\Jac}$}
\index{terms}{Jacobian} This $\underline{\Jac}$ is non-compact: it is
a $({\mathbb{C})^{\times}}^{(k-1)}$-bundle over $\Jac$.  We can take it
to be the Albanese variety of the open curve $C - \Par_C$, i.e. the
algebraic group built from \index{terms}{Albanese variety!of a
  non-compact curve} the mixed Hodge structure
$H_1(C  -  \Par_C, \mathbb{Z})$. Equivalently, we can take it to
be the connected component $\Pic^0(\underline{C})$ of the Picard of
the singular, compact curve $\underline{C}$ obtained from $C$ by gluing
the $k$ points of $\Par_C$ to each other transversally into a normal
crossings singularity.
\index{terms}{Jacobian!generalized!of a  singular curve}
\index{terms}{normal crossings!singularity}
\index{notations}{Cun@$\underline{C}$}
Either way, the Abel-Jacobi map induces an
isomorphism on first homologies.  As in the unramified case, this
implies that every abelian local system on $C  -  \Par_C$
extends uniquely to a Hecke eigensheaf on $\underline{\Jac}$.  In
particular, it follows that $\AJ$ does not extend to a morphism of $C$
to $\underline{\Jac}$.  (One can construct a natural compactification of
$\underline{\Jac}$ to which $\AJ$ does extend.)
\index{terms}{Hecke!eigensheaf} 

\subsubsection{The parabolic GLC}
This suggests the natural formulation of the tamely ramified, or
parabolic, version of geometric Langlands, which essentially goes back
to the classical work of Drinfeld \cite{drinfeld-jsm}.
\index{terms}{Geometric Langlands Conjecture!tamely ramified} There
are subtleties to the notion of moduli of parabolic bundles, since
those involve real weights.  We therefore start by considering the
moduli stack of quasiparabolic bundles.  For the parabolic version of
$\lan{\gLoc}$ we take the moduli stack of quasiparabolic
$\lan{G}$-flat connections on $(C, \Par_C)$ with at most logarithmic
poles along $\Par_{C}$.
\index{terms}{moduli stack!of quasiparabolic flat connections}
\index{terms}{flat!connection!logarithmic}
For the parabolic analogue
of $\sBun$ there are two options: the moduli stack $\sBun$
of $G$-bundles on $C$ with a reduction of their structure group
over $\Par_C$ to the Borel $B \subset G$, and the moduli stack
$\underline{\sBun}$ of $G$-bundles on $C$ with a
reduction of their structure group over $ \Par_C$ to the unipotent
radical $[B,B] \subset B \subset G$.  Note that $\underline{\sBun}$
is fibered over $\sBun$, the fiber being the Cartesian product
$T^{ (k-1)}$, where
$T := B/ [B,B]$ is the maximal torus of
$G$, generalizing the relation between $\Jac$ and
$\underline{\Jac}$ in the abelian case.  \index{terms}{moduli stack!of
  quasiparabolic bundles} \index{terms}{moduli stack!of bundles
  equipped with a reduction of their structure group to a unipotent
  radical of a Borel over $\Par_{C}$}
\index{notations}{Bunsun@$\underline{\sBun}$}
\index{notations}{T@$T$}

The parabolic GLC then predicts the existence of a fully faithful
functor of derived categories, relating coherent $\mathcal{O}$-modules
on $\lan{\gLoc}$ to coherent $\mathcal{D}$-modules on
$\underline{\sBun}$:
  \index{terms}{parabolic GLC|see{Geometric Langlands \linebreak
      Conjecture, tamely ramified}}
  \index{terms}{PGLC@$\mathsf{PGLC}$|see{Geometric Langlands
      \linebreak Conjecture, tamely ramified}}
\[
\tag{{\bfseries\sf PGLC}} \label{eq:pglc}
\boxed{
\mathfrak{c} : D_{\op{coh}}(\lan{\gLoc},\mathcal{O})
\stackrel{\cong}{\longrightarrow}
D_{\op{coh}}(\underline{\sBun},\mathcal{D}).}
\]
This $\mathfrak{c}$ needs to intertwine the action of the
tensorization functors on $D_{\op{coh}}(\lan{\gLoc},\mathcal{O})$ with the
action of the Hecke functors on
$D_{\op{coh}}(\underline{\sBun},\mathcal{D})$, uniformly over
$C$. The tensorization and Hecke functors and their intertwining
property are the straightforward analogues of the unramified case. 

The
intertwining property implies the existence of parabolic Hecke eigen
sheaves:
\index{terms}{tensorization!functors}
\index{terms}{Hecke!functors}
\index{terms}{Hecke!eigensheaves!parabolic}
\index{terms}{intertwining property}
\index{terms}{Langlands!functor}
The Langlands functor $\mathfrak{c}$ sends the structure sheaves of
points $\mathbb{V}$ in the parabolic $\lan{\gLoc}$ to Hecke eigen sheaves
(or more precisely, Hecke eigen $\mathcal{D}$-modules) on the
parabolic $\underline{\sBun}$.  As previously, this is simply the special
case of the intertwining property when $F$ is the structure sheaf
$\mathcal{O}_{\mathbb{V}}$ of a point:
\[ \tag{{\bfseries\sf PHES}} \label{eq:phes}
\boxed{H^{\mu}\left(\mathfrak{c}(\mathcal{O}_{\mathbb{V}})\right) =
\mathfrak{c}(\mathcal{O}_{\mathbb{V}})\boxtimes \rho^{\mu}(\mathbb{V}).}
\]
We see immediately that the objects and the correspondence specialize
correctly to the abelian case. In particular,
$\underline{\sBun}$ specializes to what we were calling above
$\underline{\Jac}$, while $\sBun$ specializes to $\Jac$.  As we
will see in the appendix, the {\sf PGLC} can be stated directly on
$\sBun$.  The Hecke eigen $\mathcal{D}$-module we want on
$\underline{\sBun}$ can instead be interpreted as a twisted \linebreak
$\mathcal{D}$-module on the $T^{(k-1)}$-quotient $\sBun$.  The
twisting itself is controlled by the residue of the  meromorphic
connection $\nabla$ for the eigenvlue $\mathbb{V} = (V,\nabla)$; the
untwisted case arises from $\nabla$'s that do not have a pole.
\index{notations}{Buns@$\sBun$}
\index{notations}{Bunsun@$\underline{\sBun}$}
\index{notations}{LFlats@$\lan{\gLoc}$}
\index{notations}{LT@$\lan{T}$}
\index{terms}{Geometric Langlands Conjecture!tamely ramified}
\index{terms}{flat!connection!logarithmic}
We
discuss various versions of the conjecture in the appendix, where the
emphasis is on the equivalent formulation in terms of twisted
differential operators, which is common in the literature.

The parabolic version of our program calls for a proof of the {\sf PGLC} by
reducing it to the abelian case via open NAHT. We review the basic
steps in this program in section \ref{program} of the
appendix. \index{terms}{Non-abelian Hodge theory}

\subsection{Our results}

Use of NAHT in the context of the Hecke eigensheaf condition requires
calculation of Chern classes for direct images of Higgs (and flat)
bundles.  In order for the open version of NAHT to apply to the map
$q: \Heck \to \sBun\times C$ where each of these spaces comes with the
appropriate parabolic divisor, the relevant objects must have
vanishing (parabolic) Chern classes. In \cite{dps} we developed
algebraic formulas for the direct images when the fibers are
semistable curves and the parabolic Higgs bundles satisfy a nilpotence
assumption along the parabolic divisor (cf. (2.2) in
\cite{dps}).

\index{terms}{parabolic!divisor}
\index{terms}{parabolic!Chern classes}
\index{terms}{parabolic!Higgs bundle!tame}
\index{terms}{parabolic!flat bundle!tame}
\index{terms}{nilpotence assumption}

In the Dolbeault setting the {\bfseries\em nilpotence assumption}
 for a tame parabolic Higgs bundle is
just the condition that the residues
of the Higgs field are nilpotent. In the de Rham setting this is
equivalent to Mochizuki's {\bfseries\em purely imaginary condition}
\cite{mochizuki-kh1} which requires for a tame parabolic flat vector
bundle that the eigenvalues of the residues be real and the
parabolic weights be the negative of these
eigenvalues. In the Betti setting the condition is equivalent to the
condition that the
tame parabolic local system has a trivial parabolic structure and
that the eigenvalues of the local monodromies belong to the unit circle.
\index{terms}{purely imaginary condition}
\index{terms}{Betti!version}
\index{terms}{Dolbeault!version}
\index{terms}{de Rham!version}
In the
proof of the algebraic pushforward formula in \cite{dps} the condition
is needed for a technical comparison between the growth filtration for
a harmonic metric and Kashiwara's $V$-filtration for the associated
$\mathcal{D}$-module.  \index{terms}{growth filtration}
\index{terms}{Vfiltration@$V$-filtration} We believe that the
comparison and the direct image formula hold without imposing the
nilpotence condition but we have not pursued this issue.

The natural structure \cite{eyal,bottacin} on the moduli of
meromorphic Higgs bundles is that of an (algebraic or analytic)
Poisson manifold. Such manifolds come with a natural foliation by
symplectic leaves, holomorphic submanifolds on which the Poisson
structure is non-degenerate, i.e. symplectic. In the case of the
moduli of meromorphic Higgs bundles, Markman shows that this foliation
by symplectic leaves is actually algebraic, specified by fixing the
conjugacy classes (or more precisely, the coadjoint orbits) of the
residues. This picture extends verbatim to our parabolic setup.

In the setup of the {\sf PGLC} the effect of imposing the nilpotence
condition is that it restricts the Hecke eigenvalues to be points in a
particular
symplectic leaf in the Poisson moduli of parabolic Higgs (or flat)
bundles. Again we believe that this restriction is unnecessary but we
have not developed the necessary NAHT to attempt to tackle other
symplectic leaves. \index{terms}{Hecke!eigenvalues}
\index{terms}{Poisson!moduli of parabolic Higgs bundles}
\index{terms}{Poisson!moduli of parabolic flat bundles}

One of the few cases where an explicit solution of {\sf PGLC} is known
over the complex numbers is when $G=\mathbb{P}SL_{2}(\mathbb{C})$ or
$GL_{2}(\mathbb{C})$, $C=\mathbb{P}^1$ and $\Par_{C}$ consists of $4$
points \index{notations}{SL2@$SL_{2}(\mathbb{C})$}
\index{notations}{GL2@$GL_{2}(\mathbb{C})$} \cite{arinkin-lysenko}. In
that case the moduli space itself is a curve. In fact, it can be
identified with $C$, and then the wobbly locus is identified with
$\Par_{C}$.  \index{terms}{Geometric Langlands Conjecture!tamely
  ramified} The technique used in \cite{arinkin-lysenko} is very
different than ours. They do not use NAHT and in particular they do
not need the analogue of our nilpotence assumption.  (They do impose
the analogue of the non-resonance condition, see below.)

In order to test our proposal, we work it out in full detail for the
case when $G=\mathbb{P}SL_{2}(\mathbb{C})$ or $GL_{2}(\mathbb{C})$,
$C=\mathbb{P}^{1}$, the parabolic divisor $\Par_{C}$ consists of $5$
points, and the flat bundle is purely imaginary - the condition that
translates into the nilpotence assumption of \cite{dps} for the
residues of the Higgs fields.  We will see that in this case the
moduli space is a del Pezzo quartic surface $X \subset \mathbb{P}^4$
and the wobbly locus consists of the 16 lines on that surface.  We
apply the results of \cite{dps} to our situation, where the above map
$q$ is a conic bundle over $X \times C$.  \index{notations}{X@$X$}
\index{terms}{del Pezzo surface!quartic}
\index{terms}{purely imaginary condition}

\noindent
Another (non-parabolic) instance of {\sf GLC} with $3$-dimensional
Hitchin fibers \index{terms}{Hitchin!fiber} arises when
$\lan{G}=SL_{2}(\mathbb{C})$ or $GL_{2}(\mathbb{C})$, $C$
\index{terms}{C@$C$} is a smooth curve of genus $2$, and $\Par_{C}$ is
empty. Our results in that case should appear soon\footnote{On an
  appropriately renormalized time scale} in \cite{dps2}.

\subsubsection{In more detail}

\noindent
We now describe the main results of the present work in more
detail. Let $\lan{G} = GL_{2}(\mathbb{C})$ or $SL_{2}(\mathbb{C})$ and
let $\mathbb{V} = (V,{}^{V}\bF,\nabla)$ be a flat tame quasi-parabolic
$\lan{G}$-bundle on $\mathbb{P}^1$, with quasi-parabolic structure at
$5$ distinct points $\{p_1,\dots, p_5\}$.

We will assume that $\mathbb{V}$ is a general point of
$\lan{\gLoc}$ in the following sense:
\begin{itemize}
\item Assume that the residue, $\res_{p_i}(\nabla)$, is regular
  semisimple, purely imaginary and non-resonant (i.e. its eigenvalues,
  which by assumption are real, do not differ by an integer.)
\item
Let $\blambda$ be one of the two eigenvalues of $\res_{p_i}(\nabla)$
(so $\blambda$ is real).  Endow $\mathbb{V}$ with a parabolic
structure by assigning parabolic weight = -$\blambda$ to the
corresponding eigenvector.  The NAHT converts this to a parabolic
Higgs bundle $\mathbb{E} = (E,\theta,{}^{E}\bF, -\blambda)$: the
weights don't change, while the residue of the Higgs field $\theta$ at
each $p_i$ is nilpotent.  Assume that the spectral curve for
$(E,\theta)$ is smooth.  \index{terms}{purely imaginary condition}
\index{terms}{genericity condition} This is an open condition on
$\mathbb{E}$ and hence on $\mathbb{V}$.
\end{itemize}

\index{terms}{parabolic!structure}
\index{terms}{quasi-parabolic!structure}

\

\noindent
Note that the nilpotency of the residues together with the
compatibility of the Higgs field with the parabolic filtration
together are equivalent to the condition that the action of the Higgs
field on the associated graded pieces for the parabolic filtration is
nilpotent. Parabolic Higgs bundles satisfying this later condition are
called \emph{\bfseries strongly parabolic} and we will typically say
'strongly parabolic Higgs bundle' to indicate that our Higgs bundle
satisfy the nilpotence assumption.

\index{terms}{Higgs bundle!strongly parabolic}
\index{terms}{nilpotence assumption}

\

\noindent
Using this notion of genericity our main result is the following

\begin{thm} \label{thm:MAIN}
Parabolic Hecke eigen sheaves exist for generic, purely imaginary flat
tame $\lan{G}$-bundles on $\mathbb{P}^1$ with parabolic structure at $5$
distinct points, for $\lan{G}=GL_{2}(\mathbb{C})$ or
$SL_{2}(\mathbb{C})$.

\index{notations}{GL2@$GL_{2}(\mathbb{C})$}
\index{notations}{PSL2@$\mathbb{P}SL_{2}(\mathbb{C})$}

Concretely, if $\mathbb{V}$ is a generic, purely imaginary flat
tame $\lan{G} = SL_{2}(\mathbb{C})$-bundle, then
there is a Hecke eigensheaf
$\mathfrak{c}_{\blambda}(\mathcal{O}_{\mathbb{V} })$ on
the stack $\sBun$ of quasi-parabolic $\mathbb{P}GL_{2}(\mathbb{C})$-bundles on
$\mathbb{P}^1$, with parabolic structure at $\{p_1,\dots, p_5\}$ with
eigenvalue $\mathbb{V}$.  The sheaf
$\mathfrak{c}_{\blambda}(\mathcal{O}_{\mathbb{V} })$ is a
twisted $\mathcal{D}$-module on $\sBun$ with a twist determined by
$\blambda$.  Restricted to the very stable open subset of
$\sBun$,
$\mathfrak{c}_{\blambda}(\mathcal{O}_{\mathbb{V} })$ is a
flat bundle of rank $4$. The result for
$\lan{G}=GL_{2}(\mathbb{C})$ is analogous.
\end{thm}

\index{notations}{c@$\mathfrak{c}$}
\index{notations}{Vbb@$\mathbb{V}$}
\index{notations}{SL2@$SL_{2}(\mathbb{C})$}
\index{terms}{quasi-parabolic}
\index{terms}{Hecke!eigensheaf}
\index{notations}{Buns@$\sBun$}
\index{notations}{lambdat@$\widetilde{\lambda}$}
\index{terms}{Dmodule@$\mathcal{D}$-module}
\index{terms}{residue!nilpotent}
\index{notations}{res@$\res$}
\index{terms}{twisted!$\mathcal{D}$-module}
\index{terms}{twisted!differential operators}
\index{terms}{GL2@$GL_{2}(\mathbb{C})$}

\

\noindent
For a more technical discussion of the formulation of the {\sf PGLC}
in the language of TDOs, see the appendix.  In the appendix we also
discuss, following our program in \cite{dp-jdg}, the translation of
the tamely ramified geometric Langlands correspondence into a
parabolic Hecke eigensheaf problem for tamely ramified Higgs bundles
on the moduli space of bundles. \index{terms}{TDO|see{twisted
    differential operators}} Solving this parabolic Hecke eigensheaf
problem for Higgs bundles is the main objective of our paper.
This is achieved by a detailed analysis of  the geometry of the
relevant moduli spaces which we explain next.

Let $C = \mathbb{P}^{1}$ and let
$\Par_{C} = p_{1} + p_{2} + p_{3} + p_{4} + p_{5}$ be the sum of five
distinct points.  The moduli space of rank two parabolic bundles on
$(C,\Par_{C})$ depends on a set of numerical invariants - the degree
of the level zero bundle in the parabolic family and the set of
parabolic weights.  The collection of weights has a chamber structure
and the moduli space depends only on the chamber and not on the
particular collection of weights in that chamber. As a first
observation we have the following

\begin{thm}[Lemma~\ref{lem:moduli.isos},
    Theorem~\ref{theo:shapes.of.moduli}]
  \label{thm:deg0=deg1}
  There is a dominant chamber
of parabolic weights for  rank two parabolic bundles on
$(C,\Par_{C})$  such that for all weights in this chamber:
\begin{itemize}
\item every semistable parabolic bundle is
stable; 
\item the connected components of the
moduli space corresponding to different degrees are canonically
isomorphic to the $dP_{5}$ del Pezzo surface $X$ obtained by blowing up
the $5$ points $\{p_{i}\}_{i=1}^{5}$ sitting on the diagonal
conic $C \subset S^{2}C
\cong \mathbb{P}^{2}$.
\end{itemize}
\end{thm}

\index{terms}{dominant chamber!of parabolic weights}

\

\noindent
The way in which the moduli space of parabolic bundles depends on the
collection of weights is worked out in section \ref {ssec:chambers}.
We exhibit a chamber structure on the the parameter space of all
weights, and show that the isomorphism class of the moduli space is
constant in each chamber and jumps across chamber walls.  As mentioned
in Theorem~\ref{thm:deg0=deg1} for the most interesting chamber, the
moduli space of parabolic bundles is a $dP_{5}$.  In this work we fix
the parabolic weights to be in the interior of this $dP_{5}$ chamber.
We have not checked carefully what happens if we take our moduli space
to correspond to weights in a different chamber.

\

\noindent
Equivalently: the del Pezzo
$X$ can be described in its anticanonical model as the
intersection of two quadrics in $\mathbb{P}^{4}$. The parameter space
of the pencil of quadrics $\{Q_{t}\}_{t \in C}$ vanishing on $X$ is
naturally identified with $C$ and the divisor $\Par_{C}$ corresponds
to the locus of singular quadrics in the pencil. Using this one checks

\

\begin{thm}[Proposition~\ref{prop:wobbly}]  \label{thm:wobbly}
  The wobbly locus in $X$ is the union of the 16 lines \linebreak
  $L_{I} \subset X
  \subset \mathbb{P}^{4}$ which are naturally labeled by the subsets
$I \in \evenL$, where $\evenL$ is the collection of subsets in
  $\{1,2,3,4,5\}$ of even cardinality.
\end{thm}

\

\noindent
From the point of view of the anti-canonical model of $X$ the basic
Hecke correspondence parametrizing the modifications of parabolic
bundles at a single point can be compactified and resolved (see
Chapter~\ref{sec:hecke}) to the correspondence
\begin{equation} \label{eq:parHeckdiag2}
\xymatrix@R-1pc{
& H \ar[dl]_-{p} \ar[dr]^-{q} & \\
X & & X\times C
}
\end{equation}
where
\begin{itemize}
\item $H = \op{Bl}_{\coprod_{I} \widehat{L_{I}\times L_{I}}}
    \op{Bl}_{\Delta} (X\times X)$;
\item the two maps $H \to X$  correspond to
    the blow down map $H \to X\times X$ followed by the first or
    second projection; 
\item the map $H \to C$ is the resolution of the rational map $X\times
  X \dashrightarrow C$ which sends $(x,y) \in X\times X$ to the unique
  $t \in C$ such that $Q_{t} \subset \mathbb{P}^{4}$
  contains the line through the two points $x, y \in \mathbb{P}^{4}$.
\end{itemize}
Note that $H$ is smooth by construction. The general fibers of $q$ are
smooth rational curves (Hecke lines) and the general fibers of $p$ are
smooth $dP_{6}$ del Pezzo surfaces. Furthermore, as explained in
Section~\ref{ssec:pardivisors}, all spaces in the Hecke diagram
\eqref{eq:parHeckdiag2}
are naturally equipped
with normal crossings parabolic divisors:
$\Par_{C} = \sum_{i=1}^{5} p_{i}$, $\Par_{X} = \sum_{I} L_{I}$,
\linebreak 
$\Par_{X\times C} = \Par_{X}\times C + X\times \Par_{C}$, and 
$\Par_{H} = p^{*}\Par_{X} + q^{*}\Par_{X\times C}$.
This geometry provides the setup  needed to formulate the
  parabolic version of the Higgs Hecke eigensheaf problem.

 The solution to this problem for $\lan{G} = GL_{2}(\mathbb{C})$ is derived
 in Chapters \ref{sec:eigensheaves} and \ref{sec:solution} and is
 summarized in the following theorem. The case $\lan{G} =
 \mathbb{P}GL_{2}(\mathbb{C})$ follows with obvious modifications.

\

\begin{thm}[Proposition~\ref{prop:numerics}]
  \label{thm:MAIN.Higgs}
  Fix a pair of parabolic weights on $\Par_{C}$, specified by a pair
  of real vectors $a, b \in \mathbb{R}^{5}$ that are generic
  except for the linear relation $\sum_{i = 1}^{5} (a_{i}+b_{i}) = 3$
  (which, by \eqref{eq:spectral.parch}, expresses the condition that
  $\parch_{1}(E_{\bullet}) = 0$). Also, fix a pair of integral vectors
  $n_{1}, n_{2} \in \mathbb{Z}^{5}$ such that $\sum_{i = 1}^{5}
  (n_{1,i}+ n_{2,i}) = -3$.  Then
\begin{description}
\item[Hecke kernel:] The pair $(a,b)$ determines a natural parabolic
  line bundle $\mycal{I}_{\bullet} =
  \mathcal{O}_{H}(\bzeta\Par_{H})_{\bullet}$ and $\mycal{I}_{\bullet}' =
  \mathcal{O}_{H}(\bzeta'\Par_{H})_{\bullet}$ on $(H,\Par_{H})$ with
  $\parch_{1}(\mycal{I}_{\bullet}) = \parch_{1}(\mycal{I}_{\bullet}')  = 0$.
\item[Hecke eigensheaf:] For any tame strongly parabolic rank two
  Higgs bundle $(E_{\bullet},\theta)$ on \linebreak
  $(C,\Par_{C})$ with weights
  specified by $(a,b)$ and a smooth spectral cover $\sC$, there exists
  a natural Hecke eigensheaf on the (disconnected) moduli
  space. Explicitly, this is specified by a pair of objects
  $(F_{\bullet},\varphi)$ and $(F_{\bullet}',\varphi')$ on the del
  Pezzo $X$, so that
  \begin{itemize}
  \item $(F_{\bullet},\varphi)$ and $(F_{\bullet}',\varphi')$ are
    stable, strongly parabolic, rank
    four Higgs bundles on $(X,\Par_{X})$ with $\parch_{1} = 0$,
    $\parch_{2} = 0$.
 \item The parabolic Hecke eigensheaf conditions
\[
\begin{aligned}
q_{*} (p^{*}(F_{\bullet},\varphi)\otimes (\mycal{I}_{\bullet},0)) & = 
p_{X}^{*}(F_{\bullet}',\varphi')\otimes p_{C}^{*}(E_{\bullet},\theta) \\
q_{*} (p^{*}(F_{\bullet}',\varphi')\otimes (\mycal{I}_{\bullet}',0)) & = 
p_{X}^{*}(F_{\bullet},\varphi)\otimes p_{C}^{*}(E_{\bullet},\theta) 
\end{aligned}
\]
hold.
\item The parabolic structures on $(F_{\bullet},\varphi)$ and
  $(F_{\bullet}',\varphi')$ each have two jumps, of rank one and three,
  on every component $L_{I}$ of $\Par_{X}$.
\item The parabolic weights of $(F_{\bullet},\varphi)$ and
  $(F_{\bullet}',\varphi')$ are specified by vectors $e, d, e', d' \in
  \mathbb{R}^{\evenL} = \mathbb{R}^{16}$ which can be taken to be
  \[
\left| \; 
    \begin{aligned}
e & = \frac{5}{8}\sigma - 3\acan +\frac{3}{8}(\unone  - 2\incdelta)
(r_{1} - r_{2}), \\[+0.5pc]
d & = \frac{1}{8}\sigma + \acan
- \frac{1}{8}(\unone - 2\incdelta)(r_{1} - r_{2}), \\[0.5pc]
e' & = \frac{5}{8}\sigma  -
3\acan - \frac{3}{8}(\unone  - 2\incdelta)
(r_{1} - r_{2}), \\[+0.5pc]
d' & = \frac{1}{8}\sigma + \acan
      + \frac{1}{8}(\unone - 2\incdelta)(r_{1} - r_{2}).
\end{aligned}
\right.
    \]
    where $\mathfrak{u} = 2\op{pr}_{\{5\}}(\delta_{L_{\varnothing}})$,
    $r_{1} = a + n_{1}$ and $r_{2} = b + n_{2}$, $\sigma \in
    \mathbb{R}^{\evenL}$ is a vector with all its coordinates equal to
    $1$, and $\unone$ and $\incdelta$ are $16\times 5$ matrices whose
    entries corresponding to $i \in \{1,2,3,4,5\}$ and $I \in \evenL$
    are given by $\unone_{Ii} = 1$ for all $I,i$, while
    $\incdelta_{Ii} = 1$ if $i \in I$ and $\delta_{Ii} = 0$ if $i
    \notin I$.
\item Both $(F_{\bullet},\varphi)$ and $(F_{\bullet}',\varphi')$
  correspond to spectral data specified on a spectral cover of $X$
  which is birational to the Hitchin fiber through
  $(E_{\bullet},\theta)$.
\end{itemize}
\end{description}   
\end{thm}

\

\noindent
As explained briefly above, and in detail in
Appendix~\ref{app-ram.glc}, Theorem~\ref{thm:MAIN} follows immediately
from Theorem~\ref{thm:MAIN.Higgs} by conjugating the assignment
$(E_{\bullet},\theta) \mapsto (F_{\bullet},\varphi)$ by two
non-abelian Hodge correspondences (one on $(C,\Par_{C})$ and one on
$(X,\Par_{X})$) and by taking a twisted Deligne-Goresky-MacPherson
extension
\index{terms}{DeligneGG@Deligne-Goresky-MacPherson!extension}
of the flat bundle corresponding to $(F_{\bullet},\varphi)$
from $X - \Par_{X}$ to the moduli stack $\sBun$.  Under the
non-abelian Hodge correspondence the assignment on parabolic
parameters $(a,b) \mapsto (e,d)$ becomes the Okamoto map which
converts the spectrum of the residue of the eigenvalue flat bundle
$\mathbb{V}$ into the twisting parameters for the eigensheaf twisted
$\mathcal{D}$-module $\mathfrak{c}\left(\mathcal{O}_{\mathbb{V}}\right)$. In
Section~\ref{ssec:okamoto} we carry out this translation and compute
the Okamoto map in our setting. \index{terms}{Okamoto!map}

\

\noindent
Curiously the non-abelian Hodge theory constraints on the modular
spectral data turn out to be sufficiently restrictive to rigidify the
parabolic eigensheaf problem. To illustrate this, in
Claim~\ref{claim:spectral.is.eigensheaf} we show that for any choice
of parabolic weights on $(X,\Par_{X})$ producing modular spectral data
whose Higgs bundle has vanishing $\parch_{1}$ and $\parch_{2}$ we can
find a companion Hecke eigenvalue. We are very grateful to the
anonymous referee for suggesting that such a statement might be
true. The analysis of this question lead us to a clearer and
streamlined version of the numerical equations producing the Hecke
eigensheaf in the first place.

\subsubsection{Organization}

This work is organized as follows. Parabolic objects (bundles, Higgs
bundles, flat bundles)
are reviewed in Chapter
\ref{sec:parabolic.objects}. In particular, section
\ref{ssec:operations} discusses the natural operations one can perform
on parabolic bundles or parabolic Higgs bundles: pullbacks
\cite{mochizuki-kh1,iyer-simpson-dr,iyer-simpson}, tensor products
\cite{mochizuki-kh1,iyer-simpson-dr,iyer-simpson}, and pushforwards
\cite{dps}.  Notions of stability or semistability are discussed in
section \ref{ssec:stability}.  In particular, subsection
\ref{sssec:ss.general} includes a statement of Mochizuki's non-abelian
Hodge theorem, while subsection \ref{sssec:ss.P1} begins the detailed
analysis we need for the case of rank $2$ parabolic bundles on
$\mathbb{P}^1$. \index{terms}{natural operations}
\index{terms}{parabolic!objects}
\index{terms}{parabolic!bundle}
\index{terms}{parabolic!Higgs bundle}
\index{terms}{parabolic!flat bundle}

Chapter \ref{sec:moduli.spaces} studies the quasi projective moduli
spaces determined by the above notions of stability or semistability,
focusing on $\mathbb{P}^1$ with parabolic structure at $5$ points
$\{p_1,\dots, p_5\}$. We work out the moduli spaces produced by
GIT. The main object that we will work with in the rest of the paper
turns out to be a del Pezzo quartic surface $X \subset
\mathbb{P}^{4}$.
\index{terms}{del Pezzo surface!quartic}
\index{notations}{X@$X$} \index{terms}{stability}
\index{terms}{semistability}
\index{terms}{parabolic!structure}

In Chapter \ref{sec:hecke} we construct explicit models for the
version of the Hecke correspondences relating the GIT moduli spaces in
our case.  Since a Hecke transform of a stable bundle need not be
stable, the naive Hecke correspondence
\index{terms}{Hecke!correspondence} \index{terms}{parabolic!Hecke
  correspondence} needs to be compactified and resolved.  We will see
that it is obtained from the product $X \times X$ by blowing up the
diagonal $X$, and then blowing up the proper transforms of the $16$
surfaces $L \times L$, where $L$ runs through the $16$ lines in the
del Pezzo surface $X$.  In particular, the Hecke correspondence in our
case can be described as a conic bundle over $X \times
C$. \index{terms}{line!on a quartic del Pezzo} \index{terms}{blow-up}

In Chapter \ref{sec:modular} we give a detailed description of the
rational map from the Hitchin fiber to the moduli of bundles.
\index{terms}{rational!map} \index{terms}{Hitchin!fiber}
\index{terms}{moduli space!of parabolic bundles} We show that the base
locus of this map consists of the $16$ points of order two on each
component of the Hitchin fiber and that the blow-up of these points
resolves the map and realizes a blow-up of (a connected component of)
the Hitchin fiber as a finite  degree four cover of the moduli space $X$ of
parabolic bundles. We also show that the blown-up Hitchin fiber is
naturally embedded in the bundle of logarithmic one forms on $X$ with
poles along the union of the sixteen $(-1)$-curves on $X$. Thus the
blown-up Hitchin fiber is a spectral cover for a meromorphic Higgs
bundle on $X$. We will refer to this cover as the {\em modular
  spectral cover}.  In section \ref{sssec:wobbly} we analyze the
wobbly locus in our situation and show that it is the reduced divisor
consisting of the 16 lines in $X$. \index{terms}{spectral
  cover!modular} \index{terms}{wobbly!locus}

Non-abelian Hodge theory converts the question of constructing the
geometric Langlands correspondence for tamely ramified flat bundles
on $C$ to the problem of solving an eigensheaf problem for the Hecke
action on tamely ramified semistable parabolic Higgs bundles on
$X$. To spell this problem out we need to first understand the
configuration of parabolic divisors in our birational model for the
parabolic Hecke correspondence, as well as the relevant Higgs data.
This is done in detail in Chapter \ref{sec:eigensheaves}.

After this preliminary work we are ready to solve all constraints on
the desired eigensheaf so that both Mochizuki's non-abelian Hodge
theory conditions and the Hecke eigensheaf condition are
satisfied. This is carried out in Chapter \ref{sec:solution}, where
the solution, unique up to the obvious symmetries, is found
explicitly.

\bigskip

\

\noindent {\bf Acknowledgments.} We owe special thanks to Dima
Arinkin, Carlos Simpson and Edward Witten for collaboration on related
projects, and for many illuminating explanations and comments on the
manuscript and related issues. We would also like to thank Roman
Bezrukavnikov, Pierre Deligne, Dennis Gaitsgory, Tamas Hausel, Ludmil
Katzarkov, and Ana Peon-Nieto for several interesting discussions.
Finally we would also like to thank the anonymous referee for the many
insightful comments and corrections on the first version of this
manuscript. In particular their suggestion that the non-abelian Hodge
theory constraints are already sufficiently strong to enforce
automorphic behavior prompted us to clarify and simplify the solution
to the abelianized version of the Hecke eigensheaf problem.

During the preparation of this work Ron Donagi was supported in part
by NSF grants DMS 1304962, DMS 1603526, DMS 2001673, and by Simons HMS
Collaboration grant \# 390287. Tony Pantev was supported in part by
NSF grant RTG 0636606, DMS 1302242, DMS 1601438, DMS 1901876, and by
Simons HMS Collaboration grant \# 347070. At one stage the authors
visited the Laboratory of Mirror Symmetry NRU HSE and were supported
through RF Government grant, ag. No 14.641.31.0001.

\section{Parabolic objects} \label{sec:parabolic.objects}

In this first chapter we introduce our main objects of study - parabolic
bundles and tamely ramified parabolic flat bundles and parabolic Higgs
bundles. We also recall their basic properties and fundamental
invariants. \index{terms}{parabolic!objects}

\subsection{Parabolic bundles and Higgs bundles}
\label{ssec:par.vb}

The most basic notion we will deal with is the notion of a parabolic
bundle. We will follow the treatment adopted in
\cite{carlos-nc,mochizuki-kh1,mochizuki-kh2,mochizuki-D1,mochizuki-D2}
and will focus on parabolic vector bundles rather than the more
general parabolic principal bundles.

\index{terms}{parabolic!vector bundles}
\index{terms}{parabolic!principal bundle}
\index{terms}{parabolic!Higgs bundle}

\subsubsection{Parabolic bundles} \label{sssec:par.vb} 
Fix a pair $(X,\Par_{X})$, where $X$ is a smooth complex projective
variety, $\Par_{X} \subset X$ is a reduced divisor with simple normal
crossings, and ${\displaystyle \Par_{X} = \sum_{i \in S} D_{i}}$ are
the irreducible components of $\Par_{X}$. With this notation we have
the following

 \begin{defi} \label{defi:parabolic.vb}
A {\bfseries parabolic vector bundle on $\left(X,\sum_{i \in S}
  D_{i}\right)$} is a collection of locally free coherent sheaves
$\boldsymbol{E}_{\bullet} =
\left\{ \boldsymbol{E}_{t} \right\}_{t \in \mathbb{R}^{S}}$ together
with inclusions $\boldsymbol{E}_{t} \subset \boldsymbol{E}_{s}$ of
sheaves of $\mathcal{O}_{X}$-modules, specified for all $t \leq s$,
and satisfying the conditions:
\begin{description}
\item[{\bfseries [semicontinuity]}]{ for every
  $t \in  \mathbb{R}^{S}$, there exists a real number $c
  > 0$ so that  $\boldsymbol{E}_{t +
  \boldsymbol{\varepsilon}} = \boldsymbol{E}_{t}$ for all
  functions $\boldsymbol{\varepsilon} : S \to [0,c]$.}
\item[{\bfseries [support]}]{ if
  $\boldsymbol{\delta}_{i} : S 
  \to \mathbb{R}$ is
  the  characteristic function of $i \in S$, then for all
  $t \in \mathbb{R}^{S}$ we have
  $\boldsymbol{E}_{t + \boldsymbol{\delta}_{i}} =
  \boldsymbol{E}_{t}(D_{i})$ (compatibly with the inclusion).}
\end{description}
The real vectors $t \in \mathbb{R}^{S}$ are  called  {\bfseries
  parabolic levels}, and the vector bundles $\boldsymbol{E}_{t}$ are
called {\bfseries level bundles} or {\bfseries level slices} of the
parabolic bundle $\boldsymbol{E}_{\bullet}$. 
\end{defi}

\

\noindent
While dealing with abelianization and the formalism of pullbacks and
pushforwards we will sometimes need to work with parabolic bundles for
which the parabolic structure is specified along non-reduced normal
crossing divisors. These are defined in a similar manner. To take into
account the non-reduced structure the parabolic data on non-reduced
normal crossings divisors must incorporate a choice of partitions of
the multiplicities of the divisor components.

\index{terms}{parabolic!vector bundles}
\index{terms}{parabolic!structure!along a non-reduced divisor}

Fix a pair $(X,\Par_{X})$, where $X$ is a smooth complex projective
variety, $\Par_{X} \subset X$ is a possibly non-reduced effective
divisor whose reduction has simple normal crossings, and
${\displaystyle \Par_{X} = \sum_{i \in S} m_{i}D_{i}}$, where the
$D_{i}$ are irreducible components of $\Par_{X}$. For every $i \in S$
fix a partition $m_{i} = \sum_{j \in \mathsf{J}(i)} m_{i}^{j}$ of the
multiplicity $m_{i}$, with a finite set $\mathsf{J}(i)$ labeling the
parts of the partition.  We will write $\Par_{X} = \sum_{\substack{i
    \in S\\ j \in \mathsf{J}(i)}} m_{i}^{j} D_{i}$ for the associated
partitioning of the parabolic divisor. With this notation we now
define:

\

 \begin{defi} \label{defi:nr.parabolic.vb}
A {\bfseries parabolic vector bundle on
  $\left(X,\sum_{\substack{i
      \in S\\ j \in \mathsf{J}(i)}} m_{i}^{j} D_{i}\right)$}
is a collection of locally free coherent sheaves
  $\left\{ \boldsymbol{E}_{t}
  \right\}_{t \in \prod_{i \in S} \mathbb{R}^{\mathsf{J}(i)}}$ together with
  inclusions $\boldsymbol{E}_{t} \subset
  \boldsymbol{E}_{s}$ of 
  sheaves of $\mathcal{O}_{X}$-modules, specified for all $t \leq
  s$, and satisfying the conditions:
\begin{description}
\item[{\bfseries [semicontinuity]}]{ for every
  $t \in  \prod_{i \in S} \mathbb{R}^{\mathsf{J}(i)}$, there exists a
  real number $c
  > 0$ so that  $\boldsymbol{E}_{t +
  \boldsymbol{\varepsilon}} = \boldsymbol{E}_{t}$ for all
  functions $\boldsymbol{\varepsilon} :
  \sqcup_{i \in S} \mathsf{J}(i) \to [0,c]$.}
\item[{\bfseries [support]}]{ if
  $\boldsymbol{\delta}_{i}^{j} : \sqcup_{i \in S} \mathsf{J}(i) 
  \to \mathbb{R}$ is
  the  characteristic function of $j \in \mathsf{J}(i)$, then for all
  $t \in \prod_{i \in S} \mathbb{R}^{\mathsf{J}(i)}$ we have
  $\boldsymbol{E}_{t + \boldsymbol{\delta}_{i}} =
  \boldsymbol{E}_{t}(m_{i}^{j}D_{i})$ (compatibly with the inclusion).}
\end{description}
\end{defi}

\

\begin{rem} \label{rem:nonreduced.nonseparated}
  {\bfseries (i)} \ Note that Definition~\ref{defi:parabolic.vb} is a
  special case of Definition~\ref{defi:nr.parabolic.vb}. Indeed, in
  the case when all multiplicities $m_{i}$ are equal to $1$ and the
  data of the partitions is vacuous, the data and properties in
  Definition~\ref{defi:nr.parabolic.vb} become the same as the data
  and properties in Definition~\ref{defi:parabolic.vb}. In this case
  we will say that $\boldsymbol{E}_{\bullet}$ is a parabolic bundle on
  $(X,\Par_{X})$.

\

\noindent
{\bfseries (ii)} \ Alternatively the partition data in
Definition~\ref{defi:nr.parabolic.vb} can be encoded in a slight
modification of the geometry.  Indeed, suppose we are in the setting
when $\Par_{X}$ is non-reduced and we have partitioned the
multiplicities of the components as above. Consider the nonseparated
scheme $\mathsf{Z}$ obtained by gluing copies of $X$ along $X^{o} := X
- \cup_{i} D_{i}$ so that in $\mathsf{Z}$ the component $D_{i}$ is
replaced by a collection of disjoint non-separated components labeled
by the elements in $\mathsf{J}(i)$. In other words
\[
\mathsf{Z} = \op{colim}\left[\xymatrix@1{X^{o} \ar[r] &
  \bigsqcup_{\substack{i \in S \\ j \in \mathsf{J}(i)}} X\times \{j\}} \right],
\]
where the coproduct is taken in the category of schemes.  By
construction $\mathsf{Z}$ contains a natural divisor
$\Par_{\mathbf{Z}} = \sum_{\substack{i \in S \\ j \in \mathsf{J}(i)}}
m_{i}^{j} \mathsf{D}_{i}^{j}$, where the $\mathsf{D}_{i}^{j}$
corresponds to $D_{i}\times\{j\} \subset X\times \{j\}$ and the
$D_{i}^{j}$'s are the irreducible components of $\Par_{\mathbf{Z}}$.

Note that the data of a parabolic vector bundle
$\boldsymbol{E}_{\bullet}$ on $X$ with parabolis structure on the
non-reduced partitioned divisor $\sum_{i,j} m_{i}^{j} D_{i}$ defined
as in Definition~\ref{defi:nr.parabolic.vb} is the same thing as a
collection of nested locally free sheaves on $\mathsf{Z}$ with the
usual semicontinuity condition specified on $\sum_{i,j}
\mathsf{D}_{i}^{j}$ and with the suport condition as in
Definition~\ref{defi:nr.parabolic.vb}, i.e. a jump is given by
tensoring with $\mathcal{O}_{Z}(m_{i}^{j} \mathsf{D}_{i}^{j})$. In
other words, passing to the non-separated $\mathsf{Z}$ allows us to
formulate Definition~\ref{defi:nr.parabolic.vb} only in terms of
multiplicities and without any reference to partitions.

Note also that this interpretation of
Definition~\ref{defi:nr.parabolic.vb} in terms of the non-separated
scheme $\mathsf{Z}$ is consistent with the non-separated versions of
the moduli of bundles that naturally appear in the tamely ramified GLC
(see Appendix~\ref{app-ram.glc}).
\end{rem}

\index{terms}{parabolic!vector bundles}
\index{terms}{parabolic!structure!along a non-reduced divisor}

\

\subsubsection{Flags and weights} \label{sssec:flags.and.weights}

Assume again that $X$ is smooth projective and
$\Par_{X} \subset X$ is a reduced  simple normal crossings
divisor.  Fix a parabolic bundle $\boldsymbol{E}_{\bullet}$ on $(X,\Par_{X})$
and a vector $\boldsymbol{c} \in \mathbb{R}^{S}$.

For every $i \in S$ we get an induced filtration
  $\{{}^{i}F_{a}\}_{c_{i} - 1 < a \leq c_{i}}$ of the
  restricted sheaf 
  $\boldsymbol{E}_{\boldsymbol{c}|D_{i}}$:
\[
   {}^{i}F_{a} = \bigcup_{\substack{\boldsymbol{t} \leq
    \boldsymbol{c} 
    \\ \alpha_{i} \leq a}} \boldsymbol{E}_{t}
\]
Define
\[
{}^{i}\!\op{gr}_{a}
  \boldsymbol{E}_{\boldsymbol{c}} :=
             {}^{i}F_{a}/{}^{i}F_{<a}.
\]
Now the {\bfseries  [semicontinuity]} property in
Definition~\ref{defi:parabolic.vb} implies that 
  the {\bfseries set of parabolic weights}
\[
\weights(\boldsymbol{E}_{\boldsymbol{c}},i) = \left\{ a \in
(c_{i}-1,c_{i}] \ \left| \ {}^{i}\!\op{gr}_{a} \neq 0\right.\right\}
\]
is finite.

Note that the single bundle $\boldsymbol{E}_{\boldsymbol{c}}$ together
with the flags $\{{}^{i}F_{a} | \ i \in S, a \in
\weights(\boldsymbol{E}_{\boldsymbol{c}},i)\}$ reconstruct the
parabolic bundle $\boldsymbol{E}_{\bullet}$ and so we can equivalently
describe parabolic bundles as bundles furnished with flags by
subbundles along the components of the parabolic divisor and with a
collection of parabolic weights labeling each step of each flag.

\

\subsubsection{Locally abelian parabolic
  bundles} \label{sssec:locally.abelian} 

The first and most basic examples of parabolic bundles are the
parabolic line bundles. In parallel to the usual theory parabolic line
bundles are given by $\mathbb{R}$-divisors:

\index{terms}{parabolic!line bundle}
\index{terms}{parabolic!bundle!locally abelian}

\begin{ex} \label{ex:parabolic.lb}
  A {\bfseries parabolic line bundle} on $(X,\sum_{i} D_{i})$
  is a parabolic
  vector bundle  $F_{\bullet}$ for which all sheaves $F_{t}$
  are invertible.  

Given $\boldsymbol{a} \in \mathbb{R}^{S}$ we define a
  parabolic line bundle $\mathcal{O}_{X}\left( \sum_{i
  \in S} \boldsymbol{a}_{i}D_{i}\right)_{\bullet}$ by setting 
\[
\left(\mathcal{O}_{X}\left( \sum_{i
  \in S} \boldsymbol{a}_{i}D_{i}\right)\right)_{t}
  := \mathcal{O}_{X}\left( \sum_{i
  \in S} \fl*{\boldsymbol{a}_{i}+t_{i}} D_{i}\right)
\]
Tautologically, every parabolic
line bundle $F_{\bullet}$ is isomorphic to $L\otimes
\mathcal{O}_{X}\left( \sum_{i \in S}
\boldsymbol{a}_{i}D_{i}\right)_{\bullet}$ for some $L \in
\Pic(X)$, and some $\boldsymbol{a} \in \mathbb{R}^{S}$.
\end{ex}

\

\noindent
For geometric considerations and crucially for the purposes of the
non-abelian Hodge correspondence it is useful to restrict attention to 
special types of parabolic bundles - Mochizuki's locally abelian
parabolic bundles \cite{mochizuki-kh1,mochizuki-kh2}.
\index{terms}{Non-abelian Hodge theory}

\begin{defi} \label{defi:locally.abelian}
  A parabolic bundle 
$F_{\bullet}$ on $(X,\sum_{i}D_{i})$  is a {\bfseries locally abelian
bundle}, if in a Zariski neighborhood of any point $x \in  X$ there is an
isomorphism between $F_{\bullet}$ and a direct sum of parabolic line
bundles.
\end{defi}

\

\begin{rem} \label{rem:locally.abelian}
  A parabolic bundle $\left({\boldsymbol{E}}_{\boldsymbol{c}}, \{
  {}^{i}F_{\bullet} \}_{i \in S}\right)$  specified by parabolic filtrations and
weights is locally abelian if and only if  on
  every intersection $D_{i_{1}}\cap \cdots \cap D_{i_{k}}$ the
  iterated graded ${}^{i_{1}}\!\op{gr}_{a_{1}} \cdots
  {}^{i_{k}}\!\op{gr}_{a_{k}} {\boldsymbol{E}}_{\boldsymbol{c}}$ does
  not depend on the order of the components.
\end{rem}

We can define
  similarly  locally abelian parabolic flat bundles, Higgs bundles,
  or more generally locally abelian parabolic $z$-connections. 

\subsubsection{Locally abelian parabolic
	$z$-connections} \label{sssec:loc.ab.lambda.connections}

Recall Deligne's notion \cite{carlos-naht,carlos-hf} of a $z$-connection:
      
\begin{defi} \label{defi:loc.ab.lambda.connections} Let
  $z \in \mathbb{C}$. A {\bfseries $z$-connection with tame
    ramification along a normal crossings divisor $D \subset X$}, is a
  pair $(E,\nabla^{z})$, where:
\begin{itemize}
\item {$E$ is an algebraic   vector bundle on $X$;}
\item {$\nabla^{z} : E \to E\otimes \Omega^{1}_{X}(\log
  D)$, is a $\mathbb{C}$-linear map 
  satisfying the $z$-twisted Leibnitz rule
$\nabla^{z}(f\cdot s) = f \nabla^{z}s + z
  s \otimes df$.
}
\end{itemize}
We say that $\nabla^{z}$ is flat if
$\nabla^{z}\circ \nabla^{z} = 0$.
\end{defi}
This notion naturally interpolates between flat bundles and Higgs
bundles. Indeed, by definition we have 
\[
\begin{aligned}
(\text{flat
  $(z=1)$-connection}) &  =  (\text{flat connection with regular
  singularities on }
D), \\ 
(\text{flat $(z=0)$-connection}) & =  (\text{Higgs
  bundle with logarithmic poles on } D),
\end{aligned}
\]
Moreover, for every $z \neq 0$ and a flat $z$-connection
$\nabla^{z}$, we have that $z^{-1}\nabla^{z}$ is a flat
$1$-connection.

The notion of a $z$-connection also interacts in a straightforward manner
with parabolic structures:

\index{terms}{tame parabolic!$z$-connection}
\index{terms}{strongly parabolic!$z$-connection}

\begin{defi} \label{defi:tame.parabolic.z.connection}
Let $X$ be a smooth projective variety and $\Par_{X} \subset X$ be a
reduced strict normal crossings divisor. A {\bfseries tame parabolic
  $z$-connection on $(X,\Par_{X})$} is a
pair $\left(\boldsymbol{E}_{\bullet},\nabla^{z}_{\bullet}\right)$, where
\begin{itemize}
\item {$\boldsymbol{E}_{\bullet}$ is a parabolic bundle on $(X,\Par_{X})$;}
\item {$\nabla^{z}_{t} : \boldsymbol{E}_{t} \to
  \boldsymbol{E}_{t}\otimes \Omega^{1}_{X}(\log \Par_{X})$ is a
  tame flat $z$-connection specified for all
  $t \in \mathbb{R}^{S}$ (compatibly with
  the inclusions).}
\end{itemize}
A  tame parabolic $z$-connection
  $\left(\boldsymbol{E}_{\bullet},\nabla^{z}\right)$ is {\bfseries
      locally abelian}  if
  the underlying parabolic bundle $\boldsymbol{E}_{\bullet}$
is locally abelian. It is  {\bfseries strongly
    parabolic}
if the
action of the residue of $\nabla^{z}_{t}$ on the associated
graded for the parabolic filtration is zero.
\end{defi}
In particular we get notions of tame parabolic (or locally abelian parabolic,
or strongly parabolic) flat bundles or Higgs bundles on $(X,\Par_{X})$.

\index{terms}{parabolic!structure}

\subsection{Parabolic Chern classes} \label{ssec:par.ch}

In this section we review the Iyer-Simpson formula for the parabolic
Chern character of a locally abelian parabolic vector bundle
\index{terms}{parabolic!Chern classes}
\index{terms}{parabolic!Chern character}

Let $\boldsymbol{E}_{\bullet}$ be a locally abelian parabolic vector
bundle on $\left(X,\sum_{i} D_{i}\right)$.  The {\bfseries parabolic Chern
  character} of $\boldsymbol{E}_{\bullet}$ can be defined
\cite{iyer-simpson} in the usual way by using the splitting principle
and the normalization that for a parabolic line bundle given by a real
divisor we have
$c_{1}(\mathcal{O}_{X}(\sum_{i} a_{i}D_{i})_{\bullet}) =
\sum_{i}a_{i}[D_{i}] \in H^{2}(X,\mathbb{R})$. Explicitly the
parabolic Chern character can be computed \cite{iyer-simpson} by
the {\bfseries
  Iyer-Simpson formula:}
\begin{equation} \label{Iyer-Simpson formula}
  \boxed{ \parch(\boldsymbol{E}_{\bullet}) =
    \parch({}_{\boldsymbol{c}}\boldsymbol{E}) = \frac{ \prod_{i\in
        S} \int_{c_{i}}^{c_{i}+1} dt_{i}
      \left[ch\left(\boldsymbol{E}_{t_{i}}\right)e^{-\sum_{i \in S}
          t_{i} D_{i}}\right]}{\prod_{i\in S} \int_{0}^{1} dt_{i}
      e^{-\sum_{i \in S} t_{i} D_{i}}},}
\end{equation}
where $\boldsymbol{c} \in \mathbb{R}^{S}$ is any base point.

\

\begin{rem} \label{rem:iyer-simpson} Given
  $\boldsymbol{c} \in \mathbb{R}^{S}$ we define the {\bfseries
    $\boldsymbol{c}$-truncation ${}_{\boldsymbol{c}}\boldsymbol{E}$ of
    $\boldsymbol{E}_{\bullet}$} to be the collection
  $\{ \boldsymbol{E}_{t} \}_{\boldsymbol{c} < t \leq \boldsymbol{c}
    +\boldsymbol{\delta}}$, with
  $\boldsymbol{\delta} = \sum_{i \in S} \boldsymbol{\delta}_{i}$.  The
  {\bfseries [support]} condition in the definition of a parabolic
  bundle implies that $\boldsymbol{E}_{\bullet}$ can be reconstructed
  effectively from any truncation ${}_{\boldsymbol{c}}\boldsymbol{E}$.
  Notice that the numerator of the Iyer-Simpson formula is independent
  of the choice of truncation.
\end{rem}

\

\noindent
For future reference note that we have an easy formula for the first
parabolic Chern class in terms of parabolic filtrations and weights:

\index{terms}{parabolic!first Chern class}

\begin{ex} \label{ex:parch1}
The first
  parabolic Chern class of $\boldsymbol{E}_{\bullet}$ is given by:
\[
\parch_{1}(\boldsymbol{E}_{\bullet}) = c_{1}(\boldsymbol{E}_{\boldsymbol{c}}) - \sum_{i \in
  S} \left(\sum_{a \in \weights(\boldsymbol{E}_{\boldsymbol{c}},i)}
   a\op{rank} {}^{i}\!\op{gr}_{a}
   \boldsymbol{E}_{\boldsymbol{c}}\right)\cdot [D_{i}].
   \]
   As usual here the right hand side is viewed as an element in
   $H^{2}(X,\mathbb{Z})$ 
\end{ex}

\subsection{Natural operations}
\label{ssec:operations}

In this section we review, following the work of various authors, the natural operations one can perform on
parabolic bundles or parabolic Higgs bundles: 
pullbacks \cite{mochizuki-kh1,iyer-simpson-dr,iyer-simpson}, 
tensor products \cite{mochizuki-kh1,iyer-simpson-dr,iyer-simpson},
and pushforwards \cite{dps}.

\subsubsection{Pullbacks} \label{sssec:pullbacks}

Suppose $f : Y \to X$ is a morphism of smooth projective varieties and
let $\Par_{Y} \subset Y$ and $\Par_{X} \subset X$ are reduced simple
normal crossings divisors such that $f^{*}\Par_{X}$ is a subdivisor of
$\Par_{Y}$. In other words we have
$\Par_{Y} = f^{*}\Par_{X} + \sum_{k \in \mathsf{K}}
\mathfrak{D}_{k}$. Let $\Par_{X} = \sum_{i \in S} D_{i}$ and let
$f^{*}D_{i} = \sum_{j \in \mathsf{J}(i)} \mathfrak{D}_{i}^{j}$ denote
the components of the pullback of $D_{i}$. Write
$\mathsf{J} = \sqcup_{i \in S} \mathsf{J}(i)$. Then the set
$\mathsf{J} \sqcup \mathsf{K}$ labels the components of
$\Par_{Y}$ and  a parabolic level
$u \in \mathbb{R}^{\mathsf{J}}\times \mathbb{R}^{\mathsf{K}}$ on
$(Y,\Par_{Y})$ is  given by a
collection of real numbers
$\left(\prod_{i \in S} (u_{i}^{j})_{j \in \mathsf{J}(i)}\right)\times
(u_{k})_{k \in \mathsf{K}}$.

The pullback of parabolic bundles from $X$ to $Y$ is defined
\cite{mochizuki-kh1,iyer-simpson-dr,iyer-simpson} so that on parabolic
line bundles it is given by the pullback of real divisors, it respects
the locally abelian condition, and in the case of rational weights is
compatible with the usual pullback of bundles on root stacks. These
properties characterize the pullback uniquely and we can define it
again by using the splitting principle. 

We can also give an explicit description of the pullback in terms of
the pullbacks of the various level bundles of a parabolic bundle. Indeed from
the above characterization it follows that the pullback of a parabolic
bundle will have trivial parabolic structure
\index{terms}{parabolic!structure}
 on the complementary
divisor
$\sum_{k \in \mathsf{K}} \mathfrak{D}_{k} = \Par_{Y} - f^{*}\Par_{X}$
and for parabolic levels that take constant values on each
$\mathsf{J}(i)$ will be just the sheaf theoretic pullback of the
corresponding level bundle from $X$. This gives all level bundles of the
pulled back parabolic bundle for which the level values are constant
on the preimages of the components of $\Par_{X}$. But this special  family of
level bundles uniquely reconstructs the full pull back parabolic
bundle.  The key to the reconstruction is the following  straightforward
remark.

\begin{rem} \label{rem:reconstruct}
Let $\mathsf{Z}$ be a smooth projective variety and let
$\Par_{\mathsf{Z}} = \sum_{i=1}^{s} \mathsf{D}_{i}$ be a reduced
strict normal crossings divisor. If $\boldsymbol{E}_{\bullet}$ is a
parabolic bundle on $(\mathsf{Z},\Par_{\mathsf{Z}})$, then all level
bundles $\boldsymbol{E}_{(t_{1}, \ldots, t_{s})}$ of
$\boldsymbol{E}_{\bullet}$ can be reconstruucted from the 'diagonal'
level bundles $\boldsymbol{E}_{(t, \ldots, t)}$, $t \in
\mathbb{R}^{s}$. Indeed we have:
\[
\boldsymbol{E}_{(t_{1}, \ldots, t_{s})} = \bigcap_{i =1}^{s} \left[
  \boldsymbol{E}_{(t_{i}, \ldots, t_{i})}\otimes
  \mathcal{O}_{Z}\left( \mathsf{M} \sum_{j \neq i}
  \mathsf{D}_{i} \right)\right],
\]
where $\mathsf{M}$ is an integer satisfying
$\mathsf{M}  \gg \max\{ t_{1}, \ldots, t_{s} \}$.
\end{rem}

\

\noindent
Using the procedure from this remark we can give an explicit
description of pull backs of parabolic bundles. Before we can spell
this out we will need to introduce some notation.

Let $\boldsymbol{\Sigma}$ denote the set of all sections
$\sigma : S \to \mathsf{J}$ of the natural projection
$\mathsf{J} = \sqcup_{i \in S} \mathsf{J}(i) \to S$.  For every
parabolic level
$u \in \mathbb{R}^{\mathsf{J}}\times \mathbb{R}^{\mathsf{K}}$ on
$(Y,\Par_{Y})$ and every $\sigma \in \boldsymbol{\Sigma}$ write
$u^{\sigma}$ for the vector
$(u_{i}^{\sigma(i)})_{i \in S} \in \mathbb{R}^{S}$.  Also, given a
parabolic level
$u \in \mathbb{R}^{\mathsf{J}}\times \mathbb{R}^{\mathsf{K}}$ and a
label $i \in S$ let $\mathsf{M}(u,i)$ be an integer satisfying
$\mathsf{M}(u,i) \gg \max\{ u_{i}^{j} \ | \ j \in \mathsf{J}(i) \}$.
With this notation we now have

\

\begin{defi} \label{defi:pullback}
  The {\bfseries pullback} of a parabolic bundle
  $\boldsymbol{E}_{\bullet}$ on $(X,\Par_{X})$ is the parabolic bundle
  $f^{*}\boldsymbol{E}_{\bullet}$ on $(Y,\Par_{Y})$ whose level bundle for
  a level $u \in \mathbb{R}^{\mathsf{J}}\times
  \mathbb{R}^{\mathsf{K}}$  is given by
\[
(f^{*}\boldsymbol{E}_{\bullet})_{u} = 
    \bigcap_{\sigma \in \boldsymbol{\Sigma}}
\left[f^{*}(\boldsymbol{E}_{u^{\sigma} })\otimes \mathcal{O}_{Y}\left(
  \sum_{i \in S}\mathsf{M}(u,i) \left( \sum_{j \in
      \mathsf{J}(i)-\sigma(i)}
    \mathfrak{D}_{i}^{j}\right)\right)\right] \otimes
\mathcal{O}_{Y}\left( \sum_{k \in
    \mathsf{K}}\fl*{u_{k}}\mathfrak{D}_{k}
\right).
\]
\end{defi}

\

\medskip

\begin{rem} \label{rem:pullback}
  \begin{enumerate}
  \item[{\bfseries (i)}] The definition is easily seen to be
    independent of the choice of the integers $\mathsf{M}(u,i)$ and
    for all practical purposes we can take these integers to be infinitely
    large.
 \item[{\bfseries (ii)}] The definition is tautologically compatible
   with locally abelian structures and the pullback of real divisors
   and thus is equivalent to the definition of pullbacks in
   \cite{iyer-simpson}. In particular this pullback is compatible with
   taking parabolic Chern classes, i.e. we have
   $\parch(f^{*}\boldsymbol{E}_{\bullet}) =
   f^{*}\parch(\boldsymbol{E}_{\bullet})$.
   \index{terms}{parabolic!Chern classes}
\end{enumerate}
\end{rem}

\

\subsubsection{Tensor products} \label{sssec:tensor}

Let $\boldsymbol{E}_{\bullet}$ and $\boldsymbol{F}_{\bullet}$ be
locally abelian parabolic vector bundles on $(X,\Par_{X})$, $\Par_{X}
= \sum_{i \in S} D_{i}$.  As explained in
\cite{mochizuki-kh1,iyer-simpson-dr,iyer-simpson} one can define a
canonical tensor product of such bundles which is uniquely
characterized by requiring that it preserves the locally abelian
property and that for parabolic line bundles it corresponds to the
addition of real divisors. Furthermore it is shown in
\cite{iyer-simpson} that for bundles with rational parabolic weights,
under the correspondence of \cite{borne-root}, the tensor product is
compatible with the usual tensor product of vector bundles on an
appropriate root stack. Again it is possible to express the tensor
product in terms of level bundles which leads to the following
equivalent definition.

\begin{defi} \label{defi:tensor} The {\bfseries tensor product} of
  the locally abelian parabolic bundles $\boldsymbol{E}_{\bullet}$ and
  $\boldsymbol{F}_{\bullet}$ is the locally abelian parabolic bundle
  $\boldsymbol{E}_{\bullet}\otimes \boldsymbol{F}_{\bullet}$ whose
  level bundle at level $t \in \mathbb{R}^{S}$ is given by
\[
\left(\boldsymbol{E}_{\bullet}\otimes \boldsymbol{F}_{\bullet}\right) =
\bigcup_{\xi \in \mathbb{R}^{S}} \boldsymbol{E}_{\xi}\otimes
\boldsymbol{F}_{t-\xi}.
\]
\end{defi}

\

\begin{rem} \label{rem:tensor}
\item[{\bfseries (i)}] Suppose $\boldsymbol{E}_{\bullet}$ is an
  arbitrary parabolic bundle but $\boldsymbol{F}_{\bullet}$ is a
  parabolic line bundle on $(X,\Par_{X})$. Then, as explained above,
\[
\boldsymbol{F}_{\bullet} = L\left(\sum_{i \in S} d_{i} D_{i}\right)_{\bullet},
\quad
\text{for some} \ (d_{i})_{i \in S} = d \in \mathbb{R}^{S}, \ \text{and some} \
L \in \Pic(X).
\]
In this case the formula for the level bundles of the tensor product
\[
\boldsymbol{E}_{\bullet}\otimes \boldsymbol{F}_{\bullet} = \boldsymbol{E}_{\bullet}
  \otimes L\left(\sum_{i \in S} d_{i} D_{i}\right)_{\bullet}
\]
simplifies to
\[
\left(\boldsymbol{E}_{\bullet}
\otimes L\left(\sum_{i \in S} d_{i} D_{i}\right)_{\bullet}\right)_{t} =
\boldsymbol{E}_{t+d}\otimes L.
\]

\

\noindent
{\bfseries (ii)} \ The compatibility of the tensor product with the
locally abelian structures implies again \cite{iyer-simpson} that the
parabolic Chern character is multiplicative  for the tensor product,
i.e.
\index{terms}{parabolic!Chern character}
\[
\parch(\boldsymbol{E}_{\bullet}\otimes \boldsymbol{F}_{\bullet}) =
  \parch(\boldsymbol{E}_{\bullet})\cdot\parch(\boldsymbol{F}_{\bullet})
\]
\end{rem}

\

\subsubsection{Pushforwards} \label{sssec:pushforwards}

Let $f : Y \to X$ be a morphism between smooth projective varieties.
Given a semistable tame parabolic Higgs bundle $(F_{\bullet},\varphi)$
we want to understand the derived pushforward of
$(F_{\bullet},\varphi)$ via $f$ in the category of Higgs sheaves. By
definition this pushforward is the complex of tame Higgs bundles which
is the specialization at zero of the $f$-pushforward (in the sense of
\cite{claude-twistor,mochizuki-D1}) of the twistor $\mathcal{D}$-module
corresponding to $(F_{\bullet},\varphi)$. In \cite{dps} we gave
the following direct algebraic description of such pushforwards.

Assume $Y$ and $X$ are equipped with simple normal crossings divisors
$\Par_{Y} \subset Y$ and $\Par_{X} \subset X$. Decompose $\Par_{Y}$ as
a sum $\Par_{Y} = \Par_{Y}^{\ver} + \Par_{Y}^{\hor}$ of a vertical and
horizontal part. That is $\Par_{Y}^{\ver}$ is the sum of all
components of $\Par_{Y}$ which map to proper subvarieties in $X$, and
$\Par_{Y}^{\hor}$ is the sum of all components that dominate
$X$. Additionally we will assume that
$\Par_{Y}^{\ver} = f^{*}\Par_{X}$ scheme theoretically, that $f$ is
smooth away from $\Par_{X}$, that $\Par_{Y}^{\hor}$ is smooth over
$X$, and that $\Par_{Y}^{\hor}\cap \Par_{Y}^{\ver}$ is contained in
the smooth locus of $\Par_{Y}^{\ver}$. Note that these conditions
imply that
$\Par_{Y}^{\hor} = \sum_{k \in \mathsf{K}} \mathfrak{D}_{k}$ is a sum
of disconnected smooth components $\mathfrak{D}_{k}$ each of which is
also smooth over $X$.

Given a parabolic level $t$ on $X$ we will write $\lift(t)$ for the
parabolic level on $Y$ which assigns $0$ to each component of
$\Par_{Y}^{\hor}$, while to each component of $\Par_{Y}^{\ver}$ it
assigns the value of the level $t$ on the image of that vertical
component under $f$. Now for each component $\mathfrak{D}_{k}$ of the
horizontal divisor we can consider the associated graded
$\op{gr}^{\mathfrak{D}_{k}}F_{\lift(t)}$ of
$F_{\lift(t)}$ with respect to the parabolic filtration along
$\mathfrak{D}_{k}$. By definition this a vector bundle on
$\mathfrak{D}_{k}$ given by
\[
\op{gr}^{\mathfrak{D}_{k}}F_{\lift(t)} =  F_{\lift(t)}/F_{\lift(t)
  - \varepsilon\bdelta_{k}},
\]
where $\bdelta_{k}$ is the characteristic function of $k$ (viewed as a
parabolic level on $Y$), and $\varepsilon > 0$ is a small real number.
Note that by construction $\op{gr}^{\mathfrak{D}_{k}}F_{\lift(t)}$
is a quotient of the vector bundle $F_{\lift(t)|\mathfrak{D}_{k}}$ and
that the residue $\res_{\mathfrak{D}_{k}} \varphi :
F_{\lift(t)|\mathfrak{D}_{k}} \to F_{\lift(t)|\mathfrak{D}_{k}}$
descends to an endomorphism  $\op{gr}-\res_{\mathfrak{D}_{k}}\varphi \in
\op{End}\left(\op{gr}^{\mathfrak{D}_{k}}F_{\lift(t)}\right)$. The nilpotent
  part of this endomorphism induces a monodromy weight filtration
  $W_{\bullet}\left(\op{gr}^{\mathfrak{D}_{k}}F_{\lift(t)}\right)$ on
  the vector bundle
  $\op{gr}^{\mathfrak{D}_{k}}F_{\lift(t)}$. The vector bundle
  $\op{gr}^{\mathfrak{D}_{k}}F_{\lift(t)}$ on $\mathfrak{D}_{k}$ can
  be viewed as a torsion sheaf on $Y$ supported on $\mathfrak{D}_{k}$
  and so we get a torsion sheaf
\[
  \op{gr}^{\Par_{Y}^{\hor}} F_{\lift(t)} := \bigoplus_{k \in K}
  \op{gr}^{\mathfrak{D}_{k}}F_{\lift(t)}
\]
which satisfies:
\begin{itemize}
\item $\op{gr}^{\Par_{Y}^{\hor}} F_{\lift(t)}$ is a vector bundle on
  its support $\Par_{Y}^{\hor} = \sqcup_{k \in K } \mathfrak{D}_{k}$;
\item $\op{gr}^{\Par_{Y}^{\hor}} F_{\lift(t)}$ 
  is naturally filtered by the subsheaves
\[
W_{\ell}\op{gr}^{\Par_{Y}^{\hor}} F_{\lift(t)} := \bigoplus_{k \in K}
  W_{\ell} \op{gr}^{\mathfrak{D}_{k}}F_{\lift(t)};
\]
\item there is a natural  surjective map of $\mathcal{O}_{Y}$-modules:
\begin{equation} \label{eq:surject}
F_{\lift(t)} \to \op{gr}^{\Par_{Y}^{\hor}} F_{\lift(t)}.
\end{equation}
\end{itemize}
Let $W_{\ell}(\hor,F_{\lift(t)})$ be the pullback of
$W_{\ell}\op{gr}^{\Par_{Y}^{\hor}} F_{\lift(t)} \subset
\op{gr}^{\Par_{Y}^{\hor}} F_{\lift(t)}$ via the surjetcive map
\eqref{eq:surject}. Then $W_{\bullet}(\hor,F_{\lift(t)})$ is a
filtration of $F_{\lift(t)}$ by locally free subsheaves which are
equal to $F_{\lift(t)}$ away from $\Par_{Y}^{\hor}$.

For each $i \geq 1$ we have a global and an
$f$-relative residue maps:
\[
\xymatrix@R-2pc{
\Omega^{i}_{Y}(\log \Par_{Y}) \ar[r] & 
\Omega^{i-1}_{\Par_{Y}^{\hor}} \\
\Omega^{i}_{Y/X}(\log \Par_{Y}) \ar[r] & 
\Omega^{i-1}_{\Par_{Y}^{\hor}/X}
}
\]
Tensoring with $W_{0}(\hor,F_{\lift(t)})$ gives maps
\begin{equation} \label{eq:Wres}
\xymatrix@R-2pc{
W_{0}(\hor,F_{\lift(t)})\otimes \Omega^{i}_{Y}(\log \Par_{Y}) \ar[r] & 
W_{0}(\hor,F_{\lift(t)})_{|\Par_{Y}^{\hor}}\otimes \Omega^{i-1}_{\Par_{Y}^{\hor}}, \\
W_{0}(\hor,F_{\lift(t)})\otimes \Omega^{i}_{Y/X}(\log \Par_{Y}) \ar[r] & 
W_{0}(\hor,F_{\lift(t)})_{|\Par_{Y}^{\hor}}\otimes\Omega^{i-1}_{\Par_{Y}^{\hor}/X}.
}
\end{equation}
We define locally free sheaves 
\[
  \begin{aligned}
W_{-2,0}\left(\hor, F_{\lift(t)}\otimes \Omega^{i}_{Y}(\log \Par_{Y})\right)  &
\subset W_{0}(\hor,F_{\lift(t)})\otimes \Omega^{i}_{Y}(\log \Par_{Y}) \\
W_{-2,0}\left(\hor, F_{\lift(t)}\otimes \Omega^{i}_{Y/X}(\log \Par_{Y})\right)  &
\subset W_{0}(\hor,F_{\lift(t)})\otimes \Omega^{i}_{Y/X}(\log \Par_{Y}) 
  \end{aligned}
\]
as the preimages of
\[
  \begin{aligned}
    W_{-2}(\hor, F_{\lift(t)})_{|\Par_{Y}}\otimes
    \Omega^{i-1}_{\Par_{Y}}(\log \Par_{Y}) &
\subset W_{0}(\hor,F_{\lift(t)})_{|\Par_{Y}}\otimes
\Omega^{i-1}_{\Par_{Y}}(\log \Par_{Y}) \\
W_{-2}(\hor, F_{\lift(t)})_{|\Par_{Y}}\otimes
    \Omega^{i-1}_{\Par_{Y}/X}(\log \Par_{Y}) &
\subset W_{0}(\hor,F_{\lift(t)})_{|\Par_{Y}}\otimes
\Omega^{i-1}_{\Par_{Y}/X}(\log \Par_{Y})
  \end{aligned}
\]
under the maps \eqref{eq:Wres}. These subsheaves are preserved by
$\varphi$ and so we get weight modified  global and relative Dolbeault
complexes for $(F_{\lift(t)},\varphi)$:
\begin{equation} \label{eq-absolute-general}
\Dol(Y,F_{\lift(t)})
 := 
\left[ \begin{array}{c}
    W_{0}(\hor,F_{\lift(t)}) \\[+0.3pc]
    \downarrow \wedge \varphi \\[+0.3pc]
    W_{-2,0}\left( \hor, F_{\lift(t)}\otimes \Omega^{1}_{Y}(\log \Par_{Y}) \right) \\[+0.3pc]
    \downarrow \wedge \varphi \\[+0.3pc]
    W_{-2,0}\left( \hor, F_{\lift(t)}\otimes \Omega^{2}_{Y}(\log \Par_{Y}) \right) \\[+0.3pc]
    \downarrow \wedge \varphi \\[+0.3pc]
    \vdots  \\[+0.3pc]
    \downarrow \wedge \varphi \\[+0.3pc]
    W_{-2,0}\left( \hor, F_{\lift(t)}\otimes \Omega^{\dim_{Y}}_{Y}(\log \Par_{Y}) \right) 
\end{array}\right] \  \begin{array}{c}
     0 \\[+0.3pc]
    \\[+0.3pc]
    1 \\[+0.3pc]
    \\[+0.3pc]
    2 \\[+0.3pc]
    \\[+0.3pc]
    \vdots  \\[+0.3pc]
    \\[+0.3pc]
    \dim_{Y}
\end{array}
\end{equation}
and
\begin{equation} \label{eq-relative-general}
\Dol(f,F_{\lift(t)})
 := 
\left[ \begin{array}{c}
    W_{0}(\hor,F_{\lift(t)}) \\[+0.3pc]
    \downarrow \wedge \varphi \\[+0.3pc]
    W_{-2,0}\left( \hor, F_{\lift(t)}\otimes \Omega^{1}_{Y/X}(\log \Par_{Y}) \right) \\[+0.3pc]
    \downarrow \wedge \varphi \\[+0.3pc]
    W_{-2,0}\left( \hor, F_{\lift(t)}\otimes \Omega^{2}_{Y/X}(\log \Par_{Y}) \right) \\[+0.3pc]
    \downarrow \wedge \varphi \\[+0.3pc]
    \vdots  \\[+0.3pc]
    \downarrow \wedge \varphi \\[+0.3pc]
    W_{-2,0}\left( \hor, F_{\lift(t)}\otimes
         \Omega^{\dim_{Y/X}}_{Y/X}(\log \Par_{Y}) \right) 
\end{array}\right] \  \begin{array}{c}
     0 \\[+0.3pc]
    \\[+0.3pc]
    1 \\[+0.3pc]
    \\[+0.3pc]
    2 \\[+0.3pc]
    \\[+0.3pc]
    \vdots  \\[+0.3pc]
    \\[+0.3pc]
    \dim_{Y}-\dim_{X}
\end{array}
\end{equation}
We also define inductively subcomplexes $I^{k}(F_{\lift(t)})$ in the
global Dolbeault complex by setting:
\[
\begin{aligned}
  I^{0}(F_{\lift(t)}) & := \Dol(Y,F_{\lift(t)}),
  \\ I^{k+1}(F_{\lift(t)}) & :=
  \op{im}\left[I^{k}(F_{\lift(t)})\otimes f^{*} \Omega^{1}_{X}(\log \Par_{X})
    \to
    \Dol(Y,F_{\lift(t)}) \right].
\end{aligned}
\]
With this notation we get a short exact sequence of complexes
\begin{equation} \label{eq:ses.DOL}
  \xymatrix{
    0 \ar[d] \\
    \Dol(f,F_{\lift(t)})[-1]\otimes f^{*}\Omega^{1}_{X}(\log \Par_{X})
    \ar[d] \\
    \Dol(Y,F_{\lift(t)})/I^{2}(F_{\lift(t)}) \ar[d] \\
    \Dol(f,F_{\lift(t)}) \ar[d] \\
   0
  }
\end{equation}
which we can view as a morphism
\[
\mathfrak{d}(\varphi) : \Dol(f,F_{\lift(t)}) \to  \Dol(f,F_{\lift(t)})\otimes
f^{*}\Omega^{1}_{X}(\log \Par_{X})
\]
in the derived category $D^{b}_{\op{coh}}(Y,\mathcal{O}_{Y})$.

For every parabolic level $t$ on $X$ consider the sheaf theoretic
pusforward of the pair
$(\Dol(f,F_{\lift(t)}),\mathfrak{d}(\varphi))$. This gives a parabolic
Higgs complex
\begin{equation} \label{eq:alg.push}
f_{*}\mathfrak{d}(\varphi) : f_{*}\Dol(f,F_{\lift(\bullet)}) \longrightarrow
f_{*}\Dol(f,F_{\lift(\bullet)})\otimes \Omega^{1}_{X}(\log \Par_{X})
\end{equation}
on $X$ which conjecturally coincides with Higgs pushforward
$f_{*}(F_{\bullet},\varphi)$ defined above via pushing forward the
twistor $\mathcal{D}$-module corresponding to $(F_{\bullet},\varphi)$.

In \cite{dps} we proved this conjecture in the case when $f : Y \to X$
is a split semistable family of curves and the residue of $\phi$ is
nilpotent.  We noted that the pushforward is uniquely determined by
its restriction to any open in $X$ whose complement has codimension at
least two. We also explained how to apply the formula to a general
family $f : Y \to X$ by using the Abramovich-Karu semistable reduction
theorem but we will not be needing this here since we will only use
the algebraic pushforward formula \eqref{eq:alg.push} for split
semistable families of curves.

\subsection{Stability} \label{ssec:stability}

In this section we will discuss stability of parabolic bundles and
tame parabolic Higgs bundles.

\subsubsection{Semistable parabolic bundles in general}
\label{sssec:ss.general}

Let $(X,\mathcal{O}_{X}(1))$ be a polarized smooth projective variety
of dimension $d$, and let $\Par_{X} \subset X$ be a strict normal
crossings divisor.  Let $\mycal{E}_{\bullet}$ be a torsion free
parabolic sheaf on $(X,\Par_{X})$.  Following \cite{mochizuki-kh1} we
define the rank, parabolic degree, and slope of $\mycal{E}_{\bullet}$
as follows:
\begin{itemize}
\item The  \textit{\bfseries rank of $\mycal{E}_{\bullet}$} is by
  definition the rank of the coherent sheaf $\mycal{E}_{t}$. 
  Tautologically, this is independent of the parabolic level $t$.  We will
  denote this number by $\op{rank} \mycal{E}_{\bullet}$. Note that
  from the Iyer-Simpson formula we also have that $\op{rank}
  \mycal{E}_{\bullet} = \parch_{0}(\mycal{E}_{\bullet})$.
\item The $\mathcal{O}_{X}(1)$-\textit{\bfseries parabolic degree} of
  $\mycal{E}_{\bullet}$ is
  \[
  \pardeg(\mycal{E}_{\bullet}) =
  \int_{X}\parch_{1}(\mycal{E}_{\bullet})\cdot
  c_{1}(\mathcal{O}_{X}(1))^{d-1}.
  \]
\item  The $\mathcal{O}_{X}(1)$-\textit{\bfseries parabolic slope} of
  $\mycal{E}_{\bullet}$ is
  \[
\mathsf{par}\bmu(\mycal{E}_{\bullet}) =
\frac{\pardeg(\mycal{E}_{\bullet})}{\op{rank}(\mycal{E}_{\bullet})}
\]
\end{itemize}

\
With this notation we now have the
following standard

\begin{defi} \label{defi:stable.vb}
A parabolic vector bundle $\boldsymbol{E}_{\bullet}$ on $(X,\Par_{X})$
is {\bfseries slope stable (slope semistable)} if for every saturated
parabolic subsheaf
$\mycal{E}'_{\bullet} \subset \boldsymbol{E}_{\bullet}$ we have
\begin{equation} \label{eq:stable.inequality}
\mathsf{par}\bmu(\mycal{E}'_{\bullet}) < 
\mathsf{par}\bmu(\boldsymbol{E}_{\bullet}) \qquad \left(
\mathsf{par}\bmu(\mycal{E}'_{\bullet}) \leq
\mathsf{par}\bmu(\boldsymbol{E}_{\bullet})\right).
\end{equation}

Similarly for a fixed $z \in \mathbb{C}$, a tame parabolic $z$-connection
$(\boldsymbol{E}_{\bullet},\nabla^{z})$ on $(X,\Par_{X})$ is
{\bfseries slope stable (slope semistable)} if for every
$\nabla^{z}$-invariant saturated parabolic subsheaf
$\mycal{E}'_{\bullet} \subset \boldsymbol{E}_{\bullet}$ the inequality
\eqref{eq:stable.inequality} holds. In particular a tame parabolic
Higgs bundle $(\boldsymbol{E}_{\bullet},\theta)$ is {\bfseries slope
  stable (slope semistable)} if for every saturated parabolic subsheaf
$\mycal{E}'_{\bullet} \subset \boldsymbol{E}_{\bullet}$ for which
$\theta(\mycal{E}'_{\bullet}) \subset \mycal{E}'_{\bullet}\otimes
\Omega^{1}_{X}(\log \Par_{X})$, the inequality
\eqref{eq:stable.inequality} holds.
\end{defi}

\

\noindent
Here as usual we call a parabolic subsheaf $\mycal{E}'_{\bullet}
\subset \boldsymbol{E}_{\bullet}$ \textit{\bfseries saturated} if for
every parabolic level $t$ the quotient
$\boldsymbol{E}_{t}/\mycal{E}'_{t}$ is torsion free.

\

\noindent
For us the stability of tame parabolic Higgs bundles is important
because of Mochizuki's non-abelian Hodge theorem:
\index{terms}{Non-abelian Hodge theory}

\

\smallskip

\noindent
{\bfseries Theorem.}\label{Mo}
 \ \cite{mochizuki-kh1} \textit{Let
  $(X,\mathcal{O}_{X}(1))$ be a polarized projective variety, and let
  $\Par_{X} \subset X$ be a Weil divisor. Assume that $X$ is smooth
  and $\Par_{X}$ has strict normal crossings away from a closed
  subvariety $\mathfrak{Z} \subset X$ of codimension $\geq 3$.
  Suppose $(\boldsymbol{E}_{\bullet},\theta)$ is a tame parabolic
  Higgs bundle on $X - \mathfrak{Z}$ with vanishing
  $\parch_{1}(\boldsymbol{E}_{\bullet})$ and
  $\parch_{2}(\boldsymbol{E}_{\bullet})$. If
  $(\boldsymbol{E}_{\bullet},\theta)$ is slope stable, then there
  exists a canonical analytic family
  $\left(\boldsymbol{V}^{z}_{\bullet},\nabla^{z}\right)$ of stable
  tame parabolic flat $z$-connections on $X - \mathfrak{Z}$, which
  depends analytically on $z$, and satisfies
  $\left(\boldsymbol{V}^{0}_{\bullet},\nabla^{0}\right) =
  (\boldsymbol{E}_{\bullet},\theta)$.
  }

\

\subsubsection{Semistable parabolic bundles on $\mathbb{P}^{1}$}
\label{sssec:ss.P1} 

In this subsection we specialize the foregoing to the case that  
$V$ has rank 2 on a base $X$ that is isomorphic to $\mathbb{P}^{1}$. 
We prove an explicit criterion for stability of such bundles.

Let $C \cong \mathbb{P}^{1}$ be a smooth proper rational
curve, and let $\{ p_{1}, \ldots p_{n} \}$ be a fixed collection of
$n$ distinct points. We would like to understand the stability of
rank two parabolic bundles on $C$ with a parabolic structure along the
divisor $\Par_{C} = \sum_{i = 1}^{n} p_{i}$.

\index{terms}{parabolic!structure}

Suppose $V$ is a rank two
holomorphic bundle on $C$. As discussed in
section~\ref{ssec:par.vb}, a non-degenerate parabolic
structure on $V$ along $\Par_{C}$
consists of the data:
\begin{itemize}
\item a collection $\bF = \{ F_{i} \}_{i = 1}^{n}$ of lines $F_{i}
\subset V_{p_{i}}$, and
\item a collection $(\ba,\bb)$, where $\ba = \{ a_{i} \}_{i=1}^{n}$, $\bb =
  \{ b_{i} \}_{i = 1}^{n}$ of real parabolic weights $\left\{
\left(a_{i},b_{i}\right)\right\}_{i = 1}^{n}$, where $-1 < a_{i} <
b_{i} \leq 0$. 
\end{itemize}
We can also consider degenerate parabolic structures by allowing $a_{i}
= b_{i}$ at some $p_{i}$. In this case the flag $F_{i} \subset
V_{p_{i}}$ will not be a line but will be equal to $0$. 

\index{terms}{parabolic!structure!degenerate}

The parabolic first Chern class (= parabolic determinant) of the
\index{terms}{parabolic!first Chern class}
parabolic bundle \linebreak $(V,\bF,(\ba,\bb))$ is (see
section~\ref{ssec:par.ch}) the $\mathbb{R}$-Cartier divisor
\[
\parch_{1}(V,\bF,(\ba,\bb)) := \det V\otimes
\mathcal{O}_{C}\left(-\sum_{i=1}^{5}
\left(a_{i}+b_{i}\right)p_{i}\right) \in
\Pic(C)\otimes_{\mathbb{Z}} \mathbb{R}, 
\] 
and the parabolic degree of $(V,\bF,(\ba,\bb))$ is the real number 
$\pardeg(V,\bF,(\ba,\bb)) := \deg V - \sum_{i = 1 }^{5} (a_{i} + b_{i})$.
Similarly, if we have a line bundle $L$, then a parabolic
structure on $L$ along $\Par_{C}$  is specified by a collection $\be = \{
e_{i} \}_{i =1}^{n}$, $-1 < e_{i} \leq 0$, $e_{i} \in \mathbb{Q}$ of
parabolic weights. The parabolic first Chern class of $(L,\be)$ is the
$\mathbb{R}$-Cartier divisor
\[
\parch_{1}(L,\be) := L\otimes \mathcal{O}_{C}\left( - \sum_{i = 1}^{n}
e_{i} p_{i}\right) \in \Pic(C)\otimes_{\mathbb{Z}} \mathbb{R},
\]
and the parabolic degree of $(L,\be)$ is the real number
$\pardeg(L,\be) := \deg L - \sum_{i = 1}^{n}
e_{i}$. \index{terms}{parabolic!first Chern class}
Suppose next that 
 $\mathcal{V} = (V,\bF,(\ba,\bb))$ is a rank two parabolic bundle and
$L \subset V$ 
is a line subbundle. Define the {\bfseries critical weights}
$\be^{\op{crit}}(L,\mathcal{V}) = \left\{
e_{i}^{\op{crit}}(L,\mathcal{V}) \right\}$  of
 $L$ relative to $\mathcal{V}$ by 
setting
\[
e_{i}^{\op{crit}}(L,\mathcal{V}) := \inf  \left\{ e_{i} \in \be \, \left| \,
(L,\be) \text{ is  a parabolic subsheaf of } \mathcal{V} \right. \right\}.
\]
In particular $(L,\be)$ is a parabolic subsheaf of $\mathcal{V}$ if and
only if $e_{i} \geq e_{i}^{\op{crit}}(L,\mathcal{V})$ for all $i = 1,
\ldots n$. 

From the definition (see Definition~\ref{defi:stable.vb} and the formula
in Example~\ref{ex:parch1}) of parabolic stability it is straighforward to
check that the critical weights are given by
\[
e_{i}^{\op{crit}}(L,\mathcal{V})  = \begin{cases} a_{i},  & \text{ if }
  \op{im}\left[L_{p_{i}} \to V_{p_{i}}\right] \subseteq F_{i}, \\
 b_{i},  & \text{ if }
 \op{im}\left[L_{p_{i}} \to V_{p_{i}}\right] \nsubseteq F_{i}.
\end{cases}
\]
\

\noindent
The collection of all parabolic bundles with underlying bundle $V$ and
parabolic weights $(\ba,\bb)$ is naturally parameterized by points
$\bF = (F_{1},\ldots,F_{n}) \in \mathbb{P}(V_{p_{1}})\times \cdots
\times \mathbb{P}(V_{p_{n}})$.  The projective line
$\mathbb{P}(V_{p_{i}})$ can be viewed as the fiber of the ruled surface
$\mathbb{P}(V) \to C$ over the point $p_{i} \in C$. 
Thus we can think of a parabolic bundle with a fixed
underlying bundle $V$ and fixed parabolic weights $(\ba,\bb)$ as a
collection of $n$-points $F_{1}, \ldots, F_{n} \in \mathbb{P}(V)$,
mapping to the points $p_{1}, \ldots, p_{n} \in C$.

The surface $\mathbb{P}(V)$ is a useful device for studying the moduli
of parabolic bundles.  We will write $\gamma: \mathbb{P}(V) \to C$ for
the natural projection and $\mathcal{O}_{V}(1)$ for the tautological
hyperplane line bundle, uniquely characterized by the property
$\gamma_{*}\mathcal{O}_{V}(1) = V^{\vee}$. Now given a parabolic
bundle $\mathcal{V} = (V,\bF,(\ba,\bb))$ and a $d \in \mathbb{Z}$
define $\kk_{d}(\mathcal{V})$ by
\[
\kk_{d}(\mathcal{V}) := \max \left\{ \# (\bF\cap D) \,
  \left| \, 
\begin{minipage}[c]{2in}
$D$ is an irreducible curve in the linear system $\left|
  \mathcal{O}_{V}(1)\otimes 
  \gamma^{*} \mathcal{O}_{C}(\deg V + d) 
\right|$.  
\end{minipage}\right.\right\}.
\]
With this notation we have the following 

\begin{lem} \label{lem:stability} Fix a rank two vector bundle 
$V  \cong \mathcal{O}_{C}(k_{1})\oplus \mathcal{O}_{C}(k_{2})$ on $C \cong \mathbb{P}^{1}$,
a divisor $\Par_{C}$ consisting of $n$ distinct points of $C$,
and a balanced system of weights $(\ba,\bb)$ 
with respect to the divisor $\Par_{C}$,
i.e. one in which all $a_{i}$'s
are equal to some $a \in (-1,0]\cap \mathbb{Q}$ and all $b_{i}$'s are
equal to some $b \in [a,0]\cap \mathbb{Q}$.

Let $\mathcal{V} = (V,\bF,(\ba,\bb))$ be a rank two parabolic bundle
on $C$ with parabolic structure along $\Par_{C}$.
\index{terms}{parabolic!structure}
Then $\mathcal{V}$
  is semistable if and only if for all integers $d \geq \min
  \{-k_{1}, -k_{2}\}$ we have
\begin{equation} \label{eq:stable}
(b-a) \kk_{d}(\mathcal{V})  \; \leq \;  
\frac{1}{2} \deg(V)  + (b - a) \frac{n}{2} + d.
\end{equation}
Furthermore $\mathcal{V}$ is stable if and only if \eqref{eq:stable}
holds with a strict inequality.
\end{lem}
{\bfseries Proof.} The parabolic bundle $\mathcal{V}$ is semistable if and only
if for every saturated line subbundle $L \subset V$ and every choice
of weights $\be$ such that $(L,\be)$ is a parabolic subbundle of
$\mathcal{V}$ we have $\pardeg(L,\be) <
\pardeg(\mathcal{V})/2$. 

Indeed suppose that $L \subset V$ is an arbitrary locally free
subsheaf of rank one, and $\be$ is a collection of parabolic weights
for which $(L,\be)$ is a destabilizing parabolic subbundle of
$\mathcal{V}$. First note that for every parabolic subsheaf $(L,\be)
\subset \mathcal{V}$, we have $\pardeg(L,\be) \leq
\pardeg(L,\be^{\op{crit}}(L,\mathcal{V}))$ and so
$(L,\be^{\op{crit}}(L,\mathcal{V})) \subset \mathcal{V}$ is also a
parabolic subsheaf which destabilizes $\mathcal{V}$.

Let $L \subset \widetilde{L} \subset V$ be
the saturation of $L \subset V$. If $p_{i}$ is not in the support of
the torsion subsheaf of $V/L$, then $\op{im}[L_{p_{i}} \to V_{p_{i}}]
= \op{im}[\widetilde{L}_{p_{i}} \to V_{p_{i}}]$ and so
$e_{i}^{\op{crit}}(L,\mathcal{V}) =
e_{i}^{\op{crit}}(\widetilde{L},\mathcal{V})$. If $p_{i}$ is in the
support of the torsion subsheaf of $V/L$, then $\op{im}[L_{p_{i}} \to
V_{p_{i}}] = 0$, and so $e_{i}^{\op{crit}}(L,\mathcal{V}) = a$,
whereas $e_{i}^{\op{crit}}(\widetilde{L},\mathcal{V})$ can be either
$a$ or $b$. Hence $e^{\op{crit}}_{i}(\widetilde{L},\mathcal{V})-a \geq
0$ for all such $p_{i}$ and so
\[
\begin{aligned}
\pardeg\left(\widetilde{L},\be^{\op{crit}}(\widetilde{L},\mathcal{V})\right)
& =  \pardeg(L,\be^{\op{crit}}(L,\mathcal{V}) -  \# \left(
\op{supp}(\op{Tor}(V/L)) - \Par_{C} \right) \\
& \ \hspace{9.35pc} - \sum_{p_{i} \in
  \op{supp}(\op{Tor}(V/L))\cap \Par_{C}}
\left(e^{\op{crit}}_{i}(\widetilde{L},\mathcal{V})-a\right)  \\
& \geq \pardeg(L,\be^{\op{crit}}(L,\mathcal{V}) \\
& > \pardeg(\mathcal{V})/2.
\end{aligned}
\]
In other words 
$\left(\widetilde{L},\be^{\op{crit}}(\widetilde{L},\mathcal{V})\right)$
is a saturated destabilizing parabolic line subbundle of $\mathcal{V}$. 

Suppose next that $L \subset V$ is a saturated line subbundle of
degree $-d \in \mathbb{Z}$. The inclusion $L \subset V$ corresponds to
a non-zero section in $H^{0}(C,\mathcal{O}_{C}(d+k_{1})\oplus
\mathcal{O}_{C}(d+k_{2}))$. In particular $d \geq \min \{ -k_{1},
-k_{2} \}$.  Since $L$ is saturated we have that $V/L \cong
\mathcal{O}_{C}(\deg V + d)$ and the quotient map $V \twoheadrightarrow
\mathcal{O}_{C}(\deg V + d)$ corresponds to a nonzero
section 
\[
s \in H^{0}(C,V^{\vee}\otimes \mathcal{O}_{C}(\deg V + d)) =
H^{0}(\mathbb{P}(V),\mathcal{O}_{V}(1)\otimes
\gamma^{*}\mathcal{O}_{C}(\deg V + d)).
\] 
The  zero divisor $D \subset
\mathbb{P}(V)$ of the section $s \in
H^{0}(\mathbb{P}(V),\mathcal{O}_{V}(1)\otimes 
\gamma^{*}\mathcal{O}_{C}(\deg V + d))$ is irreducible by the
saturation condition and is the section of $\gamma : \mathbb{P}(V) \to
C$ that corresponds to the line bundle $L
\subset V$. Therefore we have
\[
\pardeg\left(L,\be^{\op{crit}}(L,\mathcal{V})\right) = \deg L - a\cdot
\#(\bF\cap D) - b\cdot (n - \#(\bF\cap D)). 
\]
Thus $\mathcal{V}$ is semistable if and only if for all  $D$ we have 
\[
-d - bn + (b-a)\cdot \#(\bF\cap D) \, \leq \, \frac{1}{2}\deg V -
(a+b) \frac{n}{2}
\]
or equivalently
\begin{equation} \label{eq:stable.with.D}
(b-a) \cdot \#(\bF\cap D) \, \leq \, \frac{1}{2}\deg V + (b-a)\frac{n}{2}
+ d.
\end{equation}
Replacing $\#(\bF\cap D)$ in \eqref{eq:stable.with.D} with its maximal
possible value $\kk_{d}(\mathcal{V})$ we get the inequality 
\eqref{eq:stable}. The strict inequality statement is obvious since
the bound  $\kk_{d}(\mathcal{V})$ is always achieved for some $D$. The
lemma is proven. \ \hfill $\Box$

\section{Moduli Spaces} \label{sec:moduli.spaces}

The goal of this chapter is to give an explicit description of a basic geometric object of interest in this work, 
the moduli space of rank 2 semistable vector bundles
on $\mathbb{P}^{1}$ with parabolic structure at 5 marked points.

Let as before $C \cong \mathbb{P}^{1}$, and let $p_{1}$, $p_{2}$,
$p_{3}$, $p_{4}$, $p_{5}$ be five fixed distinct points in $C$. 
We want to describe in detail the geometry of the moduli
space of semistable rank two parabolic vector bundles on $C$ with a
specified parabolic structure along $\{p_{1}, p_{2}, p_{3}, p_{4},
p_{5} \}$.  To streamline the discussion and avoid the analysis of
several cumbersome cases we will make the following two assumptions:

\index{terms}{parabolic!structure}

\begin{itemize}
\item The weights $(\ba,\bb)$ determining the parabolic structure are 
balanced, i.e. all $a_{i}$'s
are equal to some $a \in (-1,0]\cap \mathbb{Q}$ and all $b_{i}$'s are
equal to some $b \in (a,0]\cap \mathbb{Q}$.
\item The parabolic structure is non-degenerate, i.e. $a < b$. 
\end{itemize}

\index{terms}{parabolic!structure!balanced}
\index{terms}{parabolic!structure!non-degenerate}
\index{notations}{babb@$(\ba,\bb)$}
\

\noindent
These assumptions simplify the behavior of the spaces of stable
parabolic bundles and restrict the possible shapes that these spaces
may take. Nevertheless they do give us a rich enough picture which
exhibits all subtleties of the general story without the extra burden
of complicated notation. The complete story can be worked out in a
similar way. In fact the moduli space of parabolic
$SL_{2}(\mathbb{C})$-bundles on $(\mathbb{P}^{1},p_{1},\ldots,p_{n})$
was studied by Bauer \cite{bauer} and by Loray-Saito
\cite{loray-saito} for arbitrary choices of the weights. Here we
essentially rederive and refine Bauer's results in the case $n = 5$ by
a different method which allows us to analyze the geometry of the
moduli space in more detail, allows us to identify the wobbly locus,
and provides the structural properties that we will need later on for
the description of the parabolic Hecke correspondences.

\

\noindent
With these assumptions in mind let us fix numbers $d \in \mathbb{Z}$,
and $a, b \in \mathbb{Q}$ with $-1 < a < b \leq 0$. Consider the
moduli problem $N^{\natural}_{d}(a,b)$ of non-degenerate rank two
parabolic bundles of determinant
$\mathcal{O}(d)$ on $C$, with parabolic structure along
$\Par_{C} = p_{1} +  p_{2} +  p_{3} +  p_{4} +  p_{5}$ with balanced
weights $a_{i} 
= a$ and $b_{i} = b$. In other words $N^{\natural}_{d}(a,b)$ is the functor
\[
N^{\natural}_{d}(a,b) : \left( \op{{\sf Sch}}/\mathbb{C}
\right)^{\op{op}} \to \left(\op{{\sf Set}}\right)
\]
from complex schemes to sets, which assigns to every complex scheme
$T$ the set $N^{\natural}_{d}(a,b)(T)$ of isomorphism classes of
$T$-families of semi-stable non-degenerate rank two parabolic bundles
of determinant $\mathcal{O}(d)$ on $(C,\Par_{C})$ with parabolic weights $a_{i} = a$ and $b_{i} = b$.

\index{notations}{Nnatdab@$N^{\natural}_{d}(a,b)$}
 
By the classical Mehta-Seshadri theorem \cite{mehta-seshadri,seshadri-book} 
we know that $N^{\natural}_{d}(a,b)$ is coarsely represented by a
projective moduli space $N_{d}(a,b)$. Our goal is to describe
$N_{d}(a,b)$ explicitly for the various choices of $(d,a,b)$. For this
it will be useful to understand the dependence of  $N_{d}(a,b)$ on the
parameters $(d,a,b)$.

\index{notations}{Ndab@$N_{d}(a,b)$}

\subsection{A family of moduli problems} \label{sec:family}

In this section we will augment the family of moduli problems
$N^{\natural}_{d}(a,b)$ slightly, thereby clarifying the symmetries of the situation. 

Recall
\cite{carlos-nc,borne-root,iyer-simpson}  
that a rank two parabolic bundle
$\mathcal{V}$ has two
different equivalent incarnations: 
\begin{description}
\item[$(\dagger)$] as a bundle $\mathcal{V} =
(V,\bF,(\ba,\bb))$ equipped with flags $\bF = \{ F_{i} \subset
  V_{p_{i}}\}$  and weights $\ba = \{a_{i}\}$ and $\bb = \{b_{i}\}$ at
  the points 
  $p_{i} \in \Par_{C}$, and 
\item[$(\ddagger)$] as a collection $\mathcal{V} = V_{\bullet} = \{
V_{w} \}_{w : \Par_{C} \to \mathbb{R}}$   where
\begin{itemize}
\item $V_{w}$ are locally free sheaves on $C$;
\item The collection is monotone, i.e. $V_{w} \subset V_{w'} $
  whenever $w \leq w'$ and these inclusions are all compatible with
  each other.
\item $V_{\bullet}$ is semicontinuous, i.e. for every $w$ we can find a
real  number $c > 0$ so that $V_{w + \bepsilon} = V_{w}$ for all
functions $\bepsilon : \Par_{C} \to [0,c]$;
\item  if $\bdelta_{p_{i}} : \Par_{C} \to
\mathbb{R}$ is the characteristic function of $p_{i} \in
\Par_{C}$, then $V_{w + \bdelta_{p_{i}}} = V_{w}(p_{i})$ and this
identification is compatible with the inclusions of $V_{w}$'s;
\end{itemize}
\end{description}

\

\noindent
If a parabolic rank two bundle $\mathcal{V}$ is given by the data 
$V_{\bullet}$, then we can reconstruct the data $(V,\bF,(\ba,\bb))$ by
setting $V := V_{0}$, taking $a_{i}$, $b_{i}$ to be
the points of semicontinuity jumps of $V_{t\cdot \bdelta_{p_{i}}}$ 
on the interval
$(-1,0]$, and taking $F_{i} \subset V_{p_{i}}$ to be the lines
  corresponding to these jumps. In other words  the numbers 
$-1 < a_{i} \leq b_{i}
  \leq 0$ are uniquely determined by the identification
\[
V_{t\cdot \bdelta_{p_{i}}} = \begin{cases}
V, & \text{ for } \quad  b_{i} < t \leq 0, \\
V_{b_{i}}, & \text{ for } \quad  
a_{i} < t \leq
b_{i},  \\
V(-p_{i}), & \text{ for } -1 \leq t \leq a_{i},
\end{cases}
\]
and $F_{i} =
\op{im}\left[\left(V_{b_{i}}\right)_{p_{i}} \to V_{p_{i}}\right]$. 
Conversely given data $(V,\bF,(\ba,\bb))$ of type $(\dagger)$ we can
define data $V_{\bullet}$ of type $(\ddagger)$ by setting
\[
V_{t\cdot \bdelta_{p_{i}}} := \begin{cases}
V, & \text{ for } \quad  b_{i} < t \leq 0, \\
\ker\left[V \to V_{p_{i}}/F_{i}\right], & \text{ for } \quad  
a_{i} < t \leq
b_{i},  \\
V(-p_{i}), & \text{ for } -1 \leq t \leq a_{i}.
\end{cases}
\]
We extend this to other values of $t$ by setting 
$V_{(t+n)\cdot \bdelta_{p_{i}}} :=V_{t\cdot \bdelta_{p_{i}}}(np_i)$ for integer $n$,
and define:
\[
V_{w} :=  \bigcap_{i=1}^{5} V_{w(p_{i})\cdot \bdelta_{p_{i}}}.
\]
From this point of view it is natural to extend the definition of the
data $(\dagger)$ to consist of triples $(V,\bF,(\ba,\bb))$ for which
the parabolic weights $a_{i}$ and $b_{i}$ are arbitrary real numbers
satisfying $a_{i} \leq b_{i} \leq a_{i}+1$. Again any such data 
$(V,\bF,(\ba,\bb))$ gives rise to data of type $(\ddagger)$ via the
assignment
\[
V_{t\cdot \bdelta_{p_{i}}} := \begin{cases}
V, & \text{ for } \quad  b_{i} < t \leq a_{i}+1, \\
\ker\left[V \to V_{p_{i}}/F_{i}\right], & \text{ for } \quad  
a_{i} < t \leq
b_{i},
\end{cases}
\]
We again extend this to other values of $t$ by setting 
$V_{(t+n)\cdot \bdelta_{p_{i}}} :=V_{t\cdot \bdelta_{p_{i}}}(np_i)$ for integer $n$,
and we define: 
\[
V_{w} :=  \bigcap_{i=1}^{5} V_{w(p_{i})\cdot \bdelta_{p_{i}}}.
\]
The parabolic Chern classes \index{terms}{parabolic!Chern classes}
and the stability conditions transfer
verbatim to this case and so for every triple $(d,\ba,\bb)$ with
$d \in \mathbb{Z}$, and $\ba, \bb \in \mathbb{R}^{5}$, satisfying
$a_{i} < b_{i} < a_{i}+1$ we can get a well defined moduli functor
$N^{\natural}_{d}(\ba,\bb) : \left( \op{{\sf Sch}}/\mathbb{C}
\right)^{\op{op}} \to \left(\op{{\sf Set}}\right)$, which assigns to
every complex scheme $T$ the set $N^{\natural}_{d}(\ba,\bb)(T)$ of
isomorphism classes of $T$-families of semi-stable non-degenerate rank
two parabolic bundles of determinant $\mathcal{O}(d)$ on
$(C,\Par_{C})$ and parabolic weights $(\ba,\bb)$. The Mehta-Seshadri
theorem \cite{mehta-seshadri,seshadri-book} again applies to this
situation and implies that for every $(d,\ba,\bb)$ the functor
$N^{\natural}_{d}(\ba,\bb)$ is coarsely representable by a projective
moduli space $N_{d}(\ba,\bb)$.

The family of functors $N^{\natural}_{d}(\ba,\bb)$ parametrized by
$(d,\ba,\bb) \in \mathbb{Z}\times \mathbb{R}^{5}\times \mathbb{R}^{5}$ 
with $\ba \leq \bb \leq \ba +1$
has some very natural automorphisms:

\index{notations}{Nnatdab@$N^{\natural}_{d}(\ba,\bb)$}
\index{notations}{tk@$\bt_{k}$}
\index{notations}{ws@$\bws_{\varepsilon}$}
\index{notations}{ps@$\bps$}

\

\begin{description}
\item[{\bfseries (tensorization)}] 
Tensoring
with $\mathcal{O}_{C}(k)$ transforms a parabolic bundle
$(V,\bF,(\ba,\bb))$ into a new parabolic bundle $\left(V\otimes
\mathcal{O}_{C}(k), \left\{ F_{i}\otimes
\mathcal{O}_{C}(k)_{p_{i}}\right\}_{i=1}^{5},(\ba,\bb)\right)$. This
operation transforms the associated data $(\ddagger)$ in   an obvious
manner by sending  
%$V_{\bullet}$  
$V_{w}$ to $V_{w}\otimes \mathcal{O}_{C}(k)$ for all $w$. In particular for any $k$
tensoring with $\mathcal{O}_{C}(k)$ gives rise to an isomorphism
of moduli functors:
\[
\bt_{k} : N^{\natural}_{d}(\ba,\bb) \to N^{\natural}_{d+2k}(\ba,\bb),
\]
and these isomorphisms satisfy $\bt_{k}\circ \bt_{k'} = \bt_{k+k'}$. 
\item[{\bfseries (weight shift)}] The additive group
  $\op{Fun}(\Par_{C},\mathbb{R})$ acts on parabolic bundles
  $(V,\bF,(\ba,\bb))$  by shifting the weights:  a function $\bepsilon
  : \Par_{C} \to 
  \mathbb{R}$  sends $(V,\bF,(\ba,\bb))$  to
  $(V,\bF,(\ba+\bepsilon,\bb+\bepsilon))$. Equivalently $\bepsilon$
  acts on the associated data  $(\ddagger)$ by shifting the labelling:
  it sends $\mathcal{V} = V_{\bullet}$ to a new parabolic bundle
  $\mathcal{V} + \bepsilon$ which is defined by $\mathcal{V} + \bepsilon
  = (V+\bepsilon)_{\bullet}$ with $(V+\bepsilon)_{w} :=
  V_{w+\bepsilon}$. This gives rise to an isomorphism
of moduli functors:
\[
\bws_{\bepsilon} : N^{\natural}_{d}(\ba,\bb) \to
N^{\natural}_{d}(\ba+\bepsilon,\bb+\bepsilon), 
\]
and these isomorphisms satisfy $\bws_{\bepsilon}\circ \bws_{\bepsilon'} =
\bws_{\bepsilon+\bepsilon'}$. 
\item[{\bfseries (parity shift)}] For a given $p_{i} \in \Par_{C}$ we
  can shift the $p_{i}$-discontinuity of a parabolic bundle to the
  right. Explicitly given $\mathcal{V} = (V,\bF,(\ba,\bb))$  we define the
  parity shifted bundle ${}^{p_{i}}\mathcal{V} =
  \left({}^{p_{i}}V,{}^{p_{i}}\bF,{}^{p_{i}}(\ba,\bb)\right)$ by
  setting 
\[
\begin{aligned}
{}^{p_{i}}V &  :=  \ker\left[V \to V_{p_{i}}/F_{i}\right]\otimes
\mathcal{O}_{C}(p_{i}), \\
{}^{p_{i}}F_{i} & := \ker\left[{}^{p_{i}}V_{p_{i}} \to F_{i}\otimes
  \mathcal{O}_{C}(p_{i})_{p_{i}}\right], \\
{}^{p_{i}}(\ba,\bb) & := \left(\bb,\ba + \bdelta_{p_{i}}\right).
\end{aligned}
\]
In particular this operation does not change the associated data of type
$(\ddagger)$. The parity shift across $p_{i}$ induces an isomorphism of
moduli functors
\[
\bps_{p_{i}} : N^{\natural}_{d}(\ba,\bb) \to
N^{\natural}_{d+1}\left(\bb,\ba + \bdelta_{p_{i}}\right). 
\]
\end{description}

\

\noindent
If we specialize to  the balanced non-degenerate case,
i.e. the case when the weights $(\ba,\bb)$ satisfy $a_{i} = a$, $b_{i}
= b$, with $a < b < a+1$, then we get a slightly more restricted
family of isomorphisms of moduli functors:
\begin{description}
\item[\quad {\bfseries (tensorization)}] 
$\bt_{k} : N^{\natural}_{d}(a,b) \to N^{\natural}_{d+2k}(a,b)$,
  for $k \in \mathbb{Z}$.
\item[\quad {\bfseries (weight shift)}] $\bws_{\varepsilon} :
  N^{\natural}_{d}(a,b) \to 
  N^{\natural}_{d}(a+\varepsilon,b+\varepsilon)$, for $\varepsilon \in
  \mathbb{R}$; here $\bws_{\varepsilon} := \bws_{\varepsilon\cdot
    \sum_{i = 1}^{5}\bdelta_{p_{i}}}$. 
\item[\quad {\bfseries (parity shift)}] $\bps  : N^{\natural}_{d}(a,b) \to
  N^{\natural}_{d+5}(b,a+1)$; here $\bps := \bps_{p_{1}}\circ
  \bps_{p_{2}}\circ \bps_{p_{3}}\circ \bps_{p_{4}}\circ \bps_{p_{5}}$.
\end{description}
By using these isomorphisms we get many accidental isomorphisms 
among the moduli problems $N^{\natural}_{d}(a,b)$. For future
reference we record this fact in the following

\begin{lem} \label{lem:moduli.isos}
If $a < b < a+1$ the  moduli problem $N^{\natural}_{d}(a,b)$ is
isomorphic to a moduli 
problem $N^{\natural}_{0}(0,q)$ with $0 < q < 1$.
\end{lem}
{\bfseries Proof.} If $d$ is even, then we have an isomorphism 
\[
\bt_{-\frac{d}{2}}\circ
\bws_{-a} : 
N^{\natural}_{d}(a,b) \stackrel{\cong}{\longrightarrow}
N^{\natural}_{0}(0,b-a). 
\] 
If $d$
is odd, then we have an isomorphism 
\[
\bt_{-\frac{d+5}{2}}\circ 
\bws_{-b}\circ \bps : N^{\natural}_{d}(a,b) \stackrel{\cong}{\longrightarrow}
N^{\natural}_{0}(0,a-b+1).
\]
This proves the lemma. \ \hfill $\Box$

\

\noindent
In view of this lemma it suffices to understand the moduli spaces
$N^{\natural}_{0}(0,q)$ for $0 < q < 1$.

\subsection{Birational models and chambers} \label{ssec:chambers}

The stability condition for a bundle $(V,\bF,(0,q))$ depends on the
parity of the degree of $V$ as well as the parameter $q$.  This
dependence is piecewise constant.  In this section we analyze this
stability condition.  For each parity, we find 4 chambers in $q$-space
and describe the stability condition in each.  Since the two parities
are related by the isomorphism of moduli spaces $N_{0}(0,q) \cong
N_{1}(0,1-q)$ explained in the previous section, we are free to choose
the most convenient parity for any given calculation. Most of the
analysis is done for the even parity but for one of the chambers we
switch to odd parity where the description of the geometry of the
moduli space is more straightforward.

For  degree $0$, it turns out that there are only three options for the
underlying
vector bundle $V$: it can be $\mathcal{O}_{C}(-k)\oplus
\mathcal{O}_{C}(k)$ for $k=0, 1, 2$.
We examine these separately in the three subsections below.
The results are then combined  in the next section \ref{sec:git}
(in Table~\ref{table:deg.0.stable}) 
and applied to identify the GIT moduli spaces in the four chambers.

Let  $V$ be  of degree 0, so that
$V \cong \mathcal{O}_{C}(-k)\oplus
\mathcal{O}_{C}(k)$ for some $k \geq 0$. In this case the stability
criterion of 
Lemma~\ref{lem:stability} implies that $\mathcal{V} =
(V,\bF,(0,q))$ is semistable if and only if for all $d \geq -k$ we
have
\begin{equation} \label{eq:deg0.stability}
\kk_{d}(\mathcal{V}) \leq \frac{5}{2} + \frac{d}{q},
\end{equation}
where $\kk_{d}(\mathcal{V})$ is defined immediately above
Lemma~\ref{lem:stability}.

This allows us to classify completely the semistable parabolic bundles
with balanced parabolic weights. Indeed, note that if $d = -k$, the
linear system $|\mathcal{O}_{V}(1)\otimes
\gamma^{*}\mathcal{O}_{C}(-k)|$ defining $\kk_{d}(\mathcal{V})$
contains a unique irreducible curve. So for $d = -k$ we have
$\kk_{-k}(\mathcal{V}) \geq 0$, and \eqref{eq:deg0.stability} gives
\[
0 \leq \kk_{-k}(\mathcal{V}) \leq \frac{5}{2} - \frac{k}{q}
\leq \frac{5}{2} - k,
\]
confirming that $k \leq 2$.  Thus we must
have $k = 0, 1, 2$. We analyze these three cases in detail:

\subsubsection{Semistable parabolic structures on
  $\mathcal{O}_{C}\oplus \mathcal{O}_{C}$} \label{sssec:k=0}
Suppose $k = 0$. Then $\mathbb{P}(V) = \mathbb{P}^{1}\times C$
is isomorphic to the Hirzebruch surface $\mathbb{F}_{0}$, and for $d
\geq 0$ the linear system defining $\kk_{d}(\mathcal{V})$ is the linear
system of all curves of bidegree $(1,d)$. So by
\eqref{eq:deg0.stability} we get the following 
stability constraints on the $5$-tuple of points $\bF$ in
$\mathbb{P}(V) = \mathbb{P}^{1}\times C$:
\begin{description}
\item[\framebox{\framebox{$d=0$}} \ :]  \ At most $2$ of the points $F_{1},
  \ldots, F_{5}$ can 
  lie on a $(1,0)$ ruling. In particular if we embed
  $\mathbb{P}^{1}\times C$ as a quadric in $\mathbb{P}^{3}$, no $3$ 
  of the points $F_{1}, \ldots, F_{5}$ are colinear in
  $\mathbb{P}^{3}$. This stability constraint is valid for all values
  of $0 < q < 1$.
\item[\framebox{\framebox{$d=1$}} \ :]  \ \begin{description}
\item[\fbox{$\frac{2}{3} < q < 1$}] \ No $4$ of $F_{1}, \ldots, F_{5}$ can
  lie on an irreducible $(1,1)$ curve. Combined with the $d = 0$
  constraint this implies that from 
  the point of view of the embedding of the quadric in
  $\mathbb{P}^{3}$, no $4$  
  of the points $F_{1}, \ldots, F_{5}$ are coplanar.
\item[\fbox{$\frac{2}{5} < q \leq \frac{2}{3}$}] \ The $5$ points $F_{1},
  \ldots, F_{5}$ can not lie on an irreducible $(1,1)$ curve. From the
  point of view of $\mathbb{P}^{3}$ this means that the $5$ points can
  not all  be  coplanar.
\item[\fbox{$0 < q \leq \frac{2}{5}$}] \ There is no constraint
  on the five points in this case.
\end{description}
\item[\framebox{\framebox{$d=2$}} \ :]  \  \begin{description}
\item[\fbox{$\frac{4}{5} \leq q \leq 1$}] \ The $5$ points $F_{1},
  \ldots, F_{5}$ can not lie on an irreducible $(1,2)$ curve. Combined
  with the previous constraints this shows that there is no semistable
  $5$-point configuration for such parabolic weights.
\item[\fbox{$0 < q < \frac{4}{5}$}] \  There is no constraint
  on the five points in this case.
\end{description}
\item[\framebox{\framebox{$d\geq 3$}} \ :]  \ There is no constraint
  on the five points in this case.
\end{description}
\ 

\noindent
In summary: the semistable parabolic bundles with underlying bundle
$\mathcal{O}_{C}\oplus \mathcal{O}_{C}$ and balanced parabolic weights
$(0,q)$ correspond to all choices of $5$-tuples of
points $F_{1}, \ldots, F_{5} \in \mathbb{P}^{1}\times C$, such that 
$F_{i}$ projects to $p_{i} \in C$ for all $i$ and the parabolic
weights and the  $5$-tuple are
constrained as follows:
\begin{itemize}
\item $0 < q \leq \frac{2}{5}$ and no three of the five points 
 $F_{1}, \ldots, F_{5} \in \mathbb{P}^{1}\times C$ can lie on a
 $(1,0)$ ruling.
\item $\frac{2}{5} < q \leq \frac{2}{3}$, no three of the
 five points $F_{1}, \ldots, F_{5} \in \mathbb{P}^{1}\times C$ can lie
 on a $(1,0)$ ruling, and the five points can not all lie on an
 irreducible $(1,1)$ curve.
\item $\frac{2}{3} < q <  \frac{4}{5}$, no three of the
 five points $F_{1}, \ldots, F_{5} \in \mathbb{P}^{1}\times C$ can lie
 on a $(1,0)$ ruling, and no  four of the
 five points $F_{1}, \ldots, F_{5} \in \mathbb{P}^{1}\times C$ can lie on an
 irreducible $(1,1)$ curve.
\end{itemize}

\subsubsection{Semistable parabolic structures on
 $\mathcal{O}_{C}(-1)\oplus \mathcal{O}_{C}(1)$} \label{sssec:k=1}

Suppose $k = 1$. Then $\gamma : \mathbb{P}(V) \to C$ is isomorphic to
the Hirzebruch surface $\mathbb{F}_{2}$.  In this case $d \geq -1$ and
the linear system defining $\kk_{d}(\mathcal{V})$ is
$|\mathcal{O}_{V}(1)\otimes \gamma^{*}\mathcal{O}_{C}(d)|$ and we get
the following stability constraints on the $5$-tuple of points $\bF$:

\index{terms}{parabolic!structure!semistable}

\begin{description}
\item[\framebox{\framebox{$d=-1$}} \ :] In this case the allowed
  values for $\kk_{-1}(\mathcal{V})$ are $\kk_{-1}(\mathcal{V}) = 0,
  1$ and we have 
\begin{description}
\item[\framebox{$\frac{2}{3} < q < 1$} \ :] At most one of the five
  points $F_{1}$, \ldots, $F_{5}$ lies on the unique $(-2)$ curve in
  $\mathbb{P}(V) \cong \mathbb{F}_{2}$. 
\item[\framebox{$\frac{2}{5} < q \leq \frac{2}{3}$} \ :] None of the five
  points $F_{1}$, \ldots, $F_{5}$ lies on the $-2$ curve in
  $\mathbb{P}(V) \cong \mathbb{F}_{2}$. 
\item[\framebox{$0 < q \leq \frac{2}{5}$} \ :] No collection $\bF$ of
  five points is stable.
\end{description}
\item[\framebox{\framebox{$d=0$}} \ :] At most 2 of the points
  $F_{1}$, \ldots, $F_{5}$ can lie on an irreducible curve in the
  linear system $\mathcal{O}_{V}(1)$. This stability constraint is
  valid for all values of $0 < q < 1$. Note however that the linear
  system $|\mathcal{O}_{V}(1)|$ has a 
  base component, namely the unique $(-2)$ curve in
  $\mathbb{P}(V) \cong \mathbb{F}_{2}$. Hence  $|\mathcal{O}_{V}(1)|$ does
  not have any irreducible members and so this stability condition is
  vacuous. 
\item[\framebox{\framebox{$d=1$}} \ :] In this case
  $|\mathcal{O}_{V}(1)\otimes \gamma^{*}\mathcal{O}_{C}(1)| \cong
  \mathbb{P}^{3}$ and so we have a curve in this linear system that
  passes through every three points in $\mathbb{P}(V)$. Thus
  $\kk_{1}(\mathcal{V}) \geq 3$ for all $\mathcal{V}$ and the interesting
  values for $\kk_{1}(\mathcal{V})$ are $\kk_{1}(\mathcal{V}) = 4,
  5$. We have 
\begin{description}
\item[\framebox{$\frac{2}{3} < q < 1$} \ :] At most three of the
  five points $F_{1}$, \ldots, $F_{5}$ can lie on an irreducible curve
  in the linear system $|\mathcal{O}_{V}(1)\otimes
  \gamma^{*}\mathcal{O}_{C}(1)|$.
\item[\framebox{$\frac{2}{5} < q \leq \frac{2}{3}$} \ :] At most
  four of the five points $F_{1}$, \ldots, $F_{5}$ can lie on an
  irreducible curve in the linear system $|\mathcal{O}_{V}(1)\otimes
  \gamma^{*}\mathcal{O}_{C}(1)|$.
\item[\framebox{$0 < q \leq \frac{2}{5}$} \ :] There is no
  stability constraint on the five points.
\end{description}
\item[\framebox{\framebox{$d=2$}} \ :] In this case
  $|\mathcal{O}_{V}(1)\otimes \gamma^{*}\mathcal{O}_{C}(2)| \cong
  \mathbb{P}^{5}$ and so we have a curve in this linear system that
  passes through every five points in $\mathbb{P}(V)$. Thus
  $\kk_{2}(\mathcal{V}) = 5$ and we have
\begin{description}
\item[\framebox{$\frac{4}{5} \leq q < 1$} \ :] There are no stable
  configurations of five points.
\item[\framebox{$0 < q < \frac{4}{5}$} \ :] There is no
  stability constraint on the five points in this case.
  \end{description} 
\item[\framebox{\framebox{$d\geq 3$}} \ :] There is no
  stability constraint on the five points in this case.
\end{description}
\

\noindent
In summary: the semistable parabolic bundles with underlying bundle $V
= \mathcal{O}_{C}(-1)\oplus \mathcal{O}_{C}(1)$ and balanced parabolic
weights $(0,q)$ correspond to all choices of $5$-tuples of
points $F_{1}, \ldots, F_{5} \in \mathbb{P}(V)$, such that $F_{i}$
projects to $p_{i} \in C$ for all $i$ and the parabolic weights and
the $5$-tuple are constrained as follows:
\begin{itemize}
\item $\frac{2}{5} < q \leq \frac{2}{3}$, none of the five points
  lie on the unique irreducible $-2$ curve in $\mathbb{P}(V)$, and at
  most $4$ of the five points 
  lie on an irreducible curve in the linear 
  system $|\mathcal{O}_{V}(1)\otimes \gamma^{*}\mathcal{O}_{C}(1)|$.
\item $\frac{2}{3} < q < \frac{4}{5}$, at most $1$ of the five points
  lies on the unique irreducible $-2$ curve in $\mathbb{P}(V)$, and at
  most $3$ of the five points 
  lie on an irreducible curve in the linear 
  system $|\mathcal{O}_{V}(1)\otimes \gamma^{*}\mathcal{O}_{C}(1)|$.
\item $\frac{4}{5} \leq q < 1$, there are no stable configurations of
  five points.
\end{itemize}

\

\subsubsection{Semistable parabolic structures on
 $\mathcal{O}_{C}(-2)\oplus \mathcal{O}_{C}(2)$} \label{sssec:k=2}

Suppose $k = 2$. Then $\gamma : \mathbb{P}(V) \to C$ is isomorphic to
the Hirzebruch surface $\mathbb{F}_{4}$.  In this case $d \geq -2$,
the linear system defining $\kk_{d}(\mathcal{V})$ is
$|\mathcal{O}_{V}(1)\otimes \gamma^{*}\mathcal{O}_{C}(d)|$, and we get
the following stability constraints on the $5$-tuple of points $\bF$:

\index{terms}{parabolic!structure!semistable}

\begin{description}
\item[\framebox{\framebox{$d=-2$}} \ :] In this case the only allowed
  value for $\kk_{-2}(\mathcal{V})$ is $\kk_{-2}(\mathcal{V}) = 0$,
  and we have 
\begin{description}
\item[\framebox{$0 < q < \frac{4}{5}$} \ :] There are no stable configurations
  of five points. 
\item[\framebox{$\frac{4}{5} \leq q \leq 1$} \ :] None of the five
  points $F_{1}$, \ldots, $F_{5}$ lies on the unique $-4$ curve in the surface 
  $\mathbb{P}(V) \cong \mathbb{F}_{4}$. 
\end{description}
\item[\framebox{\framebox{$d=-1,0,1$}} \ :] The linear
  system $|\mathcal{O}_{V}(1)\otimes \gamma^{*}\mathcal{O}_{C}(d)|$ has a 
  base component, namely the unique $(-4)$ curve in
  $\mathbb{P}(V) \cong \mathbb{F}_{4}$. Hence
  $|\mathcal{O}_{V}(1)\otimes \gamma^{*}\mathcal{O}_{C}(d)|$ does
  not have any irreducible members and so the stability constraint is
  vacuous in these cases.
\item[\framebox{\framebox{$d = 2$}} \ :] \ \hfill \quad
  \begin{description}
  \item[\framebox{$0 < q \leq \frac{4}{5}$} \ :] There is no stability
    constraint on the five points in this case.
\item[\framebox{$\frac{4}{5} < q \leq 1$} \ :] At most four of the
  five points $F_{1}$, \ldots, $F_{5}$ lie on the unique $-4$ curve in
  the surface $\mathbb{P}(V) \cong \mathbb{F}_{4}$. This holds
  automatically in view of the stronger condition found in the case $d
  = -2$.
  \end{description}
\item[\framebox{\framebox{$d \geq 3$}} \ :] There is no stability
    constraint on the five points in this case.
\end{description}
\

\

\noindent
In summary: the semistable parabolic bundles with underlying bundle $V
= \mathcal{O}_{C}(-2)\oplus \mathcal{O}_{C}(2)$ and balanced parabolic
weights $(0,q)$ correspond to all choices of $5$-tuples of
points $F_{1}, \ldots, F_{5} \in \mathbb{P}(V)$, such that $F_{i}$
projects to $p_{i} \in C$ for all $i$ and the parabolic weights and
the $5$-tuple are constrained as follows:
\begin{itemize}
\item $0 < q < \frac{4}{5}$, there are no stable configurations of
  five points.
\item $\frac{4}{5} \leq  q \leq 1$, none of the five points
  lie on the unique irreducible $-4$ curve in $\mathbb{P}(V)$.
\end{itemize}

\subsection{The GIT construction} \label{sec:git} 

\begin{figure}[!ht]
\begin{center}
\psfrag{q}[c][c][1][0]{{$q$}}
\psfrag{0}[c][c][1][0]{{$0$}}
\psfrag{2/5}[c][c][1][0]{{${\displaystyle \frac{2}{5}}$}}
\psfrag{2/3}[c][c][1][0]{{${\displaystyle \frac{2}{3}}$}}
\psfrag{4/5}[c][c][1][0]{{${\displaystyle \frac{4}{5}}$}}
\psfrag{1}[c][c][1][0]{{$1$}}
\epsfig{file=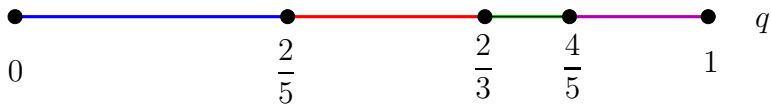,width=4in} 
\end{center}
\caption{The four chambers in the space of parabolic
  weights in even degree}\label{fig:chambers}  
\end{figure}

Recall that by the classical Mehta-Seshadri theorem
\cite{mehta-seshadri,seshadri-book} we know that
$N^{\natural}_{d}(a,b)$ is coarsely represented by a projective moduli
space $N_{d}(a,b)$.  We are now ready to construct and analyze the
structure of these projective moduli spaces $N_d(a,b)$.  This is
accomplished in Theorem \ref{theo:shapes.of.moduli} for $N_{0}(0,q)$
and in Proposition \ref{prop:deg1.moduli} for $N_{1}(0,q)$, for $0 < q
< 1$.  These spaces are of course isomorphic, via a parity shift
interchanging $q$ with $1-q$.  Nevertheless, we find it easier to
analyze some chambers in one parity and the others in the opposite
parity.  The discussion in the previous section shows that for both
moduli problems there are four natural chambers in the space of
parabolic weights according to the value of $q$ (see
Figure~\ref{fig:chambers}), and that the stability changes only after
crossing from one chamber to an adjacent one.

\

\subsubsection{The moduli spaces $N_{0}(0,q)$} \label{sssec-N0}

The moduli spaces $N_{0}(0,q)$ are canonically identified for all $q$
in the interior of any given chamber. The semi-stable parabolic data
$\mathcal{V} = (V,\bF,(\ba,\bb))$ for the various chambers is recorded
in Table~\ref{table:deg.0.stable}.

\begin{table}[H]
\begin{center}
{\small
\begin{tabular}{|c|c|c|} 
\hline
$V$ & $\bF$ & $q$ \\ 
\hline
\hline 
$\mathcal{O}_{C}\oplus \mathcal{O}_{C}$
& 
\begin{minipage}[c]{3.8in} \addtolength{\baselineskip}{-5pt}
\

\begin{itemize}
\item At most $2$ points in $\bF$ lie on a $(1,0)$ ruling in 
$\mathbb{P}(V)$.
\end{itemize}

\
\end{minipage} &
${\displaystyle 0 < q \leq \frac{2}{5}}$ \\
\hline\hline
$\mathcal{O}_{C}\oplus \mathcal{O}_{C}$ & 
\begin{minipage}[c]{3.8in} \addtolength{\baselineskip}{-5pt}
\

\begin{itemize}
\item At most $2$  points in $\bF$ lie on a $(1,0)$ ruling in 
$\mathbb{P}(V)$.
\item At most $4$  points in $\bF$ lie on an irreducible  $(1,1)$ curve in 
$\mathbb{P}(V)$.
\end{itemize}

\
\end{minipage} &
\multirow{2}{*}{${\displaystyle \frac{2}{5} < q \leq \frac{2}{3}}$} \\
\cline{1-2} 
$\mathcal{O}_{C}(-1)\oplus \mathcal{O}_{C}(1)$ &
\begin{minipage}[c]{3.8in} \addtolength{\baselineskip}{-5pt}
\

\begin{itemize}
\item No  points in $\bF$ lie on the unique rational $(-2)$ curve in 
$\mathbb{P}(V)$.
\item At most $4$ points in $\bF$ lie on an irreducible curve $D
  \subset \mathbb{P}(V)$, $D \in |\mathcal{O}_{V}(1)\otimes
  \gamma^{*}\mathcal{O}_{C}(1)|$. 
\end{itemize}

\
\end{minipage}
& \\
\hline\hline
$\mathcal{O}_{C}\oplus \mathcal{O}_{C}$ & 
\begin{minipage}[c]{3.8in} \addtolength{\baselineskip}{-5pt}
\

\begin{itemize}
\item At most $2$  points in $\bF$ lie on a $(1,0)$ ruling in 
$\mathbb{P}(V)$.
\item At most $3$  points in $\bF$ lie on an irreducible  $(1,1)$ curve in 
$\mathbb{P}(V)$.
\end{itemize}

\
\end{minipage} &
\multirow{2}{*}{${\displaystyle \frac{2}{3} < q < \frac{4}{5}}$} \\
\cline{1-2} 
$\mathcal{O}_{C}(-1)\oplus \mathcal{O}_{C}(1)$ &
\begin{minipage}[c]{3.8in} \addtolength{\baselineskip}{-4pt}
\

\begin{itemize}
\item At most $1$ point in $\bF$ lie on the unique rational $(-2)$ curve in 
$\mathbb{P}(V)$.
\item At most $3$ points in $\bF$ lie on an irreducible curve $D
  \subset \mathbb{P}(V)$, $D \in |\mathcal{O}_{V}(1)\otimes
  \gamma^{*}\mathcal{O}_{C}(1)|$. 
\end{itemize}

\
\end{minipage}
& \\
\hline
\hline
$\mathcal{O}_{C}(-2)\oplus \mathcal{O}_{C}(2)$ &
\begin{minipage}[c]{3.8in} \addtolength{\baselineskip}{-5pt}
\

\begin{itemize}
\item No  points in $\bF$ lie on the unique rational $(-4)$ curve in 
$\mathbb{P}(V)$.
\end{itemize}

\
\end{minipage}
& ${\displaystyle \frac{4}{5} \leq q < 1}$ \\
\hline\hline
\end{tabular}
}
\end{center}
\caption{Stability configurations for the degree $0$ bundles
  $\mathcal{O}_{C}\oplus \mathcal{O}_{C}$, $\mathcal{O}_{C}(-1)\oplus
  \mathcal{O}_{C}(1)$, and $\mathcal{O}_{C}(-2)\oplus
  \mathcal{O}_{C}(2)$ .}
\label{table:deg.0.stable}
\end{table}
\

\noindent
Note also that the characterization of stability in Lemma~\ref{lem:stability} 
implies that for weights $(0,q)$ in the interior of a chamber 
all semi-stable parabolic bundles are strictly stable. Thus by the
general theory of the dependence of GIT quotients on a linearization 
(see e.g. \cite{mehta-seshadri,seshadri-book,boden-hu}) it follows that
the corresponding moduli space $N_{0}(0,q)$ will be smooth and
projective. We are now ready to descibe these varieties explicitly:

\begin{theo} \label{theo:shapes.of.moduli}
The moduli spaces $N_{0}(0,q)$ in the various
chambers have the following shape:
\begin{itemize}
\item If ${0 < q \leq 2/5}$, then the moduli
  space $N_{0}(0,q)$ is a smooth del Pezzo surface
  $dP_{4}$. \linebreak Specifically $N_{0}(0,q)$ is the del Pezzo surface
  obtained by  first blowing-up the projective plane  $S^{2}C$ at the five
  points $p_{1}, p_{2}, p_{3}, p_{4}, p_{5} \in C \subset S^{2}C$, and
  then contracting the strict transform of the conic $C \subset S^{2}C$.
\item If ${2/5 < q \leq 2/3}$,
  then the moduli
  space $N_{0}(0,q)$ is a smooth del Pezzo surface
  $dP_{5}$. \linebreak Specifically $N_{0}(0,q)$ is the del Pezzo surface
  obtained by blowing-up the projective plane  $S^{2}C$ at the five
  points $p_{1}, p_{2}, p_{3}, p_{4}, p_{5} \in C \subset S^{2}C$. 
\item If ${2/3 < q < 4/5}$,
  then the moduli
  space $N_{0}(0,q)$ is a projective plane. Specifically \linebreak 
$N_{0}(0,q) =
  S^{2}C \cong \mathbb{P}^{2}$.
\item If ${4/5 \leq  q \leq 1}$, then the moduli space
  $N_{0}(0,q)$ is a single point.
  \end{itemize}

The moduli spaces $N_{1}(0,q)$ are given by $N_{1}(0,q) \cong N_{0}(0,1-q)$.

\end{theo}

\

\medskip

\

\noindent
The proof of this theorem will take up  the remainder of Section~\ref{sec:git}.
We start with an explicit description of the moduli
corresponding to $q$ in the first chamber:

\begin{prop} \label{prop:dP4} Suppose $0 < q < 2/5$. Then 
  $N_{0}(0,q)$ is naturally isomorphic to the blow-up of
  $\mathbb{P}^{1}\times \mathbb{P}^{1}$ at three points. Equivalently
  $N_{0}(0,q)$ is a del Pezzo surface $dP_{4}$ obtained by
  blowing up four general points on $\mathbb{P}^{2}$.
\end{prop}
{\bfseries Proof.} If $0 < q < 2/5$ and $\mathcal{V} =
(V,\bF,(0,q))$ is stable, then $V \cong \mathcal{O}_{C}\oplus
\mathcal{O}_{C}$ and thus we can describe the moduli space
$N_{0}(0,q)$ as the GIT quotient
\[
\begin{aligned}
N_{0}(0,q) & = \left.\left.\left(\mathbb{P}^{1}_{p_{1}}\times
\mathbb{P}^{1}_{p_{2}}\times  \mathbb{P}^{1}_{p_{3}}\times
\mathbb{P}^{1}_{p_{4}}\times
\mathbb{P}^{1}_{p_{5}}\right)\right/\right/\mathbb{P}SL_{2}(\mathbb{C}) \\
&  =
\left.\left(\mathbb{P}^{1}_{p_{1}}\times
\mathbb{P}^{1}_{p_{2}}\times  \mathbb{P}^{1}_{p_{3}}\times
\mathbb{P}^{1}_{p_{4}}\times
\mathbb{P}^{1}_{p_{5}}\right)^{\text{stable}}
\right/\mathbb{P}SL_{2}(\mathbb{C}).
\end{aligned}
\]
Here $\mathbb{P}^{1}_{p_{i}}$ denotes the fiber of the ruled surface 
$\mathbb{P}(\mathcal{O}_{C}\oplus \mathcal{O}_{C})$ over $p_{i} \in
C$, the subset
\[
\left(\mathbb{P}^{1}_{p_{1}}\times
\mathbb{P}^{1}_{p_{2}}\times  \mathbb{P}^{1}_{p_{3}}\times
\mathbb{P}^{1}_{p_{4}}\times
\mathbb{P}^{1}_{p_{5}}\right)^{\text{stable}} \subset
\left(\mathbb{P}^{1}_{p_{1}}\times 
\mathbb{P}^{1}_{p_{2}}\times  \mathbb{P}^{1}_{p_{3}}\times
\mathbb{P}^{1}_{p_{4}}\times
\mathbb{P}^{1}_{p_{5}}\right)
\]
is the open subset parameterizung stable configurations $\bF$ of five
points in $\mathbb{P}(\mathcal{O}_{C}\oplus \mathcal{O}_{C})$ over $\{
p_{i}\}_{i=1}^{5} \subset C$,  and $\mathbb{P}SL_{2}(\mathbb{C})$ is
viewed as the group of all projective bundle automorphisms of $\gamma
: \mathbb{P}(\mathcal{O}_{C}\oplus \mathcal{O}_{C}) \to C$. 

Since $\mathbb{P}(\mathcal{O}_{C}\oplus
\mathcal{O}_{C}) = C\times \mathbb{P}^{1}$, then we can identify
all $\mathbb{P}^{1}_{p_{i}}$ with the fiber $\mathbb{P}^{1}$ and write
each point $F_{i}$ as a pair: $F_{i} = (p_{i},f_{i}) \in C\times
\mathbb{P}^{1}$. Thus the collection $\bF$ is given by the point
$(f_{1},f_{2},f_{3},f_{4},f_{5}) \in \mathbb{P}^{1}\times
\mathbb{P}^{1}\times \mathbb{P}^{1}\times \mathbb{P}^{1}\times
\mathbb{P}^{1}$ and the stability condition on $\bF$ (see
\ref{sssec:k=0}) becomes the condition that no three of the $f_{i}$'s
coincide.

Let now $I \subset \{ 1, 2, 3, 4, 5 \}$ be a subset with $\# I =
3$. Consider the subset  $\widetilde{U}_{I} \subset
\left(\mathbb{P}^{1}\right)^{\times 5}$ defined by  
\[
\widetilde{U}_{I} := \left\{\left.  (f_{1},f_{2},f_{3},f_{4},f_{5}) \in
\left(\mathbb{P}^{1}\right)^{\times 5} \, \right| \,
  \text{ $f_{a} \neq f_{b}$  for
  all $a\neq b \in I$ }\right\}.
\]
Then $\widetilde{U}_{I} \subset \left( \mathbb{P}^{1}\times
\mathbb{P}^{1}\times \mathbb{P}^{1}\times \mathbb{P}^{1}\times
\mathbb{P}^{1} \right)^{\text{stable}}$ is a
$\mathbb{P}SL_{2}(\mathbb{C})$-invariant open subset and $\{
\widetilde{U}_{I} \}$ is an open covering of $\left(
\mathbb{P}^{1}\times \mathbb{P}^{1}\times \mathbb{P}^{1}\times
\mathbb{P}^{1}\times \mathbb{P}^{1} \right)^{\text{stable}}$.
$\mathbb{P}SL_{2}(\mathbb{C})$ Thus $N_{0}(0,q)$ is glued out
of the quotients $U_{I} =
\widetilde{U}_{I}/\mathbb{P}SL_{2}(\mathbb{C})$. Suppose for
concretness we take the subset $\{ 1, 2, 3 \}$. Since for every
$(f_{1},f_{2},f_{3},f_{4},f_{5}) \in \widetilde{U}_{\{ 1, 2, 3 \}}$
the three points $f_{1},f_{2},f_{3} \in \mathbb{P}^{1}$ are distinct
it follows that we can fnd a unique element in
$\mathbb{P}SL_{2}(\mathbb{C})$ which sends $(0,1,\infty)$ to
$(f_{1},f_{2},f_{3})$. So each $\mathbb{P}SL_{2}(\mathbb{C})$-orbit in
$\widetilde{U}_{\{ 1,2,3 \}}$ contains a preferred point, namely the unique
point of the form $(0,1,\infty,*,*)$. This shows that the
$\mathbb{P}SL_{2}(\mathbb{C})$ action on $\widetilde{U}_{\{ 1, 2, 3
\}}$ is free, shows that the collection of all
$(0,1,\infty,f_{4},f_{5}) \in \widetilde{U}_{\{ 1, 2, 3 \}}$ is a
slice for the action and so for all $I$ we get:
\[
\begin{aligned}
\widetilde{U}_{I} & = U_{I} \times \mathbb{P}SL_{2}(\mathbb{C}) \\
U_{I} & = \left(\mathbb{P}^{1}\times
\mathbb{P}^{1}\right) - \left\{ (0,0), (1,1), (\infty,\infty)
\right\}. 
\end{aligned}
\]
It is not hard to see that when we glue all $U_{I}$'s together we get
the blow-up of $\mathbb{P}^{1}\times \mathbb{P}^{1}$ at the three
points $(0,0)$, $(1,1)$, $(\infty,\infty)$:

\begin{lem} \label{lem:dP4.iso} Let $X := \op{Bl}_{\{ (0,0), (1,1),
    (\infty,\infty)\}}(\mathbb{P}^{1}\times \mathbb{P}^{1})$, and let 
$\varphi : N_{0}(0,q) \dashrightarrow X$
be the birational map which on the open $U_{\{1,2,3\}} =
\left(\mathbb{P}^{1}_{p_{4}}\times 
\mathbb{P}^{1}_{p_{5}}\right) - \left\{ (0,0), (1,1), (\infty,\infty)
\right\}$ is given by the natural inclusion 
\[
\left(\mathbb{P}^{1}\times 
\mathbb{P}^{1}\right) - \left\{ (0,0), (1,1), (\infty,\infty)
\right\} \subset \op{Bl}_{\{ (0,0), (1,1),
    (\infty,\infty)\}}(\mathbb{P}^{1}\times \mathbb{P}^{1}).
\]
Then $\varphi$ is an isomorphism.
\end{lem}
{\bfseries Proof.} Consider the permutation action of $S_{5}$ on
$(\mathbb{P}^{1})^{\times 5}$ and look at the action of the $5$-cycle 
$(12345) \in S_{5}$. The action of $S_{5}$ on
$(\mathbb{P}^{1})^{\times 5}$ commutes with the diagonal action of 
$\mathbb{P}SL_{2}(\mathbb{C})$ and so $(12345)$ induces a birational
automorphism of $U_{\{1,2,3\}} =
\left(\mathbb{P}^{1}\times 
\mathbb{P}^{1}\right) - \left\{ (0,0), (1,1), (\infty,\infty)
\right\}$ and hence a birational automorphism $\sigma$ of
$\mathbb{P}^{1}\times \mathbb{P}^{1}$. We will check that $\sigma$
extends to a biregular automorphism of $X$. 

\begin{figure}[!ht]
\begin{center}
\psfrag{x=0}[c][c][1][0]{{$x=0$}}
\psfrag{x=1}[c][c][1][0]{{$x=1$}}
\psfrag{x=inf}[c][c][1][0]{{$x=\infty$}}
\psfrag{y=0}[c][c][1][0]{{$y=0$}}
\psfrag{y=1}[c][c][1][0]{{$y=1$}}
\psfrag{y=inf}[c][c][1][0]{{$y=\infty$}}
\psfrag{x=y}[c][c][1][0]{{$x=y$}}
\epsfig{file=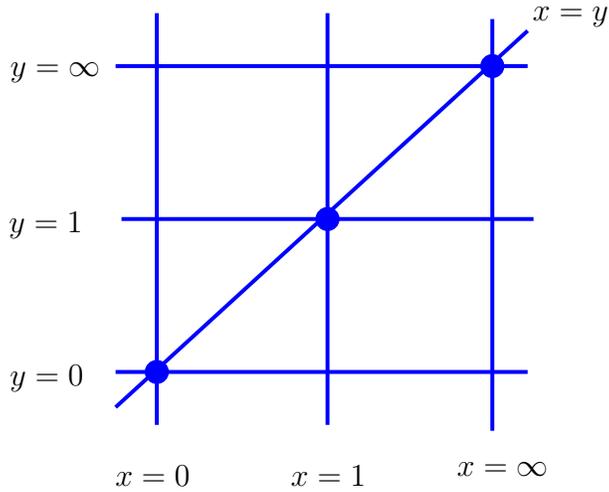,width=3in} 
\end{center}
\caption{The ten lines on $dP_{4}$.}\label{fig:tencurves} 
\end{figure}

Let
$[(0,1,\infty,f_{4},f_{5})] \in U_{\{1,2,3\}}$ be a generic
point. The $5$-cycle $(12345)$ sends $(0,1,\infty,f_{4},f_{5}) \in
\widetilde{U}_{\{1,2,3\}}$ to $(1,\infty,f_{4},f_{5},0)$, which by
genericity will again be in $\widetilde{U}_{\{1,2,3\}}$. The canonical
representative in the $\mathbb{P}SL_{2}(\mathbb{C})$-orbit of
$(0,1,\infty,f_{4},f_{5})$ is 
\[
\left(0,1,\infty,\frac{f_{5}-1}{f_{5}-f_{4}},\frac{1}{f_{4}} \right),
\]
i.e. $\sigma$ is the birational automorphism of $\mathbb{P}^{1}\times
\mathbb{P}^{1}$ given by
\begin{equation} \label{eq:sigma.formula}
\sigma(x,y) = \left( \frac{y-1}{y-x}, \frac{1}{x}\right) \text{ for } 
(x,y) \in \mathbb{P}^{1}\times 
\mathbb{P}^{1}.
 \end{equation}
The surface $X$ is the blow-up of $\mathbb{P}^{1}\times
\mathbb{P}^{1}$  at the three points $(0,0)$, $(1,1)$,
$(\infty,\infty)$ and so is a del Pezzo surface $dP_{4}$. It contains
$10$ lines ($= (-1)$-curves). From the point of view of $\mathbb{P}^{1}\times
\mathbb{P}^{1}$ these curves are the three exceptional divisors
$E_{0}$, $E_{1}$, $E_{\infty}$ corresponding to $(0,0)$, $(1,1)$,
$(\infty,\infty)$, and the strict transforms of the seven curves $x=0$,
$x=1$, $x=\infty$, $y = 0$, $y=1$, $y=\infty$, $x=y$ - see
Figure~\ref{fig:tencurves}.

\

\noindent
From the formula \eqref{eq:sigma.formula} we see that $\sigma$ acts on
$X$ by mapping
\[
\xymatrix@1@R-1pc{
E_{0} \ar[r]^-{\sigma} &  E_{\infty}   \ar[r]^-{\sigma} & (y=0)
  \ar[r]^-{\sigma} & (x=y) \ar[r]^-{\sigma} & (x=\infty)
  \ar[r]^-{\sigma} & E_{0}, \\
E_{1} \ar[r]^-{\sigma} &  (y=1)  \ar[r]^-{\sigma} & (x=0)
  \ar[r]^-{\sigma} & (y=\infty) \ar[r]^-{\sigma} & (x=1)
  \ar[r]^-{\sigma} & E_{1},
}
\]
and so $\sigma$ acts on $X$ as a biregular automorphism of order $5$.

Next observe that the flipping involution $\tau : \mathbb{P}^{1}\times
\mathbb{P}^{1} \to \mathbb{P}^{1}\times
\mathbb{P}^{1}$, given by $\tau(x,y) = (y,x)$ preserves the three
points $(0,0)$, $(1,1)$,
$(\infty,\infty)$ and so extends to a biregular involution of $X$. The
flipping involution and the $5$-cycle generate the group of $S_{5}$ of
birational automorphisms of $\mathbb{P}^{1}\times
\mathbb{P}^{1} \to \mathbb{P}^{1}$ and so their extensions $\sigma,
\tau \in \op{Aut}(X)$ generate $S_{5} \subset \op{Aut}(X)$. This shows
that all permutations of the five factors in $\left( \mathbb{P}^{1}
\right)^{\times 5}$ extend to biregular automorphisms of $X$. 

Since the rational map $\varphi$ is $S_{5}$-equivariant this implies
that for every $\xi \in S_{5}$ the restriction of $\varphi$ to
$U_{\{\xi(1),\xi(2),\xi(3)\}}$ is $\xi\circ
\varphi_{|U_{\{1,2,3\}}}\circ 
\xi^{-1}$. But when $\xi$ runs over $S_{5}$ the opens 
$U_{\{\xi(1),\xi(2),\xi(3)\}}$ cover all of $N_{0}(0,q)$ and
therefore $\varphi$ is an isomorphism. This proves the lemma and the
proposition. \ \hfill $\Box$

\

\medskip

\noindent
Next we give a description of the moduli spaces
corresponding to $q$ in the  second chamber:

\begin{prop} \label{prop:dP5} 
Suppose $2/5 < q \leq 2/3$. Then
  $N_{0}(0,q)$ is naturally isomorphic to the blow-up of
  $\mathbb{P}^{1}\times \mathbb{P}^{1}$ at four points. Equivalently
  $N_{0}(0,q)$ is a del Pezzo surface $dP_{5}$ obtained by
  blowing up five general points on $\mathbb{P}^{2}$. 
\end{prop}
{\bfseries Proof.} By the general theory of variations of GIT quotients
\cite{dolgachev-hu,thaddeus} or from the analysis in \cite{boden-hu}
we know that if we vary the parabolic weights $(0,q)$ and cross
from the chamber $0 < q < 2/5$ to the chamber $2/5 <
q < 2/3$ across the wall $q = \frac{2}{5}$, then the
moduli space $N_{0}(0,q)$ undergoes a flip. From
Table~\ref{table:deg.0.stable} we see that the classes of semi-stable
or stable
parabolic bundles do not change if $q$ varies in the interval $0 <
q \leq 2/5$ and so for  $q$ in this chamber we have 
$N_{0}(0,q) = X =
\op{Bl}_{(0,0),(1,1),(\infty,\infty)}(\mathbb{P}^{1}\times
\mathbb{P}^{1})$. In particular if $2/5 <
q \leq 2/3$ we must have a birational morphism
$\varepsilon : N_{0}(0,q) \to X$. But both
$N_{0}(0,q)$ and $X$ parametrize stable parabolic bundles. So
they are both smooth surfaces and thus $\varepsilon :
N_{0}(0,q) \to X$ must contract consequtively 
finitely many smooth rational
curves of self intersection $(-1)$.

Table~\ref{table:deg.0.stable} also shows that when $q$ crosses
  into the chamber $2/5 < q <  2/3$ the moduli space
  $N_{0}(0,q)$ becomes a union
\[
N_{0}(0,q) = N_{0}(0,q)_{k=0} \sqcup N_{0}(0,q)_{k=1}
\]
where $N_{0}(0,q)_{k=0} \subset N_{0}(0,q)$ parametrizes parabolic
bundles whose underlying bundle is $\mathcal{O}_{C}\oplus
\mathcal{O}_{C}$, and $N_{0}(0,q)_{k=1}\subset N_{0}(0,q)$
parametrizes parabolic bundles whose underlying bundle is
$\mathcal{O}_{C}(-1)\oplus \mathcal{O}_{C}(1)$.  Furthermore, from
Table~\ref{table:deg.0.stable} it follows that the only effect the
wall crossing has on parabolic bundles of the form
$(\mathcal{O}_{C}\oplus \mathcal{O}_{C}, \bF, 0,q)$ is that an
additional stability constraint is imposed: the five points in $\bF
\subset C\times \mathbb{P}^{1}$ can not lie on an irreducible $(1,1)$
curve in the surface $C\times \mathbb{P}^{1}$. This is a closed
condition on five-tuples of points and so for $2/5 < q \leq 2/3$ we
have that $N_{0}(0,q)_{k=0} \subset X$ is an open subset.

To describe this open subset explicitly we need to analyze the
additional stability constraint in more detail.  Write again $\bF = \{
(p_{i},f_{i} \}_{i = 1}^{5} \subset C\times \mathbb{P}^{1}$ and for
concreteness let us suppose that $(f_{1},f_{2},f_{3},f_{4},f_{5}) \in
\widetilde{U}_{\{1,2,3\}}$. Then the same reasoning as before shows
that without a loss of generality we may assume that $\bF = \{
(p_{1},0),(p_{2},1), (p_{3},\infty), (p_{4},f_{4}), (p_{5},f_{5})
\}$. Now consider the unique isomorphism $C \cong \mathbb{P}^{1}$
which sends $p_{1}$ to $0$, $p_{2}$ to $1$, and $p_{3}$ to
$\infty$. This identifies $C\times \mathbb{P}^{1}$ with
$\mathbb{P}^{1}\times \mathbb{P}^{1}$, and under this identification
$\bF$ becomes the five-tuple of points
\[
\{ (0,0), (1,1), (\infty,\infty), (p_{4},f_{4}), (p_{5},f_{5})\}
\subset \mathbb{P}^{1}\times \mathbb{P}^{1}.
\]
But the three points $(0,0)$, $(1,1)$, $(\infty,\infty)$ already lie
on the diagonal in $\mathbb{P}^{1}\times \mathbb{P}^{1}$, and so $\bF$
will be contained in an irreducible $(1,1)$ curve if and only if it is
contained in the diagonal. This means that $\bF$ will be unstable if
and only if $f_{4} = p_{4}$ and $f_{5} = p_{5}$. Therefore:
\begin{itemize}
\item  the open subset $U_{\{ 1,2,3 \}} \subset 
\mathbb{P}^{1}_{p_{4}}\times
\mathbb{P}^{1}_{p_{5}}$ contains a unique point
$(p_{4},p_{5}) \in \mathbb{P}^{1}_{p_{4}}\times
\mathbb{P}^{1}_{p_{5}}$ which becomes unstable for 
$q$ in the chamber $2/5 < q < 2/3$;
\item this point corresponds to a $\bF = \{ (p_{i}, f_{i}) \}_{i =
 1}^{5}$ in which all $f_{i}$ are distinct.
\end{itemize}
But if $I \subset \{1,2,3,4,5\}$ is any
subset of cardinality $3$ and if we have an unstable \linebreak 
$\bG = \{
(p_{i},g_{i}) \}_{i = 1}^{5}$ such that 
$(g_{1},g_{2},g_{3},g_{4},g_{5}) \in
\widetilde{U}_{I}$, then we can repeat the above reasoning to conlude
that all $g_{i}$ must be distinct. In particular
$(g_{1},g_{2},g_{3},g_{4},g_{5}) \in 
\cap_{J} \widetilde{U}_{J} \subset \widetilde{U}_{\{ 1,2,3\}}$ and
hence the image $\varphi(g_{1},g_{2},g_{3},g_{4},g_{5})$ will be
contained in $U_{\{ 1,2,3\}}$. This shows that
$\varphi(g_{1},g_{2},g_{3},g_{4},g_{5}) = (p_{4},p_{5})$, and so
$(p_{4},p_{5})$ will be the unique point in the
surface $\op{Bl}_{(0,0),(1,1),(\infty,\infty)}(\mathbb{P}^{1}\times
\mathbb{P}^{1})$ that becomes unstable when we cross the wall $b-a =
\frac{2}{5}$. This shows that when $b-a \in
\left(\frac{2}{5},\frac{2}{3}\right)$ we have
\[
N_{0}(0,q)_{k=0} =
\op{Bl}_{(0,0),(1,1),(\infty,\infty)}(\mathbb{P}^{1}\times 
\mathbb{P}^{1}) - \{(p_{4},p_{5})\} = X - \op{pt}.
\]
This shows that $N_{0}(0,q)$ is the blow-up of $X$ at finitely
 many smooth points which are infinitesimally close to
 $(p_{4},p_{5})$.

On the other hand the union of exceptional divisors for the birational
morphism \linebreak $N_{0}(0,q) \to X$ is precisely
$N_{0}(0,q)_{k=1}$ and so the proposition will be proven if we
can show that $N_{0}(0,q)_{k=1}$ is irreducible. To argue the
irreducibility we will describe
  $N_{0}(0,q)_{k=1}$ explicitly as a quotient.

Let $V = \mathcal{O}_{C}(-1)\oplus \mathcal{O}_{C}(1)$. Then
$N_{0}(0,q)_{k=1}$ is the quotient
\[
N_{0}(0,q)_{k=1} = \left.\left( \mathbb{P}(V)_{p_{1}}\times
\mathbb{P}(V)_{p_{2}}\times  \mathbb{P}(V)_{p_{3}}\times
\mathbb{P}(V)_{p_{4}}\times
\mathbb{P}(V)_{p_{5}}\right)^{\op{stable}}\right/\op{Aut}(\mathbb{P}(V)),
\]
where $\left( \mathbb{P}(V)_{p_{1}}\times
\mathbb{P}(V)_{p_{2}}\times  \mathbb{P}(V)_{p_{3}}\times
\mathbb{P}(V)_{p_{4}}\times
\mathbb{P}(V)_{p_{5}}\right)^{\op{stable}}$ denotes the locus of
stable five tuples of points, and $\op{Aut}(\mathbb{P}(V))$ 
is the group of global
projective bundle automorphisms  of $\mathbb{P}(V)$ acting by
evaluation on each fiber
$\mathbb{P}(V)_{p_{i}}$. 

Note that in contrast with the $\mathcal{O}_{C}\oplus \mathcal{O}_{C}$
case the group of projective bundle automorphisms
\[
\op{Aut}(\mathbb{P}(V)) = \mathbb{C}^{\times} \ltimes
H^{0}(C,\mathcal{O}_{C}(2)) 
\]
is not reductive but solvable. Nevertheless this group acts freely and
with closed orbits on $\left( \mathbb{P}(V)_{p_{1}}\times
\mathbb{P}(V)_{p_{2}}\times \mathbb{P}(V)_{p_{3}}\times
\mathbb{P}(V)_{p_{4}}\times
\mathbb{P}(V)_{p_{5}}\right)^{\op{stable}}$ which allows us to
construct $N_{0}(0,q)_{k=1}$ as the above quotient. To analyze
the structure of $N_{0}(0,q)_{k=1}$ we will again construct an
explicit slice to the action.  

From Table~\ref{table:deg.0.stable} we see that a five tuple
$(f_{1},f_{2},f_{3},f_{4},f_{5})$ will belong to  \linebreak 
$\left(
\mathbb{P}(V)_{p_{1}}\times 
\mathbb{P}(V)_{p_{2}}\times \mathbb{P}(V_{p_{3}})\times
\mathbb{P}(V)_{p_{4}}\times
\mathbb{P}(V)_{p_{5}}\right)^{\op{stable}}$ if and only if
\begin{itemize}
\item none of the five points $\{ (p_{i},f_{i})\}_{i=1}^{5}$ lie on
  the unique $(-2)$ curve in $\mathbb{P}(V) \cong \mathbb{F}_{2}$;
\item the five points $\{ (p_{i},f_{i})\}_{i=1}^{5}$ do not lie on
  any irreducible curve in the linear system
  $|\mathcal{O}_{V}(1)\otimes \gamma^{*}\mathcal{O}_{C}(1)|$.  
\end{itemize}
To understand better the first of these two conditions consider 
the identification 
$\mathbb{P}(V)$ with $\mathbb{P}(V\otimes \mathcal{O}_{C}(1)) =
\mathbb{P}(\mathcal{O}_{C}\oplus \mathcal{O}_{C}(2))$.  The
$(-2)$ curve in $\mathbb{P}(V)$ is then given as the section of
$\gamma : \mathbb{P}(V) \to C$ which to each point $p \in C$ assigns
the point of $\mathbb{P}(V)_{p}$ which corresponds to the line
$(\mathcal{O}_{C})_{p} \subset (\mathcal{O}_{C}\oplus
\mathcal{O}_{C}(2))_{p}$. Thus  $\left(
\mathbb{P}(V_{p_{1}})\times 
\mathbb{P}(V_{p_{2}})\times \mathbb{P}(V_{p_{3}})\times
\mathbb{P}(V_{p_{4}})\times
\mathbb{P}(V_{p_{5}})\right)^{\op{stable}}$ is contained in the open subset
\[
\prod_{i=1}^{5}
\mathcal{O}_{C}(2)_{p_{i}} \cong \mathbb{C}^{5} \subset
\mathbb{P}(V_{p_{1}})\times  
\mathbb{P}(V_{p_{2}})\times \mathbb{P}(V_{p_{3}})\times
\mathbb{P}(V_{p_{4}})\times
\mathbb{P}(V_{p_{5}}).
\]
Here we have embedded the total space of $\mathcal{O}_{C}(2)$ in
$\mathbb{P}(V)$ by sending a point $z \in (\mathcal{O}_{C}(2))_{p}$ to
the point $(1:z) \in \mathbb{P}((\mathcal{O}_{C}\oplus
\mathcal{O}_{C}(2))_{p}) = \mathbb{P}(V)_{p}$.

The second condition characterizes the stable points
$(f_{1},f_{2},f_{3},f_{4},f_{5}) \in \prod_{i=1}^{5}
\mathcal{O}_{C}(2)_{p_{i}} \cong \mathbb{C}^{5}$ as the points which
do not belong to the orbit of $0 \in \mathbb{C}^{5}$ under the action
of the subgroup $H^{0}(C,\mathcal{O}_{C}(2)) \subset
\op{Aut}(\mathbb{P}(V))$. But if $g \in 
H^{0}(C,\mathcal{O}_{C}(2))$ is a global section then the
corresponding automorphism of $\mathbb{P}(V)$ preserves the open
subset $\op{tot}(\mathcal{O}_{C}(2)) \subset \mathbb{P}(V)$ and acts
on each line $\mathcal{O}_{C}(2)_{p}$ as translation by the vector
$g(p) \in \mathcal{O}_{C}(2)_{p}$. Therefore such a $g$ acts on
$\prod_{i=1}^{5} 
\mathcal{O}_{C}(2)_{p_{i}} \cong \mathbb{C}^{5}$ as translation by the
vector $(g(p_{1}),g(p_{2}),g(p_{3}),g(p_{4}),g(p_{5}))$.

In particular we see that
\[
\left(
\mathbb{P}(V_{p_{1}})\times 
\mathbb{P}(V_{p_{2}})\times \mathbb{P}(V_{p_{3}})\times
\mathbb{P}(V_{p_{4}})\times
\mathbb{P}(V_{p_{5}})\right)^{\op{stable}} = \prod_{i=1}^{5}
\mathcal{O}_{C}(2)_{p_{i}} - H^{0}(C,\mathcal{O}_{C}(2)) =
\mathbb{C}^{5} - \mathbb{C}^{3},
\] 
and that the group $\op{Aut}(\mathbb{P}(V)) =
\mathbb{C}^{\times}\ltimes \mathbb{C}^{3}$ acts on this space in the
obvious manner: an element in $\mathbb{C}^{\times}$ acts on
$\mathbb{C}^{5} - \mathbb{C}^{3}$ by scaling, and a vector in
$\mathbb{C}^{3}$ acts on $\mathbb{C}^5 - \mathbb{C}^{3}$ by
translation.
Therefore we get 
\[
\begin{aligned}
N_{0}(0,q)_{k=1} & = \left.\left( \frac{\mathbb{C}^{5} -
  \mathbb{C}^{3}}{\mathbb{C}^{3}} \right)\right/ \mathbb{C}^{\times}
  \\[1pc] 
& = \left.\left( \frac{\mathbb{C}^{5}}{\mathbb{C}^{3}} - \{ 0 \}
  \right)\right/ \mathbb{C}^{\times} \\[1pc]
& \cong \mathbb{P}^{1}.
\end{aligned}
\]
This proves the irreducibility of $N_{0}(0,q)_{k=1}$ and
finishes the proof of the proposition.
\ \hfill $\Box$

\

\noindent
Since the moduli space $N_{0}(0,q)$ corresponding to a choice of
parabolic weights in the second chamber will be our main object of
study, it will be useful to have a better and intrinsic understanding of
its geometric structure. For future reference we give three equivalent
synthetic descriptions of this moduli space. First we will need some
notation. 

Consider the plane
$S^{2}C \cong \mathbb{P}^{2}$ with its natural diagonal conic $C
\subset S^{2}C$. Let $t_{i} \subset S^{2}C$ be the line tangent to the
conic $C$ at the point $p_{i} \in C$.

\begin{prop} \label{prop:synthetic.N0} Suppose $\frac{2}{5} < q <
  \frac{2}{3}$.  Then 
\begin{description}
\item[{\bfseries (i)}] The moduli space $N_{0}(0,q)$ is
  canonically isomorphic to the blow-up of the quadric surface
  $C\times C$ at the four points $(p_{1},p_{1})$, $(p_{2},p_{2})$,
  $(p_{3},p_{3})$, and $(p_{4},p_{5})$.
\item[{\bfseries (ii)}] The moduli space $N_{0}(0,q)$ is
  equipped with a natural morphism $\beta_{0} : N_{0}(0,q) \to
  S^{2}C$ which is uniquely characterized by the following properties:
\begin{itemize}
\item $\beta_{0}^{*} \mathcal{O}_{S^{2}C}(1) =
  K_{N_{0}(0,q)}^{\otimes -2}$
\item $\beta_{0}$ is a finite Galois cover with Galois group
  $(\mathbb{Z}/2)^{4}$.
\item $\beta_{0}$ is unramified over $S^{2}C - \cup_{i = 1}^{5}
  t_{i}$. 
\end{itemize}
In particular $N_{0}(0,q)$ is 
isomorphic to the normalization
of the fiber product 
\[
Y_{12}\times_{S^{2}C} Y_{23}\times_{S^{2}C}
Y_{34}\times_{S^{2}C} Y_{45},
\] 
where $Y_{ij} \to S^{2}C$
is the unique (square root) double cover branched along
$t_{i}+t_{j}$. 
\item[{\bfseries (iii)}] The moduli space $N_{0}(0,q)$ is
  isomorphic to the blow-up of $S^{2}C$ at the five points $p_{1},
  p_{2}, p_{3}, p_{4}, p_{5} \in C \subset S^{2}C$. 
\end{description}
\end{prop}
{\bf Proof.} {\bfseries (i)} This is just Proposition~\ref{prop:dP5}.

\

\noindent
{\bfseries (ii)}  Again we will choose a coordinate on $C \cong
\mathbb{P}^{1}$ so that  $p_{1}$, $p_{2}$, and $p_{3}$ are mapped to
$0$, $1$, and $\infty$. Consider the surface $X =
\op{Bl}_{(0.,0),(1,1),(\infty,\infty)}(\mathbb{P}^{1}\times
\mathbb{P}^{1})$ i.e. the moduli space $N_{0}(0,q)$ for $0< b-a \leq
\frac{2}{5}$. Recall that by its modular interpretation $X$ is
identified with the quotient
\[
X = \left.\left( \mathbb{P}^{1}\times  \mathbb{P}^{1}\times
\mathbb{P}^{1}\times  \mathbb{P}^{1}\times
\mathbb{P}^{1}\right)^{\op{stable}}\right/ \mathbb{P}SL_{2}(\mathbb{C})
\] 
where the stability corresponds to the chamber of weights $0< b-a \leq
\frac{2}{5}$. 

However we have a natural rational map
\[
\xymatrix@1{ \varphi : \mathbb{P}^{1}\times  \mathbb{P}^{1}\times
\mathbb{P}^{1}\times  \mathbb{P}^{1}\times
\mathbb{P}^{1} \ar@{-->}[r] & \mathbb{P}^{5}}
\]
defined as follows.

\begin{figure}[!ht]
\begin{center}
\psfrag{P1f}[c][c][1][0]{{$\mathbb{P}^{1}_{f}$}}
\psfrag{P1p}[c][c][1][0]{{$\mathbb{P}^{1}_{p}$}}
\psfrag{(p1,f1)}[c][c][1][0]{{$(p_{1},f_{1})$}}
\psfrag{(p4,f4)}[c][c][1][0]{{$(p_{4},f_{4})$}}
\psfrag{ramification}[c][c][0.8][0]{{\begin{minipage}[c]{0.6in}
      \addtolength{\baselineskip}{-5pt} \begin{center}
ramification \\ point \end{center} \end{minipage}}}
\psfrag{branch point}[c][c][0.8][0]{{branch point}}
\psfrag{p1}[c][c][1][0]{{$p_{1}$}}
\psfrag{p2}[c][c][1][0]{{$p_{2}$}}
\psfrag{p3}[c][c][1][0]{{$p_{3}$}}
\psfrag{p4}[c][c][1][0]{{$p_{4}$}}
\psfrag{p5}[c][c][1][0]{{$p_{5}$}}
\psfrag{(1,2) curve}[c][c][1][0]{{$(1,2)$ curve}}
\epsfig{file=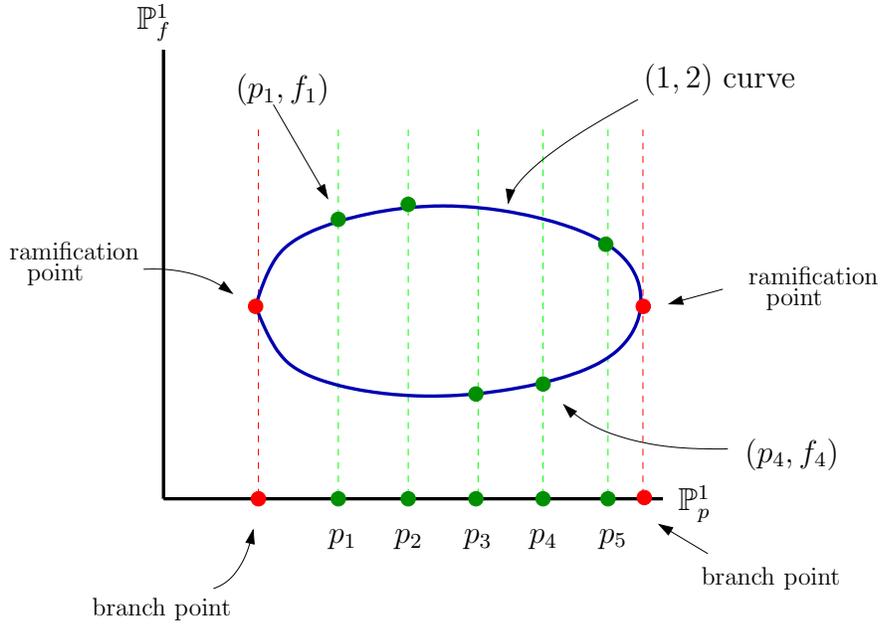,width=4.5in} 
\end{center}
\caption{The $(1,2)$ curve through a five-tuple of points}
\label{fig:twisted.cubic} 
\end{figure}

\begin{itemize}
\item For any five tuple $(f_{1},f_{2},f_{3},f_{4},f_{5}) \in
  \left(\mathbb{P}^{1}\right)^{\times 5}$ consider the five tuple of
  points $\bF = \{ (0,f_{1}), (1,f_{2}), (\infty,f_{3}),
  (p_{4},f_{4}), (p_{5},f_{5}) \}$ in the surface $C\times
  \mathbb{P}^{1} \cong \mathbb{P}^{1}\times\mathbb{P}^{1}$. 
\item Consider the linear system $|\mathcal{O}(1,2)| \cong \mathbb{P}^{5}$ on
  $\mathbb{P}^{1}\times\mathbb{P}^{1}$.  For a generic choice of
  $(f_{1},f_{2},f_{3},f_{4},f_{5})$, there will be a unique curve (see
  Figure~\ref{fig:twisted.cubic}) in
  this linear system passing through the five points $\bF$ and so we
  define 
\[
\xymatrix@R-1pc{
\varphi : & \mathbb{P}^{1}\times  \mathbb{P}^{1}\times
\mathbb{P}^{1}\times  \mathbb{P}^{1}\times
\mathbb{P}^{1} \ar@{-->}[r] & \mathbb{P}^{5}  =
\mathbb{P}(H^{0}(\mathbb{P}^{1}\times\mathbb{P}^{1},\mathcal{O}(1,2))) \\
 & (f_{1},f_{2},f_{3},f_{4},f_{5}) \ar[r] &
\left(\text{\begin{minipage}[c]{1.8in} the unique $(1,2)$ curve in
    $\mathbb{P}^{1}\times\mathbb{P}^{1}$ through the five tuple of
    points $\bF$. \end{minipage}}\right).
}
\]
\end{itemize}

\

\medskip

\noindent
It is not hard to describe the fibers of the rational map $\varphi$
geometrically. Indeed, if $\Phi \in |\mathcal{O}(1,2)|$ is any member
of the linear system, then $\varphi^{-1}(\Phi)$ consists of five tuples
of points $(f_{1},f_{2},f_{3},f_{4},f_{5}) \in
\left(\mathbb{P}^{1}\right)^{\times 5}$ such that $\bF \subset
\Phi$. But the points $p_{1}$, $p_{2}$, $p_{3}$, $p_{4}$, $p_{5}$ are
fixed, and so if $\Phi$ is a smooth curve which is not tangent to any
of the rulings $\{ p_{i} \}\times \mathbb{P}^{1}$, then there are
exactly two points in $\Phi$ mapping to each $p_{i}$. Therefore for
such $\Phi$'s we get that $\varphi^{-1}(\Phi)$ consists of the $2^{5}$
five tuples of points corresponding to choosing one of the two
preimages of each $p_{i}$ in $\Phi$.

\begin{figure}[!ht]
\begin{center}
\psfrag{P1f}[c][c][1][0]{{$\mathbb{P}^{1}_{f}$}}
\psfrag{P1p}[c][c][1][0]{{$\mathbb{P}^{1}_{p}$}}
\psfrag{ramification}[c][c][0.8][0]{{\begin{minipage}[c]{0.6in} 
\addtolength{\baselineskip}{-5pt}
\begin{center}
ramification \\ point \end{center} \end{minipage}}}
\psfrag{branch point}[c][c][0.8][0]{{branch point}}
\psfrag{Phi}[c][c][1][0]{{$\Phi$}}
\psfrag{g(Phi)}[c][c][1][0]{{$g(\Phi)$}}
\epsfig{file=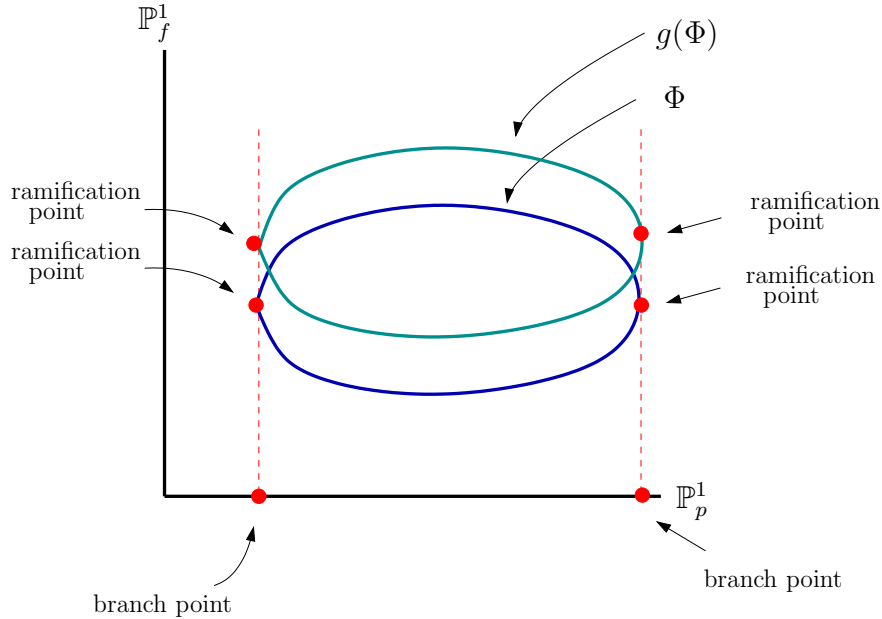,width=4.5in} 
\end{center}
\caption{The image of a $(1,2)$-curve 
under $g \in \mathbb{P}SL_{2}(\mathbb{C})$}
\label{fig:slide.cubic} 
\end{figure}

Next observe that we also have a natural rational map (see
Figure~\ref{fig:twisted.cubic}) 
\[
\xymatrix@R-1pc{
\mathfrak{b}_{0} : & \mathbb{P}^{5} = \mathbb{P}(H^{0}(\mathcal{O}(1,2)))
\ar@{-->}[r] & \mathbb{P}^{2} = S^{2}C \\
& \Phi \ar[r] & \left(\text{\begin{minipage}[c]{1.6in}
   the two branch
    points of the projection\\[-1pc]
\[
\Phi \stackrel{2;1}{\to} C =
    \mathbb{P}^{1}_{p}.
\]
\end{minipage}}\right)
}
\]
Furthermore note that the diagonal action of $\mathbb{P}SL_{2}(\mathbb{C})$ on
$\left(\mathbb{P}^{1}\right)^{5}$ corresponds to the action on
$C\times \mathbb{P}^{1} = \mathbb{P}^{1}\times \mathbb{P}^{1}$ which
is trivial on the first factor and tautological on the second
one. This action preserves the linear system $|\mathcal{O}(1,2)|$, and
so for every $\Phi \in  |\mathcal{O}(1,2)|$ and every $g \in
\mathbb{P}SL_{2}(\mathbb{C})$ the image $g(\Phi)$ will be another
$(1,2)$ curve in $\mathbb{P}^{1}\times \mathbb{P}^{1}$. But if $(p,f)
\in \Phi$ we have that $g(p,f) = (p,g(f)) \in g(\Phi)$ and so $\Phi
\to C$
and $g(\Phi) \to C$ have the same branch points  as depicted in 
Figure~\ref{fig:slide.cubic}.  

\

\noindent
In particular $\mathfrak{b}_{0}\circ g = \mathfrak{b}_{0}$ for all $g \in
\mathbb{P}SL_{2}(\mathbb{C})$. This shows that $\mathfrak{b}_{0}$ factors through
the quotient by $\mathbb{P}SL_{2}(\mathbb{C})$ and so we have a commutative
diagram
\[
\xymatrix{
\left(\mathbb{P}^{1}\right)^{5} \ar@{-->}[r]^-{\varphi} \ar@{-->}[d] &
\mathbb{P}^{5} \ar@{-->}[d]^-{\mathfrak{b}_{0}} \\
N_{0}(0,q) \ar@{-->}[r]_-{\beta_{0}} & 
\mathbb{P}^{2} = S^{2}C
}
\]
where $N_{0}(0,q) \cong dP_{4}$ is the moduli space
corresponding to the chamber $0 < b-a \leq \frac{2}{5}$ and 
$\beta_{0} : N_{0}(0,q) \dashrightarrow S^{2}C$ is the
induced rational map. 

\

\begin{lem} \label{lem:beta0.morphism}
Suppose $2/5 < q \leq 2/3$. Then the rational map $\beta_{0}$ extends
to a finite 
morphism 
\[
\beta_{0} : N_{0}(0,q) \to S^{2}C.
\]
\end{lem}
{\bfseries Proof.}
We can analyze the indeterminacy locus of $\beta_{0}$ from the
definition. First note that $\varphi$ is well defined at a point
$(f_{1},f_{2},f_{3},f_{4},f_{5}) \in \left(\mathbb{P}^{1}\right)^{5}$
if and only if the corresponding five-tuple of points $\bF = \{
(p_{1},f_{1}),(p_{2},f_{2}),(p_{3},f_{3}),(p_{4},f_{4}),(p_{5},f_{5})\}
\subset C\times \mathbb{P}^{1}$ imposes five independent conditions on
the linear system $|\mathcal{O}(1,2)|$. Now suppose that
$(f_{1},f_{2},f_{3},f_{4},f_{5})$ is $(0,q)$-stable for $0 < q \leq
2/5$. Note that all five points in $\bF$ are distinct and by
stability no two points in $\bF$ can lie on a $(1,0)$ ruling of
$C\times \mathbb{P}^{1}$. Thus there are three distinct points in
$\bF$, say $(p_{1},f_{1})$, $(p_{2},f_{2})$, $(p_{3},f_{3})$, that
impose independent conditions on $|\mathcal{O}(1,2)|$, and so we get a
well defined morphism from the blow up of $C\times \mathbb{P}^{1}$ at
$(p_{1},f_{1})$, $(p_{2},f_{2})$, $(p_{3},f_{3})$ to the projective
plane $\mathbb{P}^{2}$ given by the linear system of all $(1,2)$
curves passing through the three points $(p_{1},f_{1})$,
$(p_{2},f_{2})$, $(p_{3},f_{3})$. This map has degree $1$ and so is a
birational morphism. It contracts four disjoint $(-1)$ curves in
$\op{Bl}_{\{ (p_{1},f_{1}),(p_{2},f_{2}),(p_{3},f_{3}) \}}(C\times
\mathbb{P}^{1})$: the strict transforms of the three $(1,0)$ rulings
passing through eacho of $(p_{1},f_{1})$, $(p_{2},f_{2})$,
$(p_{3},f_{3})$, and the strict transform of the unique irreducible
$(1,1)$ curve through all three points $(p_{1},f_{1})$,
$(p_{2},f_{2})$, $(p_{3},f_{3})$. Now the $(1,2)$ curves through
$(p_{1},f_{1})$, $(p_{2},f_{2})$, $(p_{3},f_{3})$ are pullbacks of
lines in $\mathbb{P}^{2}$ and so the five points $\bF$ will impose
independent conditions if and only if $(p_{4},f_{4})$ and
$(p_{5},f_{5})$ map to distinct points in $\mathbb{P}^{2}$. But the
$(p_{4},f_{4})$ and $(p_{5},f_{5})$ can only map to the same point if
they both lie on one of the rational curves contracted by the
map. Since by stability they can not lie on any of the three
$(1,0)$-rulings it follows that  $(p_{4},f_{4})$ and $(p_{5},f_{5})$
must lie on the $(1,1)$-curve. As we saw in the proof of
Proposition~\ref{prop:dP5} this happens precisely when
$(f_{1},f_{2},f_{3},f_{4},f_{5})$ is $(0,q)$-stable for $0 < q \leq
2/5$ but is {\em not} $(0,q)$-stable for $2/5 < q \leq
2/3$. In other words, the map $\phi$ is a morphism over the
open subset of $\left(\mathbb{P}^{1}\right)^{5}$ of all points that are 
$(0,q)$-stable for $2/5 < q \leq
2/3$.

Suppose next $\Phi \subset C\times \mathbb{P}^{1}$ is a $(1,2)$-curve,
then $\mathfrak{b}_{0}$ is well defined at the point $\Phi \in
\mathbb{P}^{5}$ as long as the discriminant of the map $\Phi \to C$ is
not all of $C$. But this discriminant will be all of $C$ if and only
if $\Phi$ is of the form $2C\times \{ f \} + \{ p \}\times
\mathbb{P}^{1}$. In particular if $\Phi =
\phi(f_{1},f_{2},f_{3},f_{4},f_{5})$ we must have that at least $4$ of
the five points in $\bF$ lie on $C\times \{ f \}$. Therefore
$(f_{1},f_{2},f_{3},f_{4},f_{5})$ can not be stable for any choice of
the parabolic weights. This shows that if $\frac{2}{5} < b-a \leq
\frac{2}{3}$, then $\beta_{0}$ is a well defined morphism on
$N_{0}(0,q)_{k=0}$. Similar but easier reasoning shows that
$\beta_{0}$ is a well defined and injective morphism on the rational
curve $N_{0}(0,q)_{k=1}$.

This implies that if $\frac{2}{5} < b-a \leq \frac{2}{3}$, then
$\beta_{0} : N_{0}(0,q) \to \mathbb{P}^{2} = S^{2}C$ is a well
defined morphism.  Next we will argue that in this case $\beta_{0}$ is
a finite morphism. Indeed, since $N_{0}(0,q)$ is projective we
only need to check that $\beta_{0}$ is quasi-finite. Since $\beta_{0}$
is an inclusion of $N_{0}(0,q)_{k=1}$ in $\mathbb{P}^{2}$ this
amounts to checking that $\beta_{0}$ is quasi-finite over the open
$N_{0}(0,q)_{k=0} \subset N_{0}(0,q)$. By construction
$N_{0}(0,q)_{k=0}$ is the quotient of the open subset
$\left(\left(\mathbb{P}^{1}\right)^{5}\right)^{\op{stable}}$ in 
$\left(\mathbb{P}^{1}\right)^{5}$ consisting of $(0,q)$-stable points
for $2/5 < q \leq 2/3$ by the free diagonal action
of $\mathbb{P}SL_{2}(\mathbb{C})$. So to show that $\beta_{0}$ is
quasi-finite it suffices to check that eac positive dimensional
component of a fiber of $\mathfrak{b}_{0}\circ
\varphi_{|\left(\left(\mathbb{P}^{1}\right)^{5}\right)^{\op{stable}}}$
is a $\mathbb{P}SL_{2}(\mathbb{C})$-orbit. This is clear for the
fibers of $\mathfrak{b}_{0}$ and so the only problem can occur if
$\varphi_{|\left(\left(\mathbb{P}^{1}\right)^{5}\right)^{\op{stable}}}$
has a positive dimensional fiber. This happens only over curves $\Phi$
which are the union of a $(0,1)$-ruling of $C\times \mathbb{P}^{1}$
sitting over one of the five marked points $p_{i} \in C$, and two
distinct $(1,0)$-rulings. Moreover by stability the $\varphi$-preimage
of such a $\Phi$ in
$\left(\left(\mathbb{P}^{1}\right)^{5}\right)^{\op{stable}}$ consists
of five-tuples $\bF$ in which one point lies on the $(0,1)$-ruling and each
of the $(1,0)$-rulings contains exactly two points from $\bF$. The
positive dimensional component of the $\varphi$-fiber then traced by
moving the point on the $(0,1)$-ruling arbitrarily. But all such
choices lie in a single $\mathbb{P}SL_{2}(\mathbb{C})$-orbit since
$\mathbb{P}SL_{2}(\mathbb{C})$ will move transitively the triple of
points in the fiber. This shows that $\beta_{0}$ is quasi-finite on 
$N_{0}(0,q)_{k=0}$ and completes the proof of the lemma. \
\hfill $\Box$

\

\noindent
Next we describe the  structure of the finite map $\beta_{0}$:

\begin{lem} \label{lem:extend.dP5}
Suppose that $2/5 < q \leq 2/3$, then 
\[
\beta_{0} : N_{0}(0,q) \to S^{2}C
\] 
is a finite branched Galois cover with Galois group 
\[
\op{Gal}(N_{0}(0,q)/S^{2}C) \cong
\left(\mathbb{Z}/2\right)^{4}. 
\]
\end{lem}
{\bfseries Proof.}  From our analysis of the fibers of the rational map
$\varphi : \left(\mathbb{P}^{1}\right)^{5} \to \mathbb{P}^{5}$ it
follows that over a dense open set $U \subset \mathbb{P}^{5}$ the map
$\varphi : \varphi^{-1}(U) \to U$ is a torsor over a bundle of groups
$\mathcal{G} \to U$ which is a twisted form of
$\left(\mathbb{Z}/2\right)^{5}$. Indeed let $U \subset
|\mathcal{O}(1,2)|$  be the open subset of all smooth 
$\Phi \in |\mathcal{O}(1,2)|$ that are not
tangent to any of the rulings $\{ p_{i} \}\times \mathbb{P}^{1}$. As
we explained before, for any $\Phi \in U$ 
the fiber $\varphi^{-1}(\Phi)$ consists of the $2^{5}$ five tuples of
points corresponding to choosing one of the two preimages of each
$p_{i}$ in $\Phi$. Consider the group $\mathcal{G}_{\Phi}$ generated by the
flips of the two preimages of $p_{i}$ in $\Phi$ acts naturally and
simply transitively on the fiber $\varphi^{-1}(\Phi)$. This group is
simply the product
\[
\mathcal{G}_{\Phi} = \prod_{i=1}^{5}S_{\Phi_{p_{i}}}
\] 
of the symmetric groups
$S_{\Phi_{p_{i}}}$ of the fibers $\Phi_{p_{i}}$ of $\Phi \to C$ over
$p_{i} \in C$, and is thus isomorphic to
$\left(\mathbb{Z}/2\right)^{5}$. The bundle of groups $\mathcal{G} \to
U$ is then the natural bundle with fibers $\mathcal{G}_{\Phi}$, $\Phi
\in U$. 

In fact the bundle of groups $\mycal{G}$ is trivial and can be
described explicitly in terms of global automorphisms of the
Hirzebruch surface $C\times \mathbb{P}^{1}$. The first observation
that makes this description possible is the following
 
\begin{claim} \label{claim:global.a.Phi} Let $\Phi$ be a smooth
  $(1,2)$ curve in the surface $C\times \mathbb{P}^{1}$, then there
  exists a unique element $\alpha_{\Phi} \in
  \mathbb{P}SL_{2}(\mathbb{C})$ of order two, so that the automorphism
  $\op{id}_{C}\times \alpha_{\Phi} : C\times \mathbb{P}^{1} \to
  C\times \mathbb{P}^{1}$ preserves $\Phi$ and induces on $\Phi$ the
  covering involution of the  double cover $\Phi \to C$.
\end{claim}
{\bfseries Proof.} If we choose an inclusion of $C\times \mathbb{P}^{1}$ as
a quadric in $\mathbb{P}^{3}$, the $(1,2)$ curves in $C\times
\mathbb{P}^{1}$ become twisted cubics in $\mathbb{P}^{3}$. Note that the
linear automorphisms of $\mathbb{P}^{3}$ act transitively on the set
of pairs consisting of a smooth quadric surface in $\mathbb{P}^{3}$
and a smooth twisted cubic contained in that quadric surface. In
particular, for any smooth $(1,2)$ curve $\Phi \subset C\times
\mathbb{P}^{1}$ we can choose an embedding of $C\times \mathbb{P}^{1}$
in $\mathbb{P}^{3}$ with coordinates $(z_{0}:z_{1},z_{2}:z_{3})$ so that 
\begin{itemize}
\item $\Phi$ is the image of the map $\mathbb{P}^{1} \to
  \mathbb{P}^{3}$, $(t_{0}:t_{1}) \mapsto
  (s_{0}^{3}:s_{0}^{2}s_{1}:s_{0}s_{1}^{2}:s_{1}^{3})$. 
\item $C\times \mathbb{P}^{1}$ is the quadric with equation
  $z_{0}z_{3} = z_{1}z_{2}$.
\item The rulings $\{ p \}\times \mathbb{P}^{1}$ map to the lines 
$z_{1}\mu_{1} = z_{3}\mu_{0}$, $z_{0}\mu_{1} = z_{2}\mu_{0}$.
\end{itemize}
The covering involution for the cover $\Phi \to C$ now is given by
$(s_{0}:s_{1}) \mapsto (s_{0}:-s_{1})$, and it obviously extends to
the the global automorphism $A_{\Phi}(z_{0}:z_{1}:z_{2}:z_{3}) :=
(z_{0}:-z_{1}:z_{2}:-z_{3})$ of $\mathbb{P}^{3}$. This automorphism
preserves the quadric and preserves each ruling $\{ p \}\times
\mathbb{P}^{1}$. The assignment $p \mapsto A_{\Phi|\{ p \}\times
\mathbb{P}^{1}}$ is an regular morphism from $C$ to the affine variety
$\mathbb{P}SL_{2}$ and is hence constant. Therefore if we take  $\alpha_{\Phi}
:= A_{\Phi|\{ p \}\times \mathbb{P}^{1}}$ for some $p$, then
the restriction of $A_{\Phi}$ to the image of $C\times
\mathbb{P}^{1}$ will be the automorphism $\op{id}_{C}\times
\alpha_{\Phi}$. This proves the claim. \ \hfill $\Box$

\

\noindent
Let  $\alpha : U \to
\mathbb{P}SL_{2}(\mathbb{C})$ denote the morphism which to each $\Phi
\in U$ assigns the involution $\alpha_{\Phi}$.  The map $\alpha$ and
Claim~\ref{claim:global.a.Phi} immediately yield
a global isomorphism of bundles of groups :
\[
\xymatrix@R-1pc{
U\times \left( \mathbb{Z}/2 \right)^{5} \ar[rr]^-{\bT}_-{\cong} \ar[rd] & &
\mathcal{G} \ar[ld]  & \hspace{-1.5pc} \subset U\times
\mathbb{P}SL_{2}(\mathbb{C})^{5}\\
& U
}
\]
given by the formula
\[
\bT\left(\Phi,\varepsilon_{1},\varepsilon_{2},\varepsilon_{3},\varepsilon_{4},
\varepsilon_{5}\right) 
:=
\left(\Phi,\alpha_{\Phi}^{\varepsilon_{1}},\alpha_{\Phi}^{\varepsilon_{2}},
\alpha_{\Phi}^{\varepsilon_{3}},\alpha_{\Phi}^{\varepsilon_{4}},
\alpha_{\Phi}^{\varepsilon_{5}}\right).
\]
In other words we get an action of $\left( \mathbb{Z}/2 \right)^{5}$
on the the variety $\varphi^{-1}(U)$ given by 
\begin{equation} \label{eq:Z2^5}
(\varepsilon_{1},\varepsilon_{2},\varepsilon_{3},\varepsilon_{4},
\varepsilon_{5})\cdot (f_{1},f_{2},f_{3},f_{4},f_{5}) := 
\bT\left(\varphi(f_{1},f_{2},f_{3},f_{4},f_{5}),\varepsilon_{1},
\varepsilon_{2},\varepsilon_{3},\varepsilon_{4}, 
\varepsilon_{5}\right).
\end{equation}
From the formula \eqref{eq:Z2^5} it is clear that $\left( \mathbb{Z}/2
\right)^{5}$  preserves and acts simply
transitively on the fibers of $\varphi$. In other words
$\varphi^{-1}(U) \to U$ is an etale $\left( \mathbb{Z}/2 \right)^{5}$-Galois
cover.

Next observe that since $\varphi^{-1}(U) \subset \left(
\mathbb{P}^{1}\right)^{5}$ is open and dense, it follows the $\left(
\mathbb{Z}/2 \right)^{5}$-action \eqref{eq:Z2^5} extends to a
birational action of $\left( \mathbb{Z}/2 \right)^{5}$ on $\left(
\mathbb{P}^{1}\right)^{5}$. Note also that the proof of
Claim~\ref{claim:global.a.Phi} shows that the map $\alpha$ is
$\mathbb{P}SL_{2}(\mathbb{C})$-equivariant, i.e. for every $g \in
\mathbb{P}SL_{2}(\mathbb{C})$ and every $\Phi \in U$ we have
\[
\alpha_{g(\Phi)} = g\alpha_{\phi}g^{-1}.
\]
So, the action \eqref{eq:Z2^5} induces a birational action of $\left(
\mathbb{Z}/2 \right)^{5}$ on $N_{0}(0,q)$ which preserves 
the fibers of $\beta_{0}$ and acts transitively on the general fiber. 

Let now $g \in \left( \mathbb{Z}/2 \right)^{5}$. The graph of the
induced birational automorphism of $N_{0}(0,q)$ will be
contained in the fiber product $N_{0}(0,q)\times_{ S^{2}C}
N_{0}(0,q) \subset N_{0}(0,q)\times
N_{0}(0,q)$ and so the closure of this graph will also be
contained in $N_{0}(0,q)\times_{ S^{2}C}
N_{0}(0,q)$. But by Lemma~\ref{lem:beta0.morphism}
$N_{0}(0,q)\times_{ S^{2}C} N_{0}(0,q)$ is finite over
$N_{0}(0,q)$ and so the closure of the graph of $g$ in
$N_{0}(0,q)\times N_{0}(0,q)$ maps birationally and
finitely onto $N_{0}(0,q)$. Since $N_{0}(0,q)$ is a
smooth surface this implies that the closure of $g$ is isomorphic to
$N_{0}(0,q)$, i.e. $g$ extends to a biregular utomorphism of
$N_{0}(0,q)$.

\

\noindent
Since $(\mathbb{Z}/2)^{5}$ acts simply transitively on the general
fiber of $\beta_{0}$ it follows that the image of $(\mathbb{Z}/2)^{5}$
in the automorphism group of the surface $N_{0}(0,q)$ will act
simply transitively on the general fiber of $\beta_{0}$. But
$N_{0}(0,q)_{k=0}$ is the
$\mathbb{P}SL_{2}(\mathbb{C})$-quotient of the stable locus in
$\left(\mathbb{P}^{1}\right)^{5}$ and so the kernel of the map from
$(\mathbb{Z}/2)^{5}$ to the automorphism group of $N_{0}(0,q)$
consists of all elements
$\left(\varepsilon_{1}, \varepsilon_{2}, \varepsilon_{3},
\varepsilon_{4}, \varepsilon_{5}\right)$ such that 
for every $\boldsymbol{f} = 
(f_{1},f_{2},f_{3},f_{4},f_{5}) \in \varphi^{-1}(U)$ we
can find a $g^{\boldsymbol{f}} \in \mathbb{P}SL_{2}(\mathbb{C})$ so that 
\[
\left(\varepsilon_{1}, \varepsilon_{2}, \varepsilon_{3},
\varepsilon_{4}, \varepsilon_{5}\right)\cdot 
(f_{1},f_{2},f_{3},f_{4},f_{5})  =
\left(g^{\boldsymbol{f}}(f_{1}), g^{\boldsymbol{f}}(f_{2}),
g^{\boldsymbol{f}}(f_{3}), g^{\boldsymbol{f}}(f_{4}),
g^{\boldsymbol{f}}(f_{5})
\right).
\]
Substituting the formula \eqref{eq:Z2^5} for the action in this
condition yields 
\[
g^{\boldsymbol{f}} =
\alpha_{\varphi\left(\boldsymbol{f}\right)}^{\varepsilon_{i}}, 
\quad i = 1, \ldots, 5,
\]
which is possible if and only if $\varepsilon_{1} = \varepsilon_{2} =
\varepsilon_{3} = \varepsilon_{4} = \varepsilon_{5}$.  Therefore the
kernel of the map from $(\mathbb{Z}/2)^{5}$ to the automorphism group
of $N_{0}(0,q)$ is $\mathbb{Z}/2$ which sits diagonally in
$(\mathbb{Z}/2)^{5}$. This proves the lemma. \ \hfill $\Box$

\

\noindent
We are now ready to prove part {\bfseries (ii)} of
Proposition \ref{prop:synthetic.N0}. The proof of Lemma~\ref{lem:extend.dP5} shows that the
map $\beta_{0}$ is unramified over the open set $\mathfrak{b}_{0}(U)
\subset S^{2}C$. But the complement of this open set parameterizes the
branch divisors of all $(1,2)$ curves in $C\times \mathbb{P}^{1}$
which are tangent to one of the rulings $\{ p_{i} \}\times
\mathbb{P}^{1}$. Thus $S^{2}C - \mathfrak{b}_{0}(U) = \cup_{i =1}^{5}
(p_{i} + C)$. But the curve $p_{i} + C \subset S^{2}C$ is precisely
the line in $\mathbb{P}^{2} = S^{2}C$ which  is tangent to the
diagonal conic $C \subset S^{2}C$ at the point $p_{i}$. Thus $p_{i} +
C = t_{i}$ and $\beta_{0} : N_{0}(0,q) \to S^{2}C$ is
unramified over the complement of the five lines $t_{i}$, $i = 1,
\ldots,5$. Furthermore a general point of
each $t_{i}$ comes from a {\em smooth} $(1,2)$  curve $\Phi$ which is
tangent to $\{ p_{i} \}\times
\mathbb{P}^{1}$. Now the explicit definition of the Galois
action on $\left(\mathbb{P}^{1}\right)^{5}$ immediately shows that
$\beta_{0}$ has a simple Galois ramification over such a point. In
particular the quotient of $N_{0}(0,q)$ by an index two
subgroup in the Galois group will be a double cover of $S^{2}C$
branched at two of the lines $t_{i}$ it is clear that all the pairs of
lines arise in this way and so from the universal property of the root cover
branched at the pair of lines we conclude that for every $i \neq j \in
\{1,2,3,4,5\}$ there is a finite morphism $N_{0}(0,q) \to
Y_{ij}$ which is more over equivariant for some surjective
homomorphism $\left(\mathbb{Z}/2\right)^{5} \to \mathbb{Z}/2$ of
Galois groups. This
shows that we have a $(\mathbb{Z}/2)^{4}$-equivariant birational
morphism 
\[
N_{0}(0,q) \to Y_{12}\times_{S^{2}C} Y_{23}\times_{S^{2}C}
Y_{34}\times_{S^{2}C} Y_{45}.
\]
But $N_{0}(0,q)$ is smooth and so this birational morphism
must be the normalization morphism. 

Finally, the interpretation of $N_{0}(0,q)$ as the
normalization of the fiber product of root covers shows that the
branch divisor of $\beta_{0}$ is the divisor $\mathfrak{B} = 
\sum_{i = 1}^{5} t_{i}$
and so by the simple Galois ramification property we have that the
ramification divisor $\mathfrak{R}$ of $\beta_{0}$ satisfies 
\[
2\mathfrak{R} = \beta_{0}^{*} \mathfrak{B} = \beta_{0}^{*}
\left(\sum_{i = 1}^{5} t_{i}\right).
\] 
Now by the Hurwitz formula we get 
\[
K_{N_{0}(0,q)} =
\mathcal{O}_{N_{0}(0,q)}(\mathfrak{R})\otimes \beta_{0}^{*}
K_{S^{2}C} = \mathcal{O}_{N_{0}(0,q)}(\mathfrak{R})\otimes 
\beta_{0}^{*} \mathcal{O}_{\mathbb{P}^{2}}(-3),
\]
and so
\[
K_{N_{0}(0,q)}^{\otimes 2} =
\mathcal{O}_{N_{0}(0,q)}(2\mathfrak{R}) \otimes 
\beta_{0}^{*} \mathcal{O}_{\mathbb{P}^{2}}(-3) = \beta_{0}^{*}
\mathcal{O}_{\mathbb{P}^{2}}(5-6) = \beta_{0}^{*}
\mathcal{O}_{\mathbb{P}^{2}}(-1).
\]
This completes the proof of part {\bfseries (ii)}.

\

\noindent
The statement {\bfseries (iii)} follows immediately from the statement
{\bfseries (i)}. Indeed consider the blow-up of $C\times C$ in the
four points $(p_{1},p_{1})$, $(p_{2},p_{2})$, $(p_{3},p_{3})$, and
$(p_{4},p_{5})$ contains five natural disjoint $(-1)$-curves - the three
exceptional divisors corresponding to $(p_{1},p_{1})$,
$(p_{2},p_{2})$, $(p_{3},p_{3})$ and the strict transforms of the two
rulings of $C\times C$ passing through $(p_{4},p_{5})$. When we
contract these five curves we get a birational morphism from
$N_{0}(0,q)$ which maps the diagonal copy of $C$ in $C\times
C$ to a conic, and the five exceptional divisors to the points $p_{i}
\in C$ on that conic. Thus this $\mathbb{P}^{2}$ is canonically
identified with $S^{2}C$ and $N_{0}(0,q)$ is identified with
the blow-up of $S^{2}C$ at the five points $p_{i}$ sitting on the
conic $C \subset S^{2}C$. 

This completes the proof of Proposition \ref{prop:synthetic.N0}.
\ \hfill $\Box$

\

\bigskip

\begin{rem} \label{rem:pencil}
{\bfseries (a)} Note that the isomorphisms in parts {\bfseries (i)}
and {\bfseries (iii)} of  Proposition \ref{prop:synthetic.N0} depend on the
labelling of the points $p_{i}$, or more precisely depend on
partitioning the five points into $3$ and $2$. 

\

\noindent
{\bfseries (b)} After the fact, the $16$-sheeted cover picture from
  part {\bfseries (ii)} of the proposition can be described in very
  simple terms via the anti-canonical emebedding of
  $N_{0}(0,q) = dP_{5}$. Indeed the anticanonical embedding
  realizes a general $dP_{5}$ as a complete intersection of two quadrics in
  $\mathbb{P}^{4}$. The corresponding pencil of quadrics has five
  singular members and the vertices of these quadrics are five points
  in general position in $\mathbb{P}^{4}$. These five points give rise
  to a system of projective coordinates on $\mathbb{P}^{4}$: the
  unique projective coordinate system
  $(x_{1}:x_{2}:x_{3}:x_{4}:x_{5})$ for which the five points are the
  vertices of the standard tetrahedron. In these coordinates we can
  describe the $dP_{5}$ as the intersection of the quadrics
\[
\sum_{i=1}^{5} x_{i}^{2}  = 0, \qquad  \sum_{i=1}^{5} \lambda_{i}
 x_{i}^{2}  = 0,
\]
where the $\lambda_{i}$'s denote the positions of the singular fibers
of the pencil of quadrics.

Now note that we have a natural morphism of degree $16$
\[
\mathbb{P}^{4} \to \mathbb{P}^{4}, \quad
(x_{1}:x_{2}:x_{3}:x_{4}:x_{5}) \mapsto
(x_{1}^{2}:x_{2}^{2}:x_{3}^{2}:x_{4}^{2},x_{5}^{2}),
\]
which maps the $dP_{5}$ to the plane with equations $\sum_{i=1}^{5}
y_{i} = 0$, $\sum_{i=1}^{5}
\lambda_{i}y_{i} = 0$ in the target $\mathbb{P}^{4}$. The restriction
of this map to $dP_{5}$ is precisely the morphism $\beta_{0}$ as can
be easily seen from the branching behavior.
\end{rem}

\

\noindent
Next we should describe the moduli space $N_{0}(0,q)$ for $q$ in one
of the last two chambers. First we deal with the fourth chamber

\begin{prop} \label{prop:P2}
Suppose $4/5 \leq  q \leq 1$. Then
  $N_{0}(0,q)$ is a point. 
\end{prop}
{\bfseries Proof.} In this case we have $V = \mathcal{O}_{C}(-2)\oplus
\mathcal{O}_{C}(2)$  and $N_{0}(0,q)$ is the quotient
\[
N_{0}(0,q) = \left.\left( \mathbb{P}(V)_{p_{1}}\times
\mathbb{P}(V)_{p_{2}}\times  \mathbb{P}(V)_{p_{3}}\times
\mathbb{P}(V)_{p_{4}}\times
\mathbb{P}(V)_{p_{5}}\right)^{\op{stable}}\right/\op{Aut}(\mathbb{P}(V)).
\]
The solvable group $\op{Aut}(\mathbb{P}(V)) =
\mathbb{C}^{\times} \ltimes H^{0}(C,\mathcal{O}(4))$ acts 
with closed orbits on the stable locus. Indeed, from the last row of
Table~\ref{table:deg.0.stable} we see that a five tuple
$(f_{1},f_{2},f_{3},f_{4},f_{5})$ will belong to $\left(
\mathbb{P}(V)_{p_{1}}\times \mathbb{P}(V)_{p_{2}}\times
\mathbb{P}(V_{p_{3}})\times \mathbb{P}(V)_{p_{4}}\times
\mathbb{P}(V)_{p_{5}}\right)^{\op{stable}}$ if and only if none of the
points $f_{i}$ lie on the unique $(-4)$-curve in $\mathbb{P}(V) \cong
\mathbb{F}_{4}$. In other words we have
\[
\left( \mathbb{P}(V)_{p_{1}}\times \mathbb{P}(V)_{p_{2}}\times
\mathbb{P}(V_{p_{3}})\times \mathbb{P}(V)_{p_{4}}\times
\mathbb{P}(V)_{p_{5}}\right)^{\op{stable}} = \prod_{i = 1}^{5}
\mathcal{O}_{C}(4)_{p_{i}} \cong \mathbb{C}^{5},
\]
and an element $(z,g) \in \mathbb{C}^{\times} \ltimes
H^{0}(C,\mathcal{O}(4))$ acts on this $\mathbb{C}^{5}$ by
\[
(f_{1},f_{2},f_{3},f_{4},f_{5}) \mapsto (z\cdot f_{1} + g(p_{1}),z\cdot f_{2}
+ g(p_{2}), z\cdot f_{3} + g(p_{3}), z\cdot f_{4} + g(p_{4}), z\cdot
f_{5} + g(p_{5})).
\]
But by Lagrange interpolation any $(f_{1},f_{2},f_{3},f_{4},f_{5}) \in
\prod_{i = 1}^{5} \mathcal{O}_{C}(4)_{p_{i}}$ is equal to \linebreak
$(g(p_{1}),g(p_{2}),g(p_{3}),g(p_{4}),g(p_{5}))$ for a unique $g$ and
so in this case the stack quotient \linebreak
$[\mathbb{C}^{5}/\op{Aut}(\mathbb{P}(V))]$ is isomorphic to
$B\mathbb{C}^{\times}$ and the coarse moduli space $N_{0}(0,q)$ is a
\linebreak point. \ \hfill $\Box$

\

\medskip

Finally we have to describe $N_{0}(0,q)$ for $2/3 < q < 4/5$. However,
it turns out that the analysis of this moduli space is easier if we
recast it as a moduli space of parabolic bundles of underlying degree
one. By Lemma~\ref{lem:moduli.isos} we have a natural isomorphism of
$N_{0}(0,q)$ and $N_{1}(0,1-q)$. Thus it suffices to describe
$N_{1}(0,q)$ for $1/5 \leq q < 1/3$.

\

\subsubsection{The moduli spaces $N_{1}(0,q)$} \label{sssec-N1}

We have already dealt with all chambers but one. So here we fix a real
number $q$ with $1/5 \leq q < 1/3$ and suppose we have a parabolic
rank two vector bundles $\mathcal{V} = (V,\bF,(\ba,\bb))$ equipped
with an isomorphism $\det V \cong \mathcal{O}_{C}(1)$ and with fixed
balanced parabolic weights: $\ba = \{ a_{i} \}_{i =1}^{5}$, $\bb = \{
b_{i} \}_{i =1}^{5}$, with $a_{i} = 0$ and $b_{i} = q$ for all $i$.
In this case we will have $V \cong \mathcal{O}_{C}(-k)\oplus
\mathcal{O}_{C}(1+k)$ for some $k \geq 0$ and the stability criterion
of Lemma~\ref{lem:stability} implies that $\mathcal{V} =
(V,\bF,(\ba,\bb))$ is semistable if and only if for all $d \geq -1-k$
we have
\begin{equation} \label{eq:deg1.stability}
\kk_{d}(\mathcal{V}) \leq \frac{5}{2} + \frac{2d+1}{2q},
\end{equation}
or equivalently
\[
- d \leq \left( \frac{5}{2} - \kk_{d}(\mathcal{V})\right)q + \frac{1}{2}.
\]
Since
$0 \leq \kk_{d}(\mathcal{V})$ whenever the linear system
$|\mathcal{O}_{V}(1)\otimes \gamma^{*}\mathcal{O}_{C}(d)|$
contains an irreducible curve, and 
 $1/5 \leq q < 1/3$ by assumption, we get 
\[
- d \leq \left( \frac{5}{2} - \kk_{d}(\mathcal{V})\right)q + \frac{1}{2} <
 \left( \frac{5}{2} - 0)\right)\frac{1}{3} + \frac{1}{2}  = \frac{4}{3}.
\]
Taking $- d = 1+k$ we conclude that in this case we must have $k = 0$, 
i.e. $V \cong
\mathcal{O}_{C}\oplus \mathcal{O}_{C}(1)$. 

It will be convenient to identify the
surface $\mathbb{P}(V) \cong \mathbb{F}_{1}$ with the blow-up of
$\mathbb{P}^{2}$ at a point $\pt \in \mathbb{P}^{2}$. Under this
identification the 
five fibers of $\gamma : \mathbb{P}(V) \to C$ over the points $p_{i}$
become five distinct lines $\ell_{p_{1}}, \ell_{p_{2}},
\ell_{p_{3}},\ell_{p_{4}},\ell_{p_{5}} \subset \mathbb{P}^{2}$  all
passing through $\pt$. The quasi-parabolic structures on $V =
\mathcal{O}_{C}\oplus \mathcal{O}_{C}(1)$ then correspond to
configurations of five tuples of points
$\bF = (F_{1},F_{2},F_{3},F_{4},F_{5})$ in $\mathbb{P}^{2}$, so that $F_{i}
\in  \ell_{p_{i}}$. It is straightforward to rewrite the
semi-stability conditions on $\mathcal{V} = (V,\bF,(0,q))$  in terms
of this geometry. 

\index{terms}{quasi-parabolic!structure}

\begin{table}[!ht]
\begin{center}
{\small
\begin{tabular}{|c|c|c|} 
\hline
$V$ & $\bF$ & $q$ \\ 
\hline
\hline 
$\mathcal{O}_{C}\oplus \mathcal{O}_{C}(1)$ & 
\begin{minipage}[c]{3.8in} \addtolength{\baselineskip}{-5pt}
\

\begin{itemize}
\item $F_{i} \neq \pt$ for all $i = 1, \ldots, 5$.
\item At most $4$  points in $\bF$ are collinear.
\end{itemize}

\
\end{minipage} &
${\displaystyle \frac{1}{5} \leq q < \frac{1}{3}}$ \\
\hline\hline
\end{tabular}
}
\end{center}
\caption{Stability configurations for $\mathcal{O}_{C}\oplus
  \mathcal{O}_{C}(1)$ }
\label{table:deg.1.stable}
\end{table}
\

\vspace{1cm}
\noindent
The moduli spaces corresponding to the various chambers
can be described explicitly:

\begin{prop} \label{prop:deg1.moduli} Let $0 \leq q < 1$. Then:
\begin{description}
\item[{\bfseries (i)}] If $0 \leq q < 1/5$,
  then the moduli space $N_{1}(0,q)$ is empty.
\item[{\bfseries (ii)}]  If  $1/5 \leq q < 1/3$, 
  then the moduli space $N_{1}(0,q)$ is isomorphic to
  $\mathbb{P}^{2}$.
\item[{\bfseries (iii)}]  If  $1/3 \leq q < 3/5$, 
  then the moduli space $N_{1}(0,q)$ is isomorphic to a del
  Pezzo surface $dP_{5}$. 
\item[{\bfseries (iv)}]  If  $3/5 \leq q  < 1$, 
  then the moduli space $N_{1}(0,q)$ is isomorphic to a del
  Pezzo surface $dP_{4}$.
\end{description}
\end{prop}
{\bfseries Proof.} Part {\bfseries (i)} follows from $N_{1}(0,q) \cong
H_{0}(0,1-q)$ and Proposition~\ref{prop:P2}.

\

\noindent
For part{\bfseries (ii)} recall that $V = \mathcal{O}_{C}\oplus
\mathcal{O}_{C}(1)$ and that we have identified $\mathbb{P}(V)$ with
the blow-up of a plane $\mathbb{P}^{2}$ at a point $\pt$. Under this
identification the five fibers of $\gamma : \mathbb{P}(V) \to C$ over
the points $p_{i}$ become five distinct lines $\ell_{p_{1}},
\ell_{p_{2}}, \ell_{p_{3}},\ell_{p_{4}},\ell_{p_{5}} \subset
\mathbb{P}^{2}$ all passing through $\pt$. The moduli space 
$N_{1}(0,q)$ is constructed as a quotient
\[
\left.\left( \ell_{p_{1}}\times \ell_{p_{2}}\times \ell_{p_{3}}\times
\ell_{p_{4}}\times 
\ell_{p_{5}} \right)^{\text{stable}}\right/ \op{Aut}(\mathbb{P}(V)),
\]
where:
\begin{itemize}
\item $\left( \ell_{p_{1}}\times \ell_{p_{2}}\times \ell_{p_{3}}\times
\ell_{p_{4}}\times 
\ell_{p_{5}} \right)^{\text{stable}} \subset \ell_{p_{1}}\times
\ell_{p_{2}}\times \ell_{p_{3}}\times 
\ell_{p_{4}}\times 
\ell_{p_{5}} $ denotes the $(0,q)$-stable five-tuples of points with the
stability condition corresponding to  $1/5 <
  q < 1/3$. 
\item $\op{Aut}(\mathbb{P}(V)) = \mathbb{C}^{\times}\ltimes
  H^{0}(C,\mathcal{O}_{C}(1)) \cong  \mathbb{C}^{\times}\ltimes
  \mathbb{C}^{2}$. 
\end{itemize}
The first stability constraint in the second row of
Table~\ref{table:deg.1.stable} shows that 
\[
\left( \ell_{p_{1}}\times \ell_{p_{2}}\times \ell_{p_{3}}\times
\ell_{p_{4}}\times 
\ell_{p_{5}} \right)^{\text{stable}} \subset \prod_{i=1}^{5}
\left(\ell_{p_{i}} - \{\pt\} \right) \cong \mathbb{C}^{5}.
\]
Furthermore $H^{0}(C,\mathcal{O}_{C}(1))$ maps to $\mathbb{C}^{5}$ by
evaluation 
\[
g \in H^{0}(C,\mathcal{O}_{C}(1)) \mapsto
(g(p_{1}),g(p_{2}),g(p_{3}),g(p_{4}),g(p_{5}))
\] 
and acts on
$\prod_{i=1}^{5} \left(\ell_{p_{i}} - \{\pt\} \right) \cong
\mathbb{C}^{5}$ by translation. Thus a point in $\mathbb{C}^{5}$
satisfies the second condition from the second row of
Table~\ref{table:deg.1.stable} if and only if modulo the translation
action of $\mathbb{C}^{2} \cong H^{0}(C,\mathcal{O}_{C}(1))$ not all
five coordinates of this point are equal to zero. In other words, the
set of stable points $\left( \ell_{p_{1}}\times \ell_{p_{2}}\times
\ell_{p_{3}}\times \ell_{p_{4}}\times \ell_{p_{5}}
\right)^{\text{stable}}$ consists of all points in $\mathbb{C}^{5}$
which are {\em not}  contained in the orbit of $(0,0,0,0,0) \in
\mathbb{C}^{5}$ under the $\mathbb{C}^{2}$-action. This shows that 
\[
\begin{aligned}
\left.\left( \ell_{p_{1}}\times \ell_{p_{2}}\times \ell_{p_{3}}\times
\ell_{p_{4}}\times 
\ell_{p_{5}} \right)^{\text{stable}}\right/ \op{Aut}(\mathbb{P}(V)) &
= \left.\left( \mathbb{C}^{5} - \mathbb{C}^{2} \right)\right/
\left(\mathbb{C}^{\times} \ltimes \mathbb{C}^{2}\right) \\
& = \left.\left( \mathbb{C}^{3} -
\{0\}\right)\right/\mathbb{C}^{\times} \\
& \cong \mathbb{P}^{2}.
\end{aligned}
\]
\

\noindent
Parts {\bfseries (iii)} and {\bfseries (iv)} follow immediately from
Lemma~\ref{lem:moduli.isos}, Proposition~\ref{prop:dP4} and
Proposition~\ref{prop:dP5}. The proposition is proven.  \ \hfill
$\Box$

\

\noindent
Together Proposition~\ref{prop:dP4}, Proposition~\ref{prop:dP5},
Proposition~\ref{prop:P2}, Proposition~\ref{prop:deg1.moduli}, and the
synthetic picture from Proposition~\ref{prop:synthetic.N0} complete
the proof of Theorem~\ref{theo:shapes.of.moduli}.

\

\noindent
The description of the moduli spaces given in 
Theorem~\ref{theo:shapes.of.moduli} provides a concrete geometric
setup for trying to understand the Hecke correspondences and the
eigensheaf property in the realm of parabolic bundles.

\section{The Hecke correspondence} \label{sec:hecke}

Hecke correspondences are correspondences between moduli problems of
bundles or parabolic bundles. They parametrize elementary
modifications of bundles at points of the base curve $C$ and so have
modular interpretation themselves. Elementary modifications do not
preserve stability in general and so Hecke correspondences are most
naturally described in the context of moduli for all, not necessarily
stable, bundles.  Alternatively, one can describe the open subset of
the Hecke correspondence involving stable bundles on both sides, and
then construct a compactified and resolved version of this open Hecke
correspondence. This is the approach we pursue here.

After a brief review of Hecke correspondences of moduli stacks and
moduli spaces in general, we focus on the case at hand where the
parabolic weights are in the dominant chamber, i.e. the weights are
chosen so that the moduli space dominates by a morphism the moduli
spaces corresponding to choices in any other chamber.
\index{terms}{dominant chamber!of parabolic weights}
Specifically we
specialize to $q = 1/2$ which is in the dominant chamber. The
parabolic moduli spaces $N_d(0,\frac{1}{2})$, for $d=0$ or $1$, are
both identified with the $dP_5$ surface $X$. We describe the
compactified and resolved Hecke correspondence explicitly in this
case. We will see that it is obtained from the product $X \times X$ by
blowing up the diagonal $X$, and then blowing up the proper transforms
of the 16 surfaces $L \times L$, where $L$ runs through the 16 lines
in the del Pezzo surface $X$.

The Hecke correspondences have natural extensions, or rather
abelianizations, that act on Higgs bundles.  The above issue regarding
stability pertains to these abelianizations as well.  We will discuss
and analyze these abelianized Hecke correspondences in a later
chapter.

\subsection{Hecke correspondences of moduli stacks and 
moduli spaces}
\label{ssec-Hecke_general}

In this section we review the geometry of Hecke correspondences in a
general setting, both for moduli spaces and for stacks.

Let as before $(d,\ba,\bb)$ be a triple with $d \in \mathbb{Z}$ and
$\ba, \bb : \Par_{C} \to \mathbb{R}$. We will write $\sBun_{d}(\ba,\bb)$
for the moduli stack of all rank two parabolic bundles on $C$ with
determinant $\mathcal{O}_{C}(d)$ and parabolic weights $(\ba,\bb)$ and
we will write $\sN_{d}(\ba,\bb) \subset \sBun_{d}(\ba,\bb)$ for the
open substack of stable parabolic bundles. It is well known (see
e.g. \cite{faltings,laszlo-sorger,heinloth}) that $\sBun_{d}(\ba,\bb)$
is an Artin algebraic stack of infinite but locally finite type, and
that $\sN_{d}(\ba,\bb)$ is a substack of finite type. The moduli
problem $N^{\natural}_{d}(\ba,\bb)$ for stable\footnote{Here
  we assume
  that the weights are chosen generically so that every semi-stable
  parabolic bundle is actually stable.} parabolic bundles that we
studied in Section~\ref{sec:moduli.spaces} is the sheaf of sets
associated with $\sN_{d}(\ba,\bb)$ and so $N_{d}(a,b)$ is the coarse
moduli space of the stack $\sN_{d}(\ba,\bb)$. In particular we have a
natural morphism $\sN_{d}(\ba,\bb) \to N_{d}(\ba,\bb)$.

\index{notations}{Buns@$\sBun$}
  \index{notations}{Bunsdab@$\sBun_{d}(\ba,\bb)$}
  \index{notations}{Ndabs@$\sN_{d}(\ba,\bb)$}

Fix parabolic weights $(\ba,\bb)$.  The basic Hecke correspondence
associated with the fundamental representation of $GL_{2}(\mathbb{C})$
is (see e.g. \cite{heinloth}) the moduli stack $\sHeck(\ba,\bb)$ of triples 
\[
\left((V,\bF,(\ba,\bb)), (V',\bF',(\ba,\bb)),\beta \right),
\]
where $(V,\bF,(\ba,\bb))$ and $(V',\bF',(\ba,\bb))$ are parabolic
bundles with $\det V = \mathcal{O}_{C}$, $\det V' =
\mathcal{O}_{C}(1)$, and $\beta : V \to V'$ is an injective morphism
of parabolic locally free sheaves, such that $\op{length}\left( V'/V
\right) = 1$. The moduli stack $\sHeck(\ba,\bb)$ is a correspondence:
\[
\xymatrix@-1pc{
& \sHeck(\ba,\bb) \ar[dl]_-{p} \ar[dr]^-{q} & \\
\sBun_{0}(\ba,\bb) & & \sBun_{1}(\ba,\bb)\times C 
}
\]
where   \index{notations}{Buns@$\sBun$}
\[
\begin{aligned}
p\left((V,\bF,(\ba,\bb)), (V',\bF',(\ba,\bb)),\beta \right)  & :=
(V,\bF,(\ba,\bb)), \text{ and } \\
q\left((V,\bF,(\ba,\bb)), 
(V',\bF',(\ba,\bb)),\beta \right) & := \left( (V',\bF',(\ba,\bb)),
\op{supp}(V'/V) \right).
\end{aligned}
\]
If $x \in C$ is a point, then we can also consider the Hecke
correspondence $\sHeck_{x}(\ba,\bb) \subset \sHeck(\ba,\bb)$ 
at $x$ which is just the fiber at $x$ of the natural map $\sHeck(\ba,\bb) \to
C$. We will write 
\[
\xymatrix@-1pc{
& \sHeck_{x}(\ba,\bb) \ar[dl]_-{p_{x}} \ar[dr]^-{q_{x}} & \\
\sBun_{0}(\ba,\bb) & & \sBun_{1}(\ba,\bb)
}
\]
for the restrictions of the projections. The stacks
$\sBun_{0}(\ba,\bb)$, $\sBun_{1}(\ba,\bb)$, $\sHeck_{x}(\ba,\bb)$, and
$\sHeck(\ba,\bb)$ are smooth and the maps $p_{x}$, $q_{x}$ are
$\mathbb{P}^{1}$-bundles \cite{faltings,heinloth}.  The
correspondences $\sHeck(\ba,\bb)$ (respectively $\sHeck_{x}(\ba,\bb)$)
induce Hecke operators which send (twisted) $\mathcal{D}$-modules or Higgs
sheaves on $\sBun_{0}(\ba,\bb)$ to (twisted) $\mathcal{D}$-modules or Higgs
sheaves on $\sBun_{1}(\ba,\bb) \times C$ (respectively
$\sBun_{1}(\ba,\bb)$). These Hecke operators are defined as integral
transforms with kernels a (twisted) rank one flat bundle and a
(twisted) rank one Higgs sheaf on $\sHeck(\ba,\bb)$ (respectively
$\sHeck_{x}(\ba,\bb)$). In particular, the Hecke operators involve
pushforwards of an (untwisted) $\mathcal{D}$-module or a Higgs bundle via the
map $q$ (respectively $q_{x}$).

\

\noindent
The restrictions of the Hecke correspondences to the stable loci:
\[
\begin{aligned}
\sHeck^{\op{stable}}(\ba,\bb) & :=
\left[\sHeck(\ba,\bb)\right]
\bigtimes_{\sBun_{0}(\ba,\bb)\times \sBun_{1}(\ba,\bb)
  \times C}\displaylimits \left[ \sN_{0}(\ba,\bb)\times \sN_{1}(\ba,\bb)
  \times C\right] \\[1pc]
\sHeck_{x}^{\op{stable}}(\ba,\bb)& :=
\left[\sHeck_{x}(\ba,\bb)\right] \bigtimes_{\sBun_{0}(\ba,\bb)\times
  \sBun_{1}(\ba,\bb)}\displaylimits 
\left[ \sN_{0}(\ba,\bb)\times \sN_{1}(\ba,\bb)\right] 
\end{aligned}
\]
are no longer proper over either side and so before we can use them as
correspondences we need to compactify their fibers. Consider the
coarse moduli spaces   \index{notations}{Buns@$\sBun$}
\[
\begin{aligned}
\Heck^{\op{stable}}(\ba,\bb) & \subset N_{0}(\ba,\bb) \times
N_{1}(\ba,\bb) \times C,  \text{ and } \\
\Heck_{x}^{\op{stable}}(\ba,\bb)    & \subset N_{0}(\ba,\bb) \times
N_{1}(\ba,\bb)
\end{aligned}
\]
of the moduli stacks $\sHeck^{\op{stable}}(\ba,\bb)$ and
$\sHeck_{x}^{\op{stable}}(\ba,\bb)$. We can easily compactify these
spaces by taking their closure in $N_{0}(\ba,\bb) \times
N_{1}(\ba,\bb) \times C$ and $ N_{0}(\ba,\bb) \times N_{1}(\ba,\bb)$
respectively. We can further blow-up this closures to make them smooth
and to make sure that the complements of
$\Heck^{\op{stable}}(\ba,\bb)$ and $\Heck^{\op{stable}}_{x}(\ba,\bb)$
in the blow-up are divisors with normal crossings. These blow-up
closures will be our substitutes for the Hecke correspondences in the
parabolic context.

\

\subsection{The Hecke correspondence on $X$} \label{ssec-Hecke_X}

In this section we give explicit descriptions of the these compactified and resolved Hecke correspondences 
in the case that we are interested in.
Specifically we would like to understand the Hecke
correspondence for balanced parabolic weights $(a,b)$ such that both 
$N_{0}(a,b)$ and $N_{1}(a,b)$ are isomorphic to a $dP_{5}$. For 
instance we can take $(a,b) = (0,1/2)$. By Lemma~\ref{lem:moduli.isos}
we have a canonical isomorphism
\begin{equation} \label{eq:basic.iso}
\bt_{-3}\circ 
\bws_{-\frac{1}{2}}\circ \bps : N_{1}\left(0,\frac{1}{2}\right) 
\stackrel{\cong}{\longrightarrow}
N_{0}\left(0,\frac{1}{2}\right).
\end{equation}
We will use Proposition~\ref{prop:synthetic.N0} and the isomorphism
\eqref{eq:basic.iso} to identify $N_{0}\left(0,\frac{1}{2}\right)$ and
$N_{1}\left(0,\frac{1}{2}\right)$  with the del Pezzo surface $X$ obtained by
blowing up the projective plane $\mathbb{P}^{2} = S^{2}C$ at the five
points $p_{1}, p_{2}, p_{3}, p_{4}, p_{5}$ on the conic $C \subset
S^{2}C$.

\

\bigskip

\noindent
It turns out that our resolved and compactified Hecke correspondence
is closely related to the dual variety of $X$ in its anti-canonical
embedding. Before we give its precise description we will need some
notation

Consider the embedding $X \subset \mathbb{P}^{4}$ given by the linear
system $\left| - K_{X} \right|$. As in
Remark~\ref{rem:pencil}{\bfseries (b)} we can choose
coordinates $(x_{1}:x_{2}:x_{3}:x_{4}:x_{5})$ on $\mathbb{P}^{4}$ so
that $X$ is given as a complete intersection of two quadrics:
\[
X \; : \quad
\sum_{i=1}^{5} x_{i}^{2}  = 0, \qquad \sum_{i=1}^{5} \lambda_{i}
x_{i}^{2}  = 0.
\]
In other words $X$ is the base locus of the pencil of
quadrics $\{ Q_{\lambda} \}_{\lambda \in \mathbb{P}^{1}}$, where 
\[
Q_{\lambda} = \sum_{i = 1}^{5} (\lambda - \lambda_{i}) x_{i}^{2}.
\]
Let $\Lambda \cong  \mathbb{P}^{1}$ be the parameter line for this
pencil. Note that the pencil $\{ Q_{\lambda} \}_{\lambda \in
  \mathbb{P}^{1}}$ has exactly five singular fibers sitting over the points
$\lambda_{1}$, $\lambda_{2}$, $\lambda_{3}$, $\lambda_{4}$,
$\lambda_{5}$ in $\Lambda$. Furthermore, the set of five points
$\{\lambda_{1},\lambda_{2},\lambda_{3},\lambda_{4}, \lambda_{5}\}$ is
in a natural bijection with the set of five points
$\{p_{1},p_{2},p_{3},p_{4},p_{5}\}$ in $C$. Indeed, by
Proposition~\ref{prop:synthetic.N0}, we can also identify $X$ with the
16-sheeted Galois cover of $S^{2}C$  branched along the five lines
tangent to the conic $C \subset S^{2}C$ at the points $p_{i} \in
C$. On the other hand, as we saw in Remark~\ref{rem:pencil}{\bfseries
  (b)}, the covering map $\beta_{0} : X \to S^{2}C$ is the restriction
to $X$ of the morphism $\mathbb{P}^{4} \to \mathbb{P}^{4}$ that
squares all coordinates. In particular the set of branch lines of $\beta_{0}$
is in a bijection with the set of five coordinate functions $\{
x_{1}, x_{2},x_{3},x_{4},x_{5}\}$ which in turn is in bijection with
the set $\{\lambda_{1},\lambda_{2},\lambda_{3},\lambda_{4},
\lambda_{5}\}$. Equivalently we can say that the two sets
$\{\lambda_{i}\}$ and $\{p_{i}\}$ are in bijection since they can both
be identified with the set of vertices of the quadrics in the pencil
cutting out $X$. From now on we will assume that the coordinates in
$\mathbb{P}^{4}$ that we choose in Remark~\ref{rem:pencil}{\bfseries
  (b)} are ordered so that under the bijection matching
$\{\lambda_{1},\lambda_{2},\lambda_{3},\lambda_{4}, \lambda_{5}\}$ and
$\{p_{1},p_{2},p_{3},p_{4},p_{5}\}$ we have $\lambda_{i}
\leftrightarrow p_{i}$ for all $i = 1, \ldots, 5$. 

We now have the following 

\begin{prop} \label{prop:C=Lambda} The bijection $\lambda_{i}
\leftrightarrow p_{i}$ arising from the identification of the
sets
$\{\lambda_{i}\}$ and $\{p_{i}\}$ with the set of vertices of the
singular quadrics in the pencil $\{Q_{\lambda} \}_{\lambda \in \Lambda}$
extends to a global isomorphism $\iota : \Lambda \to C$ of
projective lines.
\end{prop}
{\bfseries Proof.} We will give two different geometric constructions of the
isomorphism $\iota$ each of which is interesting in its own right. 
Recall that in the process of proving
Proposition~\ref{prop:synthetic.N0} we argued that the preimage of the
conic $C \subset S^{2}C$ under the $16$-sheeted map $\beta_{0} : X \to
S^{2}C$ is completely reducible and becomes the union of the $16$
lines on the del Pezzo surface $X$. In particular the map $\beta_{0}$
identifies each line $L$ on $X$ with the conic $C$. Furthermore from
the basic geometry of a $dP_{5}$ we know that each line meets exactly
$5$ other lines and so the five marked points in $C$ correspond under
this identification to the five points on $L$ where $L$ meets these
five other lines. Now let us fix a line $L \subset X$. For instance we
can fix the line that gets contracted when we cross the chamber wall
$q = 2/5$. Let $\bP_{L} \subset \left(\mathbb{P}^{4}\right)^{\vee}$ be
the $\mathbb{P}^{2}$ parametrizing all hyperplanes in $\mathbb{P}^{4}$
that contain the line $L$.  There is a natural morphism
\begin{equation} \label{eq:double.cover}
\Lambda\times L \to \bP_{L}, \qquad (\lambda,x) \mapsto T_{x}Q_{\lambda}.
\end{equation}
From the definition it follows immediately that this morphism sends
curves of the form $\Lambda\times \{x\}$ and $\{\lambda\}\times L$ to
lines in $\bP_{L}$. Thus the morphism \eqref{eq:double.cover} is a
double cover branched along a conic. The ramification divisor for
\eqref{eq:double.cover} is a smooth rational curve of bidegree $(1,1)$
on $\Lambda\times L$ which can be viewed as the graph of an
isomorphism $\Lambda \to L$. 

Next suppose that $x \in L$ is the point of intersection of $L$ with
another line $L' \subset X$.  Then the $L$ and $L'$ together span the
tangent plane $T_{x}X$. Each quadric $Q_{\lambda}$ in our pencil
contains $X$ and hence contains the two lines $L$ and $L'$. Therefore
$Q_{\lambda}\cap T_{x}X \supset L\cup L'$, for all $\lambda$. This
shows that there must be a quadric in our pencil which contains the
plane $T_{x}X$. Indeed if we assume that all the $Q_{\lambda}$'s
intersect the plane $T_{x}X$ properly, then we will get a pencil of
conics in $T_{x}X$. Since each of these conics contains and is hence
equal to the reducible conic $L + L'$ we get a
contradiction. Furthermore a quadric of full rank in $\mathbb{P}^{4}$
can not contain a plane and so the quadrics in the pencil containing
$T_{x}X$ must be singular. But if we have two such quadrics, then they
will generate the pencil and so all $Q_{\lambda}$ will have to contain
$T_{x}X$ which is impossible since the general $Q_{\lambda}$ is
smooth. In other a point $x$ of intersection of $L$ with another line
$L' \subset X$ determines a unique singular quadric in the pencil or
equivalently determines one of the $\lambda_{i}$'s. Next note that
since $L \subset T_{x}X \subset Q_{\lambda_{i}}$ we will have that for
all $y \in L$ the tangent space $T_{y}Q_{\lambda_{i}}$ will contain
the plane $T_{x}X$. In other words the image of $\{\lambda_{i}
\}\times L$ under the map \eqref{eq:double.cover} is the line $\mathfrak{l}
\subset \bP_{L}$ parametrizing all hyperplanes through $T_{x}X$.
On the othe hand, for every $\lambda \in \Lambda$ we have
$T_{x}Q_{\lambda} \supset T_{x}X$ and so \eqref{eq:double.cover} sends
$\Lambda\times \{x \}$ also to $\mathfrak{l} \subset \bP_{L}$. This
shows that the preimage of the line  $\mathfrak{l} \subset \bP_{L}$
under the map \eqref{eq:double.cover} splits as $(\Lambda\times \{x
\})\cup (\{\lambda_{i}\}\times L)$, which in turn implies that
$\mathfrak{l}$ is tangent to the branch locus of
\eqref{eq:double.cover}  at the image of the point
$(\lambda_{i},x)$. Thus the point $(\lambda_{i},x)$ is a point of the
ramification divisor of \eqref{eq:double.cover}, and so our isomorpism
$\Lambda \to L$ sends $\lambda_{i}$ to $x$. Thus the composition of the map
$\Lambda \to L$ with the map $\beta_{0}$ is an isomorphism $\iota :
\Lambda \to C$ which sends the set $\{\lambda_{i}\}_{i=1}^{5}$ to the
set $\{p_{i}\}_{i=1}^{5}$. 

To finish the proof of the propositon it only remains to show that the
bijection of the set $\{\lambda_{i}\}_{i=1}^{5}$ and the set
$\{p_{i}\}_{i=1}^{5}$ that we just constructed is the same as the
identification of these two sets with the vertices of the singular
quadrics in the pencil. 

To argue that this is the case we will give an alternative
construction of $\iota$.  Consider the
squaring map $\sq : \mathbb{P}^{4} \to \mathbb{P}^{4}$
given by 
\[
\sq(x_{1}:x_{2}:x_{3}:x_{4}:x_{5}) =
(x_{1}^{2}:x_{2}^{2}:x_{3}^{2}:x_{4}^{2}:x_{5}^{2}).
\] 
This is a
branched Galois cover with Galois group $(\mathbb{Z}/2)^{4}$  which is
branched along the hyperplanes $y_{i} = 0$, for $i = 1, \ldots,
5$. Here $(y_{1}:y_{2}:y_{3}:y_{4}:y_{5})$ are the homogeneous
coordinates in the target $\mathbb{P}^{1}$. 

The squaring map $\sq$ sends each quadric $Q_{\lambda}$ to a
hyperplane $H_{\lambda}$ in the target $\mathbb{P}^{4}$ and the del
Pezzo $X \subset \mathbb{P}^{4}$ to the plane $\mathbb{P}^{2} \subset
\mathbb{P}^{4}$ cut out by the pencil $\{H_{\lambda}\}_{\lambda \in
  \Lambda}$. Now fix a $\lambda \in \Lambda$. The five coordinate
hyperplanes $H_{i} := \{y_{i} = 0\}$ intersect $H_{\lambda}$ in planes
$\bP_{i,\lambda}$. Passing to the dual projective three space
$H_{\lambda}^{\vee}$ we get six points: a point $b \in
H_{\lambda}^{\vee}$ corresponding to 
the base locus $\mathbb{P}^{2} = S^{2}C$ of 
the pencil $\{H_{\lambda}\}_{\lambda \in
  \Lambda}$, and five more points $h_{\lambda,i} \in
H_{\lambda}^{\vee}$  corresponding to the $\bP_{i,\lambda}$'s. The six
points $b$, $h_{\lambda,1}$, \ldots, $h_{\lambda,5}$ sit on a unique
twisted cubic $\mathfrak{C}_{\lambda} \subset H_{\lambda}^{\vee}$. Now
projecting out of $b$ we get a rational map  
$H_{\lambda}^{\vee} \dashrightarrow (\mathbb{P}^{2})^{\vee}$  whose
restriction to $\mathfrak{C}_{\lambda}$ gives an isomorphism from
$\mathfrak{C}_{\lambda}$  to the conic $C \subset S^{2}C =
\mathbb{P}^{2}$. Under this isomorphism each of the five points
$h_{\lambda,i}$ goes to the corresponding point $p_{i} \in C$ and the
point $b$ goes to a point $\iota(\lambda) \in C$. In this way we get
an isomorphism $\iota : \Lambda \to C$ which manifestly
sends $\lambda_{i}$ to $p_{i}$. Note that this map must be equal to
the previously constructed one since their composition is an
automorphism of $\Lambda$ which preserves the set of $\lambda_{i}$'s
and so must be the identity.  The proposition is proven. \ \hfill $\Box$

\

\bigskip

\begin{rem} \label{rem:explicit.map} The isomorphism $\iota$
  constructed in the previous proposition can be computed explicitly
  in coordinates. To simplify the calculation we will choose a better
coordinate system on the anti-canonical $\mathbb{P}^{4}$. By a linear
change of variables we can rewrite $X$ as the intersection of two
quadrics given by 
\[ 
\begin{aligned}
x_{1}^{2} & + x_{3}^{2} + x_{4}^{2} + x_{5}^{2} = 0, \\
x_{2}^{2} & + x_{3}^{2} + \lambda_{4}x_{4}^{2} + \lambda_{5}x_{5}^{2} = 0.
\end{aligned}
\]
The parameter line $\Lambda$ for the pencil generated by these two
quadrics has a natural coordinatization with an affine coordinate
$\lambda$. The quadric in the pencil corresponding to $\lambda \in
\mathbb{P}^{1}$ is given by
\[
Q_{\lambda} := \lambda x_{1}^{2} -  x_{2}^{2}  + (\lambda - 1)
x_{3}^{2} + (\lambda - \lambda_{4})x_{4}^{2} + (\lambda -
\lambda_{5})x_{5}^{2} = 0.
\]
The pencil has five singular fibers sitting over the points $0$,
$\infty$, $1$, $\lambda_{4}$, and $\lambda_{5}$ in $\mathbb{P}^{1} =
\Lambda$. 

The squaring map $\sq$ sends the del Pezzo $X \subset \mathbb{P}^{4}$ to
the plane $\mathbb{P}^{2} \subset \mathbb{P}^{4}$ which is cut by the
linear equations
\begin{equation} \label{eq:plane.in.P4}
\begin{aligned}
y_{1} & + y_{3} + y_{4} + y_{5} = 0, \\
y_{2} & + y_{3} + \lambda_{4}y_{4} + \lambda_{5}y_{5} = 0.
\end{aligned}
\end{equation}
In this $\mathbb{P}^{2}$ we have:
\begin{itemize}
\item The five branch lines for the covering map $\beta_{0} =
\sq_{|X} : X \to \mathbb{P}^{2}$. These lines are the intersection of
the $\mathbb{P}^{2}$  given by \eqref{eq:plane.in.P4} with the
coordinate hyperplanes, i.e. $L_{i} := \{ y_{i} = 0 \}\cap
\mathbb{P}^{2}$, for  $i = 1,
\ldots, 5$. 
\item The conic $C \subset \mathbb{P}^{2}$ tangent to these five lines. 
\item The five points $p_{1}, \ldots, p_{5} \in C$ such that 
\[
T_{p_{i}}C = \{ y_{i} = 0 \}\cap \mathbb{P}^{2}, \text{ for } i = 1,
\ldots, 5.
\]
\end{itemize}
We want to show that $\Lambda$ and $C$ are identified by a linear
isomorphism $C\to \Lambda$ for which
$\iota(p_{1}) = 0$, $\iota(p_{2}) = \infty$, $\iota(p_{3}) = 1$,
$\iota(p_{4}) = \lambda_{4}$, and $\iota(p_{5}) =
\lambda_{5}$. In other words we have to show that the unique
coordinate function $u : C \to
\mathbb{P}^{1}$ which satisfies $u(p_{1}) = 0$, $u(p_{2}) = \infty$,
$u(p_{3}) = 1$ also satisfies $u(p_{4}) = \lambda_{4}$ and $u(p_{5}) =
\lambda_{5}$. 

To do this we will use $y_{3}$, $y_{4}$, and $y_{5}$ as the
homogeneous coordinates on  the $\mathbb{P}^{2}$ given by 
\eqref{eq:plane.in.P4}. To simplify the notation we will relabel:
$x = y_{3}$, $y = y_{4}$, and $z = y_{5}$. In these coordinates the
five branch lines of the map $\beta_{0} : X \to \mathbb{P}^{2}$ have
equations: 
\[
\begin{aligned}
L_{1} : \quad x + &  \; y + z = 0 \\
L_{2} : \quad x + &  \; \lambda_{4}y + \lambda_{5}z = 0 \\
L_{3} : \quad x = &  \; 0 \\
L_{4} : \quad y = &  \; 0 \\
L_{5} : \quad z = &  \; 0.
\end{aligned}
\]
Let $(\mathbb{P}^{2})^{\vee}$ be the dual plane with dual
coordinates $(x^{\vee}:y^{\vee}:z^{\vee})$. In this dual plane the
lines $L_{i}$ correspond to the points 
\[
\begin{aligned}
\ell_{1} & = (1:1:1) \\
\ell_{2} & = (1:\lambda_{4}:\lambda_{5}) \\
\ell_{3} & = (1:0:0) \\
\ell_{4} & = (0:1:0) \\
\ell_{5} & = (0:0:1).
\end{aligned}
\]
The projectively dual conic $C^{\vee} \subset (\mathbb{P}^{2})^{\vee}$
is the unique conic passing through the points $\ell_{1}$, $\ell_{2}$,
$\ell_{3}$, $\ell_{4}$, $\ell_{5}$. The natural isomorphism $C \to
C^{\vee}$, given by $p \mapsto T_{p}C$ sends $p_{i}$ to
$\ell_{i}$. Therefore it suffices to show that the coordinate function
$u : C^{\vee} \to \mathbb{P}^{1}$ which satisfies $u(\ell_{1}) = 0$,
$u(\ell_{2}) = \infty$, $u(\ell_{3}) = 1$, also satisfies $u(\ell_{4})
= \lambda_{4}$ and $u(\ell_{5}) = \lambda_{5}$. 

Since $(1:0:0)$, $(0:1:0)$, $(0:0:1)$ are in $C^{\vee}$ we have that
\[
C^{\vee} : \quad \alpha y^{\vee}z^{\vee} + \beta x^{\vee}z^{\vee} +
\gamma x^{\vee}y^{\vee} = 0
\]
for some $\alpha, \beta, \gamma \in \mathbb{C}$. It is not hard to
write the desired affine coordinate $u$ for such a conic. To do this
note first that  if  $(s:t)$
are homogeneous coordinates on 
$\mathbb{P}^{1}$, then the map 
\[
\xymatrix@R-1pc{
\mathfrak{c} : & \hspace{-1pc} 
\mathbb{P}^{1} \ar[r] & (\mathbb{P}^{2})^{\vee} \\
& \hspace{-1pc} (s:t) \ar@{|->}[r] & \left(\alpha t (t-s):\beta
s(s-t): \gamma st\right)
}
\]
parametrizes $C^{\vee}$ in such a way that $\mathfrak{c}(0:1) =
\ell_{3}$, $\mathfrak{c}(1:0) = \ell_{4}$, $\mathfrak{c}(1:1) =
\ell_{5}$. In other words the ratio $v := s/t$ is an affine
coordinate on $C^{\vee}$ such that $v(\ell_{3}) = 0$, $v(\ell_{4}) =
\infty$, $v(\ell_{5}) = 1$. On the affine coordinate $u : C^{\vee} \to
\mathbb{P}^{1}$ that we want should satisfy $u(\ell_{1}) = 0$,
$u(\ell_{2}) = \infty$, $u(\ell_{3}) = 1$. Therefore $u$ and $v$ are
related by the fractional-linear transformation
\begin{equation} \label{eq:frac.lin}
u = \frac{\frac{v}{v(\ell_{1})} - 1}{\frac{v}{v(\ell_{2})} - 1}.
\end{equation}
\

\noindent
Our next task is to compute $u(\ell_{3})$ and $u(\ell_{4})$. First we
will need explicit formulas for the coefficients $\alpha$, $\beta$,
$\gamma$ of the conic $C^{\vee}$. We have that
\[
\begin{aligned}
(1:1:1) \in C^{\vee} \quad & \text{gives} \quad \alpha+\beta+\gamma = 0,
  \\
(1:\lambda_{4}:\lambda_{5}) \in C^{\vee} \quad & \text{gives} \quad
  \lambda_{4}\lambda_{5}\alpha+ \lambda_{5}\beta+\lambda_{4}\gamma = 0.
\end{aligned}
\]
So $\alpha = -\beta - \gamma$ and $\beta\lambda_{5}(1 - \lambda_{4}) +
\gamma\lambda_{4}(1-\lambda_{5}) = 0$ and up to scale we have
\[
\begin{aligned}
\alpha & = \lambda_{5} - \lambda_{4} \\
\beta & = \lambda_{4}(1 - \lambda_{5}) \\
\gamma & = \lambda_{5}(\lambda_{4} - 1).
\end{aligned}
\]
Therefore the dual conic $C^{\vee}$ is given by the equation
\[
C^{\vee} : \quad (\lambda_{5} - \lambda_{4})y^{\vee}z^{\vee} +
\lambda_{4}(1 - \lambda_{5}) x^{\vee}z^{\vee} +
\lambda_{5}(\lambda_{4} - 1) x^{\vee}y^{\vee} = 0,
\]
and so the parametrization $\mathfrak{c}$ is given by
\[
\mathfrak{c}(s:t) = \left( (\lambda_{5} - \lambda_{4})t
(t-s):\lambda_{4}(1 - \lambda_{5}) 
s(s-t): \lambda_{5}(\lambda_{4} - 1) st\right).
\]
Now solving the equations $\mathfrak{c}(s:t) = (1:1:1)$ and
$\mathfrak{c}(s:t) = (1:\lambda_{4}:\lambda_{5})$ we get
\[
\begin{aligned}
\ell_{1} & = \mathfrak{c}\left(\lambda_{4}-\lambda_{5}:\lambda_{4}(1 -
\lambda_{5})\right) \\
\ell_{2} & = \mathfrak{c}\left(\lambda_{5} - \lambda_{4}:\lambda_{5} -
1 \right).
\end{aligned}
\]
Thus we have
\[
\begin{aligned}
v(\ell_{1}) & = \frac{\lambda_{4}-\lambda_{5}}{\lambda_{4}(1 -
\lambda_{5})} \\
v(\ell_{2}) & = \frac{\lambda_{5} - \lambda_{4}}{\lambda_{5} -
1} \\
v(\ell_{3}) & = 0 \\
v(\ell_{4}) & = \infty \\
v(\ell_{5}) & = 1
\end{aligned}
\]
which after substituting in \eqref{eq:frac.lin} gives 
\[
\begin{aligned}
u(\ell_{1}) & = 0 \\
u(\ell_{2}) & = \infty \\
u(\ell_{3}) & = 1 \\
u(\ell_{4}) & = \lambda_{4} \\
u(\ell_{5}) & = \lambda_{5}.
\end{aligned}
\]
This completes the calculation of $\iota$ in coordinates.
\end{rem}

\bigskip

\

\noindent
With the isomorphism $\iota : \Lambda \to C$ at our disposal we are
now ready to describe the compactified and resolved Hecke
correspondence. Consider the blow-up $\op{Bl}_{\Delta}(X\times X)$ of
$X\times X$ along the diagonal $\Delta \subset X\times X$. Write
$\beps_{1}, \beps_{2} : \op{Bl}_{\Delta}(X\times X) \to X\times X \to
X$ for the composition of the blow-up map with the projections to the
first and second factor respectively. The fiber of the map $\beps_{1}
: \op{Bl}_{\Delta}(X\times X) \to X$ over a point $x$ is the $dP_{6}$
del Pezzo surface $\op{Bl}_{x}X$. When $x \in X - \cup_{I} L_{I}$ does
not lie on one of the $16$ lines in $X \subset \mathbb{P}^{4}$, the
anti-canonical map $\op{Bl}_{x}(X) \to \mathbb{P}^{3}$ is an embedding
and the exceptional divisor $L_{x} \subset \op{Bl}_{x}X$ corresponding
to $x$ is a line in $\mathbb{P}^{3}$. Projecting from $L_{x}$
gives a morphism from $\op{Bl}_{x}X$ to $\mathbb{P}^{1}$. The fibers
of this map are strict transforms of anti-canonical sections of $X$
having a node at $x$. But the parameter space for the pencil of such
sections can be identified canonically with $\Lambda$. Indeed, if
$\lambda \in \Lambda$ and $x \in X$ the intersection
$T_{x}Q_{\lambda}\cap X$ is an anti-canonical curve with a node at
$x$. Thus the assignment $\lambda \mapsto T_{x}Q_{\lambda}\cap X$
gives a parametrization of the pencil of anti-canonical sections of
$X$ with singularity at $x$. In this way we obtain a morphism
$\op{Bl}_{x}X \to \mathbb{P}^{1} = \Lambda$ and after composing with
$\iota : \Lambda \to C$ we get a morphism $\br_{x} : \op{Bl}_{x}X \to
C$. Note however that if $x$ lies on a line $L_{I}$ in $X \subset
\mathbb{P}^{3}$, the anti-canonical image of $\op{Bl}_{x}X$ is
singular and the map $\br_{x}$ is given by a linear pencil which has a
base component $L_{I}$, i.e. $\br_{x}$ is not a morphism in this case.

So we get a well  defined rational map 
\begin{equation} \label{eq-rationalnu}
\op{Bl}_{\Delta}(X\times X) \dashrightarrow X\times X\times C
\end{equation}
by setting
\[
a \mapsto
\left(\beps_{1}(a),\beps_{2}(a),
\br_{\beps_{1}(a)}(\beps_{2}(a))\right),
\]
which is a morphism outside of the strict transforms
$\widehat{L_{I}\times L_{I}}$ of the surfaces $L_{I}\times L_{I}
\subset X\times X$ in $\op{Bl}_{\Delta}(X\times X)$.

Since $\Delta$
intersects each surface $L_{I}\times L_{I}$ in the diagonal copy of
$L_{I}$, it follows that $\widehat{L_{I}\times L_{I}} \to L_{I}\times
L_{I}$ is an isomorphism, and that the surfaces $\widehat{L_{I}\times
  L_{I}}$ are all disjoint in $\op{Bl}_{\Delta}(X\times X)$.  
Blowing up the surfaces  $\widehat{L_{I}\times
  L_{I}}$ resolves the singularities of the map \eqref{eq-rationalnu}
and gives a morphism 
\[
\nu \; : \; \op{Bl}_{\sqcup_{I} \widehat{L_{I}\times
  L_{I}}} \op{Bl}_{\Delta}(X\times X) \longrightarrow X\times X\times C
\]
To simplify notation we set 
\begin{equation} \label{eq-H}
H := \op{Bl}_{\sqcup_{I} \widehat{L_{I}\times
  L_{I}}} \op{Bl}_{\Delta}(X\times X).
\end{equation}
The map $\nu$ allows us to identify $H$ with a resolution of the
compactified Hecke correspondence. More precisely we have the
following theorem:

\begin{theo} \label{theo:describe.hecke} The morphism 
$\nu \; : \; H  \to X\times X\times C$ 
maps $H$ birationally onto the closure of 
$\Heck^{\op{stable}}(0,1/2) \subset X\times X\times C$.
\end{theo}
{\bfseries Proof.} From the description of $\nu$ we see that the image of
$\nu$ can be described as the subvariety of $X\times X\times C$
parametrizing all triples of points $(x,y,p) \in X\times X\times C$
such that $y$ is a point on the anti-canonical curve
$\boldsymbol{\mathfrak{T}}_{x,p} := T_{x}Q_{\iota^{-1}(p)}\cap X$.
Therefore to 
prove the theorem we only need to show that for a general $x$ and $p$
the curve $\boldsymbol{\mathfrak{T}}_{x,p} \subset X$ is equal to the
Hecke line 
$\bbH_{x,p} \subset  X$
parametrizing all Hecke transforms  at $p \in C$ of the parabolic bundle
$\mathcal{V}_{x}$ corresponding to $x \in X$.

A necessary condition for the equaility
$\boldsymbol{\mathfrak{T}}_{x,p} = \bbH_{x,p}$ is that $\bbH_{x,p}
\subset X$ is an anti-canoinical section of $X$ which is singular at
$x$. To check this we will describe the Hecke line explicitly. Let $x
\in X$ and let $\mathcal{V}_{x}$ be the corresponding parabolic
bundle. Since the difference of the weights is fixed to be $q = 1/2$ the
parabolic bundle $\mathcal{V}_{x}$ is described by the data $(V,\bF)$,
where $V$ is a holomorphic rank two vector bundle with trivial
determinant on $C \cong \mathbb{P}^{1}$ and $\bF = \{ F_{i}\}_{i =
1}^{5}$, where $F_{i} \in \mathbb{P}(V_{p_{i}})$. In the proof of 
Theorem~\ref{theo:shapes.of.moduli} we showed that the moduli space
$X$ is glued together out of natural open charts $\widehat{U}_{I}$ where $I
\subset \{1,2,3,4,5\}$ is a subset with  three elements. By genericity  
we may assume that $x$ is in one fixed chart, e.g. that $x \in
\widehat{U}_{\{1,2,3\}} \subset X$. Specifically this means that
$\mathcal{V}_{x} = (V,\bF)$, the bundle $V$ is trivializable and:
\begin{itemize}
\item we identify $C$ with the standard $\mathbb{P}^{1}$ so that
  $p_{1}$, $p_{2}$, $p_{3}$ are idenitified with $0$, $1$, $\infty$
  respectively; 
\item we trivialize the bundle $V$ and identify the projectivization
  of its fiber with the  standard $\mathbb{P}^{1}$ so that lines
  $F_{1}$, $F_{2}$, $F_{3}$ are idenitified with $0$, $1$, $\infty$
  respectively; 
\item under the identification of $\mathbb{P}(V) \cong C \times
\mathbb{P}^{1}$ with 
$\mathbb{P}^{1}\times \mathbb{P}^{1}$ the parabolic data $\bF$
becomes a five tuple of points $(0,0)$, $(1,1)$, $(\infty,\infty)$,
$(p_{4},f_{4})$, $(p_{5},f_{5})$;  
\item the chart
$\widehat{U}_{\{1,2,3\}}$ is the blow-up of $\mathbb{P}^{1}\times
  \mathbb{P}^{1} - \{(0,0),(1,1),(\infty,\infty)\}$ at the point
  $(p_{4},p_{5})$; 
\item by genericity we may assume that $x$ belongs to the open subset
\[
\mathbb{P}^{1}\times 
  \mathbb{P}^{1} - \{(0,0),(1,1),(\infty,\infty),(p_{4},p_{5})\}
  \text{  of  } 
\widehat{U}_{\{1,2,3\}},
\] 
this means that $(p_{4},p_{5}) \neq
  (f_{4},f_{5})$ and 
\[
x = (f_{4},f_{5}) \in \mathbb{P}^{1}\times 
  \mathbb{P}^{1} - \{(0,0),(1,1),(\infty,\infty),(p_{4},p_{5})\}.
\]
\end{itemize}
Now suppose we are given a sixth point $p \in C$ which is distinct
from all  $p_{i}$. To specify a Hecke transform of $\mathcal{V}_{x}$ at $p \in
C$ we need to specify a flag $F \subset V_{p}$, i.e by points $f \in
\mathbb{P}^{1} = \mathbb{P}(V_{p})$. In particular we have
a natural morphism
\[
\mathfrak{hecke}_{x,p} : \mathbb{P}^{1} \to X
\]
whose image is the Hecke line $\bbH_{x,p} \subset
X$. It  is not hard to analyze this map. If
$f \in
\mathbb{P}^{1}$, then the Hecke transform of $\mathcal{V}_{x}$ centered at
$f$ is a parabolic bundle in $N_{1}(0,1/2)$ which then can be
transported to $N_{0}(0,1/2) = X$ via the isomorphism from
Lemma~\ref{lem:moduli.isos}. The resulting point in $X$ is  precisely
the point $\mathfrak{hecke}_{x,p}(f)$ of the Hecke line $\bbH_{x,p}$
corresponding to $x$. In other words $\mathfrak{hecke}_{x,p}(f)$ is
uniquely characterized by the property that if $\mathcal{V}_{x} =
(V,\bF)$, then $\mathcal{V}_{\mathfrak{hecke}_{x,p}(f)} = (V',\bF')$,
where: 
\begin{itemize}
\item $V'$ is the tensor product of the 
Hecke transform of $V$ at the points $p, p_{1},
\ldots, p_{5}$ centered at $F, F_{1}, \ldots, F_{5}$ and 
$\mathcal{O}_{C}\left(p + p_{1} + \cdots + p_{5}\right) \cong
\mathcal{O}_{C}(3)$, i.e.
\[
V' := \left(\ker\left[ V \to V_{p}/F\right]\cap \left(\bigcap_{i = 1}^{5}
\ker\left[ V \to V_{p_{i}}/F_{i}\right]\right)\right) \otimes
\mathcal{O}_{C}(3),
\]
\item $\bF'$ are the induced flags in $V'_{p_{i}}$, i.e.
\[
F_{i}' := \op{im}\left[(V(-p_{i}))_{p_{i}} \to V'_{p_{i}}\right].
\]
\end{itemize}
Recall (see e.g. \cite[Chapter~III]{beauville}, \cite{maruyama},
\cite[Appendix~A]{dopw}) that from the point of view of the ruled
surface $\mathbb{P}(V)$ the Hecke transform of $V$ centered at $f \in
\mathbb{P}(V_{p})$ amounts to blowing up the point $f \in
\mathbb{P}(V)$, and then contracting the strict transform of the fiber
$\mathbb{P}(V_{p})$. Therefore we can describe the point
$\mathfrak{hecke}_{x,p}(f)$ synthetically from the point of view of
the surface $\mathbb{P}(V) \cong \mathbb{P}^{1}\times \mathbb{P}^{1}$
as follows:
\begin{itemize}
\item Consider the pencil  $P \cong \mathbb{P}^{1}$ of all curves on
  $\mathbb{P}^{1}\times 
  \mathbb{P}^{1}$ in the linear system $\mathcal{O}(3,1)$ passing
  through the six points $(p,f)$, $(0,0)$, $(1,1)$, $(\infty,\infty)$,
  $(p_{4},f_{4})$, $(p_{5},f_{5})$.
\item Let $Z$ be the blow-up of $\mathbb{P}^{1}\times \mathbb{P}^{1}$
  at the six points $(p,f)$, $(0,0)$, $(1,1)$, $(\infty,\infty)$,
  $(p_{4},f_{4})$, $(p_{5},f_{5})$. The linear system gives a natural
  morphism $\pi : Z \to P$ with six singular fiber. Thus we have six
  natural marked points in $P$ -  the discriminant of the map $\pi$. 
\item If $x$ is generic each singular fibers of $\pi$ has two rational
  components - the strict transform of the $(1,0)$ ruling through one
  of the six points and the strict transform of the unique $(2,1)$
  curve through the remaining five points. Therefore we get a
  bijection between the six marked points in $P$ and the six points
  $p$, $0$, $1$, $\infty$, $p_{4}$, $p_{5}$ in $C \cong
  \mathbb{P}^{1}$. Therefore we can label the six points of
  $\op{Discr}(\pi)$ as $\mathfrak{d}_{p}$, $\mathfrak{d}_{0}$,
  $\mathfrak{d}_{1}$, $\mathfrak{d}_{\infty}$, $\mathfrak{d}_{p_{4}}$,
  $\mathfrak{d}_{p_{5}}$.
\item  Identify $P$ and the standard $\mathbb{P}^{1}$
  so that  $\mathfrak{d}_{0}$,
  $\mathfrak{d}_{1}$, $\mathfrak{d}_{\infty}$ are mapped to $0$, $1$, and
  $\infty$ respectively. Then $\mathfrak{hecke}_{x,p}(f)$ is the point
  in $\mathbb{P}^{1}\times 
  \mathbb{P}^{1} - \{(0,0),(1,1),(\infty,\infty),(p_{4},p_{5})\}$
  given by 
\[
\mathfrak{hecke}_{x,p}(f) = \left( \mathfrak{d}_{p_{4}},
\mathfrak{d}_{p_{5}} \right) \in \mathbb{P}^{1}\times 
  \mathbb{P}^{1} - \{(0,0),(1,1),(\infty,\infty),(p_{4},p_{5})\}.
\]
\end{itemize}
The curve $\bbH_{x,p}$ is traced by all points
$\mathfrak{hecke}_{x,p}(f) = \left( \mathfrak{d}_{p_{4}},
\mathfrak{d}_{p_{5}} \right)$ obtained by the above procedure with a
varying $f$. 

This synthetic description of the Hecke line allows us to parametrize
$\bbH_{x,p}$ explicitly in the natural coordinates on
$\mathbb{P}^{1}\times \mathbb{P}^{1} -
\{(0,0),(1,1),(\infty,\infty),(p_{4},p_{5})$. Indeed let $u$ and $v$
be the standard base and fiber affine coordinates on $\mathbb{P}(V)
\cong \mathbb{P}^{1}\times \mathbb{P}^{1}$. The pencil $P$ consists of
all polynimials 
\[
\sum_{\substack{0 \leq i \leq 3 \\ 0 \leq j \leq 1}}
a_{ij}u^{i}v^{j}
\] 
satisfying
\begin{equation} \label{eq:equations.for.P}
\sum a_{ij} p^{i}f^{j}  = 0, \;
a_{00} = 0,  \;
a_{31} = 0,  \;
\sum a_{ij} = 0, \;
\sum a_{ij} p_{4}^{i}f_{4}^{j}  = 0, \;
\sum a_{ij} p_{5}^{i}f_{5}^{j} = 0. 
\end{equation}
Furthermore:
\begin{itemize}
\item $\mathfrak{d}_{0}$ corresponds to the unique up to scale
  polynomial satisfying  \eqref{eq:equations.for.P}
  and $a_{01} = 0$, 
\item $\mathfrak{d}_{1}$ corresponds to the unique up to scale
  polynomial satisfying  \eqref{eq:equations.for.P} and $\sum_{i =
  0}^{3} a_{i0} = 0$,
\item $\mathfrak{d}_{\infty}$ corresponds to the unique up to scale
  polynomial satisfying  \eqref{eq:equations.for.P} and $a_{30} = 0$.
\end{itemize}
We can normalize these three polynomials by requiring that in each of them the
coefficient in front of $u^{2}$ is equal to $1$. After that we can
rescale the polynomials corresponding to $\mathfrak{d}_{0}$ and
$\mathfrak{d}_{\infty}$ so that the polynomial corresponding to
$\mathfrak{d}_{1}$ is equal to their sum. Let $\mathfrak{v}_{0}(u,v)$, 
$\mathfrak{v}_{\infty}(u,v)$, and $\mathfrak{v}_{1}(u,v)$ deote the
resulting rescalled and normalized polynomials. Let $\phi$ denote the affine
coordinate on $P$ which takes values $0$, $1$, $\infty$ at 
$\mathfrak{d}_{0}$, $\mathfrak{d}_{1}$, $\mathfrak{d}_{\infty}$. Then 
\[
\phi = -\frac{\mathfrak{v}_{0}(u,v)}{\mathfrak{v}_{1}(u,v)}.
\]
In terms of $\phi$ the parametrization of $\bbH_{x,p}$ is
given by the rational map
\[
\xymatrix@R-1.5pc{
\mathfrak{hecke}_{x,p} : &  \mathbb{P}^{1} \ar@{-->}[r] &
\mathbb{P}^{1}\times \mathbb{P}^{1} \\
& f \ar@{-->}[r] & \left( \phi(p_{4},f_{4}), \phi(p_{5},f_{5}) \right).
}
\]
Solving the corresponding linear systems, normalizing and rescalling
we get that
\[
\begin{aligned}
\phi(p_{4}, &  f_{4}) = \\
& \frac{(f + f_{5}p - pfp_{5} + p_{5} -
f_{5})(-pff_{4}p_{5} - p_{4}f_{4}f_{5}p + pff_{5}p_{4} +
pf_{4}p_{5}f_{5} - p_{4}fp_{5}f_{5} + p_{4}f_{4}fp_{5})}{(fp_{5}f_{5} +
f_{5}p + pfp_{5} - pp_{5}f_{5} - fp_{5} - f_{5}pf)(-f_{5}p + f_{4}p +
f_{5}p_{4} - f_{4}p_{5} - fp_{4} + fp_{5})} \\ 
\phi(p_{5}, & f_{5})  = \\
& \frac{(f_{4}p - fp_{4} - f_{4} - p + p_{4} + f)(-pff_{4}p_{5} -
p_{4}f_{4}f_{5}p + pff_{5}p_{4}+pf_{4}p_{5}f_{5}
-p_{4}fp_{5}f_{5}+p_{4}f_{4}fp_{5})}{(f_{4}p-fp_{4}-pp_{4}f_{4}+pfp_{4}
- f_{4}pf+fp_{4}f_{4})(-f_{5}p
+f_{4}p+f_{5}p_{4}-f_{4}p_{5}-fp_{4}+fp_{5})}.
\end{aligned}
\]
From these formulas it is not hard to locate the singularity of
$\bbH_{x,p}$. In general consider a rational map
$\xymatrix@1{\mathbb{P}^{1} \ar@{-->}[r] & \mathbb{P}^{1}\times
\mathbb{P}^{1}}$ given by $f \to (A(f)/B(f),C(f)/D(f))$, with $A$,
$B$, $C$, and $D$ quadratic polynomials in $f$. Reading off their
coefficients, we get a $4\times 3$ matrix. Let $(a, b, c, d)$ generate
its kernel, i.e. $aA + bB + cC + dD =0$. Then the image is a $(2,2)$
curve singular at $(-b/a, -d/c)$. Applying this to the map
$\mathfrak{hecke}_{x,p}$ we  immediately get that the image of
$\mathfrak{hecke}_{x,p}$  is
singular at the point $(f_{4},f_{5}) \in \mathbb{P}^{1}\times
\mathbb{P}^{1}$. In other words $\bbH_{x,p}$ is singular precisely at
$x$. 

Consider now the pencil of $(2,2)$ curves on $\mathbb{P}^{1}\times
\mathbb{P}^{1}$ that pass through $(0,0)$, $(1,1)$, $(\infty,\infty)$,
$(p_{4},p_{5})$, $(f_{4},f_{5})$ and in addition have a singularity at
$(f_{4},f_{5})$. We have argued that both $\bbH_{x,p}$ and
$\boldsymbol{\mathfrak{T}}_{x,p}$ belong to this linear pencil. In
other words we have two parametrizations of the parameter line for
this pencil sending $p \in C$ to the curves $\bbH_{x,p}$ and
$\boldsymbol{\mathfrak{T}}_{x,p}$ respectively. To prove the theorem
we have to show that these two parametrizations coincide (up to an
automorphism of $\mathbb{P}^{1}$). This can be seen as follows. A
natural way to label the members of this pencil is by the points of a
line on $X$. More precisely, suppose $\boldsymbol{e}_{(0,0)} \subset X
\subset \mathbb{P}^{4}$ is the line which under the natural map $X \to
\mathbb{P}^{1}\times \mathbb{P}^{1}$ contracts to the point
$(0,0)$. Since each curve in our pencil is a hyperplane section of
$X$, it follows that each such curve interects the line
$\boldsymbol{e}_{(0,0)}$ at a single point. Thus we get two morphisms
\[
\begin{aligned}
\mathfrak{h}_{x} : \; & \; C \to \boldsymbol{e}_{(0,0)}, \qquad 
\mathfrak{h}_{x}(p)
:= \bbH_{x,p}\cap  \boldsymbol{e}_{(0,0)}, \\[5pt]
\mathfrak{t}_{x} : \; & \; C \to \boldsymbol{e}_{(0,0)}, \qquad 
\mathfrak{t}_{x}(p)
:= \mathfrak{T}_{x,p}\cap  \boldsymbol{e}_{(0,0)}.
\end{aligned}
\]
The two curves $C$ and $\boldsymbol{e}_{(0,0)}$ are both isomorphic to
$\mathbb{P}^{1}$, and the theorem will follow if we can show that the
maps $\mathfrak{h}_{x}$ and $\mathfrak{t}_{x}$ are linear
isomorphisms. This is obvious for the map $\mathfrak{t}_{x}$ since the
family of hyperplanes $\{T_{x}Q_{p}\}_{p \in C}$ depends linearly on
$p$ and $\mathfrak{t}_{x}(p) = T_{x}Q_{p}\cap
\boldsymbol{e}_{(0,0)}$. To show that $\mathfrak{h}_{x} $ is linear we
will compute it explicitly. The Hecke line $\bbH_{x,p}$ is the image
of the map  $\xymatrix@1{
\mathfrak{hecke}_{x,p} : &  \mathbb{P}^{1} \ar@{-->}[r] &
\mathbb{P}^{1}\times \mathbb{P}^{1}}$, which in the affine coordinate
$f$ is given by  $\mathfrak{hecke}_{x,p}(f) = \left(
\phi(p_{4},f_{4})(f),\phi(p_{5},f_{5})(f)\right)$. 
The curve $\boldsymbol{e}_{(0,0)}$ parameterizes lines through the
origin in the $\mathbb{C}^{2}$ with coordinates $(u,v)$, and so we can
choose as an affine coordinate on $\boldsymbol{e}_{(0,0)}$ the slope
of such a line. Suppose that $f_{0}$ is the unique point on the
$f$-line which goes to $(0,0)$ under the map
$\mathfrak{hecke}_{x,p}$. Then 
\[
\mathfrak{h}_{x}(p) = \frac{\frac{\partial
    \phi(p_{5},f_{5})}{\partial f}_{|f=f_{0}}}{\frac{\partial
    \phi(p_{4},f_{4})}{\partial f}_{|f=f_{0}}}.
\] 
As we saw above the
components of $\mathfrak{hecke}_{x,p}$ factor as
\[
\begin{aligned}
\phi(p_{4},f_{4})(f) & = \bF_{1}(f)\cdot\balpha(f), \\
\phi(p_{5},f_{5})(f) & = \bF_{2}(f)\cdot\balpha(f),
\end{aligned}
\]
where $\bF_{1}$, $\bF_{2}$, $\balpha$ are the fractional linear
transformations in $f$ given by:
\[
\begin{aligned}
\bF_{1}(f) & = \frac{(1- p f_{5})f + f_{5}p + p_{5} - f_{5}}
{(p_{5}f_{5} - f_{5} - p_{5} - f_{5}p)f +  f_{5} p  - pp_{5} f_{5}}
\\[1.5pc]
\bF_{2}(f) & = \frac{(1- p f_{4})f + f_{4}p + p_{4} - f_{4}}
{(p_{4}f_{4} - f_{4} - p_{4} - f_{4}p)f +  f_{4} p  - pp_{4} f_{4}}
\\[1.5pc]
\balpha(f) & = \frac{(-pf_{4} p_{5} + pf_{5}p_{5} f_{5} + p_{4}
  f_{4}p_{5})f 
- p_{4} f_{4} f_{5} p 
  p_{4} + pf_{4} p_{5} f_{5} - p_{4} 
  }
{(p_{5}-p_{4})f -f_{5} p + f_{4} p + f_{5} p_{4} - f_{4} p_{5}}.
\end{aligned}
\]
Therefore $f_{0}$ is the unique solution of $\balpha(f) = 0$ and 
\[
\mathfrak{h}_{0} = \frac{\bF_{2}(f_{0})}{\bF_{1}(f_{0})}.
\]
Solving for $f_{0}$ and substituting in the above formulas we get that 
\[
\mathfrak{h}_{x}(p) =\frac{f_{5}\left(pf_{4}p_{5} - pf_{5}p_{4} -
pf_{4}^{2}p_{5}  + p_{4}f_{4}f_{5}p - f_{4}p_{5}f_{5} + f_{4}^{2}p_{5}
- p_{4}f_{4}p_{5}   + p_{4}p_{5}f_{5}\right)}{f_{4}\left(-p
  f_{4}p_{5}f_{5} + p_{4}pf_{5}^{2}  + pf_{4}p_{5} - p f_{5}p_{4} -
  f{5}^{2}p_{4} + p_{4}p_{5}f_{5} + f_{5}p_{4}f_{4} - p_{4}f_{4}p_{5}\right)}
\]
which is manifestly linear in $p$. 

Finally we need to check that $\mathfrak{h}_{x}(0) =
\mathfrak{t}_{x}(0)$, $\mathfrak{h}_{x}(1) =
\mathfrak{t}_{x}(1)$, $\mathfrak{h}_{x}(\infty) =
\mathfrak{t}_{x}(\infty)$.  Substituting in the formula for
$\mathfrak{h}_{x}$ we get 
\begin{equation} \label{eq:values.of.h}
\begin{aligned}
\mathfrak{h}_{x}(0) & =
\frac{(-f_{4}+p_{4})p_{5}f_{5}}{f_{4}(p_{5}-f_{5})p_{4}}, \\[0.5pc]
\mathfrak{h}_{x}(1) & = \frac{f_{5}}{f_{4}}, \\[0.5pc]
\mathfrak{h}_{x}(\infty) & = \frac{(-1+f_{4})f_{5}}{f_{4}(-1+f_{5})},
\end{aligned}
\end{equation}
\

\medskip

\noindent
To work out the images of $0$, $1$, and $\infty$ under the map
$\mathfrak{t}_{x}$ note that the pencil of curves $\{
\mathfrak{T}_{x,p}\}_{p \in C}$ has 5 singular fibers which are all
reducible. If we think of the curves $\mathfrak{T}_{x,p}$ as curves
embedded in the quadric $\mathbb{P}^{1}\times \mathbb{P}^{1}$, then
the five singular members of the pencil are $\mathfrak{T}_{x,0}$,
$\mathfrak{T}_{x,1}$, $\mathfrak{T}_{x,\infty}$,
$\mathfrak{T}_{x,p_{4}}$, $\mathfrak{T}_{x,p_{5}}$. Each of these
curves in $\mathbb{P}^{1}\times \mathbb{P}^{1}$ has two irreducible
components. The components can be described explicitly as follows:

\

\begin{itemize}
\item $\mathfrak{T}_{x,0} = \mathfrak{T}_{x,0}'\cup
\mathfrak{T}''_{x,0}$, where $\mathfrak{T}_{x,0}'$ is the unique
$(1,1)$ curve passing through the points $(f_{4},f_{5})$, $(0,0)$, and
$(p_{4},p_{5})$, and $\mathfrak{T}_{x,0}''$ is the unique $(1,1)$ curve
passing through $(f_{4},f_{5})$, $(1,1)$, and $(\infty,\infty)$;
\item $\mathfrak{T}_{x,1} = \mathfrak{T}_{x,1}'\cup
\mathfrak{T}_{x,1}''$, where $\mathfrak{T}_{x,1}'$ is the unique
$(1,1)$ curve passing through the points $(f_{4},f_{5})$, $(1,1)$, and
$(p_{4},p_{5})$, and $\mathfrak{T}_{x,1}''$ is the unique $(1,1)$ curve
passing through $(f_{4},f_{5})$, $(0,0)$, and $(\infty,\infty)$;
\item $\mathfrak{T}_{x,\infty} = \mathfrak{T}_{x,\infty}'\cup
\mathfrak{T}_{x,\infty}''$, where $\mathfrak{T}_{x,\infty}'$ is the
unique $(1,1)$ curve passing through the points $(f_{4},f_{5})$,
$(\infty,\infty)$, and $(p_{4},p_{5})$, and $\mathfrak{T}_{x,\infty}''$
is the unique $(1,1)$ curve passing through $(f_{4},f_{5})$, $(0,0)$,
and $(1,1)$;
\item $\mathfrak{T}_{x,p_{4}} = \mathfrak{T}_{x,p_{4}}'\cup
\mathfrak{T}_{x,p_{4}}''$, where $\mathfrak{T}_{x,p_{4}}'$ is the
unique $(0,1)$ curve passing through the point $(f_{4},f_{5})$, and
$\mathfrak{T}_{x,p_{4}}''$ is the unique $(2,1)$ curve passing through
$(f_{4},f_{5})$, $(0,0)$, and $(1,1)$, and $(\infty,\infty)$;
\item $\mathfrak{T}_{x,p_{5}} = \mathfrak{T}_{x,p_{5}}'\cup
\mathfrak{T}_{x,p_{5}}''$, where $\mathfrak{T}_{x,p_{5}}'$ is the
unique $(1,0)$ curve passing through the point $(f_{4},f_{5})$, and
$\mathfrak{T}_{x,p_{5}}''$ is the unique $(1,2)$ curve passing through
$(f_{4},f_{5})$, $(0,0)$, and $(1,1)$, and $(\infty,\infty)$;
\end{itemize}
Now computing the slopes of $T_{(0,0)} \mathfrak{T}'_{x,(0,0)}$,
$T_{(1,1)} \mathfrak{T}''_{x,(1,1)}$, and $T_{(\infty,\infty)}
\mathfrak{T}''_{x,(\infty,\infty)}$ we get precisely the three values
from \eqref{eq:values.of.h}. This shows that the values of
$\mathfrak{h}_{x}$ and $\mathfrak{t}_{x}$ at $0$, $1$, and $\infty$
match and completes the proof of the theorem. \hfill $\Box$

\

\bigskip

\section{The modular spectral cover} \label{sec:modular}

In this chapter we give a detailed description of the rational map from the
Hitchin fiber to the moduli of bundles. We show that a resolution of
the base locus of this map gives rise to a finite morphism and that this
realizes a blow-up of (a connected component of) the Hitchin fiber as a
degree four cover of the moduli space $X$ of parabolic bundles. We
also show that the blown-up Hitchin fiber is naturally embedded in the
bundle of logarithmic one forms on $X$ with poles along the union of
the sixteen $(-1)$-curves on $X$. Thus the blown-up Hitchin fiber is a
spectral cover for a meromorphic Higgs bundle on $X$. We will refer to
this cover as the {\em modular spectral cover}.  

\index{terms}{Hitchin!fiber}
\index{terms}{Hitchin!map}

\subsection{The fiber of the Hitchin map} \label{ssec:Hitchin.fiber}

In this section we establish some elementary facts about the geometry
of the parabolic Hitchin system in our setup.  In particular, we will
see that it induces a rational map from the Jacobian of a spectral
curve to our del Pezzo moduli space $X$, and that the degree of this
map is 4. In subsequent sections we will refine this by resolving this
map. \index{terms}{Jacobian} \index{terms}{Jacobian!of spectral curve}
\index{terms}{Hitchin!system}

For concreteness we will identify the del Pezzo surface $X$ with the
moduli space $N_{0}(0,1/2)$.  Consider also the moduli space
$\Higgs_{0}(0,1/2)$ of stable rank two degree zero parabolic Higgs
bundles on $(C,p_{1}, \ldots, p_{5})$ with balanced parabolic weights
$0 < 1/2$. (Note that we consider {\em all} stable parabolic $GL(2)$-Higgs bundles.
In particular we do not impose any restriction on the residues at the $p_i$.)
The space $\Higgs_{0}(0,1/2)$ has a number of remarkable
geometric properties that we will exploit:

\

\medskip

\noindent
{\bfseries (i)} \ $\Higgs_{0}(0,1/2)$ is a smooth connected
quasi-projective variety \cite{biswas-ramanan,yokogawa,boden-yokogawa}
of dimension 13.  (One way to count: the underlying bundle is,
generically, ${\mathcal{O}}_{C} \oplus {\mathcal{O}}_{C}$, so it has
no moduli. The Higgs field is a $2 \times 2$ matrix whose entries are
sections of $\Omega^{1}_{C}(\sum_{i=1}^{5} p_{i}) \cong
{\mathcal{O}}_{C}(3)$, so it depends on $4 \times 4 = 16$
parameters. A generic Higgs field is compatible with a finite number
(=2) of filtrations at each puncture, so still 16 parameters.  From
this we need to subtract the effect of auotomorphisms: $GL(2)$ acts on
the bundle, but with a 1-dimensional stabilizer, so 16-4+1=13. We will
see a different way to count in the next item.)  \

\smallskip

\noindent
{\bfseries (ii)} \ $\Higgs_{0}(0,1/2)$ is equipped with two natural
maps
\[
\xymatrix{
\Higgs_{0}(0,1/2) \ar[r]^-{\bh} \ar@{-->}[d]_-{\bpi} & \bB \\
X
}
\]
where
\begin{itemize}
\item $\bpi : \Higgs_{0}(0,1/2) \dashrightarrow X$ is the rational map
of forgetting the Higgs field, i.e. $\bpi(E,\theta) := E$ for a
parabolic Higgs bundle $(E,\theta)$;
\item $\bB := H^{0}(C,\Omega^{1}_{C}(\sum_{i=1}^{5} p_{i}))\oplus
  H^{0}(C,(\Omega^{1}_{C}(\sum_{i=1}^{5} p_{i}))^{\otimes 2}) \cong
  H^{0}(\mathcal{O}_{C}(3))\oplus  H^{0}(\mathcal{O}_{C}(6))$ is the
  Hitchin base (for G=GL(2)). It is a vector space of dimension $(3+1)+(6+1) = 11$;
\item $\bh : \Higgs_{0}(0,1/2) \to \bB$ is Hitchin's proper map,
  defined by $\bh(E,\theta) :=
  (-\op{tr}(\theta),\det(\theta))$. \index{terms}{Hitchin!base}
  \index{terms}{Hitchin!map}
\end{itemize}

\

\smallskip

\noindent
{\bfseries (iii)} \ The connected component of the general fiber
  of $\bh$ is a principally polarized abelian surface. More precisely
  \cite{hitchin}, \cite{eyal}, 
  \cite{donagi-msri}, \cite[chapter~4]{ron-eyal} 
  for a general point $\beta = (\beta_{1},\beta_{2}) = \bB \cong
  H^{0}(\mathcal{O}_{C}(3))\oplus  H^{0}(\mathcal{O}_{C}(6))$, each
  connected component of the Hitchin fiber $\bh^{-1}(\beta)$ is
   isomorphic to the Jacobian of a genus 2 curve $\sC_{\beta}$, 
  where $\bp_{\beta} : \sC_{\beta} \to C$ is the two sheeted
  spectral cover given by \index{terms}{Jacobian!of spectral curve}
  \index{terms}{Hitchin!fiber} 
\[
\sC_{\beta} \; : \; \lambda^{2} + \mathfrak{p}^{*}\beta_{1} \lambda +
\mathfrak{p}^{*}\beta_{2} = 0.   
\]
This equation describes $\sC_{\beta}$ as a divisor in the
surface $\op{tot}(\Omega_{C}^{1}(\sum_{i=1}^{5} p_{i})) \cong
\op{tot}(\mathcal{O}_{C}(3))$. The map $\mathfrak{p} :
\op{tot}(\mathcal{O}_{C}(3)) \to C$ is the natural projection, and 
$\lambda$ is the tautological section of
$\mathfrak{p}^{*}\mathcal{O}_{C}(3)$. 

More precisely we have the following 

\begin{lem} \label{lem:h.components} Let $\beta \in \bB$ be such that
  $\sC_{\beta}$ is smooth and such that $p_{1}$, \ldots $p_{5}$ are
  not among the six branch points of  $\bp_{\beta} : \sC_{\beta} \to
  C$. Then:
\begin{itemize}
\item[{\bfseries (a)}] The 32 connected components of $\bh^{-1}(\beta)$ are
  labelled by the quintuples of points $\tilde{p}_{1}$, \ldots, 
  $\tilde{p}_{5}$ in $\sC_{\beta}$, such that
  $\bp_{\beta}(\tilde{p}_{i}) = p_{i}$ for $i = 1, \ldots, 5$.
\item[{\bfseries (b)}] Each connected component of $\bh^{-1}(\beta)$
  is canonically isomorphic to $\Pic^{3}(\sC_{\beta})$. For every
  quintuple $\tilde{p}_{1}$, \ldots $\tilde{p}_{5}$ as in {\bfseries
  (a)}, the embedding $i_{\tilde{p}_{1},\ldots,\tilde{p}_{5}} :
  \Pic^{3}(\sC_{\beta}) \to \bh^{-1}(\beta)$ of
  $\Pic^{3}(\sC_{\beta})$ as the component labelled by
  $\tilde{p}_{1}$, \ldots $\tilde{p}_{5}$ associates to a line bundle
  $L \in \Pic^{3}(\sC_{\beta})$ a parabolic Higgs bundle
  $i_{\tilde{p}_{1},\ldots,\tilde{p}_{5}}(L) =
  ((V,\bF,(0,1/2)),\theta)$, where:
\begin{itemize}
\item[$\bullet$] $V = \bp_{\beta*}L$;
\item[$\bullet$] $\theta : V \to V\otimes \Omega^{1}_{C}(\sum_{i
  =1}^{5} p_{i})$ 
  is given by $\theta = \bp_{\beta*}((\bullet)\otimes \lambda)$;
\item[$\bullet$] $F_{i} = \ker\left[\xymatrix@1{V_{p_{i}} =
  (\bp_{\beta*}L)_{p_{i}}
  \ar[r]^-{\op{ev}_{\tilde{p}_{i}}}&  L_{\tilde{p}_{i}}}\right]$.
\end{itemize}
\end{itemize}
\end{lem}
{\bf Proof.} Every point in $\bh^{-1}(\beta)$ is a parabolic Higgs
bundle with spectral cover $\sC_{\beta}$. Note that the weights of the
parabolic structure are fixed and so do not constitute extra
data. Furthermore by definition the spectral cover associated to a
point $((V,\bF,(0,1/2)),\theta) \in
\Higgs_{0}(0,1/2)$ depends only on the meromorphic Higgs bundle
$(V,\theta)$ and is independent of the flags $\bF$.  The fact that the
meromorphic Higgs bundles underlying points in $\bh^{-1}(\beta)$ are
obtained from line bundles in $\Pic^{3}(\sC_{\beta})$ via the
prescription given in part {\bfseries (b)} follows from the general
spectral construction for Higgs bundles \cite{hitchin},
\cite{eyal}, \cite{ron-eyal}. But all possible quasi-parabolic
structures on a given meromorphic Higgs bundle $(V,\theta)$
correspond to choosing $F_{i} \subset V_{p_{i}}$ which are eigenlines
of $\op{Res}_{p_{i}} \theta$. Under our genericity assumption
\index{terms}{genericity condition}  the two
eigenlines of  $\op{Res}_{p_{i}} \theta$ are simply the fibers of
$L$ at the two preimages of $p_{i}$ in $\sC_{\beta}$. This shows that
all the parabolic structures on the meromorphic Higgs bundle 
$\bp_{\beta*} (L,(\bullet)\otimes \lambda)$ are labelled by quintuples 
of points $\tilde{p}_{1}$, \ldots $\tilde{p}_{5}$ such that 
$\bp_{\beta}(\tilde{p}_{i}) = p_{i}$. \hfill $\Box$

\

\smallskip

\noindent
{\bfseries (iv)} \ The natural structure \cite{eyal,bottacin} on the moduli space $\Higgs_{0}(0,1/2)$ of meromorphic Higgs bundles is that of an (algebraic or analytic) Poisson manifold. Such manifolds come with a natural foliation by symplectic leaves, holomorphic submanifolds on which the Poisson structure is non-degenerate, i.e. symplectic. In the case of the moduli of meromorphic Higgs bundles, Markman shows that this foliation by symplectic leaves is actually algebraic, specified by fixing the conjugacy classes (or more precisely,  the coadjoint orbits) of the residues. This picture extends verbatim to our parabolic setup:
The space $\Higgs_{0}(0,1/2)$ has a natural
  holomorphic Poisson structure and the restriction of the Hitchin map
  to a general symplectic leaf is an algebraically completely
  \index{terms}{Hitchin!map}
  \index{terms}{Hitchin!system}
  integrable system \cite{eyal,bottacin}. 

\

\smallskip

\noindent
{\bfseries (v)} \ Suppose $\beta \in \bB$ is a generic point and let 
$\tilde{p}_{1}, \ldots, \tilde{p}_{5} \in \sC_{\beta}$ be such that
$\bp_{\beta}(\tilde{p}_{i}) = p_{i}$. Let $\Par_{\, \sC_{\beta}} =
\tilde{p}_{1} + \ldots + \tilde{p}_{5}$ be the divisor given by
these five points. The composition 
\[
\bff_{\Par_{\, \sC_{\beta}}} := \bpi\circ
i_{\tilde{p}_{1},\ldots,\tilde{p}_{5}} \; : \;
\Pic^{3}(\sC_{\beta}) \dashrightarrow X
\] 
is a rational map. We have the following

\begin{lem} \label{lem:degree.four} The rational map
  $\bff_{\Par_{\, \sC_{\beta}}}$ is generically finite and of degree $4$.
\end{lem}
{\bfseries Proof.} Since $\Pic^{3}(\sC_{\beta})$ and $X$ are irreducible
surfaces, saying that the map  $\bff_{\Par_{\sC},\beta}$ is
generically finite is equivalent to saying that it is dominant.
To show that this map is dominant consider the open subset 
 $\Pic^{3}(\sC_{\beta})^{\op{stable}} \subset
\Pic^{3}(\sC_{\beta})$ parametrizing line
bundles whose pushforwards to $C$ are stable parabolic bundles. 
The map $\bff_{\Par_{\, \sC_{\beta}}} :
\Pic^{3}(\sC_{\beta})^{\op{stable}} \to X$ is a morphism of
irreducible quasi-projective surfaces and so by the fiber dimension
theorem we know that this map will be dominant if and only if it has a
zero dimensional fiber. 

Let $(V,\bF,(0,1/2)) \in X$ be a very stable parabolic bundle, i.e. a
bundle that does not admit any non-zero nilpotent Higgs fields whose
residues preserve the parabolic structure. In
Section~\ref{sssec:wobbly} we will show that such bundles do exist.
Let $\Pi \subset \Higgs_{0}(0,1/2)$ be the vector space of all
parabolic Higgs bundles with underlying parabolic bundle
$(V,\bF,(0,1/2)) \in X$. 
It is contained in the 16-dimensional $H^0(C, End(V) \otimes \mathcal{O}_{C}(3))$
as the codimension 5 subspace of $ \mathcal{O}_{C}(3)$-valued endomorphisms 
preserving $\bF$, so it is an 11 dimensional vector space.
If we restrict the Hitchin map to $\Pi$ we get a morphism $g : \Pi \to \bB$
of $11$-dimensional affine spaces. This morphism is given by $4$
linear homogeneous polynomials and $7$ quadratic homogeneous
polynomials and so descends to a rational map $\bar{g} : \mathbb{P}(\Pi)
\to \mathbb{WP}(\bB)$ between the projectivization of $\Pi$ and the
weighted projectivization of $\bB$ where $H^{0}(\mathcal{O}_{C}(3))$
has weight $1$ and $H^{0}(\mathcal{O}_{C}(6))$ has weight $2$. However
the fact that $(V,\bF,(0,1/2))$ was very stable implies that the
$g$-preimage of $0 \in \bB$ consists of the single point $\bar{0} \in \Pi$. By
the fiber dimension theorem this implies that $g$ is
dominant. Furthermore this shows that $\bar{g}$ is a morphism, which
is necessarily dominant, and hence surjective. This shows that
$\bar{g}$ is a generically finite morphism.
But since the source $\mathbb{P}(\Pi)$ is a projective space, 
so in particular it has Picard number 1, 
it follows that no  curve  in it can be contracted unless all curves are contracted,
so any generically finite morphism must be finite.
Hence $g$ too is  finite and surjective. 
This implies that the intersection 
$\Pi \cap i_{\tilde{p}_{1},\ldots,\tilde{p}_{5}}(\Pic^{3}(\sC_{\beta}))$ 
is a non-empty finite set. Therefore $\bff_{\Par_{\, \sC_{\beta}}}$ is
dominant and hence generically finite.
\index{terms}{Hitchin!fiber}

The fact that $\bff_{\Par_{\, \sC_{\beta}}}$ is generically finite implies
that the restriction of $\bpi$ to $\bh^{-1}(\beta)$ is also
generically finite. But the degree of the finite part of the map
$\bpi_{|\bh^{-1}(\beta)}$ will be equal to the number of points in the
intersection $\Pi\cap \bh^{-1}(\beta)$, where $\Pi$ is the space of
parabolic Higgs fields on a generic very stable parabolic bundle. But
by the same token the intersection of $\Pi\cap \bh^{-1}(\beta)$ is
also equal to the degree of the finite map $g : \Pi \to \bB$. Since
the latter is given by $4$ linear polynomials and $7$ quadratic ones
it follows that $g$ has degree $1^{4}\cdot 2^{7} = 2^{7}$. This
implies that the degree of  $\bpi_{|\bh^{-1}(\beta)}$ is
$2^{7}$. Since by Lemma~\ref{lem:h.components} $\bh^{-1}(\beta)$ has
$2^{5}$ connected 
components it follows that the restriction $\bff_{\Par_{\, \sC_{\beta}}}$
of $\bpi$ to a connected component of $\bh^{-1}(\beta)$ will have
degree $4$. The lemma is proven. \hfill $\Box$

\subsection{The base locus}  \label{ssec:base.locus}

In this section we will identify the base locus of the rational map
\[
\bff_{\Par_{\, \sC_{\beta}}} \; : \; \Pic^{3}(\sC_{\beta})
\dashrightarrow X.
\] 
with a set 
$\left\{ \mathfrak{P}_{I} \right\}_{I  \in \oddL}$
of $16$ points which we describe below.
In Theorem \ref{theo:map.from.jacobian} we will see that the rational map $\bff_{\Par_{\, \sC_{\beta}}}$
lifts to a morphism
\[
f_{\beta} \; : \; Y_{\beta} \to X
\subset \mathbb{P}^{4},
\]
where $Y_{\beta}$ is the blowup of  $\Pic^{3}(\sC_{\beta})$ at the  $16$ points of 
$\left\{ \mathfrak{P}_{I} \right\}_{I  \in \oddL}$.

For the duration of the section we will
fix $\beta$ and the divisor $\Par_{\sC}$. 

\begin{defi} \label{defi:16.points}  For every subset $I \subset \{
  1,2,3,4,5\}$ of odd cardinality set $\kappa(I) = (\# I - 3)/2 \in
  \{-1,0,1\}$ and define a degree three line bundle $\mathfrak{P}_{I}$ on
  $\sC_{\beta}$ by setting
\[
\mathfrak{P}_{I} = K_{\sC_{\beta}}^{\otimes - \kappa(I)}\otimes
\mathcal{O}_{\sC_{\beta}}\left( \sum_{i \in I} \tilde{p}_{i} \right)
\]
\end{defi}

\

\noindent
The collection of line bundles $\left\{ \mathfrak{P}_{I} \right\}_{I
  \in \oddL}$, where $\oddL \subset 2^{\{ 1,2,3,4,5\}}$ is the set of
all subsets of $\{ 1,2,3,4,5\}$ of odd cardinality, can be viewed as a
collection of $16$ distinct points in $\Pic^{3}(\sC_{\beta})$. Let
\[
\bmu_{\beta} : Y_{\beta} := \op{Bl}_{\{ \mathfrak{P}_{I}\}}(
\Pic^{3}(\sC_{\beta})) \to \Pic^{3}(\sC_{\beta})
\]
denote the blow-up of $\Pic^{3}(\sC_{\beta})$ at those $16$
points, and let $E_{I} \subset Y_{\beta}$ be the exceptional
divisor corresponding to the point $\mathfrak{P}_{I}$.

Consider the Abel-Jacobi map $\aj : \sC_{\beta} \to
\Pic^{3}(\sC_{\beta})$, given by $\aj(x) := K_{\sC_{\beta}}(x)$,
and the associated theta line bundle
$\mathcal{O}_{\Pic^{3}(\sC_{\beta})}\left(\aj(\sC_{\beta})\right)$
on $\Pic^{3}(\sC_{\beta})$. To simplify notation we will denote
the pull-back line bundle $\bmu_{\beta}^{*}
\mathcal{O}_{\Pic^{3}(\sC_{\beta})}\left(\aj(\sC_{\beta})\right)$
by $\btheta$, and more generally for every $\mathfrak{A} \in
\Pic^{0}(\sC_{\beta})$ we will denote by $\btheta_{\mathfrak{A}}$
the pull-back $\bmu_{\beta}^{*}t^{*}_{\mathfrak{A}}
\mathcal{O}_{\Pic^{3}(\sC_{\beta})}\left(\aj(\sC_{\beta})\right)$
of the $\mathfrak{A}$-translated theta line bundle.

With this notation we have the following 

\

\begin{theo} \label{theo:map.from.jacobian}
The rational map $\bff_{\Par_{\, \sC_{\beta}}} \; : \;
\Pic^{3}(\sC_{\beta}) 
\dashrightarrow X$ extends to a morphism
\[
f_{\beta} \; : \; Y_{\beta} \to X
\subset \mathbb{P}^{4}
\]
satisfying $f_{\beta}^{*} K_{X}^{-1} = 
\btheta_{\mathfrak{A}}
\otimes \btheta^{\otimes 3} \left(-\sum_{I \in
  \oddL} E_{I} \right)$,
where
$\mathfrak{A} := K_{\sC_{\beta}}^{\otimes 5}\left( -2
\Par_{\sC}\right)$. 
The morphism $f_{\beta} \; : \;
Y_{\beta} \to X$ is finite of degree four.
\end{theo}
{\bfseries Proof.} Recall that the $dP_{5}$ del Pezzo $X$ contains
$16$ lines. If we use the model in which $X$ is identified with the
blow-up of $\mathbb{P}^{2} = S^{2}C$ at the five points $p_{1}$,
\ldots, $p_{5}$ on the conic $C \subset S^{2}C$, then the $16$ lines
are the five exceptional divisors, the strict transforms of the ten
lines through two of the points, and the strict transform of the conic
through all five points. In particular, the $16$ can be labelled
naturally by the subsets $I \in \oddL$: if $I = \{i\}$ we will
set $L_{I} \subset X$ to be the exceptional divisor corresponding to
$p_{i}$; if $\# I = 3$, we will set $L_{I}$ to be the strict transform
of the line
through the two points $\{ p_{i} \}_{i \notin I}$; if $I = \{
1,2,3,4,5 \}$, then we will set $L_{I}$ to be the  strict transform
of the conic through all five points. 

\begin{figure}[!ht]
\begin{center}
\begin{minipage}[c]{2.8in}
\begin{center}
\psfrag{p1}[c][c][1][0]{{$p_{1}$}}
\psfrag{p2}[c][c][1][0]{{$p_{2}$}}
\psfrag{p3}[c][c][1][0]{{$p_{3}$}}
\psfrag{p4}[c][c][1][0]{{$p_{4}$}}
\psfrag{p5}[c][c][1][0]{{$p_{5}$}}
\epsfig{file=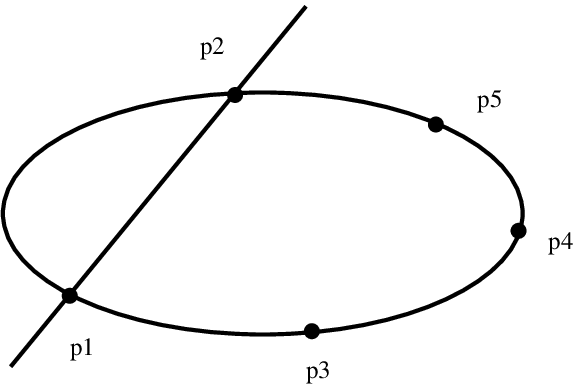,width=2.5in} 
\end{center}
\caption{The cubic giving the quadrilateral  $L_{\{1,2,3,4,5\}}\cup
  L_{\{1\}} \cup L_{\{2\}} \cup 
  L_{\{3,4,5\}}$}\label{fig:qla}  
\end{minipage}
\hspace{0.4in}
\begin{minipage}[c]{2.8in}
\begin{center}
\psfrag{p1}[c][c][1][0]{{$p_{1}$}}
\psfrag{p2}[c][c][1][0]{{$p_{2}$}}
\psfrag{p3}[c][c][1][0]{{$p_{3}$}}
\psfrag{p4}[c][c][1][0]{{$p_{4}$}}
\psfrag{p5}[c][c][1][0]{{$p_{5}$}}
\epsfig{file=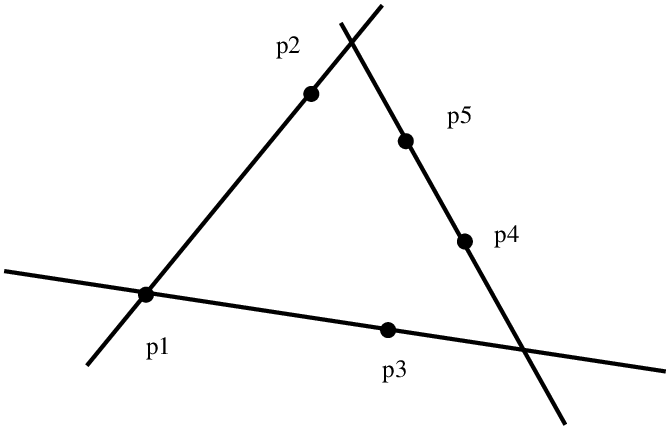,width=2.5in} 
\end{center}
\caption{The cubic giving the quadrilateral  $L_{\{1\}}\cup
  L_{\{3,4,5\}} \cup L_{\{2,4,5\}} \cup 
  L_{\{1,2,3\}}$}\label{fig:qlb}
\end{minipage}  
\end{center}
\end{figure}

There are $40$ hyperplanes in $\mathbb{P}^{4}$ that inersect $X$ in
quadrilaterals, each consisting of four of the lines. We can describe
these quadrilaterals as follows:
\begin{itemize}
\item $L_{\{1,2,3,4,5\}}\cup L_{\{1\}} \cup L_{\{2\}} \cup
  L_{\{3,4,5\}}$ and   its orbit under the action of the symmetric group 
  $S_{5}$. There are $10$ quadrilaterals of this type.
\item $L_{\{1\}}\cup L_{\{3,4,5\}}\cup L_{\{2,4,5\}} \cup
  L_{\{1,2,3\}}$ and  its orbit under the action of the symmetric group 
  $S_{5}$. There are 30 curves of this type.
\end{itemize}
In the model $X = \op{Bl}_{\{p_{1}, \ldots, p_{5}\}}(\mathbb{P}^{2})$
these quadrilaterals correspond to special reducible cubics in
$\mathbb{P}^{2}$ passing through the five points $p_{1}, \ldots,
p_{5}$. The quadrilateral $L_{\{1,2,3,4,5\}}\cup L_{\{1\}} \cup
L_{\{2\}} \cup L_{\{3,4,5\}}$ corresponds to the cubic (see
Figure~\ref{fig:qla}) which is the union of the conic through the five
points and the line through $p_{1}$ and $p_{2}$, whereas the
quadrilateral $L_{\{1\}}\cup L_{\{3,4,5\}}\cup L_{\{2,4,5\}} \cup
L_{\{1,2,3\}}$ corresponds to the triangle (see Figure~\ref{fig:qlb})
with sides the line through $p_{1}$ and $p_{2}$, the line through
$p_{1}$ and $p_{3}$, and the line through $p_{4}$ and $p_{5}$.

\

\medskip

Next given $I \in \oddL$  
consider the translated theta divisor 
\[
\bT_{I} := t_{K_{\sC_{\beta}}^{\otimes -1}(\Par_{\sC})\otimes
\mathfrak{P}_{I}^{\otimes -1}}^{*} (\aj(\sC_{\beta})) \subset
\Pic^{3}(\sC_{\beta}), 
\] 
and its strict transform $\bD_{I} \subset Y_{\beta}$ in the
blown-up surface $Y_{\beta}$.

In particular each $\bD_{I}$ is isomorphic to $\sC_{\beta}$. Consider
next the rational map $\bff_{\Par_{\, \sC_{\beta}}} :
\Pic^{3}(\sC_{\beta}) \dashrightarrow X$.  Note that
$\bff_{\Par_{\, \sC_{\beta}}}$ maps the translate $\bT_{I}$
dominantly onto the line $L_{I}$. This follows immediately from the
description of the parabolic bundles in $L_{I}$ given in the proof of
Proposition~\ref{prop:dP5}. In particular the closure of
${(\bff_{\Par_{\, \sC_{\beta}}})}^{-1}(L_{I})$ in
$\Pic^{3}(\sC_{\beta})$ is an effective divisor whose support
contains $\bT_{I}$. 

Now we have the following 

\begin{lem} \label{lem:T_I} The closure in
  $\Pic^{3}(\sC_{\beta})$ of the preimage
  ${(\bff_{\Par_{\, \sC_{\beta}}})}^{-1}(L_{I})$ is equal to 
$\bT_{I}$. 
%and the differential of $\bff_{\Par_{\, \sC_{\beta}}}$ at the generic point of $\bT_{I}$ is an isomorphism.
\end{lem}
{\bfseries Proof.} The statement is invariant under $S_{5}$ and so it
suffices to check it for $I = \{1,2,3,4,5\}, \linebreak \{3,4,5\},
\{5\}$. We consider each of these cases separately: 

\

\noindent
Suppose $I = \{1,2,3,4,5\}$. In this case the curve
$\bT_{\{1,2,3,4,5,\}}$ parametrizes all degree three line bundles on
$\sC_{\beta}$ of the form $K_{\sC_{\beta}}(x)$, with $x \in
\sC_{\beta}$. On the other hand from the proof of
Proposition~\ref{prop:dP5} we have that $L_{\{1,2,3,4,5\}} \subset X$
is the line parametrizing all parabolic bundles $(V,\bF,(0,1/2)) \in
X$ for which $V \cong \mathcal{O}_{C}(-1)\oplus
\mathcal{O}_{C}(1)$. By definition
$\bff_{\Par_{\, \sC_{\beta}}}(L)$ is a parabolic Higgs bundle whose
underlying vector bundle is $\bp_{\beta*} L$. Thus $L \in
\Pic^{3}(\sC_{\beta})$ will be in the
$\bff_{\Par_{\, \sC_{\beta}}}$-preimage of $L_{\{1,2,3,4,5\}}$ if
and only if $\bp_{\beta*} L \cong \mathcal{O}_{C}(-1)\oplus
\mathcal{O}_{C}(1)$. But if $L$ is a line bundle of degree $3$ on
$\sC_{\beta}$, then $\bp_{\beta *} L$ will be isomorphic to
$\mathcal{O}_{C}(-1)\oplus \mathcal{O}_{C}(1)$ if and only if
$\mathcal{O}_{C}(-1)\otimes \bp_{\beta *} L$ has a unique section. By
the projection formula this is equivalent to
$\bp_{\beta}^{*}\mathcal{O}_{C}(-1)\otimes L \cong
K_{\sC_{\beta}}^{-1}\otimes L$ having a unique section, i.e it is
equivalent to $L \in \bT_{\{1,2,3,4,5\}}$.

\

\noindent
Suppose $I = \{3,4,5\}$. In this case $\bT_{\{3,4,5\}}$ parametrizes
all line bundles on $\sC_{\beta}$ of the form
$\mathcal{O}_{\sC_{\beta}}(\tilde{p}_{1}+\tilde{p}_{2} + x)$ for $x
\in \sC_{\beta}$. At the same time from the proof of
Proposition~\ref{prop:dP5} we have that the generic point of
$L_{\{3,4,5\}} \subset X$ is a parabolic bundles $(V,\bF,(0,1/2)) \in
X$ for which $V \cong \mathcal{O}_{C}\oplus \mathcal{O}_{C}$ and $\bF$
corresponds to a five-tuple of points $\{(p_{i},f_{i})\}_{i = 1}^{5}$,
on the quadric $\mathbb{P}(V)$ for which $(p_{1},f_{1})$ and
$(p_{2},f_{2})$ lie on a $(0,1)$ ruling. From this description we see
that if $L \in \Pic^{3}(\sC_{\beta})$, then
$\bff_{\Par_{\, \sC_{\beta}}}(L) \in L_{\{3,4,5\}}$ if and only if
we can find a trivialization of $\bp_{\beta *}L$ so that if we
identify the fibers $(\bp_{\beta *}L)_{p_{1}}$ and $(\bp_{\beta
*}L)_{p_{2}}$ via this trivialization, then $L_{\tilde{p}_{1}}$ and
$L_{\tilde{p}_{2}}$ become the same line. Explicitly a trivialization
of $\bp_{\beta *}L$ is given by choosing two global sections $s', s''
\in H^{0}(\sC_{\beta},L)$ of $L$, and the condition on the two flags
is equivalent to saying that the meromorphic function $s'/s''$ has the
same values at the points $\tilde{p}_{1}$ and $\tilde{p}_{2}$. This is
equivalent to saying that a linear combination of $s'$ and $s''$
vanishes at both $\tilde{p}_{1}$ and $\tilde{p}_{2}$. In other words
$\bff_{\Par_{\, \sC_{\beta}}}(L) \in L_{\{3,4,5\}}$ if and only if 
there is a holomorphic section in $L$ that vanishes at $\tilde{p}_{1}$
and $\tilde{p}_{2}$, i.e. $L \in \bT_{\{3,4,5\}}$.

\

\noindent
Suppose $I = \{5\}$. The curve $\bT_{\{5\}}$ parametrizes line bundles
of the form $K_{\sC_{\beta}}^{-1}\otimes
\mathcal{O}_{\sC_{\beta}}(\tilde{p}_{1} + \tilde{p}_{2} +
\tilde{p}_{3} + \tilde{p}_{4} + x)$, with $x \in \sC_{\beta})$. Again
from the proof of Proposition~\ref{prop:dP5} we have that the generic
point of $L_{\{5\}} \subset X$ is a parabolic bundles $(V,\bF,(0,1/2))
\in X$ for which $V \cong \mathcal{O}_{C}\oplus \mathcal{O}_{C}$ and
$\bF$ corresponds to a five-tuple of points $\{(p_{i},f_{i})\}_{i =
1}^{5}$, on $\mathbb{P}(V)$ such that $(p_{1},f_{1})$,
$(p_{2},f_{2})$, $(p_{3},f_{3})$, $(p_{4},f_{4})$ lie on some curve of
bidegree $(1,1)$. From this description we see that if $L \in
\Pic^{3}(\sC_{\beta})$, then $\bff_{\Par_{\, \sC_{\beta}}}(L)
\in L_{\{5\}}$ if and only if:
\begin{itemize}
\item  we can find two sections $s'$ and $s''$
of $(\bp_{\beta *}L)\otimes \mathcal{O}_{C}(1)$ so that $s'$ and $s''$
are linearly independent at the general point of $C$, and
\item if we use $s'$ and $s''$ to trivialize $\bp_{\beta
*}L$ generically, then in this trivialization the four lines
  $L_{\tilde{p}_{1}}$, $L_{\tilde{p}_{2}}$, $L_{\tilde{p}_{3}}$,
  $L_{\tilde{p}_{4}}$ coincide.
\end{itemize}
Rewriting this as a condition on $\sC_{\beta}$ we get that we must
have two holomorphic sections $s', s'' \in
H^{0}(\sC_{\beta},\bp_{\beta}^{*} \mathcal{O}_{C}(1) \otimes L)$, such
that the meromorphic function $s'/s''$ on $\sC_{\beta}$ takes the same
values at the four points $\tilde{p}_{1}$, $\tilde{p}_{2}$,
$\tilde{p}_{3}$, $\tilde{p}_{4}$. Similarly to the previous case this
is equivalent to saying that we have a section of $\bp_{\beta}^{*}
\mathcal{O}_{C}(1) \otimes L \cong K_{\sC_{\beta}}\otimes L$ which
vanishes at the four points $\tilde{p}_{1}$, $\tilde{p}_{2}$,
$\tilde{p}_{3}$, $\tilde{p}_{4}$, i.e. this is equivalent to having 
$L \in \bT_{\{5\}}$. 
In fact, we see easily that the differential of
$\bff_{\Par_{\, \sC_{\beta}}}$  at
the generic point of $\bT_{I}$ is an isomorphism,
so the preimage equals $T_I$  scheme-theoretically.
The lemma is proven. \ \hfill $\Box$

\

\smallskip

%\begin{lem} \label{lem:diffT_I} The differential of
%  $\bff_{\Par_{\, \sC_{\beta}}}$  at
%  the generic point of $\bT_{I}$ is an isomorphism.
%\end{lem}
%{\bfseries Proof.} Let $L \in \Pic^{3}(\sC_{\beta})$ be a point at which
%$\bff_{\Par_{\, \sC_{\beta}}}$ is well defined. This means that the
%data $(V,\bF,(0,1/2))$ defined by $L$ via the procedure in
%Lemma~\ref{lem:h.components} is a stable parabolic bundle. The
%differential of $\bff_{\Par_{\, \sC_{\beta}}}$ at $L$ is a linear
%map
%\[
%d_{L}\bff_{\Par_{\, \sC_{\beta}}} : T_{L}
%\Pic^{3}(\sC_{\beta}) \to T_{(V,\bF,(0,1/2))}X.
%\]
%We will compute this map explictly by deformation theory. (A similar
%computation in the non-parabolic case was carried out in
%\cite{teleman}.) First note that by Kodaira-Spencer theory we have 

Next consider a hyperplane $H \cong \mathbb{P}^{3} \subset
\mathbb{P}^{4}$ which cuts $X$ in one of the $40$ quadrilaterals. The
closure of the preimage $(\bff_{\Par_{\, \sC_{\beta}}})^{-1}(H\cap X)$ is a
divisor on $\Pic^{3}(\sC_{\beta})$ and for a generic $\beta$
the Neron-Severi class of this divisor will be a multiple of the
principal polarization $\ell \in \op{NS}(\Pic^{3}(\sC_{\beta}))$. 
Since each $\bT_{I}$ has Neron-Severi class $\ell$, it follows that 
the closure of \linebreak $\bff_{\Par_{\, \sC_{\beta}}}^{-1}(H\cap X)$  will
have a Neron-Severi class $a\cdot \ell$ with $a \geq 4$.

Since by
Lemma~\ref{lem:degree.four} the map $\bff_{\Par_{\, \sC_{\beta}}}$
is generically finite of degree four it follows that the Neron-Severi
class of this divisor is four times the principal polarization on
$\Pic^{3}(\sC_{\beta})$. But as we just noticed, this divisor
already contains four translates of the theta divisor, and so the
closure of $\bff_{\Par_{\, \sC_{\beta}}}^{-1}(\mathbb{P}^{3})$ must
be the sum of all $t_{K_{\sC_{\beta}}^{\otimes
    -1}(\Par_{\sC})\otimes \mathfrak{P}_{I}^{\otimes -1}}^{*}
(\aj(\sC_{\beta}))$ where $I$ ranges through the four sides of the
quadrilateral corresponding to $\mathbb{P}^{3}$. The lemma is
proven. \ \hfill $\Box$

\

\noindent
{\bf Conclusion of the proof of Theorem~\ref{theo:map.from.jacobian}:}
Combining the previous two lemmas, we see that for each $I$, the
inverse image divisor ${(\bff_{\Par_{\, \sC_{\beta}}})}^{-1}(L_{I})$
is equal to $\bT_{I}$.  Pulling back to $Y_{\beta}$, we get that the
inverse image divisor $(f_{\beta})^{-1}(L_{I})$, defined modulo the
exceptional divisors $E_I$, is equal to $D_{I}$.  The inverse image of
$K_X$ can be calculated using any of our 40 quadrilaterals, and
equals, again modulo the $E_I$: $f_{\beta}^{*} K_{X}^{-1} =
\btheta_{\mathfrak{A}} \otimes \btheta^{\otimes 3}.$ So the line
bundle on $Y_{\beta}$ giving the (a priori still rational) map to
$\mathbb{P}^{4}$ must be of the form $\btheta_{\mathfrak{A}} \otimes
\btheta^{\otimes 3} \left(-\sum_{I \in \oddL} a_I E_{I} \right)$ for
some integers $a_I$. By symmetry, all the $a_I$ must be equal to some
$a$.  However, the map to $\mathbb{P}^{4}$ has degree $4 \times 4 =
16$, and this forces $a=1$ as claimed.  In particular, the rational
map on the Picard lifts to a morphism on $Y_{\beta}$, concluding the
proof of the theorem.

\ \hfill $\Box$

\smallskip

\noindent

\

\medskip

\subsection{The case of nilpotent residues} 
\label{ssec:nilpotent.residues}

\

\smallskip

\noindent

In this section we specialize the general picture of the modular
spectral cover and the abelianized parabolic Hecke correspondence to
the main case of interest for us - the case of parabolic Higgs bundles
with nilpotent residues.  (Elsewhere in the literature, these are sometimes 
called {\em strictly parabolic} Higgs bundles.)
Concretely let us write $\Higgs_{\op{nilp}}$
for the moduli space of stable rank two degree zero parabolic Higgs
bundles on $(C,p_{1},\ldots,p_{5})$ with balanced parabolic weights $0
< 1/2$ and nilpotent residues at all five parabolic points $p_{1},
\ldots, p_{5}$. The subvariety $\Higgs_{\op{nilp}} \subset
\Higgs_{0}(0,1/2)$ is a closure of a maximal dimensional symplectic
leaf inside the Poisson variety $\Higgs_{0}(0,1/2)$, namely it is the
closure of the moduli space of Higgs bundles $(E,\theta)$ such that
for all $i = 1, \ldots, 5$ the matrix $\op{res}_{p_{i}} \theta$ is a
single $2\times 2$ nilpotent Jordan block. Again the space
$\Higgs_{\op{nilp}}$ has several useful properties which follow in
straightforward way from the definition:

\

\medskip

\noindent
{\bfseries (i)} \ $\Higgs_{\op{nilp}}$ is a smooth
connected  quasi-projective variety of dimension $4$.

\

\smallskip

\noindent
{\bfseries (ii)} \ $\Higgs_{\op{nilp}}$ is equipped with two natural
maps
\[
\xymatrix{
\Higgs_{\op{nilp}} \ar[r]^-{\bh} \ar@{-->}[d]_-{\bpi} & \bB_{\op{nilp}} \\
X
}
\]
where
\begin{itemize}
\item $\bpi : \Higgs_{\op{nilp}} \dashrightarrow X$ is the rational map
of forgetting the Higgs field, i.e. $\bpi(E,\theta) := E$ for a
parabolic Higgs bundle $(E,\theta)$;
\item $\bB_{\op{nilp}}$ is the corresponding 
  Hitchin base, i.e. \index{terms}{Hitchin!base} 
\[
\begin{aligned}
\bB_{\op{nilp}} & := 
{\textstyle H^{0}\left(C,\Omega^{1}_{C}\left(\sum_{i=1}^{5} 
p_{i}\right)\right)_{\op{nilp}}\oplus
  H^{0}\left(C,\left(\Omega^{1}_{C}\left(\sum_{i=1}^{5} p_{i}\right)\right)^{\otimes
    2}\right)_{\op{nilp}}} \\[+0.5pc]
& =
  {\textstyle  H^{0}\left(\mathcal{O}_{C}(3)\left(-\sum_{i=1}^{5} p_{i}\right)\right)\oplus
  H^{0}\left(\mathcal{O}_{C}(6)\left(-\sum_{i=1}^{5} p_{i}\right)\right)} \\[+0.5pc]
& \cong 0 \oplus H^{0}(\mathcal{O}_{C}(1)) \cong {\bf C}^2;
\end{aligned}
\] 
\item $\bh : \Higgs_{\op{nilp}} \to \bB_{\op{nilp}} $ is the
  restriction of the Hitchin map to  $\Higgs_{\op{nilp}}$,
  which is the full inverse image of $\bB_{\op{nilp}}
  $. \index{terms}{Hitchin!map}
  \index{terms}{Hitchin!base}
\end{itemize}

\

\smallskip

\noindent
{\bfseries (iii)} \ In contrast with the case of a generic symplectic
leaf, the fibers of $\bh : \Higgs_{\op{nilp}} \to \bB_{\op{nilp}}$ are
all connected. 
In the general case, the multiple components arise from 
labeling the two points of the spectral cover above each $p_i$
(or equivalently, the two eigenvalues of the residue).
In the nilpotent case, the spectral cover is branched at the $p_i$
and the eigenvalues all vanish.
Again the generic fiber is a principally polarized
abelian surface. More precisely if  $\beta = (0,\beta_{2}) \in  \bB_{\op{nilp}} =
  0\oplus   H^{0}(\mathcal{O}_{C}(6)(- \sum_{i =1}^{5} p_{i}))$
  is such that $\op{div}(\beta_{2}) = p_{1} + \cdots + p_{5} + p_{6}$
  with $p_{6}$ not equal to one of the five points $p_{1}, \ldots,
  p_{5}$, then the Hitchin fiber $\bh^{-1}(\beta)$ is
   isomorphic to the Jacobian of a curve $\sC_{\beta}$, 
  where $\bp_{\beta} : \sC_{\beta} \to C$ is the smooth two sheeted
  spectral cover given by \index{terms}{Jacobian!of spectral curve}
  \index{terms}{Hitchin!fiber}
\[
\sC_{\beta} \; : \; \lambda^{2} +
\mathfrak{p}^{*}\beta_{2} = 0.   
\]
This again describes $\sC_{\beta}$ as a divisor in the
surface $\op{tot}(\Omega_{C}^{1}(\sum_{i=1}^{5} p_{i})) \cong
\op{tot}(\mathcal{O}_{C}(3))$. The map $\mathfrak{p} :
\op{tot}(\mathcal{O}_{C}(3)) \to C$ is the natural projection, and 
$\lambda$ is the tautological section of
$\mathfrak{p}^{*}\mathcal{O}_{C}(3)$.

Note that by definition the spectral cover $\sC_{\beta}$ is the
hyperelliptic curve of genus two branched over the points $p_{1}$,
\ldots, $p_{6}$. In this case we have a unique choice of a preimage
of each of the parabolic points $p_{i}$, namely the corresponding
ramification points of the cover $\sC_{\beta} \to C$. This natural
choice is the reason why the fibers of the Hitchin map are connected in
this setting.

The Stein factorization  ${\bh}' : \Higgs \to {\bB}'$
of the Hitchin fibration  $\bh : \Higgs \to \bB$
gives a 32-sheeted branched cover ${\bB}' \to \bB$
and an abelian surface fibration over ${\bB}'$. Over our locus $\bB_{\op{nilp}}$, 
the 32-sheeted cover is totally ramified. The corresponding Hitchin fibers over ${\bB}$
are non-reduced, but the fibers of $\bh : \Higgs_{\op{nilp}} \to \bB_{\op{nilp}}$ are reduced. 
They coincide with the fibers of ${\bh}'$ and, as we saw above, 
with the Jacobians of the spectral cover $\sC_{\beta}$.

For ease of reference let us write $\tilde{p}_{1}, \ldots,
\tilde{p}_{6} \in \sC_{\beta}$ for the six ramification points. The
Hitchin fiber  \index{terms}{Hitchin!fiber}
 $\bh^{-1}(\beta)$ then 
  is canonically isomorphic to $\Pic^{3}(\sC_{\beta})$. The
  isomorphism  $i_{\tilde{p}_{1},\ldots,\tilde{p}_{5}} :
  \Pic^{3}(\sC_{\beta}) \to \bh^{-1}(\beta)$ associates to a line bundle
  $L \in \Pic^{3}(\sC_{\beta})$ a parabolic Higgs bundle
  $i_{\tilde{p}_{1},\ldots,\tilde{p}_{5}}(L) =
  ((V,\bF,(0,1/2)),\theta)$, where:
\begin{itemize}
\item[$\bullet$] $V = \bp_{\beta*}L$;
\item[$\bullet$] $\theta : V \to V\otimes \Omega^{1}_{C}(\sum_{i
  =1}^{5} p_{i})$ 
  is given by $\theta = \bp_{\beta*}((\bullet)\otimes \lambda)$;
\item[$\bullet$] $F_{i} = \ker\left[\xymatrix@1{V_{p_{i}} =
  (\bp_{\beta*}L)_{p_{i}}
  \ar[r]^-{\op{ev}_{\tilde{p}_{i}}}&  L_{\tilde{p}_{i}}}\right]$
\end{itemize}

\noindent
Again the rational map 
\[
\bff_{\Par_{\, \sC_{\beta}}} := \pi\circ i_{\tilde{p}_{1},\ldots,\tilde{p}_{5}}
: \Pic^{3}(\sC_{\beta}) \dashrightarrow X
\]
can be  resolved by blowing up the sixteen points
$\{\mathfrak{P}_{I}\}_{I \in \oddL}$ in
$\Pic^{3}(\sC_{\beta})$.

Because in this case the points $\tilde{p}_{i}$ are all Weierstrass
points on the genus two curve $\sC_{\beta}$, Definition~\ref{defi:16.points} shows
that the points $\{\mathfrak{P}_{I}\}_{I \in \oddL}$ in
$\Pic^{3}(\sC_{\beta})$ are simply all points of the form
$K_{\sC_{\beta}}\otimes \kappa$ where $\kappa$ runs over the 16 theta
characteristics of $\sC_{\beta}$.

For convenience, let us use the sixth ramification point
$\tilde{p}_{6}$ to identify $\Pic^{3}(\sC_{\beta})$ with the
Jacobian $\Jac_{\beta} := \Pic^{0}(\sC_{\beta})$ of $\sC_{\beta}$ via
the map $\Pic^{3}(\sC_{\beta}) \to \Pic^{0}(\sC_{\beta})$, $L
\mapsto L(-3\tilde{p}_{6})$. Under this identification the 16 points
$\{\mathfrak{P}_{I}\}_{I \in \oddL}$  get mapped to the 16  $2$-torsion points
$\{p_{I}\}_{I \in \oddL}$ in the abelian surface $\Jac_{\beta}$. Again we will
write \index{terms}{Jacobian!of spectral curve}
\[
\bmu_{\beta} : Y_{\beta} := \op{Bl}_{\{p_{I}\}_{I \in \oddL}} \to  \Jac_{\beta}
\]
for the blow-up map, and 
\[
f_{\beta} : Y_{\beta} \to X 
\]
for the degree 4 finite morphism resolving the rational map
$\bff_{\Par_{\, \sC_{\beta}}}$, and $E_{I} \subset Y_{\beta}$ for the
exceptional divisor corresponding to the point $p_{I}$.

For future reference we will need to describe the structure of the map
$f_{\beta}$ in more detail.  Let $\Theta \subset \Jac_{\beta}$ be the theta divisor
corresponding to the point $\tilde{p}_{6}$, i.e. the image of
$\sC_{\beta}$ under the Abel-Jacobi map $\sC_{\beta} \to \Jac_{\beta}$,
$x \to \mathcal{O}_{\sC_{\beta}}(x - \tilde{p}_{6})$. Each $2$-torsion
point $p_{I} \in \Jac_{\beta}[2]$ defines a translate $\Theta_{I} =
t_{p_{I}}^{*}\Theta$ of $\Theta$ and we will write $G_{I} \subset
Y_{\beta}$ for its strict transform in $Y_{\beta}$.  
Note that, as abstract varieties, both $G_I$ and $\Theta_I$ 
are isomorphic to $\sC_{\beta}$.
With this notation we now have the following:

\begin{lem} \label{lem:bl.nilp} Let $\beta   \in
  \bB_{\op{nilp}}$ be a generic point, i.e. $\beta =
  (0,\beta_{2})$ with $\beta_{2} \in H^{0}(C,\mathcal{O}(6))$
  vanishing at six distinct points, five of which are the parabolic
  points $p_{1}, \ldots, p_{5} \in C$. Then the degree 4 finite
  morphism $f_{\beta} : 
  Y_{\beta} \to X$ satisfies $f_{\beta}^{*} L_{I} = 2E_{I} +
  G_{I}$. In particular $\sum_{I \in \oddL} E_{I}$ is contained in the
  ramification divisor of $f_{\beta}$. The projection from $G_I$ to $L_I$ is the
standard 2 to1 hyperelliptic map from a genus 2 curve to the projective line.

\end{lem} 
{\bfseries Proof.} \ By specializing the reasoning of
Theorem~\ref{theo:map.from.jacobian} we have that the degree $4$
finite morphism $f_{\beta} : Y_{\beta} \to X$ satisfies
$f_{\beta}^{*}K_{X}^{-1} = \theta^{\otimes 4}(- \sum_{I \in \oddL}
E_{I})$, where $\theta$ is the canonical theta line bundle on
$Y_{\beta}$, i.e. $\theta = \bmu_{\beta}^{*}
\mathcal{O}_{\Jac_{\beta}}(\Theta)$.

In particular the morphism $f_{\beta} : Y_{\beta} \to X \subset
\mathbb{P}^{4}$ to $\mathbb{P}^{4}$ is given by the line bundle 
$\theta^{\otimes 4}(- \sum_{I \in \oddL} E_{I})$. By construction
$f_{\beta}$ maps both  $E_{I}$ and $G_{I}$ to the line $L_{I}$. Again
$f_{\beta}$ must  map $E_{I}$ isomorphically to $L_{I}$ since 
$\theta^{\otimes 4}(- \sum_{I \in \oddL} E_{I})$ restricted to any given
$E_{J}$ gives $\mathcal{O}_{E_{J}}(1)$.   

Furthermore the restriction of $f_{\beta}$ to $G_{J}$ is given by the
linear system 
\[
\theta^{\otimes 4}\left(- \sum_{I \in \oddL}
E_{I}\right)_{|G_{J}} = K_{G_{J}}^{\otimes 4}\left( - \sum_{\substack{I \in
    \oddL \\ p_{I} \in \Theta_{J}}} p_{I}\right). 
\]
Since the configuration of points $p_{I}$ and curves $\Theta_{I}$ in
$\Jac_{\beta}$ is invariant under translations by points of order $2$, it
suffices to understand this linear system on one of the curves
$G_{J}$.

Take $J = \{ 1,2,3,4,5 \}$. Then the curve $\Theta_{\{1,2,3,4,5 \}}$
 is simply the canonical theta divisor $\Theta$ which contains exactly
 six of the
 $2$-torsion points , namely $p_{\{1\}}$, $p_{\{2\}}$, $p_{\{3\}}$,
 $p_{\{4\}}$, $p_{\{5\}}$, $p_{\{1,2,3,4,5\}}$. Using the Abel-Jacobi
 map to identify $G_{\{1,2,3,4,5 \}} = \Theta$ with $\sC_{\beta}$  we
 therefore get  
\[
\begin{aligned}
 K_{G_{\{ 1,2,3,4,5 \}}}^{\otimes 4}\left( - \sum_{\substack{I \in
    \oddL \\ p_{I} \in \Theta_{\{ 1,2,3,4,5 \}}}} p_{I}\right) & = 
K_{G_{\{ 1,2,3,4,5 \}}}^{\otimes 4}\left(- p_{\{1\}} -  p_{\{2\}} -
  p_{\{3\}} - p_{\{4\}} - p_{\{5\}} - p_{\{1,2,3,4,5\}}\right) \\
& \cong K_{\sC_{\beta}}^{\otimes 4}(- \tilde{p}_{1} - \tilde{p}_{2} -
\tilde{p}_{3} - \tilde{p}_{4} - \tilde{p}_{5} - \tilde{p}_{6} )
\\
& = K_{\sC_{\beta}}.
\end{aligned}
\]
Therefore $f_{\beta} : G_{J} \to L_{J}$ is the hyperelliptic
map. Since $f_{\beta}$ is of degree $4$ we then  must have
$f_{\beta}^{*}L_{J} = 2E_{J} + G_{J}$.
\ \hfill $\Box$

\

\begin{rem} \label{rem:GI} {\bfseries (a)} \ For the purposes of Chern
  class and Hecke calculations later on it will be useful to spell out
  the intersection patterns of the curves $E_{I}$ and $G_{I}$ on the
  smooth surface $Y_{\beta}$. To record these patterns efficiently we
  will write $E$ and $G$ for the $1\times 16$ row vectors with entries
  $E_{I}$ and $G_{I}$ respectively. The intersection patterns of the
  $E_{I}$'s and the $G_{I}$'s then are captured completely in the
  three $16\times 16$ matrices $E^{t}E$, $G^{t}E$, and $G^{t}G$ with
  entries given by the intersection numbers $E_{I}\cdot E_{J}$,
  $G_{I}\cdot E_{J}$, and $G_{I}\cdot G_{J}$ respectively.

  Since the $E_{I}$'s are the exceptional divisors blowing up 16
  distinct points on $\Jac_{\beta}$, they are disjoint and for each $I$
  we have $E_{I}^{2} = -1$.

Similarly, the curves $\Theta_{I}$ in $\Jac_{\beta}$ all intersect
transversally in the points $p_{I}$ and so the blow up separates them
from each other.  In particular the strict transforms  $G_{I}$ are all
disjoint from each other.

Since $G_{I} = \bmu_{\beta}^{*}\Theta_{I}(-\sum_{J \in \oddL, p_{J} \in
  \Theta_{I}} E_{J})$ and each $\Theta_{I}$ passes through exactly six
points of order $2$ we have that
\[
G_{I}^{2} = \Theta_{I}^{2}  + \sum_{J \in \oddL, p_{J} \in \Theta_{I}}
  E_{J}^{2} = 2 - 6 = -4. 
\]
Finally from the proof of Lemma~\ref{lem:bl.nilp} it follows that the
curve $G_{I}$ intersects transversally exactly six of the exceptional
divisors $E_{J}$: the divisor $E_{I}$ plus the five divisors $E_{J}$
corresponding to the five $J \in \oddL$ such that the line $L_{J}$
intersects the line $L_{I}$ in $X$.
To summarize we have
\begin{equation}\label{intersection pattern}
E^{t}E = - \mathbf{1}, \quad G^{t}E = 2 \cdot \mathbf{1} + \mathbb{I}, 
\quad
G^{t}G = -4\mathbf{1},
\end{equation}
where $\mathbf{1}$ denotes the $16\times 16$ identity matrix, and
$\mathbb{I} := L^{t}L$ is the intersection matrix of the $16$  lines on
the del Pezzo surface $X$. (Note that the diagonal entries in the middle equation are: $1=2 \cdot 1 - 1$.)

\

\noindent {\bfseries (b)} \ In particular, we note that the divisor $E+G$ on $Y_{\beta}$ has normal crossings, as does the divisor $L$ on $X$.

\

\noindent {\bfseries (c)} \ It is instructive to compare $Y_{\beta}$
and the configuration of curves $E$ and $G$ with the configuration of
curves on the modular spectral cover coming from the generic
symplectic leaf in $\Higgs$. As we saw in section~\ref{sec:modular},
if $\gamma \in \bB$ is a general point, and we choose a connected
component of the Hitchin fiber we again get a modular spectral cover
\index{terms}{Hitchin!fiber} $f_{\gamma} : Y_{\gamma} \to X$. If we
identify $\Pic^{3}(\sC_{\gamma})$ with
$\Jac_{\gamma} = \Pic^{0}(\sC_{\gamma})$ via the map
\[
\Jac_{\gamma} \to
  \Pic^{3}(\sC_{\gamma}), \quad A \mapsto A\otimes
  {\textstyle K_{\sC_{\gamma}}(-\sum_{i=1}^{5} \tilde{p}_{i})},
\] 
where
  $\tilde{p}_{i}$ are the five marked points on $\sC_{\gamma}$ that
  label the connected component of the Hitchin fiber, then again we
  can view  $Y_{\gamma}$ as the  blow 
up  of \index{terms}{Hitchin!fiber}
$\Jac_{\gamma}$ at $16$
distinct points, and $f_{\gamma}$ is a finite morphism of degree
$4$.

In section~\ref{sec:modular} we also argued that 
$f_{\gamma}^{*}L_{I} = E_{I} +
\bD_{I}$, where $\bD_{I} \cong \sC_{\gamma}$ is a strict transform of
a copy of the theta divisor on $\Jac_{\gamma}$ and the map $f_{\gamma} :
\bD_{I} \to L_{I}$ is a $3$-sheeted cover.

In this case we had the intersection pattern $E_{I}^{2} = -1$,
$E_{I}\bD_{I} = 0$, $\bD_{I}^{2} = -3$. In particular when we specialize
the general symplectic leaf to the nilpotent symplectic leaf and
specialize $\gamma \in \bB$ to $\beta\in \bB_{\op{nilp}}$, then 
$Y_{\gamma}$
specializes to $Y_{\beta}$ so that 
\begin{itemize}
\item $E_{I}$ in $Y_{\gamma}$ specializes to $E_{I}$ in $Y_{\beta}$,
\item $D_{I}$ in $Y_{\gamma}$ specializes to $E_{I}+G_{I}$ in
  $Y_{\beta}$. 
\end{itemize}
\end{rem}

\

\noindent
We conclude this section with  a synthetic description of the
modular spectral cover \linebreak $f_{\beta} : Y_{\beta} \to X$ for a general
$\beta \in \bB_{\op{nilp}}$. As we saw above the morphism 
\[
\xymatrix@1@C+0.5pc@M+0.5pc{
Y_{\beta} \ar[r]^-{f_{\beta}} & X \ar@{^{(}->}[r] & \mathbb{P}^{4}
}
\]
is given by a subsystem of the linear system $\left|\theta^{\otimes
    4}(-\sum_{I\in \oddL} E_{I})\right|$ on $Y_{\beta}$. It is not
hard to describe this linear system explicitly.

\

\begin{prop} \label{prop:Ybar} Let $f_{\beta} : Y_{\beta} \to X$ be the
  modular spectral cover corresponding to a general
$\beta \in \bB_{\op{nilp}}$, i.e. to a $\beta = (0,\beta_{2})$ with
$\beta_{2}$ vanishing at six distinct points $p_{1}, \ldots, p_{6}$
five of which are the parabolic points on $C$. Then:
\begin{itemize}
\item[{\em\bfseries (a)}] The inversion involution $\bimath$
  on the Abelian surface
  $\Jac_{\beta}$ lifts to an involution $\widehat{\bimath}$  on
  $Y_{\beta}$. The 
linear
  subsystem $H^{0}(Y_{\beta}, \theta^{\otimes 4}(-\sum_{I\in \oddL}
  E_{I})) \subset H^{0}(Y_{\beta}, \theta^{\otimes 4})$ coincides with
  the subsystem of $\widehat{\bimath}$-anti-invariant sections 
$H^{0}(Y_{\beta},
  \theta^{\otimes 4})^{-} \subset H^{0}(Y_{\beta}, \theta^{\otimes
    4})$. In particular 
\[
H^{0}(Y_{\beta}, \theta^{\otimes
    4}(-\sum_{I\in \oddL} E_{I})) = H^{0}(Y_{\beta}, \theta^{\otimes
    4})^{-} = \wedge^{2} H^{0}(\Jac_{\beta},\mathcal{O}(2\Theta))
\] is $6$
  dimensional.
\item[{\em\bfseries (b)}] The complete linear system 
$|\theta^{\otimes 4}(-\sum_{I\in \oddL}
  E_{I})| = |\theta^{\otimes 4}|^{-}$ gives a morphism $Y_{\beta} \to
  \mathbb{P}^{5}$ which factors throught the Kummer $K3$ surface
$\overline{Y}_{\beta} = Y_{\beta}/\widehat{\bimath}$:
\[
\xymatrix@1@C+0.5pc@M+0.5pc{ 
Y_{\beta}  \ar[r]^-{2:1} & \overline{Y}_{\beta}
  \ar@{^{(}->}[r] & \mathbb{P}^{5}
}
\]
\item[{\em\bfseries (c)}] The six points $p_{1}$, \ldots, $p_{6}$ are in
  one to one correspondence with six points $\Pi_{1}$, \ldots,
  $\Pi_{6}$ in $\mathbb{P}^{5}$ which do not lie in
  $\overline{Y}_{\beta}$ and are the vertices of a tetrahedron. The
  linear projection centered at $\Pi_{6}$ gives a $2:1$ morphism
$\bar{f}_{\beta} : \overline{Y}_{\beta} \to X \subset \mathbb{P}^{4}$
which realizes the Kummer $K3$ surface as a double cover of the del
Pezzo surface $X$ branched at a smooth bi-anti-canonical section.
The modular spectral cover map $f_{\beta} : Y_{\beta} \to X$ is then
the composition
\[
\xymatrix@1@C+0.5pc@M+0.5pc{
 Y_{\beta} \ar[r] & \overline{Y}_{\beta} \ar[r]^-{\bar{f}_{\beta}} & X.
}
\]
\end{itemize}
\end{prop}
{\bfseries Proof.} By definition $H^{0}(Y_{\beta}, \theta^{\otimes
  4}(-\sum_{I\in \oddL} E_{I})) = H^{0}(\Jac_{\beta},
\mathcal{O}(4\Theta)\otimes \mathcal{I})$, where $\mathcal{I} \subset
\mathcal{O}_{\Jac_{\beta}}$ is the ideal sheaf of the $16$ points
$\{p_{I}\}_{I \in \oddL}$. In other words we can identify
$\left|\theta^{\otimes 4}(-\sum_{I\in \oddL} E_{I})\right|$ with the
linear system of all sections of $\mathcal{O}_{\Jac_{\beta}}(4\Theta)$
that vanish at the $16$ points $\{p_{I}\}_{I \in \oddL}$. This linear
system is at least $6$ dimensional. Indeed, let
$\bimath$ be the inversion involuton on the Abelian
surface $\Jac_{\beta}$. Since the theta divisor $\Theta$ is symmetric,
the line bundle $\mathcal{O}_{\Jac_{\beta}}(4\Theta)$ has a natural
$\bimath$ equivariant structure. In particular every section of
$\mathcal{O}_{\Jac_{\beta}}(4\Theta)$ which is anti-invariant with
respect to $\bimath$ will automatically vanish at the fixed points of
$\bimath$, i.e. will belong to $H^{0}(\Jac_{\beta},
\mathcal{O}(4\Theta)\otimes \mathcal{I})$. In other words 
if we denote the space of $\bimath$-anti-invariant sections by 
$H^{0}(\Jac_{\beta}, \mathcal{O}(4\Theta))^{-}$
we have a
natural inclusion
\[
H^{0}(\Jac_{\beta}, \mathcal{O}(4\Theta))^{-} \subset H^{0}(\Jac_{\beta},
\mathcal{O}(4\Theta)\otimes \mathcal{I}).
\]
In fact one can argue that for a smooth
  $\sC_{\beta}$ this inclusion must be an isomorphism. Indeed this is
  equivalent to saying that if we look at the Kummer surface
  $\Jac_{\beta}/\bimath \subset
  \mathbb{P}(H^{0}(\Jac_{\beta},\mathcal{O}(2\Theta))^{\vee}) =
  \mathbb{P}^{3}$, there is no quadric in $\mathbb{P}^{3}$ that
  contains all $16$ nodes of  $\Jac_{\beta}/\bimath$. This can be checked by
  specializing the curve $\sC_{\beta}$ to the union of two distinct elliptic
  curves. Under this specialization the Kummer surface specializes to
  a double quadric and so in the limit there is a single quadric
  containing the images of 
  $16$ points of order $2$ and a $9$ dimensional space of quadrics not
passing through all those points. However when we move away from the
product of elliptic curves, this special quadric moves away from the
points and we get
the statement.

The $6$ dimensional linear system $H^{0}(\Jac_{\beta},
\mathcal{O}(4\Theta))^{-} =
\wedge^{2}H^{0}(\Jac_{\beta},\mathcal{O}(2\Theta))$ gives a rational map
\begin{equation} \label{eq:ai.rat}
\xymatrix@1{\Jac_{\beta} \ar@{-->}[r] & \mathbb{P}^{5} 
= \mathbb{P}( (H^{0}(\Jac_{\beta}, \mathcal{O}(4\Theta))^{-} )^{\vee}),
}
\end{equation}
with base locus the $16$ points of order $2$ in $\Jac_{\beta}$.

Blowing up this base locus resolves the rational map \eqref{eq:ai.rat}
to a morphism
\begin{equation} \label{eq:ai.mor}
Y_{\beta} \to \mathbb{P}^{5}.
\end{equation}
The inversion involution $\bimath : \Jac_{\beta} \to \Jac_{\beta}$
lifts to an involution $\widehat{\bimath} : Y_{\beta} \to Y_{\beta}$ of the
blow up and the morphism \eqref{eq:ai.mor} is given by the linear
system  $H^{0}(Y_{\beta},\theta^{\otimes 4})^{-}$ of
$\widehat{\bimath}$-anti-invariant sections in $\theta^{\otimes 4}$. In
particular factors through the Kummer $K3$ surface
$\overline{Y}_{\beta} := Y_{\beta}/\widehat{\bimath}$:
\[
\xymatrix@1@C+0.5pc@M+0.5pc{ 
Y_{\beta}  \ar[r]^-{2:1} & \overline{Y}_{\beta}
  \ar@{^{(}->}[r] & \mathbb{P}^{5}
}
\]
which embeds in $\mathbb{P}^{5}$ as a surface of degree $8$. 

Every $K3$ surface of degree $8$ in $\mathbb{P}^{5}$ is an
intersection of three quadrics. The Kummer $K3$ surfaces corresponds
to nets of quadrics $\mathbb{P}^{2} \subset
\mathbb{P}(H^{0}(\mathbb{P}^{5},\mathcal{O}(2)))$ such that the
discriminant curve inside $\mathbb{P}^{2}$ parametrizing singular
quadrics in the net is the union of six lines tangent to a conic. The
points at which these lines are tangent to the conic are exactly the
six branch points of the cover $\sC_{\beta} \to C$. This allows us to
identify the conic in parameter space $\mathbb{P}^{2}$ of the net with
$C$ and the discriminant divisor with $\cup_{i =1}^{6} \ell_{i}$,
where $\ell_{i} \subset \mathbb{P}^{2}$ is the line tangent to $C$ at
$p_{i}$. 

All quadrics in the net corresponding to points in the line $\ell_{i}$
are singular with a comon singularity at a point $\Pi_{i} \in
\mathbb{P}^{5}$. The six points $\{\Pi_{i}\}_{i = 1}^{6}$ do not lie
on the $K3$ surface $\overline{Y}_{\beta} \subset \mathbb{P}^{5}$ and
give a projective basis of $\mathbb{P}^{5}$. 

Projecting $\mathbb{P}^{5}$ from the sixth point $\Pi_{6}$ gives a
rational map 
\[
\xymatrix@1{ \mathbb{P}^{5} \ar@{-->}[r] & \mathbb{P}^{4}}
\]
which maps the $K3$ surface $\overline{Y}_{\beta} \subset
\mathbb{P}^{5}-\{\Pi_{6}\}$ to the del Pezzo surface $X \subset
\mathbb{P}^{4}$. 

The resulting morphism 
\[
\overline{f}_{\beta} : \overline{Y}_{\beta} \to X
\]
is a double cover branched along bicanonical curve $X\cap Q_{\beta}$
where $Q_{\beta} \in H^{0}(\mathbb{P}^{4},\mathcal{O}(2))$ is a
quadric which is not in the pencil defining $X$. 

In fact $\overline{Y}_{\beta}$ is the root cover
defined by  $Q_{\beta}$:
\[
\xymatrix@C-1.6pc{
\overline{Y}_{\beta}  \ar[d]_-{\bar{f}_{\beta}} & \subset &
\op{tot}(\mathcal{O}_{\mathbb{P}^{4}}(1)_{|X}) \ar[d] & \subset &
\op{tot}(\mathcal{O}_{\mathbb{P}^{4}}(1)) \ar[d] & = &
\mathbb{P}^{5}-\{\Pi_{6}\} \ar[d]  \\ 
X & = & X& \subset & \mathbb{P}^{4} & = & \mathbb{P}^{4} 
}
\]
By construction $Q_{\beta} \subset X$ is a smooth genus $5$ curve
which is tangent to all $16$ lines $L_{I} \subset X$. The preimage
$\bar{f}_{\beta}^{-1}L_{I}$ is thus a union $\overline{E}_{I}\cup
\overline{G}_{I}$ of two smooth genus zero curves meeting
transversally at a point. The curves  $\overline{E}_{I}$ and 
$\overline{G}_{I}$ are the images of $E_{I}$ and $G_{I}$ under the
double cover map $Y_{\beta} \to \overline{Y}_{\beta}$ and in fact
$Y_{\beta}$ is the root cover of $\overline{Y}_{\beta}$ branched at
the divisor $\sum_{I \in \oddL} \overline{E}_{I} \subset
\overline{Y}_{\beta}$.  \ \hfill $\Box$

\

\medskip

\subsection{Wobbly, shaky and exceptional loci} 
\label{sssec:wobbly}

In this section we will study the special loci in the moduli space $X$
where the parabolic Hecke eigensheaves will acquire a pole.  The
general analysis by Drinfeld and Laumon
\cite{laumon-nilpotent,laumon-gln} implies that the polar locus is
contained in the image of the non-zero (or equivalently non-dominant)
components of the global nilpotent cone.  Strictly speaking the
setting in \cite{laumon-nilpotent,laumon-gln} deals with the case of
holomorphic (non-ramified) and non-parabolic Higgs bundles on a
compact curve. However parabolic tamely ramified Higgs bundles with
rational parabolic weights and nilpotent residues on $C$ can be
identified with holomorphic non-parabolic Higgs bundles on an orbifold
(root stack) curve mapping to $C$ \cite{HilbertHiggs,parabolicMD}. The
reasoning and conclusions of \cite{laumon-nilpotent,laumon-gln} apply
verbatim to such orbifold curves. 

It will be useful to have a name for the points in  this polar locus. 
For concreteness we will think of $X$ as the moduli space
$N_{0}(0,1/2)$. Recall that all points of $X$ represent stable
parabolic bundles.

\begin{defi} \label{defi:vs.and.wobbly} A parabolic vector bundle
  $(V,\bF,(0,1/2)) \in X$ is {\em\bfseries very stable} if the only
  tamely ramified nilpotent Higgs field on $(V,\bF,(0,1/2)) $ is the
  zero one. A parabolic vector bundle $(V,\bF,(0,1/2)) \in X$ is
  {\em\bfseries wobbly} if it is stable but not very stable. We will
  denote the locus of wobbly bundles by $\bWob \subset X$.
\end{defi}

Thus we expect that the Hecke eigensheaves will have singularities
along the wobbly locus $\bWob \subset X$. On the other hand, from the
point of view of our approach \cite{dp-jdg} through non-abelian Hodge
theory, \index{terms}{Non-abelian Hodge theory}
we expect that a Hecke eigensheaf will be singular along the
so called {\em\bfseries shaky locus}, i.e. the locus in $X$ over which
the map from the Hitchin fiber passing through the point representing
the eigenvalue is not proper. From \cite{laumon-nilpotent,bnr} we know
that the wobbly locus is contained in the shaky locus but apriori the
two might be different. In fact by definition the shaky locus seems to
depend on the eigenvalue while the wobbly locus is fixed and does not
depend on a choice of an eigenvalue or a Hitchin fiber.  In our
setting we have already checked that the shaky locus is independent of
any choices. Indeed, by definition the shaky locus is the image of the
exceptional divisors from the modular spectral cover $Y_{\beta}$. In
the previous two sections we showed that the exceptional divisors
$\{E_{I}\}_{I \in \oddL}$ map to the $16$ lines in $X$, i.e. the shaky
locus is the normal crossings divisor $\cup_{I \in \oddL} L_{I}$ in
$X$. \index{terms}{Hitchin!fiber}
(See Lemma  \ref{lem:bl.nilp} and Remark \ref{rem:GI} {\bfseries (b)})

To fix the details of the polar geometry here we will check that this
divisor of lines describes the wobbly locus as well:

\begin{prop} \label{prop:wobbly} The wobbly locus $\bWob \subset X$  is the
  normal crossings divisor given by the $16$ lines in $X$, i.e. $\bWob
  = \cup_{I \in \oddL} L_{I}$.
\end{prop}
{\bfseries Proof.} \ 
As we saw in the proof of Proposition~\ref{prop:dP5}, the moduli space
$X$ can be identified with the blow-up of a quadric surface:
\begin{equation} \label{eq:blowup.quadric}
X = \op{Bl}_{(0,0),(1,1),(\infty,\infty),(p_{4},p_{5})}\left(
  \mathbb{P}^{1}_{p_{4}}\times \mathbb{P}^{1}_{p_{5}}\right).
\end{equation}
For ease of bookkeeping we label the exceptional curves corresponding to
the points $(0,0)$, $(1,1)$, $(\infty,\infty)$, $(p_{4},p_{5})$ by
$E_{0}$, $E_{1}$, $E_{\infty}$, and $E$ respectively .

In the proof of Proposition~\ref{prop:dP5} we also argued that the
points in $X - E$ correspond to stable parabolic bundles whose
underlying vector bundle is $\mathcal{O}_{C}\oplus \mathcal{O}_{C}$. A
traceless tamely ramified Higgs field $\theta$ on such a bundle is
given by a $2\times 2$ matrix
\[
\theta = 
\begin{pmatrix}
  \theta_{11} & \theta_{12}  \\
  \theta_{21} & -\theta_{11}
\end{pmatrix}
\]
with $\theta_{ij} \in H^{0}(C,\mathcal{O}_{C}(3))$ and such that
$F_{i} \subset \ker  \theta(p_{i})$ for $i = 1,2,3,4,5$.

If $a =
\op{gcd}(\theta_{11},\theta_{12},\theta_{21})$ then $\theta$ is nilpotent
(i.e. $\det \theta = 0$) if and only if $\theta_{11}= abc$, $\theta_{12} =
-ab^{2}$, and $\theta_{21} = ac^{2}$ for some homogeneous polynomials $a$,
$b$, $c$ on $C$ of degrees $\leq 3$. Furthermore, since $\theta_{12}$ and
$\theta_{21}$ have degree $3$, the degrees of $b$ and $c$ must be equal and
so the degree of $a$ has to be odd. Thus we have to consider two cases:

\

\noindent
{\bfseries Case 1.}  \ $\deg a = 3$. In this case $b$ and $c$ must be
constant.  If $\theta$ is not identically zero, then $b$ and $c$ are
both non-zero and the constant matrix
\[
\begin{pmatrix}
bc & -b^{2} \\ c^{2} & -bc
\end{pmatrix}
\]
is a regular nilpotent matrix with eigenvector $\begin{pmatrix} b
  \\ c \end{pmatrix}$.  Since
\begin{equation} \label{eq:nilp.theta}
\theta = a\cdot \begin{pmatrix}
bc & -b^{2} \\ c^{2} & -bc
\end{pmatrix}
\end{equation}
it follows that at any of the $p_{i}$'s which is not a zero of $a$ we
must have $F_{i} = \mathbb{C}\begin{pmatrix} b
  \\ c \end{pmatrix}$. Since $a$ can vanish at most at three of the
five points $\{p_{i}\}_{i = 1}^{5}$, it follows that at least two of
the lines $F_{i}$ must coincide. Conversely, if two of the lines
$F_{i}$ coincide, we can choose a spanning vector $\begin{pmatrix} b
  \\ c \end{pmatrix}$ for this line and then choose $a$ to be the
degree three polynomial vanishing at the other three points.  The
Higgs field $\theta$ given by \eqref{eq:nilp.theta} is a non-zero
nilpotent Higgs field which is automatically compatible with the
parabolic structure.

\

\noindent
{\bfseries Case 2.} \ $\deg a = 1$. In this case $\deg b = \deg c =
1$. Again we have a dichotomy:
\begin{itemize}
\item If $a(p_{i}) = 0$, then the flag $F_{i}$ is uncostrained.
\item If $a(p_{i}) \neq 0$, then $F_{i} = \mathbb{C}\begin{pmatrix} b(p_{i})
  \\ c(p_{i}) \end{pmatrix}$.
\end{itemize}

\

\noindent
Therefore if we encode a stable parabolic bundle
$(\mathcal{O}_{C}\oplus\mathcal{O}_{C},\bF,(0,1/2))$ by the $5$-tuple
of points $\{(p_{i},F_{i})\}_{i=1}^{5}$ in $C\times
\mathbb{P}^{1}$, then this bundle is wobbly if either two of the
five points coincide, or four of the five points lie on a $(1,1)$
curve in $C\times \mathbb{P}^{1}$.

The condition that two of the points coincide is easy to analyze in
the model \eqref{eq:blowup.quadric}.  Suppose that we have chosen
coordinates on $C$ so that the three points are $(p_{1},F_{1})$,
$(p_{2},F_{2})$, and $(p_{3},F_{3})$ do not coincide. This means that
the parabolic bundle corresponds to a point in the chart $U_{\{1,2,3\}}$
which we identified with the open set
\[
\mathbb{P}^{1}_{p_{4}}\times \mathbb{P}^{1}_{p_{5}} - \{ (0,0), (1,1),
(\infty,\infty)\},
\]
and $F_{1}$, $F_{2}$, and $F_{3}$ correspond respectively to the
points $0$, $1$, and $\infty$ in the coordinatized line
$\mathbb{P}^{1}_{p_{4}} = \mathbb{P}^{1}_{p_{5}}$.

To get two coinciding points in the $5$-tuple, then we must have one
of the following coincidences

\

\begin{center}
\begin{tabular}{ccc}
  $F_{4} = 0$ & $F_{4} = 1$ & $F_{4} = \infty$ \\
  $F_{5} = 0$ & $F_{5} = 1$ & $F_{5} = \infty$ \\
  & $F_{4} = F_{5}$
\end{tabular}
\end{center}

\

\noindent
These coincidences exactly correspond to the points on the $7$ lines
in Figure~\ref{fig:tencurves}. The remaining coincidences for two of
the five points are outside of the chart $U_{\{1,2,3\}}$ and correspond to
points on one of three exceptional divisor $E_{0}$, $E_{1}$, and
$E_{\infty}$. All together these are the ten lines on the $dP_{4}$
given as the blow up of $\mathbb{P}^{1}_{p_{4}}\times
\mathbb{P}^{1}_{p_{5}}$ at the three points $(0,0)$, $(1,1)$,
$(\infty,\infty)$.

Next we analyze the condition that four of the five points lie on a
$(1,1)$ curve in $C\times \mathbb{P}^{1}$. Let us
again choose affine coordinates $p$ on $\mathbb{P}^{1}$ and
$f$ on $\mathbb{P}^{1}$ so that the parabolic bundle is given
by the $5$-tuple of points $(0,0)$, $(1,1)$, $(\infty,\infty)$,
$(p_{4},F_{4})$, and $(p_{5},F_{5})$. If we assume that say the four
points $(0,0)$, $(\infty,\infty)$, $(p_{4},F_{4})$, and
$(p_{5},F_{5})$ lie on a $(1,1)$ curve, then this curve will give an
isomorphism $C \to \mathbb{P}^{1}$ under which we
have $0 \mapsto 0$, $\infty \mapsto \infty$, $p_{4} \mapsto F_{4}$,
and $p_{5} \mapsto F_{5}$.  However our choice of coordinates already
gives an isomorphism of $C \cong \mathbb{P}^{1}$
under which $0 \mapsto 0$ and $\infty \mapsto \infty$. Therefore in
the cooordinates $p$ and $f$ the isomorphism given by the
$(1,1)$ curve is given by $f = c\cdot p$ for some scaling
constant $c$. Therefore in these coordinates we must have
$f(F_{4}) = c\cdot p(p_{4})$ and $f(F_{5}) = c\cdot
p(p_{5})$, i.e. $f(F_{4})/f(F_{5}) =
p(p_{4})/p(p_{5})$. In the model \eqref{eq:blowup.quadric}
this corresponds to the condition that the point $(F_{4},F_{5}) \in
\mathbb{P}^{1}_{p_{4}}\times \mathbb{P}_{p_{5}}$ lies on the unique
$(1,1)$ curve passing through the points $(0,0)$, $(\infty,\infty)$,
$(p_{4},p_{5})$.

Repeating this analysis for the other sub collections of four points
we conclude that in $X-E$ the wobbly locus consists of the following
fifteen curves:
%fourteen
\begin{itemize}
\item The seven lines in Figure~\ref{fig:tencurves}.
\item The three exceptional divisors $E_{0}$, $E_{1}$, and $E_{\infty}$.
\item The strict transform of the unique $(1,1)$ curve in
  $\mathbb{P}^{1}_{p_{4}}\times \mathbb{P}^{1}_{5}$ passing through the points
  $(0,0)$, $(\infty,\infty)$, and $(p_{4},p_{5})$.
\item The strict transform of the unique $(1,1)$ curve in
  $\mathbb{P}^{1}_{p_{4}}\times \mathbb{P}^{1}_{5}$ passing through the points
  $(0,0)$, $(1,1)$, and $(p_{4},p_{5})$.
\item The strict transform of the unique $(1,1)$ curve in
  $\mathbb{P}^{1}_{p_{4}}\times \mathbb{P}^{1}_{5}$ passing through the points
  $(1,1)$, $(\infty,\infty)$, and $(p_{4},p_{5})$.
\item The strict transform of the ruling  $\{p_{4}\}\times \mathbb{P}^{1}_{5}$.
\item The strict transform of the ruling $\mathbb{P}^{1}_{p_{4}}\times
  \{p_{5}\}$.
\end{itemize}
Note that in the $dP_{5}$ model of $X$ these curves are precisely the
fourteen lines $L_{I}$ corresponding to labels $I \in \oddL$, $I \neq
\{1,2,3,4,5\}$.

In the model \eqref{eq:blowup.quadric} the last line
$L_{\{1,2,3,4,5\}}$ corresponds to the exceptional divisor $E$ so to
complete the proof of the proposition we have to show that $E$
parametrizes wobbly parabolic bundles. Again from the proof of
Proposition~\ref{prop:dP5} we know that $E$ parametrizes stable
parabolic bundles whose underlying bundle is
$V = \mathcal{O}_{C}(-1)\oplus \mathcal{O}_{C}(1)$.  Recall that if we
encode a parabolic bundle
$(\mathcal{O}_{C}(-1)\oplus \mathcal{O}_{C}(1), \bF, (0,1/2))$ by the
$5$-tuple of points $\{(p_{i},F_{i})\}_{i = 1}^{5}$ in the Hirzebruch
surface $\mathbb{P}(V) \cong \mathbb{F}_{2}$, then this parabolic
bundle is stable if and only if:
\begin{itemize}
\item the five points do not lie in the unique $-2$ curve in
  $\mathbb{P}(V)$, i.e. we can not have $F_{i} =
  \mathcal{O}(1)_{p_{i}}$ for all $i = 1, \ldots, 5$, and
\item the five points do not lie on any irreducible curve in the
  linear system $\left| \mathcal{O}_{V}(1)\otimes
    \gamma^{*}\mathcal{O}_{C}(1) \right|$. 
\end{itemize}

A traceless tamely ramified Higgs field on such
a stable bundle is thus given by a matrix
\[
\theta = \begin{pmatrix}
  \theta_{11} & \theta_{12} \\ \theta_{21} & - \theta_{11}
\end{pmatrix}
\]
with $\theta_{11} \in H^{0}(C,\mathcal{O}(3))$, $\theta_{21} \in
H^{0}(C,\mathcal{O}(5))$, $\theta_{12} \in H^{0}(C,\mathcal{O}(1))$
and such that $\theta(p_{i})(F_{i}) \subset F_{i}$.
Again there are two possibilities for $\theta$  to be nilpotent:

\

\noindent
{\bfseries Case 1:} $\theta_{12} = 0$. 

\

\noindent
{\bfseries Case 2:} $\theta_{12} \neq 0$.

\

\noindent
In the first  case we must have
$\theta_{11} = 0$ and $\theta_{21}$ is any homogeneous
polynomial of degree $5$. Notice that if $\theta_{21}$ is the unique
up to scale polynomial that vanishes at $\{p_{i}\}_{i=1}^{5}$, then
the flags $F_{i}$ can be chosen arbitrarily, and so  every stable
parabolic bundle $(V,\bF,(0,1/2)) \in L_{\{1,2,3,4,5\}} \subset X$
admits a non-zero nilpotent Higgs field, namely
\[
\theta = \begin{pmatrix}
0 & 0 \\ \theta_{21} & 0
\end{pmatrix}
\]
where $\theta_{21} \in H^{0}(C,\mathcal{O}(5))$ is a section with
divisor $p_{1} + p_{2} + p_{3} + p_{4} + p_{5}$. This shows that
$L_{\{1,2,3,4,5\}}$ is contained in the wobbly divisor and completes
the proof of the proposition.

We can also analyze the remaining cases. This is not necessary but we
record the answer here for completeness. If $\theta_{21}$ does not
vanish in all five points, then the Higgs field $\theta$ will only be
compatible with certain parabolic bundles in $L_{\{1,2,3,4,5\}}$ so
for those special bundles we will get extra nilpotent Higgs
fields. Similarly from Case 2 we will have additional nilpotent Higgs
fields on bundles in $L_{\{1,2,3,4,5\}}$ corresponding to intersection
points of $L_{\{1,2,3,4,5\}}$ with one of the other lines in $X$.  \
\hfill $\Box$

\section{Hecke eigensheaves} \label{sec:eigensheaves}

Non-abelian Hodge theory \index{terms}{Non-abelian Hodge theory}
converts the question of constructing the geometric Langlands
correspondence for tamely ramified flat bundles on $C$ to an
eigensheaf problem for the Hecke action on tamely ramified semistable
parabolic Higgs bundles on $X$. To spell this problem out we must
first understand the configuration of parabolic divisors in our
birational model for the parabolic Hecke correspondence, as well as
the relevant Higgs data.

\subsection{Parabolic divisors} \label{ssec:pardivisors}

As we saw in Section~\ref{sec:hecke} the parabolic Hecke correspondence
on the moduli space $X \cong dP_{5}$ can be compactified and resolved to the
correspondence
\[
\xymatrix@R-1pc{
& H \ar[dl]_-{p} \ar[dr]^-{q} & \\
X & & X\times C
}
\]
where: 
\begin{itemize}
\item $H = \op{Bl}_{\sqcup_{I} \widehat{L_{I}\times L_{I}}}
    \op{Bl}_{\Delta} (X\times X)$;
\item the map $p : H \to X$  is the
    composition of the blow down map $H \to X\times X$ followed by the first projection
    $\op{pr}_{1} : X\times X \to X$; 
  \item the first coordinate $p_{X}\circ q : H  \to X$ of the
    map $q : H \to X\times C$ is the composition of the blow down map
    $H \to X\times X$ followed by the second projection
    $\op{pr}_{2} : X\times X \to X$;
\item  the second coordinate
    $p_{C}\circ q : H \to C$ is the resolution of the rational
    map $X\times X \dashrightarrow C$ which for a pair of 
    points  $(x,y) \in X\times X$ assigns  $\lambda \in C$
    such that $Q_{\lambda} \subset \mathbb{P}^{4}$  is the unique
    quadric in the pencil parametrized by $C$ which contains the line
    through the two points $x, y \in \mathbb{P}^{4}$.
\end{itemize}

To formulate the parabolic Hecke eigensheaf property we will use the
maps $p$ and $q$ and the projections $p_{X} : X\times C \to X$ and
$p_{C} : X\times C \to C$ to pull back or push forward tamely
ramified parabolic Higgs bundles. For this we need to ensure that the
supports
of the parabolic structures on all participating spaces map to each
other.  The parabolic structures on $C$ are specified at the five
points $\{p_{i}\}_{i =1}^{5}$, while as we saw in section \ref{sssec:wobbly},
the parabolic structures on $X$
can only occur along the wobbly locus $\bWob = \cup_{I} L_{I}$, 
which is a normal crossing divisor in $X$ . Thus we
get the following parabolic divisors: \index{terms}{parabolic!structure}
\begin{description}
\item[$\Par_{C}$] $= p_{1} + p_{2} + p_{3} + p_{4} + p_{5}$.
\item[$\Par_{X}$] $= \sum_{I \in \oddL} L_{I}$.
\item[$\Par_{X\times C}$] $=p_{X}^{-1}(\sum_{I \in \oddL} L_{I}) +
  p_{C}^{-1}(\sum_{i =1}^{5} p_{i})$, that is
\[
\Par_{X\times C} = \sum_{I \in \oddL}L_{I}\times C + \sum_{i
    = 1}^{5} X\times p_{i}.
\]
\item[$\Par_{H}$] is made up from all components of $p^{-1}(\Par_{X})$ and
  $q^{-1}(\Par_{X\times C})$.
\end{description}

\

\noindent
It is not hard to enumerate the components of $\Par_{H}$.

From the iterated blow up description of  $H$ we have that
\[
p^{*}L_{I} = Lp_{I} + LLC_{I},
\] 
where $Lp_{I}$ is the strict transform in $H$ of the divisor
$L_{I}\times X$, while $LLC_{I}$ is the exceptional divisor in $H$
corresponding to blowing up the surface $\widehat{L_{I}\times L_{I}}$
inside $\op{Bl}_{\Delta}(X\times X)$. By symmetry between the two
factors of $X\times X$ we similarly get that
\[
q^{*}(L_{I}\times C) = Lq_{I} + LLC_{I},
\]
where $Lq_{I}$ is the strict transform in $H$ of the divisor
$X\times L_{I}$.

We can also analyze the $q$-preimages of the divisors
$X\times p_{i} \subset X\times C$. Consider the natural morphism
$p_{C}\circ q : H \to C$ and let $RR_{\lambda} \subset H$ denote the
fiber of this morphism over a point $\lambda \in C$. By definition
$RR_{\lambda} = q^{*}(X\times \lambda)$ parametrizes the (closure of
the) locus of all secant lines of $X \subset \mathbb{P}^{4}$ that are
contained in the quadric $Q_{\lambda}$.  Since $X$ is the intersection
of $Q_{\lambda}$ with another quadric, {\em all} lines in the quadric
$Q_{\lambda}$ occur: the forgetful map sends $RR_{\lambda}$ onto the
family $R_{\lambda}$ of lines in $Q_{\lambda}$. In fact,
$RR_{\lambda}$ is a small resolution of $R_{\lambda}$ at the 16 points
corresponding to the $L_I$, each of which is blown up to a line.
In particular if $\lambda \in C - \{p_{1}, \ldots, p_{5} \}$, the
divisor $RR_{\lambda}$ is smooth and irreducible. On the other hand for
every $i = 1, \ldots, 5$ the divisor $RR_{p_{i}}$ will have two
components parametrizing secant lines of $X$ contaned in the two
rulings of the quadratic cone $Q_{p_{i}}$. We will write $R_{1,i}$ and
$R_{2,i}$ for the two components of $RR_{p_{i}}$. Thus
\[
q^{*}(X\times p_{i}) = R_{1,i} + R_{2,i}, \quad \text{for all } i = 1,
\ldots, 5,
\]
and so
\[
\Par_{H} = \sum_{I \in \oddL} Lp_{I} + \sum_{I \in \oddL} Lq_{I} +
\sum_{I \in \oddL} LLC_{I} + \sum_{i = 1}^{5} (R_{1,i} + R_{2,i}).
\]
For the compatibility with the non-abelian Hodge theory we will need the
following statement.

\begin{lem} \label{lem:pearHnc}
  $\Par_{H}$ is a normal crossing divisor away from the $q$-preimage of
  finitely many  points in $X\times C$. 
\end{lem}
{\bfseries Proof.} From the iterated blow-up description of $H$ in
section \ref{ssec-Hecke_X} follows immediately that $\sum_{I \in
  \oddL} Lp_{I} + \sum_{I \in \oddL} Lq_{I} + \sum_{I \in \oddL}
LLC_{I}$ is normal crossings. Indeed, $Lp_{I}$ and $Lq_{I}$ are the
strict transforms of $L_{I}\times X$ and $X\times L_{I}$ in $H$, and
$LLC_{I}$'s are the exceptional divisors for the blow-up of
$\op{Bl}_{\Delta}(X\times X)$ of the smooth disjoint surfaces
$\widehat{L_{I}\times L_{I}}$. The divisor $\sum_{i = 1}^{5} (R_{1,i}
+ R_{2,i})$ is also clearly normal crossing from the secant line
description given above. This shows that $\Par_{H}$ is normal
crossings away from the $q$-preimages of the points of the form
$(L_{I}\cap L_{J})\times p_{i}$ for some intersecting lines $L_{I},
L_{J} \subset X$. Since $q$ is a flat map of relative dimension one
this proves the lemma. In fact it is not hard to analyze the points of
$\Par_{H}$ that lie on these curves and check directly that they are
also normal crossing points. Since this analysis is tedious and is not
needed for our arguments, we will not carry it out here. \ \hfill
$\Box$

\

\noindent
In order to perform computations with parabolic pullbacks and
pushforwards efficiently we will need to order the two components of
the divisors $RR_{p_{i}}$ consistently with the relative position of
the point $p_{i}$ to each of the lines $L_{I}$.

\subsection{Consistent labeling} \label{ssec:ordering}
The labels $I$ we use to enumerate the lines in $X$ or the $2$-torsion
points in $\Jac_{\beta}$ can be thought of as members of any one of the
following three 
bijective sets:
\[
\begin{aligned}
  \oddL & = \left(\text{\begin{minipage}[c]{1.7in} subsets of
      $\{1,2,3,4,5\}$ of odd cardinality
  \end{minipage}}\right), \\[+0.5pc]
  \evenL & = \left(\text{\begin{minipage}[c]{1.7in} subsets of
      $\{1,2,3,4,5\}$ of even cardinality
  \end{minipage}}\right), \\[+0.5pc]
\boldsymbol{\sf{P}} & = \left(\text{\begin{minipage}[c]{3in} even cardinality
    subsets of $\{1,2,3,4,5,6\}$ modulo complementation
  \end{minipage}}\right).
\end{aligned}
\]
The bijections are given by
\[
\xymatrix@R-1pc{
\evenL \ar[r] & \oddL \\
A \subset \{1,2,3,4,5\}  \ar@{|->}[r] & A^{c}\subset \{1,2,3,4,5\} 
}
\]
and
\[
\xymatrix@R-1pc{
\evenL \ar[r] &  \boldsymbol{\sf{P}}\\
A \subset \{1,2,3,4,5\}  \ar@{|->}[r] & A \subset  \{1,2,3,4,5,6\},
}
\]
where $A^{c}$ is the complement of $A$ in $\{1,2,3,4,5\}$.

\

\noindent
Note that any of these indexing sets can be viewed as a
$4$-dimensional vector space over $\mathbb{F}_{2}$. For instance
we have an embedding 
\[
  \evenL = \op{Fun}_{0}\left( \{1,2,3,4,5\},\mathbb{F}_{2} \right)
  \subset \op{Fun}\left( \{1,2,3,4,5\},\mathbb{F}_{2} \right)
\]
as the subset of all functions whose values add up to zero. Similarly 
$\boldsymbol{\sf{P}}$ is a subquotient of $\op{Fun}\left(
  \{1,2,3,4,5,6\},\mathbb{F}_{2} \right)$, namely 
\begin{equation} \label{eq:4dim}
\boldsymbol{\sf{P}} = \op{Fun}_{0}\left(
  \{1,2,3,4,5,6\},\mathbb{F}_{2} \right)/\langle\boldsymbol{\sf{1}}\rangle,
\end{equation}
where again $\op{Fun}_{0}$ denotes the functions whose values add up
to zero, and $\boldsymbol{\sf{1}}$ is the constant function which
takes value one at each point of $\{1,2,3,4,5,6\}$.

\

\noindent
Furthermore, we can map each element in $\{1, 2, 3, 4, 5, 6\}$ to a
natural vector in the $\mathbb{F}_{2}$-vector space
$\boldsymbol{\sf{P}} = \evenL = \oddL$. Namely, in the interpretation
given in equation \eqref{eq:4dim}, we set
\[
\xymatrix@M+0.5pc@R-1pc{
  6 \ar@{|->}[r] & (\text{the zero function}) \\ i \neq 6 \ar@{|->}[r] &
  \left(\text{\begin{minipage}[c]{2in} the function taking value $0$
      at $i$ and $6$ and value $1$ everywhere
      else \end{minipage}}\right).}
\]
We will write $\underline{i}$ for the vector in $\boldsymbol{\sf{P}} =
\evenL$ corresponding to $i \in \{1, 2, 3, 4, 5, 6\}$.  In the identification
$\evenL = \op{Fun}_{0}(\{1, 2, 3, 4, 5\},\mathbb{F}_{2})$, we have that
\[
\begin{aligned}
  \underline{6} & = (\text{the zero function}) \\
  \underline{i} \neq \underline{6}  & = 
  \left(\text{\begin{minipage}[c]{1.8in} the function taking value $0$
      at $i$ and value $1$ everywhere
      else \end{minipage}}\right).
\end{aligned}
\]

\

\noindent
{\bfseries Notation:} For any $I \in \boldsymbol{\sf{P}} = \evenL =
\oddL$ and any $i \in \{1,2,3,4,5,6\}$ we will write $I+i$ for the sum
of the vectors $I$ and $\underline{i}$.

\

\noindent
With this notation we now have the following algebraic characterization of incidences.

\begin{rem} \label{rem:incidences}
\begin{description}
\item[(a)] If $i \in \{1,2,3,4,5\}$ and $I \subset \{1,2,3,4,5\}$ is a
  subset of even cardinality, then
  \[
  \left( i \in I \right) \quad \Longleftrightarrow \quad
  \left(\text{\begin{minipage}[c]{1.8in} the $i$-th coordinate of the
      vector $I$ is $1$
  \end{minipage}}\right)
  \]
\item[(b)] Given $I \in \boldsymbol{\sf{P}} = \evenL = \oddL$, the
  line $L_{I} \subset X$ intersects exactly five other lines in $X$,
  namely the lines $\{L_{I+i}\}_{i = 1}^{5}$.
\end{description}
\end{rem}

\

\noindent
After these preliminaries we can now explain our rule for consistent
ordering of the two components of the divisor $RR_{p_{i}}$.

\begin{lem} \label{lem:ordering}
  For all $i \in \{1,2,3,4,5\}$ we can order the
  two components  $R_{1,i}$ and $R_{2,i}$ of $RR_{p_{i}}$ so that:
  \begin{description}
\item[{\bfseries (i)}] For all $I \in \evenL$ we have that $i \in I$ if and
  only if the line $L_{I} \subset X$ belongs to the ruling $R_{1,i}$ of $Q_{p_{i}}$.
\item[{\bfseries (ii)}] Given a point $\op{pt} \in X$ not lying on any
  of the lines write $L_{I}\times \op{pt}$ for the curve $L_{I}\times
  \op{pt} \subset X\times X$ and also for the strict transform of this
  curve in $H$. Then
  \[
  \begin{aligned}
    \left( \text{\begin{minipage}[c]{1.3in}
        $L_{I}$ belongs to the ruling $R_{1,i}$
    \end{minipage}}\right) & \Longleftrightarrow
    \left(
    \text{
      \begin{minipage}[c]{2in}
     $L_{I}\times \op{pt} \subset H$ is disjoint from the divisor $R_{1,i} \subset H$
      \end{minipage}
    } 
    \right) \\[+0.5pc]
  & \Longleftrightarrow
    \left(
    \text{
      \begin{minipage}[c]{2.5in}
        $L_{I}\times \op{pt} \subset H$ intersects  the divisor $R_{2,i} \subset H$
        transversally at a single point
      \end{minipage}
    } 
    \right)
  \end{aligned}
  \]
  \end{description}
\end{lem}

\

\noindent
{\bfseries Proof.} {\bfseries (i)} We need to specify how we
label the two rulings of each singular quadric (as $R_{1,i}$ and
$R_{2,i}$), as well as the set of 16 lines in $X$ (by labels
$I \in \evenL$). Note that each nodal quadric $Q_{p_{i}}$ contains two
rulings, each 3-dimensional, and these intersect in the surface of
lines through the vertex of $Q_{p_{i}}$. A line $L \subset X$ cannot
contain the vertex, so it belongs to precisely one of these
rulings. We choose one of the lines and label it
$00000 \in \evenL = \op{Fun}_{0}(\{1,2,3,4,5\},\mathbb{F}_{2})$. We
denote the ruling of $Q_{p_{i}}$ containing this line  by $R_{2,i}$, and we let
$R_{1,i}$ denote the other ruling of $Q_{p_{i}}$. We then label each of the
$16$ lines by the $5$-tuple $I$ whose $i$-th entry $r_{I,i}$ is 1 if the
line is in $R_{1,i}$ and 0 otherwise. This labeling has the following
properties:
\begin{itemize}
\item By definition, $i \in I$ (or equivalently, $r_{I,i}=1$) if and
only if the line $L_I \subset X$ belongs to the ruling $R_{1,i}$ of $Q_{p_{i}}$ .
\item In particular, the original line is indeed labelled $00000$.
\item If line $L_I$ meets line $L_J$, the plane they span is contained
in a unique quadric of our pencil. This quadric must be singular, say
$Q_{p_{i}}$. Then $r_{I,i}=r_{J,i}$ but $r_{I,j} \neq r_{J,j}$ for $j \neq
i$.
\item It follows that for intersecting lines $L_{I}$, $L_{J}$ as above,
$J=I+i$.
\end {itemize}
 Our labeling system depends on the choice of the initial line
$L_{00000}$. If we replace it by an intersecting line $L_{\underline{i}}$,
the rulings $R_{1,i}$, $R_{2,i}$ are preserved while the rulings
$R_{1,j}$, $R_{2,j}$ are switched, for $i \neq j$. The effect is then
to perform a translation on all the labels, $J \mapsto J+
i$. Now each of the 16 lines meets either $L_{00000}$ or another
line meeting $L_{00000}$, so iterating the above argument shows that the
effect of changing our initial line from $L_{00000}$ to $L_I$ simply
translates all the line labels, $J \mapsto J+I$, and simultaneously it
flips those ruling labels for which $i \in I$ (or equivalently
$r_{I,i}=1$).

    \

    \noindent
  {\bfseries (ii)} Choose the consistent ordering of the components of
  $RR_{p_{i}}$ from {\bfseries (i)}.  Fix a \linebreak $\op{pt} \in X
  - \cup_{I \in \evenL} L_{I}$ and consider the curve $L_{I}\times
  \op{pt} \subset H$. It is embedded as (the closure of) the set of
  points $(x,\op{pt},p)$ where $x \in L_{I}$ and $Q_{p}$ is the unique
  quadric in the pencil containing the line $\overline{x\, \op{pt}}$ spanned by $x$ and
  $\op{pt}$.

  The intersection of the plane $\op{span}(L_{I},\op{pt}) \subset
  \mathbb{P}^{4}$ with $Q_{p_{i}}$ is a conic which contains the line
  $L_{I}$ and the points $x$ and $\op{pt}$. Thus this intersection is
  $L_{I} \cup \overline{x\, \op{pt}}$. Furthermore, the plane
  $\op{span}(L_{I},\op{pt})$ can not pass through the vertex $v_{i}$
  of the cone $Q_{p_{i}}$. Indeed, the cone $Q_{p_{i}}$ contains the
  plane $\op{span}(L_{I},v_{i})$ and so $\op{span}(L_{I},v_{i})$
  intersects any other quadric $Q_{p}$ in the pencil in $L_{I}$ plus a
  residual line $\ell_{p}$. But then
  \[
  \ell_{p} \subset X = Q_{p_{i}} \cap Q_{p}
  \]
so $\ell_{p} = L_{J}$ for some $J \in \evenL$. Since $\ell_{p}$
depends algebraically on $p \in C$ it follows that there exists a $J
\in \evenL$ so that $\ell_{p} = L_{J}$ for all $p \neq p_{i}$. Since
$X$ is smooth we can not have $L_{J} = L_{I}$ and so $L_{J}$ must be a
line intersecting $L_{I}$. In particular $L_{I}$ and $L_{J}$ belong to
the same ruling of $Q_{p_{i}}$.  By the properties of the consistent
labeling it follows that $i \in I\cap J$, and so $J = I +i$.

Therefore:
\begin{itemize}
\item For any $I$ the five lines intersecting $L_{I}$ obey the following pattern
  \begin{itemize}
  \item $L_{I+i}$ and $L_{I}$ belong to the same ruling of $Q_{p_{i}}$;
  \item for $j \in \{1,2,3,4,5\}-\{i\}$ the lines $L_{I+j}$ and
    $L_{I}$ belong to different rulings of $Q_{p_{i}}$.
  \end{itemize}
\item For any $\op{pt} \in X - \cup_{I\in \evenL}L_{I}$ and any $x \in
  L_{I}$ the lines $\overline{x\, \op{pt}}$ and $L_{I}$ belong to
  different rulings of $Q_{p_{i}}$
\end{itemize}

\

\noindent
This shows that ($L_{I}$ belongs to the ruling $R_{1,i}$)
$\Longleftrightarrow$ ($i \in I$) $\Longleftrightarrow$ ($(L_{I}\times
\op{pt})\cdot R_{1,i} =0$) $\Longleftrightarrow$ ($(L_{I}\times
\op{pt})\cdot R_{2,i} =1$). \ \hfill $\Box$

\

\begin{rem} \label{rem:consistent}
  For future reference note that the same analysis implies also that  \linebreak
($L_{I}$ belongs to the ruling $R_{1,i}$) $\Longleftrightarrow$ ($i
  \in I$) $\Longleftrightarrow$ ($(\op{pt}\times L_{I})\cdot R_{1,i}
  =1$) $\Longleftrightarrow$ ($(\op{pt}\times L_{I})\cdot R_{2,i}
  =0$).
\end{rem}

\subsection{Parabolic Hecke data} \label{ssec:hecke.data}

The parabolic Hecke eigensheaf property whose solution we will feed into
the non-abelian Hodge correspondence involves semistable parabolic
Higgs bundles on the varieties $C$, $X$, and $H$. These enter the
Hecke condition as the Hecke eigenvalue, the Hecke eigenvector, and
the integral kernel for the Hecke operator respectively.

\subsubsection{Hecke eigenvalue} \label{sssec:eigenvalue}

\noindent
For $\lan{G} = GL_{2}(\mathbb{C})$ the eigenvalue for our problem is a
tame parabolic rank two Higgs bundle
\[
\left(E_{\bullet},\theta : E_{\bullet} \to E_{\bullet}\otimes
\Omega^{1}_{C}\left(\log \Par_{C}\right)\right)
\]
on $(C,\Par_{C})$, satisfying the following conditions:
\begin{itemize}
\item $\parch_{1}(E_{\bullet}) = 0$.
\item $\op{res}_{\Par_{C}}(\theta)$ is nilpotent.
\item $(E_{\bullet},\theta)$ is generic in the sense that it
  corresponds to a smooth genus 2 spectral curve:
  \[
  \xymatrix@R-1pc@C-0.5pc@M+0.5pc{
    \sC \ar@{^{(}->}[r] \ar[dr]_-{\pi} &
    \op{tot}\left( \Omega^{1}_{C}(\log \Par_{C})\right) \ar[d] \\
    & C
  }
  \]
Note that the condition that $\op{res}_{\Par_{C}}(\theta)$ is
nilpotent implies that $\pi : \sC \to C$ must be ramified over the
five parabolic points $\{ p_{1}, p_{2}, p_{3}, p_{4}, p_{5} \} =
\Par_{C} \subset C$.  In particular $\sC$ is uniquely determined by
its $6$-th branch point $p_{6} \in C$.
\end{itemize}

\

\index{notations}{Omega1logParC@$\Omega^{1}_{C}(\log \Par_{C})$}
\index{notations}{GL2@$GL_{2}(\mathbb{C})$}
\index{notations}{SL2@$SL_{2}(\mathbb{C})$}
\index{notations}{PSL2@$\mathbb{P}SL_{2}(\mathbb{C})$}

\noindent
For $\lan{G} = SL_{2}(\mathbb{C})$ the eigenvalue is given by the same
exact data but obeys the additional requirement that $E_{\bullet}$ is
an $SL_{2}$-parabolic bundle. It is well known, see e.g.
\cite[Section~2.1]{teleman.woodward} or \cite[Remark~3.7]{bbp}, that
this is equivalent to requiring that for each point in $\Par_{C}$ the
sum of the two parabolic weights assigned to the point is an integer.

To carry out our abelianization calculation we will also need a
spectral description of $(E_{\bullet},\theta)$. Since $\pi : \sC \to
C$ is ramified over the points of $\Par_{C}$, the spectral data for
$(E_{\bullet},\theta)$ will be a parabolic line bundle on $\sC$ with
parabolic structure along the non-reduced divisor $\pi^{*}
\Par_{C}$. As explained in section~\ref{sssec:par.vb}
parabolic bundles with parabolic structures on non-reduced normal
crossings divisors depend on a choice of partitions of the
multiplicities of the divisor components and are given by
semicontinuous families of nested locally free sheaves with labels and
jumps controlled by the parts of the partitions.  

\index{terms}{parabolic!structure!along a non-reduced divisor}

Let
$\tilde{p}_{1},\tilde{p}_{2},\tilde{p}_{3},\tilde{p}_{4},\tilde{p}_{5}$
  be the five ramification points of the cover $\pi : \sC \to C$
  sitting over the points of $\Par_{C}$. Specializing from a general 
  spectral cover for a Hitchin fiber in a symplectic leaf
  \index{terms}{Hitchin!fiber}
  corresponding to regular  semislimple residues one checks that the
  multiplicity $2$ of each component of the non-reduced divisor $\pi^{*}
  \Par_{C} = \sum_{i =1}^{5} 2\tilde{p}_{i}$ has to be partitioned as
  \linebreak
$2 = 1+1$. In other words, the partition decorated parabolic divisor
  on $\sC$ is $\Par_{\sC} = \sP+\sP$ with $\sP = \sum_{i =1}^{5}
  \tilde{p}_{i}$.  The spectral data for $(E_{\bullet},\theta)$ is then
  given by a parabolic line bundle $\mathfrak{a}_{\bullet,\bullet}$ on
  $(\sC,\sP+\sP)$, i.e. a family of line bundles $\{\mathfrak{a}_{s,t}\}_{s,t :
    \sP \to \mathbb{R}}$ on $\sC$ such that
\begin{itemize}
\item $\mathfrak{a}_{s,t} \subset \mathfrak{a}_{s',t'}$ for
  $s \leq s'$ and $t \leq t'$.
\item $\mathfrak{a}_{\bullet,\bullet}$ is semicontinuous for each partial level,
  i.e. for every $s$ and $t$ we can find a real number $c > 0$ so that
  $\mathfrak{a}_{s+\bepsilon,t} = \mathfrak{a}_{s,t+\bepsilon} =
  \mathfrak{a}_{s,t}$ for all
  $\bepsilon : \sP \to [0,c]$.
\item If $\bdelta_{\tilde{p}_{i}} : \sP \to \mathbb{R}$ is the
  characteristic function of $\tilde{p}_{i} \in \sP$, then
  $\mathfrak{a}_{s+\bdelta_{\tilde{p}_{i}},t} = \mathfrak{a}_{s,t+\bdelta_{\tilde{p}_{i}}} =
  \mathfrak{a}_{s,t}(\tilde{p}_{i})$
   and these identifications are compatible with the inclusions of the
  $\mathfrak{a}_{s,t}$'s.
\end{itemize}
The fact that $(E_{\bullet},\theta)$ is given by the parabolic
spectral data $\mathfrak{a}_{\bullet,\bullet}$ means that for every \linebreak
$t : \Par_{C}
\to \mathbb{R}$ we have
\[
(E_{t},\theta) = \left(\pi_{*}\mathfrak{a}_{t,t},\pi_{*}\lambda\right),
\]
where $\lambda$ is the restriction of the tautological section in
$\Gamma(\op{tot}(\Omega^{1}_{C}(\log \Par_{C})), \op{pr}^{*}
\Omega^{1}_{C}(\log \Par_{C}))$ to $\widetilde{C} \subset
\op{tot}(\Omega^{1}_{C}(\log \Par_{C}))$, and we slightly abuse
notation and write $t$ both for the map $t : \Par_{C} \to \mathbb{R}$
and for the composition $\sP \to \Par_{C} \stackrel{t}{\to}
\mathbb{R}$.

\index{notations}{adotfr@$\mathfrak{a}_{\bullet,\bullet}$}

For future reference note that every parabolic line bundle
$\mathfrak{a}_{\bullet,\bullet}$ on $(\sC,\sP+\sP)$ is specified by a
choice of a fixed ordinary line bundle and two vectors in
$\mathbb{R}^{5}$ specifying the jumps of the parabolic levels of
$\mathfrak{a}_{\bullet,\bullet}$ along each part of the partition.  In
other words, given $\mathfrak{a}_{\bullet,\bullet}$ we can always
find:
\begin{itemize}
\item a degree zero line bundle $\mathfrak{a} \in \Pic^{0}(\sC)$;
\item two vectors $a, b \in \mathbb{R}^{5}$;
\end{itemize}
so that $\mathfrak{a}_{\bullet,\bullet} =
\mathfrak{a}(a\sP+b\sP)_{\bullet,\bullet}$.
Explicitly this means that for all $s,
t \in \mathbb{R}^{5}$ we have
\begin{equation} \label{eq:parline}
  \mathfrak{a}_{s,t} = \mathfrak{a}(a\sP+b\sP)_{s,t} =
  \mathfrak{a}(\fl*{s+a}\sP + \fl*{t+b}\sP). 
\end{equation}
Note however that the parabolic line bundle
$\mathfrak{a}_{\bullet,\bullet}$ does not determine the triple
$(\mathfrak{a},a,b)$ uniquely. The formula \eqref{eq:parline} defines
a surjective group homomorphism to the Picard group of all parabolic
line bundles on $(\sC,\sP+\sP)$:
\[
\xymatrix@R-1pc{
\Pic^{0}(\sC)\times \mathbb{R}^{5}\times \mathbb{R}^{5}
\ar[r] & \op{parPic}(\sC,\sP+\sP) \\
(\mathfrak{a},a,b) \ar[r]  & \mathfrak{a}(a\sP+b\sP)_{\bullet,\bullet},
}
\]
but this  homomorphism is not injective. Its kernel consists of all triples
\[
\left(\mathcal{O}_{\sC}(-(m_{1}+m_{2})\sP), m_{1}, m_{2} \right),
\]
where $m_{1}, m_{2} \in \mathbb{Z}^{5}$ are such that $\deg
((m_{1}+m_{2})\sP) = 0$ (and so $\mathcal{O}_{\sC}(-(m_{1}+m_{2})\sP)$
is a point of order $2$ in $\Pic^{0}(\sC)$).

\

\noindent
Furthermore, consider the subgroup $\mathfrak{S} \cong
(\mathbb{Z}/2)^{5} \subset GL(\mathbb{R}^{5}\oplus \mathbb{R}^{5})$
  consisting of all linear transformations $(a,b) \to
  (\mathfrak{s}_{1}(a,b),\mathfrak{s}_{2}(a,b))$
  that swap some coordinates of $a$
  with the corresponding coordinates of $b$. It is straighforward to
  check that for any $\mathfrak{s} \in \mathfrak{S}$ the formula
  \[
  \mathfrak{a}(a\sP + b\sP)_{\bullet,\bullet} \to
  \mathfrak{a}(\mathfrak{s}_{1}(a,b)\sP +
  \mathfrak{s}_{2}(a,b)\sP)_{\bullet,\bullet}
  \]
  gives a well defined action
  of $\mathfrak{S}$ on $\op{parPic}(\sC,\sP+\sP)$. This action is non
  trivial. However, since $E_{t} = \pi_{*}\mathfrak{a}((\fl*{t+a} +
  \fl*{t+b})\sP) = \pi_{*}\mathfrak{a}((\fl*{t+\mathfrak{s}_{1}(a,b)} +
  \fl*{t+\mathfrak{s}_{2}(a,b)})\sP)$ it follows that the two spectral
  parabolic line bundles $\mathfrak{a}(a\sP + b\sP)_{\bullet,\bullet}$
  and $\mathfrak{a}(\mathfrak{s}_{1}(a,b)\sP +
  \mathfrak{s}_{2}(a,b)\sP)_{\bullet,\bullet}$ define the same strongly
  parabolic rank two Higgs bundle $(E_{\bullet},\theta)$ on
  $(C,\Par_{C})$.

  \

  \index{notations}{Sfrak@$\mathfrak{S}$}
  \index{notations}{sfrak@$\mathfrak{s}$}
    \index{notations}{sfraki@$\mathfrak{s}_{i}$}
\index{terms}{Weyl group}

  \begin{rem} \label{rem:swap} Intrinsically the group $\mathfrak{S}$ is
the product of several copies of the Weyl group of $\lan{G} =
GL_{2}(\mathbb{C})$ or $SL_{2}(\mathbb{C})$, one for each point in
$\Par_{C}$. The ambiguity introduced by the action of this group is a
familiar ambiguity in the formulation of the tamely ramified GLC. In
the setup discussed in section~\ref{ssec:setup.glc} this ambiguity
appears in the condition that the relevant symplectic leaf of the
moduli of flat $\lan{G}$-bundles is labeled by a map $\bLambda :
\Par_{C} \to \lan{\mathfrak{t}}/\!/W$.  Concretely the numbers $a_{i}$,
$b_{i}$ correspond to the eigenvalues of the residue at $p_{i}$ of the
meromorphic connection corresponding to $(E_{\bullet},\theta)$ under
the non-abelian Hodge correspondence. These eigenvalues are unordered
and choosing an ordering introduces the mentioned Weyl group
ambiguity.
\end{rem}

\

\noindent
Finally note that we can also express the condition that
$\parch_{1}(E_{\bullet}) = 0$ in terms of the spectral data
$\mathfrak{a}_{\bullet,\bullet}$. Since $E_{t} =
\pi_{*}\mathfrak{a}_{t,t}$ we have that the parabolic weights of
$E_{\bullet}$ at $p_{i}$ are $-a_{i} + \fl*{a_{i}}$ and $-b_{i} +
\fl*{b_{i}}$ and
\[
E_{0} = \pi_{*} \left(\mathfrak{a}((\fl*{a} + \fl*{b})\sP)\right).
\]
Since there are six branch points, Grothendieck-Riemann-Roch gives
\[
\deg E_{0} = \sum_{i=1}^{5} \left( \fl*{a_{i}} + \fl*{b_{i}}\right) - 3,
\]
and so
\[
\pardeg(E_{\bullet}) = -3 + \sum_{i=1}^{5} \left( \fl*{a_{i}} +
\fl*{b_{i}}\right) - \sum_{i=1}^{5}\left(-a_{i} + \fl*{a_{i}} -b_{i} +
\fl*{b_{i}}\right) = \sum_{i=1}^{5}(a_{i}+b_{i}) - 3.
\]
Thus the condition $\parch_{1}(E_{\bullet}) = 0$ is equivalent to
\begin{equation} \label{eq:spectral.parch}
\boxed{
  \sum_{i=1}^{5}(a_{i}+b_{i}) = 3.
}
\end{equation}
Finally, as mentioned at the beginning of the section, the condition
that $E_{\bullet}$ is an $SL_{2}(\mathbb{C})$ parabolic bundle is
expressed numerically as the requirement that at each parabolic point
the sum of the corresponding parabolic weights is an integer. But as
we just saw, the parabolic weights of $E_{\bullet}$ at the point
$p_{i}$ are the fractional parts $-a_{i} + \fl*{a_{i}}$ and $-b_{i} +
\fl*{b_{i}}$. Therefore $E_{\bullet}$ is an $SL_{2}(\mathbb{C})$
parabolic bundle if in addition we have
\begin{equation} \label{eq:spectral.SL2}
\boxed{
  a+b \in \mathbb{Z}^{5}.
}
\end{equation}

\subsubsection{Hecke eigensheaf} \label{sssec:eigensheaf}

As we saw in section~\ref{ssec:nilpotent.residues} the assumption that
Hecke eigenvalue $(E_{\bullet},\theta)$ has nilpotent residues and a
smooth spectral curve $\sC$ implies that the map from the fiber of the
Hitchin map $h^{-1}(\sC) \cong \Jac(\sC)$ to $X$ has degree four over
the very stable locus $X - \Par_{X}$ and that the associated modular
spectral cover $f : Y \to X$ is a finite morphism of degree
four. \index{terms}{Hitchin!fiber}

Suppose $\lan{G} = GL_{2}(\mathbb{C})$ so that also $G =
GL_{2}(\mathbb{C})$. In this case the eigensheaf problem is most
naturally formulated (see \cite[Section 2]{dp-jdg}) on the rigidified
stack of bundles in which the connected component of the generic
stabilizer is removed. In section \ref{sec:family} we checked that
when the weights $(\ba,\bb)$ belong to the dominant chamber, there are
  no strictly semistable bundles, and so the rigidified moduli stack
  is just the moduli space $\Bun(\ba,\bb)$ of stable parabolic rank
  two vector bundles. In section \ref{sec:moduli.spaces} we proved
  that $\Bun(\ba,\bb)$ is a disjoint union $\Bun(\ba,\bb) = \sqcup_{d
    \in \mathbb{Z}} N_{k}(\ba,\bb)$ with the connected component
  $N_{k}(\ba,\bb)$ parametrizing rank two stable parabolic vector
  bundles $V_{\bullet}$ with parabolic weights $(\ba,\bb)$, and such
  that the level zero bundle $V_{0}$ has degree $k \in \mathbb{Z}$. In
  section \ref{sec:family} we also saw that when the weights
  $(\ba,\bb)$ belong to the dominant chamber, each component
  $N_{k}(\ba,\bb)$ is naturally isomorphic to the $dP_{5}$ surface
  $X$.

  \

  \begin{rem} \label{rem:boldface.aa}
{\bfseries (i)} \
Here and later, e..g. in Section
  \ref{sssec:kernel}, we use our standard notational convention that
  the boldface pair $(\ba,\bb)$ denotes vectors of real parabolic
  weights $\ba = \{a_{i}\}_{i = 1}^{5}$, $\bb = \{b_{i}\}_{i =1}^{5}$,
  with $a_{i}, b_{i} \in (-1,0]$. In contrast the non-boldface letters
$a$ and $b$ used in the spectral data description in the previous
section denote vectors $a, b \in \mathbb{R}^{5}$ whose entries are
arbitrary real numbers. The relationship of the two notations is that
for parabolic bundles given by spectral data with parameters $a$ and
$b$ the parabolic weights are given by $\ba = - a + \fl{a}$ and $\bb =
-b + \fl{b}$.

\

\noindent
{\bfseries (ii)} We described the dominant chamber for the case of
balanced weights.  For general weights the dominant chamber is the
unique open chamber in the space of weights which contains the points
$(\ba,\bb)$ with $b_{i} = a_{i} + 1/2$ for all $i$.x
 \end{rem}

\index{terms}{dominant chamber!of parabolic weights}
\index{notations}{GL2@$GL_{2}(\mathbb{C})$}
\index{notations}{SL2@$SL_{2}(\mathbb{C})$}
\index{notations}{PSL2@$\mathbb{P}SL_{2}(\mathbb{C})$}

\

Similarly if $\lan{G} = SL_{2}(\mathbb{C})$ so $G =
\mathbb{P}GL_{2}(\mathbb{C})$, the Hecke eigensheaf problem is
formulated on the moduli stack $\sBun(\ba,\bb)$ of semistable
parabolic $G$ bundles with parabolic weights $(\ba,\bb)$. This stack
is a disjoint union $\sBun(\ba,\bb) = \sM_{0}(\ba,\bb)\sqcup
\sM_{1}(\ba,\bb)$ with the connected component $\sM_{k}(\ba,\bb)$
parametrizing $\mathbb{P}SL_{2}(\mathbb{C})$ semistable parabolic
bundles which are of the form $\op{ad}(V_{\bullet}) =
\op{End}_{0}(V_{\bullet})$ for a rank two parabolic vector bundle
$V_{\bullet}$ with parabolic weights $(\ba,\bb)$, and such that the
level zero bundle $V_{0}$ has degree $k \in \mathbb{Z}$. Explicitly,
when $(\ba,\bb)$ are in the dominant chamber $\sM_{k}(\ba,\bb)$ is a
smooth Deligne-Mumford stack obtained as the quotient of
$N_{k}(\ba,\bb)$ by the action of the group of $2$-torsion parabolic
line bundles on $C$ which is naturally isomorphic to
$(\mathbb{Z}/2)^{4}$. Also in section \ref{sec:git} we saw that when
$(\ba,\bb)$ are in the dominant chamber and we identify
$N_{k}(\ba,\bb)$ with $X$, this action of $(\mathbb{Z}/2)^{4}$ becomes
the action which flips the signs of subsets of the natural coordinates
of the $\mathbb{P}^{4}$ where $X$ is embedded as the intersection of
two quadrics. Thus for weights in the dominant chamber \linebreak
$\sM_{k}(\ba,\bb) = \left[N_{k}(\ba,\bb)/(\mathbb{Z}/2)^{4}\right]
\cong \left[X/(\mathbb{Z}/2)^{4}\right]$ and as we explained in
\ref{sec:git} its coarse moduli space $M_{k}(\ba,\bb)$ is naturally
identified with projective plane $X/(\mathbb{Z}/2)^{4} = S^{2}C$.  For
GLC we need to work on the stack $\sM_{k}(\ba,\bb)$ (see
\cite{dp-jdg}[Section 1]) so our eigensheaf objects will live on $X$ and
will be $(\mathbb{Z}/2)^{4}$-equivariant.

\index{notations}{Bun@$\Bun$}
\index{notations}{Nkab@$N_{k}(\ba,\bb)$}
\index{notations}{Mkab@$M_{k}(\ba,\bb)$}
\index{notations}{Mkabs@$\sM_{k}(\ba,\bb)$}
\index{terms}{componnet!connected}

Thus for $\lan{G} = GL_{2}(\mathbb{C})$ the Hecke eigensheaf
corresponding to the eigenvalue $(E_{\bullet},\theta)$ will be a
collection $\left\{(F_{\bullet}^{{k}},\varphi^{{k}})\right\}_{k \in
  \mathbb{Z}}$ of tame parabolic rank $4$ Higgs bundles on
$(X,\Par_{X})$, where $(F_{\bullet}^{{k}},\varphi^{{k}})$ is viewed as
the piece of the eigensheaf on the $k$-th connected component of the
moduli space. Similarly, for $\lan{G} = SL_{2}(\mathbb{C})$ the Hecke
eigensheaf corresponding to the eigenvalue $(E_{\bullet},\theta)$ will
be given by a pair
$\left\{(F_{\bullet}^{{k}},\varphi^{{k}})\right\}_{k \in
  \mathbb{Z}/2}$ of tame parabolic rank $4$ Higgs bundles on
$(X,\Par_{X})$, where each $(F_{\bullet}^{{k}},\varphi^{{k}})$ is
$(\mathbb{Z}/2)^{4}$-equivariant and its descent to
$\left[X/(\mathbb{Z}/2)^{4}\right]$ is viewed as the piece of the
eigensheaf on the $k$-th connected component of the moduli space.

\

\index{notations}{GL2@$GL_{2}(\mathbb{C})$}
\index{notations}{SL2@$SL_{2}(\mathbb{C})$}
\index{notations}{PSL2@$\mathbb{P}SL_{2}(\mathbb{C})$}

For both $\lan{G} = GL_{2}(\mathbb{C})$ and $\lan{G} =
SL_{2}(\mathbb{C})$ the Hecke eigensheaves are encoded in the same
type of data on each connected component of the moduli of parabolic
bundles. Following  our strategy, in order for these data to
determine an actual eigensheaf flat bundle or a twisted
$\mathcal{D}$-module, the pieces $(F_{\bullet}^{{k}},\varphi^{{k}})$
should satisfy the conditions of the construction coming from the
Fourier-Mukai transform along the Hitchin fibers and the non-abelian
Hodge correspondence. Before we start analyzing these conditions let
us simplify notation and write $(F_{\bullet},\varphi)$ for any of the
putative eigensheaf pieces $(F_{\bullet}^{{k}},\varphi^{{k}})$. Then
the Fourier-Mukai and non-abelian Hodge correspondence conditions on
such an $(F_{\bullet},\varphi)$ are:

\

\begin{itemize}
\item $\parch_{1}(F_{\bullet}) = 0$, $\parch_{2}(F_{\bullet}) = 0$.
\item $\res_{\Par_{X}}(\varphi)$ is nilpotent.
\item $(F_{\bullet},\varphi)$ is a stable parabolic Higgs bundle
with spectral cover $Y \to X$.
\end{itemize}

\

\noindent
In fact the requirement that the residues of $\varphi$ are nilpotent
is redundant since it is a consequence of the fact that
$(F_{\bullet},\varphi)$ is  given by spectral data on $Y$.

\

To spell out the spectral construction of $(F_{\bullet},\varphi)$
recall (see section~\ref{ssec:nilpotent.residues}) that:
\begin{itemize}
\item[(a)] The modular spectral cover $Y$ comes equipped with a
  natural morphism
\[
  \xymatrix@R-1pc{
    Y \ar[r]^-{\alpha} \ar[rd]_-{f} &
    \op{tot}\left(\Omega^{1}_{X}(\log L)\right) \ar[d] \\
& X
  }
  \]
\item[(b)] The parabolic divisor $\Par_{Y}$ on $Y$ is the non-reduced
  divisor $f^{*}L = 2E + G$ partitioned as $E + (E+G)$.
\end{itemize}

\index{notations}{alpha@$\alpha$}
\index{notations}{Omega1logL@$\Omega^{1}_{X}(\log L)$}
\index{notations}{Omega1logParC@$\Omega^{1}_{C}(\log \Par_{C})$}

\

\noindent
Statement (b) follows from the specialization picture described in
Remark~\ref{rem:GI} {\bfseries (b)}. To explain statement (a) we have
to analyze in more detail the tautological one form on $Y$. Consider
the spectral cover $\pi : \sC \to C$ with its natural inclusion $\sC
\hookrightarrow \op{tot}(\Omega^{1}_{C}(\log \Par_{C}))$.  This
inclusion corresponds to the tautological section
\[
\lambda \in 
\pi^{*}(\Omega^{1}_{C}(\log \Par_{C})) =
\pi^{*}(\Omega^{1}_{C}(p_{1}+p_{2}+p_{3}+p_{4}+p_{5})).
\]
By Hurwitz we have
$\pi^{*}(\Omega^{1}_{C}(p_{1}+p_{2}+p_{3}+p_{4}+p_{5})) =
\Omega^{1}_{\sC}(\tilde{p}_{1}+\tilde{p}_{2}+\tilde{p}_{3}
+\tilde{p}_{4}+\tilde{p}_{5} - \tilde{p}_{6})$. Since
$\theta$ has nilpotent residues at
$p_{1},\ldots, p_{5}$ it follows that $\lambda$
viewed as a meromorphic one form on $\sC$ has vanishing residues at
$\tilde{p}_{1}, \ldots, \tilde{p}_{5}$ and thus comes from the
unique up to scale holomorphic one form on $\sC$ that vanishes twice
at $\tilde{p}_{6}$.

In particular $\lambda$ corresponds to a holomorphic one form on the
Jacobian $\Jac = \Pic^{0}(\sC)$ of $\sC$.  Pulling back this one form
by the blow up map $\bmu : Y \to \Jac$ we obtain a canonical
holomoprphic one form $\alpha \in H^{0}(Y,\Omega^{1}_{Y})$ on $Y$. At
the same time the codifferential of $f$ gives an inclusion of
sheaves\footnote{Here we write $L$ as a shortcut notation for
  $\Par_{X} = \cup_{I} L_{I}$.}  $f^{*}(\Omega^{1}_{X}(\log L))
\subset \Omega^{1}_{Y}(\log f^{*}L) = \Omega^{1}_{Y}(\log (E +
G))$. Statement (a) is the statement that the form $\alpha$ aligns
with this inclusion i.e. $\alpha$ is actually a section in the
subsheaf
  \[
  f^{*}(\Omega^{1}_{X}(\log L))\cap \Omega^{1}_{Y} \subset
  \Omega^{1}_{Y}(\log(E + G)).
  \]
In particular $\alpha$ can be viewed as a map from $Y$ to
$\op{tot}(\Omega^{1}_{X}(\log L))$.

To see this alignment recall that the four sheeted cover $f : Y \to X$
decomposes as an iterated double cover $Y \to \overline{Y} \to X$
where $\overline{Y} = Y/(-1)$ is the Kummer $K3$ for $\Jac$, and $\bar{f} :
\overline{Y} \to X$ is the double cover branched along a smooth curve
$Q \subset X$ which is tangent to each line $L_{I}$ at a single point
corresponding to the point $p_{6} \in C$ under the natural
identification of $L_{I}$ with $C$. 

Denote the ramification divisor of the cover $\overline{Y} \to
X$ by $\overline{Q}$, denote the preimage of $\overline{Q}$
in $Y$ by $R$. We will write $\iota : R \to Y$ for
the natural inclusion and
$q : R \to \overline{Q}$ for the double cover map. 

The iterated double cover picture gives a short exact sequence
\[
0 \to f^{*}(\Omega^{1}_{X}(\log L)) \to \Omega^{1}(\log(E +G)) \to
\iota_{*} (q^{*}N^{\vee}_{\overline{Q}/\overline{Y}})\otimes
\mathcal{O}_{Y}(E+G) \to 0.
\]
But $\alpha$ is the equation of $Q \subset X$ and so at the points of
$R$ which are away from $E$ and $G$ the projection of $\alpha$ to the
conormal line of $Q$ vanishes. In other words, when viewed as a
section in $ \Omega^{1}_{Y}(\log(E +G))$, $\alpha$ belongs to the subsheaf
$f^{*}(\Omega^{1}_{X}(\log L))$. Note also that by construction
$\alpha \in \Omega^{1}_{Y}(\log(E +G))$ has zero residues along the
$E_{I}$'s and the $G_{I}$'s. Thus if we take any line bundle $M$ on
$Y$, then the $f$-pushforward of the Higgs field
\[
(-) \otimes \alpha : M \to M\otimes f^{*}(\Omega^{1}_{X}(\log L))
\]
will yield a rank four  meromorphic Higgs field
\[
f_{*}((-) \otimes \alpha ) : f_{*}M \to f_{*}M \otimes \Omega^{1}_{X}(\log L)
\]
with nilpotent residues along $L$. In particular the stable parabolic
Higgs bundle $(F_{\bullet},\varphi)$ with spectral cover $f : Y \to X$
will have nilpotent residues automatically.

\

\noindent
The modular spectral cover
  $f : Y \to X$ is not embedded in the total space of
  $\Omega^{1}_{X}(\log L)$ but rather maps birationally there. Indeed
  from the definition of $\alpha$ it follows that when viewed as a map
\[
\xymatrix@R-1pc@C-1pc{Y \ar[rr]^-{\alpha} \ar[dr] & &
  \op{tot}\left(\Omega^{1}_{X}(\log L)\right) \ar[dl] \\
& X &
}
\]
it glues the pairs of points in $Y$ sitting over the intersection
points $L_{I}\cap L_{J}$ of a pair of lines in $X$. As explained in
Proposition~\ref{prop:Ybar} (c) we have $f^{*}L_{I} = 2E_{I} + G_{I}$
and $f : E_{I} \to L_{I}$ is an isomorphism, while
$f : G_{I} \to L_{I}$ is a double cover that is naturally identified
with the cover $\widetilde{C} \to C$. Note also that if two lines
$L_{I}$ and $L_{J}$ intersect, then $E_{I}$ intersects $G_{J}$
transversally at a point and is disjoint from $E_{J}$ and
$G_{I}$. Similarly $E_{J}$ intersects $G_{I}$ transversally at a point
and is disjoint from $E_{I}$ and $G_{J}$. The set theoretic preimage
of the intersection point $L_{I}\cap L_{J}$ in $Y$ consists of the two
points $E_{I}\cap G_{J}$ and $E_{J}\cap G_{I}$. The map $\alpha$ glues
together all such pairs of points and is an isomorphism everywhere
else.  The local shape  of the cover $f: Y \to X$ over a neighborhood
of an intersection point $L_{I}\cap L_{J}$ is illustrated in
Figure~\ref{fig:localY} below.

\begin{figure}[!ht]
\begin{center}
\psfrag{LI}[c][c][1][0]{{$L_{I}$}}
\psfrag{LJ}[c][c][1][0]{{$L_{J}$}}
\psfrag{EI}[c][c][1][0]{{$E_{I}$}}
\psfrag{EJ}[c][c][1][0]{{$E_{J}$}}
\psfrag{GI}[c][c][1][0]{{$G_{I}$}}
\psfrag{GJ}[c][c][1][0]{{$G_{J}$}}
\psfrag{X}[c][c][1][0]{{$X$}}
\epsfig{file=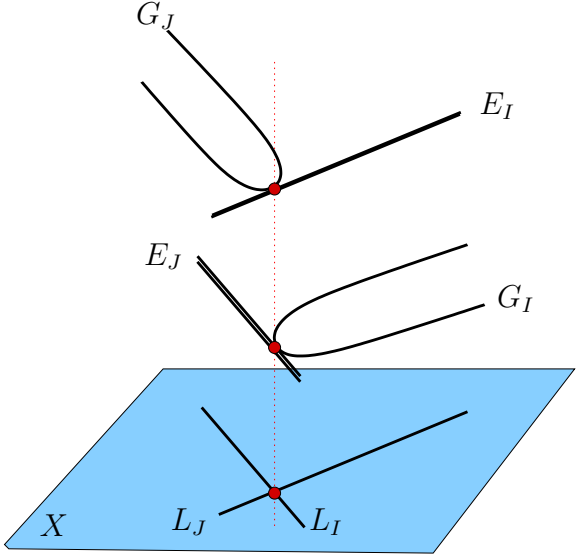,width=3in} 
\end{center}
\caption{Local picture for the cover $f : Y \to X$}\label{fig:localY}  
\end{figure}

\

\noindent
In particular, locally on $X$ near the intersection point of $L_{I}$
and $L_{J}$ the meromorphic Higgs field $\varphi$ looks like
\[
  \varphi(x,y) = \begin{pmatrix} 0 & 0 &  &  \\ 0 & 0 &  &   \\  &
     & 0 & y \\  &  & 1 & 0 \end{pmatrix} \frac{dx}{x}
  + \begin{pmatrix} 0 & x &  &  \\ 1 & 0 &  &   \\  &  & 0 & 0 \\
     &  & 0 & 0 \end{pmatrix} \frac{dy}{y},
\]
where $x = 0$ is the local equation of $L_{I}$ and $y = 0$ is the
local equation of $L_{J}$. If $u$ and $v$ are the coordinates along
the fibers of $\Omega^{1}_{X}(\log L)$ corresponding to the local
frame $\{dx/x,dy/y\}$, then $\alpha(Y) \subset \op{tot}(
\Omega^{1}_{X}(\log L))$ is locally given as:
\[
\alpha(Y) \ : \ \left\{ \begin{aligned} x & = u^{2} \\ v & = 0 \end{aligned}
  \right\} \, \cup \, \left\{ \begin{aligned} y & = v^{2} \\ u & = 0
  \end{aligned} \right\}.
\]
\

\

\noindent
Thus the Higgs field $\varphi$ in the parabolic Hecke eigensheaf
$(F_{\bullet},\varphi)$ is just the pushforward $f_{*}((-)\otimes
\alpha)$ while the parabolic bundle $F_{\bullet}$ is the pushforward
of a parabolic line bundle $M_{\bullet,\bullet}$ on $Y$ with parabolic
structure on the non-reduced divisor $E + (E+G)$.

Every such parabolic line bundle is of the form $M_{\bullet,\bullet} =
\mycal{L}(e E + d(E+G))_{\bullet,\bullet}$  for some line bundle
$\mycal{L}$ on $Y$ and some row vectors of real numbers $e =
(e_{I})_{I \in \oddL}$ and $d = (d_{I})_{I\in \oddL}$, i.e. for any
parabolic levels $S = (S_{I})_{I \in \oddL}$ and $T = (T_{I})_{I \in
  \oddL}$ we have
\begin{equation}
\label{eq:MST}
M_{S,T} = \mycal{L}(\fl{e + S}E + \fl{d +T}(E+G)).
\end{equation}

\

\noindent
For $\lan{G} = GL_{2}(\mathbb{C})$ or $\lan{G} = SL_{2}(\mathbb{C})$
the Hecke eigenvalue $(E_{\bullet},\theta)$ is given by spectral data
$(\sC,\mathfrak{a}(a\sP+b\sP)_{\bullet,\bullet})$ where $\mathfrak{a}
\in \Pic^{0}(\sC)$, and $a, b \in \mathbb{R}^{5}$ are chosen so that
$\sum_{i = 1}^{5} (a_{i}+ b_{i} ) = 3$, and in the case $\lan{G} =
SL_{2}(\mathbb{C})$ aditionally satisfy $a + b \in \mathbb{Z}^{5}$. In
order to use the compatibility with the action of the abelianized
Hecke correspondences we will look for a Hecke eigensheaf whose piece
$(F_{\bullet},\varphi)$ on a particular component of the moduli of
parabolic $G$ bundles is given by modular spectral data
$(Y,\mycal{L}_{\mathfrak{a}}(e E + d (E+G))_{\bullet,\bullet}))$,
where $\mycal{L}_{\mathfrak{a}} =
\bmu^{*}\op{FM}(\mathcal{O}_{\mathfrak{a}}) \in \Pic^{0}(Y)$. In other
words for every parabolic level $T \in \mathbb{R}^{16}$ we will have
\begin{equation}
\label{eq:FT}
 \begin{aligned}
   F_{T} & = f_{*}\mycal{L}_{\mathfrak{a}}(\fl{e + T}E + \fl{d
     +T}(E+G)),   \\
   \varphi & = f_{*}(-\otimes\alpha). 
 \end{aligned}
\end{equation}
The vectors $e, d \in \mathbb{R}^{16}$ need to be chosen so that
$\parch_{1}(F_{\bullet}) = 0$ and $\parch_{2}(F_{\bullet}) = 0$ and so
that all the pieces $(F_{\bullet},\varphi)$ together are a Hecke
eigensheaf with eigenvalue $(E_{\bullet},\theta)$. Note that for the
Hecke eigensheaf property we have to consider all pieces together
since the basic Hecke correspondence relates different connected
components of the moduli of stable parabolic $G = GL_{2}(\mathbb{C})$
or $\mathbb{P}SL_{2}(\mathbb{C})$ bundles.

\index{terms}{Hecke!eigenvalue}
\index{terms}{Hecke!eigensheaf}
\index{terms}{component!connected}
\index{notations}{GL2@$GL_{2}(\mathbb{C})$}
\index{notations}{SL2@$SL_{2}(\mathbb{C})$}
\index{notations}{PSL2@$\mathbb{P}SL_{2}(\mathbb{C})$}

\subsubsection{Hecke kernel} \label{sssec:kernel}

The basic Hecke functor that we need to understand acts on parabolic
Higgs sheaves on the moduli of parabolic bundles with structure group
$G = GL_{2}(\mathbb{C})$ or $G = \mathbb{P}SL_{2}(\mathbb{C})$. The
Higgs version of the Hecke functor is defined as an integral transform
with kernel a parabolic line Higgs bundle
$\left( \boldsymbol{\mathcal{I}}_{\bullet},\bxi\right)$
on a compactified and resolved Hecke correspondence between moduli
spaces.

For $\lan{G} = GL_{2}(\mathbb{C})$ and $G = GL_{2}(\mathbb{C})$ the
compactified and resolved Hecke correspondence is a correspondence
between moduli spaces:
\begin{equation} \label{eq:Heckfr}
\xymatrix@-1pc{
& \bH(\ba,\bb) \ar[dl]_-{p} \ar[dr]^-{q} & \\
\Bun(\ba,\bb) & & \Bun(\ba,\bb)\times C 
}
\end{equation}
Since both the moduli space $\Bun(\ba,\bb)$ and the Hecke
correspondence $\bH(\ba,\bb)$ are disconnected, the kernel object
$\left( \boldsymbol{\mathcal{I}}_{\bullet},\bxi\right)$ is a collection
of different pieces - one on each connected component of
$\bH(\ba,\bb)$.

\index{notations}{Heckabb@$\bH(\ba,\bb)$}
\index{terms}{Hecke!correspondence!compactified and resolved}

In more detail, if $\lan{G} = GL_{2}(\mathbb{C})$ and $G =
GL_{2}(\mathbb{C})$ we have $\bH(\ba,\bb) = \sqcup_{k \in
  \mathbb{Z}} \bH_{k}(\ba,\bb)$, where
\begin{equation} \label{eq:Heckkfr}
\xymatrix@-1pc{
& \bH_{k}(\ba,\bb) \ar[dl]_-{p} \ar[dr]^-{q} & \\
N_{k}(\ba,\bb) & & N_{k+1}(\ba,\bb)\times C 
}
\end{equation}
is a compactification and a resolution of the moduli space
\[
\Heck_{k}^{\op{stable}}(\ba,\bb) \subset N_{k}(\ba,\bb)\times
N_{k+1}(\ba,\bb)\times C
\]
parametrizing triples
$\left((V,\bF,(\ba,\bb)), (V',\bF',(\ba,\bb)),\beta \right)$,
where $(V,\bF,(\ba,\bb))$ and $(V',\bF',(\ba,\bb))$ are stable
parabolic bundles on $(C,\Par_{C})$ with $\deg V = k$, $\deg V' =
k+1$, and $\beta : V \to V'$ is an injective morphism of parabolic
locally free sheaves with $\op{length}(\op{coker}(\beta)) = 1$.

\index{notations}{Heckkabb@$\bH_{k}(\ba,\bb)$}

Suppose $(\ba,\bb)$ are in the dominant chamber. Then for every $k
\in \mathbb{Z}$ the natural isomorphisms $N_{k}(\ba,\bb) \cong
N_{k+1}(\ba,\bb) \cong X$ from section~\ref{sec:family} and the
description of the compactification and resolution from section
\ref{ssec-Hecke_X} identify the diagram \eqref{eq:Heckkfr} with the
diagram
\begin{equation} \label{eq:incidenceHecke}
\xymatrix@-1pc{
& H \ar[dl]_-{p} \ar[dr]^-{q} & \\
X & & X\times C 
}
\end{equation}
where $H = \op{Bl}_{\sqcup \widehat{L_{I}\times
    L_{I}}}\op{Bl}_{\Delta} (X\times X)$. Thus $\bH(\ba,\bb)$ is the
disjoint union of countably many copies of $H$, one for each degree $k
\in \mathbb{Z}$. In particular the Hecke kernel
$(\boldsymbol{\mathcal{I}}_{\bullet},\bxi)$ is a parabolic line Higgs
bundle given by $\sqcup_{k \in \mathbb{Z}}
(\mycal{I}^{k}_{\bullet},\xi)$ where $(\mycal{I}^{k}_{\bullet},\xi^{k})$
is a parabolic line Higgs bundle on $H$ residing on the $k$-th
component of $\bH(\ba,\bb)$.

\index{notations}{H@$H$}
\index{notations}{X@$X$}
\index{notations}{GL2@$GL_{2}(\mathbb{C})$}
\index{notations}{SL2@$SL_{2}(\mathbb{C})$}
\index{notations}{PSL2@$\mathbb{P}SL_{2}(\mathbb{C})$}

The setup for $\lan{G} = SL_{2}(\mathbb{C})$, $G =
\mathbb{P}GL_{2}(\mathbb{C})$ is very similar. In this case the
rigidified moduli stack $\sBun(\ba,\bb)$ of stable parabolic $G$
bundles is a smooth Deligne-Mumford stack which has two connected
components $\sBun(\ba,\bb) = \sM_{0}(\ba,\bb)\sqcup \sM_{1}(\ba,\bb)$ and
the Hecke correspondence
\begin{equation} \label{eq:sHeck.psl2}
\xymatrix@-1pc{
& \sH(\ba,\bb) \ar[dl]_-{p} \ar[dr]^-{q} & \\
\sBun(\ba,\bb) & & \sBun(\ba,\bb)\times C 
 }
\end{equation}
also splits into two disjoint
components:
\begin{equation} \label{eq:Heckabfr.psl2}
  \begin{aligned}
    \left[\begin{minipage}[c]{2.3in}
 $\xymatrix@R-1pc@C-3.5pc{
& \sH(\ba,\bb) \ar[dl]_-{p} \ar[dr]^-{q} & \\
\sBun(\ba,\bb) & & \sBun(\ba,\bb)\times C 
 }$\end{minipage}\right] & = \\[+2pc]
& \hspace{-1.6in}  \left[\begin{minipage}[c]{2.4in}\xymatrix@R-1pc@C-3.5pc{
& \sH_{0}(\ba,\bb) \ar[dl]_-{p} \ar[dr]^-{q} & \\
\sM_{0}(\ba,\bb) & & \sM_{1}(\ba,\bb)\times C 
      }\end{minipage}\right]  \bigsqcup
    \left[\begin{minipage}[c]{2.4in}\xymatrix@R-1pc@C-3.5pc{
& \sH_{1}(\ba,\bb) \ar[dl]_-{p} \ar[dr]^-{q} & \\
   \sM_{1}(\ba,\bb) & & \sM_{0}(\ba,\bb)\times C
 }\end{minipage}\right].
\end{aligned}
\end{equation}
By construction each of the compactified and resolved Hecke
correspondences appearing in the right hand side of
\eqref{eq:Heckabfr.psl2} is isomorphic to the stack quotient of
\eqref{eq:incidenceHecke} by the natural action of
$(\mathbb{Z}/2)^{4}$, i.e. to the diagram
\[
\xymatrix@R-1pc@C-3.5pc{
& \left[H/(\mathbb{Z}/2)^{4}\right] \ar[dl]_-{p} \ar[dr]^-{q} & \\
  \left[X/(\mathbb{Z}/2)^{4}\right] & &
  \left[X/(\mathbb{Z}/2)^{4}\right] \times C 
 }
\]
in which the action of $(\mathbb{Z}/2)^{4}$ on $H = \op{Bl}_{\sqcup
  \widehat{L_{I}\times L_{I}}}\op{Bl}_{\Delta} (X\times X)$ is the one
induced from the diagonal action of $(\mathbb{Z}/2)^{4}$ on $X\times
X$.

Thus for $\lan{G} = SL_{2}(\mathbb{C})$, $G =
\mathbb{P}GL_{2}(\mathbb{C})$ we can view the Hecke kernel
$\left(\boldsymbol{\mathcal{I}}_{\bullet},\bxi\right)$ as a pair
$\left\{ \left(\mycal{I}_{\bullet}^{k},\xi^{k}\right)\right\}_{k \in
  \mathbb{Z}/2}$ of $(\mathbb{Z}/2)^{4}$-equivariant parabolic line
Higgs bundles on \linebreak $(H,\Par_{H} = Lp+Lq+LLC+R_{1}+R_{2})$ which have
trivial first parabolic Chern classes.

\

\noindent
To understand these kernel objects better let us again temporarily
simplify notation and  write $\left(\mycal{I}_{\bullet},\xi\right)$ for any
of $\left(\mycal{I}_{\bullet}^{k},\xi^{k}\right)$ for either $k \in
\mathbb{Z}$ or $k \in \mathbb{Z}/2$.
Via Mochizuki's non-abelian Hodge and extension theorems
\index{terms}{Non-abelian Hodge theory} the parabolic Higgs line
bundle $\left(\mycal{I}_{\bullet},\xi\right)$ should correspond to the
appropriate connected component piece of the Hecke kernel definining
the Hecke functor on the twisted $\mathcal{D}$-modules appearing on
the automorphic side of the tamely ramified GLC (see
Appendix~\ref{app-ram.glc}).

In other words the non-abelian Hodge correspondence should convert
$\left(\mycal{I}_{\bullet},\xi\right)$ to a tamely ramified rank one
parabolic flat bundle $(\mathfrak{I}_{\bullet},\nabla_{\bullet})$ on
$(H,\Par_{H})$. The restriction of $\mathfrak{I}_{0}$ to the open $H -
\Par_{H}$ should be trivial and the Hecke kernel for the tamely
ramified GLC will be a twisted $\mathcal{D}$-module
Deligne-Goresky-MacPherson extension of the holomorphic flat bundle
$(\mathfrak{I}_{0},\nabla_{0})_{|H - \Par_{H}}$ corresponding to the
inclusion $H - \Par_{H} \hookrightarrow \sHeck$ of this open in the
Hecke stack.
\index{terms}{DeligneGG@Deligne-Goresky-MacPherson!extension} The
twisting needed for this Deligne-Goresky-MacPherson extension is
determined by the condition of compatibility with the twisting of the
eigensheaf in the GLC. In our setup we work with eigen twisted
$\mathcal{D}$-modules for which the twistings are real multiples of
the components of the wobbly divisor in $\sBun$
\index{notations}{Buns@$\sBun$} and correspond to the eigenvalues of
the residues of the flat bundle with regular singularities defined on
the complement of the wobbly locus. Furthermore we restricted to the
purely imaginary case of Mochizuki's non-abelian Hodge correspondence
and for simplicity chose to work with the closure of the symplectic
leaf selected on the Higgs side by requiring nilpotent residues.
\index{terms}{purely imaginary condition}
\index{terms}{conversion table!of the ramified nah}
By the conversion table of
\cite{carlos-nc} or \cite{dps} this means that on the flat bundle side
the parabolic weights are equal to the negative eigenvalues of the
residues and that the parabolic weights on the flat bundle side are
equal to the parabolic weights on the Higgs bundle side. In particular
the twistings for the GLC Hecke kernel will be equal to the residues
of $\nabla_{0}$, or equivalently to the negative of the parabolic
weights for $\mathfrak{I}_{\bullet}$, or equivalently equal to the
negative of the parabolic weights for $\mycal{I}_{\bullet}$. Also, the
logarithmic one form $\xi$ in the Higgs bundle
$\left(\mycal{I}_{\bullet},\xi\right)$  has
to have zero residues. Since by construcion $H$ is a smooth rational
variety, $H$ can not have any global holomorphic one forms and so we
must have that $\xi$ is identically equal to zero.

This implies that the Hecke functor on parabolic Higgs sheaves must be
given by a kernel which is a parabolic rank one Higgs bundle
\begin{equation} \label{eq:kernelfinal}
\boxed{
  \left(\mycal{I}_{\bullet},\xi\right) =
\left(\mathcal{O}(\bzeta\Par_{H})_{\bullet}, 0\right) 
}
\end{equation}
where
\begin{equation} \label{eq:zeta.basic}
  \bzeta = (lp,lq,llc,r_{1},r_{2}) \in \mathbb{R}^{16}\times
  \mathbb{R}^{16}\times \mathbb{R}^{16}\times \mathbb{R}^{5}\times
  \mathbb{R}^{5}
\end{equation}
is chosen so that $\parch_{1}(\mathcal{O}(\bzeta\Par_{H})_{\bullet}) =
[\bzeta\Par_{H}] = 0$ in $H^{2}(H,\mathbb{R})$ and so that for the
given $\bzeta$ the parabolic Hecke eigensheaf property holds for all
possible eigenvalue Higgs bundles on $(C,\Par_{C})$. 

\subsection{The eigensheaf property} \label{ssec:eigenproperty}

The parabolic Hecke eigensheaf condition is the non-abelian Hodge
translation of the ordinary Hecke eigensheaf condition that one
imposes on the automorphic side of GLC.  Explicitly we have

\index{terms}{Non-abelian Hodge theory}
\index{terms}{Hecke!eigensheaf}

\

\medskip

\begin{problem} \label{problem:Hecke.eigensheaf} Show that there exists a
natural $\bzeta = (lp,lq,llc,r_{1},r_{2})$ so that the Hecke operator
determined by the kernel
$\left(\mathcal{O}(\bzeta\Par_{H})_{\bullet},0\right)$ on $H$ has a
unique irreducible eigensheaf for any given eigenvalue. Specifically,
show that for any choice of $(a,b) \in \mathbb{R}^{5}\times
\mathbb{R}^{5}$ we can choose an apporpriate $\bzeta$ with
$\left[\bzeta \Par_{H}\right] = 0$ in $H^{2}(H,\mathbb{R})$ and so that
for any stable $(E_{\bullet},\theta)$ on $(C,\Par_{H})$ with parabolic
weights determined by $(a,b)$, we can find stable parabolic
Higgs bundles $(F_{\bullet},\varphi)$ and $(F_{\bullet}',\varphi')$ on
$(X,\Par_{X})$ so that the folllowing {\em\bfseries parabolic Hecke
  eigensheaf condition} holds
\begin{equation} \label{eq-par.hecke}
\boxed{
q_{*}\left(p^{*}(F_{\bullet},\varphi)\otimes
(\mathcal{O}(\bzeta\Par_{H})_{\bullet},0)\right) =
p_{X}^{*}(F_{\bullet}',\varphi')\otimes p_{C}^{*}(E_{\bullet},\theta).}
\end{equation}
\end{problem}

\

\bigskip

\begin{rem} \label{rem:basic.gives.all}
A solution to Problem~\ref{problem:Hecke.eigensheaf} will provide the
crucial basic building block for solving the full Hecke eigensheaf
problem for parabolic Higgs bundles for $\lan{G} = GL_{2}(\mathbb{C})$
or $\lan{G} = SL_{2}(\mathbb{C})$. Indeed, when $\lan{G} =
GL_{2}(\mathbb{C})$ solving the Hecke kernel and Hecke eigensheaf
problem on $\Bun(\ba,\bb)$ amounts to solving
Problem~\ref{problem:Hecke.eigensheaf} iteratively once for each $k
\in \mathbb{Z}$ and producing collections $\sqcup_{k \in \mathbb{Z}}
(\mathcal{O}(\bzeta^{k}\Par_{H})_{\bullet},0)$ and $\sqcup_{k \in
  \mathbb{Z}} (F_{\bullet}^{k},\varphi^{k})$ such that
\[
q_{*}\left(p^{*}(F_{\bullet}^{k},\varphi^{k})\otimes
(\mathcal{O}(\bzeta^{k}\Par_{H})_{\bullet},0)\right) =
p_{X}^{*}(F_{\bullet}^{k+1},\varphi^{k+1})\otimes p_{C}^{*}(E_{\bullet},\theta).
\]
Similarly, when $\lan{G} = SL_{2}(\mathbb{C})$ solving the Hecke
kernel and Hecke eigensheaf problem on $\sBun(\ba,\bb)$ amounts to
solving Problem~\ref{problem:Hecke.eigensheaf} twice producing
$(\mathcal{O}(\bzeta\Par_{H}),0)\sqcup
(\mathcal{O}(\bzeta'\Par_{H})_{\bullet},0)$ and a
$(\mathbb{Z}/2)^{4}$-equivariant pair $(F_{\bullet},\varphi)\sqcup
(F_{\bullet}',\varphi')$ satisfying
\[
\begin{aligned}
q_{*}\left(p^{*}(F_{\bullet},\varphi)\otimes
(\mathcal{O}(\bzeta\Par_{H})_{\bullet},0)\right) & =
p_{X}^{*}(F_{\bullet}',\varphi')\otimes p_{C}^{*}(E_{\bullet},\theta), \\[0.5pc]
q_{*}\left(p^{*}(F_{\bullet}',\varphi')\otimes
(\mathcal{O}(\bzeta'\Par_{H})_{\bullet},0)\right) & =
p_{X}^{*}(F_{\bullet},\varphi)\otimes p_{C}^{*}(E_{\bullet},\theta).
\end{aligned}
\]
\end{rem}

\

\bigskip

\noindent
The pullbacks, pushforwards, and tensor products in this formula are
all taken in the dg category of semistable parabolic Higgs bundles
with vanishing parabolic Chern classes (see
section~\ref{ssec:operations}).  \index{terms}{parabolic!Chern
  classes} That is, we use the pullback, pushforward, and tensoring
operations that are intertwined with the respective operations on
twisted $\mathcal{D}$-modules via Mochizuki's non-abelian Hodge and
extension theorems. In particular the pushforward is the one defined
via the algebraic $L^{2}$-Dolbeault complex of a Higgs bundle as in
\cite{dps}. We will review all these operations next when we use
abelianization to convert the Hodge and Hecke conditions to numerical
equations. This is a long calculation using a sequence of geometric
reductions, and will take the remainder of this section and the  entirety
of Section~\ref{ssec:ab.eigenproperty}. The final result of the
calculation is the abelianized Hecke condition which reduces
Problem~\ref{problem:Hecke.eigensheaf} to an equation of line bundles
and is summarized in Section~\ref{sssec:ab.Hecke.conditions}.

\index{terms}{Non-abelian Hodge theory}

Since $\Par{X\times C} = \sum_{I \in \oddL} L_{I}\times C + \sum_{i
  =1}^{5} X\times p_{i}$. a parabolic level on $(X\times C,
\Par_{X\times C})$ is labeled by a vector $(T,t) \in
\mathbb{R}^{16}\times \mathbb{R}^{5}$. Both sides of the condition
\eqref{eq-par.hecke} are parabolic tame Higgs bundles on $X\times
C$. We will write $\lhs_{T,t}$ and $\rhs_{T,t}$ for the corresponding
meromorphic Higgs bundles at level $(T,t)$. With this notation we have
\[
\rhs_{T,t} = \left( p_{X}^{*}(F_{\bullet}',\varphi')\otimes
p_{C}^{*}(E_{\bullet},\theta)\right)_{T,t}
\]
Taking into account that $\Par_{X}\times C = p_{X}^{*}\Par_{X}$ and
$X\times \Par_{C} = p_{C}^{*}\Par_{C}$  and applying 
the pullback formulas in Section~\ref{ssec:operations} we have
\[
\begin{aligned}
  (p_{X}^{*}F_{\bullet}')_{T,t} & = p_{X}^{*}F'_{T}\otimes
  \mathcal{O}\left(\fl{t}(X\times \Par_{C})\right), \\
  (p_{C}^{*}E_{\bullet})_{T,t}
  & = \mathcal{O}\left(\fl{T} (\Par_{X}\times C)\right)\otimes p_{C}^{*}E_{t}.
\end{aligned}
\]
Furthermore by the tensor product formulas in
Section~\ref{ssec:operations}  we get
\[
\begin{aligned}
  (p_{X}^{*}F_{\bullet}' & \otimes p_{C}^{*} E_{\bullet})_{T,t}
  = \bigcup_{(S'+S'',s'+s'') \leq (T,t)}
  \left(p_{X}^{*}F_{\bullet}'\right)_{(S',s')}\bigotimes 
  \left( p_{C}^{*} E_{\bullet}\right)_{(S'',s'')} \\[+0.5pc]
  & = \left(\bigcup_{S'+S''\leq T }
  p_{X}^{*}F_{S'}'\otimes
  \mathcal{O}\left(\fl{S''}(\Par_{X}\times C)\right)\right)
  \bigotimes
  \left(\bigcup_{s'+s" \leq t} p_{C}^{*}E_{s''}\otimes
  \mathcal{O}\left((\fl{s'}(X\times \Par_{C})\right)\right) \\[+0.5pc]
  & = p_{X}^{*} F_{T}' \otimes p_{C}^{*} E_{t}.
\end{aligned}
\]
So altogether we get
\begin{equation} \label{eq:rhs}
\boxed{
  \rhs_{T,t} = \left( p_{X}^{*} F_{T}' \otimes p_{C}^{*} E_{t},
  p_{X}^{*}\varphi\otimes \mathbf{1} +
  \mathbf{1}\otimes p_{C}^{*}\theta \right)
 }
\end{equation}
Similarly on the left hand side we need to understand the complex of
meromorphic Higgs sheaves 
\[
\lhs_{T,t} = \left( q_{*} \left(p^{*}(F_{\bullet},\varphi)\otimes
  \left(\mathcal{O}(\bzeta\Par_{H})_{\bullet},0\right)
\right)\right)_{T,t}.
\]
The first step is to describe the parabolic bundle $p^{*}F_{\bullet}$
on $(H,\Par_{H})$. Since $\Par_{H} = Lp + Lq + LLC + R_{1} + R_{2}$ we
can label the parabolic level  on $H$ by a five tuple of vectors
\[
  \left(T_{Lp},T_{Lq},T_{LLC},t_{R_{1}},t_{R_{2}}\right) \in
  \mathbb{R}^{16}\times \mathbb{R}^{16}\times \mathbb{R}^{16}\times
  \mathbb{R}^{5}\times \mathbb{R}^{5}.
\]
Since for all $I \in \oddL$ we have $p^{*}L_{I} = Lp_{I} + LLC_{I}$,
the definition of a parabolic pullback in
Section~\ref{ssec:operations} gives  
\[
  (p^{*}F_{\bullet})_{T,T_{Lq},T,t_{R_{1}},t_{R_{2}}} =
  p^{*}F_{T}\otimes \mathcal{O}(\fl{T_{Lq}}Lq +
  \fl{t_{R_{1}}}R_{1} + \fl{t_{R_{2}}}R_{2}),
\]
for every $T \in
\mathbb{R}^{16}$. Thus 
for a level $\left(T_{Lp},T_{Lq},T_{LLC},t_{R_{1}},t_{R_{2}}\right)$
we get
\[
(p^{*}F_{\bullet})_{T_{Lp},T_{Lq},T_{LLC},t_{R_{1}},t_{R_{2}}} =
\left(p^{*}F_{T_{Lp}}(\infty LLC)\cap p^{*}F_{T_{LLC}}\right)
(\infty Lp)\otimes \mathcal{O}(\fl{T_{Lq}}Lq +
  \fl{t_{R_{1}}}R_{1} + \fl{t_{R_{2}}}R_{2}).
\]
Next recall that the general parabolic tensor product formula from
Section~\ref{ssec:operations} simplifies when one of the bundles in
the tensor product is a parabolic line bundle. If
$\mycal{F}_{\bullet}$ is a parabolic bundle on $(Z,D)$ and
$\mycal{E}_{\bullet} = \mycal{L}(dD)_{\bullet}$ is a parabolic line
bundle on $(Z,D)$, then the level $w$ of the tensor product
$\mycal{F}_{\bullet}\otimes \mycal{E}_{\bullet}$ is given simply by
\[
\left(\mycal{F}_{\bullet}\otimes \mycal{E}_{\bullet}\right)_{w} =
\mycal{F}_{w+d}\otimes \mycal{L}. 
\]
Applied to our setting this remark implies that
\[
  \begin{aligned}
  \big(p^{*}F_{\bullet}\otimes &
    \mathcal{O}(\bzeta\Par_{H})_{\bullet}
  \big)_{T_{Lp},T_{Lq},T_{LLC},t_{R_{1}},t_{R_{2}}} = \\
& = \left(p^{*}F_{T_{Lp}+lp}(\infty LLC)\cap p^{*}F_{T_{LLC}+llc}
  (\infty Lp)\right) \\
& \qquad\quad \otimes \mathcal{O}(\fl{T_{Lq+lq}}Lq +
\fl{t_{R_{1}}+r_{1}}R_{1} + \fl{t_{R_{2}}+r_{2}}R_{2}).
\end{aligned}
\]
The last step is to compute the parabolic Higgs bundle pushforward of 
\[
  (U_{\bullet},\psi) :=
  \left(p^{*}F_{\bullet}\otimes \mathcal{O}(\bzeta\Par_{H})_{\bullet},
    p^{*}\varphi\otimes \mathbf{1}\right)
\]
via the map $q : H \to
X\times C$.

As explained in Section~\ref{ssec:operations} this amounts to
computing the sheaf theoretic hyperderived image  of an appropriately chosen
parabolic Dolbeault complex for $(U_{\bullet},\psi)$. Concretely
note that $q : H \to X\times C$ is a semistable proper flat morphism of relative
dimension one, and $q^{*}(Par_{X}\times C) = Lq+LLC$, and
$q^{*}(X\times \Par_{C})  = R_{1} + R_{2}$. Therefore if 
we write $\Par_{H} = \Par_{H}^{\hor}
+ \Par_{H}^{\ver}$ for the horizontal and vertical parts of the
parabolic divisor on $H$ we have
\[
  \begin{aligned}
    \Par_{H}^{\hor} & = Lp, \\
    \Par_{H}^{\ver} & = Lq+LLC + R_{1} + R_{2} = q^{*}(\Par_{X\times C}).
  \end{aligned}
\]
As usual (see Section~\ref{ssec:operations}) given a parabolic level $(T,t)$ on
$(X\times C,\Par_{X\times C})$ let us  write $\lift(T,t)$ for the
parabolic level on $(H,\Par_{H})$ which assigns $0$ to each horizontal
component of $\Par_{H}$, while to a vertical component of $\Par_{H}$
it assigns the value of $(T,t)$ on the image of that vertical
component. Thus we have
\[
\lift(T,t) = (T_{Lp},T_{Lq},T_{LLC},t_{R_{1}},t_{R_{2}}) = (0,T,T,t,t).
\]
The level $(T,t)$ of the $q$-relative parabolic Dolbeault complex of
$(U_{\bullet},\psi)$ is the two term complex
\[
\Dol(q,U_{\lift(T,t)}) = 
\xymatrix@R-2pc{
\Big[ W_{0}(\hor,U_{\lift(T,t)}) 
\ar[r]^-{\verpsi} & 
W_{-2}(\hor,U_{\lift(T,t)})
\otimes \omega_{q}(Lp)\Big]. \\
0 & 1 }
\]
Here $\verpsi$ is the composition of $\psi$ and the
projection onto the relative logarithmic forms, while $\omega_{q}$
denotes the relative dualizing sheaf of $q$. In the formula for
$\Dol(q,U_{\lift(T,t)})$ we have tacitly used the fact that because
$q$ is semistable and of relative dimension one, there are natural
identifications
\[
\begin{aligned}
\Omega^{1}_{H/X\times C}(\log \Par_{H}) & =
\Omega^{1}_{H}(\log \Par_{H})/q^{*}\Omega^{1}_{X\times
  C}(\log \Par_{X\times C}) \\
& = \omega_{q}(\Par_{H}^{\hor}) \\
& = \omega_{q}(Lp).
\end{aligned}
\]

Similarly, the level $(T,t)$ of the absolute parabolic Dolbeault
complex is
\[
\Dol(H,U_{\lift(T,t)}) = \left[ \begin{array}{c}
    W_{0}(\hor,U_{\lift(T,t)}) \\[+0.3pc]
    \downarrow \wedge \psi \\[+0.3pc]
    W_{-2,0}\left(\hor, U_{\lift(T,t)}\otimes \Omega^{1}_{H}(\log \Par_{H}) \right) \\[+0.3pc]
    \downarrow \wedge \psi \\[+0.3pc]
W_{-2,0}\left(\hor, U_{\lift(T,t)}\otimes \Omega^{2}_{H}(\log \Par_{H}) \right) \\[+0.3pc]
                                  \downarrow \wedge \psi \\[+0.3pc]
W_{-2,0}\left(\hor, U_{\lift(T,t)}\otimes \Omega^{3}_{H}(\log \Par_{H}) \right) \\[+0.3pc]
                                  \downarrow \wedge \psi \\[+0.3pc]
W_{-2,0}\left(\hor, U_{\lift(T,t)}\otimes \Omega^{4}_{H}(\log \Par_{H}) \right)
\end{array}\right] \  \begin{array}{c}
     0 \\[+0.3pc]
    \\[+0.3pc]
    1 \\[+0.3pc]
    \\[+0.3pc]
    2 \\[+0.3pc]
    \\[+0.3pc]
    3  \\[+0.3pc]
    \\[+0.3pc]
    4
\end{array}
\]
These Dolbeault complexes fit in a short exact sequence of complexes
\[
 \xymatrix{
   0 \ar[d] \\
   \Dol(q,U_{\lift(T,t)})[-1]\otimes q^{*}\Omega^{1}_{X\times
  C}(\log \Par_{X\times C}) \ar[d] \\
\Dol(H,U_{\lift(T,t)})/I^{2}(U_{\lift(T,t)}) \ar[d] \\
\Dol(q,U_{\lift(T,t)}) \ar[d] \\
0
}
\]
which can be viewed as a morphism
\[
  \mathfrak{d}(\psi) : \Dol(q,U_{\lift(T,t)}) \longrightarrow
  \Dol(q,U_{\lift(T,t)})\otimes q^{*}\Omega^{1}_{X\times
    C}(\log \Par_{X\times C})
\]
in the derived category $D^{b}_{\op{coh}}(H,\mathcal{O})$.

The $(T,t)$-level of the Higgs pushforward of $(U_{\bullet},\psi)$
is then given by the sheaf theoretic pushforward of
$\Dol(q,U_{\lift(T,t)})$ and $\mathfrak{d}(\psi)$:
\[
  \left(q_{*}(U_{\bullet},\psi)\right)_{T,t} = \left(q_{*}\Dol(q,U_{\lift(T,t)}),q_{*}
    \mathfrak{d}(\psi)\right).
\]
To streamline the calculation of this pushforward it will be useful to
introduce the following notation. Given a parabolic level $(T,t)$ on
$(X\times C,\Par_{X\times C})$ we fix $m_{I} \geq
\max(T_{I}+llc_{I},lp_{I}) + 1$ for all $I \in \oddL$ and we  set
\begin{equation}
\label{eq:btau}
\begin{aligned}
\btau_{T,t} & = \mathcal{O}_{H}(\fl{T+lq}Lq + 
\fl{t+r_{1}}R_{1} + \fl{t+r_{2}}R_{2}), \\
\bbF_{T,t} & = p^{*}F_{m}\otimes \btau_{T,t}, \\
\bbF'_{T,t} & = p^{*}F_{T+llc}\left( \fl{m-T-llc}Lp\right) \otimes
\btau_{T,t}, \\
\bbF''_{T,t} & = p^{*}F_{lp}\left( \fl{m - lp}LLC\right)\otimes
\btau_{T,t}.
\end{aligned}
\end{equation}
Then by definition we have
\[
U_{\lift(T,t)} = U_{0,T,T,t,t} = \bbF'_{T,t}\cap
\bbF''_{T,t} \ \text{inside} \ \bbF_{T,t}.
\]
Altogether this gives
\begin{equation} \label{eq:lhs}
\boxed{
\lhs_{T,t} = \left( q_{*}\Dol(q,\bbF'_{T,t}\cap
\bbF''_{T,t}), q_{*}\mathfrak{d}(\psi)\right).
}
\end{equation}
To compute the steps of the weight filtration appearing in the
definition of this Dolbeault complex we will also need to pass to the
associated graded of the parabolic filtration on $U_{\bullet}$ along
the $q$-horizontal divisor $Lp$. For this we need to  choose small positive
real numbers $\varepsilon_{I} > 0$ for all $I \in \oddL$ and
consider the subsheaves
\[
  U_{\lift(T,t)}(-\varepsilon) =
  U_{-\varepsilon,T,T,t,t} \subset U_{0,T,T,t,t} = U_{\lift(T,t)}
\]
and
the action of the residue of $\psi$ on the  quotient sheaves 
\[
U_{\lift(T,t)}/U_{\lift(T,t)}(-\varepsilon) = U_{0,T,T,t,t}/U_{-\varepsilon,T,T,t,t}.
\]
To compute $U_{\lift{T,t}}(-\varepsilon)$ note that by definition small
decreases in the $Lp$ parabolic level do not affect the family
$\bbF'_{T,t}$ and so we only need to trace how subtracting
$\varepsilon$ from the $Lp$ parabolic level will affect
$\bbF''$. In particular setting
\[
  \bbF''_{T,t}(-\varepsilon) =
  p^{*}F_{-\varepsilon + lp}\left( \fl{m - lp}LLC\right)\otimes
\btau_{T,t}
\]
gives the identification
\[
U_{\lift(T,t)}(-\varepsilon) = \bbF'_{T,t}\cap \bbF''_{T,t}(-\varepsilon)
\quad \text{inside \ } \bbF_{T,t}.
\]

\subsection{Abelianization} \label{ssec:ab.eigenproperty}

\

\noindent
Abelianization allows us to express non-abelian data on $C$ and $X$ in
terms of equivalent abelian data on their sepctral covers $\sC$ and
$Y$, respectively.  Since $(E_{\bullet},\theta)$ is described by
spectral data on $\sC$ and $(F_{\bullet},\varphi)$ and
$(F_{\bullet}',\varphi')$ are described by spectral data on $Y$, we
can rewrite all our conditions on $Y\times \sC$.  We start here by
doing this for the Hecke condition \eqref{eq-par.hecke}.

\

\smallskip

\subsubsection{Abelianization of the right hand side of
  \eqref{eq-par.hecke}} \

As a first step we can rewrite the right hand side
\eqref{eq:rhs} as data on $Y\times \sC$. This is straightforward.
First note that we have a commutative diagram
\[
\xymatrix@R-0.5pc@C-0.5pc{ & Y\times\sC \ar[dl]_-{p_{Y}} \ar[dr]^-{p_{\sC}}
  \ar[dd]^-{f\times \pi} & \\
Y \ar[dd]_-{f} & & \sC \ar[dd]^-{\pi} \\
& X\times C \ar[dl]_-{p_{X}} \ar[dr]^-{p_{C}}& \\
X & & C
}
\]
and these are fiber squares. As explained in sections
\ref{sssec:eigenvalue} and \ref{sssec:eigensheaf} for any level $(T,t)$
we have
\[
\begin{aligned}
  (E_{t},\theta) & = \left(\pi_{*}\mathfrak{a}_{t,t},
  \pi_{*}(-\otimes \lambda)\right) \\
& = 
  \left(\pi_{*}
  \mathfrak{a}((\fl{t+a}+\fl{t+b})\widetilde{P}),
  \pi_{*}(-\otimes \lambda)\right), \\
  (F_{T}',\varphi') & = \left(f_{*}M_{T,T}', f_{*}(-\otimes \alpha)\right)  \\
  & = \left(f_{*}\mycal{L}_{\mathfrak{a}}((\fl{T+e'} + \fl{T+d'})E +
  \fl{T+d'}G),f_{*}(-\otimes \alpha)\right).
\end{aligned}
\]
Combined with the formula \eqref{eq:rhs} this gives
\[
\rhs_{T,t} = \Big( (f\times \pi)_{*}(p_{Y}^{*}M_{T,T}'\otimes
p_{\sC}^{*}\mathfrak{a}_{t,t}),
(f\times
  \pi)_{*}(p_{Y}^{*}\alpha\otimes \mathbf{1} + \mathbf{1}\otimes
  p_{\sC}^{*}\lambda)\Big),
  \]
or more explicitly after expanding our shortcut notation:

\

\begin{equation} \label{eq:rhs-ab}
  \boxed{
    \begin{aligned}
      \rhs_{T,t} & = \Big((f\times \pi)_{*}(p_{Y}^{*}
      \mycal{L}_{\mathfrak{a}}((\fl{T+e'}
+ \fl{T+d'})E + \fl{T+d'}G) \\
& \hspace{6pc} \otimes p_{\sC}^{*}
  \mathfrak{a}((\fl{t+a}+\fl{t+b})\widetilde{P}), \\
  & \hspace{11pc} (f\times
  \pi)_{*}(p_{Y}^{*}\alpha\otimes \mathbf{1} + \mathbf{1}\otimes
  p_{\sC}^{*}\lambda)\Big),
\end{aligned}
}
\end{equation}
\

\smallskip

\subsubsection{Abelianization of the left hand side of
  \eqref{eq-par.hecke}} \ \label{sssec:ab.lhs}

Next we need to rewrite $\lhs_{T,t}$ in a similar manner. To do this
efficiently we need to first rewrite the parabolic Higgs bundle
$(U_{\bullet},\psi)$ on $H$ in terms of spectral data.

Since up to tensoring with a parabolic line bundle $U_{\bullet}$ is a
pull back of a parabolic Higgs bundle on $X$ whose spectral cover is
$f : Y \to X$, it follows that the spectral cover $f_{Z} : Z \to H$
for $(U_{\bullet},\psi)$ will be the fiber product of $Y$ and $H$:
\[
\xymatrix@-0.5pc{
  Z \ar[r]^-{f_{Z}} \ar[d]_-{p_{Z}} & H \ar[d]^-{p} \\
Y \ar[r]_-{f} & X
}
\]
By the same token the corresponding spectral sheaf on $Z$ will be the
tensor product of a pullback of a parabolic line bundle on $H$ with
the pullback with the spectral parabolic line bundle on
$Y$. Concretely, for a parabolic level $(T,t)$ on $(X\times
C,\Par_{X\times C})$ we fix \linebreak 
$m_{I} \geq \max(T_{I}+llc_{I},lp_{I}) +
1$ and small $\varepsilon_{I} > 0$ for all $I \in \oddL$ and using the
previously defined  line bundles $M_{T,T}$  (see \eqref{eq:MST})
and $\btau_{T,t}$ (see \eqref{eq:btau}) we set
\begin{equation} \label{eq:frakM}
\begin{aligned}
\bQ_{T,t} & = p_{Z}^{*}M_{m,m}\otimes f_{Z}^{*}\btau_{T,t}, \\
\bQ'_{T,t} & = p_{Z}^{*}M_{T+llc,T+llc}\left( \fl{m-T-llc}f_{Z}^{*}Lp\right) \otimes
f_{Z}^{*}\btau_{T,t}, \\
\bQ''_{T,t}(-\varepsilon)  & = p_{Z}^{*}M_{-\varepsilon+lp,\varepsilon+lp}
\left( \fl{m - lp}f_{Z}^{*}LLC\right)\otimes
f_{Z}^{*}\btau_{T,t}, \\
\bQ_{\lift(T,t)}(-\varepsilon)  & = \bQ'_{T,t}\cap \bQ''_{T,t}(-\varepsilon)
\quad \text{inside \ }  \bQ_{T,t}.
\end{aligned}
\end{equation}
With this notation we now have
\begin{equation} \label{eq:UupviaM}
  U_{\lift(T,t)}(-\varepsilon) =
  U_{-\varepsilon,T,T,t,t} = f_{Z*}\bQ_{\lift(T,t)}(-\varepsilon) .
\end{equation}
Furthermore, since $Z = Y\times_{X} H$ and  $p^{*}L_{I} = Lp_{I} +
LLC_{I}$ and $f^{*}L_{I} = 2E_{I} + G_{I}$ it follows that
$p_{Z}^{*}E_{I}$ and $p_{Z}^{*}G_{I}$ each have two components which
we will denote by $LpE_{I}$ and $LLE_{I}$, and $LpG_{I}$ and $LLG_{I}$
respectively. Explicitly these are Weil divisors on $Z$ defined by 
\begin{equation} \label{eq:LpEetc}
\begin{aligned}
LpE_{I} & = E_{I}\times_{L_{I}} Lp_{I}, \\
LLE_{I} & = E_{I}\times_{L_{I}} LLC_{I}, \\ 
LpG_{I} & = G_{I}\times_{L_{I}} Lp_{I}, \\
LLG_{I} & = G_{I}\times_{L_{I}} LLC_{I},
\end{aligned}
\end{equation}
 and we have 
\[
\begin{aligned}
p_{Z}^{*}E & = LpE + LLE, \\
p_{Z}^{*}G & = LpG + LLG.
\end{aligned}
\]
Tracing through the definitions we see that in terms of these divisors 
$\bQ_{\lift(T,t)}(-\varepsilon)$ is the torsion free rank one divisorial
sheaf on $Z$ given by the formula
\begin{equation} \label{eq:spectral.sheaf}
\boxed{
\begin{aligned}
\bQ_{\lift(T,t)}(-\varepsilon) & = p_{Z}^{*}\mycal{L}_{\mathfrak{a}} \otimes
\mathcal{O}_{Z}\Big(
\left(\fl{-\varepsilon + lp + e} + \fl{-\varepsilon + lp +
    d}\right)LpE \\
& \qquad\quad  + 
\left(\fl{T + llc + e} + \fl{T + llc +
    d}\right)LLE \\
& \qquad\quad  +
\fl{-\varepsilon + lp + d}LpG
+ \fl{T + llc + d}LLG
\Big)\otimes f_{Z}^{*}\btau_{T,t}.
\end{aligned}
}
\end{equation}

\

\smallskip

\

\noindent
The spectral data $(Z,\bQ_{\lift(T,t)}(-\varepsilon))$ for
$(U_{\lift(T,t)}(-\varepsilon),\psi)$ should be treated somewhat
carefully since $Z$ is not locally factorial in codimension two and
the spectral sheaf $\bQ_{\lift(T,t)}(-\varepsilon)$ will not in general
be a line bundle.  For the calculation of the Dolbeault complexes we
will need to analyze this in more detail.

\

\begin{lem} \label{lem:sing.of.Z} \ {\bfseries (a)} \ The singular locus
  $\op{Sing}(Z)$ of $Z$ is the union of sixteen disjoint surfaces
 $\Xi_{I}$,  $I \in \oddL$ each of which is
 naturally isomorphic to $C\times C = \mathbb{P}^{1}\times \mathbb{P}^{1}$.
Each $\Xi_{I}$ contains a marked ruling $\rho_{I} =
 C\times p_{6}$ and the singularity type of $Z$ is locally constant along the
 strata of the stratification $\cup_{I} \rho_{I} \subset \cup_{I}
 \Xi_{I} \subset Z$.

 \

 \noindent
 {\bfseries (b)} \ Near each $\Xi_{I}$ the space $Z$ is isomorphic to
 a family of surface singularities parametrized by $\Xi_{I}$, with
 fibers of type $A_{1}$ at points in $\Xi_{I} - \rho_{I}$ and fibers
 of type $A_{3}$ near points in $\rho_{I}$.

 \

 \noindent
 {\bfseries (c)} \ In $Z - \cup_{I} \rho_{I}$ the divisors $LpE_{I}$,
 $LLE_{I}$
 are not Cartier and are transversally represented by a ruling.
 \end{lem}

 \

 \noindent
 {\bfseries Proof.} \ The fiber product $Z = Y\times_{X} H$ is smooth
 at points where at least one of the maps $f$ and $p$ is
 submersive. In other words we have
 $\op{Sing}(Z) \subset \op{crit}(f)\times_{X} \op{crit}(p)$.

 The blowup description of $H$ in Section~\ref{ssec-Hecke_X} implies
 that the discriminant of the projection $p : H \to X$ is the wobbly
 divisor in $X$: $\op{Discr}(p) = \bWob = \cup_{I} L_{I}$. Furthermore
 for all $I \in \oddL$ we have $p^{*}L_{I} = Lp_{I} + LLC_{I}$ and the
 part of the critical locus of $p$ mapping to the line $L_{I}$ is
 scheme theoretically identified with the
 smooth surface
\[
\op{crit}(p)_{I} = Lp_{I} \cap LLC_{I}.
\]
The blowdown map $H \to X\times X$ maps $\op{crit}(p)_{I}$
isomorphically to $L_{I}\times L_{I} \subset X\times X$.  But recall that the line
  $L_{I}$ meets exactly five other lines, namely the lines $L_{I+i}$
  for $i = 1, \ldots 5$. As explained in Section~\ref{ssec-Hecke_X}
  the cross-ratios of the intersection points
  $p_{I,I+i} := L_{I}\cap L_{I+i}$ are equal to the crooss-ratios of
  the points $p_{i}$ in $C$ and so there is a unique isomorphism
  $L_{I} \cong C$ under which $p_{I,I+i}$ is mapped to $p_{i}$. Using
  this isomorphism we obtain a canonical identification of
$\op{crit}(p)_{I}$ with $C\times C$, so that the map $p :
\op{crit}(p)_{I} \to L_{I}$ corresponds to the projection of $C\times
C$ onto the first factor.

In Proposition~\ref{prop:Ybar} we saw that the four sheeted modular
spectral cover $f : Y \to X$ arises naturally as an iterated double
cover $Y \to \overline{Y} \to X$, where $\overline{Y}$ is the Kummer
$K3$ surface for the Jacobian of $\sC$. In fact 
\begin{itemize}
\item $\overline{Y} \to X$ is the double cover branched along the
  unique smooth bi-anticanonical curve $Q_{\sC}$ in $X$ which is
  tangent to each line $L_{I}$ at the point $p_{I,I} \in L_{I}$
  corresponding to $p_{6} \in C$ under the isomorphism
  $L_{I} \cong C$.
\item The preimage of each $L_{I}$ in $\overline{Y}$ is the union
  $\overline{E}_{I}\cup \overline{G}_{I}$ of 
  two smooth genus zero curves meeting at a point and $Y \to
  \overline{Y}$ is the double cover branced at the divisor $\sum_{I}
  \overline{E}_{I}$. 
\end{itemize}
Let $\overline{\op{Ram}}_{\sC} \subset \overline{Y}$ be the
ramification divisor for the double cover $\overline{Y} \to Y$, and
let $\op{Ram}_{\sC} \subset Y$ be the preimage of
$\overline{\op{Ram}}_{\sC}$. Then from the iterated double cover
description we get that
\[
\op{crit}(f) = \left(\bigcup_{I\in \oddL} E_{I}\right)\cup \op{Ram}_{\sC}
\]
and so scheme theoretically 
\[
\op{crit}(f) \times_{X} \op{crit}(p) = (E_{I}\cup
\op{Ram}_{\sC})\times_{X} \op{crit}(p)_{I}  = (E_{I}\cup
\op{Ram}_{\sC})\times_{L_{I}\cup Q_{\sC}} \op{crit}(p)_{I}.
\]
Let $\Xi_{I}$ denote the surface $E_{I}\times_{L_{I}}\op{crit}(p)_{I}
= \op{crit}(p)_{I} \cong L_{I}\times L_{I}$ and let $\rho_{I} \subset
\Xi_{I}$ be the preimage of the point $p_{I} \in L_{I}$. Then as a
scheme $\op{crit}(f) \times_{X} \op{crit}(p)$ is the disjoint union of
sixteen components $\widetilde{\Xi}_{I}$ where
$(\widetilde{\Xi}_{I})_{\op{red}} = \Xi_{I}$ and $\widetilde{\Xi}_{I}$
has an embedded length two nilpotent structure along $\rho_{I} \subset
\Xi_{I}$.

This shows that the singiularity type of $Z$ is locally constant along
the strata of the stratification induced by the inclusions
$\cup_{I} \rho_{I} \subset \cup_{I} \Xi_{I} \subset Z$ and essentially
proves {\bfseries (a)}. The part of {\bfseries (a)} that is still
missing is the statement that the singularity type actually jumps from
stratum to stratum. This is a local statement which again can be
deduced from the iterated double cover  description of $Y$ and the
blowup description of $H$.

From these descriptions it immediately follows that locally near a
curve $\rho_{I}$ the maps $f$ and $p$ are given by
\[
\xymatrix@R-1pc@C-1pc{ & H \ar[d] & (x,u,v,s) \ar[d] \\
  Y \ar[r] & X & (x,uv) \\
  (x,y) \ar[r] & (x,2xy^{2} - y^{4})
}
\]
and so $Z$ is locally isomorphic to the hypersurface in
$\mathbb{C}^{5}$ with coordinates $(x,y,u,v,s)$ given by the equation
$uv = xy^{2} - y^{4}$.

In this local model the surface $\Xi_{I}$ corresponds to the plane $y
= u = v = 0$ while the curve $\rho_{I}$ corresponds to the line $x = y
= u = v = 0$. Thus the singularity type indeed jumps from stratum to
stratum: we have a natural projection $(x,y,u,v,s) \mapsto (x,s)$ onto
$\Xi_{I}$ and we get fibers with $A_{1}$ surface singularities for $x
\neq 0$ and fibers with $A_{3}$ singularities for $x = 0$. This
completes the proof of {\bfseries (a)} and {\bfseries (b)}.

Finally for  the proof of {\bfseries (c)} note that in this local
model the divisors $Lp_{I}$ and $LLC_{I}$ are given by $u = 0$ and $v
= 0$ in the coordinates $(x,u,v,s)$ on $H$. So the Weil divisors $LpE_{I}$
and $LLE_{I}$ will be given by the equations $u = y  = 0$ and $v = y =
0$ respectively and so
are not Cartier and represented by transversal rulings. \ \hfill $\Box$

\

\begin{rem} \label{rem:Weil} Note that $\bQ_{T,t}$, $\bQ'_{T,t}$,
  $\bQ''_{T,t}(-\varepsilon)$ are line bundles on $Z$ being
  combinations of pullbacks of Cartier divisors in $Y$ or $H$. However
  the previous lemma shows that the sheaf
  $\bQ_{\lift(T,t)}(-\varepsilon)$ is in general only defined by Weil
  divisors and so is a torsion free divisorial sheaf which in not
  locally free in general. 
\end{rem}

\

\noindent
After these preliminaries we are now ready to compute the
abelianization of the $q$-relative parabolic
Dolbeault complex of $(U_{\bullet},\psi)$:
\[
\Dol(q,U_{\lift(T,t)}) = 
\xymatrix@1{
\Big[ W_{0}(\hor,U_{\lift(T,t)}) 
\ar[r]^-{\verpsi} & 
W_{-2}(\hor,U_{\lift(T,t)})
\otimes \omega_{q}(Lp)\Big].
}
\]
We have already described $U_{\lift(T,t)}$ as a push-forward of a
spectral sheaf on $Z$. To undersand $W_{0}(\hor,U_{\lift(T,t)})$ and
$W_{-2}(\hor,U_{\lift(T,t)})$ we will need to undersand the maps
$\psi$ and $\verpsi$ in the language of spectral data.

By definition $\psi = p^{*}\varphi$ and $\verpsi$ is the composition
of $\psi$ with the projection
\begin{equation} \label{eq:project.on.qvert}
  \Omega^{1}_{H}(\log \Par_{H}) \to
  \Omega^{1}_{H/X\times C}(\log \Par_{H}) = \omega_{q}(Lp).
\end{equation}
But we
already have a description of $\varphi : F_{\bullet} \to
F_{\bullet}\otimes \Omega^{1}_{X}(\log \Par_{X})$ in terms of spectral
data. By construction $\varphi = f_{*}(\alpha\otimes -)$ where $\alpha
\in H^{0}(Y,f^{*}\Omega^{1}_{X}(\log \Par_{X}))$ was the tautological
one form resolving the rational map $Y \to \Higgs_{\op{nilp}}
\dashrightarrow T^{\vee}_{X}$. Therefore
\[
\psi = p^{*}\varphi = p^{*}f_{*}(\alpha\otimes -) =
f_{Z*}\left((p^{*}_{Z}\alpha)\otimes-\right),
\]
where $p^{*}_{Z}\alpha$ is viewed as a section in
\[
p^{*}_{Z}f^{*}\Omega^{1}_{X}(\log \Par_{X})  =
f_{Z}^{*}p^{*}\Omega^{1}_{X}(\log \Par_{X}) \subset
f_{Z}^{*}\Omega^{1}_{H}(\log \Par_{H}).
\]
Pulling back the epimorphism \eqref{eq:project.on.qvert}
to $Z$ we get an epimorphism of vector bundles
\begin{equation} \label{eq:epi.for.beta}
f_{Z}^{*} \Omega^{1}_{H}(\log \Par_{H})  \to f_{Z}^{*}\omega_{q}(Lp). 
\end{equation}
Define
\[
\mybeta \in H^{0}(Z, f_{Z}^{*}\omega_{q}(Lp))
\]
to be the image of $p_{Z}^{*}\alpha \in H^{0}(Z,
f_{Z}^{*}\Omega^{1}_{H}(\log \Par_{H}))$ under the map
\eqref{eq:epi.for.beta}. Tensorization with $\mybeta$ gives a natural
map of divisorial sheaves
\[
\mybeta\otimes - : \bQ_{\lift(T,t)}(-\varepsilon) \to
\bQ_{\lift(T,t)}(-\varepsilon)\otimes f_{Z}^{*}\omega_{q}(Lp),
\]
and by the see-saw principle we get isomorphisms
\[
\begin{aligned}
\left( U_{\lift(T,t)}(-\varepsilon),\psi\right) & = \left(f_{Z*}\bQ_{\lift(T,t)}(-\varepsilon),
  f_{Z*}((p_{Z}^{*}\alpha)\otimes -)\right) \\
\left( U_{\lift(T,t)}(-\varepsilon),\verpsi\right) & = \left(f_{Z*}\bQ_{\lift(T,t)}(-\varepsilon),
  f_{Z*}(\mybeta\otimes -)\right)
\end{aligned}
\]
which provide the desired spectral description of the absolute and
$q$-relative versions of $(U_{\bullet},\psi)$.

To understand the steps $W_{\bullet}(\hor,U_{\lift(T,t)})$ we have to
analyze the action of the residue of $\psi$, or equivalently of
$\verpsi$, on the associated graded $\op{gr}^{\hor}_{0}U_{\lift(T,t)}$
of $U_{\lift(T,t)}$ with respect to the parabolic filtration along the
$q$-horizontal divisor $Lp \subset H$.  By definition
$\op{gr}^{\hor}_{0}U_{\lift(T,t)}$ is a rank four vector bundle on
$Lp$ and 
 the action of $\op{gr}-\res_{Lp} \verpsi$ is induced from the
 action of $\op{res}_{Lp}\verpsi$ on ${U_{\lift(T,t)}}_{|Lp}$.
 So to  compute
$W_{\bullet}(\hor,U_{\lift(T,t)})$ we have to rewrite in
terms of spectral data the basic commutative diagram:
\begin{equation} \label{eq:res.diagram}
\xymatrix@C+1pc{ U_{\lift(T,t)} \ar[r]^-{\verpsi}
  \ar[d]_-{\text{restriction}} & U_{\lift(T,t)} \otimes \omega_{q}(Lp)
  \ar[d]^-{\op{res}_{Lp}} 
  \\
    {U_{\lift(T,t)}}_{|Lp} \ar[r]^-{\op{res}_{Lp}\verpsi} \ar[d] &
      {U_{\lift(T,t)}}_{|Lp} \ar[d] \\
  \op{gr}^{\hor}_{0}U_{\lift(T,t)} \ar[r]_-{\op{gr-res}_{Lp} \verpsi} &
  \op{gr}^{\hor}_{0}U_{\lift(T,t)}  
}
\end{equation}
For this it will be convenient to pull back the natural residue map 
$\op{res}_{Lp} : \omega_{q}(Lp) \to \mathcal{O}_{Lp}$ on $H$ to a map
of sheaves on the spectral cover $f_{Z} : Z \to H$. By construction
$Lp$ pulls back to the non-reduced  and reducible divisor $2LpE + LpG$
and so we the pullback of $\op{res}_{Lp}$ is naturally a surjective
sheaf homomorphism 
\[
f_{Z}^{*}\op{res}_{Lp} : f_{Z}^{*}\omega_{q}(Lp) \to \mathcal{O}_{2LpE
  + LpG}. 
\]
Therefore we have
\[
\op{res}_{Lp} \verpsi =  f_{Z*}\left((f_{Z}^{*}\op{res}_{Lp}
  )(\mybeta)\otimes -\right)
\]
Furthermore we have 
\[
  \begin{aligned}
\op{gr}^{\hor}_{0} U_{\lift(T,t)} & =
U_{\lift(T,t)}/U_{\lift(T,t)}(-\varepsilon) \\
& =
f_{Z*}\bQ_{\lift(T,t)}/f_{Z*}\bQ_{\lift(T,t)}(-\varepsilon) \\
& =
f_{Z*}(\bQ_{\lift(T,t)}/\bQ_{\lift(T,t)}(-\varepsilon)),
\end{aligned}
\]
where as usual we write $\bQ_{\lift(T,t)}$ for $\bQ_{\lift(T,t)}(0)$.

In particular the diagram \eqref{eq:res.diagram} can be rewritten as the
$f_{Z}$ push forward of the following diagram of sheaves on  $Z$:
\begin{equation} \label{eq:spectralres.diagram}
\xymatrix@C+5pc{\bQ_{\lift(T,t)} \ar[r]^-{\mybeta\otimes -}
  \ar[d]_-{\text{restriction}} & \bQ_{\lift(T,t)} \otimes f_{Z}^{*}\omega_{q}(Lp)
  \ar[d]^-{f_{Z}^{*}\op{res}_{Lp}} 
  \\
    {\bQ_{\lift(T,t)}}_{|2LpE+LpG}
    \ar[r]^-{f_{Z}^{*}\op{res}_{Lp}(\mybeta)\otimes -} \ar[d] &
      {\bQ_{\lift(T,t)}}_{|2LpE+LpG} \ar[d] \\
      \bQ_{\lift(T,t)}/\bQ_{\lift(T,t)}(-\varepsilon)
      \ar[r]_-{\op{gr-}(f_{Z}^{*}\op{res}_{Lp}(\mybeta)\otimes -)} &
  \bQ_{\lift(T,t)}/\bQ_{\lift(T,t)}(-\varepsilon).
}
\end{equation}

Therefore $W_{\bullet}(\hor,U_{\lift(T,t)})$ depends on order of
vanishing of $\mybeta$ on $LpE$ and $LpG$ and on the length of the
scheme theoretic support of
$\bQ_{\lift(T,t)}/\bQ_{\lift(T,t)}(-\varepsilon)$ which we proceed to
analyze.

\subsubsection{Abelianization of the Hecke
  correspondence} \label{ssec:abH}

To understand the divisor of $\mybeta$ on $Z$ we will need to understand
the 
abelianized parabolic Hecke correspondence
\[
\xymatrix@R-1pc{
& \abH \ar[dl]_-{\abp} \ar[dr]^-{\abq} & \\
Y & & Y\times \sC
}
\]
which similarly to $H$ is a compactification and a resolution of the
usual abelian Hecke correspondence acting along the fibers of the
Hitchin map. In other words $\abH$ is a birational lift of the graph
of the  Abel-Jacobi map \index{terms}{Hitchin!map}
\[
\AJ : \sC\times \Jac \to \sC, \qquad \AJ(x,\mycal{L}) = \mycal{L}\otimes
\mathcal{O}_{\sC}(x - \tilde{p}_{6})
\]
As a first approximation of $\abH$ we can take the 
fiber
product
\begin{equation} \label{eq:defabH}
  \xymatrix@R-0.5pc@C-0.5pc{
 Y\times_{\Jac} (Y\times \sC) \ar[dd]_-{\abp}
  \ar[r]^-{\abq} & Y\times \sC \ar[d]^-{\bmu\times
      \mathbf{1}}\\
 & \Jac\times\sC \ar[d]^-{\AJ} \\
    Y \ar[r]_-{\bmu} & \Jac 
    }
  \end{equation}
  The geometry of this fiber product is easy to understand.

  \begin{lem} \label{lem:fiber.product}
 \begin{itemize}
\item[(a)] The fiber product $Y\times_{\Jac} (Y\times \sC)$ is a
  compact threefold with $96$ conifold singularities.
\item[(b)] The image of $Y\times_{\Jac} (Y\times \sC)$ in $X \times
  X\times C$ is contained in the image of $H$:
 \begin{equation} \label{eq:include.in.H}
(f\times f\times \pi)\circ (\tilde{p}\times \tilde{q})\left(Y\times_{\Jac}
    (Y\times \sC)\right) \subset (p\times q)(H).
\end{equation}
In particular we have a natural rational map
$Y\times_{\Jac} (Y\times \sC) \dashrightarrow H$.
\end{itemize}
\end{lem}

\

\begin{figure}[!ht]
\begin{center}
\psfrag{GammaI}[c][c][1][0]{{$\Gamma_{I}$}}
\psfrag{0}[c][c][1][0]{{$0$}}
\psfrag{EI+1}[c][c][1][0]{{$E_{I+1}\times \tilde{p}_{1}$}}
\psfrag{EI+2}[c][c][1][0]{{$E_{I+2}\times \tilde{p}_{2}$}}
\psfrag{EI+3}[c][c][1][0]{{$E_{I+3}\times \tilde{p}_{3}$}}
\psfrag{EI+4}[c][c][1][0]{{$E_{I+4}\times \tilde{p}_{4}$}}
\psfrag{EI+5}[c][c][1][0]{{$E_{I+5}\times \tilde{p}_{5}$}}
\psfrag{EI}[c][c][1][0]{{$E_{I}\times \tilde{p}_{6}$}}
\epsfig{file=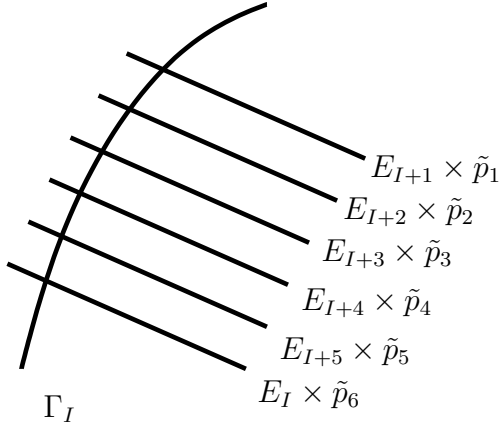,width=2.5in} 
\end{center}
\caption{The curve $\hatAJ^{-1}(p_{I})$}\label{fig:comb}  
\end{figure}

\noindent
{\bfseries Proof.} (a) \ Indeed, let
$\hatAJ : Y\times \sC \to \Jac$ denote the composition
$\hatAJ = \AJ\circ (\bmu\times \mathsf 1)$. Then
$Y\times_{\Jac} (Y\times \sC)$  is a blowup of $Y\times \sC$ in the 16
disjoint one dimensional subschemes $ \hatAJ^{-1}(p_{I})$:
\[
  Y\times_{\Jac} (Y\times \sC) = \op{Bl}_{\underset{I \in \oddL}{\sqcup}
    \hatAJ^{-1}(p_{I})}(Y\times \sC).
\]
For a point of order two $p_{I} \in \Jac$ the subscheme
$\hatAJ^{-1}(p_{I}) \subset Y\times \sC$ is a reduced and
reducible curve which has the shape of a comb (see
Figure~\ref{fig:comb}).
This curve has seven components:
\begin{equation}
\hatAJ^{-1}(p_{I}) = \Gamma_{I} \bigcup \left( \bigcup_{i =
    1}^{6} E_{I+i}\times \tilde{p}_{i} \right),
\end{equation}
where the curves $E_{I+i}\times \tilde{p}_{i}$ are the disjoint teeth
of the comb and the spine of the comb is the graph
$\Gamma_{I} = \op{Graph}(\sC \to G_{I})$ of the canonical
identification of $\sC$ with $G_{I}$. 

\
 
Thus $Y\times_{\Jac} (Y\times \sC)$ is locally isomorphic to the
blow-up of $\mathbb{C}^{3}$ at the union of two coordinate
axes. Explicitly the local picture is
\[
\xymatrix@C-1pc@R-1pc{
  & Y\times \sC \ar[d] & (x_{1},u_{1},p) \ar[d]  \\
  & \Jac\times \sC \ar[d] & (x_{1},x_{1}u_{1},p) \ar[d] \\
  Y \ar[r] & \Jac & (x_{1}+p,x_{1}u_{1}) \\
  (x_{2},u_{2}) \ar[r] & (x_{2},x_{2}u_{2})
} 
\]
and so locally the fiber product is defined by the
equations
\[
  \left| \
    \begin{aligned}
      x_{2} & = x_{1}+p \\
      x_{2}u_{2} & = x_{1}u_{1}
     \end{aligned}
  \right.
\]
or in $\mathbb{C}^{4}$ with coordinates $(x_{1},x_{2},u_{1},u_{2})$ by
the single equation $x_{2}u_{2} = x_{1}u_{1}$.

\

\noindent
For part (b) note that the standard fact \cite{dp-Langlands,dp-jdg}
that the Abel-Jacobi map on Hitchin fibers lifts the Hecke \index{terms}{Hitchin!fiber}
correspondence on the moduli of bundles gives the desired inclusion on
Zariski open sets. Since both the fiber product and $H$ are
irreducible this implies the global inclusion \eqref{eq:include.in.H}.
\ \hfill $\Box$

\

\smallskip

\noindent
Because of these conifold singularities the fiber product is not
locally factorial. To remedy this and avoid working with Weil divisors
we can pass to a small resolution of $Y\times_{\Jac}(Y\times \sC)$.
The fiber product has two natural small resolutions that are related by
flops, namely:
\[
 \op{Bl}_{\underset{I,i}{\sqcup} \widehat{E_{I}\times
   \tilde{p}_{i}}}\op{Bl}_{\Gamma_{I}}(Y\times \sC) \qquad
 \text{and} \qquad 
  \op{Bl}_{\widehat{\Gamma_{I}}}\op{Bl}_{\underset{I,i}{\sqcup} E_{I}\times
   \tilde{p}_{i}}(Y\times \sC).
\]
For concreteness we will choose the first small resolution as  our
model for the abelianized parabolic Hecke correspondence. We set
\[
\abH :=  \op{Bl}_{\underset{I,i}{\sqcup} \widehat{E_{I}\times
   \tilde{p}_{i}}}\op{Bl}_{\Gamma_{I}}(Y\times \sC),
\]
and by a slight abuse of notation will write again $\abp$ and $\abq$
for the natural maps:
\[
\xymatrix@R-1pc{
& \abH \ar[dl]_-{\abp} \ar[dr]^-{\abq} & \\
Y & & Y\times \sC
}
\]
We will denote by  $g$  the natural rational map 
\[
g : \abH \to H
\]
which is the composition of the small blowdown map $\abH \to
Y\times_{\Jac}(Y\times \sC)$ and the map described in
Lemma~\ref{lem:fiber.product}(b). Finally note that by construction
the map $\abp\times g : \abH \to Y\times H$ factors through the fiber
product $Z = Y\times_{X} H$ and so we get another rational map $g_{Z}
: \abH \to Z$ which is a closed immersion away from a closed subvariety
in $\abH$ which maps to a curve in $Y\times \sC$.

  \

  \begin{rem} \label{rem:away.from.curves}
    Since we chose $(E_{\bullet},\theta)$, $(F_{\bullet},\varphi)$ and
    $(F_{\bullet}',\varphi')$ to satisfy Mochizuki's non-abelian Hodge
    theorem conditions, by construction both the left hand side and
    the right hand side of the Hecke conditions are complexes of
    parabolic Higgs sheaves that are specializations at zero of
    twistor families of parabolic $z$-connections. Such parabolic
    objects are uniquely determined by their restrictions to any
    sufficiently ample hyperplane section curve inside $X\times
    C$. Since such a curve can be chosen to be disjoint from any
    finite collection of curves in $X\times C$, it suffices to to
    compare $\lhs_{T,t}$ and $\rhs_{T,t}$ on the complement of some
    finite collection of curves in $X\times C$.  This means that in
    the geometric considerations above we can freely excise curves in
    $X\times C$ (and their preimages in $H$, $Y\times \sC$, $\abH$,
    etc.)  over which something pathological happens.

In the considerations of the Hecke condition that follow we will tacitly assume
that the following curves in $X\times C$ (and their preimages in $H$, $Y\times
\sC$, $\abH$, etc.) were excised:
\begin{itemize}
\item all curves over which the map $g : \abH \to H$ is not proper and
  over which the map $g_{Z} : \abH \to Z$ is not a closed immersion;
\item all curves over which the maps $Lp_{I}
  \to X\times C$ is not finite;
\item all curves over which the divisors $Lp_{I}$ intersect with each other.
\end{itemize}

We will not indicate that these curves are excised in the notation but
will treat all spaces and maps as if these sources of pathological
behavior are not present. In particular we will treat $g$ and $g_{Z}$
as proper morphisms with the unspoken understanding that sufficient chunks of
the spaces under consideration have been excised for this to be indeed
the case.
  \end{rem}

Together the maps $g$ and $g_{Z}$ fit in a
commutative diagram
\[
  \xymatrix{
& \abH \ar[dl]_-{g_{Z}} \ar[dr]^-{g} & \\
Z \ar[d]_-{p_{Z}} \ar[rr]^-{f_{Z}} & & H \ar[d]^-{p} \\
Y \ar[rr]_-{f} & & X }
\]
and we can view $\abH$ as a Cartier divisor in $H$ and a Weil divisor
in $Z$ via the maps $g$ and $g_{Z}$. This is significant since the
hypersurface $\abH \subset Z$ is naturally a component of the zero
locus of $\mybeta$. More precisely we have the following

\begin{lem} \label{lem:divisor.of.beta}
  The divisor of $\mybeta \in H^{0}(Z,f_{Z}^{*}\omega_{q}(Lp))$ in $Z$ is 
\[
(\mybeta) = \abH + LpE + LpG.  
\]
\end{lem}

\

\noindent
{\bfseries Proof.} \ Since by definition $Z$ is swept by preimages of
Hecke curves in $H$ it suffices to understand the  restriction of $\mybeta$ to
the preimage of a Hecke curve. Fix a point $(x,c) \in X\times C$ and
let $H_{x,c} = q^{-1}(x,c)$ be the Hecke curve corresponding to
$(x,c)$. For a generic $(x,c)$ the curve $H_{x,c}$ is a smooth
rational curve. The image $X_{x,c}$ of $H_{x,c}$ under $p : H \to X$
is a nodal $\mathbb{P}^{1}$. The preimage $Y_{x,c} = f^{-1}(X_{x,c})
\subset Y$ is a curve of geometric genus $13$ with four nodes while
the preimage $Z_{x,c} = f_{Z}^{-1}(H_{x,c})$ is the smooth genus $13$
curve  normalizing $Y_{x,c}$. Thus we have a commutative diagram of curves 
\[
\xymatrix{
  Z_{x,c} \ar[r]^-{f_{Z}} \ar[d]_-{p_{Z}} & H_{x,c} \ar[d]^-{p} \\
  Y_{x,c} \ar[r]_-{f} & X_{x,c}
  }
\]
where the vertical maps are normalization maps.

From the blowup description of $H$ it is immediate that $LLC$ does not
intersect $H_{x,c}$. Furthermore by definition $X_{x,c}$ is the
intersection of $X$ with the projective tangent space
$\overline{T_{x}Q_{c}} \cong \mathbb{P}^{3} \subset \mathbb{P}^{4}$ of
the quadric $Q_{c}$ at the point $x \in X \subset Q_{c} \subset
\mathbb{P}^{4}$. In particular $X_{x,c}$ intersects each line $L_{I}$
transversally at a single point. Since $p^{*}(L_{I}) = Lp_{I} +
LLC_{I}$ it follows that $Lp_{x,c} := Lp\cap H_{x,c} = p^{*}(X_{x,c}\cap L)$,
i.e. the intersection $Lp_{x,c}$ of $Lp$ with the curve $H_{x,c}$ consists of
$16$ distinct points.

In particular the line bundle $\omega_{q}(Lp)_{|H_{x,c}} =
\omega_{H_{x,c}}(Lp_{x,c})$ has degree $14$ and so
$f_{Z}^{*}\omega_{q}(Lp))_{|Z_{x,c}}$ has degree $56$.

From the modular description of $H$ we have a different but equally
usefuld description of $H_{x,z}$. If we view $X$ as the moduli space
$N_{1}(a,b)$ we can interpret the point $x \in X$ as a
parabolic rank two bundle $(V,\bF)$ on $(C,\Par_{C})$ with
$\deg V = 1$. By definition the  Hecke curve $H_{x,c}$ is 
the projective space  $\mathbb{P}(V_{c})$ parametrizing all lines in the
fiber $V_{c}$ of $V$ at $c$. From this point of view  the map
\[
p : H_{x,c} \to X_{x,c} \subset X
\]
assigns to each line $\boldsymbol{L} \subset V_{c}$ the (down) Hecke
transform $({}^{\boldsymbol{L}}V,{}^{\boldsymbol{L}}\bF)$ centered at
$\boldsymbol{L}$, i.e. assuming $c \in C - \Par_{C}$ we have 
\[
{}^{\boldsymbol{L}}V = \ker\left[ V \to V_{c}/\boldsymbol{L}\right], \
\text{and}  \ {}^{\boldsymbol{L}}\bF = \bF.
\]
The section
\[
\mybeta_{|Z_{x,c}} \in H^{0}\left(Z_{x,c},f_{Z}^{*}\omega_{H_{x,c}}(Lp_{x,c})\right)
\]
is the image of
$p_{Z}^{*}\alpha \in H^{0}(Z,f_{Z}^{*}\Omega^{1}_{H}(\log \Par H))$ to
logarithmic one forms on $H_{x,c}$. But $L$ is part of the branch
divisor of $f$, the preimage of $L$ in $Y$ is $f^{*}L = 2E+G$, and as
we saw in section~\ref{sssec:eigensheaf} the residues of $\alpha$ on
$E$ and $G$ have one dimensional kernel. Therefore
$\mybeta_{|Z_{x,c}}$ must vanish along the divisor
$LpE_{x,c} + LpG_{x,c}$ which has degree $48$.

Furthermore if $z \in Z_{x,c}$ is a smooth point which is not in the
ramification locus of $f_{Z} : Z_{x,c} \to H_{x,c}$, then by
definition $\mybeta(z) \in \Omega^{1}_{Z_{x,c},z}$ is the restriction
of $\alpha(p_{Z}(z))$ on the tangent space of $Y_{x,z}$ at
$p_{Z}(z)$. The assumption that $f_{Z}$ is unramified at $z$ is
equivalent to $f$ being unramified at $p_{Z}(z)$ and so we can use the
differential of $f$ to identify $T_{Y_{x,c},p_{Z}(z)}$ with
$\op{Hom}(\boldsymbol{L},\boldsymbol{Q})$ where
$\boldsymbol{L} \subset V_{c}$ is
the line corresponding to $f_{Z}(z) \in \mathbb{P}(V_{c}) = H_{x,c}$
and $\boldsymbol{Q} = V_{c}/\boldsymbol{L}$.

If we use the codifferential of $f$ to identify
$\Omega^{1}_{Y,p_{Z}(z)}$ with
$\Omega^{1}_{X,p(\boldsymbol{L})}$, and the codifferential of $f_{Z}$
to identify  $\Omega^{1}_{Z_{x,c},z}$ with
  $\Omega^{1}_{H_{x,c},\boldsymbol{L}}$ we see that the element 
\begin{equation} \label{eq:Z.with.H}
  \mybeta(z) \in \Omega^{1}_{Z_{x,c},z} \cong
  \Omega^{1}_{H_{x,c},\boldsymbol{L}}
\end{equation}
is simply the image of 
the element
\begin{equation} \label{eq:Y.with.X}
\alpha(p_{Z}(z)) \in \Omega^{1}_{Y,p_{Z}(z)} \cong
\Omega^{1}_{X,p(\boldsymbol{L})}
\end{equation}
under the codifferential
\[
dp^{\vee}_{\boldsymbol{L}} :
\Omega^{1}_{X,p(\boldsymbol{L})} \to \Omega^{1}_{H_{x,c},\boldsymbol{L}}
\]
of the map $p : H_{x,c}
\to X$ at the point $\boldsymbol{L} \in H_{x,c}$.

Under our genericity assumption on $z$ we
can describe the point $p_{Z}(z)$ in modular terms as a parabolic rank
two Higgs bundle on $(C,\Par_{C})$, namely
\[
p_{Z}(z) = \left( {}^{\boldsymbol{L}}V,{}^{\boldsymbol{L}}\bF; \ {}^{z}\theta \right) 
\]
where $({}^{\boldsymbol{L}}V,{}^{\boldsymbol{L}}\bF) = p(\boldsymbol{L})$ is
the parabolic bundle defined above and 
\[
  {}^{z}\theta :
{}^{\boldsymbol{L}}V \to {}^{\boldsymbol{L}}V\otimes \Omega^{1}_{C}(\log
\Par_{C})
\]
is one of the four Higgs fields with nilpotent residues on
${}^{\boldsymbol{L}}V$ whose spectral cover is precisely
$\sC$. Furthermore in modular terms we have the identification
\begin{equation} \label{eq:cotangent.as.higgs}
  \Omega^{1}_{X,p(\boldsymbol{L})} =
H^{0}(C,\op{End}_{0}({}^{\boldsymbol{L}}V)\otimes \Omega^{1}_{C}(\log
\Par_{C}))_{\op{nilp}}
\end{equation}
where $H^{0}(C,\op{End}_{0}({}^{\boldsymbol{L}}V)\otimes
\Omega^{1}_{C}(\log \Par_{C}))_{\op{nilp}}$ is the vector space of all
traceless Higgs fields ${}^{\boldsymbol{L}}V \to {}^{\boldsymbol{L}}V\otimes
\Omega^{1}_{C}(\log \Par_{C})$ whose residues at the points $p_{i} \in
\Par_{C}$ are nilpotent and preserve the line $F_{i} \subset
    {}^{\boldsymbol{L}}V_{p_{i}} = V_{p_{i}}$.

Now note that by definition of the modular spectral cover $Y$, the
section $\alpha$ is simply the natural map from $Y$ to the total space
of $\Omega^{1}_{X}(\log L)$. Thus in the identifications
\eqref{eq:Y.with.X} and

\eqref{eq:cotangent.as.higgs} the element $\alpha(p_{Z}(z))$
corresponds to the element
\[
{}^{z}\theta \in H^{0}(C,\op{End}_{0}({}^{\boldsymbol{L}}V)\otimes
\Omega^{1}_{C}(\log \Par_{C}))_{\op{nilp}},
\]
and so we get that under the identification \eqref{eq:Z.with.H}
$\mybeta(z)$ corresponds to the element
\[
dp^{\vee}_{\boldsymbol{L}}({}^{z}\theta) \in \Omega^{1}_{H_{x,c},\boldsymbol{L}}. 
\]
The codifferential $dp^{\vee}_{\boldsymbol{L}}$ can be computed
explicitly. Indeed, since $H_{x,c} =
\mathbb{P}(V_{c})$  we have a canonical identification
\[
\Omega^{1}_{H_{x,c},\boldsymbol{L}} = \op{Hom}(\boldsymbol{Q},\boldsymbol{L}).
\]
By definition we have a short exact sequence of vector spaces
\begin{equation} \label{eq:center.of.Hecke}
\xymatrix@1{
  0 \ar[r] &  \boldsymbol{L} \ar[r] &  V_{c} \ar[r] & \boldsymbol{Q}
  \ar[r] &  0.}
\end{equation}
Since ${}^{\boldsymbol{L}}V$ is the Hecke transform of $V$ centered at
$\boldsymbol{L}$ we have that the fiber ${}^{\boldsymbol{L}}V_{c}$ of
${}^{\boldsymbol{L}}V$ fits canonically in a short exact sequence
\begin{equation} \label{eq:flipped.center}
\xymatrix@1{
0 \ar[r] & \boldsymbol{Q}\otimes \mathcal{O}_{C,c}(-c) \ar[r] & 
{}^{\boldsymbol{L}}V_{c} \ar[r] & \boldsymbol{L} \ar[r] & 0.}
\end{equation}
If $\eta \in H^{0}(C,\op{End}_{0}({}^{\boldsymbol{L}}V)\otimes
\Omega^{1}_{C}(\log \Par_{C}))_{\op{nilp}}$ is any Higgs field, then evaluating
$\eta$ at $c$ gives a linear map
\[
\eta(c) : {}^{\boldsymbol{L}}V_{c} \to {}^{\boldsymbol{L}}V_{c}\otimes \Omega^{1}_{C,c}
\]
Pre composing $\eta(c)$ with the inclusion map in
\eqref{eq:flipped.center} and post composing with the projection map
in \eqref{eq:flipped.center} yields a linear map
\[
\boldsymbol{Q}\otimes \mathcal{O}_{C,c}(-c) \to \boldsymbol{L}\otimes
\Omega^{1}_{C,c}
\]
or equivalently a linear map
\[
\overline{\eta(c)} : \boldsymbol{Q} \to \boldsymbol{L}\otimes
\Omega^{1}_{C}(c)_{c}.
\]
Now standard deformation theory identifies  the codifferential
$dp^{\vee}_{\boldsymbol{L}}$ with the map:
\[
\xymatrix@R-3.2pc@C+1pc@M+1pc{ dp^{\vee}_{\boldsymbol{L}} : \hspace{-5pc} &
  H^{0}(C,\op{End}_{0}({}^{\boldsymbol{L}}V)\otimes
  \Omega^{1}_{C}(\log \Par_{C}))_{\op{nilp}} \ar[r] &
\op{Hom}(\boldsymbol{Q},\boldsymbol{L}) \\
& \eta \ar@{|->}[r] & \mathsf{res}_{c}\circ \overline{\eta(c)}
}
\]
where $\mathsf{res}_{c} : \Omega^{1}_{C}(c)_{c} \to \mathcal{O}_{C,c}$
is the residue isomorphism at $c$.

In other words, up to the residue isomorphism,
$dp^{\vee}_{\boldsymbol{L}}$ assigns to $\eta$ the off-diagonal matrix
entry $\overline{\eta(c)}$ of the map $\eta(c)$.  Finally, note that
by definition $\abH_{x,c} \subset Z_{x,c}$ is precisely the
subset of all points $z \in Z_{x,c}$ for which ${}^{z}\theta(c)$ is
diagonal, i.e. is a linear map preserving $\boldsymbol{Q}\otimes
\mathcal{O}_{C,c}(-c)$ and inducing a linear endomorphism of
$\boldsymbol{L}$.

This proves that $\mybeta_{|Z_{x,c}}$ vanishes at the $8$ points of
the divisor $\abH_{x,c} \subset Z_{x,c}$. This accounts for
$56$ zeroes of the section $\mybeta_{|Z_{x,c}}$. But this is a nonzero
section in a line bundle of degree $56$ and so the divisor of
$\mybeta_{|Z_{x,c}}$ is precisely $\abH_{x,c} + LpE_{x,c} +
LpG_{x,c}$.

Finally note that locus of zeroes of $\mybeta_{|Z_{x,c}}$ for
non-generic $(x,c)$ has codimension two in $Z$ and so the divisor of
$\mybeta$ is
\[
(\mybeta) = \abH + LpE + LpG.
\]
This concludes the proof of the lemma.
\ \hfill $\Box$

\bigskip

\

\noindent
Using the formula for the divisor of $\mybeta$ and the spectral
description \eqref{eq:spectralres.diagram} of the horizontal weight
filtration of $(U_{\lift(T,t)},\verpsi)$ we can now derive a spectral
description of the left hand side entry $\lhs_{T,t}$of the Hecke
condition. To focus on cases of interest, note that by definition, the
right hand side entry $\rhs_{T,t}$ of the Hecke condition is a
parabolic rank $8$ Higgs bundle on $X\times C$. Thus we only need to
consider data for which the left hand side is also of rank $8$. Using
the definition of $\bQ_{\lift(T,t)}$ (see \eqref{eq:frakM}), $LpE$
(see \eqref{eq:LpEetc}), and the calculation of the \cite{dps} formula
for the pushforward of the relative parabolic Dolbeault complex of
$(U_{\lift(T,t)},\verpsi)$ (see \eqref{eq:lhs}) together with the
calculation of $U_{\lift{T,t}}(-\varepsilon)$ (see \eqref{eq:UupviaM}),
we can now summarize the relevant abeliazation result in the
following proposition.

\

\pagebreak

\begin{prop} \label{pro:abelianize.lhs} Let
  $(F_{\bullet},\varphi)$ and the Hecke kernel be given by vectors $e,
  d, lp, lq, llc \in \mathbb{R}^{16}$ and $r_{1}, r_{2} \in
  \mathbb{R}^{5}$ as above, and let
  \[
\lhs_{T,t} = \left( \, q_{*}\Dol(q,U_{\lift(T,t)}),
q_{*}\mathfrak{d}(\psi) \, \right)
  \]
  be the associated left hand side of the Hecke condition.  Then
  \begin{itemize}
\item For the complex $q_{*}\Dol(q,U_{\lift(T,t)})$ to be a vector
  bundle of rank $8$ it is necessary that $lp + d \in
  \mathbb{Z}^{16}$.
\item If $lp + d \in \mathbb{Z}^{16}$, then the two term complex
  $\Dol(q,U_{\lift(T,t)})$ is quasi-isomorphic to a single pure
  codimension one sheaf on $H$ which is a direct image of a sheaf
  on $\abH$ under the natural map $g : \abH \to H$.
  That is
  \begin{equation} \label{eq:push.formula}
\begin{aligned}
  \Dol(q,U_{\lift(T,t)}) & \cong
  f_{Z*}\left(\bQ_{\lift(T,t)}\left(-\sum_{\substack{\left\{I \, | \,
      lp_{I} + e_{I} \in \mathbb{Z}\right\}}}LpE_{I}
  + \abH\right)_{|\abH}\right) \\[+0.5pc]
& \cong
  g_{*}\left(g_{Z}^{*}\left(\bQ_{\lift(T,t)}\left(-
  \sum_{\substack{\left\{I \, | \,
      lp_{I} + e_{I} \in \mathbb{Z}\right\}}}LpE_{I} +
  \abH\right)\right)\right).
\end{aligned}
 \end{equation}
\end{itemize}
 In other words $\Dol(q,U_{\lift(T,t)})$ is quasi-isomorphic to a sheaf
 supported on $g(\abH) \subset H$ which is
 given by the spectral datum $(\abH,
\bQ_{\lift(T,t)}(-\sum_{\substack{\left\{I \, | \,
      lp_{I} + e_{I} \in \mathbb{Z}\right\}}}LpE_{I} + \abH)_{|\abH})$. 
\end{prop}

\

\noindent
    {\bfseries Proof.} \ Write $W_{\bullet}(\hor,\bQ_{\lift(T,t)})
    \subset \bQ_{\lift(T,t)}$ for the divisorial subsheaves defined as
    the pullbacks of the steps of the monodromy weight filtration on
    the sheaf $\bQ_{\lift(T,t)}/\bQ_{\lift(T,t)}(-\varepsilon)$
    associated to the nilpotent endomorphism
    $\op{gr-}(f_{Z}^{*}\op{res}_{Lp}(\mybeta)\otimes -)$.

 By construction (see section~\ref{sssec:ab.lhs}) the pair
 $(Z,W_{\bullet}(\hor,\bQ_{\lift(T,t)}))$ is the spectral datum for
 the sheaves $W_{\bullet}(\hor,U_{\lift(T,t)})$, that is
 $W_{\bullet}(\hor,U_{\lift(T,t)}) =
 f_{Z*}W_{\bullet}(\hor,\bQ_{\lift(T,t)})$. In particular we have
 \[
 \Dol(\hor,U_{\lift(T,t)}) = f_{Z*}\left[ \xymatrix@1{
     W_{0}(\hor,\bQ_{\lift(T,t)}) \ar[r]^-{\mybeta\otimes -} &
     W_{-2}(\hor,\bQ_{\lift(T,t)})\otimes f_{Z}^{*}\omega_{q}(Lp)
   }\right].
 \]
In other words the complex
 \begin{equation} \label{eq:actual.spectral}
\xymatrix@1{
     W_{0}(\hor,\bQ_{\lift(T,t)}) \ar[r]^-{\mybeta\otimes -} &
     W_{-2}(\hor,\bQ_{\lift(T,t)})\otimes f_{Z}^{*}\omega_{q}(Lp)}
 \end{equation}
 of divisorial sheaves on $Z$ is the spectral datum for the complex
$\Dol(\hor,U_{\lift(T,t)})$.
 It is obtained from the complex
 \begin{equation} \label{eq:naive.spectral}
\xymatrix@1{
     \bQ_{\lift(T,t)} \ar[r]^-{\mybeta\otimes -} &
     \bQ_{\lift(T,t)}\otimes f_{Z}^{*}\omega_{q}(Lp)}
 \end{equation}
by restricting the map in \eqref{eq:naive.spectral} to the torsion
free subsheaf $W_{0}(\hor,\bQ_{\lift(T,t)}) \subset \bQ_{\lift(T,t)}$
and observing that the image of the restricted map lands in the
torsion free subsheaf $W_{-2}(\hor,\bQ_{\lift(T,t)})\otimes
f_{Z}^{*}\omega_{q}(Lp) \subset \bQ_{\lift(T,t)}\otimes
f_{Z}^{*}\omega_{q}(Lp)$. Since the map in \eqref{eq:naive.spectral}
is given by tensoring with the non-zero section $\mybeta$ in the line
bundle $f_{Z}^{*}\omega_{q}(Lp)$ it follows that the map in
\eqref{eq:naive.spectral} is injective and therefore the map in
\eqref{eq:actual.spectral} is injective.

This implies that the length two complex \eqref{eq:actual.spectral}
has no cohomology in degree zero and so is quasi isomorphic to its
degree one cohomology sheaf $\mycal{H}^{1}_{T,t}$. From the formula
\eqref{eq:spectral.sheaf} we see that the support of the quotient
$\bQ_{\lift(T,t)}/\bQ_{\lift(T,t)}(-\varepsilon)$ is contained in the
divisor $LpE+LpG$ and the precise shape of this quotient depends on
whether $lp_{I} + d_{I}$ and $lp_{I} + e_{I}$ are
integers. Specifically we see that the sheaf
$\bQ_{\lift(T,t)}/\bQ_{\lift(T,t)}(-\varepsilon)$ has rank
\begin{itemize}
\item $0$ or $1$ on $LpG_{I} =$ the number of integers in the set $\{
  lp_{I} + d_{I} \}$;
\item $0$, $1$, or $2$ on $LpE_{I} =$ the number of integers in the set
$\{ lp_{I} + d_{I}, lp_{I} + e_{I} \}$.  
\end{itemize}
Furthermore from the definition of $\mybeta$ (or equivalently
$\verpsi$) we know that the action of
$\op{gr}-(f_{Z}^{*}\op{res}_{Lp_{I}}(\mybeta)\otimes -)$ on the
restriction of $\bQ_{\lift(T,t)}/\bQ_{\lift(T,t)}(-\varepsilon)$ on
either $LpG_{I}$ of $2LpE_{I}$ is given by a single Jordan block of
maximal size.

Thus we get the following possibilities for the terms and the degree
one cohomology sheaf of the complex \eqref{eq:actual.spectral}:

\

\bigskip

\begin{center}
{\scriptsize
  \begin{tabular}{|l||c|c|c|}
    \hline
    & $W_{0}\bQ_{\lift(T,t)}$ & $W_{-2}\bQ_{\lift(T,t)}$ & $\mycal{H}^{1}_{T,t}$ \\
    \hline\hline
    $\begin{aligned}
      lp_{I} + d_{I} & \in \mathbb{Z} \\[-0.5pc]
      lp_{I} + e_{I} & \in \mathbb{Z}
    \end{aligned}$
    & $\bQ_{\lift(T,t)}(-LpE_{I})$ & $\bQ_{\lift(T,t)}(-2LpE_{I}
    - LpG_{I})$ & $\bQ_{\lift(T,t)}(- LpE_{I} + \underline{LpE}_{I} +
    \underline{LpG}_{I} + \abH)_{|
      \underline{LpE}_{I} + \underline{LpG}_{I} + \abH}$ \\
    \hline
    $\begin{aligned}
      lp_{I} + d_{I} & \in \mathbb{Z} \\[-0.5pc]
      lp_{I} + e_{I} & \notin \mathbb{Z}
    \end{aligned}$
    & $\bQ_{\lift(T,t)}$ & $\bQ_{\lift(T,t)}(-LpE_{I}
    - LpG_{I})$ & $\bQ_{\lift(T,t)}(\underline{LpE}_{I} + \underline{LpG}_{I} +
    \abH)_{|
      \underline{LpE}_{I} + \underline{LpG}_{I}  +
      \abH}$ \\
    \hline
  $\begin{aligned}
      lp_{I} + d_{I} & \notin \mathbb{Z} \\[-0.5pc]
      lp_{I} + e_{I} & \in \mathbb{Z}
    \end{aligned}$
    & $\bQ_{\lift(T,t)}$ & $\bQ_{\lift(T,t)}(-LpE_{I})$ &
    $\bQ_{\lift(T,t)}(\underline{LpE}_{I} + LpG+ \abH)_{|
      \underline{LpE}_{I} + LpG + \abH}$ \\
    \hline
    $\begin{aligned}
      lp_{I} + d_{I} & \notin \mathbb{Z} \\[-0.5pc]
      lp_{I} + e_{I} & \notin \mathbb{Z}
    \end{aligned}$
    & $\bQ_{\lift(T,t)}$ & $\bQ_{\lift(T,t)}$ &
    $\bQ_{\lift(T,t)}(LpE + LpG+ \abH)_{|
      LpE + LpG + \abH}$ \\
    \hline\hline 
  \end{tabular}
  }
\end{center}

\

\bigskip

\noindent
Here we use the notation $\underline{LpE}_{I}$ and $\underline{LpG}_{I}$
for the the divisors $\sum_{J \neq I} LpE_{J}$ and $\sum_{J\neq
  I} LpG_{J}$ respectively.

Now note that if we view the threefolds $LpE_{I}$, $LpG_{I}$, and
$\abH$ with their natural maps to $X\times C$, then by definition we
have that $q\circ f_{Z} : LpE_{I} \to \abH$ is birational, $q\circ f_{Z} :
LpG_{I} \to X\times C$ is generically of degree two, and $q\circ f_{Z}
: \abH \to X\times C$ is genericall of degree eight.  Since $\abH$ is
always a component of the support of $\mycal{H}^{1}_{T,t}$ we see that
the only way the pushforward of $\mycal{H}^{1}_{T,t}$ to $X\times C$
can have rank eight is if we are in one of the first two rows of the
above table and in fact if we are in the first two rows for all values
of $I \in \oddL$. This proves the formula \eqref{eq:push.formula}
\ \hfill $\Box$

\

\subsubsection{Abelianization of the Hecke conditions}
\label{sssec:ab.Hecke.conditions}

\noindent
In summary we have reduced the parabolic Hecke eigensheaf condition
\eqref{eq-par.hecke} to the following statement:

\

\noindent
For any choice of $a, b \in \mathbb{R}^{5}$ with $\sum_{i= 1}^{5}
(a_{i} + b_{i}) = 3$, there exists a natural $\bzeta =
(lp,lq,llc,r_{1},r_{2})$ and natural vectors $e, d, e', d' \in
\mathbb{R}^{16}$ so that:
\begin{itemize}
\item The vectors $lp$ and $d$ satisfy $lp + d \in \mathbb{Z}^{16}$.
\item For any smooth spectral curve $\sC \to C$
branched over $\Par_{C}$ and any line bundle $\mathfrak{a} \in \Pic^{0}(\sC)$,
the parabolic rank two Higgs bundle $(E_{\bullet},\theta)$ on
$(C,\Par_{C})$ and the parabolic rank four Higgs bundles
$(F_{\bullet},\varphi)$ and $(F_{\bullet}',\varphi')$ on $(X,\Par_{X})$
defined in sections
\ref{sssec:eigenvalue} and \ref{sssec:eigensheaf} respectively  satisfy
the condition:  for all
  levels $(T,t)$ we have
\begin{equation} \label{eq:ab.condition}
{\textstyle q_{*}g_{*}g_{Z}^{*}(\bQ_{\lift(T,t)}(-
\sum_{\substack{\left\{ I \, | \, lp_{I} + e_{I} \in \mathbb{Z}\right\}}}
  LpE_{I} + \abH)) =
p_{X}^{*}F_{T}'\otimes p_{C}^{*}E_{t},}
\end{equation}
or equivalently, in terms of spectral data,
\begin{equation} \label{eq:ab.condition.spectral}
{\textstyle g_{Z}^{*}(\bQ_{\lift(T,t)}(-
\sum_{\substack{\left\{ I \, | \, lp_{I} + e_{I} \in \mathbb{Z}\right\}}}
LpE_{I}  + \abH)) = \abq^{*}(p_{Y}^{*}M_{T,T}'\otimes
p_{\sC}^{*}\mathfrak{a}_{t,t}),}
\end{equation}
where as usual $M_{T,T}' = \mycal{L}_{\mathfrak{a}}(e'E + d'(E+G))$.
\end{itemize}

\

\noindent
Note that:
\begin{itemize}
\item In this formulation we deliberately omitted the Higgs
fields on the two sides of the Hecke condition. The reason is that
from the spectral description we have already argued that the Higgs
fields will match generically, i.e. away from the exceptional loci of
the birational map $\abH \to \Jac\times \sC$. Thus if we know that the
underlying bundles match everywhere, then the Higgs fields will have
to match automatically.
\item As explained in Remark~\ref{rem:away.from.curves} it suffices to
  check that condition \eqref{eq:ab.condition} or equivalently condition 
  \eqref{eq:ab.condition.spectral} holds away from finitely many
  curves in $X\times C$ and their preimages in $Y\times \sC$, $H$, and
  $\abH$. Since $\abH$ and $Y\times \sC$ are isomorphic away from
  finitely many curves in $Y\times \sC$ this means that it is
  sufficient to verify the spectral version \eqref{eq:ab.condition.spectral}
  of the abelianized Hecke condition
  away from a codimension two locus in $Y\times \sC$.
\item Since the Fourier-Mukai bundle $\mycal{L}_{\mathfrak{a}}$ on the
  Jacobian satisfies the eigensheaf property
  $\aj^{*}\mycal{L}_{\mathfrak{a}} = \mycal{L}_{\mathfrak{a}}\boxtimes
  \mathfrak{a}$, it follows that the dependence on $\mathfrak{a}$
  decouples from \eqref{eq:ab.condition} and
  \eqref{eq:ab.condition.spectral}. In other words, the condition
  \eqref{eq:ab.condition} on $E_{t}$, $F_{T}$, $F_{T}'$, and
  $\bQ_{\lift(T,t)}$ (equivalently the condition
  \eqref{eq:ab.condition.spectral} on $\mathfrak{a}_{t,t}$, $M_{T,T}$,
  $M_{T,T}'$) correspondoing to a degree zero line bundle
  $\mathfrak{a}$ on $\sC$ is equivalent to the same condition on
  $E_{t}$, $F_{T}$, $F_{T}'$ and $\bQ_{\lift(T,t)}$ (equivalently the
  condition \eqref{eq:ab.condition.spectral} on $\mathcal{O}_{t,t}$,
  $M_{T,T}$, $M_{T,T}'$, and $\bQ_{\lift(T,t)}$) corresponding to the
  trivial line bundle $\mathcal{O}_{\sC}$.
\end{itemize}

\

\noindent
To summarize: the Hecke condition is equivalent to
\eqref{eq:ab.condition.spectral} which in turn is equivalent to linear
equivalence relations on  divisors on the threefold $Y\times
\sC$. Furthermore since the dependence on $\mathfrak{a} \in
\Pic^{0}(\sC)$ cancels in the equation
\eqref{eq:ab.condition.spectral}, this linear equivalence condition
becomes simply an equation on the unknown real variables $lp$, $lq$,
$llc$, $r_{1}$, $r_{2}$, $e$, $d$, $e'$, and $d'$ which we proceed to
describe.

\subsection{Linear equivalence conditions} \label{ssec:Pic.conditions}

To write down the equation in $\Pic(Y\times \sC)$ equivalent to
\eqref{eq:ab.condition.spectral} we will need notation for the natural
divisors on the smooth threefolds $Y\times \sC$ and $\abH$. Note that
the natural divisors on $\abH$ fall into two classes: the exceptional
divisors for the blowup map $\abq : \abH \to Y\times\sC$, and the
pullbacks of the natural divisors on $Y\times\sC$:
\begin{itemize}
\item For every $I \in \evenL$ and every  $i = 1, 2, \ldots, 6$ we will
  write $N_{I}$ and $M_{I.i}$ for the exceptional divisors
  corresponding to the curves $\Gamma_{I}$ and
  $E_{I}\times \tilde{p}_{i}$ respectively.
\item For every $I \in \evenL$ we will write
  $\abq^{*}(E_{I}\times \sC)$,  $\abq^{*}(G_{I}\times \sC)$ for the
  $\abq$-pullbacks of the divisors
$E_{I}\times \sC$ and $G_{I}\times \sC$, and will write
$\hatECI$ and $\hatGCI$ for the
  $\abq$-strict transforms of
$E_{I}\times \sC$ and $G_{I}\times \sC$.
\item For every $I \in \evenL$ we will write $\hatAJ^{*}(\Theta_{I})$
  for the pullback via $\hatAJ = \AJ\circ (\bmu\times \mathbf{1}) :
  Y\times \sC \to \Jac$ of the theta divisor $\Theta_{I}$ introduced
  in Section~\ref{ssec:nilpotent.residues}, and will write $S_{I}$ for
  the $\abq$-strict transform of the $\bmu\times \mathbf{1}$-strict
  transform of the divisor $\AJ^{*}(\Theta_{I})$.
\item For every $i = 1, \ldots, 5$ we will write 
  $\abq^{*}(Y\times \tilde{p}_{i})$ for the $\abq$-pullback of the divisor
  $Y\times \tilde{p}_{i}$ and $\hatYpi$ for
  the $\abq$-strict transform of $Y\times \tilde{p}_{i}$.
\end{itemize}

From the definition of these divisors, the description of the basic
divisors on $H$,  and the key diagram 
\[
  \xymatrix{
& \abH \ar[dl]_-{g_{Z}} \ar[dr]^-{g} & \\
Z \ar[d]_-{p_{Z}} \ar[rr]^-{f_{Z}} & & H \ar[d]^-{p} \\
Y \ar[rr]_-{f} & & X }
\]
we can compute the pullbacks of the standard divisors on $H$ and
$Y\times \sC$ in terms of the natural divisors on $\abH$. A
straightforward calculation, which we leave to the reader, gives:

\

\begin{equation} \label{eq:abp.and.abq.pullbacks}
  \boxed{
\begin{array}{lcll}
  \abp^{*}E_{I} & = &  N_{I} + \sum_{i = 1}^{6} M_{I,i},
  & \text{for } I \in \evenL,
  \\[+0.5pc]
  \abp^{*}G_{I} & = & S_{I} + \hatECI,
  & \text{for } I \in \evenL, \\[+0.5pc]
  \abq^{*}(E_{I}\times \sC) & =  & \hatECI + \sum_{i=1}^{6}
  M_{I,i},
& \text{for } I \in \evenL,
  \\[+0.5pc]
  \abq^{*}(G_{I}\times \sC) & = & \hatGCI + N_{I},
  & \text{for } I \in \evenL,
  \\[+0.5pc]
  \abq^{*}(Y\times \tilde{p}_{i}) & = & \hatYpi +
  \sum_{I\in \evenL}  M_{I,i},
  & \text{for } i  = 1, \ldots, 5.
\end{array}
  }
\end{equation}

\

and

\begin{equation}
  \label{eq:g.pullbacks}
  \boxed{
\begin{array}{lcll}
  g^{*}Lp_{I} & = & S_{I} + N_{I} + 2 \sum_{i = 1}^{5} M_{I+i,i},
& \text{for } I \in \evenL,
  \\[+0.5pc]
  g^{*}Lq_{I} & = & \hatECI + \hatGCI + 2 \sum_{i=1}^{5} M_{I,i},
& \text{for } I \in \evenL,
  \\[+0.5pc]
  g^{*}LLC_{I} & = & \hatECI + N_{I} +  2M_{I,6},
& \text{for } I \in \evenL,
  \\[+0.5pc]
  g^{*}R_{1,i} & = & \hatYpi +
  2\sum_{I\in \evenL} \incdelta_{I,i} M_{I,i},
& \text{for } i = 1, \ldots, 5,
  \\[+0.5pc]
  g^{*}R_{2,i} & = & \hatYpi +
  2\sum_{I\in \evenL} (1-\incdelta_{I,i}) M_{I,i},
& \text{for } i = 1, \ldots, 5.
\end{array}
  }
\end{equation}

\

\noindent
Here we use the notation $\incdelta_{I,i}$ for the incidence function
\[
\incdelta_{I,i} = \begin{cases} 1 & \text{if } i \in I, \\
0 & \text{if } i \notin I,
\end{cases}
\]
where $i \in \{1,2,3,4,5\}$ and $I \subset \{1,2,3,4,5\}$ is a subset
of even cardinality.

\

\noindent
Computing the $g_{Z}$-pullbacks of the standard divisors on $Z$ is
slightly more subtle. The key observation is that because we are
pulling back Weil divisors which in general are not Cartier, the
natural pullbacks are well defined by intersections and so 
make sense as elements in the rational Chow group of codimension one
cycles or equivalently as elements in $\Pic(\abH)\otimes
\mathbb{Q}$.

\

\begin{ex} \label{ex:cone.Weil}
As a model example consider a two dimensional quadratic cone
\[
\mathsf{Cone} : z^{2} = xy,
\]
in a three dimensional complex projective space with homogeneous
coordinates $(x:y:z:w)$.  Let $\mathsf{Rul}, \mathsf{Rul}' \subset
\mathsf{Cone}$ be two (possibly coincident) rulings, i.e.
$\mathsf{Rul}$ and $\mathsf{Rul}'$ are straingt lines in
$\mathbb{P}^{3}$, contained in the cone, and passing through the
vertex. The rulings $\mathsf{Rul}$ and $\mathsf{Rul}'$ are Weil
divisors on $\mathsf{Cone}$ with the property that $2\mathsf{Rul}$,
$2\mathsf{Rul}'$, and $\mathsf{Rul} + \mathsf{Rul}'$ are all Cartier
and belong to the linear system
$\mathcal{O}_{\mathbb{P}^{3}}(1)_{|\mathsf{Cone}}$. Since
$\mathcal{O}_{\mathbb{P}^{3}}(1)_{|\op{Rul}}$ is the hyperplane line
bundle on the projective line $\mathsf{Rul}$, it follows that the
intersections of $2\mathsf{Rul}$ and $2\mathsf{Rul}'$ with
$\mathsf{Rul}$ are linearly equivalent and are both represented by a
class of a point.  In particular both the self intersection
$\mathsf{Rul}\cdot \mathsf{Rul}$ and the intersection
$\mathsf{Rul}'\cdot \mathsf{Rul}$ belong to the rational Chow group of
codimension one cycles on $\mathsf{Rul}$ and are both represente by
the class $(1/2)\pt$. In summary, if $\iota : \mathsf{Rul}
\hookrightarrow \mathsf{Cone}$ is the inclusion, we get that the
$\iota$-pullbacks of the Weil divisors $\mathsf{Rul}, \mathsf{Rul}'
\subset \mathsf{Cone}$ are the rational Cartier divisors
\[
\begin{aligned}
\iota^{*}\mathsf{Rul} & = \frac{1}{2}\pt \in
\Pic(\mathsf{Rul})\otimes \mathbb{Q}, \\[+0.5pc]
\iota^{*}\mathsf{Rul}'
& = \frac{1}{2}\pt \in \Pic(\mathsf{Rul})\otimes \mathbb{Q}.
\end{aligned}
\]
\end{ex}

\

\medskip

\

\noindent
Going back to our situation in which we want to compute the $g_{Z} :
\abH \to Z$ pullback of the various Weil divisors on $Z$ is a relative
version of this model example. Indeed by Lemma~\ref{lem:sing.of.Z} we
know that away from finitely curves the fourfold $Z$ has sixteen
disjoint surfaces $\Xi_{I}$ of transversal $A_{1}$-surface
singularities, and that $LpE_{I}$ and $LLE_{I}$ are all Weil divisors
on $Z$ which are Cartier away from $\Xi_{I}$, each of them passes
through $\Xi_{I}$, and transversally to $\Xi_{I}$ is represented by a
ruling. Also the divisor $\abH \subset Z$ is a Weil divisor in $Z$
which modulo Cartier divisors is transversally represented by a
ruling. This can be seen locally from the definition of $Z$ and $\abH$
or by noting that by Lemma~\ref{lem:divisor.of.beta} the sum $LpE+\abH$ is
Cartier.

Taking this into account and tracing the relations among the $g$ and
$g_{Z}$ pullbacks in the basic diagram we get:

\begin{equation}
  \label{eq:gZ.pullbacks}
  \boxed{
\begin{array}{lcll}
  g_{Z}^{*}LpE_{I} & = & \frac{1}{2}N_{I}
  + \sum_{i = 1}^{5} M_{I+i,i},
  & \text{for } I \in \evenL, 
  \\[+0.5pc]
  g_{Z}^{*}LLE_{I} & = & \frac{1}{2}N_{I} + M_{I,6},
  & \text{for } I \in \evenL, 
  \\[+0.5pc]
  g_{Z}^{*}LpG_{I} & = & S_{I},
  & \text{for } I \in \evenL, 
  \\[+0.5pc]
  g_{Z}^{*}LLG_{I} & = & \hatECI,
    & \text{for } I \in \evenL.
\end{array}
  }
\end{equation}

\

\noindent
Note also that in the basic diagram we have $g = f_{Z}\circ g_{Z}$ and
so $g^{*} = g_{Z}^{*}\circ f_{Z}^{*}$.  Combined with
\eqref{eq:g.pullbacks} this gives formulas for
$g_{Z}^{*}(f_{Z}^{*}Lq_{I})$, $g_{Z}^{*}(f_{Z}^{*}R_{1,i})$, and
$g_{Z}^{*}(f_{Z}^{*}R_{2,i})$.

\

\noindent
Next observe that by definition the pullbacks
$g_{Z}^{*}\bQ_{\lift(T,t)}$ and $g_{Z}^{*}LpE$, as well as
$\abq^{*}(p_{Y}^{*}M_{T,T})$ and $\abq^{*}(p_{\sC}^{*}\mathfrak{a}_{t,t})$ can
all be written uniquely as a sum: combination of
$\abq^{*}(\text{divisors on } Y\times \sC)$ plus a combination of
$\abq$-exceptionals. Furthermore, in this expression we can take the
first summand to be integral combinations of the $\abq$-pullbacks of
$\{E_{I}\times \sC\}_{I \in \evenL}$, $\{G_{I}\times \sC\}_{I \in
  \evenL}$, $\{\hatAJ^{*}\Theta_{I}\}_{I \in \evenL}$, and $\{Y\times
\tilde{p}_{i}\}_{i=1}^{5}$.

In particular, to get the desired rewrite of
\eqref{eq:ab.condition.spectral} as linear equivalence conditions we
only need to express $g_{Z}^{*}\mathcal{O}_{Z}(\abH)$ as a linear
combination of the divisors $\{E_{I}\times \sC\}_{I \in \evenL}$,
\linebreak
$\{G_{I}\times \sC\}_{I \in \evenL}$, $\{\hatAJ^{*}\Theta_{I}\}_{I \in
  \evenL}$, and $\{Y\times \tilde{p}_{i}\}_{i=1}^{5}$ plus
$\abq$-exceptionals. The formula for $g_{Z}^{*}\mathcal{O}_{Z}(\abH)$
is complicated and will not be symmetric in $I$ and $i$ if we insist
on combinations with integral coefficients. However if use the natural
linear equivalence relations among the divisors $E_{I}$ and $G_{I}$ on
$Y$, and among the points $\tilde{p}_{i}$ viewed as divisors on $\sC$,
we can write a symmetric formula for this pullback.

\begin{lem} \label{lem:pullback.of.abH}
  As an element in $\Pic(\abH)\otimes \mathbb{Q}$ the
  $\mathbb{Q}$-Cartier divisor $g_{Z}^{*}\abH$ is given by
  the following
  combination:
  \begin{equation} \label{eq:pullback.of.abH}
    \boxed{
   \begin{aligned}
    g_{Z}^{*}\abH & = -\frac{1}{8}\sum_{I\in \evenL} S_{I} +
        \frac{5}{8}\sum_{I \in \evenL} \hatECI +
      \frac{1}{8} \sum_{I \in \evenL} \hatGCI +
      \frac{3}{5}\sum_{i = 1}^{5} \hatYpi \\[+0.5pc]
      & \qquad\qquad +  \frac{3}{8}
      \sum_{I\in \evenL} N_{I} +
      \frac{3}{5}\sum_{\substack{I\in \evenL \\ i = 1, \ldots 5}}
      M_{I,i} + M_{I,6}.
   \end{aligned}
    }
  \end{equation}
\end{lem}

\

\noindent
{\bfseries Proof.} \ By Lemma~\ref{lem:divisor.of.beta} we have
$\mathcal{O}_{Z}(LpE + LpG + \abH) = f_{Z}^{*}\left(
\omega_{q}(Lp)\right)$.  Since $f_{Z}^{*}Lp = 2LpE + LpG$ this implies
that $g_{Z}^{*}\mathcal{O}_{Z}(\abH) =
g^{*}\omega_{q}(g_{Z}^{*}LpE)$. Since we already have a formula
\eqref{eq:gZ.pullbacks} for the divisors $g_{Z}^{*}LpE_{I}$, we only
need to express the line bundle $g^{*}\omega_{q} \in \Pic(\abH)$
as a combination of the standard divisors on $\abH$.

By definition we have $\omega_{q} = \omega_{H} \otimes q^{*}
\omega_{X\times C}^{-1}$. But all points are linearly equivalent on
$C$, and the sum of all lines in $X$ is linearly equivalent to four
times the anti-canonical class of $X$. Therefore we can express
$\omega_{X\times C}^{-1}$ as a symmetric linear combination of the
components of $\Par_{X}$, i.e.
\[
\omega_{X\times C}^{-1} = p_{X}^{*}\omega_{X}^{-1}\otimes
p_{C}^{*}\omega_{C}^{-1} = \mathcal{O}_{X\times C}\left( \frac{1}{4}
\sum_{I \in \evenL} L_{I}\times C + \frac{2}{5} \sum_{i =1}^{5}
X\times p_{i}\right).
\]
But as we saw in section~\ref{ssec:pardivisors} we have
$q^{*}(L_{I}\times C) = Lq_{I} + LLC_{I}$, and $q^{*}(X\times p_{i}) =
RR_{p_{i}} = R_{1,i} + R_{2,i}$. Setting $RR := \sum_{i= 1}^{5}RR_{p_{i}}$ we get
\[
q^{*}\omega_{X\times C}^{-1} = \mathcal{O}_{H}\left( \frac{1}{4} Lq +
\frac{1}{4}LLC + \frac{2}{5} RR\right).
\]
To compute $\omega_{H}$ we will use the definition of $H$ as an
iterated blow-up:
\[
H = \op{Bl}_{\sqcup \widehat{L_{I}\times L_{I}}} \op{Bl}_{\Delta}(X\times X).
\]
Write $\varpi : H \to X\times X$ for the total blow-up map, and
denote by $\widehat{\Delta}$ the strict transform in $H$ of the
exceptional divisor of the first blow-up map
$\op{Bl}_{\Delta}(X\times)$. Then
\[
\omega_{H} = \left(\varpi^{*}\omega_{X\times X}\right)(\widehat{\Delta} + LLC).
\]
But
\[
\omega_{X\times X} = \mathcal{O}_{H}\left( -\frac{1}{4} \sum_{I}
L_{I}\times X - \frac{1}{4} \sum_{I} X \times L_{I}\right),
\]
and since $\varpi^{*}(L_{I}\times X) = Lp_{I}+LLC_{I}$ and $\varpi^{*}(X\times
L_{I}) = Lq_{I}+LLC_{I}$ we get
\[
\omega_{H} = \mathcal{O}_{H}\left( -\frac{1}{4} Lp - \frac{1}{4}Lq +
\frac{1}{2} LLC + \widehat{\Delta}\right),
\]
and so for the $q$-relative dualizing sheaf we get
\[
\omega_{q} = \mathcal{O}_{H}\left( -\frac{1}{4}Lp + \frac{3}{4} LLC +
\frac{2}{5} RR + \widehat{\Delta} \right).
\]
Since $\Pic(H)$ is spanned by the divisor components of $\Par_{H}$ we
can also express $\widehat{\Delta}$ as a combination of these divisor
components up to linear equivalence.

To figure out this expression it suffices to intersect
$\widehat{\Delta}$ and the $58$ divisor components of $\Par_{H}$ with a
collection of curves in $H$ that span $H_{2}(H,\mathbb{Z})$. We will
use the follwing curves:
\begin{itemize}
\item $L_{I}\times \pt$, where $\pt \in X - \cup_{I} L_{I}$.
\item $\pt\times L_{I}$, where $\pt \in X - \cup_{I} L_{I}$.
\item $\widehat{\Delta}_{\pt}$: the fiber of $\widehat{\Delta}
  \to \Delta$ over a point $\pt \in X - \cup_{I} L_{I}$.
\item $\Sigma_{I}$: the fiber of $LLC_{I} \to L_{I}\times L_{I}$ over
  a general point in $L_{I}\times L_{I}$.
\end{itemize}
The intersections of these curves with the standard divisors are
easily extracted from the definitions and can be organized in the
following table:

\

\begin{center}
\begin{tabular}{|c||c|c|c|c|c|c|}
\hline \hline
$ { } $ & $\widehat{\Delta}$ & $Lp_{J}$ & $Lq_{J}$ & $LLC_{J}$ & $R_{1,i}$
& $R_{2,i} $\\ 
\hline \hline
$L_I \times \pt$ & $0$ & $\bI_{I,J}$ & $0$ & $0$ & $1-\incdelta_{I,i}$
& $\incdelta_{I,i} $ \\
\hline
$\pt \times L_I$ & $0$ & $0$ & $\bI_{I,J}$ & $0$ & $1-\incdelta_{I,i}$
& $\incdelta_{I,i} $ \\
\hline
$\widehat{\Delta}_{\pt}$ & $-1$ & $0$ & $0$ & $0$ & $1$ & $1$  \\ 
\hline
$\Sigma_I$ & $0$ &  $\one_{I,J}$ &  $\one_{I,J}$ &  $-\one_{I,J}$ & $
\incdelta_{I,i}$ & $1-\incdelta_{I,i} $ \\
\hline 
\end{tabular}
\end{center}

\

\noindent
where as above $\incdelta_{I,i}$ denotes the incidence function between a
label and subset of even cardinality, and  $\one$ denotes  the identity
$16\times 16$ matrix, while $\bI$ denotes the intersection matrix of
the configuration of $16$ lines in $X$.

From these intersection patterns it immediately follows that in
$\Pic(H)$ we have
\[
\widehat{\Delta} = \frac{1}{8} Lp + \frac{1}{8}Lq - \frac{1}{4} LLC -
\frac{1}{10} RR.
\]
Substituting in the above formula for $\omega_{q}$ we get
\[
\omega_{q} = \mathcal{O}_{H}\left(-\frac{1}{8} Lp +
\frac{1}{8}Lq + \frac{1}{2} LLC +
\frac{3}{10} RR\right),
\]
which together with the formulas \eqref{eq:g.pullbacks} and
\eqref{eq:gZ.pullbacks} gives the claimed expression for
$g_{Z}^{*}\abH$ and proves the lemma.
\ \hfill $\Box$

\

\medskip

\noindent
With this formula at our disposal we are now ready to write the linear
equivalence version of the  abelianized Hecke
condition \eqref{eq:ab.condition.spectral}.

\begin{prop} \label{prop:Pic.Hecke}
For any choice of $a, b \in \mathbb{R}^{5}$ with $\sum_{i= 1}^{5}
(a_{i} + b_{i}) = 3$, there exists natural vectors $lp, lq, llc, e, d,
e', d' \in \mathbb{R}^{16}$, and $r_{1},r_{2} \in \mathbb{R}^{5}$ satisfying:

\

\noindent
For any smooth spectral curve $\sC \to C$ branched over
  $\Par_{C}$ and any line bundle $\mathfrak{a} \in \Pic^{0}(\sC)$, the
parabolic rank two Higgs bundle $(E_{\bullet},\theta)$ on $(C,\Par_{C})$
and the parabolic rank four Higgs bundles
  $(F_{\bullet},\varphi)$  and $(F_{\bullet}',\varphi')$ on $(X,\Par_{X})$
  defined by
\[
\begin{aligned}
  (E_{t},\theta) & = \left(\pi_{*}\mathfrak{a}_{t,t},
  \pi_{*}(-\otimes \lambda)\right) \\
& = 
  \left(\pi_{*}
  \mathfrak{a}((\fl{t+a}+\fl{t+b})\widetilde{P}),
  \pi_{*}(-\otimes \lambda)\right), \\[+0.7pc]
(F_{T},\varphi) & = \left(f_{*}M_{T,T}, f_{*}(-\otimes \alpha)\right)  \\
  & = \left(f_{*}\mycal{L}_{\mathfrak{a}}((\fl{T+e'} + \fl{T+d'})E +
  \fl{T+d'}G),f_{*}(-\otimes \alpha)\right), \\[+0.7pc]
(F_{T}',\varphi') & = \left(f_{*}M_{T,T}', f_{*}(-\otimes \alpha)\right)  \\
  & = \left(f_{*}\mycal{L}_{\mathfrak{a}}((\fl{T+e'} + \fl{T+d'})E +
  \fl{T+d'}G),f_{*}(-\otimes \alpha)\right).
\end{aligned}
\] 
satisfy the parabolic Hecke eigensheaf condition \eqref{eq-par.hecke}
(equivalently the abelianized Hecke condition
\eqref{eq:ab.condition.spectral}) if and only if the vectors $lp$,
$lq$, $llc$, $e$, $d$, $r_{1}$ and $r_{2}$ are such that for all
levels $(T,t)$ the divisor
  \begin{equation}  \label{eq:Hecke.divisorial}
    \boxed{
{
\begin{aligned} 
  & \sum_{I\in \evenL}  \left(\fl*{lp_{I} +d_{I}} 
  -\frac{1}{8}\right)\hatAJ^{*}\Theta_{I} \\ 
  & \quad + \sum_{I \in \evenL}
  \left(\fl*{T_{I} +llc_{I} + d_{I}} + \fl*{T_{I} +lq_{I}} +
  \frac{5}{8} - \fl*{T_{I} +e_{I}'}
  - \fl*{T_{I} +d_{I}'}\right)E_{I}\times \sC \\
& \quad + 
\sum_{I\in \evenL} \left(\fl*{T_{I} +lq_{I}} + \frac{1}{8} -
\fl*{T_{I} +d_{I}'}\right)G_{I}\times \sC \\
& \quad + 
\sum_{i=1}^{5} \left(\fl*{t_{i} +r_{1,{i}}} + \fl*{t_{i} +r_{2,{i}}}
+ \frac{3}{5} - \fl*{t_{i} +a_{i}}
- \fl*{t_{i} +b_{i}}\right)Y\times \tilde{p}_{i}
\end{aligned}
}}
  \end{equation}
  is linearly equivalent to zero on $Y\times \sC$.
\end{prop}

\

\noindent
{\bfseries Proof.} The abelianized (spectral cover) Hecke condition
\eqref{eq:ab.condition.spectral} is an equation in $\Pic(\abH)\otimes
\mathbb{Q}$. As explained in section~\ref{sssec:ab.Hecke.conditions}
it is enough to check this equation for $\mathfrak{a}
=\mathcal{O}_{\sC}$ and away from closed subvarieties in $\abH$ that
map to codimension two loci in $Y\times \sC$. In particular, since the
exceptional divisors $M_{I,i}$ and $N_{I}$ for the birational morphism
$\abq : \abH \to Y\times \sC$ contract to curves we can check the
respective equations modulo the span of those divisors.

Since by
\eqref{eq:gZ.pullbacks} $g_{Z}^{*}LpE_{I}$ is a combination of $M_{I+i,i}$'s
and $N_{I}$'s we see that
\[
{\textstyle g_{Z}^{*} (\bQ_{\lift(T,t)}(-
\sum_{\substack{\{I \, | \, lp_{I} + e_{I} \in \mathbb{Z}\}}}LpE_{I} +\abH)) = \\[+0.5pc]
   = g_{Z}^{*}(\bQ_{\lift(T,t)}(\abH))} 
\]
modulo the span the $\abq$-exceptional
divisors. Finally note that modulo the $\abq$-exceptional divisors the
$\abq$-pullback of any divisor on $Y\times \sC$ is equal to its $\abq$
strict transform. Therefore working modulo the $M_{I,i}$'s and $N_{I}$'s we get:
\[
\begin{aligned}
  g_{Z}^{*}& {\textstyle (\bQ_{\lift(T,t)}(-\sum_{\substack{\{I \, |
        \, lp_{I} + e_{I} \in \mathbb{Z}\}}}LpE_{I} +\abH))} =
  \\[+0.5pc] & = g_{Z}^{*}(\bQ_{\lift(T,t)}(\abH)) \\[+0.5pc] & =
  (g_{Z}^{*}\bQ_{\lift(T,t)})\otimes
  \mathcal{O}_{\abH}(g_{Z}^{*}\abH)) \\[+0.5pc] & =
  (g_{Z}^{*}\bQ_{\lift(T,t)})\otimes \mathcal{O}_{\abH}\left(
  -\frac{1}{8} \sum_{I\in \evenL}S_{I} + \frac{5}{8}\sum_{I\in
    \evenL}\hatECI + \frac{1}{8}\sum_{I\in \evenL}\hatGCI +
  \frac{3}{5}\sum_{i =1}^{5} \hatYpi \right) \\[+0.5pc] & =
  (g_{Z}^{*}\bQ_{\lift(T,t)})\otimes \\[+0.5pc] & \qquad\quad
  \abq^{*}\mathcal{O}_{Y\times \sC}\left( -\frac{1}{8} \sum_{I\in
    \evenL}\hatAJ^{*}\Theta_{I} + \frac{5}{8} \sum_{I\in \evenL}
  E_{I}\times \sC + \frac{1}{8}\sum_{I\in \evenL} G_{I}\times \sC +
  \frac{3}{5} \sum_{i=1}^{5} Y\times \tilde{p}_{i} \right).
\end{aligned}
\]
Also, specializing the formula \eqref{eq:spectral.sheaf} to
$\varepsilon = 0$ and $\mathfrak{a} = \mathcal{O}_{\sC}$ and taking into account
the pullback formulas \eqref{eq:gZ.pullbacks} we get that modulo
$M_{I,i}$'s and $N_{I}$'s we have
\[
\begin{aligned}
  g_{Z}^{*}\bQ_{\lift(T,t)} & = \mathcal{O}_{\abH}\left( \fl*{T +
    llc + d}g_{Z}^{*}LpG + \fl*{lp + d} g_{Z}^{*}LLG \right. \\[+0.5pc] &
  \qquad\quad \left. + \fl*{T+llc+d} g^{*}Lq  +
  \fl*{t+r_{1}} g^{*}R_{1} + \fl*{t+r_{2}}g^{*}R_{2} \right)
  \\[+0.5pc] & = \mathcal{O}_{\abH}\left( \fl*{lp + d}S + \fl*{T + llc +
    d}\widehat{\mathsf{EC}} \right. \\[+0.5pc] & \qquad\quad \left. +
  \fl*{T+lq}(\widehat{\mathsf{EC}} + \widehat{\mathsf{GC}}) +
 (\fl*{t+r_{1}} +
  \fl*{t+r_{2}})\widehat{\mathsf{Yp}} \right) \\[+0.5pc] & =
  \abq^{*}\mathcal{O}_{Y\times \sC}\left(\fl*{lp+d}\hatAJ^{*}\Theta +
  (\fl*{T+lq} +  \fl*{T + llc +
    d}) E\times \sC \right. \\[+0.5pc]
  & \qquad\qquad\quad \left. + \fl*{T+lq}  G\times \sC +
  (\fl*{t+r_{1}} +
  \fl*{t+r_{2}}) Y \times \tilde{P}  \right) .
\end{aligned}
\]
Similarly specializing the formulas \eqref{eq:parline} and
\eqref{eq:FT}  to $\mathfrak{a} = \mathcal{O}_{\sC}$ and taking into account
the pullback formulas \eqref{eq:abp.and.abq.pullbacks}
 we compute that modulo
$M_{I,i}$'s and $N_{I}$'s we have
\[
\begin{aligned}
  \abq^{*}\left(p_{Y}^{*}M_{T,T}'\otimes p_{\sC}^{*}\mathfrak{a}_{t,t}\right) & =
 \abq^{*}\mathcal{O}_{Y\times \sC}\left( \left( \fl*{T + e'} + \fl*{T +
     d'}\right) E\times  \sC  + \fl*{T + d'} G\times \sC
\right. \\[+0.5pc]
& \qquad\qquad\quad \left. + \left( \fl*{t + a} + \fl*{t + b}\right)
  Y\times \widetilde{P} \right).
\end{aligned}
\]
Now setting the expressions for $g_{Z}^{*}\bQ_{\lift(T,t)}(\abH)$ and
$\abq^{*}\left(p_{Y}^{*}M_{T,T}'\otimes p_{\sC}^{*}\mathfrak{a}_{t,t}\right)$ to
be equal to each other and taking into account that $\abq^{*}$ gives
an isomorphism between $\Pic(Y\times \sC)$ and
$\Pic(\abH)/\langle M_{I,i}, N_{I} \rangle$ we arrive at the condition
\eqref{eq:Hecke.divisorial}. This concludes the proof of  the
proposition. 
\ \hfill $\Box$

\section{Solving the constraints} \label{sec:solution}

After all this preliminary work we are now ready to solve all
constraints on $e$, $d$, $e'$, $d'$, and $\bzeta =
(lp,lq,llc,r_{1},r_{2})$ so that both Mochizuki's non-abelian Hodge
theory condition and the Hecke eigensheaf condition are satisfied.
\index{terms}{Non-abelian Hodge theory}

\subsection{Chern characters}

In this section we calculate the Chern characters of $\LL$ and $F_T$.
\index{terms}{parabolic!Chern character}
Recall our notation:
$d,e,d',e',T$ are 16-dimensional, real-valued column vectors, as is
$\sigma$ which is the vector with all entries equal to $1$.  To record
the calculations more easily we will also introduce the shortcut
notation
\[
\bba := \fl*{T+e}, \ \bbb := \fl*{T+d},
\]
The vectors $\bba,\bbb,d,e,T,\sigma$ can all be paired with row
vectors $L,E,G$ whose entries are cohomology (or rational equivalence)
classes on $X,Y,Y$ respectively.

On $Y$ we have the line bundle
\[ \LL:=\OO_Y(E\bba+(E+G)\bbb),
\] and we let $F_T$ be its direct image on $X$:
\[ F_T := f_* \LL,
\]
where $f:Y \to X$ is our degree 4 map.  We want to calculate the
Chern classes of $F_T$. A straighforward calculation on the del Pezzo
surface $X$ and the blown-up abelian surface $Y$  gives the following
formulas for the Todd
classes:
\[
\begin{aligned} \op{Td}(Y) & = \textsf{1}_{Y} - \frac{1}{2} E \sigma,
\\ \op{Td}(X)^{-1} & = \textsf{1}_{X} - \frac{1}{8} L \sigma,
\end{aligned}
\]
and the intersection pairings on $Y$, as given in formula
\eqref{intersection pattern}:
\[
  E^{t} E =- \one, \quad G^{t} E =2\one + \bI, \quad G^{t} G =-4\one,
\]
expressed in terms of the $16\times 16$ identity matrix $\one$ and
the intersection matrix
\[ \bI := L^{t} L
\]
of the 16 lines on the del Pezzo $X$.

\

\begin{rem} \label{rem-degree-convention} In the formulas for the Todd
classes and in all Chern class calculations below we use the
convention in which a cohomology class on an algebraic surface $S$ is
written as a sum of three pieces which are classes in
$H^{0}(S,\mathbb{Q})$, $H^{2}(S,\mathbb{Q})$, and
$H^{4}(S,\mathbb{Q})$ respectively. The first piece in
$H^{0}(S,\mathbb{Q})$ will always be written as an integer multiple of
the Poincare dual $\textsf{1}_{S} \in H^{0}(S,\mathbb{Q})$ of the
fundamental class. The second piece in $H^{2}(S,\mathbb{Q})$ will be a
rational linear combination of divisors, and the third piece in
$H^{4}(S,\mathbb{Q})$ will be a rational multiple of the class $\pt
\in H^{4}$ of a point. So for instance in the formula for $\op{Td}(Y)$
above, the first piece is just the fundamental class $\textsf{1}_{Y}$,
the second piece is the rational divisor $-(1/2)E\sigma = -
(1/2)\sum_{I} E_{I}$, and the third piece which is supressed from the
notation is the zero multiple of $\pt$.
\end{rem}

\

\medskip

\noindent Upstairs we find:
\[
\begin{aligned} ch(\LL) & = e^{E\bba+(E+G)\bbb} \\ & = \textsf{1}_{Y} +
  \left(E\bba+(E+G)\bbb\right) \\ & \hspace{1in} + \frac{1}{2}
  \left(-\bba^t \bba -
  2 \bba^t \bbb + 2 \bba^t (2\cdot\one + \bI) \bbb + \bbb^t (-\one
  -4\cdot\one +4\cdot\one +
2\bI)\bbb\right)\pt \\ & = \textsf{1}_{Y} + \left(E\bba+(E+G)\bbb\right) +
\frac{1}{2} \left(-(\bba-\bbb)^2 +2 (\bba+\bbb)^t \bI \bbb\right)\pt.
\end{aligned}
\] Therefore downstairs:
\[
\begin{aligned} ch(F_T) & = \op{Td}(X)^{-1} f_* \left[ch(\LL)
\op{Td}(Y)\right] \\[+0.5pc] & = \Big(\textsf{1}_{X} - \frac{1}{8} L
\sigma\Big) f_*\Big[\Big(\textsf{1}_{Y} - \frac{1}{2} E \sigma\Big)
\Big(\textsf{1}_{Y} + E\bba+(E+G)\bbb \\[+0.5pc] & \hspace{2.6in} +
\frac{1}{2} \Big(-(\bba-\bbb)^2 +2 (\bba+\bbb)^t \bI \bbb\Big) \pt \Big)\Big]
\\[+0.5pc] & = \Big(\textsf{1}_{X} - \frac{1}{8} L \sigma\Big)
f_*\Big[ \textsf{1}_{Y} + E\bba+(E+G)\bbb + \frac{1}{2} \Big(-(\bba-\bbb)^2 +2
(\bba+\bbb)^t \bI \bbb\Big) \pt \\[+0.5pc] & \hspace{2.6in} - \frac{1}{2} E
\sigma (\textsf{1}_{Y} + E\bba+(E+G)\bbb) \Big] \\[+0.5pc] & =
\Big(\textsf{1}_{X} - \frac{1}{8} L \sigma\Big) f_*\Big[
\textsf{1}_{Y} + E(\bba+\bbb- \frac{\sigma}{2}) +G\bbb \\[+0.5pc] &
\hspace{2.2in} + \frac{1}{2} \Big(-(\bba-\bbb)^2 + 2(\bba+\bbb)^t\bI \bbb +
\sigma^t
(\one (\bba -\bbb) - \bI \bbb )\Big) \pt\big] \\[+0.5pc] & = \Big(\textsf{1}_{X}
- \frac{1}{8} L \sigma\Big) \Big[4\cdot \textsf{1}_{X} + L\Big(\bba+3\bbb-
\frac{\sigma}{2}\Big) \\[+0.5pc] & \hspace{2.2in} + \frac{1}{2}
\Big((\bba-\bbb)^t \sigma - (\bba-\bbb)^2 +(2\bba+2\bbb-\sigma)^t
\bI \bbb\Big) \pt \Big]
\\[+0.5pc] & = 4\cdot \textsf{1}_{X} + L\Big(\bba+3\bbb -\sigma\Big) +
\frac{1}{2} \Big(-8\bbb^t \sigma - (\bba-\bbb)^2 +2(\bba+\bbb)^t\bI
\bbb +8\Big) \pt,
\end{aligned}
\]
where in the last step we used the facts that $\sigma$ is an
eigenvector of $\bI$ with eigenvalue $4$ and that $\sigma^2 =16$.

In summary we have

\begin{equation} \label{eq-ch.of.FT} \boxed{
\begin{aligned} ch(\LL) & = \textsf{1}_{Y} + \left(E\bba+(E+G)\bbb\right) +
\frac{1}{2} \left(-(\bba-\bbb)^2 +2 (\bba+\bbb)^t \bI \bbb\right)\pt,
\\[+0.5pc] ch(F_T) & = 4\cdot \textsf{1}_{X} + L\Big(\bba+3\bbb
-\sigma\Big) + \frac{1}{2} \Big(-8\bbb^t \sigma - (\bba-\bbb)^2
+2(\bba+\bbb)^t\bI \bbb +8\Big) \pt.
\end{aligned} }
\end{equation}

\

\noindent Next we will plug this into the Iyer-Simpson formula
\eqref{Iyer-Simpson formula} from \cite{iyer-simpson} to compute the
parabolic Chern character of $F_{\bullet}$.
\index{terms}{parabolic!Chern character}

\subsection{Parabolic Chern characters}

By \eqref{Iyer-Simpson formula}, the parabolic Chern character of $F_{\bullet}$
is the ratio of two 16-dimensional integrals:
\[ \frac{\left[\prod_I \int_{-d_I}^{1-{d_I}} dT_I\right]
\left(ch(F_T)\exp(-LT)\right)} {\left[\prod_I \int_{0}^{1} dT_I\right]
\exp(-LT)}.
\]
We recall that in this formula, the denominator involves
integration from $0$ to $1$, while the numerator is periodic and can
be integrated over any unit length interval. We choose $[-d_I,1-d_I]$
since this is the interval on which $\bbb= \fl*{T+d} \equiv 0$, which
will be quite convenient in the next step.

In order to be able to apply the non abelian Hodge correspondence, we
need to set this parabolic Chern character equal to $4\cdot
\textsf{1}_{X}$. We are led to the condition:
\index{terms}{parabolic!Chern character}

\begin{equation} \label{ch} \left[\prod_I \int_{-d_I}^{1-{d_I}}
dT_I\right] \left(ch(F_T)\exp(-LT)\right) = 4 \left[\prod_I
\int_{0}^{1} dT_I\right] \exp(-LT).
\end{equation}

\noindent Since
\[ \exp(-LT) = \textsf{1}_{X} - LT + \frac{1}{2} T^t \bI T \pt,
\] we find that the RHS of \eqref{ch} is:

\[ \text{RHS} = 4\left(\textsf{1}_{X}-\frac{1}{2}L\sigma +
\frac{22}{3} \pt\right) = 4\cdot \textsf{1}_{X} - {2}L\sigma +
\frac{88}{3} \pt.
\]
Indeed, the definite integral over $[0,1]$ of $T_I$ with respect to
$dT_I$ is $\frac{1}{2}$, and integration with respect to the remaining
15 variables does not change this value. Integration of the various
terms in $T^t \bI T$ is carried out similarly: the $i$-th row contains
one diagonal term $-T_I^2$, which integrates to $-\frac{1}{3}$, and 5
non-zero off-diagonal entries $T_I T_J$ (occuring whenever $\bI_{I,J}
= 1$) which integrate to $\frac{1}{4}$.  It follows that the total
integral of $T^t \bI T$ is then $-\frac{1}{3} \cdot 16 + \frac{1}{4 }
\cdot 16 \cdot 5 = \frac{44}{3},$ completing the verification of the
RHS formula.

On the left hand side we have:
\[
\begin{aligned} ch(F_T)\exp(-LT) & = \\[+0.5pc] & \hspace{-0.5in} =
  \Big( 4\cdot \textsf{1}_{X} + L(\bba+3\bbb -\sigma) + \frac{1}{2}
  \Big(-8\bbb^t
\sigma - (\bba-\bbb)^2 +2(\bba+\bbb)^t \bI \bbb +8\Big) \pt\Big) \\[+0.5pc] &
\hspace{2in} \cdot \Big(\textsf{1}_{X} -LT + \frac{1}{2} T^t \bI T
\pt\Big) \\[+0.5pc] & \hspace{-0.5in} = 4\cdot \textsf{1}_{X} + L(\bba+3\bbb
-\sigma) + \frac{1}{2} \Big(-8\bbb^t \sigma - (\bba-\bbb)^2 +2(\bba+\bbb)^t\bI
\bbb
+8\Big) \pt \\[+0.5pc] & \hspace{1in} -4LT - (\bba+3\bbb -\sigma)^t \bI T
\pt +2 T^t \bI T \pt \\[+0.5pc] & \hspace{-0.5in} = 4\cdot
\textsf{1}_{X} + L(\bba+3\bbb -\sigma -4T) \\[+0.5pc] & + \Big[\frac{1}{2}
  \Big(-8b^t \sigma - (\bba-\bbb)^2 +2(\bba+\bbb)^t\bI \bbb +8\Big) -
  (\bba+3\bbb -\sigma -
2T)^t \bI T\Big] \pt,
\end{aligned}
\]
which simplifies when we substitute $\bbb=0$ to:
\[
(ch(F_T)\exp(-LT))_{|\bbb=0} = 4\cdot \textsf{1}_{X} + L(\bba -\sigma
-4T) + \left[4-\frac{\bba^2}{2} - (\bba -\sigma - 2T)^t \bI T\right] \pt.
\]
In order to carry out the integral on the LHS, we use the
additional observation that the integral $\int_{-d_I}^{1-{d_I}} dT_I
\bba_J = \int_{-d_I}^{1-{d_I}} dT_I \fl*{T_J+e_J} $ equals $\fl*{T_J+e_J}
$ if $I \neq J$, while it equals $e_J-d_J$ when $I=J$. The multiple
integral is therefore
\[
\left[\prod_I \int_{-d_I}^{1-{d_I}} dT_I\right] a_J = \left[\prod_I
\int_{-d_I}^{1-{d_I}} dT_I\right] \fl*{T_J+e_J} = e_J-d_J.
\]
Let $\cube(d) \subset \mathbb{R}^{16}$ denote the cube $\cube(d) =
\prod_{I} [-d_{I},1-d_{I}]$. To simplify notation we will write
$\int_{\cube(d)}$ for the integral $\prod_{I}
\int_{-d_{I}}^{1-d_{I}}dT_{I}$. With this notation the result of the
previous calculation can be written more succinctly as:
\[ \int_{\cube(d)} \bba=e-d.
\] Similarly
\[ \int_{\cube(d)} T =\frac{\sigma}{2} -d.
\]
These formulas allow us to integrate the degree $0$ and degree $2$
pieces of $ch(F_{T})\exp(-LT)$ over the cube $\cube(d)$. Taking
into account the fact that $\bbb = 0$ over the cube $\cube(d)$ we get
\[
\int_{\cube(d)} \left[ch(F_{T})\exp(-LT)\right]_{0} = \left(
\int_{\cube(d)} 4\right) \textsf{1}_{X} = 4\cdot \textsf{1}_{X},
\]
and
\[
\begin{aligned} \int_{\cube(d)} \left[ch(F_{T})\exp(-LT)\right]_{2} &
= \int_{\cube(d)} L(\bba - \sigma - 4T) = L \int_{\cube(d)} (\bba - \sigma -
4T) \\[+0.5pc] & = L \left( e - d -\sigma - 2\sigma + 4d \right)
\\[+0.5pc] & = L (e + 3d - 3\sigma).
\end{aligned}
\]
In particular the degree zero pieces of the LHS and the RHS are
equal, and equating the degree two pieces gives $L (e + 3d - 3\sigma)
= - 2L\sigma$ or equivalently $L (e + 3d - \sigma) = 0$ in
$H^{2}(X,\mathbb{R})$.

In other words $\parch_{0}(F_{\bullet}) = 4$ and
$\parch_{1}(F_{\bullet}) = 0$ if and only if the vectors $e$ and
$d$ satisfy the equation $L (e + 3d - \sigma) = 0$ in
$H^{2}(X,\mathbb{R})$. The latter is a relation among the lines on $X$
and so equivalent to requiring that vectors $e, d \in \mathbb{R}^{16}$
satisfy the matrix equation
\begin{equation} \label{eq-parch1} \boxed{ \bI (e + 3d - \sigma) = 0}
\end{equation}

\

\noindent To compute the degree four piece of the LHS we need to
compute additional integrals whose verification is somewhat tedious
but straightforward.  The first is:
\[
\begin{aligned} \int _{\cube(d)} T^t \bI T & = - \sum_I
\left(d_I^2-d_I+\frac{1}{3}\right) + \sumij
\left(d_I-\frac{1}{2}\right)\left(d_J-\frac{1}{2}\right) \\[+0.5pc] &
= - \sum_I \left(d_I^2-d_I+\frac{1}{3}\right) + \sum_{I}\left( d_I
-\frac{1}{2} \right)^{2} \\[+0.5pc] & \hspace{0.5in} - \sum_{I}\left(
d_I -\frac{1}{2} \right)^{2} + \sumij
\left(d_I-\frac{1}{2}\right)\left(d_J-\frac{1}{2}\right) \\[+0.5pc] &
= \sum_{I}\left(-d_I^2 + d_I - \frac{1}{3} + d_{I}^{2} -d_{I} +
\frac{1}{4}\right) + \left( d - \frac{\sigma}{2}\right)^{t}\bI \left(
d - \frac{\sigma}{2}\right) \\[+0.5pc] & = 16\cdot
\left(-\frac{1}{12}\right) + \left( d - \frac{\sigma}{2}\right)^{t}\bI
\left( d - \frac{\sigma}{2}\right) \\[+0.5pc] & = -\frac{4}{3} +
d^{t}\bI d - \sigma^{t}\bI d + \frac{1}{4} \sigma^{t}\bI\sigma
\\[+0.5pc] & = -\frac{4}{3} + d^{t}\bI d - 4\sigma^{t}d + 16
\\[+0.5pc] & = \frac{44}{3} + d^{t}\bI d - 4\sigma^{t}d.
\end{aligned}
\]
Thus
\begin{equation} \label{eq-TIT} \int_{\cube(d)} \left(2T^{t}\bI
T\right) = \frac{88}{3} - 8\sigma^{t}d + 2d^{t}\bI d.
\end{equation}
Note that when the $d_I$ vanish, this reduces to the
same value $88/3$ which we found on the RHS.  Similarly we verify:
\begin{equation} \label{eq-a2}
  \int_{\cube(d)}\left(-\frac{1}{2}
\bba^2\right) = \frac{1}{2} \fl*{e-d}^2 - \frac{1}{2}\fl*{e-d}^t (2e -2d
-\sigma) - \frac{1}{2}(e-d)^t\sigma.
\end{equation}
Finally we compute the term
\[
\begin{aligned} \int_{\cube(d)} \left( - (\bba -\sigma)^{t}\bI T\right) &
\\[+0.5pc] & \hspace{-0.5in} = \sum_{I} \int_{-d_{I}}^{1-d_{I}}
dT_{I}((\bba_{I} -1)T_{I}) - \sumij \left(\int_{-d_{I}}^{1-d_{I}}
dT_{I}(\bba_{I} - 1)\right) \left(\int_{-d_{J}}^{1-d_{J}} dT_{J}
T_{J}\right) \\[+0.5pc] & \hspace{-0.5in} = \sum_{I}
\int_{-d_{I}}^{1-d_{I}} dT_{I}((\bba_{I} -1)T_{I}) - \sumij \left(e_{I} -
d_{I} - 1\right)\left( \frac{1}{2} - d_{J}\right).
\end{aligned}
\]
The second summand simplifies:
\[
\begin{aligned} \sumij \left(e_{I} - d_{I} - 1\right)\left(
\frac{1}{2} - d_{J}\right) & = \sumij \left( \frac{1}{2} e_{I} -
\frac{1}{2} d_{I} - \frac{1}{2} - e_{I}d_{J} + d_{I}d_{J} + d_{J}
\right) \\[+0.5pc] & = \frac{5}{2} \sum_{I} e_{I} + \frac{5}{2}
\sum_{I} d_{I} - 40 - \sumij (e_{I} - d_{I})d_{J} \\[+0.5pc] & =
\frac{5}{2} \sum_{I} e_{I} + \frac{5}{2} \sum_{I} d_{I} - 40 - \sumij
(e_{I} - d_{I})d_{J}\\[+0.5pc] & \hspace{+0.3in} + \sum_{I}(e_{I} -
d_{I})d_{I} - \sum_{I}(e_{I} - d_{I})d_{I} \\[+0.5pc] & = \frac{5}{2}
e^{t}\sigma + \frac{5}{2} d^{t}\sigma -40 - (e-d)^{t}\bI d - e^{t}d +
d^{2} \\[0.5pc] & = - e^{t}\bI d + d^{t}\bI d -e^{t}d + d^{2} +
\frac{5}{2} e^{t}\sigma + \frac{5}{2} d^{t}\sigma - 40.
\end{aligned}
\]
For the first summand we compute:
\[
\begin{aligned} \sum_{I} \int_{-d_{I}}^{1-d_{I}} dT_{I}((\bba_{I}
-1)T_{I}) & = \frac{1}{2}\sum_{I} \Big[ (\fl*{e_{I} - d_{I}}
-1)\Big((1 - e_{I} + \fl*{e_{I} - d_{I}})^{2} - d_{I}^{2}\Big)
\\[+0.5pc] & \hspace{0.5in} + \fl*{e_{I} - d_{I}}\Big( (1-d_{I})^{2} -
(1 - e_{I} + \fl*{e_{I} - d_{I}})^{2}\Big) \\[+0.5pc] &
\hspace{-0.5in} = \frac{1}{2}\sum_{I} \Big[ - (1 - e_{I} + \fl*{e_{I}
- d_{I}})^{2} + d_{I}^{2} + \fl*{e_{I} - d_{I}}(1 - 2d_{I})\Big]
\\[+0.5pc] & \hspace{-0.5in} = \frac{1}{2}\sum_{I} \Big[
-\fl*{e_{I}-d_{I}}^{2} + \fl*{e_{I} - d_{I}}(2e_{I} - 2d_{I} -1) -
(1-e_{I})^{2} + d_{I}^{2} \Big]\\[+0.5pc] & \hspace{-0.5in} = -
\frac{1}{2}\fl*{e-d}^{2} + \frac{1}{2}\fl*{e-d}^{t}(2e -2d -\sigma) -
\frac{1}{2} (\sigma - e)^{2} + \frac{1}{2} d^{2} \\[+0.5pc] &
\hspace{-0.5in} = - \frac{1}{2}\fl*{e-d}^{2} +
\frac{1}{2}\fl*{e-d}^{t}(2e -2d -\sigma) + e^{t}\sigma -
\frac{1}{2}e^{2} + \frac{1}{2} d^{2} - 8.
\end{aligned}
\]
Together the two summands give
\[
\begin{aligned} \int_{\cube(d)} \left( - (\bba-\sigma)^{t}\bI T\right) &
= - \frac{1}{2}\fl*{e-d}^{2} + \frac{1}{2}\fl*{e-d}^{t}(2e -2d
-\sigma) + e^{t}\sigma - \frac{1}{2}e^{2} + \frac{1}{2} d^{2} - 8
\\[+0.5pc] & \hspace{0.5in} - \left( - e^{t}\bI d + d^{t}\bI d -e^{t}d
+ d^{2} + \frac{5}{2} e^{t}\sigma + \frac{5}{2} d^{t}\sigma - 40
\right)
\end{aligned}
\]
or
\begin{equation} \label{eq-aminussigma}
  \begin{aligned} \int_{\cube(d)} \left( - (\bba-\sigma)^{t}\bI T\right)
& = - \frac{1}{2}\fl*{e-d}^{2} + \frac{1}{2}\fl*{e-d}^{t}(2e -2d
-\sigma) \\[+0.5pc] & \hspace{0.5in} -\frac{1}{2}(e-d)^{2} + e^{t}\bI
d - d^{t}\bI d - \frac{3}{2} e^{t}\sigma - \frac{5}{2} d^{t}\sigma +
32
\end{aligned}
\end{equation}
Now combining the formulas \eqref{eq-TIT},
\eqref{eq-a2}, and \eqref{eq-aminussigma} we get an expression for the
degree four piece of the LHS:
\[
\begin{aligned} \int_{\cube(d)} \left[ ch(F_{T})e^{-LT}\right]_{4} & =
\int_{\cube(d)} \left( 4 - \frac{1}{2}\bba^{2} + 2T^{t}\bI T -
(\bba-\sigma)^{t}\bI T \right) \\[+0.5pc] & \hspace{-0.3in} = 4
\\[+0.5pc] & \hspace{-0.3in}+ \frac{1}{2} \fl*{e-d}^2 -
\frac{1}{2}\fl*{e-d}^t (2e -2d -\sigma) - \frac{1}{2}(e-d)^t\sigma
\\[+0.5pc] & \hspace{-0.3in} + \frac{88}{3} - 8\sigma^{t}d + 2d^{t}\bI
d \\[+0.5pc] & \hspace{-0.3in} - \frac{1}{2}\fl*{e-d}^{2} +
\frac{1}{2}\fl*{e-d}^{t}(2e -2d -\sigma) \\[+0.5pc] & \hspace{-0.3in}
-\frac{1}{2}(e-d)^{2} + e^{t}\bI d - d^{t}\bI d - \frac{3}{2}
e^{t}\sigma - \frac{5}{2} d^{t}\sigma + 32 \\[+0.5pc] &
\hspace{-0.5in} = 36 + \frac{88}{3} -\frac{1}{2}(e-d)^{2} + e^{t}\bI d
+ d^{t}\bI d - (2e + 10d)^{t}\sigma
\end{aligned}
\]
Equating this to the degree four piece of the RHS, i.e. the
constant $88/3$ we conclude that the condition
$\parch_{2}(F_{\bullet}) = 0$ is equivalent to the quadratic
equation
\begin{equation} \label{eq-parch2} \boxed{ 36 -\frac{1}{2}(e-d)^{2} +
e^{t}\bI d + d^{t}\bI d - (2e + 10d)^{t}\sigma = 0 }
\end{equation} This can be simplified somewhat further. Taking into
account the $\parch_{1}(F_{\bullet}) = 0$ condition
\eqref{eq-parch1} we have that
\[ \bI e = - 3\bI d + 4\sigma.
\] Substituting this into \eqref{eq-parch2} we get the quadratic
equation
\[ 36 -\frac{1}{2}(e-d)^{2} - 2d^{t}\bI d - (2e + 6d)^{t}\sigma = 0.
  \] Therefore we have proven the following

  \

  \begin{lem} \label{lem-parch} The conditions
$\parch_{1}(F_{\bullet}) = 0$ and $\parch_{2}(F_{\bullet}) =
0$ are equivalent to the following equations on the vectors $e, d \in
\mathbb{R}^{16}$:
 \[ \boxed{
   \begin{aligned} \bI (e + 3d - \sigma) & = 0, \\[+0.5pc] 36
-\frac{1}{2}(e-d)^{2} - 2d^{t}\bI d - (2e + 6d)^{t}\sigma & = 0.
   \end{aligned} }
 \]
  \end{lem} %{\bfseries Proof.} See the calculations above.\ \hfill  $\Box$

      \

\subsection{Killing the Chern classes} \label{ssec:kill.parch}
To solve the conditions in
Lemma~\ref{lem-parch}, we start by introducing the vector:
\index{terms}{parabolic!Chern classes}
\[ v:= e+3d-\sigma.
\] In terms of $v$ the linear equation from Lemma~\ref{lem-parch} now
says that $\bI v=0$, while the quadratic equation from
Lemma~\ref{lem-parch} becomes:
\begin{equation} \label{eq-parch2viav} -\frac{1}{2}v^2 - 8d^2 +4vd
-2d^t \bI d +4 d^t \sigma -4 =0.
\end{equation} Recall that the intersection matrix $\bI$ has three
eigenspaces, of dimension $1,5,10$, corresponding to eigenvalues
$4,-4,0$ respectively. We decompose the unknown vector $d$ into its
components in the eigenspaces of $\bI$:
\[
\begin{aligned} d & = d_{\{1\}} + d_{\{5\}} + d_{\{10\}} \\ \bI
d_{\{1\}} & = 4 d_{\{1\}}, \ d_{\{1\}} = x\frac{\sigma}{8}, \ x^2 = 4
{d_{\{1\}}}^2 \\ \bI d_{\{5\}} & = -4 d_{\{5\}} \\ \bI d_{\{10\}} & =
0.
\end{aligned}
\]
Equation \eqref{eq-parch2viav} simplifies to:
\begin{equation} \label{ch2''} -\frac{1}{2}(v-4d_{\{10\}})^2 -
4(x-1)^2 =0.
\end{equation} All solutions of the Chern class constraints are therefore
given by:
\begin{equation} \label{chsoln}
\begin{aligned} d_{\{5\}}& {\text{ and }} d_{\{10\}} {\text{ arbitrary
}} \\ x &=1, {\text{ so }} \ d_{\{1\}} = \frac{\sigma}{8} \\ v
&=4d_{\{10\}}, {\text{ so: }} \\ e &=4d_{\{10\}} - 3d + \sigma \\
&=d_{\{10\}} -3d_{\{5\}} +\frac{5\sigma}{8}
\end{aligned}
\end{equation}

\subsection{Hecke conditions}

\noindent
We start with two pairs of vectors $e, d$ and $e',d'$ each obeying
\eqref{chsoln}. In section~\ref{ssec:kill.parch} we proved that for
every $\mathfrak{a} \in \Pic^{0}(\sC)$ the modular spectral data
$\left(Y,\mycal{L}_{\mathfrak{a}}(eE +
d(E+G))_{\bullet,\bullet}\right)$ and
$\left(Y,\mycal{L}_{\mathfrak{a}}(e'E +
d'(E+G))_{\bullet,\bullet}\right)$ define stable rank four parabolic
Higgs bundles $(F_{\bullet},\varphi)$ and $(F_{\bullet}',\varphi')$ on
$(X,\Par_{X})$ with vanishing first and second parabolic Chern class.
Here, as in section~\ref{ssec:hecke.data}, we write
$\mycal{L}_{\mathfrak{a}}$ for the line bundle on $Y$ which is the
pullback of the Fourier-Mukai transform of the skyscraper sheaf
$\mathcal{O}_{\mathfrak{a}}$ on $\Pic^{0}(\sC)$.

Fix $a, b \in \mathbb{R}^{5}$ such that $\sum_{i =1}^{5} (a_{i}+b_{i})
= 3$.  In this and the next section we derive and solve the additional
constraints on $e$, $d$, $e'$, $d'$ and the coefficients $\bzeta$ (see
\eqref{eq:zeta.basic}) of the Hecke kernel so that for every
$\mathfrak{a} \in \Pic^{0}(\sC)$, the parabolic Higgs bundles
$(F_{\bullet},\varphi)$, $(F_{\bullet}',\varphi')$ on $(X,\Par_{X})$
satisfy the parabolic Hecke eigensheaf condition \eqref{eq-par.hecke}
with eigenvalue $(E_{\bullet},\theta)$ defined by the spectral data
$(\sC,\mathfrak{a}(a\widetilde{P}
+b\widetilde{P})_{\bullet,\bullet})$.

Interestingly, as we will see in
Claim~\ref{claim:spectral.is.eigensheaf}, the correspondence between
eigenvalues and eigensheaves works in both directions.  In other
words, given vectors $e, d$ that obey \eqref{chsoln}, we
can also solve for and find vectors $a$ and $b$ for the eigenvalue, 
$\bzeta$ for the Hecke kernel, and vectors $e', d'$ that  obey \eqref{chsoln},
so that $(F_{\bullet},\varphi)$ and
$(F_{\bullet}',\varphi')$ given by
$\left(Y,\mycal{L}_{\mathfrak{a}}(eE +
d(E+G))_{\bullet,\bullet}\right)$ and 
$\left(Y,\mycal{L}_{\mathfrak{a}}(e'E +
d'(E+G))_{\bullet,\bullet}\right)$ satisfy the parabolic Hecke
eigensheaf condition \eqref{eq-par.hecke} with Hecke eigenvalue
$(E_{\bullet},\theta)$ correspinding to the spectral data
$(\sC,\mathfrak{a}(a\widetilde{P} +
b\widetilde{P}))_{\bullet,\bullet}$.

Let us first show that given an eigenvalue parabolic Higgs bundle on
$\sC$ we can solve for $e, d$ and $e', d'$ so that the Chern class
conditions, the Hecke condition, and the Hecke kernel condition are
all satisfied.  We are given $a_i, b_i$, $i =1, \ldots, 5$.  The Hecke
condition is a condition on the coefficients $e_I, d_I, lp_I, lq_I,
llc_I, r_{1,i},r_{2,i}$.  As in the previous section, our convention
is that these are column vectors, where $I$ denotes a $16$-dimensional
index while $i$ is $5$ dimensional. The corresponding vectors of
divisors $\Theta_I$, $Lp_I$, $Lq_I$, $LLC_I$, $R_{1,i}$, $R_{2,i}$ are
row vectors. We showed in Proposition~\ref{pro:abelianize.lhs} and
Proposition~\ref{prop:Pic.Hecke} that the Hecke condition is
equivalent to requiring that $lp + d \in \mathbb{Z}^{16}$ together
with the vanishing in $\Pic(Y\times \sC)$ of:
  \begin{equation}  \label{Heck}
    \boxed{
{
\begin{aligned} 
  & \sum_{I\in \evenL}  \left(\fl*{lp_{I} +d_{I}} 
  -\frac{1}{8}\right)\hatAJ^{*}\Theta_{I} \\ 
  & \quad + \sum_{I \in \evenL}
  \left(\fl*{T_{I} +llc_{I} + d_{I}} + \fl*{T_{I} +lq_{I}} +
  \frac{5}{8} - \fl*{T_{I} +e_{I}'}
  - \fl*{T_{I} +d_{I}'}\right)E_{I}\times \sC \\
& \quad + 
\sum_{I\in \evenL} \left(\fl*{T_{I} +lq_{I}} + \frac{1}{8} -
\fl*{T_{I} +d_{I}'}\right)G_{I}\times \sC \\
& \quad + 
\sum_{i=1}^{5} \left(\fl*{t_{i} +r_{1,{i}}} + \fl*{t_{i} +r_{2,{i}}}
+ \frac{3}{5} - \fl*{t_{i} +a_{i}}
- \fl*{t_{i} +b_{i}}\right)Y\times \tilde{p}_{i}.
\end{aligned}
}}
  \end{equation}

A necessary condition for \eqref{Heck} to hold for all $T,t$ is that
it must remain unchanged as either $T$ or $t$ goes through a value at
which one of the floor functions above is discontinuous: the positive
jump must be offset by a negative one. This implies that (up to
interchanging the $a$'s with the $b$'s if necessary) if we define new
column vectors:
\[
\begin{aligned} & A_I := lp_I + d_I\\ & B_I := lq_I - d_I'\\ & C_I :=
llc_I +d_I -e_I'\\ & n_{1,i} := r_{1,i} - a_i \\ & n_{2,i} := r_{2,i} -
b_i
\end{aligned}
\]
we must have that all of $B_I, C_I, n_{1,i}, n_{2,i}$ are integers.
Since from Proposition~\ref{pro:abelianize.lhs} we already had the
condition that $A_{I} = lp_{I} + d_{I}$ are integers we conclude that
for Hecke condition to hold it is necessary for $A_{I}, B_I, C_I,
n_{1,i}, n_{2,i}$ to be all integers.

Once these integralities are imposed, the behavior of the Hecke
condition as $T,t$ vary is under control. We still need to make sure
the equations hold for some specific values of $T$ and $t$. For this
it is convenient to take $T_I := -d_I'$ and $t_{i} = - a_{i}$ so
\eqref{Heck} becomes:
\[
\begin{aligned} &\sum_{I \in \evenL} \hatAJ^{*}\Theta_I \left(\fl*{lp_I
    +d_I}-\frac{1}{8}\right) + \\ &\sum_{I \in \evenL}
  E_I\times \sC \left(\fl*{-d_{I}' + llc_I + d_{I}} + \fl*{-d_{I}' + lq_{I}}
  + \frac{5}{8} - \fl*{-d_I' + e_{I}'} \right) + \\
  &\sum_{I \in \evenL} G_I\times \sC \left(\fl*{-d_I'+lq_I} 
  + \frac{1}{8} \right) + \\ &
  \sum_{i=1}^{5} Y \times \widetilde{p}_i \left(\fl*{-a_{i}
+r_{1,i}} + \fl*{-a_{i} +r_{2,i}} + \frac{3}{5} - 
\fl*{-a_{i} +b_{i}}\right)
\end{aligned}
\]
or equivalently:
\[
\begin{aligned}
  &\sum_{I \in \evenL} \hatAJ^{*}\Theta_I
  \left(A_I-\frac{1}{8}\right) + \\
  &\sum_{I \in \evenL} E_I\times \sC \left(C_{I}  + B_I + \frac{5}{8}\right) + \\
  &\sum_{I \in \evenL} G_I\times \sC \left(B_I + \frac{1}{8}\right) + \\ &
\sum_{i=1}^{5} Y \times \widetilde{p}_i \left(\fl*{-a_{i} +r_{1,i}} + \fl*{-a_{i}
+r_{2,i}} + \frac{3}{5} -  \fl*{-a_{i}  +b_{i}}\right).
\end{aligned}
\]
This needs to vanish in $\Pic(Y\times \sC) \otimes \bR$. It is
easier to solve first a slightly weaker condition: the vanishing of
the same expression in cohomology $H^2(Y\times \sC, \bR)$. A full set of
cohomological relations among the natural divisors in $Y \times
\sC$ is:
\begin{itemize}
\item[(i)] $G_I\times \sC + E_I\times \sC +
  \sum_{\substack{\{ J \, | \,  \bI_{IJ}=1\}}}
  E_J\times \sC$ \ is independent of $I\in \evenL$;
\item[(ii)] $\hatAJ^{*}\Theta_I + E_I \times
  \sC$ \ is independent of $I \in \evenL$; 
\item[(iii)]  $Y \times
\widetilde{p}_i$  is independent of $i \in \{1, 2, 3, 4, 5\}$.
\end{itemize}
We use these relations to simplify \eqref{Heck}. Since $Y \times
\widetilde{p}_i$ occurs only in the  relation (iii), its coefficent in 
\eqref{Heck} must vanish modulo the  relation (iii). Explicitly this means
that we must have
\[
\begin{aligned}
0 & = \sum_{i = 1}^{5} \left( \fl*{-a_{i} +r_{1,i}} + \fl*{-a_{i} 
  +r_{2,i}} + \frac{3}{5} - \fl*{-a_{i} +b_{i}}\right) \\[+0,5pc]
& = \sum_{i = 1}^{5} \left( \fl*{n_{1,i}} + \fl*{-a_{i} + b_{i} 
  + n_{2,i}} + \frac{3}{5} - \fl*{-a_{i} +b_{i}}\right) \\[+0,5pc]
& = \sum_{i = 1}^{5} \left( n_{1,i} + n_{2, i} +  \fl*{-a_{i} + b_{i}}
  - \fl*{-a_{i} +b_{i}}\right)  + 3 \\[+0,5pc]
& = \sum_{i = 1}^{5} \left( n_{1,i} + n_{2, i} \right)  + 3.
\end{aligned}
\]
Thus the last line in \eqref{Heck} vanishes in $H^{2}(Y\times
\sC,\mathbb{R})$ if and only if
\[
\sum_{i = 1}^{5} (n_{1,i} + n_{2,i}) = - 3.
\]
Furthermore since the relations (iii) are independent of relations (i)
and (ii)  we must have that the sum of the remaining  terms in \eqref{Heck}:
\[
\begin{aligned}
  &\sum_{I \in \evenL} \hatAJ^{*} \Theta_I
  \left(A_I-\frac{1}{8}\right)  + \\
  &\sum_{I \in \evenL} E_I\times \sC \left(C_{I} + B_I + \frac{5}{8}
  \right)  + \\
  &\sum_{I \in \evenL} G_I\times \sC \left(B_I + \frac{1}{8}\right)
\end{aligned}
\]
must vanish modulo relations (i) and (ii). Relations (i) and (ii)
in turn allow
us to eliminate $\hatAJ^{*}\Theta_I$ and $ G_I\times \sC $, leading
to the vanishing of the new coefficients of $E_I\times \sC$:
\[
\begin{aligned} 0 & = C_{I} + B_I  + \frac{5}{8} - A_I
  +\frac{1}{8} - B_I -\frac{1}{8} - \sum_{\substack{\{ J \, | \,  \bI_{I,J}=1\}}}
    \displaylimits \left(B_J
    +\frac{1}{8}\right) \\[+0.5pc]
    & = C_I -A_I -\sum_{\substack{\{ J \, | \,  \bI_{I,J}=1\}}}
    \displaylimits  B_J.
\end{aligned}
\]
More concisely, the Hecke condition is that:
\begin{equation} \label{Hecksoln}
  \begin{aligned}
& {\text {
all of } } A_I, B_I, C_I, n_{1,i}, n_{2,i} {\text { are integers}} \\
    &C_I = A_I + \sum_{\substack{\{ J \, | \,  \bI_{I,J}=1\}}}
    \displaylimits
    B_J \\[+0.5pc]
& \sum_{i = 1}^{5} (n_{1,i} + n_{2,i}) = - 3.
\end{aligned}
\end{equation}
To summarize, solutions to the equations considered so far, i.e. the
non-abelian Hodge theory constraints \eqref{ssec:kill.parch} and the
parabolic Hecke eigensheaf condition \eqref{eq-par.hecke}, are
determined uniquely in terms of the given vectors $a,b \in
\mathbb{R}^{5}$ and the vectors of unknowns $A,B,d,d',n_1,n_2$. In
terms of these vectors the coefficients for the Hecke kernel and the
remaining coefficients for the parabolic structure of the eigensheaf
are given by
\begin{equation} \label{eq:kernel.coefficients}
\left| \quad
  \begin{aligned} & lp := A-d\\ & lq:=B+d'\\ & llc := C+e'-d \\ &
  r_1 := a +n_1\\ & r_2 := b +n_2
  \end{aligned}
\right.
\end{equation}
and the non-abelian Hodge theory constraints \eqref{ssec:kill.parch} are
given by
\begin{equation} \label{eq:nah.constraints}
  \left| \quad
  \begin{aligned}
    & e, d, e', d' \in \mathbb{R}^{16} \\
    & d_{\{1\}} = d_{\{1\}}' = \frac{\sigma}{8} \\
     & e =d_{\{10\}} -3 d_{\{5\}} +\frac{5\sigma}{8} \\
  & e' =d_{\{10\}}' -3 d_{\{5\}}' +\frac{5\sigma}{8} 
  \end{aligned}
\right.
\end{equation}
whle the parabolic Hecke eigensheaf constraints \eqref{eq-par.hecke} are
equivalent to
\begin{equation} \label{eq:ABCn12}
  \left| \quad
  \begin{aligned}
& A,B,C \in \mathbb{Z}^{16}, \;  n_1,n_2 \in \mathbb{Z}^{5} \\
& C = A +
    (\bI +\one)B \\
    & \Sigma_{i=1}^5
(n_{1,i} + n_{2,i}) = -3
  \end{aligned}
\right.
\end{equation}

\subsection{The class of the kernel} \label{ssec:class.of.the.kernel}

\noindent
The remaining condition is that the first parabolic Chern class 
of the kernel $\mycal{I}_{\bullet} =
\mathcal{O}_{H}\left(\bzeta\Par_{H}\right)_{\bullet}$ vanishes in
$\Pic(H)\otimes \mathbb{R} = H^{2}(H,\mathbb{R})$:
\begin{equation} \label{Phi}
  Lp \cdot lp + Lq \cdot lq + LLC \cdot llc
+ R_1 \cdot r_1 + R_2 \cdot r_2 = 0 \in H^{2}(H,\mathbb{R}).
\end{equation}
Since the rank of $\Pic(H)$ equals $2\times 6 +1 +16 = 29$, this
condition is equivalent to 29 numerical equations on the
coefficients. Such conditions can be written explicitly by
intersecting the divisor in the left hand side of \eqref{Phi} with any
curve in $H$.

\

\begin{ex} \label{ex:hecke.curve.condition}
For instance we can intersect the divisor $\bzeta\Par_{H}$ with a
generic Hecke curve i.e. the $q$-preimage $H_{x,c} = q^{-1}(x,c)$ of a
general point $(x,c) \in X\times C$. Since the divisors $Lq$, $LLC$,
$R_{1}$, and $R_{2}$ are all vertical for $q$, they all have
intersection number $0$ with $H_{c,x}$ and so we have
\[
(\bzeta\Par_{H})\cdot H_{x,c} = \sum_{I \in \evenL} lp_{I}(Lp_{I}\cdot
H_{x,c}).
\]
But in the proof of Lemma~\ref{lem:divisor.of.beta} we argued that
$Lp_{I} \cdot H_{x,c} = 1$ for all $I$.  So this bit of the condition for the first parabolic Chern class of the kernel to vanish becomes: $\sum_{I
  \in \evenL} lp_{I} = 0$. Using the symmetry of the two copies of $X$
in the blowup description of $H$ this implies that also $\sum_{I \in
  \evenL} lq_{I} = 0$.  In other words, intersecting the homology
class $[\bzeta\Par_{H}]$ with general Hecke curves
imposes  the natural necessary conditions
\begin{equation} \label{eq:{1}component}
lp^{t}\cdot \sigma = lq^{t}\cdot \sigma = 0.
\end{equation}
\end{ex}

\

\noindent
To convert \eqref{Phi} into a complete set of equations on $\bzeta$,
we need to intersect $[\bzeta\Par_{H}]$ with a collection of curves
that span $H_2(H,\bZ)$. We use the following redundant\footnote{We
  will see below that the conditions imposed by the last bunch, the
  $\Sigma_I$, are linearly dependent on the others.} collection of 49
curves:
\[
L_I \times \pt, \ \pt \times L_I, \ \widehat{\Delta}_{\pt}, \ \Sigma_I.
\]
We tabulate the
intersection numbers of these curves with the various divisors
below. Here $\incdelta_{I,i} $ equals 1 if $ i \in I$ and 0 if $i \not\in
I$.

\

\medskip

\begin{center}
\begin{tabular}{| l || l | l | l | l | l | } \hline $ { } $ &
$Lp_J$ & $Lq_J$ & $LLC_J$ & $R_{1,i}$ & $R_{2,i} $\\ \hline \hline $
  L_I \times \pt$ & $\bI_{I,J}$ &0 & 0 & $ 1 -\incdelta_{I,i}$ &
  $\incdelta_{I,i} $
  \\ \hline $ \pt \times L_I$ & 0 & $\bI_{I,J}$ & 0 &
  $1 -\incdelta_{I,i}$ &
  $\incdelta_{I,i}$ \\ \hline $ \widehat{\Delta}_{\pt}$ & 0 & 0 & 0 &
  1 & 1
\\ \hline $ \Sigma_I$ & $\one_{I,J}$ & $\one_{I,J}$ & $-\one_{I,J}$ &
$ \incdelta_{I,i}$ & $1 -\incdelta_{I,i} $ \\ \hline
\end{tabular}
\end{center}

\

\medskip

The 49 conditions we get are therefore:
\begin{align}
 \label{eq:Lpt} 0 & = \ \sum_{J \in \evenL}
  \bI_{I,J} lp_J + \sum_{i = 1}^{5}\left((1-\incdelta_{I,i}) r_{1,i} +
  \incdelta_{I,i} r_{2,i}\right) \\ 
\label{eq:ptL} 0 & = \ \sum_{J\in \evenL}
\bI_{I,J} (lp_J - lq_J) \\
\label{eq:Deltapt} 0 & = \ \sum_{i=1}^{5} ( r_{1,i} +
  r_{2,i} ) \\
\label{eq:SigmaI} 0 & = \ (lp_I + lq_I - llc_I) + \sum_{i=1}^{5}
  ((1-\incdelta_{I,i}) r_{2,i} + \incdelta_{I,i} r_{1,i}).
\end{align}
This system of linear equations is redundant. In fact the equations
\eqref{eq:SigmaI} are linear combinations of the equations
\eqref{eq:Lpt}, \eqref{eq:ptL}, \eqref{eq:Deltapt}. For ease of
reference we record this in the following lemma.

\

\begin{lem} \label{lem:redundant}
  Fix the vectors $a$ and $b$ in $\mathbb{R}^{5}$.  
Assume the vectors $e$, $d$, $e'$, $d'$, $A$, $B$, $C$, $n_{1}$, and
$n_{2}$ satisfy the non-abelian Hodge theory constraints
\eqref{eq:nah.constraints} and the Hecke eigensheaf constraints
\eqref{eq:ABCn12} we derived in the previsous two sections. Then the
Hecke kernel $\mathcal{O}_{H}(\bzeta\Par_{H})_{\bullet}$ satisfies the
non-abelian Hodge theory condition
$\parch_{1}(\mathcal{O}_{H}(\bzeta\Par_{H})_{\bullet}) =
[\bzeta\Par_{H}] = 0$ in $H^{2}(H,\mathbb{R})$ if and only if the
following additional set of equations on $\bzeta =
(lp,lq,llc,r_{1},r_{2}) = (A-d,B+d',C+e'-d,a+n_{1},b+n_{2})$ holds:
  \begin{equation} \label{eq:33cond}
\left| \quad
\begin{aligned}
0 & = \bI \, lp + \frac{1}{2}(\unone - 2\incdelta)(r_{1} - r_{2}), \\[+0.5pc]
0 & = \bI \, lq + \frac{1}{2}(\unone - 2\incdelta)(r_{1} - r_{2}), \\[+0.5pc]
0 & =  \sum_{i=1}^{5} ( r_{1,i} + r_{2,i}).
\end{aligned}
\right.
\end{equation}
Here $\unone$ is the $16\times 5$ matrix all of whose entries are
equal to $1$, and $\incdelta$ is the $16\times 5$ matrix whose entries
are labeled by $(I,i)$ with $I \in \evenL$, $i = 1, \ldots 5$, and
$\incdelta_{I,i} = 1$ if $i \in I$ and $\incdelta_{I,i} = 0$ if $i
\not\in I$.
\end{lem}
{\bfseries Proof.} Rewriting \eqref{eq:Lpt}, \eqref{eq:ptL},
\eqref{eq:Deltapt}, and \eqref{eq:SigmaI} in matrix form we get a
complete set of equations on $\bzeta$:
\begin{equation} \label{eq:49cond}
\begin{aligned}
  0 & = 
  \bI \, lp + (\unone-\incdelta) r_{1} + 
  \incdelta \, r_{2}, \\
  0 & =  \bI \, (lp - lq), \\
  0 & = \sum_{i=1}^{5} ( r_{1,i} +
  r_{2,i} ), \\
  0 & = lp + lq - llc + (\unone-\incdelta) r_{2} + \incdelta \, r_{1}.
\end{aligned}
\end{equation}
First note that
\[
(\unone-\incdelta) r_{1} + \incdelta \,  r_{2} = \frac{1}{2}(\unone -
2\incdelta)(r_{1} - r_{2}) + \frac{1}{2}\unone(r_{1} + r_{2}).
\]
But $\unone(r_{1} + r_{2}) = (\sum_{i = 1}^{5}(r_{1,i}+r_{2,i}))\cdot
\sigma$.  Therefore the equations  \eqref{eq:49cond} are equivalent to 
\begin{equation} \label{eq:49cond.simpler}
\begin{aligned}
  0 & = 
  \bI \, lp + \frac{1}{2}(\unone - 2\incdelta)(r_{1} - r_{2}), \\ 
  0 & =   \bI \, lq + \frac{1}{2}(\unone - 2\incdelta)(r_{1} - r_{2}), \\ 
  0 & = \sum_{i=1}^{5} ( r_{1,i} +
  r_{2,i} ), \\
  0 & =  lp + lq - llc - \frac{1}{2}(\unone - 2\incdelta)(r_{1} - r_{2}).
\end{aligned}
\end{equation}
So to prove the lemma we need to show that the last line in
\eqref{eq:49cond.simpler} follows from the first three lines and the
non-abelian Hodge and Hecke eigensheaf conditions we derived in the
previous two sections. Indeed from \eqref{eq:kernel.coefficients},
\eqref{eq:ABCn12}, and \eqref{eq:nah.constraints}
we have
\[
\begin{aligned}
  lp + lq - llc & = (A - d) + (B + d') - (C + e' - d),  \\
  & = A + B - (A + (\bI + \one)B) - e' + d', \\
  & = - \bI\, B +  4 d_{\{5\}}' - \frac{4}{8} \sigma.
\end{aligned}
\]
Taking into account that 
$d_{\{1\}}' = \sigma/8$   is an eigenvector of $\bI$ of
eigenvalue $4$, $d_{\{5\}}'$ is an eigenvector of $\bI$ of
eigenvalue $-4$, and $d_{\{10\}}' \in \ker \bI$  we conclude that
\[
lp + lq - llc = - \bI (B+d'). 
\]
But by \eqref{eq:kernel.coefficients} we have $B + d' = lq$ and from
the second equation in \eqref{eq:49cond.simpler} we have \linebreak
$- \bI\, lq =
\frac{1}{2}(\unone - 2\incdelta)(r_{1} - r_{2})$. This concludes the
proof of the lemma. \ \hfill $\Box$

\

\noindent
To find the solutions of the conditions \eqref{eq:33cond} let us note
that from the definition of $\incdelta$ and $\unone$ we have
\[
\bI\, \incdelta = 4(\unone - \incdelta) \quad \text{and} \quad
\bI\, \unone  = 4 \unone
\]
and so
\[
\bI (\unone - 2\incdelta) = -4 (\unone - 2\incdelta).
\]
In particular, the image of the matrix $\frac{1}{2}(\unone -
2\incdelta)$ is contained in the $5$-dimensional eigenspace of $\bI$
corresponding to the eigenvalue $-4$. Since $\bI\, lp = 4lp_{\{1\}} -
4lp_{\{5\}}$ and $\bI\, lq = 4lq_{\{1\}} - 4lq_{\{5\}}$ we conclude
that the first two conditions in \eqref{eq:33cond} are equivalent to
\begin{equation} \label{eq:33cond.simpler}
\begin{aligned}
  0 = & lp_{\{1\}} = lq_{\{1\}}, \\
  0 = & -4lp_{\{5\}} + \frac{1}{2}(\unone -
  2\incdelta)(r_{1} - r_{2}), \\
  0 = & -4lq_{\{5\}} + \frac{1}{2}(\unone -
  2\incdelta)(r_{1} - r_{2}), \\
\end{aligned}
\end{equation}
Note also that
\[
lp_{\{1\}} = \frac{lp^{t}\cdot \sigma}{16} \sigma, \quad \text{and} \quad
lq_{\{1\}} = \frac{lq^{t}\cdot \sigma}{16} \sigma
\]
and so the conditions $ lp_{\{1\}} = lq_{\{1\}} = 0$ recover the
conditions $lp^{t}\cdot \sigma = lq^{t}\cdot \sigma = 0$ that we
derived independently in Example~\ref{ex:hecke.curve.condition}.

Using the simplification  \eqref{eq:33cond.simpler} we can now rewrite
\eqref{eq:33cond} as equations on the vectors
$A,B,d,d',n_1,n_2$:
\[
\begin{aligned}
  0 & = \ A_{\{1\}} - d_{\{1\}} = B_{\{1\}} + d_{\{1\}}' \\
  0 & =  \ -4(A_{\{5\}} -d_{\{5\}}) +
  \frac{1}{2}(\unone-2\incdelta)((a+ n_{1}) - (b+n_{2})), \\
  0  & = \ -4(B_{\{5\}} + d_{\{5\}}') +
  \frac{1}{2}(\unone-2\incdelta)((a+ n_{1}) - (b+n_{2})), \\
  0 & = \ \sum_i ( a_i + n_{1,i} + b_i + n_{2,i} ).
\end{aligned}
\]
Therefore if we write
$\mathsf{x} := \frac{1}{8}(\unone - 2\incdelta)(r_{1} - r_{2}) =
\frac{1}{8}(\unone - 2\incdelta)((a+ n_{1}) - (b+n_{2}))$
the  general solution of  \eqref{eq:33cond}  consists of vectors
$A,B \in \mathbb{Z}^{16}$,
$n_1,n_2 \in \mathbb{Z}^{5}$, $d, d' \in \mathbb{R}^{16}$ satisfying:
\begin{equation} \label{eq:fixedab.cond}
\left| \quad 
\begin{aligned} & d_{\{1\}} = A_{\{1\}}=  \frac{\sigma}{8}, \\
 &  d_{\{1\}}' = - B_{\{1\}} =  
  \frac{\sigma}{8} \\
  &  d_{\{5\}}  = A_{\{5\}} - \mathsf{x} \\
  & d_{\{5\}}' = - B_{\{5\}}  + \mathsf{x} \\
  & \Sigma_{i=1}^5
(n_{1,i} + n_{2,i}) = -3 \\ & A_{\{10\}}, B_{\{10\}}, d_{\{10\}},
d_{\{10\}}' {\mbox{ are unconstrained.}}
\end{aligned}
\right.
\end{equation}

\

\begin{rem} \label{eq:final.cond}
At this stage it is important to emphasize a subtlety in the way we
have parametrized our unknowns. Note that in deriving the combined
constraints \eqref{eq:fixedab.cond} we used fixed vectors $a, b \in
\mathbb{R}^{5}$ to describe the eigenvalue parabolic Higgs bundle
$(E_{\bullet},\theta)$ on $(C,\Par_{C})$. However, as explained in
section \ref{sssec:eigenvalue} this description is
overparametrized. That is, different pairs $(a,b)$ can give rise to
the same $(E_{\bullet},\theta)$. Indeed, $(E_{\bullet},\theta)$
corresponds to a parabolic spectral line bundle $\mathfrak{a}(a\sP +
b\sP)_{\bullet,\bullet}$ on $(\sC,\sP+\sP)$, and two triples
$(\mathfrak{a},a,b)$ and $(\mathfrak{a}',a',b')$ produce the same
parabolic Higgs bundle $(E_{\bullet},\theta)$ if and only if
\[
(\mathfrak{a}',a',b') = \left( \mathfrak{a}\otimes\mathcal{O}_{\sC}(-(m_{1} +
m_{2})\sP),\mathfrak{s}_{1}(a,b) + m_{1},\mathfrak{s}_{2}(a,b) + m_{2}\right)
\]
where $m_{1}, m_{2} \in \mathbb{Z}^{5}$, satisfying
$\deg((m_{1}+m_{2})\sP) = 0$ and $\mathfrak{s} \in \mathfrak{S} \cong
(\mathbb{Z}/2)^{5}$, where $\mathfrak{S}$ is the group generated by
all possible swaps $a_{i} \leftrightarrow b_{i}$.

In other words, the conditions \eqref{eq:fixedab.cond} are slightly
more restrictive than what is dictated by the parabolic Hecke
eigensheaf problem since we have suppressed the extra freedom
contained in the vectors $m_{1}, m_{2} \in \mathbb{Z}^{5}$ and in the
group elements $\mathfrak{s} \in \mathfrak{S}$. In fact the freedom 
created by introducing the integral vectors $m_{1}$ and $m_{2}$ does
not affect the shape of the conditions \eqref{eq:fixedab.cond} since
we can just translate $n_{1}$ and $n_{2}$ by $m_{1}$ and $m_{2}$
respectively. On the other hand, the action of $\mathfrak{S}$ does
provide additional freedom in formulating the combined constraints.

This freedom can be summarized as follows. Suppose
$(E_{\bullet},\theta)$ is a strongly parabolic rank two Higgs bundle
with parabolic structure determined by some pair of vectors $a,b \in
\mathbb{R}^{5}$.  Assume the vectors $e$, $d$, $e'$, $d'$, $A$, $B$,
$C$, $n_{1}$, and $n_{2}$ satisfy the non-abelian Hodge theory
constraints \eqref{eq:nah.constraints} and the Hecke eigensheaf
constraints \eqref{eq:ABCn12} with eigenvalue
$(E_{\bullet},\theta)$. Then the constraints
\eqref{eq:nah.constraints} and \eqref{eq:ABCn12} are
$\mathfrak{S}$-invariant and for every $\mathfrak{s} \in \mathfrak{S}$
the vector
\[
\mathfrak{s}(\bzeta) :=
\left(lp,lq,llc,\mathfrak{s}_{1}(r_{1},r_{2}),
\mathfrak{s}_{2}(r_{1},r_{2})\right)
\]
defines a Hecke kernel
$\mathcal{O}_{H}(\mathfrak{s}(\bzeta)\Par_{H})_{\bullet}$ which
satisfies the non-abelian Hodge theory condition
$\parch_{1}(\mathcal{O}_{H}(\mathfrak{s}(\bzeta)\Par_{H})_{\bullet}) =
[\mathfrak{s}(\bzeta)\Par_{H}] = 0$ in $H^{2}(H,\mathbb{R})$ if and
only if $(lp,lq,llc,r_{1},r_{2}) = \linebreak
(A-d,B+d',C+e'-d,a+n_{1},b+n_{2})$
obey the equations \eqref{eq:fixedab.cond} but with \linebreak
$\mathsf{x} = \frac{1}{8}(\unone -
2\incdelta)\left(\mathfrak{s}_{1}(r_{1},r_{2}) -
\mathfrak{s}_{2}(r_{1},r_{2})\right)$
\end{rem}

\index{notations}{Sfrak@$\mathfrak{S}$}
\index{notations}{sfrak@$\mathfrak{s}$}
\index{notations}{sfraki@$\mathfrak{s}_{i}$}

\

\bigskip

\noindent
It is now straightforward to check that these equations have a unique
solution modulo the obvious overparametrizations, e.g. modulo the
appropriate integral shifts in $a$, $b$, $d$, and $d'$ and modulo the
unconstrained components $d_{\{10\}}, d_{\{10\}}'$ and $A_{\{10\}},
B_{\{10\}}$, which vanish when mapped to the parabolic Picard
varieties of $X$ and $H$.

\

To explain this properly consider the free
abelian groups
\[
\begin{aligned}
  \lat_{C} & = \op{Fun}(\Par_{C},\mathbb{Z}) = \bigoplus_{i = 1}^{5}
  \mathbb{Z}\delta_{p_{i}} \\ \lat_{X} & =
  \op{Fun}(\Par_{X},\mathbb{Z}) = \bigoplus_{I \in \evenL}
  \mathbb{Z}\delta_{L_{I}}
\end{aligned}
\]
with their natural bases of delta functions. The matrices $\bI$,
$\one$, $\unone$, and $\incdelta$ can then be viewed as group homomorphisms
\[
\bI, \; \one : \lat_{X} \to \lat_{X}, \qquad \unone, \; \incdelta :
\lat_{C} \to \lat_{X}.
\]
Consider the vector $\sigma \in \lat_{X}$ with all coordinates in the
standard basis equal to $1$ and let
\[
\sigma^{\perp} = \left\{ v \in \lat_{X} \; \left| \; \sigma^{t}v = 0
\right.\right\}.
\]
Let $\lat_{X}^{0} \subset (\lat_{X}\otimes\mathbb{R})_{\{5\}} \subset
\lat_{X}\otimes\mathbb{R}$ be the projection of the free abelian group
$\sigma^{\perp}$ to the $5$-dimensional eigenspace of $\mathbb{I}$
corresponding to eigenvalue $-4$. With this notation we now have the
folllowing

\

\begin{lem} \label{lem:cosets} 
  \begin{itemize}
  \item[(a)]  $\lat_{X}^{0}$ is a free abelian group of rank $5$.
  \item[(b)] The linear operator
    \[
    \frac{1}{4}(\unone - 2\incdelta) :
    \lat_{C}\otimes \mathbb{R} \to (\lat_{X}\otimes\mathbb{R})_{\{5\}} \subset
\lat_{X}\otimes\mathbb{R}
\]
sends the rank $5$ lattice $\lat_{C}$ to a finite index sublattice of
$\lat_{X}^{0}$.
\item[(c)] Fix $n_{1}, n_{2} \in \lat_{C}$ vectors satisfying $\sum_{i
  = 1}^{5} (n_{1,i} + n_{2,i}) = -3$, and let $r_{1} = a + n_{1}$ and
  $r_{2} = b + n_{2}$.  Let $\acan \in (\lat_{X}\otimes
  \mathbb{R})_{\{5\}}$ be the projection of any vector in
  $\lat_{X}\cap (\frac{\sigma}{8} + \sigma^{\perp})$, e.g. $\acan$ can be the
  projection of the vector $2\delta_{L_{\varnothing}}$. Then the two
  conditions on the components $d_{\{5\}}$ and $d_{\{5\}}'$ in
  \eqref{eq:final.cond} are equivalent to the conditions that there
  exists an $\mathfrak{s} \in \mathfrak{S}$ so that if
  $\mathsf{x}^{\mathfrak{s}} = \frac{1}{8}(\unone -
2\incdelta)\left(\mathfrak{s}_{1}(r_{1},r_{2}) -
\mathfrak{s}_{2}(r_{1},r_{2})\right)$ we have 
  \begin{equation} \label{eq:cosets}
 \begin{aligned} 
   d_{\{5\}}  & \in  \acan - \mathsf{x}^{\mathfrak{s}}  +
     \frac{1}{2}\lat_{X}^{0}, \\
     d_{\{5\}}'  & \in  \acan +  \mathsf{x}^{\mathfrak{s}} + 
      \frac{1}{2}\lat_{X}^{0}.
 \end{aligned}
\end{equation}
  \end{itemize}
\end{lem}
{\bfseries Proof.} Let $\op{pr}_{\{5\}} : \lat_{X}\otimes \mathbb{R}
\to \lat_{X}\otimes \mathbb{R}$ be the projection on the
$5$-diemnsional eigenspace for $\bI$ corresponding to eigenvalue $-4$.
Since $\sigma$ is an $\bI$-eigenvector with eigenvalue $4$, it
follows, that the restrction of $\op{pr}_{\{5\}}$ to
$\sigma^{\perp}\otimes \mathbb{R}$ is equal to the restriction of
$-\frac{1}{4}\bI$ to $\sigma^{\perp}\otimes \mathbb{R}$. In particular
$\op{pr}_{\{5\}} : \sigma^{\perp}\otimes \mathbb{R} \to
(\lat_{X}\otimes \mathbb{R})_{\{5\}}$ is surjective. This proves (a).

\

\noindent
An integral basis of $\sigma^{\perp}$ is given by the vectors
$\delta_{L_{\varnothing}} - \delta_{L_{I}}$, for $I \neq \varnothing$.
Let $v_{I} = -\frac{1}{4}\bI(\delta_{L_{\varnothing}} -
\delta_{L_{I}})$. The $15$ vectors $v_{I}$ span $\lat_{X}^{0}$ over
$\mathbb{Z}$.  The first five vectors $v_{1^{c}}$, $v_{2^{c}}$,
$v_{3^{c}}$, $v_{4^{c}}$, $v_{5^{c}}$ span $\lat_{X}^{0}$ over
$\mathbb{Q}$, and the remaining $10$ vectors are $v_{\{ij\}} =
\frac{1}{2} (\sum_{k=1}^{5} v_{k^{c}}) - v_{i^{c}} - v_{j^{c}}$.  Thus
$\lat_{X}^{0}$ is spanned over $\mathbb{Z}$ by the vector
$\mathfrak{v} = \frac{1}{2} (\sum_{k=1}^{5} v_{k^{c}})$ and any four
of the five vectors $v_{1^{c}}$, $v_{2^{c}}$, $v_{3^{c}}$,
$v_{4^{c}}$, $v_{5^{c}}$.  On the other hand, by definition the image
of the lattice $\lat_{C}$ under $\frac{1}{4}(\unone - 2\incdelta)$ is
spanned by the five vectors $\frac{1}{4}(\unone -
2\incdelta)(\delta_{p_{i}}) = 2v_{i^{c}} - \mathfrak{v}$.  This shows
that $\frac{1}{4}(\unone - 2\incdelta)(\lat_{C})$ is a sublattice of
finite index in $\lat_{X}^{0}$ and proves (b).

\

\noindent
Part (c) follows from (b), the two equations involving $d_{\{5\}}$ and
$d_{\{5\}}'$ in \eqref{eq:fixedab.cond}, and \linebreak
Remark~\ref{eq:final.cond}.  \ \hfill $\Box$

\

\noindent
We have found a solution of all constraints in cohomology which is
unique modulo the obvious redundancies.  Uniqueness implies that these
numerics also solve the slightly stronger abelianized Hecke eigensheaf
constraint in the Picard. In summary we have the following
proposition.

\begin{prop} \label{prop:numerics} Fix a pair of parabolic weights
  on $\Par_{C}$, specified by a pair
  of real vectors $a, b \in \mathbb{R}^{5}$ that are generic
  except for the linear relation $\sum_{i = 1}^{5} (a_{i}+b_{i}) = 3$
  (which, by \eqref{eq:spectral.parch}, expresses the condition that
  $\parch_{1}(E_{\bullet}) = 0$). Also, fix a pair of integral vectors
  $n_{1}, n_{2} \in \mathbb{Z}^{5}$ such that $\sum_{i = 1}^{5}
  (n_{1,i}+ n_{2,i}) = -3$.  Then
\begin{description}
\item[Hecke kernel:] The pair $(a,b)$ determines a natural parabolic
  line bundle $\mycal{I}_{\bullet} =
  \mathcal{O}_{H}(\bzeta\Par_{H})_{\bullet}$ on $(H,\Par_{H})$ with
  $\parch_{1}(\mycal{I}_{\bullet}) = 0$.
\item[Hecke eigensheaf:] For any tame strongly
  parabolic rank two Higgs bundle $(E_{\bullet},\theta)$ on
  $(C,\Par_{C})$ with weights specified by $(a,b)$
  and a smooth spectral cover $\sC$, there exists
  a natural Hecke eigensheaf on the (disconnected) moduli
  space. Explicitly, this is specified by a pair of objects  
  $(F_{\bullet},\varphi)$ and $(F_{\bullet}',\varphi')$ on the del
  Pezzo $X$, so that
  \begin{itemize}
  \item $(F_{\bullet},\varphi)$ and $(F_{\bullet}',\varphi')$ are
    stable, strongly parabolic, rank
    four Higgs bundles on $(X,\Par_{X})$ with $\parch_{1} = 0$,
    $\parch_{2} = 0$.
 \item The parabolic Hecke eigensheaf condition (see \eqref{eq-par.hecke})
\[
\boxed{q_{*} (p^{*}(F_{\bullet},\varphi)\otimes (\mycal{I}_{\bullet},0)) = 
p_{X}^{*}(F_{\bullet}',\varphi')\otimes p_{C}^{*}(E_{\bullet},\theta)}
\]
holds.
\item The parabolic structures on $(F_{\bullet},\varphi)$ and
  $(F_{\bullet}',\varphi')$ each have two jumps, of rank one and three,
  on every component $L_{I}$ of $\Par_{X}$.
\item The parabolic weights of $(F_{\bullet},\varphi)$ and
  $(F_{\bullet}',\varphi')$ are specified by vectors $e, d, e', d' \in
  \mathbb{R}^{\evenL} = \mathbb{R}^{16}$ which can be taken to be
  \begin{equation} \label{eq:basic.values}
\left| \; 
    \begin{aligned}
e & = \frac{5}{8}\sigma - 3\acan +\frac{3}{8}(\unone  - 2\incdelta)
(r_{1} - r_{2}), \\[+0.5pc]
d & = \frac{1}{8}\sigma + \acan
- \frac{1}{8}(\unone - 2\incdelta)(r_{1} - r_{2}), \\[0.5pc]
e' & = \frac{5}{8}\sigma  -
3\acan - \frac{3}{8}(\unone  - 2\incdelta)
(r_{1} - r_{2}), \\[+0.5pc]
d' & = \frac{1}{8}\sigma + \acan
      + \frac{1}{8}(\unone - 2\incdelta)(r_{1} - r_{2}).
\end{aligned}
\right.
    \end{equation}
    where we have used the normalizations
    $\mathfrak{u} = 2\op{pr}_{\{5\}}(\delta_{L_{\varnothing}})$,
    $r_{1} = a + n_{1}$ and $r_{2} = b + n_{2}$.
\item Both $(F_{\bullet},\varphi)$ and $(F_{\bullet}',\varphi')$
  correspond to spectral data specified on a spectral cover of $X$
  which is birational to the Hitchin fiber through
  $(E_{\bullet},\theta)$.
\end{itemize}
\end{description}
\end{prop}
{\bfseries Proof.} Let $e$, $d$, $e'$, $d'$ be given by the formulas
\eqref{eq:basic.values}. Set $A = 2\delta_{L_{\varnothing}}$, $B =
-2\delta_{L_{\varnothing}}$, $C = A + (\bI + \one)B$, $r_{1} =
a+n_{1}$, $r_{2} = b + n_{2}$ and $\bzeta =
(A-d,B+d',C+e'-d,r_{1},r_{2})$. Suppose $\pi : \sC \to C$ is a smooth
spectral cover and let $f : Y \to X$ be the corresponding modular spectral cover. Given a line bundle  $\mathfrak{a} \in \Pic^{0}(\sC)$,  let
$(E_{\bullet},\theta)$ be the strongly parabolic stable rank two Higgs bundle
given by the spectral parabolic line bundle $\mathfrak{a}(a\sP +
b\sP)_{\bullet,\bullet}$ on $(\sC,\sP+\sP)$. 

Define $(F_{\bullet},\varphi)$ and $(F_{\bullet}',\varphi')$ to be the
rank four parabolic Higgs bundles on $X$ given by the modular spectral
cover $f : Y \to X$ and the spectral parabolic line bundles \linebreak
$\mycal{L}_{\mathfrak{a}}(e E + (e+d)G)_{\bullet,\bullet}$ and
$\mycal{L}_{\mathfrak{a}}(e' E + (e'+d')G)_{\bullet,\bullet}$
respectively. The spectral correspondence combined with the fact that
$Y$ is smooth and irreducible immediately implies that the Higgs
bundles on $X$ are stable. Also, in section \ref{sssec:eigensheaf} we
checked that any parabolic Higgs bundle on $X$ given by spectral data
on $Y$ is automatically strongly parabolic.

Finally we chose the formulas defining $e$, $d$, $e'$, $d'$ and
$\bzeta$ so that they solve the non-abelian Hodge condition
\eqref{eq:nah.constraints}, the Hecke eigensheaf condition
\eqref{eq:ABCn12}, and the Hecke kernel condition
\eqref{eq:final.cond} (with $\mathfrak{s} = \op{id}$). Therefore we
have
\[
q_{*} (p^{*}(F_{\bullet},\varphi)\otimes (\mycal{I}_{\bullet},0)) = 
p_{X}^{*}(F_{\bullet}',\varphi')\otimes p_{C}^{*}(E_{\bullet},\theta)
\]
which completes the proof of the proposition. \ \hfill $\Box$

\

\noindent
Once the basic Hecke relation in Proposition~\ref{prop:numerics}(b) is
established, it is straighforward to obtain the relations proving the
parabolic Higgs version of the geometric Langlands correspondence on
$(C,\Par_{C})$ in the cases when $\lan{G} = GL_{2}(\mathbb{C})$ and
$\lan{G} = SL_{2}(\mathbb{C})$.

\

\begin{prop} \label{prop:gl2.and.sl2}
Fix generic parabolic weights on $\Par_{C}$ specified by a pair of
real vectors \linebreak $a, b \in \lat_{C}\otimes \mathbb{R} \cong
\mathbb{R}^{5}$ s.t. $\sum_{i = 1}^{5} (a_{i} + b_{i}) = 3$. Then

\

\noindent
\textup{\bfseries (i)} \mbox{\boldmath{$[\lan{G} = GL_{2}(\mathbb{C})]$}}
\ There exists a collection
    of parabolic line bundles
$\left\{\mycal{I}_{\bullet}^{k}\right\}_{k\in \mathbb{Z}}$ on
$(H,\Par_{H})$ with $\parch_{1}(\mycal{I}_{\bullet}^{k}) = 0$ for all
$k$ and such that for any tame strongly parabolic rank two Higgs
bundle $(E_{\bullet},\theta)$ on $(C,\Par_{C})$ with parabolic weights
given by $(a,b)$ and a smooth spectral
cover, there exists a collection $\left\{
(F_{\bullet}^{k},\varphi^{k}) \right\}_{k \in \mathbb{Z}}$ so that for
all $k \in \mathbb{Z}$ we have
\begin{itemize}
  \item $(F_{\bullet}^{k},\varphi^{k})$ is a stable, strongly
    parabolic, rank four Higgs bundle on $(X,\Par_{X})$ with
    $\parch_{1} = 0$, $\parch_{2} = 0$.
 \item The parabolic Hecke eigensheaf condition
\[
q_{*} (p^{*}(F_{\bullet}^{k},\varphi^{k})\otimes
  (\mycal{I}_{\bullet}^{k},0)) = 
p_{X}^{*}(F_{\bullet}^{k+1},\varphi^{k+1})\otimes p_{C}^{*}(E_{\bullet},\theta)
\]
holds on $X\times C$.
\end{itemize}

\

\noindent
\textup{\bfseries (ii)}
\mbox{\boldmath{$[\lan{G} = SL_{2}(\mathbb{C})]$}} \ If
in addition $a + b \in \mathbb{Z}^{5}$ there exist natural parabolic
line bundles $\mathcal{O}(\bzeta\Par_{H})$ and 
$\mathcal{O}(\bzeta'\Par_{H})_{\bullet}$ on $H$ and   
$(\mathbb{Z}/2)^{4}$-equivariant stable strongly parabolic Higgs bundles
$(F_{\bullet},\varphi)$ and 
$(F_{\bullet}',\varphi')$  on $X$ satisfying
\[
\begin{aligned}
q_{*}\left(p^{*}(F_{\bullet},\varphi)\otimes
(\mathcal{O}(\bzeta\Par_{H})_{\bullet},0)\right) & =
p_{X}^{*}(F_{\bullet}',\varphi')\otimes p_{C}^{*}(E_{\bullet},\theta), \\[0.5pc]
q_{*}\left(p^{*}(F_{\bullet}',\varphi')\otimes
(\mathcal{O}(\bzeta'\Par_{H})_{\bullet},0)\right) & =
p_{X}^{*}(F_{\bullet},\varphi)\otimes p_{C}^{*}(E_{\bullet},\theta).
\end{aligned}
\]

\end{prop}
{\bfseries Proof.} This follows directly from
Proposition~\ref{prop:numerics}. Explicitly fix two integral vectors
$n_{1}, n_{2} \in \lat_{C} \cong \mathbb{Z}^{5}$ such that $\sum_{i =
  1}^{5} (n_{1,i} + n_{2,i}) = -3$ and let $r_{1} = a + n_{1}$, $r_{2}
= b + n_{2}$ and $\mathfrak{u} =
\op{pr}_{\{5\}}\left(\delta_{L_{\varnothing}}\right)$.
Now for every $k \in
\mathbb{Z}$ define
\[
A^{k} = 2\delta_{L_{\varnothing}}, \
B^{k} = -2\delta_{L_{\varnothing}}, \
C^{k} = A^{k} + (\bI + \one)B^{k}, \
n_{1}^{k} = n_{1 + (k\op{mod} 2)} , \
n_{2}^{k} = n_{2 - (k\op{mod} 2)} 
\]
and also $r_{1}^{k} = r_{1 + (k\op{mod} 2)}$,
$r_{2}^{k} = r_{2 - (k\op{mod} 2)}$, and 
\[
\left| \;
\begin{aligned}
  e^{k} & = \frac{5}{8}\sigma  - 3u + 
  \frac{3}{8}(\unone - 2\incdelta)(r_{1}^{k} - r_{2}^{k}) \\[+0.5pc] 
  d^{k} & = \frac{1}{8}\sigma + u -
  \frac{1}{8}(\unone - 2\incdelta)(r_{1}^{k} - r_{2}^{k})
\end{aligned}
\right.
\]
Finally set
\[
\bzeta^{k}  = (A^{k} - d^{k},B^{k} + d^{k+1},C^{k} + e^{k+1} - d^{k},
r_{1}^{k}, r_{2}^{k}),
\]
Now applying Lemma~\ref{lem:cosets}(c) with $\mathfrak{s} \in
\mathfrak{S}$ taken to be the total swap
$\mathfrak{s}(\mathsf{y}_{1},\mathsf{y}_{2}) =
(\mathsf{y}_{2},\mathsf{y}_{1})$ we see that
$\bzeta^{k}$ defines a Hecke kernel
parabolic line bundle $\mycal{I}^{k}_{\bullet} =
\mathcal{O}_{H}(\bzeta^{k}\Par_{H})_{\bullet})$ with
$\parch_{1}(\mycal{I}^{k}_{\bullet}) = 0$.

Now suppose  $(E_{\bullet},\theta)$ is a strongly parabolic rank
two Higgs bundle  on $(C,\Par_{C})$ with
parabolic weights given by $(a,b)$ and a smooth spectral cover
$\pi : \sC \to C$. Then $(E_{\bullet},\theta)$ corresponds to a line bundle
$\mathfrak{a} \in \Pic^{0}(\sC)$, i.e. for every parabolic level $t
\in \mathbb{R}^{5}$ we have $(E_{t},\theta) =
(\pi_{*}\mathfrak{a}((\fl*{t+a} + \fl*{t+b})\sP),\pi_{*}(\lambda\otimes
-))$.

Let $f : Y \to X$ be the modular spectral cover of $X$ corresponding
to $\pi : \sC \to C$.  Then by Proposition \ref{prop:numerics} the
spectral data $\mycal{L}_{\mathfrak{a}}(e^{k}E +
d^{k}(E+G))_{\bullet,\bullet}$ 
on $Y$ defines a rank four strongly parabolic Higgs bundles
$(F^{k},\varphi^{k})$ on $(X,\Par_{X})$ with
vanishinig $\parch_{1}$ and $\parch_{2}$, and such that the parabolic
Hecke eigensheaf condition
\[
q_{*} (p^{*}(F_{\bullet}^{k},\varphi^{k})\otimes
  (\mycal{I}_{\bullet}^{k},0)) = 
p_{X}^{*}(F_{\bullet}^{k+1},\varphi^{k+1})\otimes p_{C}^{*}(E_{\bullet},\theta)
\]
holds on $X\times C$. This proves part {\bfseries (i)} of the proposition.

\index{terms}{total swap}
\index{notations}{ifrakb@$\bfri$}
\index{notations}{zeta@$\bzeta$}
\index{notations}{e@$e$}
\index{notations}{d@$d$}
\index{terms}{Hecke!kernel}

\

\noindent
For part {\bfseries (ii)} take $\bzeta = \bzeta^{0}$, $\bzeta' =
\bzeta^{1}$, $e = e^{0}$, $e' = e^{1}$, $d = d^{0}$, $d' = d^{1}$, and
also $(F_{\bullet},\varphi) = (F_{\bullet}^{0},\varphi^{0})$ and
$(F_{\bullet}',\varphi) = (F_{\bullet}^{1},\varphi^{1})$. As we just
saw in the proof of part {\bfseries (i)} these choices automatically
satisfy
\[
\begin{aligned}
q_{*}\left(p^{*}(F_{\bullet},\varphi)\otimes
(\mathcal{O}(\bzeta\Par_{H})_{\bullet},0)\right) & =
p_{X}^{*}(F_{\bullet}',\varphi')\otimes p_{C}^{*}(E_{\bullet},\theta), \\[0.5pc]
q_{*}\left(p^{*}(F_{\bullet}',\varphi')\otimes
(\mathcal{O}(\bzeta'\Par_{H})_{\bullet},0)\right) & =
p_{X}^{*}(F_{\bullet},\varphi)\otimes p_{C}^{*}(E_{\bullet},\theta).
\end{aligned}
\]
The $(\mathbb{Z}/2)^{4}$-equivariance of all objects is checked by an
explicit computation which is very similar to the computations
simplifying the Hecke kernel conditions. \ \hfill $\Box$

\

\noindent
Finally let us note that Mochizuki's non-abelian Hodge theory
constraints on the modular spectral data and on the parabolic Hecke
kernel line bundle are the most restrictive ones. Indeed, as the
following claim shows, imposing non-abelian Hodge theory constraints
already ensures automorphicity. We thank the referee for
suggesting this might be true.

\

\begin{claim} \label{claim:spectral.is.eigensheaf}
Fix a smooth spectal cover $\pi :\sC \to C$ ramified over $\Par_{C}$,
and let $f : Y \to X$ be the corresponding modular spectral
cover. Given vectors $e, d$  in $\lat_{X}\otimes
\mathbb{R}$ that obey the non-abelian Hodge theory constraints
\eqref{chsoln}, there exist vectors $\bzeta = (lp,lq,llc,r_{1},r_{2})$,
$e', d' \in \lat_{X}\otimes
\mathbb{R}$
and vectors $a, b \in \lat_{C}\otimes \mathbb{R}$ so that $e', d'$ also obey
\eqref{chsoln}, and if
$(F_{\bullet},\varphi)$ and $(F_{\bullet}',\varphi')$ are the strongly
parabolic rank four Higgs bundles on $(X,\Par_{X})$ given by modular
spectral data $\left(Y,\mycal{L}_{\mathfrak{a}}(eE +
d(E+G))_{\bullet,\bullet}\right)$ and
$\left(Y,\mycal{L}_{\mathfrak{a}}(e'E +
d'(E+G))_{\bullet,\bullet}\right)$ with $\mathfrak{a} \in
\Pic^{0}(\sC)$, then
\[
q_{*}\left(p^{*}(F_{\bullet},\varphi)\otimes
(\mathcal{O}(\bzeta\Par_{H})_{\bullet},0)\right)  =
p_{X}^{*}(F_{\bullet}',\varphi')\otimes p_{C}^{*}(E_{\bullet},\theta),
\]
where $(E_{\bullet},\theta)$ is the strongly parabolic rank two Higgs
bundle on $(C,\Par_{C})$ corresponding to the spectral data
$(\sC,\mathfrak{a}(a\widetilde{P} +
b\widetilde{P}))_{\bullet,\bullet}$.
\end{claim}
{\bfseries Proof.} \ This follows immediately from the Hecke kernel
condition \ref{eq:fixedab.cond}. Indeed fix again the normalizations
$\acan = \op{pr}_{\{5\}}(2\delta_{L_{\varnothing}})$ and $A =
2\delta_{L_{\varnothing}}$, $B = -2\delta_{L_{\varnothing}}$, $C = A +
(\bI + \one)B$. Now let  $r_{1}, r_{2} \in
\lat_{C}\otimes{\mathbb R}$ be such that 
\[
(\unone - 2\incdelta)(r_{1} - r_{2}) = 8(\acan - d_{\{5\}}), \quad
  r_{1} + r_{2} = 0.
\]
Note that by Lemma~\ref{lem:cosets}(b) the map $(\unone - 2\incdelta) :
\lat_{C}\otimes \mathbb{R} \to (\lat_{X}\otimes \mathbb{R})_{\{5\}}$ is a
linear isomorphism so these equations have a unique solution.  Now
choose any pair of integral vectors $n_{1}, n_{2} \in \lat_{C}$
satisfying $\sum_{i = 1}^{5}(n_{1,i} + n_{2,i}) = -3$ and set $a =
r_{1} - n_{1}$ and $b = r_{2} - n_{2}$. With these choices we can now define
\[
\begin{aligned}
  e' & = \frac{5}{8}\sigma - 3\acan + \frac{3}{8} (\unone -
  2\incdelta)(r_{1} - r_{2}), \\[0.5pc]
  d' & = \frac{1}{8}\sigma +
  \acan - \frac{1}{8} (\unone - 2\incdelta)(r_{1} - r_{2}),
  \\[0.5pc]
  \bzeta & = (A - d, B+ d', C + e' -d, r_{1}, r_{2}).
\end{aligned}
\]
By construction these choices produce spectral line bundles on $\sC$
and on $Y$ so that the corresponding Higgs bundles satisfy the
non-abelian Hodge, Hecke kernel, and Hecke eigensheaf conditions.
\ \hfill $\Box$

\subsection{The Okamoto map} \label{ssec:okamoto}

In section~\ref{ssec:class.of.the.kernel} we used the Fourier-Mukai
transform along the Hitchin fibers and solved the vanishing of
parabolic Chern classes constraints and the Hecke eigensheaf
conditions to construct a canonical $(\text{eigenvalue}) \to
(\text{eigensheaf})$ assignment in the context of parabolic Higgs
bundles: \index{terms}{Hecke!eigensheaf}
\index{terms}{Fourier-Mukai!transform} \index{terms}{parabolic!Chern
  classes}

\begin{spacing}{1.0}
\begin{equation} \label{eq:parHeckeclassical}
  \left(\, \text{\begin{minipage}[c]{2in}
  stable purely imaginary rank two parabolic Higgs bundles
  $(E_{\bullet},\theta)$ on $(C,\Par_{C})$
  \end{minipage}}\,\right) \quad \longrightarrow \quad
  \left(\,\text{\begin{minipage}[c]{2in}
     stable purely imaginary rank four
  parabolic Higgs bundles $(F_{\bullet},\varphi)$ on
  $(X,\Par_{X})$
  \end{minipage}} \,\right).
\end{equation}
\end{spacing}

\vspace{1pc}

\noindent
One essential numerical component of this construction was the
computation of the assignment of $(a,b) \mapsto (e,d)$ which explains
how the parabolic weights of $(E_{\bullet},\theta)$ on the curve $C$
are converted into parabolic weights of $(F_{\bullet},\varphi)$ on the
moduli space $X$. Non-abelian Hodge theory identifies the assignment
$(a,b) \mapsto (e,d)$ with the \emph{\bfseries Okamoto map}, which is
an important ingredient of the tamely ramified GLC. As explained in
Appendix~\ref{app-ram.glc}, the Okamoto map records the way the
eigenvalues of the residue of the parabolic flat bundle on $C$ are
converted into the twisting parameters for the associated twisted
$\mathcal{D}$-module on $\sBun$.  In this section we work out the
non-abelian Hodge theory translation and compute the Okamoto map in
our setting.

\index{terms}{Okamoto!map}

In more detail (see also Appendix~\ref{app-ram.glc}), the inputs and
outputs of the assignment \eqref{eq:parHeckeclassical} can be
converted by application of Mochizuki's tame non-abelian Hodge
correspondence into stable tame parabolic flat bundles $\mathbb{V} =
(V_{\bullet},\nabla)$ on \index{terms}{Non-abelian Hodge theory}
$(C,\Par_{C})$ and $(W_{\bullet},\nabla)$ on $(X,\Par_{X})$
respectively. Furthermore using the twisted Deligne-Goresky-MacPherson
extension from $X - \Par_{X}$ to the stack of quasi-parabolic bundles
on $(C,\Par_{C})$, we can rewrite $(W_{\bullet},\nabla)$ as a twisted
$\mathcal{D}$-module.
\index{terms}{Dmodule@$\mathcal{D}$-module!twisted} Under these
identifications the assignment \eqref{eq:parHeckeclassical} becomes
the evaluation
\[
\mathcal{O}_{\mathbb{V}} \mapsto
\mathfrak{c}_{-\blambda}(\mathcal{O}_{\mathbb{V}}) 
\]
of the Langlands correspondence functor (see
section~\ref{ssec:setup.glc} in Appendix~\ref{app-ram.glc})
\[
\mathfrak{c}_{-\blambda} :
D(\lan{\gLoc}_{-\blambda},\mathcal{O})
\longrightarrow
D(\sBun,\mathcal{D}_{\Aom(-\blambda)})
\]
on the skyscraper sheaf supported at the point $\mathbb{V} \in
\lan{\gLoc}_{-\blambda}$.
\index{terms}{sheaf!skyscraper}

For concreteness let us focus on the setup when $\lan{G} =
SL_{2}(\mathbb{C})$, $G = \mathbb{P}SL_{2}(\mathbb{C})$, \linebreak
$\blambda :
\Par_{C} \to \lan{\mathfrak{t}}^{\op{cpt}}_{\mathbb{R}} \subset
\lan{\mathfrak{t}}$ are the parabolic weights for
$(E_{\bullet},\theta)$, $\lan{\gLoc}_{-\blambda}$ is the moduli stack
of quasi-parabolic flat $SL_{2}(\mathbb{C})$ bundles
$(V,\nabla,\mathsf{r})$ on $(C,\Par_{C})$ with $\res^{\mathsf{r}}\,
\nabla = - \blambda$, and $\sBun$ is the moduli stack of quasi
parabolic $\mathbb{P}SL_{2}(\mathbb{C})$ bundles on $(C,\Par_{C})$.

We want to compute the Okamoto map
\index{terms}{Okamoto!map}
\index{notations}{Ok@$\Aom$}
\[
\Aom : \Gamma(\Par_{C},\lan{\mathfrak{t}}) \longrightarrow
\op{Pic}(\sBun)\otimes \mathbb{C}.
\]
As explained in Appendix~\ref{app-ram.glc}, the condition that
$(E_{\bullet},\theta)$ and $(F_{\bullet},\varphi)$ are purely
imaginary translates under Simpson's conversion table \cite{carlos-nc}
into the statement that the Okamoto map is the projection to 
$\op{Pic}(\sBun)\otimes \mathbb{C}$ of the map

\begin{spacing}{1.0}
\begin{equation} \label{eq:parHeckeclassical.weights}
  \left(\, \text{\begin{minipage}[c]{1.3in}
  parabolic weights of 
  $(E_{\bullet},\theta)$
  \end{minipage}}\,\right) \quad \longrightarrow \quad
  \left(\,\text{\begin{minipage}[c]{1.3in}
     parabolic weights of  $(F_{\bullet},\varphi)$
  \end{minipage}} \,\right).
\end{equation}
\end{spacing}

\vspace{1pc}

\noindent
which is an essential part of  the assignment \eqref{eq:parHeckeclassical}. 

\index{notations}{Vnabla@$(V,\nabla)$}
\index{notations}{LFlats@$\lan{\gLoc}$}
\index{notations}{C@$C$}
\index{notations}{r@$\mathsf{r}$}
\index{notations}{ParC@$\Par_{C}$}
\index{terms}{residue!of a logarithmic connection}

In our setup the parabolic weights of the rank two bundle
$(E_{\bullet},\theta)$ were encoded in the pair $\blambda = (a,b)$, where 
\[
(a,b) \in \bigoplus_{i = 1}^{5}\displaylimits
(\mathbb{R}\delta_{p_{i}})^{\oplus 2} =
\Gamma(\Par_{C},\mathfrak{t}_{U(2)}) \subset
\Gamma(\Par_{C},\mathfrak{t}_{GL_{2}(\mathbb{C})}) = 
\bigoplus_{i = 1}^{5}\displaylimits
(\mathbb{\mathbb{C}}\delta_{p_{i}})^{\oplus 2}.
\]
Similarly, the parabolic weights of the rank four bundle
$(F_{\bullet},\varphi)$
were encoded in the pair
\[
\xymatrix@C-0.5pc@R-3pc@M+0.5pc{
  (e,d) \in \hspace{-2pc} &
  \bigoplus_{I \in \evenL}\displaylimits
  (\mathbb{R}\delta_{L_{I}})^{\oplus 2} \ar@{^{(}->}[r] &
 \bigoplus_{I \in \evenL}\displaylimits
  (\mathbb{\mathbb{C}}\delta_{L_{I}})^{\oplus 2}. \\
  & \cap & \cap \\
  & \Gamma(\Par_{X},\mathfrak{t}_{U(4)}) \ar@{^{(}->}[r] &
  \Gamma(\Par_{X},\mathfrak{t}_{GL_{4}(\mathbb{C})})}
\]
Thus in our setting the map \eqref{eq:parHeckeclassical.weights}
becomes the assignment $(a,b) \mapsto (e,d)$ which is the
$\mathbb{C}$-affine linear extension of the $\mathbb{R}$-affine linear
map found in  section~\ref{ssec:class.of.the.kernel}.

To compute the Okamoto map we need to
combine  \eqref{eq:parHeckeclassical.weights} with two more ingredients:
\begin{itemize}
\item[(i)] The condition on $(a,b)$ that ensures that $E_{\bullet}$ is an
  $SL_{2}(\mathbb{C})$ parabolic bundle.
\item[(ii)] The projection from $\oplus_{I \in \evenL}
  (\mathbb{C}\delta_{L_{I}})^{\oplus 2}$ to
  $\op{Pic}(\sBun)\otimes\mathbb{C}$. 
\end{itemize}
Step (i) is straightforward. A direct calculation with the Weyl
alcoves (see e.g. \cite[Section~2.1]{teleman.woodward} or
\cite[Remark~3.7]{bbp}) shows that a parabolic bundle $E_{\bullet}$
corresponding to a pair
$(a,b) \in \oplus_{i =1}^{5} (\mathbb{R}\delta_{p_{i}})^{\oplus 2}$ will be an
$SL_{2}(\mathbb{C})$ parabolic bundle if and only if we have
$a+b \in \oplus_{i = 1}^{5} \mathbb{Z}\delta_{p_{i}}$.  Recall also (see
equation \eqref{eq:spectral.parch}) that the $\parch_{1}(E_{\bullet})$
condition imposed by non-abelian Hodge theory gives the further
restriction $\sum_{i =1}^{5} (a_{i}+b_{i}) = 3$

Step (ii) is more delicate. First recall (see
section~\ref{sssec:assumption}) that under the purely imaginary
condition, the non-abelian Hodge correspondence identifies $\oplus_{I
  \in \evenL} (\mathbb{R}\delta_{L_{I}})^{\oplus 2}$ with the space
parametrizing eigenvalues of residues of the parabolic flat bundle
$(W_{\bullet},\nabla)$. Next, as we explain in
section~\ref{program}{\bfseries (vii)}, we can use the eigenvalues of
the residues of $\nabla$ on the various components $L_{I}$ of
$\Par_{X}$ to form a twisted Deligne-Goresky-MacPherson extension from
the very stable locus in $X$. By construction this extension is
defined intrinsically on a non-separated version of $X$ in which every
component $L_{I}$ of $\Par_{X}$ is replaced by a number of
non-separated copies of itself, where the copies are labeled by the
distinct eigenvalues of $\res_{L_{I}}\! \nabla$. Since in our setting
the eigenvalues are the negatives of the parabolic parameters $e$ and
$d$, the intrinsic twisted Deligne-Goresky-MacPherson extension of
$(W_{\bullet},\nabla)$ is defined on the space $\mathfrak{X}$ obtained
from $X$ by doubling each $L_{I}$, thus replacing it by two
non-separated divisors $L_{I}^{e}, L_{I}^{d} \subset \mathfrak{X}$.
\index{terms}{DeligneGG@Deligne-Goresky-MacPherson!extension}
Concretely $\mathfrak{X}$ is obtained as the gluing of $X$ and sixteen
surfaces $X\times \{I\}$, $I \in \evenL$, where $X\times \{I\}$ is
glued to $X$ along the naturally identified opens $X - L_{I}$ and
$(X-L_{I})\times \{I\}$. To keep track of divisors we will denote by
$L_{I}^{d}$ the copy of $L_{I}$ sitting inside the chart $X \subset
\mathfrak{X}$ and will write $L_{I}^{e}$ for the copy of $L_{I}$ in
the chart $X\times \{I\}$.  By definition then the intrinsic twisted
Deligne-Goresky-MacPherson extension of $(W_{\bullet},\nabla)$ is a
module over the sheaf of TDO on $\mathfrak{X}$ corresponding to the
element $\mathcal{O}_{\mathfrak{X}}(\sum_{I \in \evenL}
(e_{I}L_{I}^{e} + d_{I}L_{I}^{d}) \in \op{Pic}(\mathfrak{X})\otimes
\mathbb{C}$.  In other words the natural twisted minimal extension
from $X - \Par_{X}$ to $\mathfrak{X}$ corresponds to the map on
twistings given by
\[
\bigoplus_{I \in \evenL} (\mathbb{C}\delta_{L_{I}})^{\oplus 2} \to
\op{Pic}(\mathfrak{X})\otimes \mathbb{C}, \quad (e,d) \mapsto
\mathcal{O}_{\mathfrak{X}}( L^{e}e + L^{d}d).
\]
Thus to understand the map in step (ii) we need to relate the twisted
$\mathcal{D}$-module on the space $\mathfrak{X}$ to a twisted
$\mathcal{D}$-module on $\sBun$.  \index{notations}{Buns@$\sBun$} To
that end recall \cite{arinkin-4points} that for non-resonant twisting
the twisted $\mathcal{D}$-modules on $\sBun$ are all extensions by
zero of the twisted $\mathcal{D}$-modules on the open $\sBun' \subset
\sBun$ parametrizing indecomposable quasi-parabolic
$\mathbb{P}SL_{2}(\mathbb{C})$ bundles on
$(C,\Par_{C})$. \index{notations}{Bunsprime@$\sBun'$} In
\cite{arinkin-lysenko} it is shown that the stack $\sBun'$ is
representable by a non-separated algebraic space and that the
restriction of line bundles from $\sBun$ to $\sBun'$ induces an
isomorphism on Picard groups.

In fact, some time ago, Deligne noted that for
$\mathbb{P}^{1}$ with any number of points, a quasi-parabolic
$\mathbb{P}SL_{2}(\mathbb{C})$ bundle is indecomposable if and only if
it is stable for some choice of parabolic weights. Thus the moduli
space $\sBun'$ is glued from GIT moduli spaces of stable
quasi-parabolic bundles for various choices of weights. In a private
communication \cite{arinkin-letter} Arinkin pointed out that our
analysis of the moduli of stable parabolic bundles for the various
chambers of weights leads to an explicit description of $\sBun'$. A
similar description  was also carried out in the recent work of
Komyo-Saito \cite{komyo-saito}.  Concretely $\sBun'$ contains the del
Pezzo surface $X$ so that the complement $\sBun' - X$ is a finite
set. The entire space $\sBun'$ can be constructed
\cite{arinkin-letter,komyo-saito} from $X$ as follows. For every $I
\in \evenL$ denote by $X_{I}$ the blow down of $L_{I} \subset X$, and
let $x_{I} \in X_{I}$ denote the image of $L_{I}$. Then $\sBun'$ is
obtained as the gluing of $X$ and all the $X_{I}$ where $X$ is glued
to $X_{I}$ along the open $X - L_{I} = X_{I} - x_{I}$.  Therefore the
closure of each line $L_{I} \subset X$ in $\sBun'$ is isomorphic to
$\mathbb{P}^{1}$ with five points doubled.  In particular we see that
restrictions of line bundles to the two opens $X \subset \sBun'
\subset \sBun$ induce identifications of Picard groups $\op{Pic}
(\sBun) = \op{Pic}(\sBun') = \op{Pic}(X)$.

The space $\mathfrak{X}$ maps naturally to $\sBun'$ where the map
contracts all divisors $L_{I}^{e}$ to the non-separated points $x_{I}
\in \sBun'$ and maps the divisors $L_{I}^{d}$ isomorphically to $L_{I}
\subset X \subset \sBun'$. Thus for $\blambda = (a,b)$ the
automorphic TDO $\mathfrak{c}_{-\blambda}(\mathcal{O}_{\mathbb{V}})$ is
given by the $\mathcal{O}_{\sBun'}(dL)$-twisted
Deligne-Goresky-MacPherson extension of $(W_{\bullet},\nabla)_{|U}$ to
$\sBun'$ where extending along $L_{I}$ corresponds to the eigenvalue
$d_{I}$ of $\res_{L_{I}}\! \nabla$.
\index{terms}{DeligneGG@Deligne-Goresky-MacPherson!extension}

\

\begin{rem} \label{rem:Dmodule.is.a.pullback}
  The previous discussion tacitly includes the statement that the
  intrinsic twisted minimal extension of $(W_{\bullet},\nabla)_{|U}$
  to $\mathfrak{X}$ (corresponding to \emph{\bfseries all} eigenvalues
  $e_{I}$ and $d_{I}$) and the twisted minimal extension
  $(W_{\bullet},\nabla)_{|U}$ to $\sBun'$ (corresponding to the
  eigenvalues $d_{I}$ only) are compatible. Concretely this means that
  the twisted $\mathcal{D}$-module $(j_{U \subset
    \mathfrak{X}})_{!*}^{\mathcal{O}_{\mathfrak{X}}(eL^{e}+dL^{d})}
  \left( (W_{\bullet},\nabla)_{|U}\right)$ on $\mathfrak{X}$ is equal
  to the pullback of the twisted $\mathcal{D}$-module $(j_{U \subset
    \sBun'})_{!*}^{\mathcal{O}_{\sBun'}(dL^{d})} \left(
  (W_{\bullet},\nabla)_{|U}\right)$ on $\sBun'$ under the blow down
  map $\mathfrak{X} \to \sBun'$. This feature is expected from the
  Hecke eigen $\mathcal{D}$-modules arising in the geometric Langlands
  correspondence and it can be checked directly \cite{arinkin-letter}
  for the automorphic $\mathcal{D}$-modules obtained by the
  Drinfeld-Laumon construction.  In our setting the compatibility
  amounts to checking that for each divisor $L_{I}^{e} \subset
  \mathfrak{X}$ the sheaf of nearby cycles of $(j_{U \subset
    \mathfrak{X}})_{!*}^{\mathcal{O}_{\mathfrak{X}}(eL^{e}+dL^{d})}
  \left( (W_{\bullet},\nabla)_{|U}\right)$ at $L_{I}^{e}$ is
  trivial. This follows from the compatibility of the six operations with
  the non-abelian Hodge correspondence, and the fact that the
  eigenvalues $e_{I}$ correspond to the exceptional divisors $E_{I}
  \subset Y$ on the modular spectral cover and that the tautological
  one form defining the Higgs field $(F_{\bullet},\varphi)$ does not
  have poles along the $E_{I}$.
\end{rem}

\

\noindent
With these preliminaries in place we can now assemble everything
together and compute the Okamoto map. As we saw, we can view the
Okamoto map as a map that takes values in $H^{2}(X,\mathbb{C}) =
\op{Pic}(X)\otimes \mathbb{C} = \op{Pic}(\sBun')\otimes \mathbb{C} =
\op{Pic}(\sBun)\otimes \mathbb{C}$. It is instructive to note that the
map on parabolic weights computed in
section~\ref{ssec:class.of.the.kernel} is affine linear while the
Okamoto symmetry relating the residues of the connection on $C$ to the
twisting of the corresponding twisted $\mathcal{D}$-module on $\sBun$
is linear.  The discrepancy arises because in our discussion of the
ramified geometric Langlands problem we did not trace carefully the
normalization on twistings coming from trivializing the
non-commutative deformation direction in the quantum Langlands
correspondence \cite{dennis-quantum,travkin}. For a simple group $G$ this
normalization is always a shift by a multiple of the canonical class
of $\sBun$. To avoid subtleties of discussing such normalization we
will consider only the values of the Okamoto map modulo the span of
the canonical class, i.e. we will consider the reduced linear map
\[
\underline{\Aom} : \mathbb{C}^{5} \to H^{2}(X,\mathbb{C})/\mathbb{C}K_{X}
\]
which on the nose will be equal to the linear part of the affine
linear map on parabolic weights.

Let $a,b \in \oplus_{i = 1}^{5} \mathbb{C}\delta_{p_{i}}$ be the
two vectors defining the purely imaginary
$SL_{2}(\mathbb{C})$ parabolic Higgs bundle $(E_{\bullet},\theta)$
with $\parch_{1} = 0$. As explained above this means that the pair
$(a,b) \in \oplus_{i = 1}^{5} (\mathbb{R}\delta_{p_{i}})^{\oplus 2}$ satisfies
$a + b \in \oplus_{i = 1}^{5}\mathbb{Z}\delta_{p_{i}}$, and
$\sum_{i = 1}^{5} (a_{i} + b_{i}) = 3$. Let $(e,d)$ be the
corresponding parabolic weights for the Hecke eigensheaf
$(F_{\bullet},\varphi)$ found in
section~\ref{ssec:class.of.the.kernel}. Explicitly $e$ and $d$ break
uniquely into a sum of components
\[
  \begin{aligned}
    e & = e_{\{1\}} + e_{\{5\}} + e_{\{10\}} \\
    d & = d_{\{1\}} + d_{\{5\}} + d_{\{10\}} 
    \end{aligned}
\]
which are eigenvectors of the intersection matrix $\mathbb{I}$ for eigenvalues 
$4$, $-4$, and $0$ respectively. In
section~\ref{ssec:class.of.the.kernel} we found that these components
can be expressed in terms of $(a,b)$ as 
\[
  \left| \begin{aligned} e_{\{1\}} & = \frac{5}{8}\sigma \\
      e_{\{5\}} & = -3 d_{\{5\}} \\
      e_{\{10\}} & = d_{\{10\}} 
    \end{aligned} \right.  \qquad \text{and} \qquad
  \left| \begin{aligned}
      d_{\{1\}} & = \frac{1}{8}\sigma \\
      d_{\{5\}} & = \acan - \frac{1}{8}\left[ (\unone -
        2\incdelta)((a+n_{1}) - (b+n_{2}))\right] \\
      d_{\{10\}} & = \text{unconstrained}
    \end{aligned} \right.,
\]
where
\begin{itemize}
\item $\sigma$ is a column vector with $16$ components all of which
  are equal to $1$.
\item $n_{1}, n_{2} \in \mathbb{Z}^{5}$ are
  arbitrary vectors satisfying $\sum_{i = 1}^{5} (n_{1,i} + n_{2,i}) =
  -3$.
\item $\acan = \op{pr}_{\{5\}}(2\delta_{L_{\varnothing}})$.
\end{itemize}

\

\noindent
The vectors $n_{1}$ and $n_{2}$ appear in these formulas to account
for the redundancy arising when parametrizing the parabolic weights
on $(C,\Par_{C})$ in terms of the pair $(a,b)$. Since we have the
freedom of choosing $n_{1}$ and $n_{2}$ we can simplify the formulas by taking
\[
\left| \begin{aligned}
  n_{1} & = 0, \\ n_{2} & = - (a+b). \end{aligned}\right.
\]
Thus $a+n_{1} = a$ and $b+n_{2} = -a$ and hence
\[
\left|
\begin{aligned}
  d_{\{1\}} & = \frac{1}{8}\sigma \\
  d_{\{5\}} & = \acan - \frac{1}{4}
  \left(\unone  - 2\incdelta\right)a \\
  d_{\{10\}} & = \text{ uncostrained}
\end{aligned}
\right.
\]
Since $d_{\{10\}}$ is in the kernel of the intersection matrix
$\mathbb{I}$ we have that $Ld_{\{10\}} = 0$ in
$H^{2}(X,\mathbb{C})$. Also since $\sigma/8  =
\op{pr}_{\{1\}}(2\delta_{L_{\varnothing}})$ and
$\acan = \op{pr}_{\{5\}}(2\delta_{L_{\varnothing}})$ we have that\footnote{Recall
  our convention that $L$ denotes the row vector with entries $L_{I}$ and $d$
  denotes  the column vector with entries $d_{I}$.}
\begin{equation} \label{eq:dL}
  Ld = 2L_{\varnothing}  - \frac{1}{4} L\left( \unone - 2\incdelta \right) a
  \quad \text{in} \ H^{2}(X,\mathbb{C}).
\end{equation}
We want to express this sum in an appropriate basis of
$H^{2}(X,\mathbb{C})$. From the modular point of view we have a
canonical basis of the Picard group of $\sBun$ which consists of line
bundles $\btheta$, $F_{1}$, \ldots $F_{5}$, where $\btheta$ is the
determinant of the cohomology of the universal bundle, and $F_{i}$ is
the line bundle corresponding to the quasi-parabolic structure on the
restriction of the universal bundle on $\sBun\times \{p_{i}\}$.
\index{terms}{quasi-parabolic!structure}
To
transport this basis to $H^{2}(X,\mathbb{C}) = \op{Pic}(X)\otimes
\mathbb{C}$ we need to also pay attention to the choice of an
embedding of $X$ in $\sBun$.  This embedding depends on the choice of
normalized universal rank two bundle  $(\mycal{V},\mathsf{r})$ on
$X\times C$.  In our setting there are natural ways to choose such a
universal bundle coming from the geometry of the Hecke
correspondence. As we saw before, the natural map $q : H \to X\times C$
has $16$ distinguished rational sections given by the divisors $Lp_{I}
\subset H$. For any $I \in \evenL$ we have a universal bundle
$({}^{I}\!\mycal{V},{}^{I}\!F_{\bullet})$ where the rank two vector bundle
is ${}^{I}\!\mycal{V} = q_{*}\mathcal{O}_{H}(Lp_{I})$ and the
quasi-parabolic structure given by line subbundles ${}^{I}\!F_{i}
\subset {}^{I}\!\mycal{V}_{|X\times \{p_{i}\}}$ coming from imposing
vanishing of the fiberwise sections of $\mathcal{O}_{H}(Lp_{I})$ at
the locus of intersection of the two components $R_{1,i}$ and
$R_{2,i}$ of the preimage $RR_{i} = q^{-1}(X\times \{ p_{i}\})$.  Once
we fix the parabolic weights, e.g. by fixing the pair $(a,b)$, we can
view the pair $({}^{I}\!\mycal{V},{}^{I}\!F_{\bullet})$ as a family of
stable rank two parabolic bundles on $(C,\Par_{C})$ and hence we get a
morphism ${}^{I}\!\jmath : X \hookrightarrow \sBun' \subset \sBun$ for
which ${}^{I}\!\jmath^{*}$ induces an isomorphism on complexified Picard
groups.

We can use any of these isomorphisms to transport the standard basis
$\btheta$, $F_{1}$, \ldots, $F_{5}$ of $\op{Pic}(\sBun)$ to a basis
${}^{I}\!\jmath^{*}\btheta$, ${}^{I}\!\jmath^{*}F_{1}$, \ldots,
${}^{I}\!\jmath^{*}F_{5}$. Note that by construction we have that
${}^{I}\!\jmath^{*}\btheta = \det \left( Rq_{*} {}^{I}\!\mycal{V}
\right)^{-1}$, and ${}^{I}\!\jmath^{*}F_{i} = {}^{I}\!F_{i}$ for all $i =
1, \ldots, 5$. To be specific we will choose the normalized universal
bundle corresponding to $\varnothing \in \evenL$ and will write
\[
\bvartheta := {}^{\varnothing}\!\jmath^{*}\btheta, \quad \bFF_{1} :=
           {}^{\varnothing}\! F_{1}, \quad \ldots, \quad \bFF_{5} :=
           {}^{\varnothing}\! F_{5}
\]
for the corresponding basis of $H^{2}(X,\mathbb{C})$. One can express
the elements of this basis as combinations of lines directly from the
defnition. A straightforward computation gives
\[
\left| \;
\begin{aligned}
  \bvartheta & = - L_{\varnothing}, \\
  \bFF_{1} & = L_{\varnothing} - L_{1^{c}}, \\
   & \cdots \\
  \bFF_{5} & = L_{\varnothing} - L_{5^{c}},
\end{aligned}
\right.
\]
where $i^{c} \in \evenL$ denotes the complement of $\{ i \}$ in $\{ 1,
2, \ldots, 5 \}$. Using these expressions and the relations among the
lines in $H^{2}(X,\mathbb{C})$ we can compute the coordinates of all
lines in the given basis:
\[
\left| \;
\begin{aligned}
  L_{\varnothing} & = - \bvartheta \\
  L_{i^{c}} & = - \bvartheta -
  \bFF_{i}, \qquad i = 1, \ldots, 5 \\
  L_{\{ij\}} & = - \bvartheta +
  \bFF_{i} + \bFF_{j} - \frac{1}{2} \sum_{k = 1}^{5} \bFF_{k},
  \qquad
  i\neq j, \ i, j \in \{1, \ldots, 5\}.
\end{aligned}
\right.
\]
Write $\underline{\bvartheta}, \ubf_{1}$, $\ubf_{2}$, \ldots,
$\ubf_{5}$ for the images of $\bvartheta$, $\bFF_{1}$,
$\bFF_{2}$, \ldots,
$\bFF_{5}$ in $H^{2}(X,\mathbb{C})/\mathbb{C}K_{X}$.  Then \linebreak 
$\ubf_{1}$, $\ubf_{2}$, \ldots,
$\ubf_{5}$ is a basis of $H^{2}(X,\mathbb{C})/\mathbb{C}
K_{X}$ and our question becomes the question of computing the matrix
of the reduced
Okamoto map
\[
\xymatrix@R-4pc@M+1pc{ \underline{\Aom} : \hspace{-3pc} &
  \lat_{C}\otimes \mathbb{C} \ar[r] &
  H^{2}(X,\mathbb{C})/\mathbb{C}K_{X} \\ & a \ar[r] & Ld \ \mod
  \ \mathbb{C}K_{X} }
\]
using the basis $(\delta_{p_{1}}, \ldots, \delta_{p_{5}})$ for the
source and the basis
$(\ubf_{1}, \ldots, \ubf_{5})$ for the target.

In view of \eqref{eq:dL} we have that for every $a \in \mathbb{C}^{5}$
the value of the Okamoto map is given by
\[
\underline{\Aom}(a) = 2 L_{\varnothing} -  \frac{1}{4} L
\left( \unone - \incdelta \right) a \ \mod \ \mathbb{C}K_{X}.
\]
Evaluating on a standard basis vector $\delta_{p_{i}} \in
\lat_{C}\otimes \mathbb{C} \cong
\mathbb{C}^{5}$
we then get
\[
\begin{aligned}
  \underline{\Aom}(\delta_{p_{i}}) & = 2 L_{\varnothing} - \frac{1}{4} L\left(
    \unone -
    2\incdelta \right) \delta_{p_{i}}\ \mod \ \mathbb{C}K_{X} \\
  & =  2 L_{\varnothing} - \frac{1}{4} \sum_{I \in \evenL} \left( 1 -
    2\incdelta_{Ii} \right) L_{I}\ \mod \ \mathbb{C}K_{X} \\
    & =  2 L_{\varnothing} - \frac{1}{4} \left(L_{\varnothing} + L_{i^{c}}
    -  \sum_{j \neq i
    } L_{j^{c}} -  \sum_{j \neq i} L_{\{ij\}} + \sum_{\substack{k <
        l\\ k, l \neq i}}
    L_{\{kl\}} \right)   \ \mod \ \mathbb{C}K_{X}  \\
    & =  2 L_{\varnothing} - \frac{1}{2} \sum_{j\neq i}
    \left(L_{j^{c}} + L_{\{ij\}} \right)
  \ \mod \
  \mathbb{C}K_{X} \\
  & = -2\bvartheta - 2\left(-2\bvartheta + \bFF_{i} -
  \frac{1}{2}\sum_{k =
    1}^{5}\bFF_{k}\right) \ \mod \
  \mathbb{C}K_{X}  \\
  & = 2\bvartheta - 2\bFF_{i} + \sum_{k =
    1}^{5}\bFF_{k} \ \mod \
  \mathbb{C}K_{X}
\end{aligned}
\]
Taking into account that
\[
  \begin{aligned}
    -K_{X} & = \frac{3}{2}L_{\varnothing} +
    \frac{1}{2} \sum_{k = 1}^{5} L_{k^{c}}  \\
    & = -\frac{3}{2}\bvartheta + \frac{1}{2} \sum_{k =
      1}^{5}\left(-\bvartheta -
      \bFF_{k}\right) \\
    & = -4\bvartheta - \frac{1}{2}\sum_{k=1}^{5} \bFF_{k},
  \end{aligned}
\]
we see that
$\underline{\bvartheta} = -\frac{1}{8}\sum_{k =1}^{5}\ubf_{k}$ in
$H^{2}(X,\mathbb{C})/\mathbb{C}K_{X}$, and so for every
$i = 1, \ldots, 5$ we have
\[
\underline{\Aom}(\delta_{p_{i}}) = -2 \ubf_{i} + \frac{3}{4}
\sum_{j =1}^{5}\ubf_{j}.
%\sum_{j \neq i} \ubf_{j}.
\]
In other words the matrix of the reduced Okamoto map 
$\underline{\Aom} : \lat_{C}\otimes \mathbb{C} \to
H^{2}(X,\mathbb{C})/\mathbb{C}K_{X}$ in the bases
$(\delta_{p_{1}},\ldots, \delta_{p_{5}})$ and $(\ubf_{1},\ldots,\ubf_{5})$ is 
\[
\boxed{
\frac{1}{4}\begin{pmatrix}
  -5 & 3 & 3 & 3 & 3 \\
  3 & -5 & 3 & 3 & 3 \\
  3 & 3 & -5 & 3 & 3 \\
  3 & 3 & 3 & -5 & 3 \\
  3 & 3 & 3 & 3 & -5
\end{pmatrix}.
}
\]

\section{Summary} \label{sec:summary}

For the convenience of the reader we conclude our case study of the
tamely ramified parabolic geometric Langlands correspondence by
summarizing again the main results and giving a brief road map to the
proofs.

\noindent
In our setup we start with a flat tame quasi-parabolic $\lan{G} =
GL_{2}(\mathbb{C})$ or $SL_{2}(\mathbb{C})$-bundle $\mathbb{V}$ on
$\mathbb{P}^1$, with quasi-parabolic structure at $5$ distinct points
$\{p_1,\dots, p_5\}$. Additionally we assume that $\mathbb{V}$ is a
general point of the moduli stack $\lan{\gLoc}$ of rank two parabolic
flat bundles, where genericity means:
\begin{itemize}
\item The residue,
$\res_{p_i}(\nabla)$, is assumed to be  regular semisimple, purely
  imaginary, and non-resonant, i.e. its eigenvalues, which by assumption are
  real, do not differ by an integer.
\item Let $\blambda = (\lambda_{1},\ldots,\lambda_{5})$ be such that for
  each $i$ the coordinate $\lambda_{i}$ is one of the two eigenvalues
  of $\res_{p_i}(\nabla)$ (so $\blambda$ is real).  Endow $\mathbb{V}$
  with a parabolic structure by assigning parabolic weight =
  $-\blambda$ to the corresponding eigenvector.  NAHT converts this to
  a parabolic Higgs bundle $\mathbb{E} = (E,\theta,{}^{E}\bF,
  -\blambda)$: the weights don't change, while the residue of the
  Higgs field $\theta$ at each $p_i$ is nilpotent.  We assume that the
  spectral curve for $(E,\theta)$ is smooth.  \index{terms}{purely
    imaginary condition} \index{terms}{genericity condition} This is
  an open condition on $\mathbb{E}$ and hence on $\mathbb{V}$.
\end{itemize}

\index{terms}{parabolic!structure}
\index{terms}{quasi-parabolic!structure}

\

\noindent
Using this notion of genericity our main result is the following

\

\noindent
    {\bfseries Theorem \ref{thm:MAIN}.} \
    {\em Parabolic Hecke eigensheaves exist for
  generic, purely imaginary flat tame $\lan{G} =GL_{2}(\mathbb{C})$ or
  $SL_{2}(\mathbb{C})$-bundles on $\mathbb{P}^1$ with parabolic
  structure at $5$ distinct points.

Concretely, if  $\mathbb{V}$ is a generic, purely imaginary flat tame
$\lan{G} = SL_{2}(\mathbb{C})$-bundle, then
there is a Hecke eigensheaf
$\mathfrak{c}_{\blambda}(\mathcal{O}_{\mathbb{V} })$ on the stack
$\sBun$ of quasi-parabolic $\mathbb{P}GL_{2}(\mathbb{C})$-bundles on
$\mathbb{P}^1$, with parabolic structure at $\{p_1,\dots, p_5\}$ with
eigenvalue $\mathbb{V}$.  The sheaf
$\mathfrak{c}_{\blambda}(\mathcal{O}_{\mathbb{V} })$ is a twisted
$\mathcal{D}$-module on $\sBun$ with a twist determined by $\blambda$.
Restricted to the very stable open subset of $\sBun$,
$\mathfrak{c}_{\blambda}(\mathcal{O}_{\mathbb{V} })$ is a flat bundle
of rank $4$.}

\index{notations}{c@$\mathfrak{c}$}
\index{notations}{Vbb@$\mathbb{V}$}
\index{notations}{SL2@$SL_{2}(\mathbb{C})$}
\index{terms}{quasi-parabolic}
\index{terms}{Hecke!eigensheaf}
\index{notations}{Buns@$\sBun$}
\index{notations}{lambdat@$\widetilde{\lambda}$}
\index{terms}{Dmodule@$\mathcal{D}$-module}
\index{terms}{residue!nilpotent}
\index{notations}{res@$\res$}
\index{terms}{twisted!$\mathcal{D}$-module}
\index{terms}{twisted!differential operators}
\index{terms}{GL2@$GL_{2}(\mathbb{C})$}

\

\

\noindent
For a more technical discussion of the formulation of the {\sf PGLC}
in the language of TDOs, see Appendix~\ref{app-ram.glc}. Also, in
Appendix~\ref{app-ram.glc} we review how, following our program in
\cite{dp-jdg}, the tamely ramified geometric Langlands conjecture in
general and the statement of Theorem~\ref{thm:MAIN} in particular
translates into a parabolic Hecke eigensheaf problem for tamely
ramified Higgs bundles on the moduli space of bundles.
\index{terms}{TDO|see{twisted differential operators}} Solving this
parabolic Hecke eigensheaf problem for $SL_{2}(\mathbb{C})$ or
$SL_{2}(\mathbb{C})$ Higgs bundles on $\mathbb{P}^{1}$ is achieved by
a detailed analysis of the geometry of the relevant moduli spaces
which we summarize next.

Let $C = \mathbb{P}^{1}$ and let
$\Par_{C} = p_{1} + p_{2} + p_{3} + p_{4} + p_{5}$ be the sum of five
distinct points.  The moduli space of rank two parabolic bundles on
$(C,\Par_{C})$ depends on a set of numerical invariants - the degree
of the level zero bundle in the parabolic family and the set of
parabolic weights.  The collection of weights has a chamber structure
and the moduli space depends only on the chamber and not on the
particular collection of weights in that chamber. As a first
observation we have the following

\

\noindent
{\bfseries Theorem  \ref{thm:deg0=deg1} (Lemma~\ref{lem:moduli.isos},
    Theorem~\ref{theo:shapes.of.moduli}])} \ {\em
  There is a dominant chamber
of parabolic weights for  rank two parabolic bundles on
$(C,\Par_{C})$  such that for all weights in this chamber:
\begin{itemize}
\item every semistable parabolic bundle is
stable; 
\item the connected components of the
moduli space corresponding to different degrees are canonically
isomorphic to the $dP_{5}$ del Pezzo surface $X$ obtained by blowing up
the $5$ points $\{p_{i}\}_{i=1}^{5}$ on the diagonal  conic $C \subset S^{2}C
\cong \mathbb{P}^{2}$.
\end{itemize}
}

\

\noindent
Equivalently: $X$ can be described in its anticanonical model as the
intersection of two quadrics in $\mathbb{P}^{4}$. The parameter space
of the pencil of quadrics $\{Q_{c}\}_{c \in C}$ vanishing on $X$ is
naturally identified with $C$ and the divisor $\Par_{C}$ corresponds
to the locus of singular quadrics in the pencil. Using this, one checks

\

\noindent
{\bfseries Theorem \ref{thm:wobbly} (Proposition~\ref{prop:wobbly})}
\ {\em The wobbly locus in $X$ is the union of the 16 lines \linebreak
  $L_{I} \subset X \subset \mathbb{P}^{4}$ which are naturally labeled
  by $I \in \evenL$, where $\evenL$ is the collection of subsets in
  $\{1,2,3,4,5\}$ of even cardinality.  }

\

\

\noindent
In Chapter~\ref{sec:hecke} we show that from the point of view of the
anti-canonical model of $X$ the basic Hecke correspondence
parametrizing the modifications of parabolic bundles at a single point
can be compactified and resolved to the
correspondence
\begin{equation} \label{eq:parHeckdiag}
\xymatrix@R-1pc{
& H \ar[dl]_-{p} \ar[dr]^-{q} & \\
X & & X\times C
}
\end{equation}
where
\begin{itemize}
\item $H = \op{Bl}_{\coprod_{I} \widehat{L_{I}\times L_{I}}}
    \op{Bl}_{\Delta} (X\times X)$;
\item the two maps $H \to X$  correspond to
    the blow down map $H \to X\times X$ followed by the first or
    second projection; 
\item the map $H \to C$ is the resolution of the rational map $X\times
  X \dashrightarrow C$ which sends $(x,y) \in X\times X$ to the unique
  $t \in C$ such that $Q_{t} \subset \mathbb{P}^{4}$
  contains the line through the two points $x, y \in \mathbb{P}^{4}$.
\end{itemize}
Note that $H$ is smooth by construction. The general fibers of $q$ are
smooth rational curves (Hecke lines) and the general fibers of $p$ are
smooth $dP_{6}$ del Pezzo surfaces. Furthermore, as explained in
Section~\ref{ssec:pardivisors}, all spaces in the Hecke diagram
\eqref{eq:parHeckdiag}
are naturally equipped
with normal crossings parabolic divisors:
$\Par_{C} = \sum_{i=1}^{5} p_{i}$, $\Par_{X} = \sum_{I} L_{I}$,
\linebreak 
$\Par_{X\times C} = \Par_{X}\times C + X\times \Par_{C}$, and 
$\Par_{H} = p^{*}\Par_{X} + q^{*}\Par_{X\times C}$.
This geometry provides the setup  needed to formulate the
  parabolic version of the Higgs Hecke eigensheaf problem.

 The solution to this problem for $\lan{G} = GL_{2}(\mathbb{C})$ or
 $SL_{2}(\mathbb{C})$ is derived in Chapters \ref{sec:eigensheaves}
 and \ref{sec:solution} and is summarized in
 Proposition~\ref{prop:gl2.and.sl2}.  The case $G =
 \mathbb{P}GL_{2}(\mathbb{C})$ follows with obvious modifications.
 These statements are all straightforward corollaries from the
 following theorem.
 
\

\

\noindent
    {\bfseries Theorem \ref{thm:MAIN.Higgs} (Proposition~\ref{prop:numerics})}
    \ {\em
Fix a pair of parabolic weights
  on $\Par_{C}$, specified by a pair
  of real vectors $a, b \in \mathbb{R}^{5}$ that are generic
  except for the linear relation $\sum_{i = 1}^{5} (a_{i}+b_{i}) = 3$
  (which, by \eqref{eq:spectral.parch}, expresses the condition that
  $\parch_{1}(E_{\bullet}) = 0$). Also, fix a pair of integral vectors
  $n_{1}, n_{2} \in \mathbb{Z}^{5}$ such that $\sum_{i = 1}^{5}
  (n_{1,i}+ n_{2,i}) = -3$. Then
\begin{description}
\item[Hecke kernel:] The pair $(a,b)$ determines natural parabolic
  line bundles $\mycal{I}_{\bullet} =
  \mathcal{O}_{H}(\bzeta\Par_{H})_{\bullet}$ and $\mycal{I}_{\bullet}' =
  \mathcal{O}_{H}(\bzeta'\Par_{H})_{\bullet}$  on $(H,\Par_{H})$ with
  $\parch_{1}(\mycal{I}_{\bullet}) =
  \parch_{1}(\mycal{I}_{\bullet}') = 0$.
\item[Hecke eigensheaf:] For any tame strongly
  parabolic rank two Higgs bundle $(E_{\bullet},\theta)$ on
  $(C,\Par_{C})$ with weights specified by $(a,b)$
  and a smooth spectral cover $\sC$, there exists
  a natural Hecke eigensheaf on the (disconnected) moduli
  space. Explicitly, this is specified by a pair of objects  
  $(F_{\bullet},\varphi)$ and $(F_{\bullet}',\varphi')$  on the del
  Pezzo $X$, so that
  \begin{itemize}
  \item $(F_{\bullet},\varphi)$ and $(F_{\bullet}',\varphi')$ are
    stable, strongly parabolic, rank
    four Higgs bundles on $(X,\Par_{X})$ with $\parch_{1} = 0$,
    $\parch_{2} = 0$.
 \item The parabolic Hecke eigensheaf conditions
\[
\begin{aligned}
q_{*}(p^{*}(F_{\bullet},\varphi)\otimes (\mycal{I}_{\bullet},0)) & = 
p_{X}^{*}(F_{\bullet}',\varphi')\otimes p_{C}^{*}(E_{\bullet},\theta) \\
q_{*}(p^{*}(F_{\bullet}',\varphi')\otimes (\mycal{I}_{\bullet}',0)) & = 
p_{X}^{*}(F_{\bullet},\varphi)\otimes p_{C}^{*}(E_{\bullet},\theta) 
\end{aligned}
\]
hold.
\item The parabolic structures on $(F_{\bullet},\varphi)$ and
  $(F_{\bullet}',\varphi')$ each have two jumps, of rank one and three,
  on every component $L_{I}$ of $\Par_{X}$.
\item The parabolic weights of $(F_{\bullet},\varphi)$ and
  $(F_{\bullet}',\varphi')$ are specified by vectors $e, d, e', d' \in
  \mathbb{R}^{\evenL} = \mathbb{R}^{16}$ which can be taken to be
  \[
\left| \; 
    \begin{aligned}
e & = \frac{5}{8}\sigma - 3\acan +\frac{3}{8}(\unone  - 2\incdelta)
(r_{1} - r_{2}), \\[+0.5pc]
d & = \frac{1}{8}\sigma + \acan
- \frac{1}{8}(\unone - 2\incdelta)(r_{1} - r_{2}), \\[0.5pc]
e' & = \frac{5}{8}\sigma  -
3\acan - \frac{3}{8}(\unone  - 2\incdelta)
(r_{1} - r_{2}), \\[+0.5pc]
d' & = \frac{1}{8}\sigma + \acan
      + \frac{1}{8}(\unone - 2\incdelta)(r_{1} - r_{2}).
\end{aligned}
\right.
    \]
   where $\mathfrak{u} = 2\op{pr}_{\{5\}}(\delta_{L_{\varnothing}})$,
    $r_{1} = a + n_{1}$ and $r_{2} = b + n_{2}$, $\sigma \in
    \mathbb{R}^{\evenL}$ is a vector with all its coordinates equal to
    $1$, and $\unone$ and $\incdelta$ are $16\times 5$ matrices whose
    entries corresponding to $i \in \{1,2,3,4,5\}$ and $I \in \evenL$
    are given by $\unone_{Ii} = 1$ for all $I,i$, while
    $\incdelta_{Ii} = 1$ if $i \in I$ and $\delta_{Ii} = 0$ if $i
    \notin I$.
\item Both $(F_{\bullet},\varphi)$ and $(F_{\bullet}',\varphi')$
  correspond to spectral data specified on a spectral cover of $X$
  which is birational to the Hitchin fiber through
  $(E_{\bullet},\theta)$.
\end{itemize}
\end{description}   
}

\

\noindent
As explained briefly above, and in detail in
Appendix~\ref{app-ram.glc}, Theorem~\ref{thm:MAIN} follows
from Theorem~\ref{thm:MAIN.Higgs} by conjugating the assignment
$(E_{\bullet},\theta) \mapsto (F_{\bullet},\varphi)$ by two
non-abelian Hodge correspondences (one on $(C,\Par_{C})$ and one on
$(X,\Par_{X})$) and by taking a twisted Deligne-Goresky-MacPherson
extension
\index{terms}{DeligneGG@Deligne-Goresky-MacPherson!extension}
of the flat bundle corresponding to $(F_{\bullet},\varphi)$
from $X - \Par_{X}$ to the moduli stack $\sBun$.  Under the
non-abelian Hodge correspondence the assignment on parabolic
parameters $(a,b) \mapsto (e,d)$ becomes the Okamoto map which
converts the spectrum of the residue of the eigenvalue flat bundle
$\mathbb{V}$ into the twisting parameters for the eigensheaf twisted
$\mathcal{D}$-module $\mathfrak{c}\left(\mathcal{O}_{\mathbb{V}}\right)$. In
Section~\ref{ssec:okamoto} we carry out this translation and compute
the Okamoto map in our setting. \index{terms}{Okamoto!map}

\

\noindent
Curiously the non-abelian Hodge theory constraints on the modular
spectral data turn out to be sufficiently restrictive to rigidify the
parabolic eigensheaf problem. To illustrate this, in
Claim~\ref{claim:spectral.is.eigensheaf} we show that for any choice
of parabolic weights on $(X,\Par_{X})$ producing modular spectral data
whose Higgs bundle has vanishing $\parch_{1}$ and $\parch_{2}$ we can
find a companion Hecke eigenvalue (we thank the
referee for suggesting that this statement might be true):

\

\noindent
{\bfseries Claim~\ref{claim:spectral.is.eigensheaf}} \ {\em
Fix a smooth spectal cover $\pi :\sC \to C$ ramified over $\Par_{C}$,
and let $f : Y \to X$ be the corresponding modular spectral
cover. Given vectors $e, d$  in $\lat_{X}\otimes
\mathbb{R}$ that obey the non-abelian Hodge theory constraints
\eqref{chsoln}, there exist vectors $\bzeta = (lp,lq,llc,r_{1},r_{2})$,
$e', d' \in \lat_{X}\otimes
\mathbb{R}$
and vectors $a, b \in \lat_{C}\otimes \mathbb{R}$ so that $e', d'$ also obey
\eqref{chsoln}, and if
$(F_{\bullet},\varphi)$ and $(F_{\bullet}',\varphi')$ are the strongly
parabolic rank four Higgs bundles on $(X,\Par_{X})$ given by modular
spectral data $\left(Y,\mycal{L}_{\mathfrak{a}}(eE +
d(E+G))_{\bullet,\bullet}\right)$ and
$\left(Y,\mycal{L}_{\mathfrak{a}}(e'E +
d'(E+G))_{\bullet,\bullet}\right)$ with $\mathfrak{a} \in
\Pic^{0}(\sC)$, then
\[
q_{*}\left(p^{*}(F_{\bullet},\varphi)\otimes
(\mathcal{O}(\bzeta\Par_{H})_{\bullet},0)\right)  =
p_{X}^{*}(F_{\bullet}',\varphi')\otimes p_{C}^{*}(E_{\bullet},\theta),
\]
where $(E_{\bullet},\theta)$ is the strongly parabolic rank two Higgs
bundle on $(C,\Par_{C})$ correspinding to the spectral data
$(\sC,\mathfrak{a}(a\widetilde{P} +
b\widetilde{P}))_{\bullet,\bullet}$.
}

\newpage
\appendix

\Appendix{TDO and the tamely ramified GLC}
\label{app-ram.glc}

\setcounter{equation}{0}

\

\noindent
In this appendix we briefly discuss the translation of the ramified
geometric Langlands correspondence into a parabolic Hecke eigensheaf
problem for ramified Higgs bundles on the moduli space of
bundles. Since the Langlands transformation of Stokes data is not very
well understood, we focus on the case of tame ramification with generic
residues where the geometric Langlands conjecture has a clear
formulation.

\index{terms}{Geometric Langlands Conjecture!tamely ramified}
\index{terms}{Higgs!bundle!tamely ramified}
\index{terms}{Hecke!eigensheaf}

\subsection{Setup for the tamely ramified GLC}
\label{ssec:setup.glc}

\

\noindent One way to formulate the tamely ramified geometric Langlands
correspondence, perhaps the most natural, is as a family of
dualities - one for each closure of a regular symplectic leaf in the stack
$\lan{\gLoc}$ of logarithmic quasi-parabolic flat connections on
principal $\lan{G}$-bundles \cite{arinkin-lysenko,arinkin-4points}.
\index{terms}{flat!bundle!quasi-parabolic}
\index{terms}{flat!connection!logarithmic!quasi-parabolic}

Specifically $\lan{\gLoc}$ is the moduli stack of triples
$(V,\nabla,\mathsf{r})$ where $V$ is a principal $\lan{G}$ bundle
on $C$, $\nabla$ is a logarithmic flat connection on $V$, and
$\mathsf{r} \in \Gamma(\Par_{C}, V/\lan{B})$ is a Borel reduction
of the structure group of $V$ at the points of $\Par_{C}$ which is
furthermore preserved by the residue of $\nabla$.

\index{notations}{Vnabla@$(V,\nabla)$}
\index{notations}{LFlats@$\lan{\gLoc}$}
\index{notations}{C@$C$}
\index{notations}{LB@$\lan{B}$}
\index{notations}{r@$\mathsf{r}$}
\index{notations}{ParC@$\Par_{C}$}
\index{terms}{residue!of a logarithmic connection}
\index{terms}{Borel!reduction!of the structure group}
\index{terms}{Borel!subgroup}

Next, fix a closure of a regular adjoint orbit for each point in the
parabolic divisor, i.e. a map $\bLambda : \Par_{C} \to
\lan{\mathfrak{g}}/\!/\lan{G} = \lan{\mathfrak{t}}/\!/W$. 
The map $\bLambda$
labels a closure $\lan{\gLoc}_{\bLambda}$ of a symplectic leaf in
$\lan{\gLoc}$. This closure has also a modular interpretation as the
closed substack $\lan{\gLoc}_{\bLambda} \subset \lan{\gLoc}$
parametrizing quasi-parabolic flat bundles with residues in
$\bLambda$.  On the other hand, using the isomorphism
$\mathfrak{t}^{\vee} = \lan{\mathfrak{t}}$ we can interpret $\bLambda$
as an assignment of a $W$-orbit of linear functions on $\mathfrak{t}$
to each point in the parabolic divisor $\Par_{C}$.

\index{terms}{adjoint orbit!regular}
\index{terms}{symplectic leaf!in $\lan{\gLoc}$}
\index{notations}{LFlats@$\lan{gLoc}_{\bLambda}$}
\index{notations}{Lambdab@$\bLambda$}
\index{notations}{tfrak@$\mathfrak{t}$}
\index{notations}{Ltfrak@$\lan{\mathfrak{t}}$}

When $\bLambda$ is generic (i.e. has values not contained in any
reflection hyperplane) it gives rise (see Remark~\ref{rem-tdo}) to a
derived category ${}_{\bLambda}D(\sBun,\mathcal{D})$ of
$\bLambda$-twisted $\mathcal{D}$-modules on the stack of quasi-parabolic
principal $G$-bundles on $(C,\Par_{C})$. Ignoring subtleties (see
\cite{dp-Langlands,dp-jdg,ag-sing,bznp}) arising from choices of
connected components, twists by gerbes, and singular supports of
coherent sheaves, the $\bLambda$-layer of the geometric Langlands
conjecture in this case predicts the existence of an equivalence of
categories
\begin{equation} \label{eq-clambda} \mathfrak{c}_{\bLambda} :
D(\lan{\gLoc}_{\bLambda},\mathcal{O}) \stackrel{\cong}{\to}
{}_{\bLambda}D(\sBun,\mathcal{D})
\end{equation} which as usual intertwines tensorization and Hecke
operators (see \cite{dp-Langlands,dp-jdg}).

\index{terms}{Dmodule@$\mathcal{D}$-module!$\bLambda$-twisted}
\index{notations}{Buns@$\sBun$}
\index{notations}{cbiglambda@$\mathfrak{c}_{\bLambda}$}

\

\begin{rem} \label{rem-tdo} Suppose $(V,\nabla,\mathsf{r}) \in
\lan{\gLoc}$. The residue of the connection $\nabla$ at the parabolic
divisor $\Par_{C}$ is a section $\res \nabla \in
\Gamma(\Par_{C},\op{ad}(V))$. Since $\op{ad}(V)//\lan{G} =
V\times_{\op{ad}} \lan{\mathfrak{g}}/\!/\lan{G}$, the bundle
$\op{ad}(V)/\!/\lan{G}$ is naturally trivialized:
$\op{ad}(V)/\!/\lan{G} \cong (\lan{\mathfrak{g}}\otimes
\mathcal{O}_{C})/\!/\lan{G}.$ 
Composing with this trivialization we get a natural map
\[ \xymatrix{
\Par_{C} \ar[r]^-{\res \nabla} \ar@/_1pc/[rrd]_-{\underline{\res}\,
\nabla} & \op{ad}(V)_{|\Par_{C}} \ar[r] & \op{ad}(V)_{|\Par_{C}}
/\!/\lan{G} \ar[d]^-{\cong} \\ & &
(\lan{\mathfrak{g}}/\!/\lan{G})\times \Par_{C}.  }
\]
By definition, the closed substack $\lan{\gLoc}_{\bLambda} \subset
\lan{\gLoc}$ parametrizes $(V,\nabla,\mathsf{r}) \in \lan{\gLoc}$
with $\underline{\res}\, \nabla = \bLambda$.

\index{notations}{res@$\res$}
\index{notations}{LG@$\lan{G}$}
\index{notations}{LFlats@$\lan{\gLoc}$}

  On the other hand, for any $(V,\nabla,\mathsf{r}) \in
\lan{\gLoc}$, the quasi-parabolic structure $\mathsf{r}$ gives a
reduction of the structure group of $V_{|\Par_{C}}$ to $\lan{B}$.
\index{terms}{quasi-parabolic!structure}
Write $V^{\mathsf{r}} \subset V_{|\Par_{C}}$ for the corresponding
$\lan{B}$-bundle. The compatibility of $\nabla$ and $\mathsf{r}$
implies that $\res \nabla \in \Gamma(\Par_{C}, \op{ad}(V^{\text{{\sf
r}}}))$ is a section in the adjoint bundle of the $\lan{B}$-bundle
$V^{\mathsf{r}}$. Identifying $\lan{T}$ with the quotient
$\lan{B}/[\lan{B},\lan{B}]$, we get a linear map $\lan{\mathfrak{b}}
\to \lan{\mathfrak{t}}$ which intertwines the adjoint action of
$\lan{B}$ on $\lan{\mathfrak{b}}$ with the trivial action of $\lan{B}$
on $\lan{\mathfrak{t}}$. (In other words, the map is $\lan{B}$-invariant.)
Thus we get a natural map of Lie algebra
bundles
 \[ \op{ad}(V^{\mathsf{r}}) \to V^{\text{{\sf
r}}}\times_{\op{trivial}} \lan{\mathfrak{t}} = \lan{\mathfrak{t}}
\otimes \mathcal{O}_{\Par_{C}}
\] on $\Par_{C}$.  Composing with this map we get a different shadow
of the residue:
\[ \xymatrix@1{
\Par_{C} \ar[r]^-{\res \nabla} \ar@/_1pc/[rd]_-{\res^{\mathsf{r}}
\nabla} & \op{ad}(V^{\mathsf{r}}) \ar[d] \\ &
\lan{\mathfrak{t}}\times \Par_{C}.  }
\]
Note that by construction, the section $\res^{\mathsf{r}} \,
\nabla \in \Gamma(\Par_{C}, \lan{{\mathfrak t}}\otimes \mathcal{O})$
lifts the section \linebreak $\underline{\res}\, \nabla \in
\Gamma(\Par_{C}, (\lan{{\mathfrak t}}\otimes \mathcal{O})/\!/W)$. In
particular, if we choose any lift $\blambda \in
\Gamma(\Par_{C}, \lan{{\mathfrak t}}\otimes \mathcal{O})$ of
$\bLambda$, we will get a connected component
$\lan{\gLoc}_{\blambda}$ of $\lan{\gLoc}_{\bLambda}$
parametrizing all triples $(V,\nabla,\mathsf{r})$ with
$\res^{\mathsf{r}} \, \nabla = \blambda$. The symplectic leaf
$\lan{\gLoc}_{\bLambda}$ thus is the disjoint union of the components
$\lan{\gLoc}_{\blambda}$ where $\blambda$ ranges over
all lifts of $\bLambda$.

On the other hand we consider the stack $\sBun$ of quasi-parabolic
$G$-bundles, i.e. $G$-bundles on $C$ equipped with a reduction of the
structure group to $B$ over the parabolic divisor $\Par_{C} \subset
C$. In particualr, applying this reduction to the universal bundle
$\mycal{V} \to \sBun\times C$, we obtain a well defined $B$-bundle
$\mycal{V}^{\mathsf{r}}$ over $\sBun\times \Par_{C}$.  The quotient
$\mycal{V}^{\mathsf{r}}/[B,B]$ is then a well defined $T$-bundle on
$\sBun\times \Par_{C}$ or equivalently a collection $\{
\mycal{T}_{x}\}_{x \in \Par_{C}}$ of $T$-bundles on $\sBun$, where
\[ \mycal{T}_{x} = \mycal{V}^{\mathsf{r}}/[B,B]_{\left|
\sBun\times\{ x\}\right.}.
\]
The complexified character lattice of $T$ is naturally
identified with the dual space $\mathfrak{t}^{\vee}$ of the Cartan
algebra $\mathfrak{t} = \op{Lie}(T)$. Thus we get a natural linear map
\[
\Xi : \Gamma(\Par_{C},\mathfrak{t}^{\vee}\otimes \mathcal{O}) \to
\Pic(\sBun)\otimes \mathbb{C}, \quad c \mapsto \Xi_{c},
\]
where
\[
\Xi_{c} = \bigotimes_{x \in \Par_{C}}
c(x)\left(\mycal{T}_{x}\right) \; \in
\Pic(\sBun)\otimes \mathbb{C}
\]
and $c(x)\left(\mycal{T}_{x}\right)$ is the element in the
complexified Picard of $\sBun$ which is associated to the $T$-bundle
$\mycal{T}_{x}$ via the complexified character
$c(x) \in \mathfrak{t}^{\vee}$. The map $\Xi$ is known to be injective
\cite{ls}.  In fact when $G$ is simple it is well known \cite{ls} that
the Picard group of $\sBun$ is finitely generated and that the finite
dimensional vector space $\Pic(\sBun)\otimes \mathbb{C}$ decomposes as
$\Pic(\sBun)\otimes \mathbb{C} = \mathbb{C}\cdot \btheta \oplus
\Gamma(\Par_{C},\mathfrak{t}^{\vee}\otimes \mathcal{O})$.  Here as
usual $\btheta$ denotes the non-abelian theta line bundle,
i.e. $\btheta := \det(Rp_{\sBun *} \op{ad}(\mycal{V}))^{\vee}$, where
$p_{\sBun} : \sBun\times C \to \sBun$ is the natural projection.
\index{terms}{character!complexified}
\index{terms}{Picard!group!complexified}
\index{notations}{Vunivs@$\mycal{V}$}
\index{notations}{Lmap@$\mycal{L}$}
\index{notations}{Lmapc@$\mycal{L}_{c}$} \index{terms}{theta!line
  bundle!nonabelian} \index{notations}{thetab@$\btheta$}
\index{notations}{Buns@$\sBun$}
\index{notations}{PicBunotC@$\Pic(\sBun)\otimes\mathbb{C}$}

As usual, each element $\xi \in \Pic(\sBun)\otimes
\mathbb{C}$ gives rise \cite{kashiwara-tdo},
\cite[Section~2]{bb-jantzen}, \cite[Section~I.9]{bjork-analyticD} to a
sheaf of twisted differential operators (TDO)
$\mathcal{D}_{\xi}$ and to the associated derived
category $D(\sBun,\mathcal{D}_{\xi})$ of twisted
$\mathcal{D}$-modules.

\index{terms}{twisted differential operators}
\index{notations}{DM@${\mathcal{D}_{\Xi}}$}
\index{terms}{Dmodule@$\mathcal{D}$-module!twisted}
\index{terms}{Dmodule@$\mathcal{D}$-module!$\bLambda$-twisted}

But the Langlands duality isomorphism
$\lan{\mathfrak{t}} = \mathfrak{t}^{\vee}$ allows us to view each
$\blambda : \Par_{C} \to \lan{\mathfrak{t}}$ as an element in
$\Gamma(\Par_{C},\mathfrak{t}^{\vee}\otimes \mathcal{O})$.  The
geometric Langlands conjecture now predicts the existence of the so
called {\bfseries\em Okamoto map} which is a natural map
\index{terms}{Okamoto!map} \index{notations}{Ok@$\Aom$}
\[
\Aom : \Gamma(\Par_{C},\mathfrak{t}^{\vee}\otimes \mathcal{O}) \to
\mathbb{C}\cdot \btheta\oplus
\Gamma(\Par_{C},\mathfrak{t}^{\vee}\otimes \mathcal{O}) =
\Pic(\sBun)\otimes \mathbb{C}
\]
with the property that for each $\blambda \in
\Gamma(\Par_{C},\lan{\mathfrak{t}}\otimes\mathcal{O})$ there exists an
equivalence of categories
\begin{equation} \mathfrak{c}_{\blambda} :
D(\lan{\gLoc}_{\blambda},\mathcal{O}) \stackrel{\cong}{\to}
D(\sBun,\mathcal{D}_{\Aom(\blambda)}).   \label{eq-clambdatilde}
\end{equation}
These can be put together to give the Langlands correspondence
\eqref{eq-clambda}. Indeed, since $\lan{\gLoc}_{\bLambda}$ is a disjoint
union of connected components $\lan{\gLoc}_{\blambda}$, it follows
that the derived category $D(\lan{\gLoc}_{\bLambda},\mathcal{O})$ decomposes
into a direct sum
\[
D(\lan{\gLoc}_{\bLambda},\mathcal{O}) =
\bigsqcup_{\substack{\blambda \in
\Gamma(\Par_{C},\lan{\mathfrak{t}}\otimes \mathcal{O}) \\ \text{ a
lift of } \bLambda}} D(\lan{\gLoc}_{\blambda},\mathcal{O}).
\]
Similarly we can define the category
${}_{\bLambda}D(\sBun,\mathcal{D})$ of $\Aom(\blambda)$-twisted
$\mathcal{D}$-modules on $\sBun$ as the direct sum
\[ {}_{\bLambda}D(\sBun,\mathcal{D}) =
\bigsqcup_{\substack{\blambda \in
\Gamma(\Par_{C},\lan{\mathfrak{t}}\otimes \mathcal{O}) \\ \text{ a
lift of } \bLambda}} D(\sBun,\mathcal{D}_{\Aom(\blambda)}).
\]
Now the conjectural Langlands correspondence
$\mathfrak{c}_{\bLambda}$ is just the assembly $\sqcup\,
\mathfrak{c}_{\blambda}$ of the Langlands correspondences for the
individual components.
\index{notations}{clambda@$\mathfrak{c}_{\blambda}$}
\index{notations}{cbiglambda@$\mathfrak{c}_{\bLambda}$}
\index{terms}{Langlands!functor}

\end{rem}

\subsection{The non-abelian Hodge theory
  approach} \label{ssec:nah.strategy}

The images of structure sheaves of points under
$\mathfrak{c}_{\blambda}$ are $\blambda$-twisted $\mathcal{D}$-modules
which satisfy the Hecke eigensheaf property.  As explained in
\cite{dp-jdg}, our strategy for constructing such Hecke eigensheaves
is to conjugate a Fourier-Mukai transform with two non-abelian Hodge
correspondences: one for semistable parabolic flat logarithmic
$\lan{G}$ connections on $(C,\Par_{C})$ and the other for semistable
parabolic vector bundles on $(\Bun,\Par_{\Bun})$ equipped with flat
logarithmic connections. Here $\Bun$ is the projective moduli space of
semistable parabolic $G$ bundles on $(C,\Par_{C})$, and $\Par_{\Bun}$
is the divisor in $\Bun$ parametrizing wobbly bundles.

To carry this out we must choose parabolic weights for quasi-parabolic
$\lan{G}$ bundles on $(C,\Par_{C})$ and compatible parabolic weights
for the quasi-parabolic vector bundles on $(\Bun,\Par_{\Bun})$. The
weights have to be chosen so that the respective bundles have zero
parabolic Chern classes and so can be inserted in the non-abelian
Hodge correspondence.

\subsubsection{An assumption} \label{sssec:assumption}
In order to proceed, we make an assumption on our objects, namely that
they fit in the ``purely imaginary'' case of the non-abelian Hodge
correspondence.
\index{terms}{purely imaginary condition}
This assumption is needed in order for the results of
\cite{dps} to apply. In terms of the flat $\lan{G}$ local systems, the
assumption is that the eigenvalues of the local monodromies are in the
unit circle and that the associated local systems are not filtered,
i.e. they have trivial parabolic structure along $\Par_{C}$.
 According to the conversion table from
\cite{carlos-nc} or
 \cite{dps} (we use the sign convention of \cite{dps} here), this means that 
the residues of the Higgs field of the corresponding Higgs bundle are nilpotent.
\index{terms}{conversion table!of the ramified nah}

In terms of our original bundle with flat $\lan{G}$ connection, the
condition is:

\begin{itemize}
\item The map $\blambda : \Par_{C} \to \lan{\mathfrak{t}}$ takes
  values in the subspace
  $\lan{\mathfrak{t}}^{\op{cpt}}_{\mathbb{R}}$, where
    $\lan{\mathfrak{t}}^{\op{cpt}}_{\mathbb{R}}$ denotes the real form
      of $\lan{\mathfrak{t}}$ corresponding to the compact real form
      of $\lan{G}$.
    \item The parabolic weights for the quasi-parabolic flat
      $\lan{G}$-bundles, which in general are specified by a map
      \cite[Definition~2.1]{teleman.woodward}, \cite[Section~6.1 and
      Remark~1.0.4(2)]{balaji.seshadri}
      $ \Par_{C} \to
      \lan{\mathfrak{t}}^{\op{cpt}}_{\mathbb{R}}/\!/W^{\mathsf{aff}}
      \subset \lan{\mathfrak{t}}^{\op{cpt}}_{\mathbb{R}} =
      \chr(\lan{T})\otimes \mathbb{R}$, must be equal  to
      $-\blambda$.  Here
        the group $\lan{G}$ is assumed to be simple and the quotient
        $\lan{\mathfrak{t}}^{\op{cpt}}_{\mathbb{R}}/\!/W^{\mathsf{aff}}$
        is embedded in $\lan{\mathfrak{t}}^{\op{cpt}}_{\mathbb{R}}$
        via its natural identification with the (real) Weyl alcove
        \[
        \left\{ x \in \lan{\mathfrak{t}}^{\op{cpt}}_{\mathbb{R}} \;
          \left| \;  \alpha_{\op{max}}(x) \leq 1, \ \alpha(x) > 0,
            \text{ for all $\lan{B}$-positive roots } \ \alpha \right.\right\}.
         \]
    \end{itemize}
    \index{terms}{quasi-parabolic!flat bundle}
    \index{notations}{charLT@$\chr(\lan{T})$}
    \index{terms}{Weyl!alcove}
    \index{terms}{Weyl!group}
    \index{terms}{Weyl!group!affine}
    \index{notations}{Ltcpt@$\lan{\mathfrak{t}}^{\op{cpt}}_{\mathbb{R}}$}

    \

 \noindent
 In terms of the table in \cite[page 720]{carlos-nc}, where the three
 real parameters are denoted $\alpha, b, c$, the purely imaginary
 condition is $b=c=0$, so $(E,\theta)$ has weights $\alpha$ and
 eigenvalues 0; $(V,\nabla)$ has weights $\alpha$ and eigenvalues
 $-\alpha$; and $L$ has weights 0 and eigenvalues
 $\exp(-2\pi\sqrt{-1}\alpha)$. The second table in \cite[page 12]{dps}
 specifies our condition. In terms of the parameters $r_i, \beta_i$ in
 the more general first table there, the condition is that $r_i=0$.

 Now fix
 $\blambda : \Par_{C} \to \lan{\mathfrak{t}}^{\op{cpt}}_{\mathbb{R}}
 \subset \lan{\mathfrak{t}}$.  Under Simpson's non-abelian Hodge
 theorem, \cite{carlos-nc} tame semistable parabolic flat $\lan{G}$
 bundles $(V,\nabla,\mathsf{r})$ on $(C,\Par_{C})$ that have parabolic
 weights $\blambda$ and conjugacy classes of residues
 $\res^{\mathsf{r}}\, \nabla = - \blambda \in
 \Gamma(\Par_{C},(\lan{\mathfrak{t}}\otimes \mathcal{O}))$ correspond
 to tame semistable parabolic $\lan{G}$ Higgs bundles
 $(E,\theta,\mathsf{r})$ with parabolic weights $\blambda$ and
 nilpotent residues, i.e.
 $\res^{\mathsf{r}} \, \theta = 0 \in
 \Gamma(\Par_{C},(\lan{\mathfrak{t}}\otimes \mathcal{O}))$.  Using
 this choice of $\blambda$ and following the strategy from
 \cite{dp-jdg} we can now formulate the parabolic Hecke eigensheaf
 problem that is equivalent to the existence of
 $\mathfrak{c}_{-\blambda}$ in \eqref{eq-clambdatilde}.  As in
 \cite{dp-Langlands,dp-jdg} we fix once and for all a $G$-invariant
 bilinear pairing on $\mathfrak{g}$ and we use it to identify the
 Cartan subalgebras of $G$ and $\lan{G}$:
 $\lan{\mathfrak{t}} = \mathfrak{t}^{\vee} \ \widetilde{\to} \
 \mathfrak{t}$. Using this identification we will also view $\blambda$
 as a map $\blambda : \Par_{C} \to \mathfrak{t}$.
 \index{terms}{parabolic!weights}
 \index{notations}{lan@$\mathsf{lan}$}

\

\noindent
Next consider 
\begin{description}
\item[$\lan{\Higgs}_{0}(\blambda)$] - the moduli space of
  semistable $\lan{G}$ parabolic Higgs bundles on $(C,\Par_{C})$ with
  parabolic weights $\blambda : \Par_{C} \to
  \lan{\mathfrak{t}}^{\op{cpt}}_{\mathbb{R}}\subset
  \lan{\mathfrak{t}}$
  and nilpotent residues.
\item[$X$] - an appropriate resolution of the moduli space
  $\Bun(\blambda)$ of semistable parabolic
  $G$ bundles on $(C,\Par_{C})$ with parabolic weights
  $\blambda : \Par_{C} \to
  \mathfrak{t}^{\op{cpt}}_{\mathbb{R}}\subset \mathfrak{t}$.
\end{description}
\index{terms}{Mochizuki's!theorem}

\

\noindent
Here ``appropriate'' means that $X$ should be chosen so that (see
\cite{dp-jdg}) the hypotheses of Mochizuki's non-abelian Hodge theorem
hold. Concretely this means that $X$ should be a proper
Deligne-Mumford stack whose coarse moduli space maps birationally onto
the projective variety $\Bun(\blambda)$ and
$X$ contains a closed substack $Z$ of codimension at least three so
that $X$ is smooth away from $Z$, and the preimage $\Par_{X} \subset
X$ of the wobbly locus in $\Bun(\blambda)$
is a (possibly non-reduced) normal croosing divisor away from $Z$.

\index{notations}{Bunlanlambda@$\Bun(-\mathsf{lan}(\blambda))$}
\index{notations}{ParX@$\Par_{X}$}
\index{terms}{wobbly!locus}
\index{notations}{LHiggs0lambda@$\lan{\Higgs}_{0}(\blambda)$}
\index{notations}{Higgs0lanlambda@$\Higgs_{0}(\blambda)$}

By definition the moduli space $\lan{\Higgs}_{0}(\blambda)$ is the
closure of a symplectic leaf in the Poisson space
$\lan{\Higgs}(\blambda)$ of all semistable tame parabolic Higgs
bundles (with no constrains on the residues). Similarly the analogous
moduli space $\Higgs_{0}(\blambda)$ of semistable tame parabolic $G$
Higgs bundles with nilpotent residues is the closure of a symplectic
leaf. The reasoning in \cite{dp-Langlands} generalizes to show that
the spaces $\lan{\Higgs}_{0}(\blambda)$ and $\Higgs_{0}(\blambda)$ are
dual integrable systems, i.e. their generic Hitchin fibers are dual
abelian varieties.  The strategy of \cite{dp-jdg} is to use this
duality to construct the $\Aom(-\blambda)$-twisted
$\mathcal{D}$-modules that correspond to structure sheaves of points
under $\mathfrak{c}_{-\blambda}$. \index{terms}{Hitchin!fiber}
\index{notations}{Ok@$\Aom$}
    
\subsubsection{The parabolic program}\label{program}
The construction is done in several steps, reviewed here from
\cite[section 3.2 (page 101)]{dp-jdg}, with a slight twist explained
at the end of the present section:
 \begin{list}{{\bfseries (\roman{num})}}{\usecounter{num}} 
 \item Start with a stable tame parabolic flat $\lan{G}$ bundle
   $\mathbb{V} = (V,\nabla,\mathsf{r})$ on $(C,\Par_{C})$ with
   parabolic weights $\blambda$ and
   $\res^{\mathsf{r}}\, \nabla = - \blambda$, and a vanishing
   parabolic first Chern class.
   \index{terms}{parabolic!first Chern class}
 \item Let
   $\mathbb{E} = (E,\theta,\mathsf{r}) =
   \mathsf{nah}_{C,\Par_{C}}\mathbb{V}$ be the stable tame parabolic
   Higgs $\lan{G}$ bundle corresponding to $\mathbb{V}$ by the
   non-abelian Hodge theorem \cite{carlos-nc,mochizuki-kh2}. Note that
   since by {\bfseries (i)} the flat bundle $\mathbb{V}$ satisfies
   Mochizuki's purely imaginary condition, the non-abelian Hodge
   theory conversion table from \cite[page 720]{carlos-nc} implies
   that $\mathbb{E}$ has parabolic weights $\blambda$ and nilpotent
   residues.
   \index{terms}{purely imaginary condition}
   \index{terms}{conversion table!of the ramified nah}

   We will furthermore assume that $\mathbb{V}$ is general in the sense of
   \cite{dp-Langlands}.  \index{terms}{parabolic!weights}
   \index{notations}{nahpar@$\mathsf{nah}_{C,\Par_{C}}$}
   \index{terms}{genericity condition}
   Concretely
   this means that under the Hitchin map \index{terms}{Hitchin!map}
   $\lan{\boldsymbol{h}} : \lan{\Higgs}_{0}(\blambda) \to
   \boldsymbol{B}$ the point $\mathbb{E}$ maps to a point
   $\lan{\boldsymbol{h}}(\mathbb{E}) \in \boldsymbol{B}$ which
   corresponds to a smooth cameral cover $\sC \to C$ which has simple
   Galois ramification. In particular $\mathbb{E}$ belongs to a smooth
   Hitchin fiber in $\lan{\Higgs}_{0}(\blambda)$.
 \item Use the Fourier-Mukai transform along the Hitchin fibers to
   \index{terms}{Hitchin!fiber} \index{terms}{sheaf!skyscraper}
   convert the skyscraper sheaf $\mathcal{O}_{\mathbb{E}}$ on
   $\lan{\Higgs}_{0}(\blambda)$ to a coherent complex
   \[
     \op{FM}(\mathcal{O}_{\mathbb{E}}) \in
       D^{b}_{\op{coh}}(\Higgs_{0}(\blambda))). 
     \]
     By the genericity assumption on $\mathbb{V}$ and $\mathbb{E}$,
     this coherent complex is just a translation invariant
     line bundle on the Prym
     variety $P$ for the cameral cover $\sC$ over $C$ which
     under our assumptions can be identified with
     the Hitchin fiber
     \index{terms}{Hitchin!map}
     \index{terms}{Hitchin!fiber}
     \index{terms}{Fourier-Mukai!transform}
     \index{terms}{Prym!variety}
     \index{notations}{hmap@$\boldsymbol{h}$}
     \index{notations}{Lhmap@$\lan{\boldsymbol{h}}$}
     \index{terms}{genericity condition}
     $\boldsymbol{h}^{-1}(\sC) \subset
     \Higgs_{0}(\blambda)$. In other words,
     if we write
     $i_{P} : P \hookrightarrow
     \Higgs_{0}(\blambda)$ for the
     inclusion map, we will have $\op{FM}(\mathcal{O}_{\mathbb{E}}) =
     i_{P*} \mycal{L}_{\mathbb{E}}$ for some translation invariant line bundle
     $\mycal{L}_{\mathbb{E}} \in \Pic(P)$ on the abelian variety
     $P$.
   \item The cotangent bundle to the smooth locus of
     $\Bun(\blambda)$ embeds as a
     Zariski open subset in
     $\Higgs_{0}(\blambda)$. This allows us
     to treat the pair $(P,\mycal{L}_{\mathbb{E}})$ as spectral data for a
     meromorphic Higgs bundle on
     $\Bun(\blambda)$, and use it to
     construct a tame parabolic Higgs bundle on $X$ as follows.

 Consider the 'forgetting of the Higgs field'  map
\[ \xymatrix@R-3pc@C+1pc@M+1pc{
   \Higgs_{0}(\blambda)
   \ar@{-->}[r]^-{\boldsymbol{\pi}} &
   \Bun(\blambda) \\
   (E,\theta,\mathsf{r}) \ar[r] & (E,\mathsf{r})
 }
\]
This is a rational map which is not defined at points which are stable as
Higgs bundles but whose underlying bundles are unstable.

Composing this map with the inclusion $i_{P}$ gives a quasi-finite
rational map
\[
\boldsymbol{\pi}\circ i_{P} : P \dashrightarrow
\Bun(\blambda)
\]
which we expect to be proper over the open set of very stable bundles
(see \cite{pp-n} for a proof in the non-parabolic case). Choose a
birational modification
\begin{equation}
  \bmu : Y \to P \label{eq:modify.Y}
\end{equation}
of $P$ so that the rational
map $\boldsymbol{\pi}\circ i_{P}$ becomes a generically finite
morphism $f : Y \to X$ with $\Par_{Y} := f^{*}\Par_{X}$ being a normal
crossing divisor away from codimension three,  and such that
the natural rational map $P \dashrightarrow T^{\vee}_{\Bun(\blambda)}$
factors through a morphism  $\alpha : Y \to
T^{\vee}_{X}\left(\log \Par_{X}\right)$. Thus the birational
modification $Y$ of a Hitchin fiber, \index{terms}{Hitchin!fiber}
together with the map $\alpha$, is a spectral cover of the resolution
$X$ of the moduli space $\Bun(\blambda)$. We call it the
\textit{\bfseries modular spectral cover}.  The fact that the input
parabolic Higgs bundle $\mathbb{E}$ on $C$ had nilpotent residues at
$\Par_{C}$ now translates into the statement that the map $\alpha$ and
the map $f$ align in the sense that $\alpha$ is actually a section of
$f^{*}\left(\Omega^{1}_{X}(\log \Par_{X})\right) \cap \Omega^{1}_{Y}
\subset \Omega^{1}_{Y}(\log f^{*}\Par_{X})$. In other words, viewed as
a logarithmic form on $Y$, $\alpha$ has zero residues, i.e. is
holomorphic. (In our setting we verify all these properties in detail
in chapter \ref{sec:modular}.)

In particular for
any line bundle $M \in \Pic(Y)$ the pushforward
$(f_{*}M,f_{*}(\alpha\otimes -))$ will be a meromorphic Higgs sheaf on
$X$ with logarithmic poles on $\Par_{X}$ and nilpotent residues.
\index{terms}{spectral cover!modular}
\index{terms}{logarithmic!poles}
\index{terms}{logarithmic!cotangent bundle}
\index{terms}{nilpotent residues}
\index{notations}{TvXlogParX@$T^{\vee}_{X}\left(\log \Par_{X}\right)$}
\index{notations}{X@$X$}

Hence for any assignment of real coefficients to the components of the
parabolic divisor on $Y$, i.e. any locally constant map:
\[
\bzeta_{Y} : \Par_{Y} \to \mathbb{R},
\]
we get a tame
parabolic Higgs bundle $(F_{\bullet},\varphi)$ on $(X,\Par_{X})$ with
nilpotent residues, defined as
\index{terms}{Higgs!bundle!tame}
\begin{equation}\label{eq:Fphi}
\begin{aligned}
F_{\bullet} & = f_{*}\left(\bmu^{*}\mycal{L}_{\mathbb{E}}(\bzeta_{Y}\cdot
\Par_{Y})_{\bullet}\right), \\
\varphi & = f_{*}(\alpha\otimes -),
\end{aligned}
\end{equation}
where $\mu$ is the map defined in \eqref{eq:modify.Y} and
$\mycal{L}_{\mathbb{E}}$ is the translation invariant line bundle
defined in step {\bfseries (iii)}. Since the spectral cover $Y$ is
irreducible, the pair $(F_{\bullet},\varphi)$ will be automatically
stable as a a parabolic Higgs bundle.  (Any potentially
destabilizing Higgs subsheaf would have to be supported on a component
of $Y$.)  Therefore it can be fed into Mochizuki's non-abelian Hodge
correspondence provided that it satisfies the vanishing of Chern
classes needed for the correspondence.
\index{terms}{parabolic!Higgs bundle}
\index{notations}{zetaY@$\bzeta_{Y}$}
\index{terms}{parabolic!Chern classes}

Therefore we must solve the overdermined system of affine linear and
quadratic equations on $\bzeta_{Y}$ given by
\index{terms}{Mochizuki's!conditions}
\begin{equation}
\label{eq:Chern.class.conditions}
\parch_{1}(F_{\bullet}) = 0 \quad \text{and} \quad \parch_{2}(F_{\bullet}) = 0.
\end{equation}
\item We need to impose a further constraint on $\bzeta_{Y}$ 
in order for $(F_{\bullet},\varphi)$ to be a parabolic Hecke eigensheaf.
This will then allow us to convert $(F_{\bullet},\varphi)$ to an automorphic
  twisted $\mathcal{D}$-module on
  $\Bun_{0}(\blambda)$ which is a Hecke
  eigensheaf with eigenvalue $\mathbb{V}$.
Concretely this means the following. Given a representation $\rho :
\lan{G} \to GL_{N}(\mathbb{C})$ we consider the moduli $\Heck^{\rho}$ of
quadruples $((E,\mathsf{r}),(E',\mathsf{r}'),x,\mathsf{elm})$ where
$(E,\mathsf{r}), (E',\mathsf{r}') \in
\Bun(\blambda)$, $x \in C$, and
\[
\mathsf{elm} : (E,\mathsf{r})_{|C - \{x\}} \stackrel{\cong}{\longrightarrow}
(E',\mathsf{r}')_{|C - \{x\}}
\]
is an isomorphism of parabolic bundles away from $x$ whose pole at $x$
is bounded by the highest weight $\mathsf{high}^{\rho}$ of $\rho$,
i.e. for every finite dimensional irreducible complex
\index{terms}{highest weight}
\index{notations}{highrho@$\mathsf{high}^{\rho}$}
\index{notations}{rho@$\rho$}
\index{notations}{sigma@$\sigma$}
\index{terms}{parabolic!bundles}
\index{notations}{elm@$\mathsf{elm}$}
representation $\sigma  : G \to GL_{R}(\mathbb{C})$, the
isomorphism $\mathsf{elm}$ induces an injective map of locally free
sheaves on $C$:
\[
\sigma(\mathsf{elm}) : \sigma(E) \hookrightarrow \sigma(E')\otimes
\mathcal{O}_{C}\left(\left\langle
\mathsf{high}^{\sigma},\mathsf{high}^{\rho}\right\rangle\cdot x \right).
\]
Here $\sigma(E)$ and $\sigma(E')$ denote the vector bundles associated
with $E$ and $E'$ via $\sigma$, $\mathsf{high}^{\sigma}$ and
$\mathsf{high}^{\rho}$ denote the highest weights of $\sigma$ and
$\rho$ respectively, and $\langle-,-\rangle$  is the Langlands duality
pairing between characters of $G$ and $\lan{G}$.

The moduli space $\Heck^{\rho}$ has natural morphisms
\begin{equation} \label{eq:maps.from.Heck}
\xymatrix@R-2pc{ \Heck^{\rho} \ar[r] &
  \Bun(\blambda)
  \\ ((E,\mathsf{r}),(E',\mathsf{r}'),x,\mathsf{elm}) \ar[r] &
  (E,\mathsf{r}) \\
& \\
  \Heck^{\rho} \ar[r] &
  \Bun(\blambda)\times C
  \\ ((E,\mathsf{r}),(E',\mathsf{r}'),x,\mathsf{elm}) \ar[r] &
  ((E',\mathsf{r}'),x) }
\end{equation}
and so gives a correspondence between
$\Bun(\blambda)$ and
$\Bun(\blambda)\times C$.

We have to choose a compactification and a 
resolution $H$ of the Hecke correspondence $\Heck^{\rho}$ so that
\begin{itemize}
\item The maps \eqref{eq:maps.from.Heck} extend to proper morhisms
  \index{notations}{H@$H$}
  \index{notations}{ParH@$\Par_{H}$}
 \index{notations}{zetaH@$\bzeta_{H}$}
\[
\xymatrix@R-1pc{
& H \ar[dl]_-{p} \ar[dr]^-{q} \\
X & & X\times C
}  
\]
\item If we set $\Par_{H} = p^{*}\Par_{X} + q^{*}(\Par_{X}\times C) +
  q^{*}(X\times \Par_{C})$, then away from some subvariety of
  codimension three in $H$ we have that $H$ is smooth, and $\Par_{H}$ is a
  normal crossings divisor.
\item Away from subvarieties of codimension three in $H$ and $X\times
  C$ the map \linebreak $q : H \to X\times C$ is a proper semistable
  morphism with discriminant contained in $\Par_{X}\times C +
  X\times\Par_{C}$ and such that the horizontal part of the divisor
  $\Par_{H}$ is smooth over $X\times C$. In particular the algebraic
  formula from \cite{dps} for the $L^{2}$ pushforward via $q$ of tame
  parabolic Higgs bundles on $H$ will
  hold. \index{terms}{pushforward!of parabolic Higgs bundles}
\end{itemize}

\

\noindent
For this $H$ we need next to choose an assignment $\bzeta_{H} :
\Par_{H} \to \mathbb{R}$ of real coefficients for the components of
$\Par_{H}$ so that $(F_{\bullet},\varphi)$ is a Hecke eigensheaf with
eigenvalue $\mathbb{E}$ for the Hecke operator given by a kernel
parabolic Higgs line bundle on $H$. This parabolic Higgs line bundle
is the non-abelian Hodge counterpart of the intersection cohomology
sheaf on the Hecke stack and so \cite[Section~2]{dp-Langlands} its
restriction to the complement of the parabolic divisor is equal to the
trivial line bundle equipped with zero Higgs field. Therefore in the
Higgs setting the Hecke kernel is given by a parabolic Higgs line
bundle of the form
$(\mathcal{O}_{H}(\bzeta_{H}\Par_{H})_{\bullet},0)$ on $H$. In other
words we have to choose the real vectors $\bzeta_{Y}$ and
$\bzeta_{H}$ so that $(F_{\bullet},\varphi)$ satisfies the parabolic
Hecke eigensheaf condition
\index{notations}{Fphi@$(F_{\bullet},\varphi)$}
\index{terms}{Hecke!eigensheaves!parabolic}
\index{terms}{Hecke!kernel}
\begin{equation}
  \label{eq:parabolic.Hecke.eigensheaf}
\boxed{
q_{*}\left(p^{*}(F_{\bullet},\varphi)\otimes
(\mathcal{O}_{H}\left(\bzeta_{H}\Par_{H})_{\bullet},0\right)\right) =
p_{X}^{*}(F_{\bullet},\varphi)\otimes
p_{C}^{*}\left(\rho(\mathbb{E})_{\bullet}\right).}
\end{equation}
Note that $\rho$ appears implicitly on the left hand side of
\eqref{eq:parabolic.Hecke.eigensheaf} since $H$ is a resolution of
$\Heck^{\rho}$ and so $p$, $q$, and $\bzeta_{H}$ all depend on
$\rho$. Similarly $\bzeta_{Y}$ appears implicitly on both sides of
\eqref{eq:parabolic.Hecke.eigensheaf} becuase $F_{\bullet}$ is defined
in terms of $\bzeta_{Y}$ through the formula \eqref{eq:Fphi}.

Also, since the Hecke operators (and similarly the abelianized Hecke
operators) form a commutative algebra \cite[Section~2]{dp-Langlands}
it is enough to verify condition \eqref{eq:parabolic.Hecke.eigensheaf}
for any coleection of $\rho$'s generating the representation ring of
$\lan{G}$. Thus in our case where $\lan{G}$ is $SL_{2}(\mathbb{C})$ it
suffices to check \eqref{eq:parabolic.Hecke.eigensheaf} only for the basic
Hecke correspondence corresponding to the fundamental representation
of $SL_{2}(\mathbb{C})$.

\item Once we have ensured that  $\bzeta_{Y}$  is
  chosen so that \eqref{eq:Chern.class.conditions}
  holds for all $\rho$, apply Mochizuki's non-abelian
  Hodge correspondence on $(X,\Par_{X})$ and convert
  $(F_{\bullet},\varphi)$ to a  tame purely imaginary  parabolic flat bundle
\[
(W_{\bullet},\nabla) = \mathsf{nah}_{(X,\Par_{X})} (F_{\bullet},\varphi).
\]
If in addition $\bzeta_{H}$ is chosen so that
\eqref{eq:parabolic.Hecke.eigensheaf} holds, the parabolic flat bundle
$(W_{\bullet},\nabla)$ will satisfy the de Rham version
\index{terms}{de Rham!version} of the
parabolic Hecke eigensheaf property. This is automatic since
\cite{dps} ensures that the algebraic de Rham and Dolbeault $L^{2}$
pushforwards by $q$ will be compatible with the non-abelian Hodge
corresponcences on $H$ and $X\times C$.
\index{notations}{nahparX@$\mathsf{nah}_{X,\Par_{X}}$}

The choice of $\bzeta_{Y}$ determines the parabolic weights of
$F_{\bullet}$. Using the conversion table from \cite[page
720]{carlos-nc} or the equivalent table from \cite[page 12]{dps},
these in turn determine the eigenvalues of the residues of the purely
imginary tame parabolic flat bundle $(W_{\bullet},\nabla)$. Taking the
sum of copies of the components of the parabolic divisor with
coefficients the eigenvalues of the residues of $(W_{\bullet},\nabla)$
yields an element in the complexified Picard of $\sBun$. This
element depends only on $\bzeta_{Y}$ which in turn is a solution of a
set of equations determined by $\blambda$ so we will define
$\Aom(\blambda)$ to be this element in the complexified Picard.
\index{notations}{Ok@$\Aom$}

\index{terms}{conversion table!of the ramified nah}

\item Finally note that the very stable locus $U = X - \Par_{X}$ is a
  Zariski open subset in the moduli stack $\sBun$ of
  quasi-parabolic $G$ bundles on $(C,\Par_{C})$. Using the  open immersion
$j : U
\hookrightarrow \sBun$
we can construct an $\Aom(\blambda)$-twisted $\mathcal{D}$-module
$\mathfrak{c}_{\blambda}\left(\mathcal{O}_{\mathbb{V}}\right)$
on $\sBun$ by setting
\[
\mathfrak{c}_{\blambda}\left(\mathcal{O}_{\mathbb{V}}\right) :=
j_{!*}^{\Aom(\blambda)}\left((W_{\bullet},\nabla)_{|U}\right).
\]
Here $(W_{\bullet},\nabla)_{|U}$ is the holomorphic flat bundle on $U$
obtained by restricting any level of the parabolic flat bundle
$(W_{\bullet},\nabla)$ and $j_{!*}^{\Aom(\blambda)}$ is the
twisted middle perversity Deligne-Goresky-MacPherson extension to $\sBun$. 
\index{terms}{DeligneGG@Deligne-Goresky-MacPherson!extension}
  \index{notations}{Buns@$\sBun$}

A twisted reformulation of Mochizuki's extension theorem then implies
that the twisted $\mathcal{D}$-module
$\mathfrak{c}_{\blambda}\left(\mathcal{O}_{\mathbb{V}}\right)$
will be  a Hecke eigensheaf with eigenvalue $\mathbb{V}$.
 \end{list}

 \

 \noindent
The steps {\bfseries (i)-(vii)} follow the strategy from \cite{dp-jdg}
almost verbatim. The only difference is in step {\bfseries (vii)}
where instead of of the standard Deligne-Goresky-MacPherson extension
we have to use a twisted Deligne-Goresky-MacPherson extension. This is
necessary since in the setting of the unramified GLC from
\cite{dp-jdg} the Hecke eigensheaf is an ordinary $\mathcal{D}$-module
while in our setting of the tamely ramified GLC, the Hecke eigensheaf
is a twisted $\mathcal{D}$-module.  For completeness we review twisted
Deligne-Goresky-MacPherson extensions and the twisted reformulation of
Mochizuki's extension theorem in the next two sections.
\index{terms}{DeligneGG@Deligne-Goresky-MacPherson!extension}
 
 \

 \subsection{Twisted Deligne-Goresky-MacPherson extensions}
 \label{ssec:twisted.GM}

\noindent
In the last step of the strategy explained in the previous section we
needed to construct a minimal extension of a holomorphic
flat bundle defined on a quasi-projective variety $U$ to a twisted
$\mathcal{D}$-module on a compactification of $U$. This is a minor
modification of the Deligne-Goresky-MacPherson extension functor which
we review next.

\index{terms}{DeligneGG@Deligne-Goresky-MacPherson!extension}

Suppose $X$ is a complex smooth irreducible variety, let
$j : U \hookrightarrow X$ be the inclusion of a Zariski open subset,
and $i: Z = X - U \hookrightarrow X$ the complementary closed
embedding.  These give (see e.g. \cite{htt}) the standard functors/adjunctions
$i^{*} \dashv i_{*} = i_{!} \dashv i^{!}$ and
$j_{!} \dashv j^{!} = j^{*} \dashv j_{*}$ of the recollement picture
relating the derived categories of holonomic $\mathcal{D}$-modules on
$X$, $U$, and $Z$. We also have the classical notion of a minimal (or middle
perversity) Deligne-Goresky-MacPherson extension $j_{!*}$ for holonomic
$\mathcal{D}$-modules:
\index{terms}{DeligneGG@Deligne-Goresky-MacPherson!extension}
\index{terms}{recollement}
\index{terms}{recollement!functors}
\index{terms}{recollement!adjunctions}
\index{notations}{isubstar@$i_{*}$}
\index{notations}{isupstar@$i^{*}$}
\index{notations}{isubshriek@$i_{\shriek}$}
\index{notations}{isupshriek@$i^{\shriek}$}
\index{notations}{jsubstar@$j_{*}$}
\index{notations}{jsupstar@$j^{*}$}
\index{notations}{jsubshriek@$j_{\shriek}$}
\index{notations}{jsupshriek@$j^{\shriek}$}
\index{terms}{Dmodule@$\mathcal{D}$-module!holonomic}

\

\noindent
For every
  holonomic $\mathcal{D}$-module $N$ on $U$ there exists a unique
  extension of $N$ to a holonomic $\mathcal{D}$-module $j_{!*}N$ on $X$
  characterized completely (see \cite{bbd,gm1,htt})
  by the properties:
\begin{itemize}
\item $j^{*}j_{!*}N = N$, i.e. $j_{!*}N$ is an extension of $N$;
\item $j_{!*}N$ has no sub or quotient $\mathcal{D}_{X}$-module
  supported on $Z$.
\end{itemize}

\index{terms}{Dmodule@$\mathcal{D}$-module!holonomic}
\index{notations}{jshriekstar@$j_{\shriek *}$}

\

\noindent 
The $\mathcal{D}$-module $j_{!*}N$ is selected by cohomological
bounds.  We have a natural morphism $j_{!}N \to j_{*}N$ in the derived
category of holonomic $\mathcal{D}$-modules.  Since for an open
immersion, pushforward of $\mathcal{D}$-modules is right exact, while
compactly supported pushforward of $\mathcal{D}$-modules is left
exact, we will have $j_{!}N \in D_{\op{hol}}^{\leq 0}(X,\mathcal{D})$ and
$j_{*}N \in D_{\op{hol}}^{\geq 0}(X,\mathcal{D})$.
  
Thus the morphism $j_{!}N \to j_{*}N$ has to
factor through the degree zero cohomology sheaves:
\[
j_{!}N \to \mycal{H}^{0}(j_{!}N) \to \mycal{H}^{0}(j_{*}N) \to j_{*}N
\]
and $j_{!*}N$ is defined as the image
\[
j_{!*}N := \op{im}\left[\mycal{H}^{0}(j_{!}N) \to \mycal{H}^{0}(j_{*}N)\right].
\]
Checking the support condition is straightforward from this
construction.

The \emph{\bfseries Deligne-Goresky-MacPherson functor} $j_{!*}$ has
many useful features.  It commutes with Verdier duality by definition
and sends irreducible $\mathcal{D}$-modules to irreducible
$\mathcal{D}$-modules. Also, all irreducible holonomic
$\mathcal{D}$-modules are minimal extensions of irreducible flat
bundles from irreducible locally closed subvarieties.

\index{terms}{DeligneGG@Deligne-Goresky-MacPherson!functor}
\index{notations}{jshriekstar@$j_{\shriek *}$}

In particular  the category of irreducible holonomic
$\mathcal{D}$-modules is equivalent to the category of equivalence
classes of pairs $(U,N)$ where $U \subset X$ is a locally closed
subvariety, $N$ is an irreducible flat bundle on $U$, and the
equivalence relation is generated by the elementary equivalences
$(U,N)\sim (U',N_{|U'})$ for all open immersions $U' \subset U$.

Finally, recall that the Deligne-Goresky-MacPherson extension is the
de Rham counterpart of the middle perversity extension on the Betti
side. That is, under the Riemann-Hilbert equivalence the action of
$j_{\shriek *}$ on regular holonomic $\mathcal{D}$-modules corresponds
to the middle perversity extension functor for complexes of
constructible sheaves \cite{bbd}.

\subsubsection{Deligne-Goresky-MacPherson extensions of
  flat bundles} \label{ssec:gm.flat}

Suppose that $D = X-U$ is a strict normal crossings divisor. Let
$D = D_{1}\cup D_{2} \cup \ldots \cup D_{k}$ be the smooth irreducible
components of $D$ and let $(V_{U},\nabla_{U})$ be a flat bundle on $U$
with regular singularities along $D$. For the generalization to the
TDO setting we will need to have a more explicit model for the
Deligne-Goresky-MacPherson extension of $(V_{U},\nabla_{U})$ to
$X$. We will use the following  construction explained to us by
D. Arinkin.

\index{terms}{DeligneGG@Deligne-Goresky-MacPherson!extension}

\

\bigskip

\noindent
Choose a meromorphic flat bundle $(V,\nabla)$ on $X$ such
that\footnote{It is easy to construct such an extension of
  $(V_{U},\nabla_{U})$, e.g. we could take Deligne's canonical
  extension.  }
\begin{itemize}
\item $V$ is a vector bundle on $X$ satisfying $V_{|U} = V_{U}$ and
 $\nabla : V \to V\otimes \Omega^{1}_{X}(\log D)$ is a
  logarithmic connection with $\nabla_{|U} = \nabla_{U}$.
\item For any $i$ the residue $\op{res}_{D_{i}}(\nabla) \in
  H^{0}(D_{i},\op{End}(V))$ has no eigenvalues that are negative integers.
\end{itemize}

\noindent
Next consider the quasicoherent sheaf $V(*D)$ of sections in $V$ with
poles of arbitrary order along $D$. The connection $\nabla$ endows
$V(*D)$ with a natural $\mathcal{D}_{X}$-module structure - the unique
action generated by multiplication by functions and the action of 
vector fields given by $\nabla_{\xi}$ for $\xi \in T_{X}$. By
definition $V$ is a subsheaf of $V(*D)$ and we can construct the sub
$\mathcal{D}_{X}$-module $V_{!*}$ of $V(*D)$ generated by $V$.  The
$\mathcal{D}$-module $V_{!*}$ is the Deligne-Goresky-MacPherson
extension of $(V_{U},\nabla_{U})$, i.e. we have
\[
V_{!*} = j_{!*}(V_{U},\nabla_{U}).
\]
Indeed the question is local near each component since by the normal
crossings property the local fundamental group of $U$ near $D$ is
abelian and so the local contributions of the components of $D$
commute in the Deligne-Goresky-MacPherson extension. But locally near a
component of $D$ the comparison question is essentially one
dimensional since only the transversal direction to the component
matters.

This reduces the question to the following

\index{notations}{Vshriekstar@${V_{\shriek *}}$}
\index{notations}{VUnablaU@$(V_{U},\nabla_{U})$}
\index{notations}{VstarD@$V(*D)$}
\index{notations}{jshriekstar@$j_{\shriek *}$}

\

\noindent
{\bfseries Example:} Let $X$ be a smooth curve and let $p\in X$ be a
point with $U = X - \{p\}$. Let $(V_{U},\nabla_{U})$ and $(V,\nabla)$
be as above.  Now if $z$ is a local coordinate centered at $p$ we
have by definition that $V_{!*}$ is  the subsheaf in $V(*p)$ given by 
\[
V_{!*} = V + \nabla_{\frac{d}{dz}} V + \left(\nabla_{\frac{d}{dz}}\right)^{2} V +
\left(\nabla_{\frac{d}{dz}}\right)^{3} V + \cdots \ \subset V(*p).
\]
On the other hand $j_{!*}(V_{U},\nabla_{U})$ in this case is very
simple to describe. Since for an open immersion $j_{*}$ is exact and
is the same for $\mathcal{D}$-modules and $\mathcal{O}$-modules we get
\[
R^{0}j_{*}(V_{U},\nabla_{U}) = j_{*}(V_{U},\nabla_{U}) = j_{*}V_{U} = V(*p).
\]
Similarly we have $R^{0}j_{!}(V_{U},\nabla_{U}) = j_{!}(V_{U},\nabla_{U}) $, and 
since $ j_{!}(V_{U},\nabla_{U})$ is a subsheaf in $j_{*}V_{U},\nabla_{U}) $ we have
that in this case
\[
j_{!*}(V_{U},\nabla_{U}) = j_{!}(V_{U},\nabla_{U}).
\]
To compare with $V_{!*}$ note that the submodule
$j_{!}(V_{U},\nabla_{U}) \subset j_{*}(V_{U},\nabla_{U}) = j_{*}V_{U}$
can be understood as the kernel of the map from $ j_{*}V_{U}$ onto its
``$\mathcal{D}$-module fiber at $p$''. 
  
Explicitly from the standard $\mathcal{D}$-module recollement (see
\cite{bbd}) for the inclusions \linebreak
$U\stackrel{j}{\hookrightarrow} X \stackrel{i}{\hookleftarrow} p$ we
have that
\[
j_{!*}(V_{U},\nabla_{U}) =j_{!}(V_{U},\nabla_{U}) = \ker \left[
  j_{*}V_{U} \to i_{*}i^{*}j_{*}V_{U}\right] = \ker \left[ V(*p) \to
  i_{*}i^{*}V(*p)\right],
\]
where $i^{*}$ assigns the $\mathcal{D}$-module fiber at $p$ and
$i_{*}$ assigns this fiber as a $\mathcal{D}$-module with support at
$p$ (necessarily a sum of delta-function $\mathcal{D}$-modules).
Equivalently $i_{*}i^{*}V(*p)$ is the maximal quotient of $V(*p)$ that
is a direct sum of delta-function $\mathcal{D}$-modules at $p$ and so
$j_{!*}(V_{U},\nabla_{U}) \subset V(*p)$ is the minimal submodule of
$V(*p)$ that has no nonzero maps to any delta-function
$\mathcal{D}$-module at $p$. Because of our condition that
$\op{res}_{p} \nabla$ has no negative integer eigenvalues this
extension is minimal and so must be equal to
$j_{!*}(V_{U},\nabla_{U})$.  Also by definition $V_{!*}$ has no nozero
maps to the delta-function $\mathcal{D}$-module at $p$ and again the
condition on $\op{res}_{p} \nabla$ ensures that $V_{!*}$ is minimal
with this property yielding the equality $V_{!*} =
j_{!*}(V_{U},\nabla_{U}) $ as promised.

\index{terms}{deltafun@delta function!$\mathcal{D}$-module}

As a consistency check note that if we choose a different extension
$(V',\nabla')$ of $(V_{U},\nabla_{U})$ satisfying
$\op{spectrum}(\op{res}_{p} \nabla') \cap \mathbb{Z}_{< 0} =
\varnothing$, then $V(-Mp) \subset V' \subset V(-Np)$ for some
$M \geq N > 0$.  In particular $V(*p)$ can be identified with $V'(*p)$
and since $(\nabla_{d/dz})^{N} V' \subset V$ and
$(\nabla_{d/dz})^{M}V \subset V'$ we get $V'_{!*} = V_{!*}$ as we
should.

\

\bigskip

\noindent
In low brow terms we can construct the minimal extension $V_{!*}$ via
Hecke modifications of $V$. As above the question reduces readily to the
curve case so we will only spell out the construction in this setting.

The subsheaf $V_{!*} \subset V(*p)$ can be described as the result of
an infinite sequence of iterated up Hecke transforms of $V$ at
$p$. Indeed, we defined $V_{!*}$ as the sub $\mathcal{D}$-module of
$V(*p)$ generated by $V$, i.e. we had
\[
V_{!*} = V + \nabla_{\frac{d}{dz}} V + \left(\nabla_{\frac{d}{dz}}\right)^{2} V +
\left(\nabla_{\frac{d}{dz}}\right)^{3} V + \cdots \ \subset V(*p).
\]
Explicitly we can think of $V_{!*}$ as a union (colimit) of successive
up Hecke transforms.  Start with  ${}^{0}V = V \subset V(*p)$ and consider
\[
{}^{1}V = V + \nabla_{\frac{d}{dz}} V \subset  V(*p).
\]
${}^{1}V$ is a locally free sheaf which is the up Hecke modification of
${}^{0}V$ such that the kernel of the map
\[
{}^{0}V_{p} \to {}^{1}V_{p}
\]
between fibers is the image of the residue of $\nabla$ at $p$.

This process iterates in an obvious manner. The $\mathcal{D}$-module
structure on $V(*p)$ induces a logarithmic connection ${}^{1}\nabla$ on
$^{1}V$ and we can consider
\[
{}^{2}V = {}^{1}V + {}^{1}\nabla_{\frac{d}{dz}} {}^{1}V \subset V(*p)
\]
which is the up Hecke modification of ${}^{1}V$ such that the kernel
of the map
\[
{}^{1}V_{p} \to {}^{2}V_{p}
\]
between fibers is the image of the residue of ${}^{1}\nabla$ at $p$.
Continuing in this manner we get a nested sequence of logarithmic connections
\[
(V,\nabla) = \left({}^{0}V,{}^{0}\nabla\right) \subset
\left({}^{1}V,{}^{1}\nabla\right) \subset
\left({}^{2}V,{}^{2}\nabla\right) \subset \cdots \subset
\left(V(*p),\nabla\right)
\]
and we have $V_{!*} = 
\left(\bigcup_{i} {}^{i}V,\nabla\right)$.
\index{terms}{Hecke!modification}
\index{terms}{flat!bundle!logarithmic}
\index{notations}{Vshriekstar@${V_{\shriek *}}$}

\

\noindent
The above descriptions of $V_{!*}$ are tailor made for generalization
to the twisted setting. Consider again the case when $X$ is a smooth
compact curve, $p \in X$ is a point, $U = X - \{p\}$, $j : U \hookrightarrow
X$ and $(V_{U},\nabla_{U})$ is a flat bundle on $U$ with at worst
regular singularity at $p$. Fix $c \in \mathbb{C}$ and let
$\mathcal{D}_{c}$ be the sheaf of TDO associated to the element
$\mathcal{O}_{X}(c\cdot p) \in \Pic(X)\otimes \mathbb{C}$ in the
complexified Picard of $X$.

\index{terms}{Picard!group!complexified}
\index{terms}{twisted!differential operators}
\index{notations}{Dc@$\mathcal{D}_{c}$}
\index{notations}{PicXotC@$\Pic(X)\otimes \mathbb{C}$}
\index{notations}{res@$\res$}

As before we choose an extension $(V,\nabla)$ of $(V_{U},\nabla_{U})$
to $X$ for which:
\begin{itemize}
\item $V$ is a vector bundle on $X$ with $V_{|U} = V_{U}$ and $\nabla
  : V \to V\otimes \Omega^{1}_{X}(\log p)$ is a logarithmic connection
  with $\nabla_{|U} = \nabla_{U}$.
\item The residue $\res_{p} \nabla \in \op{End}(V_{p})$ has no
  eigenvalue in $c + \mathbb{Z}_{<0} \subset \mathbb{C}$.
\end{itemize}

Now define a subsheaf $V_{!*}^c \subset V(*p)$ 
by setting
\[
V_{!*}^{c} = V + \left(\nabla_{\frac{d}{dz}} - \frac{c}{z}\cdot \op{id}\right)V +
\left(\nabla_{\frac{d}{dz}} - \frac{c}{z}\cdot \op{id}\right)^{2}V +
\left(\nabla_{\frac{d}{dz}} - \frac{c}{z}\cdot \op{id}\right)^{3}V + \cdots
\quad V(*p).
\]
This subsheaf $V_{!*}^{c} \subset V(*p)$ carries a natural structure of a
$\mathcal{D}_{c}$-module. Indeed since $\mathcal{D}_{c|U}$ is
naturally isomorphic to the ordinary differential operators
$\mathcal{D}_{U}$, we can view $\mathcal{D}_{c}$ as a subsheaf of
rings in the sheaf $\mathcal{D}_{X}(*p)$ of differential operators on
$X$ whose coefficients are meromorphic functions with poles at most at
$p$. In fact $\mathcal{D}_{c} \subset \mathcal{D}_{X}(*p)$ is the
subsheaf generated by multiplication by sections in $\mathcal{O}_{X}$
and the  meromorphic differential operator $d/dz - c/z \in
\mathcal{D}^{\leq 1}(*p)$ .

\index{notations}{Vshriekstarc@${V_{\shriek *}^{c}}$}

As before the connnection $\nabla$ endows $V(*p)$ with a
$\mathcal{D}_{X}(*p)$-module structure and by construction the
subsheaf $V_{!*}^{c} \subset V(*p)$ is preserved by the subring
$\mathcal{D}_{c} \subset \mathcal{D}_{X}(*p)$. Thus $V_{!*}^{c}$ is a
well defined $\mathcal{D}_{c}$-module which again is easily seen to be
independent of the particular extension $(V,\nabla)$ of
$(V_{U},\nabla_{U})$.

Alternatively we can define $V_{!*}^{c}$ as a union of a sequence of
up Hecke transforms ${}^{k}V^{c}$ of $V$ by setting ${}^{0}V^{c} = V$,
${}^{0}\nabla = \nabla$, and defining inductively ${}^{k+1}V^{c}$ to
be
\[
  {}^{k+1}V^{c} = {}^{k}V^{c} + \left( {}^{k}\nabla_{\frac{d}{dz}}
    - \frac{c}{z}\cdot \op{id}\right){}^{k}V^{c} \quad V(*p),
\]
that is - ${}^{k+1}V^{c}$ is the up Hecke transform of ${}^{k}V^{c}$
for which the kernel of the natural map
$\left({}^{k}V^{c}\right)_{p} \to \left({}^{k+1}V^{c}\right)_{p}$ is
the image of $\res_{p} {}^{k}\nabla - c\cdot \op{id}$. To continue the
iteration we set ${}^{k+1}\nabla$ to be the connection on
${}^{k+1}V^{c}$ induced from the $\mathcal{D}_{X}(*p)$ module
structure on $V(*p)$. Similarly to the untwisted case we get
$V_{!*}^{c} = \left( \bigcup_{i} {}^{i}V^{c}, {}^{i}\nabla\right)$

Now suppose we have a smooth (space or stack) $X$ and a strict normal
crossing divisor $D = \sum_{i} D_{i}$ in $X$. Fix a collection of
complex numbers $c = \{c_{i}\}_{i}$ and let $\mathcal{D}_{c}$ be the
sheaf of twisted differential operators corresponding to the element
$\mathcal{O}_{X}(c\cdot D) := \mathcal{O}_{X}(\sum_{i} c_{i} D_{i})
\in \Pic(X)\otimes \mathbb{C})$.  If $U = X - D$ and
$(V_{U},\nabla_{U})$ is a flat bundle on $U$ with at worst regular
singularities along $D$, we can apply the above construction of the
twisted Deligne-Goresky-MacPherson extension for each component of the
divisor. As usual the commutativity of the local fundamental group of
$U$ near $D$ implies that the order in which we do the extensions
across the different $D_{i}$ does not matter.

\index{terms}{DeligneGG@Deligne-Goresky-MacPherson!extension!twisted}
\index{notations}{VUnablaU@$(V_{U},\nabla_{U})$}
\index{notations}{VstarD@$V(*D)$}

Explicitly, we choose an extension $(V,\nabla)$ of
$(V_{U},\nabla_{U})$ so that $V$ is a vector bundle, and $\nabla$ is a
flat logarithmic connection such that $\res_{D_{k}} \nabla$ does not
have any eigenvalues in $c_{k} + \mathbb{Z}_{<0}$.

If $p \in D \subset X$ and $D$ has local equation
$z_{1}\cdots z_{s} = 0$ for some locally defined coordinate functions
centered at $p$, then $\mathcal{D}_{c}$ can be described as the
subring in the sheaf $\mathcal{D}_{X}(*D)$ of differential operators
with coefficients meromorphic functions with poles at $D$, generated
by multiplication by holomorphic functions and the first order
meromorphic differential operators $d/d_{z_{i}} - c_{i}/z_{i}$. Note
that even though these generators depend on the choice of local
equations of the components of $D$, the subring $\mathcal{D}_{c}$
depends only on the twisting parameters $c$ and not on the particular
functions $z_{i}$.

Following the above prescription we can now describe the $c$-twisted
Deligne-Goresky-MacPherson extension of $(V_{U},\nabla_{U})$ to be the
subsheaf $j_{!*}^{c}(V_{U},\nabla_{U})$ in $V(*D)$ given locally by
\[
j_{!*}^{c}(V_{U},\nabla_{U})  = \sum_{n = 0}^{\infty} \sum_{i =
  1}^{s}\left(\nabla_{\frac{\partial}{\partial z_{i}}} -
  \frac{c_{i}}{z_{i}} \right)^{n} V \quad  \subset V(*D).
\]
By the same token this subsheaf does not depend on the choice of local
equations $z_{i}$ and is preserved by the TDO $\mathcal{D}_{c}$. 

\index{notations}{Dc@$\mathcal{D}_{c}$}
  
\

\subsection{Remarks on untwisting} \label{ssec:untwisting}

The final input needed in step {\bfseries (vii)} of the strategy from
section~\ref{ssec:nah.strategy} is a twisted reformulation of
Mochizuki's extension theorem. In the context of the unramified
Langlands correspondence the extension theorem says that the (usual) 
Deligne-Goresky-MacPherson extension gives an equivalence of the category
of tame stable parabolic flat bundles on $(X,\Par_{X})$ with fixed
parabolic weights and the category of simple regular holonomic
$\mathcal{D}$-modules on $X$ which are smooth on $X - \Par_{X}$. This gives a
compatibility between the Hecke functors acting on parabolic flat
bundles and the Hecke functors acting on $\mathcal{D}$-modules.

\index{terms}{untwisting}
\index{terms}{ParX@$\Par_{X}$}

For the tamely ramified Langlands correspondence we need a variant of
this statement in which the ordinary Deligne-Goresky-MacPherson
extension is replaced by a twisted Deligne-Goresky-MacPherson
extension for a twisting coming from a complex linear combination of
the components of the boundary divisor. This variant is in fact not a
new statement but rather an application of the standard Mochizuki
extension theorem in the equivariant setting. This follows immediately
from the observation that the category of modules over TDO coming from
elements in the complexified Picard variety can be ``untwisted'',
i.e. can be viewed as ordinary untwisted $\mathcal{D}$-modules on a different
space satisfying an appropriate equivariance condition. This
interpretation of twisted $\mathcal{D}$-modules is well known to experts but
does not seem to appear in the literature. For completeness we recall
it here in the form explained to us by D.~Arinkin.

Given a complex number $c$ let $\mathbb{K}_{c} =
\left(\mathcal{O}_{\mathbb{G}_{m}}, d + c\cdot d\log z\right)$ denote
the Kummer rank one flat bundle on $\mathbb{G}_{m}$, where $z$ is the
coordinate on $\mathbb{G}_{m}$. Thus $\mathbb{K}_{c}$ has residue $c$
at $z=0$ and residue $-c$ at $z = \infty$.  The Kummer flat bundle
$\mathbb{K}_{c}$ is a character $\mathcal{D}$-module on the group
$\mathbb{G}_{m}$, i.e. on $\mathbb{G}_{m}\times \mathbb{G}_{m}$ we
have a canonical isomorphism $\mathsf{mult}^{*}\mathbb{K}_{c} \cong
p_{1}^{*}\mathbb{K}_{c}\otimes p_{2}^{*}\mathbb{K}_{c}$ satisfying the
obvious  cocycle condition on $\mathbb{G}_{m}\times
\mathbb{G}_{m}\times \mathbb{G}_{m}$.

Let now $X$ be a smooth space or stack and let $L$ be a line bundle on
$X$. Let \linebreak $\mathsf{Y} = \op{tot}(L^{\times}) = \op{tot}(L) -
(0-\text{section})$ denote the total space of its frame bundle,
i.e. the total space of the principlal $\mathbb{G}_{m}$-bundle
$L^{\times} = \op{Isom}_{X}(\mathcal{O}_{X},L)$. The group
$\mathbb{G}_{m}$ acts freely on $\mathsf{Y}$ and we have
$\mathsf{Y}/\mathbb{G}_{m} = X$.  Write $\mathsf{p} : \mathsf{Y} \to
X$ for the natural projection and $\mathsf{act} : \mathbb{G}_{m}\times
\mathsf{Y} \to \mathsf{Y}$ for the action map.

Consider the category $D(\mathsf{Y},\mathcal{D})^{\mathbb{G}_{m}}_{c}$
of $\mathbb{K}_{c}$-twisted $\mathbb{G}_{m}$-equivariant $\mathcal{D}$-modules
on $\mathsf{Y}$. By definition, the objects of
$D(\mathsf{Y},\mathcal{D})^{\mathbb{G}_{m}}_{c}$ are $\mathcal{D}$-modules $M$
on $\mathsf{Y}$ equipped with an isomorphism
\[
\gamma_{M} : \mathsf{act}^{*}M \stackrel{\cong}{\longrightarrow}
p_{\mathbb{G}_{m}}^{*}\mathbb{K}_{c}\otimes p_{\mathsf{Y}}^{*}M \quad \text{on}
\quad \mathbb{G}_{m}\times \mathsf{Y}.
\]
The isomorphism $\gamma_{M}$ has to satisfy the natural cocycle
condition on $\mathbb{G}_{m}\times \mathbb{G}_{m}\times \mathsf{Y}$
corresponding to the character $\mathcal{D}$-module structure on
$\mathbb{K}_{c}$.  The morphisms in
$D(\mathsf{Y},\mathcal{D})^{\mathbb{G}_{m}}_{c}$ are simply morphisms
$M_{1} \to M_{2}$ of $\mathcal{D}$-modules on $\mathsf{Y}$ intertwining the
$\mathbb{K}_{c}$-twisted equivariant structures $\gamma_{M_{1}}$ and
$\gamma_{M_{2}}$.

Now the untwisting statement is the claim that the category
$D(\mathsf{Y},\mathcal{D})^{\mathbb{G}_{m}}_{c}$ is equivalent to the
category $D(X,\mathcal{D}_{c})$ of modules over the TDO
$\mathcal{D}_{c}$ acting on $L^{\otimes c} \in \Pic(X)\otimes
\mathbb{C}$. Indeed suppose $\{ U_{i} \}$ is an open cover of $X$
trivializing $L$ and let $\ell_{i} \in \Gamma(U_{i},L)$ be a local
frame for $L_{|U_{i}}$.  Let $g_{ij} = \ell_{j}/\ell_{i} \in
\Gamma(U_{ij},\mathcal{O}^{\times})$.  Then a section of $L^{\otimes c}$
over say some open $U \subset X$ is a collection of functions $s_{i}
\in \Gamma(U\cap U_{i},\mathcal{O})$ such that $t_{ij} =
s_{i}/s_{j}$ is a solution to the differential equation $(d + c\cdot
d\log g_{ij})(t_{ij}) = 0$. In $\mathcal{D}$-module language this means that
a $\mathcal{D}_{c}$-module on $X$ is simply a collection $\{ M_{i} \}$
of $\mathcal{D}$-modules on $U_{i}$, together with compatible $\mathcal{D}$-module
isomorphisms $\alpha_{ij} : M_{i} \widetilde{\to} M_{j} \otimes
(\mathcal{O}, d + c\cdot d\log g_{ij})$ over $U_{ij}$ satisfying the
natural cocycle condition on triple overlaps.

On the other hand a $\mathbb{K}_{c}$-twisted equivariant $\mathcal{D}$-module
$M$ on $\mathsf{Y}$ gives rise to $\mathcal{D}$-modules $M_{i} = \ell_{i}^{*}M$
where we view the local section $\ell_{i}$ as maps $\ell_{i} : U_{i}
\to \mathsf{Y}$. Since $\ell_{j} = \ell_{i}g_{ij}$, 
the twisted
equivariant structure on $M$ gives isomorphisms $M_{i} \cong
g_{ij}^{*}\mathbb{K}_{c}\otimes M_{j}$ satisfying the cocycle
condition. But $g_{ij}^{*}\mathbb{K}_{c} = (\mathcal{O},d + c\cdot
d\log g_{ij})$ so this is precisely the structure of a
$\mathcal{D}_{c}$-module. It is straightforward to check that this
construction on objects extends to an equivalence of categories
$D(X,\mathcal{D}_{c}) \ \widetilde{\to} \
D(\mathsf{Y},\mathcal{D})^{\mathbb{G}_{m}}_{c}$.

More generally if $L_{1}$, $L_{2}$, \ldots, $L_{k}$ is a collection of
line bundles on $X$, and $\mathsf{c} = \{ c_{i} \}_{i =1}^{k}$ is a
collection of complex numbers we can form the TDO
$\mathcal{D}_{\mathsf{c}}$ of differential operators acting on $\otimes_{i
  = 1}^{k} L_{i}^{\otimes c_{i}} \in \Pic(X)\otimes \mathbb{C}$.  We
can also consider the principal $\mathbb{G}_{m}^{\times k}$-bundle
$\mathsf{Y} \to X$ defined as the fiber product over $X$ of the
principal $\mathbb{G}_{m}$-bundles $\mathsf{Y}_{i} =
\op{tot}(L_{i}^{\times})$.

Now the above construction gives an equivalence between the category
$D(X,\mathcal{D}_{\mathsf{c}})$ of \linebreak
$\mathcal{D}_{\mathsf{c}}$-modules and the category
$D(\mathsf{Y},\mathcal{D})^{\mathbb{G}_{m}^{\times k}}_{\mathsf{c}}$
of twisted $\mathbb{G}_{m}^{\times k}$-equivariant $D$
modules on $\mathsf{Y}$, i.e. $\mathcal{D}$-modules $M$ on $\mathsf{Y}$
equipped with isomorphisms
\[
      \mathsf{act}^{*} M \cong p_{\mathbb{G}_{m}^{\times
          k}}^{*}\left(\mathbb{K}_{c_{1}}\boxtimes \cdots \boxtimes 
      \mathbb{K}_{c_{k}}\right) \otimes p_{\mathsf{Y}}^{*} M.
\]
This allows us to recast statements about twisted
$\mathcal{D}$-modules on a space $X$ as statements about quasi
equivariant ordinary $D$ modules on the principal torus bundle
$\mathsf{Y}$ over $X$. In particular conjugating the ordinary
Deligne-Goresky-MacPherson extension functor on $Y$ with this
untwisting equivalence and applying Mochizuki's extension theorem on
$\mathsf{Y}$ will give us the twisted version of Mochizuki's extension
theorem on $X$.

\subsection{Geometric class field theory revisited}
\label{ssec-ramified-cft}

We conclude this appendix and the entire work by returning to the abelian case 
presented in sections \ref{abelian} and \ref{parabelian}. The untwisting functor
allows us to show that in the case when $G = \lan{G} = GL_{1}$ the
ramified Langlands correspondence \eqref{eq-clambdatilde} is
equivalent to the familiar character property from the tamely ramified
geometric class field theory.

Geometric class field theory is the
  abelian case of geometric Langlands. The unramified case is the
  elementary statement that every rank one flat bundle on a compact
  curve $C$ extends, uniquely, to a rank one flat bundle on the
  Jacobian $\Jac = \Jac(C) := \Pic^0(C)$. Equivalently, pulling back
  by the Abel-Jacobi map $\aj : C \to \Jac$ gives an isomorphism
  \linebreak
  $\aj^{*} : \op{\bf Flat }_{\Jac,G} \to \op{\bf Flat}_{C,G}$ for
  $G=GL_{1}(\mathbb{C}) = \mathbb{C}^{\times}$, or more generally for
  any affine abelian $G$.  \index{terms}{Jacobian}
  \index{terms}{Abel-Jacobi map} \index{terms}{flat!bundle}
  \index{notations}{aaj@$\aj$} \index{notations}{J@$\Jac$}
  \index{notations}{FlatCG@$ \op{\bf Flat}_{C,G}$}
  \index{notations}{FlatJG@$\op{\bf Flat }_{\Jac,G}$} Of course, the
  map $\aj$ depends on a base point in $C$. The more natural map is
  \linebreak $\AJ : C \times \Pic(C) \to \Pic(C)$ sending
  $x \in C, \ L \in \Pic^{d}(C)$ to
  $L \otimes \mathcal{O}_C(x) \in \Pic^{d+1}(C)$. The graph of this
  $\AJ$ is the abelian Hecke correspondence
  $H \subset C \times \Pic(C) \times \Pic(C)$. Geometric class field
  theory is the statement that every rank one flat bundle $\mathbb{L}$
  on $C$ determines a flat bundle, say
  $\mathfrak{c}(\mathcal{O}_{\mathbb{L}})$, on $\Pic(C)$, which is an
  eigensheaf for $\AJ^{*}$ with eigenvalue $\mathbb{L}$:
  \index{terms}{Hecke!eigensheaf}
\[
  \AJ^{*}(\mathfrak{c}(\mathcal{O}_{\mathbb{L}})) = \op{pr}^{*}_{C}
  \mathbb{L} \otimes \op{pr}^{*}_{\Pic(C)}
  \mathfrak{c}(\mathcal{O}_{\mathbb{L}}).
\]
We want to describe the ramified version of this.  \index{terms}{Class
  Field Theory} \index{terms}{Class Field Theory!ramified} We give
ourselves the compact curve $C$ of genus $g$ and a divisor
$D=p_1 + \cdots + p_k$ consisting of $k$ distinct points of $C$. These
determine a curve $\overline{C}$ of arithmetic genus $g+k-1$, obtained
from $C$ by gluing the $k$ points of $D$ to each other, transversally
into a normal crossings singularity. (In particular, for $k \geq 3$,
this curve is not locally planar.) We replace the Jacobian $\Jac$ of $C$
and its Picard $\Pic(C)$ by the generalized Jacobian
$\Pic^0(\overline{C})$ and the Picard variety
$\Pic(\overline{C})$. There is still an Abel-Jacobi
map
$\AJ : (C - D) \times \Pic(\overline{C}) \to \Pic(\overline{C})$
sending $x \in C - D, \ L \in \Pic^d(\overline{C})$ to
$L \otimes \mathcal{O}_{\overline{C}}(x) \in \Pic^{d+1}(\overline{C})$ as
before. But now this map does not extend from $C - D$ to 
$\overline{C}$ nor to $C$.

One way to see this is to construct a compactification of (each
component of) $\Pic(\overline{C})$. Note that $\Pic(\overline{C})$
fibers over $\Pic(C)$, the fibers being torsors for the
$k-1$-dimensional affine algebraic torus
$T := \op{Hom}(K, \mathbb{G}_{m})$, where
$K := \ker\left[\mathbb{Z}^d \stackrel{\mathsf{sum}}{\longrightarrow}
  \mathbb{Z}\right]$. The fiber of $\Pic(\overline{C})$ over
$L \in \Pic(C)$ is canonically identified with the quotient by the
diagonal action of $\mathbb{G}_{m}$ of the product of punctured
fibers $\prod_{i =1}^{k} L_{p_{i}}^{\times}$ of $L$. More globally
$\Pic(\overline{C})$ is the quotient by the diagonal action of
$\mathbb{G}_{m}$ of the direct image under the projection
$\Pic(C) \times D \to \Pic(C)$ of the restriction to
$\Pic(C) \times D$ of the universal $\mathbb{G}_{m}$-bundle on
$\Pic(C) \times C$.)  The compactification $\Pic^{+}(\overline{C})$ is
constructed as the associated fiber bundle fiber the toric variety
$\mathbb{P}^{k-1} $. Fiber by fiber, this replaces the fiber, which is
the $T$-torsor
$\left(\prod_{i=1}^k L_{p_{i}}^{\times} \right) /
\mathbb{G}_{m}$, by the projective space
${\mathbb{P}} \left(\oplus_{i=1}^k L_{p_{i}}\right)$. Now the
Abel-Jacobi map extends to a rational map
\[
  \AJ^{+}: C \times \Pic^{+}(\overline{C}) \to \Pic^{+}(\overline{C})
\]
whose indeterminacy locus has codimension $k-1$ and is the union of
$k$ sections of \linebreak $\Pic^{+}(\overline{C}) \to
\Pic(\overline{C})$ - one for each $p_{i} \in D \subset C$. The
inverse image under $\AJ^{+}$ of the boundary $\boldsymbol{\partial}
:= \Pic^{+}(\overline{C}) - \Pic(\overline{C}) $ is the union $(D
\times \Pic^{+}(\overline{C})) \cup (C \times \boldsymbol{\partial})$. In
particular, this verifies that the map $\AJ : (C - D) \times
\Pic(\overline{C}) \to \Pic(\overline{C})$ does not extend from $C -
D$ to all of $\overline{C}$ nor to $C$.

Tamely ramified geometric class field theory can now be stated
analogously to the unramified case: every meromorphic rank one flat
bundle $\mathbb{L}$ on $C$ with tame ramification along $D$ determines
a meromorphic rank one flat bundle
$\mathfrak{c}(\mathcal{O}_{\mathbb{L}})$, on $\Pic^{+}(\overline{C})$,
with tame ramification along $\boldsymbol{\partial}$, which is an
eigensheaf for $\AJ^{+*}$ with eigenvalue $\ell$:
\[
  \AJ^{+*}(\mathfrak{c}(\mathcal{O}_{\mathbb{L}})) = \op{pr}^{*}_{C}
  \mathbb{L} \otimes \op{pr}^{*}_{\Pic^{+}(\overline{C})}
  \mathfrak{c}(\mathcal{O}_{\mathbb{L}}).
\]
By construction the flat bundle
$\mathfrak{c}(\mathcal{O}_{\mathbb{L}})$ is holomorphic on the
$\Pic(\overline{C})$ and is twisted equivariant for the natural
$\mathbb{G}^{\times (k-1)}_{m}$ action. The descent of
$\mathfrak{c}(\mathcal{O}_{\mathbb{L}})$ to $\Pic(C)$ will then be the
twisted $\mathcal{D}$-module which is the Hecke eigensheaf in the sense of
\eqref{eq-clambdatilde}.

\newpage

\addcontentsline{toc}{section}{References}

\bibliographystyle{halpha}
\bibliography{langlands}

\begin{thebibliography}{DOPW01}

\bibitem[AG15]{ag-sing}
D.~{Arinkin} and D.~{Gaitsgory}.
\newblock {Singular support of coherent sheaves, and the geometric {L}anglands
  conjecture}.
\newblock {\em Selecta Math. (N.S.)}, 21(1):1--199, 2015.

\bibitem[AL97]{arinkin-lysenko}
D.~Arinkin and S.~Lysenko.
\newblock Isomorphisms between moduli spaces of {${\rm SL}(2)$}-bundles with
  connections on {$\mathbb{P}^1\setminus \{x_1,\cdots, x_4\}$}.
\newblock {\em Math. Res. Lett.}, 4(2-3):181--190, 1997.

\bibitem[Ari01]{arinkin-4points}
D.~Arinkin.
\newblock Orthogonality of natural sheaves on moduli stacks of {$\rm
  SL(2)$}-bundles with connections on {$\mathbb{P}^1$} minus 4 points.
\newblock {\em Selecta Math. (N.S.)}, 7(2):213--239, 2001.

\bibitem[Ari19]{arinkin-letter}
D.~Arinkin.
\newblock letter to the authors, 2019.
\newblock February 13.

\bibitem[Bau91]{bauer}
S.~Bauer.
\newblock Parabolic bundles, elliptic surfaces and {${\rm
  SU}(2)$}-representation spaces of genus zero {F}uchsian groups.
\newblock {\em Math. Ann.}, 290(3):509--526, 1991.

\bibitem[BB93]{bb-jantzen}
A.~Beilinson and J.~Bernstein.
\newblock A proof of {J}antzen conjectures.
\newblock In {\em I. {M}. {G}el$\prime$fand {S}eminar}, volume~16 of {\em Adv.
  Soviet Math.}, pages 1--50. Amer. Math. Soc., Providence, RI, 1993.

\bibitem[BB07]{bb}
A.~Braverman and R.~Bezrukavnikov.
\newblock Geometric {L}anglands correspondence for {$\mathcal{D}$}-modules in
  prime characteristic: the {${\rm GL}(n)$} case.
\newblock {\em Pure Appl. Math. Q.}, 3(1, part 3):153--179, 2007.

\bibitem[BBD82]{bbd}
A.~Beilinson, J.~Bernstein, and P.~Deligne.
\newblock Faisceaux pervers.
\newblock In {\em Analysis and topology on singular spaces, I (Luminy, 1981)},
  volume 100 of {\em Ast\'erisque}, pages 5--171. Soc. Math. France, Paris,
  1982.

\bibitem[BBP17]{bbp}
V.~Balaji, I.~Biswas, and Y.~Pandey.
\newblock Connections on parahoric torsors over curves.
\newblock {\em Publ. Res. Inst. Math. Sci.}, 53(4):551--585, 2017.

\bibitem[BD03]{beilinson-drinfeld-langlands}
A.~Beilinson and V.~Drinfeld.
\newblock Quantization of {H}itchin's integrable system and {H}ecke
  eigensheaves.
\newblock {B}ook, in preparation, 2003.

\bibitem[Bea83]{beauville}
A.~Beauville.
\newblock {\em Complex algebraic surfaces}, volume~68 of {\em London
  Mathematical Society Lecture Note Series}.
\newblock Cambridge University Press, Cambridge, 1983.
\newblock Translated from the French by R. Barlow, N. I. Shepherd-Barron and M.
  Reid.

\bibitem[BH95]{boden-hu}
H.~Boden and Y.~Hu.
\newblock Variations of moduli of parabolic bundles.
\newblock {\em Math. Ann.}, 301(3):539--559, 1995.

\bibitem[Bj{\"{o}}93]{bjork-analyticD}
J.-E. Bj{\"{o}}rk.
\newblock {\em Analytic {$\mathcal{D}$}-modules and applications}, volume 247
  of {\em Mathematics and its Applications}.
\newblock Kluwer Academic Publishers Group, Dordrecht, 1993.

\bibitem[BNP13]{bznp}
D.~{Ben-Zvi}, D.~{Nadler}, and A.~{Preygel}.
\newblock {Integral transforms for coherent sheaves}.
\newblock {\em ArXiv e-prints}, December 2013, 1312.7164.

\bibitem[BNR89]{bnr}
A.~Beauville, M.S. Narasimhan, and S.~Ramanan.
\newblock Spectral curves and the generalised theta divisor.
\newblock {\em J. Reine Angew. Math.}, 398:169--179, 1989.

\bibitem[Bor07]{borne-root}
N.~Borne.
\newblock Fibr\'es paraboliques et champ des racines.
\newblock {\em Int. Math. Res. Not. IMRN}, 16:38, 2007.

\bibitem[Bot95]{bottacin}
F.~Bottacin.
\newblock Symplectic geometry on moduli spaces of stable pairs.
\newblock {\em Ann. Sci. \'Ecole Norm. Sup. (4)}, 28(4):391--433, 1995.

\bibitem[BR94]{biswas-ramanan}
I.~Biswas and S.~Ramanan.
\newblock An infinitesimal study of the moduli of {H}itchin pairs.
\newblock {\em J. London Math. Soc. (2)}, 49(2):219--231, 1994.

\bibitem[BS15]{balaji.seshadri}
V.~Balaji and C.~S. Seshadri.
\newblock Moduli of parahoric {$\mathcal{G}$}-torsors on a compact {R}iemann
  surface.
\newblock {\em J. Algebraic Geom.}, 24(1):1--49, 2015.

\bibitem[BY96]{boden-yokogawa}
H.~Boden and K.~Yokogawa.
\newblock Moduli spaces of parabolic {H}iggs bundles and parabolic {$K(D)$}
  pairs over smooth curves: {I}.
\newblock {\em I. Internat. J. Math}, 7:573--598, 1996.

\bibitem[CDDP15]{parabolicMD}
{W.-y.} {Chuang}, D.-E. Diaconescu, R.~Donagi, and T.~Pantev.
\newblock {Parabolic refined invariants and {M}acdonald polynomials}.
\newblock {\em Commun. Math. Phys.}, 335(3):1323--1379, 2015.

\bibitem[DH98]{dolgachev-hu}
I.~Dolgachev and Y.~Hu.
\newblock Variation of geometric invariant theory quotients.
\newblock {\em Inst. Hautes \'Etudes Sci. Publ. Math.}, 87:5--56, 1998.
\newblock With an appendix by Nicolas Ressayre.

\bibitem[DM96]{ron-eyal}
R.~Donagi and E.~Markman.
\newblock Spectral covers, algebraically completely integrable, {H}amiltonian
  systems, and moduli of bundles.
\newblock In {\em Integrable systems and quantum groups ({M}ontecatini {T}erme,
  1993)}, volume 1620 of {\em Lecture Notes in Math.}, pages 1--119. Springer,
  Berlin, 1996.

\bibitem[Don95]{donagi-msri}
R.~Donagi.
\newblock Spectral covers.
\newblock In {\em Current topics in complex algebraic geometry (Berkeley, CA,
  1992/93)}, volume~28 of {\em Math. Sci. Res. Inst. Publ.}, pages 65--86.
  Cambridge Univ. Press, Cambridge, 1995.

\bibitem[DOPW01]{dopw}
R.~Donagi, B.~Ovrut, T.~Pantev, and D.~Waldram.
\newblock Standard-model bundles.
\newblock {\em Adv. Theor. Math. Phys.}, 5(3):563--615, 2001.

\bibitem[DP09]{dp-jdg}
R.~Donagi and T.~Pantev.
\newblock Geometric {L}anglands and non-abelian {H}odge theory.
\newblock In {\em Surveys in differential geometry. {V}ol. {XIII}. {G}eometry,
  analysis, and algebraic geometry: forty years of the {J}ournal of
  {D}ifferential {G}eometry}, volume~13 of {\em Surv. Differ. Geom.}, pages
  85--116. Int. Press, Somerville, MA, 2009.

\bibitem[DP12]{dp-Langlands}
R.~Donagi and T.~Pantev.
\newblock Langlands duality for {H}itchin systems.
\newblock {\em Invent. Math.}, 183(3):653--735, 2012.

\bibitem[DPS16]{dps}
R.~{Donagi}, T.~{Pantev}, and C.~{Simpson}.
\newblock {Direct Images in Non Abelian {H}odge Theory}.
\newblock {\em ArXiv e-prints}, December 2016, 1612.06388.

\bibitem[DPS19]{dps2}
R.~{Donagi}, T.~{Pantev}, and C.~{Simpson}.
\newblock {Hecke eigensheaves in genus two.}, 2019.
\newblock in preparation.

\bibitem[Dri80]{drinfeld-icm}
V.~Drinfeld.
\newblock Langlands' conjecture for {${GL}(2)$} over functional fields.
\newblock In {\em Proceedings of the {I}nternational {C}ongress of
  {M}athematicians ({H}elsinki, 1978)}, pages 565--574, Helsinki, 1980. Acad.
  Sci. Fennica.

\bibitem[Dri83]{drinfeld-ajm}
V.~Drinfeld.
\newblock Two-dimensional {$l$}-adic representations of the fundamental group
  of a curve over a finite field and automorphic forms on {$GL(2)$}.
\newblock {\em Amer. J. Math.}, 105(1):85--114, 1983.

\bibitem[Dri87]{drinfeld-jsm}
V.~Drinfeld.
\newblock Two-dimensional {$l$}-adic representations of the {G}alois group of a
  global field of characteristic $p$ and automorphic forms on {$GL(2)$}.
\newblock {\em J. of Soviet Math.}, 36:93--105, 1987.

\bibitem[Fal93]{faltings}
G.~Faltings.
\newblock Stable {$G$}-bundles and projective connections.
\newblock {\em J. Algebraic Geom.}, 2(3):507--568, 1993.

\bibitem[FGKV98]{fgkv}
E.~Frenkel, D.~Gaitsgory, D.~Kazhdan, and K.~Vilonen.
\newblock Geometric realization of {W}hittaker functions and the {L}anglands
  conjecture.
\newblock {\em J. Amer. Math. Soc.}, 11(2):451--484, 1998.

\bibitem[FGV01]{fgv}
E.~Frenkel, D.~Gaitsgory, and K.~Vilonen.
\newblock Whittaker patterns in the geometry of moduli spaces of bundles on
  curves.
\newblock {\em Ann. of Math. (2)}, 153(3):699--748, 2001.

\bibitem[FW08]{frenkel-witten}
E.~Frenkel and E.~Witten.
\newblock Geometric endoscopy and mirror symmetry.
\newblock {\em Commun. Number Theory Phys.}, 2(1):113--283, 2008.

\bibitem[Gai01]{dennis-nearby}
D.~Gaitsgory.
\newblock Construction of central elements in the affine {H}ecke algebra via
  nearby cycles.
\newblock {\em Invent. Math.}, 144(2):253--280, 2001.

\bibitem[{Gai}16]{dennis-quantum}
D.~{Gaitsgory}.
\newblock {Quantum {L}anglands Correspondence}.
\newblock {\em arXiv e-prints}, page arXiv:1601.05279, Jan 2016, 1601.05279.

\bibitem[GM80]{gm1}
M.~{Goresky} and R.~{MacPherson}.
\newblock Intersection homology theory.
\newblock {\em Topology}, 19(2):135--162, 1980.

\bibitem[Gro12]{HilbertHiggs}
M.~Groechenig.
\newblock Hilbert schemes as moduli of {H}iggs bundles and local systems, 2012.
\newblock arXiv:1206.5116.

\bibitem[GW06]{gukov-witten}
S.~Gukov and E.~Witten.
\newblock Gauge theory, ramification, and the geometric {L}anglands program,
  2006.
\newblock arXiv.org:hep-th/0612073.

\bibitem[Hei04]{heinloth}
J.~Heinloth.
\newblock Coherent sheaves with parabolic structure and construction of {H}ecke
  eigensheaves for some ramified local systems.
\newblock {\em Ann. Inst. Fourier (Grenoble)}, 54(7):2235--2325 (2005), 2004.

\bibitem[Hit87]{hitchin}
N.~Hitchin.
\newblock Stable bundles and integrable systems.
\newblock {\em Duke Math. J.}, 54(1):91--114, 1987.

\bibitem[HNY13]{hny}
J.~Heinloth, B.-C. Ng\^{o}, and Z.~Yun.
\newblock Kloosterman sheaves for reductive groups.
\newblock {\em Ann. of Math. (2)}, 177(1):241--310, 2013.

\bibitem[HTT08]{htt}
R.~{Hotta}, K.~{Takeuchi}, and T.~{Tanisaki}.
\newblock {\em {$D$}-modules, perverse sheaves, and representation theory},
  volume 236 of {\em Progress in Mathematics}.
\newblock Birkh\"{a}user Boston, Inc., Boston, MA, 2008.
\newblock Translated from the 1995 Japanese edition by Takeuchi.

\bibitem[IS07]{iyer-simpson-dr}
J.~Iyer and C.~Simpson.
\newblock A relation between the parabolic {C}hern characters of the de {R}ham
  bundles.
\newblock {\em Math. Ann.}, 338(2):347--383, 2007.

\bibitem[IS08]{iyer-simpson}
J.~Iyer and C.~Simpson.
\newblock The {C}hern character of a parabolic bundle, and a parabolic
  corollary of {R}eznikov's theorem.
\newblock In {\em Geometry and dynamics of groups and spaces}, volume 265 of
  {\em Progr. Math.}, pages 439--485. Birkh\"auser, Basel, 2008.

\bibitem[Kas89]{kashiwara-tdo}
M.~Kashiwara.
\newblock Representation theory and {$D$}-modules on flag varieties.
\newblock {\em Ast\'erisque}, 173-174:9, 55--109, 1989.
\newblock Orbites unipotentes et repr\'esentations, III.

\bibitem[KS19]{komyo-saito}
A.~Komyo and M.-H. Saito.
\newblock Explicit description of jumping phenomena on moduli spaces of
  parabolic connections and {H}ilbert schemes of points on surfaces.
\newblock {\em Kyoto J. Math.}, 59(3):515--552, 2019.

\bibitem[KW06]{kw}
A.~Kapustin and E.~Witten.
\newblock Electric-magnetic duality and the geometric {L}anglands program.,
  2006.
\newblock hep-th/0604151.

\bibitem[Laf02]{lafforgue}
L.~Lafforgue.
\newblock Chtoucas de {D}rinfeld et correspondance de {L}anglands.
\newblock {\em Invent. Math.}, 147(1):1--241, 2002.

\bibitem[Laf18]{vlafforgue}
V.~Lafforgue.
\newblock Chtoucas pour les groupes r\'{e}ductifs et param\'{e}trisation de
  {L}anglands globale.
\newblock {\em J. Amer. Math. Soc.}, 31(3):719--891, 2018.

\bibitem[Lau87]{laumon-langlands}
G.~Laumon.
\newblock Correspondance de {L}anglands g\'eom\'etrique pour les corps de
  fonctions.
\newblock {\em Duke Math. J.}, 54(2):309--359, 1987.

\bibitem[Lau88]{laumon-nilpotent}
G.~Laumon.
\newblock Un analogue global du c\^one nilpotent.
\newblock {\em Duke Math. J.}, 57(2):647--671, 1988.

\bibitem[Lau95]{laumon-gln}
G.~Laumon.
\newblock Faisceaux automorphes pour {GL(n)}: la premiere construction de
  drinfeld, 1995.
\newblock arXiv.org:alg-geom/9511004.

\bibitem[Lau03]{laumon-fgv}
G.~Laumon.
\newblock Travaux de {F}renkel, {G}aitsgory et {V}ilonen sur la correspondance
  de {D}rinfeld-{L}anglands.
\newblock {\em Ast\'erisque}, 290:Exp. No. 906, ix, 267--284, 2003.
\newblock S{\'e}minaire Bourbaki. Vol. 2001/2002.

\bibitem[LS97a]{laszlo-sorger}
Y.~Laszlo and C.~Sorger.
\newblock The line bundles on the moduli of parabolic {$G$}-bundles over curves
  and their sections.
\newblock {\em Ann. Sci. \'Ecole Norm. Sup. (4)}, 30(4):499--525, 1997.

\bibitem[LS97b]{ls}
Y.~Laszlo and C.~Sorger.
\newblock The line bundles on the moduli of parabolic {$G$}-bundles over curves
  and their sections.
\newblock {\em Ann. Sci. \'{E}cole Norm. Sup. (4)}, 30(4):499--525, 1997.

\bibitem[LS15]{loray-saito}
F.~Loray and M.-H. Saito.
\newblock Lagrangian fibrations in duality on moduli spaces of rank 2
  logarithmic connections over the projective line.
\newblock {\em Int. Math. Res. Not. IMRN}, 4:995--1043, 2015.

\bibitem[Mar82]{maruyama}
M.~Maruyama.
\newblock Elementary transformations in the theory of algebraic vector bundles.
\newblock In {\em Algebraic geometry (La R\'abida, 1981)}, pages 241--266.
  Springer, Berlin, 1982.

\bibitem[Mar94]{eyal}
E.~Markman.
\newblock Spectral curves and integrable systems.
\newblock {\em Compositio Math.}, 93(3):255--290, 1994.

\bibitem[Moc06]{mochizuki-kh1}
T.~Mochizuki.
\newblock Kobayashi-{H}itchin correspondence for tame harmonic bundles and an
  application.
\newblock {\em Ast\'erisque}, 309:viii+117, 2006.

\bibitem[Moc07a]{mochizuki-D1}
T.~Mochizuki.
\newblock Asymptotic behaviour of tame harmonic bundles and an application to
  pure twistor {$D$}-modules. {I}.
\newblock {\em Mem. Amer. Math. Soc.}, 185(869):xii+324, 2007.

\bibitem[Moc07b]{mochizuki-D2}
T.~Mochizuki.
\newblock Asymptotic behaviour of tame harmonic bundles and an application to
  pure twistor {$D$}-modules. {II}.
\newblock {\em Mem. Amer. Math. Soc.}, 185(870):xii+565, 2007.

\bibitem[Moc09]{mochizuki-kh2}
T.~Mochizuki.
\newblock Kobayashi-{H}itchin correspondence for tame harmonic bundles. {II}.
\newblock {\em Geom. Topol.}, 13(1):359--455, 2009.

\bibitem[MS80]{mehta-seshadri}
V.B. Mehta and C.S. Seshadri.
\newblock Moduli of vector bundles on curves with parabolic structures.
\newblock {\em Math. Ann.}, 248(3):205--239, 1980.

\bibitem[PP17]{pp-n}
C.~{Pauly} and A.~{Pe{\'o}n-Nieto}.
\newblock {Very stable bundles and properness of the {H}itchin map}.
\newblock {\em ArXiv e-prints}, October 2017, 1710.10152.

\bibitem[Sab05]{claude-twistor}
C.~Sabbah.
\newblock Polarizable twistor {$\mathcal{D}$}-modules.
\newblock {\em Ast\'erisque}, 300:vi+208, 2005.

\bibitem[Ses82]{seshadri-book}
C.S. Seshadri.
\newblock {\em Fibr\'es vectoriels sur les courbes alg\'ebriques}, volume~96 of
  {\em Ast\'erisque}.
\newblock Soci\'et\'e Math\'ematique de France, Paris, 1982.
\newblock Notes written by J.-M. Drezet from a course at the {\'E}cole Normale
  Sup{\'e}rieure, June 1980.

\bibitem[Shi18]{shen}
S.~Shiyu.
\newblock Tamely ramified geometric {L}anglands correspondence in positive
  characteristic, 2018, 1810.12491.

\bibitem[Sim90]{carlos-nc}
C.~Simpson.
\newblock Harmonic bundles on noncompact curves.
\newblock {\em J. Amer. Math. Soc.}, 3(3):713--770, 1990.

\bibitem[Sim91]{carlos-naht}
C.~Simpson.
\newblock Nonabelian {H}odge theory.
\newblock In {\em Proceedings of the {I}nternational {C}ongress of
  {M}athematicians, {V}ol.\ {I}, {II} ({K}yoto, 1990)}, pages 747--756, Tokyo,
  1991. Math. Soc. Japan.

\bibitem[Sim97]{carlos-hf}
C.~Simpson.
\newblock The {H}odge filtration on nonabelian cohomology.
\newblock In {\em Algebraic geometry---{S}anta {C}ruz 1995}, volume~62 of {\em
  Proc. Sympos. Pure Math.}, pages 217--281. Amer. Math. Soc., Providence, RI,
  1997.

\bibitem[Tha96]{thaddeus}
M.~Thaddeus.
\newblock Geometric invariant theory and flips.
\newblock {\em J. Amer. Math. Soc.}, 9(3):691--723, 1996.

\bibitem[Tra16]{travkin}
R.~Travkin.
\newblock Quantum geometric {L}anglands correspondence in positive
  characteristic: the {${\rm GL}_N$} case.
\newblock {\em Duke Math. J.}, 165(7):1283--1361, 2016.

\bibitem[TW03]{teleman.woodward}
C.~Teleman and C.~Woodward.
\newblock Parabolic bundles, products of conjugacy classes and
  {G}romov-{W}itten invariants.
\newblock {\em Ann. Inst. Fourier (Grenoble)}, 53(3):713--748, 2003.

\bibitem[Yok95]{yokogawa}
K.~Yokogawa.
\newblock Infinitesimal deformation of parabolic {H}iggs sheaves.
\newblock {\em Internat. J. Math.}, 6(1):125--148, 1995.

\bibitem[Yun17]{yun-rigidity}
Z.~Yun.
\newblock Rigidity in the {L}anglands correspondence and applications.
\newblock In {\em Proceedings of the {S}ixth {I}nternational {C}ongress of
  {C}hinese {M}athematicians. {V}ol. {I}}, volume~36 of {\em Adv. Lect. Math.
  (ALM)}, pages 199--234. Int. Press, Somerville, MA, 2017.

\end{thebibliography}

\

\newpage

\printindex{notations}{List of notations}
\printindex{terms}{Index}

\newpage

 \

 \hfill

\noindent
Ron Donagi, {\sf Department of Mathematics, University of
  Pennsylvania, David Rittenhouse \linebreak
Laboratory, 209 South 33rd Street,
  Philadelphia, PA 19104, USA}, donagi@math.upenn.edu

\

\noindent
Tony Pantev, {\sf Department of Mathematics, University of
  Pennsylvania, David Rittenhouse \linebreak
Laboratory, 209 South 33rd Street,
  Philadelphia, PA 19104, USA}, tpantev@math.upenn.edu

\end{document}